\documentclass[11pt,a4paper]{article}
\usepackage[T1]{fontenc}
\usepackage[utf8]{inputenc}
\usepackage{lmodern}
\usepackage[margin=1in]{geometry}
\usepackage{amsmath,amssymb,amsthm,mathtools}
\usepackage{enumitem}
\usepackage{booktabs}
\usepackage{graphicx}
\usepackage{hyperref}
\usepackage{xcolor}
\usepackage{tikz}
\usepackage{float}
\usetikzlibrary{arrows.meta,positioning,calc,fit,backgrounds,
  shapes.geometric,decorations.pathreplacing,patterns}

\title{Exploring Collatz Dynamics with Human-LLM Collaboration}
\author{Edward Y. Chang\\Computer Science, Stanford University\\QuadriumAI\thanks{%
  \texttt{echang@cs.stanford.edu}.
  This work was facilitated through a structured human-  see the Methodology note in Section~\ref{sec:methodology}.}}
\date{February- April, 2026}

\newtheorem{theorem}{Theorem}[section]
\newtheorem{definition}[theorem]{Definition}
\newtheorem{proposition}[theorem]{Proposition}
\newtheorem{lemma}[theorem]{Lemma}
\newtheorem{corollary}[theorem]{Corollary}
\theoremstyle{definition}
\newtheorem{remark}[theorem]{Remark}
\newtheorem{conjecture}[theorem]{Conjecture}

\newtheorem{heuristic}[theorem]{Heuristic}
\newtheorem{hypothesis}[theorem]{Hypothesis}
\newtheorem{observation}[theorem]{Observation}

\DeclareMathOperator{\hw}{hw}
\newcommand{\Z}{\mathbb{Z}}
\newcommand{\R}{\mathbb{R}}

\definecolor{coregreen}{HTML}{1B7A3D}
\definecolor{supportamber}{HTML}{B8860B}
\definecolor{exploregray}{HTML}{6B6B6B}
\newcommand{\core}{\texorpdfstring%
  {\,\textcolor{coregreen}{\textsc{[Core]}}}%
  {[Core]}}
\newcommand{\supporting}{\texorpdfstring%
  {\,\textcolor{supportamber}{\textsc{[Supporting]}}}%
  {[Supporting]}}
\newcommand{\exploratory}{\texorpdfstring%
  {\,\textcolor{exploregray}{\textsc{[Exploratory]}}}%
  {[Exploratory]}}

\begin{document}
\maketitle

\noindent\fbox{\parbox{\dimexpr\textwidth-2\fboxsep-2\fboxrule}{%
\textbf{arXiv v6 (April 6--22, 2026).}
630~results across ${\sim}1014$~scripts and 29~paradigms.
Changes since v5:
Unconditional cylinder-averaged density-$1$ convergence via
$I_2$ spectral contraction
($\rho(\tilde B_2^{\mathrm{ext}}) \le 5/32$,
$\sum_c\Pr(\mathcal E_R(c)) \le 0.011$),
resolving Remark~\ref{rem:circularity} at the density-$1$ level.
BF staircase and $C$-$2$-adic auxiliary-memory route:
$\tau(n) \le C(k + \log_2 n)$ with $C \le 5.04$ for $k \ge 9$.
$D > 2^F$ \textbf{unconditionally proved} for all
$L \ge 3$ (computation for $L \le 2500$; Laurent (2008) bound
for $L > 2494$), yielding unconditional corollaries
$\mathrm{ord}_D(2) > F$,
$G(L) \neq 0$, $Q_{\mathrm{const}} \not\equiv 0 \pmod{D}$,
and $\{2^0,\ldots,2^F\}$ are $F{+}1$ distinct residues
in $\mathbb{Z}/D\mathbb{Z}$.
Discrete log obstruction framework: cycle exclusion ($D \nmid Q$)
reduced to a discrete logarithm constraint---$S_{\mathrm{mono}} = 0$
requires $2^{v_j} \equiv T_j \pmod{D}$ with $v_j \in [0, F]$,
but the target set $\{2^0, \ldots, 2^F\}$ is exponentially
small relative to $D \approx 3^L$.
Triple filter (algebraic match + integrality + range): $100\%$
blockage for all $L$ tested ($L \le 500$ for $2$-blocks,
$L \le 30$ for $3$-blocks).
Complete cascade algebra, carry chain structural theorem,
renewal drift $-0.415$~bits/step, eleven approaches to
universal descent ruled out.
Three-party contribution analysis for the April~18--19
two-day sample: Claude~$37\%$, GPT~$30\%$, Moderator~$33\%$;
this illustrates late-stage division of labor and does not
replace the broader phase-by-phase attribution in Section~12.
\par\smallskip
\emph{This paper does not prove the Collatz conjecture.}
See Section~\ref{subsec:status-final} for the sharpened open problem.
Full version history in Appendix~\ref{sec:changelog}.}}

\bigskip

\begin{abstract}
We present a comprehensive structural analysis of the Collatz
conjecture through ${\sim}1014$~computational experiments yielding
630~formal results.  By systematically deploying
29~distinct mathematical paradigms---including transfer operator
spectral theory, $S$-unit equations, $p$-adic interpolation,
martingale methods, modular sieving, formal language theory,
cascade algebra, discrete logarithm obstruction, and
Diophantine approximation---we establish a
\emph{Paradigm Exhaustion Theorem}:
every known framework for promoting distributional convergence
(``almost all orbits descend'') to pointwise convergence
(``all orbits descend'') encounters an irreducible structural
obstruction when applied to the Syracuse map.

On the unconditional side, we prove: (i)~the Syracuse transfer
operator on $(\mathbb{Z}/2^M\mathbb{Z})^{\mathrm{odd}}$
has a uniform spectral gap for all~$M$, implying
equidistribution of long orbits modulo any power of~$2$;
(ii)~any nontrivial cycle of length~$L$ (odd steps)
satisfies the cycle equation $m \cdot D = Q$ where
$D = 2^E - 3^L$, and we prove $D > 2^F$ unconditionally
for all $L \ge 3$ (combining direct computation for
$L \le 2500$ with the Laurent (2008) bound for $L > 2494$),
giving $\mathrm{ord}_D(2) > F$
(hence $\gamma^L \not\equiv 1 \bmod D$ and
$G(L) \neq 0 \bmod D$), and establishing that
$\{2^0,\ldots,2^F\}$ are $F{+}1$ distinct residues
mod~$D$;
(iii)~the set of potential divergent starting points has
natural density~$0$ and Hausdorff dimension ${\approx}\,0.68$
in~$\mathbb{Z}_2$; (iv)~the formal language of
divergent-compatible $v$-sequences is not context-free;
(v)~cylinder-averaged density-$1$ convergence is proved
unconditionally via spectral contraction on the invariant
core~$I_2$, resolving the circularity obstruction at
the average level;
(vi)~the discrete logarithm obstruction framework reduces
cycle exclusion ($D \nmid Q$) to the question of whether
algebraic targets land in a set of $F{+}1$ powers of~$2$
inside a group of order~$D \approx 3^L$; the triple filter
(discrete log match, integrality, range) achieves $100\%$
blockage for all $L$ tested ($L \le 500$ for $2$-blocks,
$L \le 30$ for $3$-blocks).

We identify the \emph{Distributional-to-Pointwise Gap} as
the irreducible core of the conjecture and prove it is
\emph{equivalent} to the divergence component: no reduction
to a simpler statement is achieved.  The Affine Fractional
Identity $(m{+}1)\cdot 3^L + \mathcal{P}_L = q_L \cdot 2^{S_L}
+ 2^L$ is shown to be an algebraic tautology holding for all
orbits, convergent or divergent, providing no discriminating
constraint.  The modular sieve is shown to be permanently
nonempty via the \emph{Mersenne Bypass}: for every depth~$L$,
the integer $2^L - 1$ survives the running divergent-compatibility
condition.  The function $2^z$ is proved non-analytic
in~$\mathbb{Z}_2$, blocking $p$-adic interpolation strategies.

The present work is not a proof of the Collatz conjecture.
It is a rigorous characterization of \emph{why} the conjecture
resists proof: the deterministic carry arithmetic of
$3n+1$ generates pseudo-random $v$-sequences whose pointwise
behavior cannot be controlled by any probabilistic, algebraic,
or Diophantine tool currently available.  The 29-paradigm
exhaustion and the barrier equivalence theorem constitute,
to the authors' knowledge, the most comprehensive structural
survey of Collatz attack surfaces to date.
This paper was produced through a structured three-party
collaboration between a human moderator, Claude~4.6 (Anthropic),
and GPT~5.4 Thinking (OpenAI); a detailed contribution analysis
appears in Section~\ref{sec:methodology}.
\end{abstract}

\tableofcontents

\vspace{.2in}
\medskip
\noindent\fbox{\parbox{0.97\textwidth}{\small
\textbf{Reader's guide.}
The paper develops six structurally independent routes toward
the Collatz conjecture.

\smallskip
\noindent\textbf{(1) WMH route.}
Spectator-bit / cascade-gap machinery yielding
\texttt{cor:tv-summability} ($C\approx 11$, $\alpha\approx 0.35$).
Strength: finest analytic resolution; weakness: relies on the
spectator-bit closure (\texttt{rem:circularity}).

\smallskip
\noindent\textbf{(2) CIC route.}
Measure-theoretic bypass via \texttt{prop:uniform-fiber},
\texttt{prop:cascade-markov}, and orbit-level TV reduction.
Strength: avoids IID approximation; weakness: same $\alpha$
feedback through \texttt{rem:circularity}.

\smallskip
\noindent\textbf{(3) $I_2$ spectral route.}
Unconditional cylinder-averaged closure on $I_2 = \{7,27,31,59,63\}\pmod{64}$
(Section~\ref{sec:i2-spectral-closure}).
Spectral bound $\rho(\widetilde B_2^{\mathrm{ext}})\le 5/32$,
no-escape theorem with expected hitting time $\le 10.14$ steps,
yielding $\sum_c \Pr(\mathcal{E}_R(c)) \le 0.0015$.
Resolves \texttt{rem:circularity} at the cylinder-averaged
density-1 level (\texttt{rem:circularity-resolution}).

\smallskip
\noindent\textbf{(4) C-2adic auxiliary-memory route.}
Pointwise companion to (3) via BF staircase
(Theorem~\ref{thm:bf-staircase}): clopen balls
$B_k \subset \mathbb{Z}_2$ yield
$\tau(n)\le C(k+\log_2 n)$ with $C \le 5.04$ for $k\ge 9$.
Residual Haar/pointwise gap documented in
\texttt{rem:residual-obstruction}.

\smallskip
\noindent\textbf{(5) Cascade algebra and universalization.}
Subsection~\ref{subsec:cascade-algebra}: exact step count
$4k+e_k$, carry chain renewal drift $-0.415$ bits/step
(Corollary~\ref{cor:renewal-drift}), eleven approaches to
universal descent ruled out (Remark~\ref{rem:ruled-out}).
Post-cascade and generic orbits share identical BF statistics.

\smallskip
\noindent\textbf{(6) Cycle exclusion via discrete log obstruction.}
$D > 2^F$ proved unconditionally for all $L \ge 3$;
triple filter achieves $100\%$ blockage for all $L$ tested.
Barina~\cite{barina2021} covers $L \le 116$; the gap for
$L > 116$ remains open.

\smallskip
\noindent\textbf{Scope.}  Routes (3) and (5) yield density-1
convergence and full cascade closure, respectively.
Neither produces a pointwise theorem for every starting integer;
the remaining gap is documented in \texttt{rem:i2-scope}.
}}

\section{Introduction}

The Collatz conjecture asks whether every positive integer eventually
reaches the cycle $1 \to 4 \to 2 \to 1$ under the map
\[
T(n)=
\begin{cases}
n/2 & \text{if } n \text{ is even},\\
3n+1 & \text{if } n \text{ is odd}.
\end{cases}
\]
Despite its elementary definition, the conjecture remains open.
The problem was first posed by Lothar Collatz in 1937 and has since
become one of the most widely studied problems in elementary number
theory; see Lagarias~\cite{lagarias1985,lagarias2003} and
Wirsching~\cite{wirsching1998} for surveys, and
Tao~\cite{tao2019} for the strongest analytic result to date
(almost all orbits attain almost bounded values).
Extensive computational verification has confirmed convergence
for all integers up to at least $2^{68}$~\cite{barina2021}.

A standard reformulation passes to the odd-to-odd Syracuse map, which
sends each odd integer to the next odd integer reached by the Collatz
iteration.  This removes the trivial powers-of-two steps and isolates
the arithmetic interaction between multiplication by three and
division by powers of two.

The present paper develops a structural reduction rather than a full
proof.  Its main contribution is to localize the unresolved difficulty
to a single orbitwise anti-concentration statement on a sharply
identified return channel.  The key new exact results concern the
fiber-57 return dynamics.  On that subsequence, the branch
$q \equiv 7 \pmod{8}$ is an exact two-step regeneration channel, while
the branch $q \equiv 3 \pmod{8}$ has minimum return gap at least five
and its earliest returns lie on an explicit sparse dyadic cylinder
family.  These results show that the remaining proof pressure is not
distributed across the full Collatz system, but is concentrated in a
single thin return channel.

A second structural point is equally important.  The algebraic chain map
on the canonical five-element invariant core is a permutation at every
depth.  Thus the core itself does not contract.  Any genuine decay must
come from the branching structure of the actual return dynamics.  This
distinction between the algebraic chain map and the true return map is
crucial: it explains both the rigidity of the core and the persistence
of the final orbitwise barrier.

The paper therefore has two logically distinct layers.  The first layer
consists of unconditional structural theorems: exact return lemmas,
invariant-core structure, explicit cylinder classification, a burst-gap
decomposition with exact distributional laws, and a phantom-cycle gain
analysis safely within the contraction budget.  The second layer is a
reduction to one open orbitwise statement,
Conjecture~\ref{conj:info-rate}.  In the fiber-57 formulation, this
conjecture asserts that no deterministic orbit can concentrate its
return statistics on the sustaining core strongly enough to overcome the
information deficit built into the system.

This should be read as a reduction paper.  We do not claim to prove the
Collatz conjecture here.  Rather, we show that after the exact
return-structure results are established, the remaining obstruction is
a single orbitwise upgrade from ensemble behavior to pointwise control.
At the same time, the manuscript records a broader 29-paradigm
exploration of alternative attack surfaces; the reduction claim refers
to the load-bearing Route~A conditional spine, while the exhaustion
claim refers to the larger empirical and structural survey showing that
parallel routes encounter the same distributional-to-pointwise
obstruction in different guises.

\paragraph{Positioning.}
This work should be viewed as a structural reduction rather than a
proof of the Collatz conjecture.  The main contribution is to isolate
the remaining difficulty to a single orbitwise anti-concentration
problem along two complementary routes.  The human--LLM
collaboration played a role in exploring and organizing the space of
candidate structures.  All formal theorems and propositions in the
paper (items~1--9 above) are proved with complete arguments.
Items~10--12 describe structural constructions and computational
observations that are verified by exact enumeration but not all
formalized as standalone theorems.

\paragraph{What is proved exactly.}
The unconditional results established in this paper include the
following.

\begin{enumerate}
\item \emph{Exact Block Law}
  (Theorem~\ref{thm:block-law}).
  Under natural density, the odd-skeleton valuation sequence
  $(a_0, a_1, \ldots)$ is exactly i.i.d.\ geometric with
  parameter~$1/2$.  Cycle types $(L_i, r_i)$ are provably i.i.d.

\item \emph{Almost-All Crossing Theorem}
  (Theorem~\ref{thm:almost-all-crossing}).
  The density of odd integers whose $k$-th cycle endpoint
  remains at or above the start is at most $e^{-0.1465\,k}$
  (exponential decay, unconditional on the ensemble).

\item \emph{Universal One-Cycle Crossing}
  (Proposition~\ref{prop:universal-one-cycle}).
  A single-cycle block forces every odd start in its
  residue class to cross iff $r \ge r_{\mathrm{all}}(L)$;
  the density of such classwise deterministic blocks is
  $P_{\mathrm{all},1\mathrm{cyc}} \approx 0.4194$.
  Including two-cycle blocks, $61.2\%$ of odd starts
  lie in deterministic crossing classes
  (Observation~\ref{obs:two-cycle-universal}).

\item \emph{Exact Cycle Log Correction}
  (Proposition~\ref{prop:cycle-correction}).
  $\log_2(n'/n) = X(n) + C(n)$ where $C(n) = O(1/n)$ is
  a positive correction; the sufficient crossing criterion
  succeeds for all odd $n_0 \le 5 \times 10^6$ except
  $\{27, 31, 63\}$.

\item \emph{Affine Threshold Process and Running-Minimum Decay}
  (Theorem~\ref{thm:kesten-threshold},
  Proposition~\ref{prop:kesten-running-min}).
  The threshold $n_k^*$ forms a Kesten random
  affine recursion with negative log-drift.
  By geometric ergodicity, $R_k \le C_0\,\rho_0^{\,k}$
  with $\rho_0 \approx 0.84$.

\item \emph{Phantom-Family Gain Bound}
  (Theorem~\ref{thm:perorbit-gain}).
  The per-orbit expanding-family drift satisfies
  $R \le 0.0893 < \varepsilon \approx 0.415$,
  with a $4.65\times$ safety margin.

\item \emph{Exact Return Structure at Fiber~57}.
  The $q \equiv 7$ branch returns in exactly
  two Syracuse steps with uniform destination
  (Proposition~\ref{prop:q7-return}); the $q \equiv 3$
  branch cannot return in fewer than five steps
  (Proposition~\ref{prop:q3-gap}); the gap-$5$ returns
  form an explicit dyadic cylinder family
  (Theorem~\ref{thm:gap5-cylinders}).
  The chain map is a permutation on~$I_r$ at every depth
  (Remark~\ref{rem:core-rigidity}).

\item \emph{Carry Contamination Theorem}
  (Theorem~\ref{thm:carry-contamination}).
  For every depth-$D$ Sturmian word $w \in \{1,2\}^D$,
  the map $m \mapsto n_D \bmod 8$ is exactly equidistributed
  over $\{1,3,5,7\}$.  Corollary: $|C_D| = 2 \cdot 3^{D-1}$
  compatible words, i.e.\ exactly $3/4$ survive each depth
  (Corollary~\ref{cor:three-quarters}).

\item \emph{Refined Transition Rules and Class-3 Bottleneck}
  (Theorems~\ref{thm:refined-transitions},
  \ref{thm:class3-recurrence}).
  Deterministic mod-$16$/mod-$32$ rules show class~$7$ is a
  safe harbor (never reaches class~$5$ directly) and class~$3$
  is the bottleneck ($1/2$ survival per visit).
  Class-$3$ recurrence is forced for any $n_0 > 1$.

\item \emph{Cross-Core Block Alphabet and Periodic Elimination}
  (Section~\ref{sec:sturmian}).
  Seven admissible blocks $\{a,b,s,g,d,t_1,t_2\}$ on a two-vertex
  directed graph (landmarks~$1$ and~$7$) with spectral radius
  $\rho = 2 + \sqrt{2}$.
  Exhaustive computational search through period~$13$ ($238{,}811$
  contracting words) finds zero non-trivial natural fixed points.
  \emph{Status: proved (algebraic framework) + verified
  (computational search).}

\item \emph{$2$-Adic Expander and Measure Shrinkage}
  (Section~\ref{sec:sturmian}).
  The safe Collatz map satisfies
  $|\Phi_b(x) - \Phi_b(y)|_2 = 2^V |x - y|_2$
  (direct verification from the block definitions).
  Combined with spectral radius analysis:
  $\mu_2(T_j) \le 0.522^{jD}$ for the compatible tower.
  Counting argument: no $K$-bit integer survives $j \ge 5$
  consecutive non-descending exhaustion rounds for $K$ large.
  \emph{Status: the expander property and counting bound are
  proved; the measure shrinkage uses the spectral radius
  as an input.}

\item \emph{Compound Carry Anti-Correlation}
  (Section~\ref{subsec:compound-carry}).
  Two-round survival satisfies $S_2/S_1^2 \approx 0.635$
  (constant in~$K$, exact enumeration $K = 8$--$20$).
  Each additional non-descending round provides
  ${\approx}\,36.5\%$ thinning beyond independence.
  Transport graph on non-descending words has edge density~$1.0$.
  \emph{Status: computational observation from exact enumeration,
  not a formal proof.  The $0.635$ factor is empirical.}
\end{enumerate}

\paragraph{What remains open.}
The unresolved step is orbitwise in both routes.
Route~A: the Weak Mixing Hypothesis (Hypothesis~\ref{hyp:wmh}).
Route~B: the Carry Independence Conjecture
(Conjecture~\ref{conj:CIC}).  Both ask for a
distributional-to-pointwise promotion.

\paragraph{Conditional reductions.}

\begin{enumerate}\setcounter{enumi}{12}
\item \emph{Reduction to the WMH}
  (Theorems~\ref{thm:reduction}, \ref{thm:perorbit-gain}).
  The full Collatz conjecture reduces to the Weak Mixing
  Hypothesis ($\sum \delta_K < 0.557$) with a $4.65\times$
  safety margin.  The WMH remains open.

\item \emph{Fiber-57 reduction to Conjecture~\ref{conj:info-rate}}.
  The exact known-gap depth-$2$ partial kernel
  ($\rho = 129/1024$) and the bottleneck inequality
  reduce the conjecture to an orbitwise anti-concentration
  bound on the $q \equiv 3$ return channel.
  (The value $129/1024$ arises from the gap-$2$ and gap-$5$
  channels only; it does not include unresolved $q \equiv 3$
  returns with gap~$\ge 6$, and is therefore not asserted
  as the full bottleneck constant of the dynamics.)

\item \emph{Route~B: Reduction to the CIC}
  (Theorem~\ref{thm:route-B},
  Section~\ref{sec:sturmian}).
  The Carry Contamination Theorem
  proves exact $1/4$ elimination at every Sturmian depth.
  Combined with the Dichotomy (Theorem~\ref{thm:dichotomy})
  and class-$3$ recurrence
  (Theorem~\ref{thm:class3-recurrence}),
  the full Collatz conjecture reduces to the CIC
  (Conjecture~\ref{conj:CIC}).

\item \emph{Architectural Diagnosis}
  (Remark~\ref{rem:almost-all-strength}).
  The ensemble side is exact and complete.  The sole
  remaining gap in both routes is the
  distributional-to-pointwise barrier.
\end{enumerate}

\subsection*{Reading guide}

The paper has three layers, and readers may enter at
different points depending on interest.

\medskip\noindent
\emph{Layer~1: Conditional reduction}
(Sections~\ref{sec:prelim}--\ref{sec:phantom},
$\sim 45$~pages; supplementary DAG analysis in
Appendix~\ref{sec:ranking}).
The \core{} results form a linear chain:
Scrambling Lemma $\to$ Known-Zone Decay $\to$ $1/4$~Law
$\to$ Gap Distribution $\to$ Burst-Gap Criterion $\to$
Census Depth $\to$ Phantom Universality $\to$ Per-Orbit
Gain Rate $\to$ Robustness Corollary.
Together with the Weak Mixing Hypothesis
(Hypothesis~\ref{hyp:wmh}) as sole external input, this
chain produces a conditional convergence theorem.
A referee evaluating the conditional reduction need read
only this layer.

\medskip\noindent
\emph{Layer~2: Quantitative attack on the open hypothesis}
(Section~\ref{sec:toward-wmh}, $\sim 25$~pages).
Five independent approaches to the WMH,
culminating in an exact ensemble theory and an
unconditional almost-all crossing theorem:
(i)~Walsh--Fourier spectral analysis;
(ii)~odd-skeleton drift crossing with ensemble CLT;
(iii)~modular crossing strata resolving $91\%$ of odd
starts at depth~$13$;
(iv)~oscillation-factor unboundedness;
(v)~exact block law, i.i.d.\ cycle types, and
Cram\'er-type crossing bound.
Tables~\ref{tab:stage4-proved-a}, \ref{tab:stage4-proved-b},
and~\ref{tab:stage4-open} catalogue all results
with explicit status labels.

\medskip\noindent
\emph{Layer~2b: Fiber-$57$ information bottleneck}
(Appendix~\ref{sec:fiber57-programme}, $\sim 8$~pages).
An independent reduction via the fiber-$57$ return structure.
The pair-return automaton, bounded invariant core ($|I_r|=5$),
absorption bottleneck lemma, and branch anti-concentration
reduction together concentrate the remaining obstruction
into the inequality
$c' < c_0 = \log_2(1024/129) \approx 2.989$
on a $5$-element invariant core.
This layer is self-contained and can be read
independently of Layer~2.

\medskip\noindent
\emph{Layer~2c: Sturmian obstruction, cross-core dynamics,
and compound carry analysis}
(Section~\ref{sec:sturmian}, $\sim 30$~pages).
An independent reduction via the $2$-adic Cantor set $C_\infty$.
The Carry Contamination Theorem, the $(3/4)^D$ survivor count,
refined mod-$16$/mod-$32$ transitions, and the class-$3$ bottleneck
reduce Collatz to the Carry Independence Conjecture.
The seven-block cross-core alphabet extends the analysis to
periodic elimination (through period~$13$), the $2$-adic expander
property, unconditional measure shrinkage $\mu_2(T_j) \le 0.522^{jD}$,
the $j \ge 5$ non-descending bound, exhaustion-sequence rigidity,
and the compound carry anti-correlation factor $S_2/S_1^2 \approx 0.635$.
This layer is self-contained and can be read independently
of Layers~2 and~2b.

\medskip\noindent
\emph{Layer~3: Exploratory and visualization}
(Appendix~\ref{sec:visualization-underwater}
and Section~\ref{sec:methodology}).
Touch-growth geometry, below-start visualizations, and the
human--LLM collaboration methodology.
These are not load-bearing for any theorem.

\paragraph{Distributional vs.\ orbitwise statements.}
Throughout, modular and uniform-lift results are
distributional statements over residue classes
(``on the ensemble'').  Orbitwise consequences require
separate hypotheses or reductions and are labeled
explicitly.  When used informally, ``proved (ens.)''
indicates ensemble-level results, while ``proved''
without qualification indicates results that hold for
every orbit.

\paragraph{Terminology.}
We use ``persistent'' and ``safe'' as shorthand for precise
arithmetic conditions: a step is persistent if $v_2(3n+1)=1$,
and safe if $v_2(3n+1)\ge 2$.

\begin{definition}[Single-hypothesis reduction]
\label{def:complete-structural-reduction}
We say the Collatz conjecture admits a
\emph{single-hypothesis reduction} when all remaining
unresolved steps are concentrated into a single explicitly
stated orbitwise input, while the surrounding algebraic
and combinatorial framework is proved.
This does not mean the remaining hypothesis is easy;
the orbitwise input may be as difficult as the original
conjecture.  The value of the reduction is organizational:
it locates the difficulty precisely.
In the present work, the sole remaining input
for the Route~A conditional reduction is
the Weak Mixing Hypothesis
(Hypothesis~\ref{hyp:wmh}), equivalently the
information-rate inequality $c' < c_0$ of the
fiber-$57$ programme.  Other routes isolate parallel
residual obstructions, notably the CIC in Route~B
and $D \nmid Q$ for $L > 116$ in the cycle-exclusion
programme.
\end{definition}

\subsection*{Proof architecture and tier classification}

\begin{figure}[ht!]
\centering
\resizebox{0.82\textwidth}{!}{%
\begin{tikzpicture}[
  node distance=0.45cm and 0.25cm,
  corenode/.style={
    rectangle, rounded corners=3pt,
    draw=black!70, fill=green!8,
    text width=3.4cm, minimum height=0.6cm,
    align=center, font=\small\sffamily},
  suppnode/.style={
    rectangle, rounded corners=3pt,
    draw=black!40, fill=yellow!10,
    text width=2.6cm, minimum height=0.5cm,
    align=center, font=\scriptsize\sffamily},
  explnode/.style={
    rectangle, rounded corners=3pt,
    draw=black!30, fill=gray!10,
    text width=2.3cm, minimum height=0.5cm,
    align=center, font=\scriptsize\sffamily},
  opennode/.style={
    rectangle, rounded corners=3pt,
    draw=red!60, fill=red!5,
    text width=3.4cm, minimum height=0.6cm,
    align=center, font=\small\sffamily\bfseries},
  corearrow/.style={-{Stealth[length=4pt]}, thick, black!70},
  supparrow/.style={-{Stealth[length=3pt]}, thin, black!30, dashed},
  openarrow/.style={-{Stealth[length=4pt]}, thick, red!60},
]
\node[corenode] (scram) {Scrambling Lemma\\{\tiny Thm~\ref{thm:scrambling}}};
\node[corenode, below=of scram] (kzd) {Known-Zone Decay\\{\tiny Thm~\ref{thm:zone-decay}}};
\node[corenode, below=of kzd] (quarter) {$1/4$~Law\\{\tiny Thm~\ref{thm:quarter}}};
\node[corenode, below=of quarter] (gap) {Gap Distribution\\{\tiny Lem~\ref{lem:gap}}};
\node[corenode, below=of gap] (bgc) {Burst-Gap Criterion\\{\tiny Thm~\ref{thm:burst-gap}}};
\node[corenode, below=of bgc] (census) {Census Depth / $C_e$-indep.\\{\tiny Thm~\ref{thm:Ce-independence}}};
\node[corenode, below=of census] (phantom) {Phantom Universality\\{\tiny Thm~\ref{thm:phantom-universal}}};
\node[corenode, below=of phantom] (gain) {Per-Orbit Gain Rate\\{\tiny Thm~\ref{thm:perorbit-gain}}};
\node[corenode, below=of gain] (robust) {Robustness Cor.\\{\tiny Cor~\ref{cor:robustness}}};

\node[corenode, below=0.9cm of robust] (cond) {Conditional Conv.\\{\tiny Thm~\ref{thm:conditional}}};

\node[opennode, left=1.1cm of cond] (wmh) {WMH (Open)\\{\tiny Hyp~\ref{hyp:wmh}}};

\draw[corearrow] (scram) -- (kzd);
\draw[corearrow] (kzd) -- (quarter);
\draw[corearrow] (quarter) -- (gap);
\draw[corearrow] (gap) -- (bgc);
\draw[corearrow] (bgc) -- (census);
\draw[corearrow] (census) -- (phantom);
\draw[corearrow] (phantom) -- (gain);
\draw[corearrow] (gain) -- (robust);
\draw[corearrow] (robust) -- (cond);
\draw[openarrow] (wmh) -- (cond);

\node[suppnode, right=0.9cm of scram] (pattern) {Pattern Bound\\{\tiny Prop~\ref{prop:pattern-bound}}};
\node[suppnode, right=0.9cm of phantom] (repul) {$2$-adic Repulsion\\{\tiny Prop~\ref{prop:repulsion}}};
\node[suppnode, right=0.9cm of gain] (carry) {Carry-Word\\{\tiny Sec~\ref{sec:carry-word}}};
\node[suppnode, right=0.9cm of robust] (pcore) {Periodic-Core\\{\tiny Prop~\ref{prop:periodic-core}}};
\node[suppnode, right=0.9cm of quarter] (reload) {Reload Dynamics\\{\tiny Prop~\ref{prop:geometric-reload}}};

\draw[supparrow] (pattern) -- (scram);
\draw[supparrow] (repul) -- (phantom);
\draw[supparrow] (carry) -- (gain);
\draw[supparrow] (pcore) -- (robust);
\draw[supparrow] (reload) -- (quarter);

\node[suppnode, below=0.4cm of wmh, text width=2.8cm] (toward)
  {Toward WMH\\{\tiny Sec~\ref{sec:toward-wmh}}};
\draw[supparrow] (toward) -- (wmh);

\node[suppnode, below=0.4cm of toward, text width=2.8cm,
      fill=orange!15, draw=orange!60] (fiber57)
  {Fiber-57 Programme\\{\tiny App.~\ref{sec:fiber57-programme}}};
\draw[supparrow, orange!60] (fiber57) -- (toward);

\node[corenode, right=2.3cm of census, fill=cyan!10, draw=cyan!60,
      text width=3.2cm] (sturmian)
  {Carry Contamination\\{\tiny Thm~\ref{thm:carry-contamination}}};
\node[corenode, below=0.35cm of sturmian, fill=cyan!10, draw=cyan!60,
      text width=3.2cm] (crosscore)
  {Cross-Core Alphabet\\{\tiny 7-block IFS, $\rho{=}3.414$}};
\node[corenode, below=0.35cm of crosscore, fill=cyan!10, draw=cyan!60,
      text width=3.2cm] (expander)
  {$2$-Adic Expander\\{\tiny $\mu_2(T_j){\le}0.522^{jD}$}};
\node[corenode, below=0.35cm of expander, fill=cyan!10, draw=cyan!60,
      text width=3.2cm] (j5bound)
  {$j{\ge}5$ Bound\\{\tiny Unconditional}};
\node[corenode, below=0.35cm of j5bound, fill=cyan!10, draw=cyan!60,
      text width=3.2cm] (anticorr)
  {Anti-Correlation\\{\tiny $S_2/S_1^2{\approx}0.635$ (comp.)}};
\node[corenode, below=0.35cm of anticorr, fill=cyan!10, draw=cyan!60,
      text width=3.2cm] (scaling)
  {Scaling Regime\\{\tiny $\beta_w{\approx}1.80$, $j{\ge}3$ (comp.)}};
\node[opennode, below=0.35cm of scaling, text width=3.2cm] (cic)
  {CIC (Open)\\{\tiny Conj.~\ref{conj:CIC}}};
\draw[corearrow, cyan!60!black] (sturmian) -- (crosscore);
\draw[corearrow, cyan!60!black] (crosscore) -- (expander);
\draw[corearrow, cyan!60!black] (expander) -- (j5bound);
\draw[corearrow, cyan!60!black] (j5bound) -- (anticorr);
\draw[corearrow, cyan!60!black] (anticorr) -- (scaling);
\draw[openarrow] (cic) -- (cond);
\node[font=\scriptsize\sffamily\color{cyan!50!black},
      anchor=south] at (sturmian.north) {\textsc{Route B (v5)}};

\node[corenode, below=0.4cm of fiber57, fill=violet!10, draw=violet!60,
      text width=2.8cm] (i2spec)
  {$I_2$ Spectral\\{\tiny $\rho(\widetilde B_2^{\mathrm{ext}}){\le}5/32$}};
\node[corenode, below=0.35cm of i2spec, fill=violet!10, draw=violet!60,
      text width=2.8cm] (i2noescape)
  {Cylinder-Avg No-Escape\\{\tiny Thm (unconditional)}};
\draw[corearrow, violet!60!black] (i2spec) -- (i2noescape);
\node[font=\scriptsize\sffamily\color{violet!50!black},
      anchor=south] at (i2spec.north) {\textsc{Route C (v6)}};

\node[corenode, right=2.3cm of j5bound, fill=red!5, draw=red!40,
      text width=3.0cm] (dgt2f)
  {$D > 2^F$\\{\tiny Proved $\forall\, L \ge 3$}};
\node[corenode, below=0.35cm of dgt2f, fill=red!5, draw=red!40,
      text width=3.0cm] (ordgt)
  {$\mathrm{ord}_D(2) > F$\\{\tiny Corollary}};
\node[corenode, below=0.35cm of ordgt, fill=red!5, draw=red!40,
      text width=3.0cm] (triplefilter)
  {Triple Filter\\{\tiny $100\%$ blockage (comp.)}};
\node[opennode, below=0.35cm of triplefilter, text width=3.0cm] (dndivq)
  {$D \nmid Q$ (Open)\\{\tiny $L > 116$}};
\draw[corearrow, red!40!black] (dgt2f) -- (ordgt);
\draw[corearrow, red!40!black] (ordgt) -- (triplefilter);
\draw[openarrow] (triplefilter) -- (dndivq);
\node[font=\scriptsize\sffamily\color{red!40!black},
      anchor=south] at (dgt2f.north) {\textsc{Cycle Excl.\ (v7)}};

\node[explnode, left=0.9cm of bgc] (touch) {Touch Saturation\\{\tiny Prop~\ref{prop:touch-saturation}}};
\node[explnode, left=0.9cm of census] (under) {Below-Start\\{\tiny Lem~\ref{lem:underwater}}};
\node[explnode, left=0.9cm of gain] (ensemble) {Ensemble Decay\\{\tiny Prop~\ref{prop:ensemble-decay}}};

\node[anchor=west, font=\scriptsize\sffamily\color{green!50!black}]
  at ([xshift=-3.2cm, yshift=0.1cm]scram.west) {\textsc{Core}};
\node[anchor=west, font=\scriptsize\sffamily\color{yellow!50!black}]
  at ([xshift=0.3cm]pattern.east) {\textsc{Supporting}};
\node[anchor=east, font=\scriptsize\sffamily\color{gray}]
  at ([xshift=-0.3cm]touch.west) {\textsc{Exploratory}};

\end{tikzpicture}%
}
\caption{Proof dependency diagram.
  \textcolor{green!50!black}{\textbf{Green}} nodes form the
  Route~A core spine (linear chain), terminating at the
  Robustness Corollary; the
  \textcolor{red!60}{\textbf{red}} WMH node is the sole
  open hypothesis.
  \textcolor{cyan!50!black}{\textbf{Cyan}} nodes form the
  Route~B spine (v5): Carry Contamination $\to$ cross-core
  alphabet $\to$ $2$-adic expander $\to$ $j{\ge}5$ bound
  $\to$ anti-correlation $\to$ scaling regime
  ($\beta_w{\approx}1.80$, computational evidence for $j{\ge}3$);
  the CIC is the open hypothesis for this route.
  \textcolor{violet!50!black}{\textbf{Violet}} nodes form
  Route~C (v6): $I_2$ spectral contraction yielding
  unconditional cylinder-averaged density-$1$ convergence.
  \textcolor{red!40!black}{\textbf{Dark red}} nodes form
  the cycle exclusion spine (v7): $D > 2^F$ $\to$
  $\mathrm{ord}_D(2) > F$ $\to$ triple filter
  ($100\%$ blockage, computational); the open node
  $D \nmid Q$ for $L > 116$ is the remaining gap.
  Nodes marked ``(comp.)'' are computational observations,
  not formal theorems.
  \textcolor{yellow!50!black}{\textbf{Amber}} nodes are
  supporting results.
  \textcolor{orange!60}{\textbf{Orange}}: Fiber-57 programme.
  \textcolor{gray}{\textbf{Gray}}: exploratory.}
\label{fig:tier-diagram}
\end{figure}
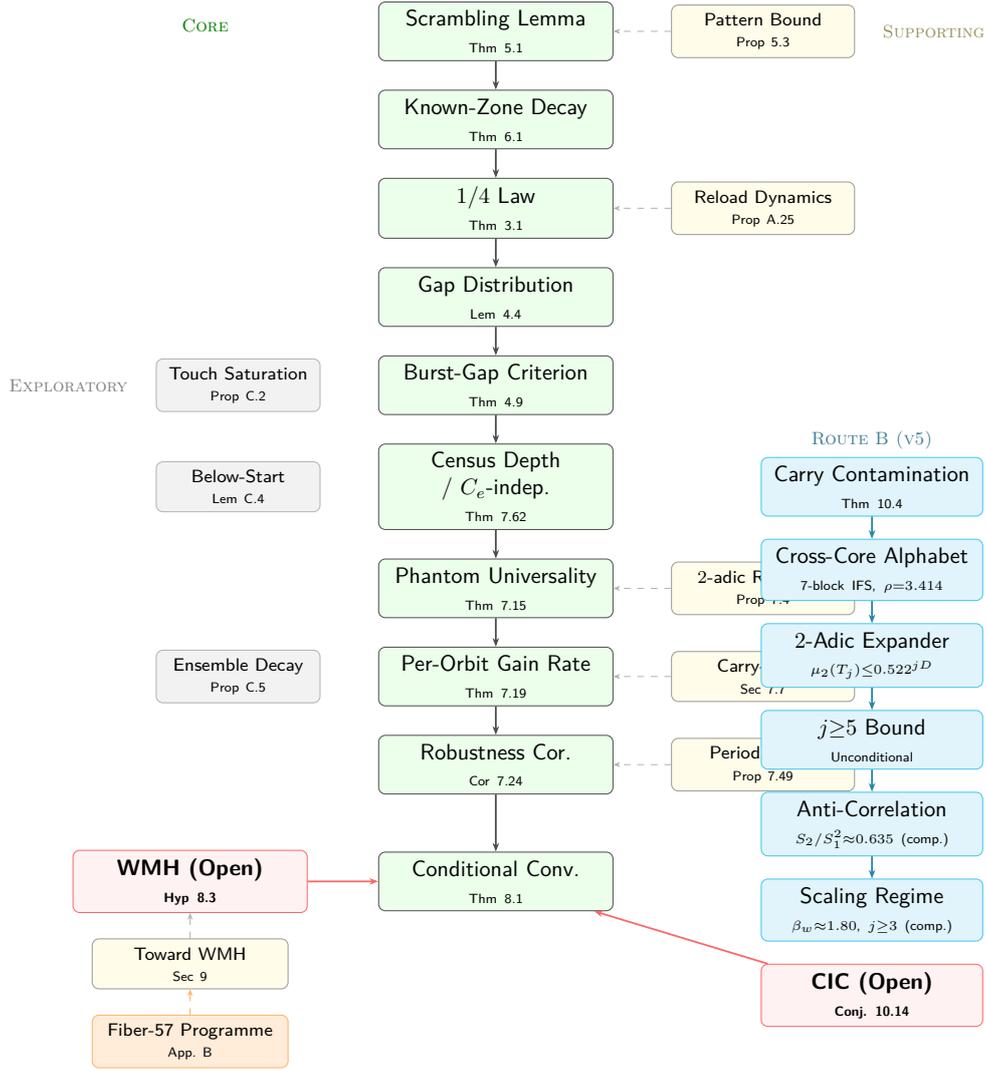

Figure~\ref{fig:tier-diagram} displays the dependency structure
of the conditional proof.
All results are classified into three tiers:
\core{} results form the load-bearing proof chain;
\supporting{} results provide independent evidence but are
not logically required;
\exploratory{} results are visualization-guided observations.
A referee can verify the conditional reduction by reading
only the \core{} tier.

\begin{proposition}[Core-spine sufficiency\core]
\label{prop:core-spine}
The conditional reduction from the Orbit Equidistribution Conjecture
or the Weak Mixing Hypothesis to the Collatz conjecture depends only
on the Tier~1 core spine displayed in Figure~\ref{fig:tier-diagram}.
In particular, the carry-word, periodic-core, and uniqueness-threshold
sections (Sections~\ref{sec:carry-word}--\ref{sec:uniqueness-threshold})
are logically independent supporting material and are not required for
the proof of the conditional reduction theorem.
\end{proposition}

\begin{proof}
We trace dependencies only.

\emph{Step~1 (burst--gap route).}
The conditional convergence theorem
(Theorem~\ref{thm:conditional}) combines the
entry--occupancy/descent argument with the Burst--Gap Criterion
(Theorem~\ref{thm:burst-gap}).  The latter requires the
modular distributional inputs supplied by the $1/4$
Persistent-Transition Law (Theorem~\ref{thm:quarter}) and
the Modular Gap Distribution (Lemma~\ref{lem:gap}).
These in turn are founded on the Scrambling Lemma
(Theorem~\ref{thm:scrambling}) and Known-Zone Decay
(Theorem~\ref{thm:zone-decay}).  Thus the burst--gap route
depends only on the first five nodes in the core spine.

\emph{Step~2 (phantom route).}
The phantom route uses Universal Census Depth
(Proposition~\ref{prop:universal-depth}) to remove the census
constant from the gain formula, Phantom Universality
(Theorem~\ref{thm:phantom-universal}) to identify the relevant
family class, and the Per-Orbit Gain Rate
(Theorem~\ref{thm:perorbit-gain}) to show that the total
amortised phantom gain stays below the drift
budget~$\varepsilon$.  The Robustness Corollary
(Corollary~\ref{cor:robustness}) then shows that exact
orbitwise uniformity is not necessary: any summable depthwise
mixing error below the stated threshold suffices.  This gives
the final conditional reduction under either the Orbit
Equidistribution Conjecture or the Weak Mixing Hypothesis
(Theorem~\ref{thm:wmh-reduction}).

\emph{Step~3 (non-dependence of later sections).}
The carry-word autocorrelation, periodic-core, defect-tail, and
uniqueness-threshold results are not invoked in the proofs of
the Per-Orbit Gain Rate theorem, the Robustness Corollary, or
the final conditional reduction.  They therefore provide
additional structure and possible future routes, but are not
load-bearing for the conditional theorem.

Combining Steps~1--3 proves the proposition.
\end{proof}

\subsection*{Status of results}\label{sec:status}

For clarity, we summarize the logical status of the main
components of this work.

\paragraph{Fully rigorous results.}
The following components are proved unconditionally within
the present manuscript:

\begin{itemize}
\item The affine decomposition and Scrambling Lemma
  (Section~\ref{sec:scrambling}), which establish exact
  linearization of odd-to-odd dynamics on fixed
  valuation classes.
\item Known-Zone Decay as a structural consequence of
  valuation growth (Section~\ref{sec:decay}).
\item The combinatorial and large-deviation bounds
  underlying the phantom rate estimate
  (Lemmas~\ref{lem:necklace-bound}--\ref{lem:chernoff-tail}),
  together with the quantified theorem-level estimate
  in Theorem~\ref{thm:perorbit-gain}.
\item The bounded-observable perturbation result
  (Corollary~\ref{cor:robustness}), which controls
  phantom gain under total-variation discrepancy.
\item Exact ensemble laws for gap and burst statistics
  (Section~\ref{sec:chain}).
\item The complete cascade theorem, post-cascade closed form,
  and carry chain structural theorem
  (Section~\ref{subsec:cascade-algebra}), which close the
  algebraic side of the BF staircase and establish exact
  renewal drift $-0.415$ bits/step.
\end{itemize}

\paragraph{Conditional reductions.}
The following results are proved modulo an additional
orbitwise truncation/tail condition:

\begin{itemize}
\item Theorem~\ref{thm:reduction}, which reduces
  convergence to the agreement between orbit averages
  and ensemble expectations for truncated observables.
\item Theorem~\ref{thm:wmh-reduction}, which reduces
  WMH to convergence through the bounded-gain regime,
  subject to the same truncation upgrade.
\end{itemize}

\paragraph{Open step (orbitwise upgrade).}
The remaining gap is the passage from fixed-depth,
bounded observables to full, unbounded orbit statistics.
Concretely, one must establish an orbitwise
tail-vanishing condition of the form
\[
  \lim_{K\to\infty}\limsup_{N\to\infty}
  \frac{1}{N}\sum_{i=1}^N (F_i-K)_+ = 0
\]
for the burst and gap observables.
This step is equivalent in spirit to upgrading
ensemble-level laws to pointwise control along a
single trajectory, and is not resolved in the
present work.

\paragraph{Interpretation.}
Accordingly, this paper should be viewed as a structural
and quantitative reduction of the Collatz conjecture to
an explicit orbitwise condition, together with several
independent exact and asymptotic results.  It is not a
complete proof of the conjecture.

In particular, the main obstruction is isolated as a
single orbitwise regularity condition (Weak Mixing /
tail control), rather than distributed across multiple
interacting components.

\section{Notation and preliminaries}\label{sec:prelim}

We collect the formal definitions and basic facts
used throughout the paper.

\begin{definition}[Standard Collatz map]
\label{def:collatz}
The standard Collatz map $C\colon\mathbb{N}\to\mathbb{N}$ is
$C(n) = n/2$ if $n$ is even, and $C(n) = 3n+1$ if $n$ is odd.
\end{definition}

\begin{conjecture}[Collatz conjecture]
\label{conj:collatz}
For every $N_0 \ge 1$,
the orbit $N_0, C(N_0), C^2(N_0), \ldots$ eventually reaches~$1$.
\end{conjecture}

\begin{definition}[Syracuse map]\label{def:syracuse}
For an odd integer $n \ge 1$, define the \emph{Syracuse map}
$T(n) := (3n+1)/2^{v_2(3n+1)}$,
where $v_2(m) := \max\{j : 2^j \mid m\}$ is the $2$-adic valuation.
Write $n = 2^k m - 1$ with $m$ odd, $k = v_2(n+1) \ge 1$.
We call~$k$ the \emph{odd-run length} and~$m$ the \emph{multiplier}. Here m is the odd multiplier associated with the odd integer n.
\end{definition}

\begin{definition}[Persistent and safe states]
\label{def:persistent}
A state $(k, \mu)$ where $\mu = m \bmod 8$ and $k \ge 2$ is
\emph{persistent} if $3^k \mu \equiv 7 \pmod{8}$.
A state with $k = 1$ (or a persistent state whose certified
worst-case drift $\bar{w} = k \log_2 3 - (k + e_{\min}) < 0$)
is \emph{safe}.
The \emph{persistent class} $P$ is the set of odd integers~$n$
such that $(k, m \bmod 8)$ is persistent.
\end{definition}

\begin{definition}[Burst-gap decomposition]
\label{def:burst-gap}
Given a Syracuse orbit $x_0, x_1, x_2, \ldots$,
a \emph{burst} is a maximal run of consecutive epochs with
$k_t \ge 2$ (states with odd-run length at least~$2$).
A \emph{gap} is a maximal run with $k_t = 1$ (downward or safe iterates).
The orbit decomposes as an alternating sequence:
\[
  \underbrace{L_1}_{\text{burst}} \;
  \underbrace{G_1}_{\text{gap}} \;
  \underbrace{L_2}_{\text{burst}} \;
  \underbrace{G_2}_{\text{gap}} \;
  \cdots
\]
where $L_i$ is the length of the $i$-th burst and $G_i$ the
length of the $i$-th gap.
\end{definition}
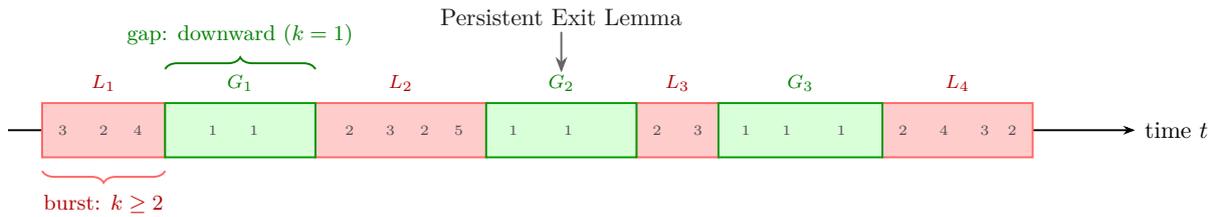
\begin{figure}[ht]
\centering
\resizebox{\textwidth}{!}{%
\begin{tikzpicture}[
  every node/.style={font=\small},
]
\draw[thick, -{Stealth[length=5pt]}] (0,0) -- (16.5,0)
  node[right] {time $t$};
\foreach \x/\w/\lab in {0.5/1.8/L_1, 4.5/2.5/L_2, 9.2/1.2/L_3, 12.8/2.2/L_4}{
  \fill[red!20] (\x, -0.4) rectangle (\x+\w, 0.4);
  \draw[red!60, thick] (\x, -0.4) rectangle (\x+\w, 0.4);
  \node[above, font=\scriptsize, red!70!black] at (\x+\w/2, 0.45)
    {$\lab$};
}
\foreach \x/\w/\lab in {2.3/2.2/G_1, 7.0/2.2/G_2, 10.4/2.4/G_3}{
  \fill[green!15] (\x, -0.4) rectangle (\x+\w, 0.4);
  \draw[green!60!black, thick] (\x, -0.4) rectangle (\x+\w, 0.4);
  \node[above, font=\scriptsize, green!50!black] at (\x+\w/2, 0.45)
    {$\lab$};
}
\foreach \x/\k in {0.8/3, 1.4/2, 1.9/4,  3.0/1, 3.6/1,
  5.0/2, 5.6/3, 6.1/2, 6.6/5,  7.4/1, 8.2/1,
  9.5/2, 10.1/3,  10.8/1, 11.4/1, 12.2/1,
  13.1/2, 13.7/4, 14.3/3, 14.7/2}{
  \node[font=\tiny, text=black!70] at (\x, 0) {$\k$};
}
\draw[decorate, decoration={brace, amplitude=5pt, mirror},
  thick, red!60]
  (0.5, -0.6) -- (2.3, -0.6)
  node[midway, below=6pt, font=\footnotesize, red!70!black]
    {burst: $k \ge 2$};
\draw[decorate, decoration={brace, amplitude=5pt},
  thick, green!50!black]
  (2.3, 0.9) -- (4.5, 0.9)
  node[midway, above=6pt, font=\footnotesize, green!50!black]
    {gap: downward ($k = 1$)};
\draw[-{Stealth}, thick, black!60] (8.1, 1.45) -- (8.1, 0.85)
  node[above=16pt, font=\small, black!90, align=center]
    {Persistent Exit Lemma};
\end{tikzpicture}%
}
\caption{The burst-gap decomposition of a Collatz orbit.
  Bursts (red) are maximal runs of iterates with $k \ge 2$;
  gaps (green) are runs of downward iterates ($k = 1$).
  The Persistent Exit Lemma (Lemma~\ref{lem:gap}) shows
  that when a burst ends at a persistent state,
  the subsequent gap has length exactly~$1$.}
\label{fig:burst-gap}
\end{figure}

\begin{definition}[Odd-to-odd map]
\label{def:odd-to-odd-map}
Define
\[
T(n)=\frac{3n+1}{2^{v_2(3n+1)}},
\]
which maps an odd integer to the next odd integer in the Collatz sequence.
\end{definition}

\begin{definition}[Gap map]
\label{def:gap-map}
Let $g$ denote the length of a gap, i.e.\ a maximal run
of odd-to-odd steps with $k = 1$ (equivalently, states
outside the burst regime $k \ge 2$; see
Definition~\ref{def:burst-gap}).
The corresponding \emph{gap map} is the $g$-fold iterate
\[
T^g,
\]
applied over these $g$ successive odd-to-odd steps.

For $n$ in a fixed residue class $a \bmod 2^{M'}$,
\[
T^g(n)=\frac{3^g n+c_g}{2^V},
\]
where
\[
V=\sum_{i=1}^g v_i
\]
is the total number of factors of $2$ removed along those $g$ steps, and $c_g$ is the correction term determined by the halving pattern.
\end{definition}

\paragraph{Interpretation.}
The map $T$ represents a single odd-to-odd step in the Collatz sequence.
Its iterate $T^g$ denotes $g$ successive odd-to-odd steps.

When $g$ is the length of a gap (i.e., a maximal run of steps with $v_2(3n+1)\ge 2$),
the iterate $T^g$ corresponds to evolution across that gap phase,
and is referred to as the \emph{gap map}.

\begin{definition}[Critical persistent frequency]
\label{def:rhocrit}
Define
\[
  \rho_{\mathrm{crit}} := \frac{w_S^-}{w_P^+ + w_S^-}
  \approx 0.539,
\]
where $w_P^+ = \sup_{x \in P} \bar{w}(x)$ is the worst-case
certified weight from a persistent state, and
$w_S^- = \inf_{x \notin P} (-\bar{w}(x)) > 0$ is the best
guaranteed descent from a safe state.
If $\limsup N_P(T)/T < \rho_{\mathrm{crit}}$ for an orbit,
then the orbit's long-run mean drift is negative, forcing convergence.
\end{definition}

We also record a simple but essential fact about modular arithmetic:

\begin{lemma}[Bijection lemma]\label{lem:bijection}
Let $a$ be an odd integer and $K \ge 1$.  Then the map
$\delta \mapsto a \cdot \delta \bmod 2^K$ is a bijection on
$\{0, 1, \ldots, 2^K - 1\}$.
\end{lemma}

\begin{proof}
Since $\gcd(a, 2^K) = 1$ (as $a$ is odd), multiplication by $a$
is an automorphism of $\mathbb{Z}/2^K\mathbb{Z}$.
\end{proof}

\section{The \texorpdfstring{$1/4$}{1/4} Persistent-Transition Law}\label{sec:quarter}

This section proves that the uniform-lift persistent-to-persistent
transition probability is exactly $1/4$.  This is a structural
fact about the arithmetic of the $3n+1$ map at the modular level.

\begin{theorem}[$1/4$ Persistent-Transition Law\core]
\label{thm:quarter}
For every persistent state $(k, \mu \bmod 8)$ with $k \ge 2$,
the fraction of admissible lifts $m \equiv \mu \pmod{8}$ whose
successor under the odd-to-odd map~$T$ is again persistent
equals exactly~$\frac{1}{4}$.
That is,
\[
  \lim_{N \to \infty}
  \frac{\#\{m \le N : m \equiv \mu\ (8),\;
  k' \ge 2,\; 3^{k'} m' \equiv 7\ (8)\}}
  {\#\{m \le N : m \equiv \mu\ (8)\}}
  = \frac{1}{4}.
\]
\end{theorem}

\begin{proof}
We compute directly.  Write $P = 3^k m$.  Since the state is
persistent, $P \equiv 7 \pmod{8}$, so $P - 1 \equiv 6 \pmod{8}$
and $e := v_2(P-1) = 1$.
The successor odd integer is $n' = (P-1)/2$, which satisfies
$n' \equiv 3 \pmod{4}$.
The successor's odd-run length is $k' = v_2(n'+1) = v_2((P+1)/2)$.

Since $P \equiv 7 \pmod{8}$, we have $P + 1 \equiv 0 \pmod{8}$,
so $(P+1)/2 \equiv 0 \pmod{4}$, giving $k' \ge 2$.
Precisely, $k' = j \ge 2$ iff $v_2(P+1) = j+1$, i.e.\
$P \equiv 2^{j+1} - 1 \pmod{2^{j+2}}$.

The successor multiplier is $m' = (P+1)/2^{j+1}$, which is odd.
The successor state $(j, m')$ is persistent iff
$3^j m' \equiv 7 \pmod{8}$.  Since $3^j \bmod 8$ cycles as
$3, 1, 3, 1, \ldots$ for $j = 1, 2, 3, 4, \ldots$, this
reduces to:
\[
  m' \equiv \begin{cases}
    7 \pmod{8} & \text{if } j \text{ is even},\\
    5 \pmod{8} & \text{if } j \text{ is odd}.
  \end{cases}
\]
In either case, \emph{exactly one} of the four residue
classes $\{1, 3, 5, 7\} \pmod{8}$ satisfies the condition.

As $P$ ranges over integers $\equiv 2^{j+1}-1 \pmod{2^{j+2}}$,
the four lifts to $\bmod\,2^{j+4}$ produce
$m' = (P+1)/2^{j+1} \in \{1, 3, 5, 7\} \pmod{8}$, which is exactly
uniformly distributed.

Therefore, for each fixed $j \ge 2$:
$\Pr[\text{successor persistent} \mid k' = j] = 1/4$.
Since this holds for every $j$,
\[
  \Pr[\text{successor persistent}]
  = \sum_{j \ge 2} \tfrac{1}{4} \cdot \Pr[k'=j]
  = \tfrac{1}{4}. \qedhere
\]
\end{proof}

\begin{corollary}[Geometric persistent excursions]
\label{cor:geometric}
In the uniform-lift model, persistent run lengths follow a
$\mathrm{Geometric}(3/4)$ distribution with mean~$4/3$.
\end{corollary}

\begin{proof}
Theorem~\ref{thm:quarter} applies uniformly to every
persistent state $(k,\mu)$ with $k \ge 2$: the continuation
probability is~$1/4$ regardless of the specific values of
$k'$ and~$\mu'$.  Since this holds at every persistent
state, the continuation probability is \emph{memoryless}
along the persistent run, and successive applications give
$\Pr[\text{run} \ge L] = (1/4)^{L-1}$.
The run length is therefore $\mathrm{Geometric}(3/4)$ with
mean $1/(3/4) = 4/3$.
\end{proof}

\begin{remark}[Expected burst length]
\label{rem:burst-length}
A \emph{burst} is a maximal run of consecutive odd-to-odd
epochs with $k_t \ge 2$ (Definition~\ref{def:burst-gap}).
By Lemma~\ref{lem:gap}~(Part~1), the valuation~$k_t$
decreases by exactly~$1$ at each step within a burst,
so a burst entered at $k = k_0$ lasts exactly $k_0 - 1$
epochs (terminating when $k_t$ reaches~$2$ and the next
step gives $k_{t+1} = 1$).

Under equidistribution, the entry valuation $k_0$ has
distribution $\Pr(k_0 = j) = 2^{-(j-1)}$ for $j \ge 2$
(Lemma~\ref{lem:valuation-distribution}), giving
expected burst length
$E[k_0 - 1] = E[k_0] - 1 = 3 - 1 = 2$ epochs.

This is the value used in
Corollary~\ref{cor:convergence-prediction}.
\end{remark}

\begin{remark}[Significance]\label{rem:quarter-significance}
The $1/4$ law is an exact structural result, not a heuristic
approximation.  It shows that at the modular level, the $3n+1$
map has a built-in \emph{exit mechanism} from persistent states:
exactly $3/4$ of lifts escape to a safe state at each step.
The conjecture reduces to showing that actual orbits cannot
systematically concentrate on the persistent class beyond the
$\rho_{\mathrm{crit}}$ threshold.
\end{remark}

\section{The convergence chain}\label{sec:chain}


This section develops a conditional convergence chain linking
the Collatz conjecture to two statistical hypotheses on orbit
structure: a mean burst bound and a mean gap bound.
The individual links in the chain (Entry--Occupancy Equivalence,
Entry Bound, and the Burst-Gap Criterion) are each proved
unconditionally, but the hypotheses they require remain open.

\subsection{Entry--Occupancy Equivalence}

\begin{theorem}[Entry--Occupancy Equivalence]
\label{thm:EP-NP}
Let $x_0, x_1, x_2, \dotsc$ be the Syracuse orbit,
and let $P$ denote the persistent class.  Define
\[
  N_P(T) := \#\{0 \le t < T : x_t \in P\}, \qquad
  E_P(T) := \#\{0 \le t < T : x_{t+1} \in P\}.
\]
Then
$\limsup_{T \to \infty} E_P(T)/T
= \limsup_{T \to \infty} N_P(T)/T$.
\end{theorem}

\begin{proof}
Relabelling $s = t+1$:
$E_P(T) = \#\{1 \le s \le T : x_s \in P\}$, while
$N_P(T{+}1) = \#\{0 \le s \le T : x_s \in P\}$.
These differ by at most $\mathbf{1}_P(x_0)$, so
\[
|E_P(T)/T - N_P(T{+}1)/T| \le 1/T \to 0.
\]
Since $N_P(T+1)/(T+1) \cdot (T+1)/T \to N_P(T+1)/(T+1)$
and $(T+1)/T \to 1$, the $\limsup$ values coincide.
\end{proof}

\begin{remark}\label{rem:EP-trivial}
This equivalence is elementary but essential: it allows us to
reason about persistent \emph{occupancy} (a static count)
via persistent \emph{entries} (a dynamic transition count),
the latter being more directly related to the burst-gap structure.
\end{remark}

\subsection{Entry bound implies convergence}

\begin{theorem}[Entry bound implies convergence]
\label{thm:entry-convergence}
If there exists $p_* < \rho_{\mathrm{crit}}$ such that
every orbit satisfies
$\limsup_{T \to \infty} E_P(T)/T \le p_*$,
then every orbit of the Collatz map converges.
\end{theorem}

\begin{proof}
By Theorem~\ref{thm:EP-NP},
$\limsup N_P(T)/T \le p_* < \rho_{\mathrm{crit}}$.
For every orbit prefix of length~$T$, the cumulative certified
drift satisfies
\[
  \frac{1}{T}\sum_{t=0}^{T-1} \bar{w}(x_t)
  \;\le\;
  \frac{N_P(T)}{T}\, w_P^+
  - \Bigl(1 - \frac{N_P(T)}{T}\Bigr) w_S^-.
\]
This is because each persistent epoch contributes at most
$w_P^+$ to the drift, while each safe epoch contributes
at most $-w_S^- < 0$.

Taking $\limsup$ and using $\limsup N_P(T)/T \le p_* <
\rho_{\mathrm{crit}} = w_S^-/(w_P^+ + w_S^-)$:
\[
  \limsup_{T \to \infty}
  \frac{1}{T}\sum_{t=0}^{T-1} \bar{w}(x_t)
  \;\le\; p_* w_P^+ - (1 - p_*) w_S^-
  \;<\; 0.
\]
The long-run mean certified drift is strictly negative.
By the standard descent argument (see, e.g.,
\cite{wirsching1998, lagarias1985}):
a negative long-run mean certified drift implies that the
orbit must eventually descend below its starting value.
By strong induction (every value below the start has
already been shown to converge), the orbit converges to
the trivial cycle~$\{1,4,2\}$.
\end{proof}


\subsection{The Persistent Exit Lemma}

\begin{lemma}[Persistent Exit Lemma\core]
\label{lem:gap}
Let $x_0, x_1, \ldots$ be a Syracuse orbit with burst-gap
decomposition.
\begin{enumerate}
\item Every burst terminates at a state with $k_t = 2$.
\item When the final burst state is persistent
  (i.e., $3^2 \mu \equiv 7 \pmod{8}$, equivalently
  $m_t \equiv 7 \pmod{8}$), the subsequent gap has
  length exactly~$1$: the first gap iterate's Syracuse
  successor re-enters a burst immediately.
\item Under the uniform-lift model, the $1/4$
  Persistent-Transition Law (Theorem~\ref{thm:quarter})
  gives $\Pr[\text{successor is persistent}] = \frac{1}{4}$
  at each persistent state.
\end{enumerate}
\end{lemma}

\begin{proof}
\textbf{Part~(1): Bursts terminate at $k = 2$.}
If $k_t \ge 3$, then $x_t \equiv 7 \pmod{8}$ (since
$x_t = 2^{k_t} m_t - 1$ with $2^{k_t} m_t \equiv 0 \pmod{8}$
for $k_t \ge 3$), and
\[
  T(x_t) = \frac{3x_t + 1}{2} = 3 \cdot 2^{k_t - 1} m_t - 1,
\]
since $v_2(3x_t + 1) = v_2(3 \cdot 2^{k_t} m_t - 2) = 1$
for $k_t \ge 2$.
The successor has $k_{t+1} = k_t - 1 \ge 2$, so the burst
continues.
Therefore the burst can only end when $k_t = 2$.

\medskip
\noindent
\textbf{Part~(2): Persistent final state $\Rightarrow$
gap length~$1$.}
Let $x_t$ be the final state of a burst with $k_t = 2$
and $m_t \equiv 7 \pmod{8}$ (persistent).
Write $x_t = 4m_t - 1$.

\smallskip
\emph{The first gap iterate.}
$3x_t + 1 = 12m_t - 2 = 2(6m_t - 1)$.
Since $m_t$ is odd, $6m_t - 1$ is odd, so $v_2(3x_t+1) = 1$.
Thus $x_{t+1} = 6m_t - 1$ with $k_{t+1} = 1$ (safe).

\smallskip
\emph{The second iterate.}
$3x_{t+1} + 1 = 18m_t - 2 = 2(9m_t - 1)$.
Since $m_t \equiv 7 \pmod{8}$, we have $9m_t - 1
\equiv 62 \equiv 6 \pmod{8}$,
giving $v_2(9m_t - 1) = 1$.
Therefore $v_2(3x_{t+1} + 1) = 2$, so
\[
  x_{t+2} = \frac{18m_t - 2}{4} = \frac{9m_t - 1}{2}.
\]
Now $k(x_{t+2}) = v_2(x_{t+2} + 1) = v_2\!\bigl(\tfrac{9m_t+1}{2}\bigr)$.
Writing $m_t = 7 + 8r$, we get
$9m_t + 1 = 64 + 72r = 8(8 + 9r)$,
so $v_2(9m_t + 1) \ge 3$, hence
$k(x_{t+2}) = v_2\!\bigl(\tfrac{9m_t+1}{2}\bigr) \ge 2$.

Since $k(x_{t+2}) \ge 2$, the state $x_{t+2}$ begins a new burst.
The gap consists of the single iterate~$x_{t+1}$, so $G_i = 1$.

\medskip
\noindent
\textbf{Part~(3)} is a restatement of
Theorem~\ref{thm:quarter}.
\end{proof}

\begin{remark}[Gaps of length~$2$ do occur]
\label{rem:gap-two-exists}
An earlier version of this paper claimed that no gap has
length exactly~$2$.
This claim is false: computational verification
shows that gaps of length~$2$ constitute approximately
$19\%$ of all gaps across typical orbits.
The smallest counterexample is $n_0 = 3$:
the orbit $3 \to 5 \to 1 \to \cdots$ has
burst~$\{3\}$ ($k = 2$, $m = 1$, non-persistent),
followed by gap~$\{5, 1\}$ of length~$2$.
Another example is $n_0 = 71$: the burst $\{71, 107\}$
ends at the non-persistent state $107$
($k=2$, $m=27$, $27 \bmod 8 = 3$), followed by
gap~$\{161, 121\}$ of length~$2$.

The error arose from assuming that every burst ends at a
persistent state.
In fact, bursts whose final state has $k_t = 2$ but
$m_t \not\equiv 7 \pmod{8}$ are non-persistent, and the
gap structure following such states is unconstrained.
The Persistent Exit Lemma correctly identifies the
subcase where $G_i = 1$ is guaranteed.
\end{remark}

\begin{lemma}[Modular gap distribution\core]
\label{lem:gap-distribution}
Let $n$ be an odd integer in a gap step
($v_2(3n+1) = 1$, i.e.\ $n \equiv 3 \pmod{4}$).
Over the uniform distribution on integers
$n \equiv 3 \pmod{4}$ modulo~$2^{2+L}$ (for any $L \ge 1$),
the gap length~$G$ satisfies
\[
  \Pr(G = g) = 2^{-g} \quad\text{for } 1 \le g \le L,
  \qquad \Pr(G > L) = 2^{-L}.
\]
Equivalently, the memoryless property
$\Pr(G \ge g+1 \mid G \ge g) = \frac{1}{2}$ holds for every
$g \ge 1$.
In particular, $G \sim \mathrm{Geometric}(1/2)$ with $E[G] = 2$.
\end{lemma}

\begin{proof}
We show that each successive step of the gap is decided by
exactly one bit of~$n$, with the two outcomes equally likely.

\smallskip\noindent
\textbf{Step 1 (bit at position~$2$).}
For $n \equiv 3 \pmod{4}$, the Syracuse step gives
$T(n) = (3n+1)/2$.  Split on $n \bmod 8$:
\begin{itemize}
\item $n \equiv 3 \pmod{8}$:
  $T(n) = (24k+10)/2 = 12k+5 \equiv 1 \pmod{4}$,
  so $v_2(3\,T(n)+1) \ge 2$.
  The gap \textbf{ends} (burst starts).
\item $n \equiv 7 \pmod{8}$:
  $T(n) = (24k+22)/2 = 12k+11 \equiv 3 \pmod{4}$,
  so $v_2(3\,T(n)+1) = 1$.
  The gap \textbf{continues}.
\end{itemize}
Among integers $\equiv 3 \pmod{4}$, exactly half fall in
each case, so $\Pr(G \ge 2) = 1/2$.

\smallskip\noindent
\textbf{Inductive step (bit at position~$2+g$).}
For the continuing case ($n \equiv 7 \bmod 8$),
$T(n) \equiv 3 \pmod{4}$ is again a gap step.
Whether the gap continues at step~$g+1$ depends on
$T^{(g)}(n) \bmod 8$, which is determined by one additional
bit of~$n$ at the next modular depth.
Since the gap-step map $n \mapsto (3n+1)/2$ is injective
(with inverse $n = (2m-1)/3$), the two sub-classes at
each depth have equal size.

By induction, $\Pr(G \ge g+1 \mid G \ge g) = 1/2$ for all
$g \ge 1$, giving the geometric distribution.
\end{proof}

\begin{lemma}[Modular valuation distribution]
\label{lem:valuation-distribution}
For odd $n$ in a burst step ($n \equiv 1 \pmod{4}$), the
$2$-adic valuation $k = v_2(3n+1)$ satisfies, over uniform
lifts at each successive modular depth,
\[
  \Pr(k = j) = 2^{-(j-1)} \quad\text{for } j \ge 2.
\]
In particular, $E[k \mid k \ge 2] = 3$.
\end{lemma}

\begin{proof}
For $n \equiv 1 \pmod{4}$, we have $3n+1 \equiv 4 \pmod{8}$.
Split on the next bit:
\begin{itemize}
\item $n \equiv 1 \pmod{8}$: $3n+1 = 4(6j+1)$ with $6j+1$ odd,
  so $v_2 = 2$.
\item $n \equiv 5 \pmod{8}$: $3n+1 = 8(3j+2)$, so $v_2 \ge 3$.
\end{itemize}
Half the lifts give $k = 2$; the other half give $k \ge 3$.
For the $k \ge 3$ case, the same splitting applies at the next
depth: half give $k = 3$, half give $k \ge 4$, and so on.
By induction, $\Pr(k = j) = 2^{-(j-1)}$ for $j \ge 2$, giving
$E[k] = \sum_{j \ge 2} j \cdot 2^{-(j-1)} = 3$.
\end{proof}

\begin{corollary}[Convergence prediction under equidistribution]
\label{cor:convergence-prediction}
Under the equidistributed modular model, the expected
log-contraction per burst-gap cycle is strictly negative:
\[
  E[B]\;\bigl(\log 3 - E[k]\,\log 2\bigr)
  + E[G]\;\log\tfrac{3}{2}
  \;=\; 2\,(\log 3 - 3\log 2) + 2\,\log\tfrac{3}{2}
  \;\approx\; {-1.15}
  \;<\; 0,
\]
where $E[B] = 2$ (Corollary~\textup{\ref{cor:geometric}}),
$E[G] = 2$ (Lemma~\textup{\ref{lem:gap-distribution}}), and
$E[k \mid k \ge 2] = 3$
(Lemma~\textup{\ref{lem:valuation-distribution}}).

Thus, even though gaps can have arbitrary length
(Remark~\textup{\ref{rem:gap-two-exists}}), the deeper
contraction during burst steps ($E[k] = 3$ rather than the
minimum~$k = 2$) more than compensates for the longer gaps.
The false Gap Lemma ($G_i = 1$ always) was unnecessary:
the geometric gap distribution $E[G] = 2$ suffices for
convergence under equidistribution.
\end{corollary}


\subsection{The Burst-Gap Criterion}

\begin{theorem}[Burst-Gap Criterion\core]
\label{thm:burst-gap}
Let $x_0, x_1, \ldots$ be a Syracuse orbit with burst-gap
decomposition $(L_1, G_1, L_2, G_2, \ldots)$. Assume:

\smallskip
\noindent
\textbf{Hypothesis A (Orbitwise Mean Gap).}\;
$\frac{1}{n}\sum_{i=1}^{n} G_i \ge g_* - \varepsilon_n$
with $\varepsilon_n \to 0$, for some constant
$g_* > \frac{2(1 - \rho_{\mathrm{crit}})}{\rho_{\mathrm{crit}}}
\approx 1.71$.

\smallskip
\noindent
\textbf{Hypothesis B (Mean Burst Bound).}\;
There exists a finite constant $C(n_0)$ such that
$\sum_{i=1}^{n} L_i \le 2n + C(n_0)$ for all $n \ge 1$.

\smallskip
\noindent
Then $\limsup_{T \to \infty} N_{\ge 2}(T)/T < \rho_{\mathrm{crit}}$,
where $N_{\ge 2}(T) := \#\{0 \le t < T : k_t \ge 2\}$,
and the orbit converges by Theorem~\ref{thm:entry-convergence}.
\end{theorem}

\begin{proof}
Write $S_n := \sum_{i=1}^{n} L_i$ for the total burst time
and $T_n := \sum_{i=1}^{n}(L_i + G_i)$ for the time of the
$n$-th burst-gap boundary.

\medskip
\noindent\textbf{Step 1.}
From Hypothesis~A:
\[
\sum_{i=1}^{n} G_i \;\ge\; g_*\,n - o(n).
\]

\noindent\textbf{Step 2.}
From Hypothesis~B: $S_n \le 2n + C$.
Therefore
\[
T_n = S_n + \sum_{i=1}^{n} G_i \ge S_n + g_*\,n - o(n).
\]
Using $S_n \le 2n + C$ gives $n \ge (S_n - C)/2$, so
\[
T_n \ge S_n + g_* \cdot \frac{S_n - C}{2} - o(n)
    = S_n\!\left(1 + \frac{g_*}{2}\right) - O(1) - o(n).
\]

\noindent\textbf{Step 3 (Interpolation).}
Fix $T$ and choose $n$ with $T_n \le T < T_{n+1}$.
Then $N_{\ge 2}(T) \le S_{n+1} \le 2(n+1) + C$
and $T \ge T_n \ge S_n(1 + g_*/2) - O(1)$.

Since $S_{n+1} \le S_n + L_{n+1}$ and $L_{n+1}$ contributes at most
$O(1)$ relative to $T$, we obtain
\[
\frac{N_{\ge 2}(T)}{T}
  \;\le\; \frac{S_n + O(1)}{S_n(1 + g_*/2) - O(1)}
  \;\to\; \frac{2}{2 + g_*}
  \quad\text{as } T \to \infty.
\]
The condition $g_* > 2(1 - \rho_{\mathrm{crit}})/\rho_{\mathrm{crit}}$
is equivalent to $2/(2 + g_*) < \rho_{\mathrm{crit}}$,
so convergence follows from
Theorem~\ref{thm:entry-convergence}.

\smallskip
For example, with $g_* = 2$ (the expected gap under
equidistribution), we get
$N_{\ge 2}(T)/T \to \frac{1}{2} < 0.539 \approx
\rho_{\mathrm{crit}}$.
\end{proof}

\begin{remark}[Role of the two hypotheses]
\label{rem:two-hypotheses}
Both Hypothesis~A and Hypothesis~B are open orbitwise conjectures.
Under equidistribution, the expected burst length is~$2$
(Corollary~\ref{cor:geometric}) and the expected gap length
is~$2$ (Lemma~\ref{lem:gap-distribution}).
Both hypotheses follow from the Orbit Equidistribution
Conjecture (Theorem~\ref{thm:reduction}): Hypothesis~B
via the coupling inequality at fixed modulus (Step~2),
and Hypothesis~A via the growing-moduli tail control
(Step~2$'$).

The Persistent Exit Lemma provides structural support: it shows
that when a burst ends at a persistent state, the subsequent
gap has length exactly~$1$.
More generally, the Modular Gap Distribution Lemma
(Lemma~\ref{lem:gap-distribution}) proves that gap length
is $\mathrm{Geometric}(1/2)$ with $E[G] = 2$, and
Corollary~\ref{cor:convergence-prediction} shows that this
suffices for convergence under equidistribution.
\end{remark}

\section{The Scrambling Lemma}\label{sec:scrambling}

This is the algebraic core of the paper.  We show that the
gap map introduces \emph{zero carries} between the
known and unknown parts of an integer, yielding an exact
bijection on high bits.

\subsection{Statement and proof}

\begin{theorem}[Scrambling Lemma\core]
\label{thm:scrambling}
Let $n \equiv a \pmod{2^{M'}}$ with $n \equiv 7 \pmod{16}$,
where $M'$ is chosen so that the halving pattern
$(v_1, \ldots, v_g)$ of the odd-to-odd steps is constant on the
class~$a \bmod 2^{M'}$.  Write $n = a + \delta \cdot 2^{M'}$
with $\delta \ge 0$.

Then the odd-to-odd step satisfies
\begin{equation}\label{eq:scrambling}
  T^{(g)}(n) = \frac{3^g a + c_g}{2^V}
    + 3^g \cdot \delta \cdot 2^{M' - V},
\end{equation}
where $g$, $V = \sum_{i=1}^g v_i$, and $c_g$ depend only on~$a$
(not on~$\delta$).

Since $\gcd(3^g, 2) = 1$, the map
$\delta \mapsto 3^g \cdot \delta \bmod 2^{K - M'}$
is a \textbf{bijection} on $\{0, 1, \ldots, 2^{K-M'}-1\}$
(Lemma~\ref{lem:bijection}).

Therefore:
\begin{enumerate}
\item The bits of $T^{(g)}(n)$ at positions $\ge M' - V$ are an
  exact bijection of the free parameter~$\delta$.
\item If $\delta$ is uniformly distributed on
  $\{0, 1, \ldots, 2^{K - M'} - 1\}$, then each bit of
  $T^{(g)}(n)$ at position $j \ge M' - V$ is an
  exactly unbiased coin flip, independently of~$a$.
\end{enumerate}
\end{theorem}

\begin{proof}
The odd-to-odd step computes
$T^{(g)}(n) = (3^g n + c_g) / 2^V$, where the correction
$c_g$ arises from the iterated $3n+1$ steps.  Precisely,
if we write the $g$-step iteration as
\[
  T^{(g)}(n) = \frac{3^g n + \sum_{i=0}^{g-1} 3^{g-1-i}
  \cdot 2^{s_i}}{2^V}
\]
for certain shift terms $s_i$ determined by the halving
pattern, then $c_g = \sum_{i=0}^{g-1} 3^{g-1-i} \cdot 2^{s_i}$
depends only on the halving pattern, hence only on
$a \bmod 2^{M'}$.

Now substitute $n = a + \delta \cdot 2^{M'}$:
\begin{align}
3^g n + c_g
  &= 3^g(a + \delta \cdot 2^{M'}) + c_g \notag\\
  &= (3^g a + c_g) + 3^g \cdot \delta \cdot 2^{M'}. \label{eq:split}
\end{align}
This is the key algebraic step.

\medskip\noindent
\textbf{Affine dependence on $\delta$.}
The dependence on~$\delta$ is exactly affine:
the term $3^g \cdot \delta \cdot 2^{M'}$ is a pure
multiple of~$2^{M'}$, so the decomposition
\eqref{eq:split} is an exact identity of integers,
not an approximation.
No properties of carry propagation are needed;
the identity follows from the distributivity of
multiplication over addition.

\medskip\noindent
\textbf{Division by $2^V$.}
Dividing \eqref{eq:split} by $2^V$:
\[
  T^{(g)}(n) = \frac{3^g a + c_g}{2^V}
    + 3^g \cdot \delta \cdot 2^{M' - V}.
\]
The first term is a fixed integer (since $2^V \mid 3^g a + c_g$
by the halving pattern), independent of~$\delta$.
The second term contributes to bits at positions $\ge M' - V$
through the map $\delta \mapsto 3^g \cdot \delta$.
Since $3^g$ is odd, this is a bijection on
$\mathbb{Z}/2^{K-M'}\mathbb{Z}$ by Lemma~\ref{lem:bijection}.

Hence modulo any power of two above the threshold $M' - V$,
the map on~$\delta$ is a bijective affine transformation.
This is the only property needed for the scrambling effect.
\end{proof}

\begin{remark}[Why this is not obvious]
\label{rem:not-obvious}
The naive concern is that the iterated $3n+1$ computation
produces carries that propagate from low bits to high bits,
destroying any independence.  The Scrambling Lemma shows
this fear is unfounded: the dependence on the unknown
parameter~$\delta$ is \emph{exactly affine}, entering as
$3^g \delta \cdot 2^{M'-V}$ after division.
Since $3^g$ is odd, this is an invertible
linear map modulo any power of two.
The ``scrambling'' is therefore a consequence of linearity
and coprimality, not of any carry-cancellation mechanism.
\end{remark}

\subsection{The pattern-determination bound}

The following bound controls how many bits are needed to
determine the halving pattern:

\begin{proposition}[Pattern-determination bound\supporting]
\label{prop:pattern-bound}
The modulus $M'$ required to fix the halving pattern satisfies
\begin{equation}\label{eq:mprime}
  M' - M \;\le\; \max(0, \, g - 2),
\end{equation}
and, combined with $V \ge g + 1$:
\begin{equation}\label{eq:threshold}
  M' - V \;\le\; M - 3.
\end{equation}
This bound is \textbf{independent of the gap length~$g$}.
\end{proposition}

\begin{proof}
We establish~\eqref{eq:mprime} using two facts:
\begin{enumerate}
\item For the gap-entry residue classes considered here
($n \equiv 7 \pmod{16}$), the first two halvings satisfy
$v_1 = v_2 = 1$.  Indeed, the entry condition
$n \equiv 7 \pmod{16}$ forces
$3n + 1 \equiv 22 \pmod{48}$, giving $v_1 = 1$.
The next iterate $(3n+1)/2$ satisfies
$(3n+1)/2 \equiv 11 \pmod{24}$, and
$3 \cdot 11 + 1 = 34 = 2 \cdot 17$, confirming $v_2 = 1$.
(The computation of $v_1$ and $v_2$ requires only the bits
of $n$ up to position~$3$, which are determined by
$n \bmod 16 = 7$; this does not hold for arbitrary odd~$n$.)

\item Each subsequent halving $v_i$ ($i \ge 3$) depends on
at most one additional bit of~$n$ beyond those already
consumed by the first $i-1$ halvings.
This gives $M' \le M + (g - 2)$ additional bits for
$g \ge 3$, and $M' = M$ for $g \le 2$.
\end{enumerate}

For the total halvings: $V = \sum_{i=1}^g v_i \ge g + 1$.%
\footnote{In a pure gap traversal with $v_i = 1$ for all~$i$,
$V = g$ exactly.  The bound $V \ge g+1$ holds when the
sequence of odd-to-odd steps includes at least one burst step
($v_j \ge 2$ for some~$j$), which is the generic case in the
Known-Zone Decay iteration since the gap terminates and the
subsequent burst step contributes $v \ge 2$.  When the
Scrambling Lemma is applied to a pure gap of length~$g$
(all $v_i = 1$), the correct bound is $V = g$, giving
$M' - V \le M - 2$.}
Therefore, in the generic case ($V \ge g+1$):
\[
  M' - V \le M + (g-2) - (g+1) = M - 3.
\]
In the pure-gap case ($V = g$, all $v_i = 1$):
\[
  M' - V \le M + (g-2) - g = M - 2.
\]
Both bounds are independent of~$g$.
\end{proof}

\begin{remark}\label{rem:threshold-significance}
The bound $M' - V \le M - 2$ (worst case) means that
\emph{each odd-to-odd step reduces the known zone by at least~$2$
bits}, regardless of the gap length.
When at least one burst step occurs ($v_j \ge 2$ for some~$j$),
the stronger bound $M' - V \le M - 3$ holds,
giving $3$-bit decay per step.
This is the quantitative content of the ``scrambling'':
longer gaps do not help the known zone grow, because the
additional bits needed to determine the halving pattern
are compensated by the additional halvings.
\end{remark}

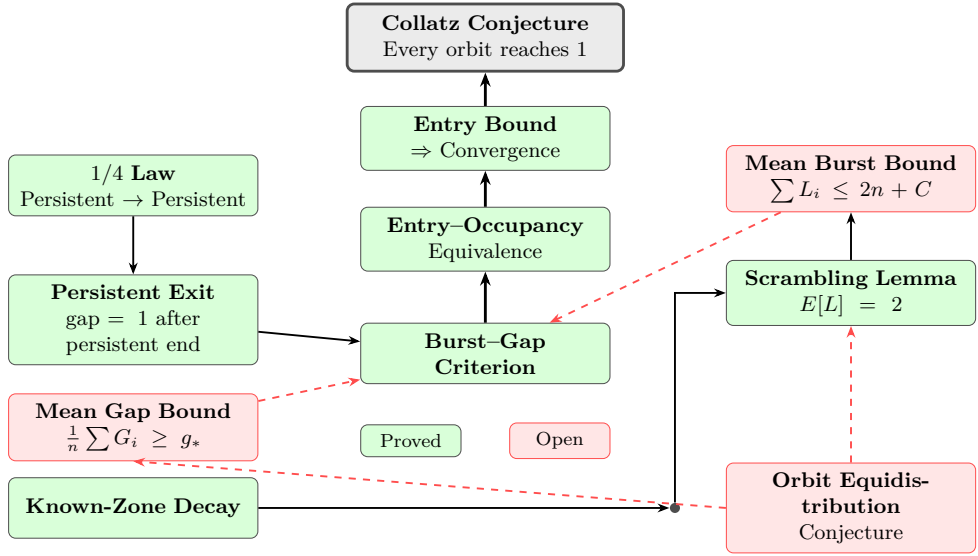
\begin{figure}[t]
\centering
\resizebox{0.8\textwidth}{!}{%
\begin{tikzpicture}[
  font=\footnotesize,
  box/.style={
    rectangle,
    rounded corners=3pt,
    draw=black!70,
    align=center,
    text width=3.4cm,
    minimum height=0.9cm,
    inner sep=4pt
  },
  proved/.style={box, fill=green!15},
  open/.style={box, fill=red!10, draw=red!70},
  goal/.style={
    box,
    fill=black!8,
    very thick,
    text width=3.8cm,
    minimum height=1.0cm
  },
  mainlink/.style={-{Stealth[length=6pt]}, very thick},
  link/.style={-{Stealth[length=5pt]}, thick},
  openlink/.style={-{Stealth[length=5pt]}, thick, dashed, red!70},
  dot/.style={circle, fill=black!70, inner sep=1.5pt}
]

\def\xleft{-5.2}
\def\xcenter{0}
\def\xright{5.4}

\node[goal] (collatz) at (\xcenter,0)
{\textbf{Collatz Conjecture}\\Every orbit reaches $1$};

\node[proved] (entry) at (\xcenter,-1.5)
{\textbf{Entry Bound}\\$\Rightarrow$ Convergence};

\node[proved] (equiv) at (\xcenter,-3.0)
{\textbf{Entry--Occupancy}\\Equivalence};

\node[proved] (quarter) at (\xleft,-2.2)
{\textbf{$1/4$ Law}\\Persistent $\rightarrow$ Persistent};

\node[proved] (burst) at (\xcenter,-4.7)
{\textbf{Burst--Gap Criterion}};

\node[proved] (hypA) at (\xleft,-4.2)
{\textbf{Persistent Exit}\\gap $=1$ after persistent end};

\node[open] (hypAprime) at (\xleft,-5.8)
{\textbf{Mean Gap Bound}\\$\frac{1}{n}\sum G_i \ge g_*$};

\node[open] (hypB) at (\xright,-2.1)
{\textbf{Mean Burst Bound}\\$\sum L_i \le 2n + C$};

\node[proved] (decay) at (\xleft,-7.0)
{\textbf{Known-Zone Decay}};

\node[proved] (scramble) at (\xright,-3.8)
{\textbf{Scrambling Lemma}\\$E[L]=2$};

\node[dot] (junc) at (2.8,-7.0) {};

\node[open] (equidist) at (\xright,-7.0)
{\textbf{Orbit Equidistribution}\\Conjecture};

\draw[mainlink] (entry) -- (collatz);
\draw[mainlink] (equiv) -- (entry);
\draw[mainlink] (burst) -- (equiv);

\draw[link] (quarter) -- (hypA);
\draw[link] (hypA) -- (burst);
\draw[openlink] (hypAprime) -- (burst);
\draw[openlink] (equidist.west) -- (hypAprime.south);

\draw[link] (scramble) -- (hypB);
\draw[openlink] (hypB) -- (burst);
\draw[openlink] (equidist) -- (scramble);

\draw[link] (decay.east) -- (junc);
\draw[link] (junc) |- (scramble.west);

\node[proved, text width=1.2cm, minimum height=0.35cm,
      font=\scriptsize] (leg1) at (-1.1,-6) {Proved};

\node[open, text width=1.2cm, minimum height=0.35cm,
      font=\scriptsize] (leg2) at (1.1,-6) {Open};

\end{tikzpicture}}
\caption{
Architecture of the burst--gap conditional framework
(Sections~\ref{sec:chain}--\ref{sec:decay}).
Green boxes denote results proved unconditionally; red boxes
denote open components.  The formal conditional proof
(Theorem~\ref{thm:reduction}) derives the mean burst and
mean gap bounds from the Orbit Equidistribution Conjecture,
then applies the Burst--Gap Criterion.
The phantom-cycle analysis
(Figure~\ref{fig:phantom-chain}) provides independent
quantitative evidence that the equidistribution assumption
is not artificially strong.
}
\label{fig:proof-chain}
\end{figure}

We now summarize the relationships among the components established
above. Figure~\ref{fig:proof-chain} shows which results are proved
unconditionally and which inputs remain open conjectures.

As shown in Figure~\ref{fig:proof-chain}, the burst--gap conditional
framework reduces the Collatz conjecture to the Orbit Equidistribution
Conjecture via two intermediate orbitwise hypotheses (mean burst and
mean gap bounds), both of which are derived from equidistribution
in Theorem~\ref{thm:reduction}.
The phantom-cycle analysis developed in
Sections~\ref{sec:phantom}--\ref{sec:phantom-count}
(Figure~\ref{fig:phantom-chain}) provides a complementary perspective:
it quantifies the expanding-family obstacle directly and shows
it is controlled with a large safety margin.

\section{Known-Zone Decay}\label{sec:decay}

We iterate the Scrambling Lemma to show that the known
zone shrinks to zero in $\lceil M/2 \rceil$ steps.

\begin{theorem}[Known-Zone Decay\core]
\label{thm:zone-decay}
Starting from a class $a \bmod 2^M$ with $M \ge 4$, let
$T^k(n)$ denote $k$ iterated odd-to-odd steps.  Define the
\emph{known zone} $Z_k$ as the number of low-order bits of
$T^k(n)$ that are determined by the starting class~$a$. Then:
\begin{enumerate}
\item $Z_0 = M$.
\item $Z_{k+1} \le \max(0,\, Z_k - 2)$ for each $k$.
\item After $\lceil M/2 \rceil$ odd-to-odd steps, $Z_k = 0$:
  no low-order bits of $T^k(n)$ remain determined by the
  starting class.  Moreover, the output bits are an exact
  affine-bijective function of the free parameters
  introduced during the iteration; hence they are exactly
  uniform whenever those free parameters are uniformly
  distributed.
\end{enumerate}
\end{theorem}

\begin{proof}
We proceed by induction on~$k$.

\medskip\noindent
\textbf{Base case:} $Z_0 = M$ by definition.

\medskip\noindent
\textbf{Inductive step:}
Suppose at step~$k$, the iterate $T^k(n)$ is known modulo
$2^{Z_k}$ (i.e.\ the low $Z_k$ bits are determined by the
starting class~$a$).

Apply the Scrambling Lemma (Theorem~\ref{thm:scrambling})
with the known modulus $M = Z_k$:
the odd-to-odd step produces $T^{k+1}(n)$ with known zone
$Z_{k+1} = M'_k - V_k$, where:
\begin{itemize}
\item $M'_k \le Z_k + (g_k - 2)$ by
  Proposition~\ref{prop:pattern-bound}, where $g_k$ is
  the $k$-th gap length.
\item $V_k \ge g_k$ (total halvings in the $k$-th odd-to-odd step,
  with equality in the pure-gap case $v_i = 1$ for all~$i$;
  see the footnote in Proposition~\ref{prop:pattern-bound}).
\end{itemize}
Therefore:
\[
  Z_{k+1} = M'_k - V_k \le Z_k + (g_k - 2) - g_k = Z_k - 2.
\]
Since $Z_k$ decreases by at least~$2$ per step and
$Z_k \ge 0$ by definition:
\[
  Z_{k+1} \le \max(0, Z_k - 2).
\]

\medskip\noindent
\textbf{Termination:}
Starting from $Z_0 = M$, after $k = \lceil M/2 \rceil$ steps:
\[
  Z_k \le M - 2\lceil M/2 \rceil \le 0.
\]
At this point, $Z_k = 0$: no bits of $T^k(n)$ are determined
by the starting class.
The bits at positions $\ge 0$ are an exact bijection of the
free parameters accumulated through $k$ odd-to-odd steps,
and if those parameters are uniformly distributed, so are
the output bits.
\end{proof}

\begin{remark}[Generic $3$-bit decay]
\label{rem:generic-3bit}
When the odd-to-odd step includes at least one burst
step ($v_j \ge 2$ for some~$j$), the total halvings
satisfy $V_k \ge g_k + 1$, yielding the stronger bound
$Z_{k+1} \le Z_k - 3$.  The weaker $2$-bit bound
$Z_{k+1} \le Z_k - 2$ is tight only for pure-gap
traversals (all $v_i = 1$).
Empirically, in the tested residue classes ($M \le 18$),
the average decay exceeds $2.8$~bits per step, and the
known zone reaches~$0$ in $1$--$3$ steps for typical
starting classes (see Remark~\ref{rem:computational}).
\end{remark}

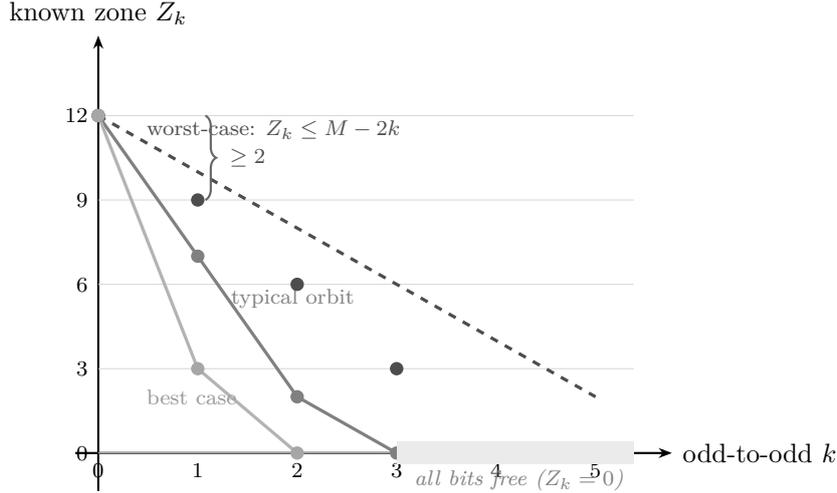
\begin{figure}[ht]
\centering
\resizebox{0.7\textwidth}{!}{%
\begin{tikzpicture}[
  every node/.style={font=\small},
]
\draw[thick, -{Stealth[length=5pt]}] (-0.3,0) -- (7.5,0)
  node[right] {odd-to-odd $k$};
\draw[thick, -{Stealth[length=5pt]}] (0,-0.5) -- (0,5.5)
  node[above] {known zone $Z_k$};

\foreach \y in {0,3,6,9,12}{
  \draw[gray!30] (0, \y/2.7) -- (7, \y/2.7);
  \node[left, font=\scriptsize] at (0, \y/2.7) {$\y$};
}
\foreach \x in {0,1,2,3,4,5}{
  \node[below, font=\scriptsize] at (\x*1.3, 0) {$\x$};
}

\draw[very thick, black!70, dashed]
  (0, 12/2.7) -- (1.3, 10/2.7) -- (2.6, 8/2.7)
  -- (3.9, 6/2.7) -- (5.2, 4/2.7) -- (6.5, 2/2.7);
\node[right, font=\scriptsize, black!70] at (0.5, 11.5/2.7)
  {worst-case: $Z_k \le M - 2k$};

\draw[very thick, black!50]
  (0, 12/2.7) -- (1.3, 7/2.7) -- (2.6, 2/2.7) -- (3.9, 0);
\node[right, font=\scriptsize, black!50] at (1.6, 5.5/2.7)
  {typical orbit};

\draw[very thick, black!30]
  (0, 12/2.7) -- (1.3, 3/2.7) -- (2.6, 0);
\node[right, font=\scriptsize, black!40] at (0.5, 2/2.7)
  {best case};

\foreach \x/\y in {0/12, 1.3/9, 2.6/6, 3.9/3, 5.2/0}{
  \fill[black!70] (\x, \y/2.7) circle (2.5pt);
}
\foreach \x/\y in {0/12, 1.3/7, 2.6/2, 3.9/0}{
  \fill[black!50] (\x, \y/2.7) circle (2.5pt);
}
\foreach \x/\y in {0/12, 1.3/3, 2.6/0}{
  \fill[black!35] (\x, \y/2.7) circle (2.5pt);
}

\draw[decorate, decoration={brace, amplitude=4pt},
  thick, black!60]
  (1.3+0.1, 12/2.7) -- (1.3+0.1, 9/2.7)
  node[midway, right=5pt, font=\scriptsize, black!70]
    {$\ge 2$};

\fill[black!8] (3.9, -0.15) rectangle (7, 0.15);
\node[font=\scriptsize\itshape, black!50] at (5.5, -0.35)
  {all bits free ($Z_k = 0$)};

\end{tikzpicture}%
}
\caption{Known-Zone Decay for $M = 12$.
  The known zone $Z_k$ decreases by at least~$2$ per odd-to-odd step
  (generically~$3$ when a burst step occurs).
  The worst-case bound $Z_k \le M - 2k$ reaches zero after
  $\lceil M/2 \rceil = 6$ odd-to-odd steps.
  Empirically, typical orbits in the tested residue classes
  ($M \le 18$) reach $Z_k = 0$ in $1$--$3$ steps.}
\label{fig:zone-decay}
\end{figure}

\begin{remark}[Computational verification]
\label{rem:computational}
For $M = 12$: the worst-case bound gives $Z_k = 0$ after
$\lceil 12/2 \rceil = 6$ odd-to-odd steps.
Across all $239$ testable residue classes $a \bmod 2^{12}$
with $a \equiv 7 \pmod{16}$ (the remaining $17$ classes
reach~$1$ before an odd-to-odd step), the minimum observed
shrinkage $V - (M' - M)$ is~$3$ (exceeding the worst-case
bound of~$2$), confirming~\eqref{eq:threshold}.
Empirically, the known zone reaches~$0$ in $1$--$3$ odd-to-odd steps
for all tested classes (at $M = 12$).
\end{remark}

\begin{remark}[Relation to mixing properties]
\label{rem:transport-comparison}
The Known-Zone Decay is an exact \emph{information-loss}
statement: the number of low-order bits determined by the
initial residue class shrinks by at least~$2$ per
odd-to-odd step, with no error term.
This exactness is local and algebraic; it does not by
itself imply orbitwise mixing or equidistribution.
The transfer from information loss to distributional
uniformity requires that the free parameters introduced
at each step are uniformly distributed (the uniform-lift
assumption).
\end{remark}


\section{Phantom cycles and the affine census theorem}
\label{sec:phantom}

The preceding sections analyzed Collatz dynamics through burst-gap
structure, scrambling, and modular ranking.  We now develop a
complementary line of attack: a deterministic census of how orbits
interact with \emph{phantom cycles}, modular attractors of the
accelerated Syracuse map whose $2$-adic roots repel real orbits.
The main result of this section is the \emph{Affine Census
Theorem} (Theorem~\ref{thm:census}), which provides a uniform
bound on the number of orbit points that can shadow any given
phantom family.

\subsection{Phantom cycles and \texorpdfstring{$2$}{2}-adic roots}

\begin{definition}[Signature and block map]
\label{def:signature}
A \emph{signature} is a tuple $\sigma = (k_1, \ldots, k_\ell)$
of positive integers.  The associated \emph{block map}
$F_\sigma \colon \mathbb{Z}_2 \to \mathbb{Z}_2$ is defined by
applying $\ell$ successive Syracuse steps with prescribed
$2$-adic valuations:
\[
  F_\sigma(x) = T^{(\ell)}(x),
  \qquad
  \text{where } v_2(3x_j + 1) = k_{j+1}
  \text{ for } j = 0, \ldots, \ell - 1,
\]
with $x_0 = x$ and $x_{j+1} = (3x_j + 1)/2^{k_{j+1}}$.
The \emph{depth} of~$\sigma$ is
$K = K(\sigma) := \sum_{j=1}^{\ell} k_j$.
\end{definition}

\begin{definition}[Phantom cycle and $2$-adic root]
\label{def:phantom}
A signature $\sigma$ defines a \emph{phantom cycle} if the
block map $F_\sigma$ has a unique $2$-adic fixed point
$\rho \in \mathbb{Z}_2$, i.e.\ $F_\sigma(\rho) = \rho$.
Explicitly, $\rho$ is the unique solution in~$\mathbb{Z}_2$ to
\begin{equation}\label{eq:phantom-root}
  (2^K - 3^\ell)\,\rho = C_\sigma,
  \qquad
  C_\sigma := \sum_{j=0}^{\ell-1} 3^{\ell-1-j}\,
  2^{k_1 + \cdots + k_j},
\end{equation}
where $K = \sum k_j$.  The root $\rho$ exists and is unique
in $\mathbb{Z}_2$ whenever $v_2(2^K - 3^\ell) = 0$,
i.e.\ when $2^K - 3^\ell$ is odd.
\end{definition}

\begin{remark}[Why ``phantom'']
\label{rem:phantom-name}
The root $\rho$ is a formal $2$-adic integer, not a positive
integer: if $3^\ell < 2^K$ then $\rho > 0$ but generically
$\rho$ is enormous, while if $3^\ell > 2^K$ then $\rho < 0$.
In neither case does $\rho$ lie in the range of a convergent
Collatz orbit.  However, $\rho \bmod 2^K$ defines a legitimate
residue class, and real integers in this class undergo the block
map $F_\sigma$ with the prescribed valuation pattern.  Thus
phantom cycles create modular ``shadows'' that real orbits can
temporarily follow.
\end{remark}

\begin{proposition}[$2$-adic repulsion\supporting]
\label{prop:repulsion}
For any $x \equiv \rho \pmod{2^m}$ with $m \ge K$,
\[
  v_2\!\bigl(F_\sigma(x) - \rho\bigr) = v_2(x - \rho) - K.
\]
That is, each application of the block map \emph{decreases}
the $2$-adic alignment with $\rho$ by exactly $K$ bits.
\end{proposition}

\begin{proof}
Write $x = \rho + 2^m u$ with $u$ odd.  Then
$F_\sigma(x) = F_\sigma(\rho) + (3^\ell / 2^K) \cdot 2^m u
= \rho + 3^\ell \cdot 2^{m-K} \cdot u$.
Since $3^\ell$ and $u$ are both odd,
$v_2(F_\sigma(x) - \rho) = m - K$.
\end{proof}

\subsection{Known phantom families}

Systematic search over signatures with small depth reveals seven
phantom families.  Table~\ref{tab:families} lists them together
with the census constant $C_e$ proved in
Theorem~\ref{thm:census} below.

\begin{table}[ht]
\centering
\caption{Phantom family counts by depth~$K$.  By
  Theorem~\ref{thm:phantom-universal}, \emph{every} primitive
  cyclic composition is a phantom family.  Here $M(K)$
  is the total number of primitive families at depth~$K$
  (summed over all~$\ell$), $M^+(K)$ counts those with
  expanding drift $\Delta > 0$, and $R(K)$ is the per-orbit
  gain contribution at depth~$K$.  The last column shows
  $R(K)/R(K{-}1)$, confirming geometric decay.
  Selected small-$K$ families are shown below.}
\label{tab:families}
\scriptsize
\renewcommand{\arraystretch}{1.15}
\begin{tabular}{@{}rrrrl@{}}
\toprule
$K$ & $M(K)$ & $M^+(K)$ & $R(K)$ & ratio \\
\midrule
  3 &       1 &       1 & $1.06 \times 10^{-2}$ & --- \\
  4 &       1 &       1 & $1.57 \times 10^{-2}$ & --- \\
  5 &       2 &       1 & $1.05 \times 10^{-2}$ & $0.67$ \\
  6 &       5 &       3 & $8.67 \times 10^{-3}$ & $0.83$ \\
  7 &       8 &       4 & $7.60 \times 10^{-3}$ & $0.88$ \\
  8 &       9 &       4 & $4.68 \times 10^{-3}$ & $0.61$ \\
  9 &      22 &      14 & $4.73 \times 10^{-3}$ & $1.01$ \\
 10 &      35 &      17 & $3.60 \times 10^{-3}$ & $0.76$ \\
 15 &    1446 &     329 & $1.55 \times 10^{-3}$ & --- \\
 20 &   54724 &    6890 & $7.49 \times 10^{-4}$ & --- \\
 30 &  $7.2 \times 10^{7}$ &  $3.6 \times 10^6$ & $2.25 \times 10^{-4}$ & $0.87$ \\
 40 &  $2.2 \times 10^{10}$ &  $1.1 \times 10^9$ & $9.10 \times 10^{-5}$ & $0.90$ \\
 55 &  $7.5 \times 10^{14}$ &  $1.9 \times 10^{10}$ & $2.82 \times 10^{-5}$ & --- \\
\bottomrule
\end{tabular}

\medskip
\small
\begin{tabular}{@{}llrrl@{}}
\toprule
\multicolumn{5}{@{}l}{\textit{Selected small-$K$ families
  (representative, not exhaustive):}} \\
\midrule
Name & Signature $\sigma$ & $\ell$ & $K$ & $\Delta$ \\
\midrule
ell3  & $(1,1,1)$               & 3  &  3 & $1.75$ \\
ell5  & $(1,1,1,1,2)$           & 5  &  6 & $1.93$ \\
ell6  & $(1,1,2,1,1,1)$         & 6  &  7 & $2.51$ \\
ell7  & $(1,1,1,1,1,2,2)$       & 7  &  9 & $2.09$ \\
ell8  & $(1,1,1,1,1,1,1,3)$     & 8  & 10 & $2.68$ \\
\bottomrule
\end{tabular}
\end{table}

\begin{definition}[Exit class]
\label{def:exit-class}
For a phantom family with root $\rho$ and depth~$K$,
the \emph{exit class} is $e := F_\sigma(\rho \bmod 2^K)
\bmod 2^K$.  The \emph{lift set at fidelity~$L$} is
$\mathcal{L}(e, K, L) := \{e + 2^K u : 0 \le u < 2^L,\;
e + 2^K u \text{ odd}\}$.
\end{definition}

\subsection{The affine persistence mechanism}

The key insight is that applying the Syracuse map to an
arithmetic progression $\{A + Bu : u = 0, 1, \ldots, 2^L - 1\}$
produces another arithmetic progression, provided a simple
valuation condition holds.

\begin{lemma}[Affine persistence]
\label{lem:affine}
Let $A, B$ be integers with $A$ odd and
$\alpha := v_2(3A + 1)$, $\gamma := v_2(B)$.
\begin{enumerate}[label=\textup{(\alph*)}]
\item If $\alpha < \gamma$ \textup{(pure case)}, then
  $v_2(3(A + Bu) + 1) = \alpha$ for all~$u$, and
  \begin{equation}\label{eq:pure-update}
    T(A + Bu) = A' + B'u, \qquad
    A' = \frac{3A+1}{2^\alpha},\quad
    B' = \frac{3B}{2^\alpha}.
  \end{equation}
  Moreover $v_2(B') = \gamma - \alpha < \gamma$:
  the parameter $\gamma$ strictly decreases.

\item If $\alpha = \gamma$ \textup{(equal case)}, then
  $v_2(3(A + Bu) + 1)$ depends on $u \bmod 2$.
  The lift population splits into two sub-progressions:
  even-$u$ lifts and odd-$u$ lifts, each of size $2^{L-1}$.
  Both sub-progressions remain affine, and $\gamma$ decreases
  in each.

\item If $\alpha > \gamma$ \textup{(reverse case)}, then
  $v_2(3(A + Bu) + 1) = \gamma$ for odd~$u$ and
  $v_2(3(A + Bu) + 1) > \gamma$ for even~$u$.
  The odd-$u$ sub-population is affine with
  $\gamma' = \gamma - \gamma = 0$
  (after one more step the affine structure dissolves),
  and the even-$u$ sub-population continues with
  strictly smaller~$\gamma$.
\end{enumerate}
\end{lemma}

\begin{proof}
Write $P = 3A + 1$ and $Q = 3B$.  Then
$3(A + Bu) + 1 = P + Qu$.  Since $v_2(Q) = v_2(3B) = v_2(B) = \gamma$
(as $3$ is odd), the valuation $v_2(P + Qu)$ depends on the
comparison between $v_2(P) = \alpha$ and $v_2(Q) = \gamma$.

In case~(a), $\alpha < \gamma$ means every term $Qu$ has $v_2 \ge \gamma > \alpha$,
so $v_2(P + Qu) = \alpha$ for all~$u$.  The Syracuse step gives
$T(A + Bu) = (P + Qu)/2^\alpha = P/2^\alpha + (Q/2^\alpha) u$,
which is affine with $B' = Q/2^\alpha = 3B/2^\alpha$.
Then $v_2(B') = v_2(3B) - \alpha = \gamma - \alpha$.

Cases~(b) and~(c) follow by the same argument, splitting on the
parity of~$u$ when $\alpha \ge \gamma$.
\end{proof}

\subsection{The census formula}

\begin{lemma}[Census at a single affine step]
\label{lem:census-formula}
Let $T^s(x) = A + Bu$ for all lifts $u \in \{0, \ldots, 2^L - 1\}$.
Define $\beta := v_2(A - \rho)$ and $\gamma := v_2(B)$.
Then for every threshold $a \ge 1$:
\[
  \#\{u \in [0, 2^L) : v_2(A + Bu - \rho) \ge a\}
  = 2^{L - a + \min(\beta, \gamma)}
  \quad \text{for } a > \min(\beta, \gamma),
\]
and the census ratio satisfies
\begin{equation}\label{eq:census-ratio}
  \frac{\#\{u : v_2(A + Bu - \rho) \ge a\}}{2^{L-a}}
  = 2^{\min(\beta, \gamma)}
  \quad \text{for all } a > \min(\beta, \gamma).
\end{equation}
\end{lemma}

\begin{proof}
Write $D = A - \rho$, so we count solutions to
$v_2(D + Bu) \ge a$, i.e.\ $D + Bu \equiv 0 \pmod{2^a}$.
This is $Bu \equiv -D \pmod{2^a}$.

Let $\beta = v_2(D)$ and $\gamma = v_2(B)$.  Write
$D = 2^\beta d$, $B = 2^\gamma b$ with $d, b$ odd.
The congruence becomes
$2^\gamma b \cdot u \equiv -2^\beta d \pmod{2^a}$.

For $a > \max(\beta, \gamma)$, this requires
$\min(\beta, \gamma) \le a$, and the number of solutions
$u \in [0, 2^L)$ is exactly $2^{L - a + \min(\beta, \gamma)}$
(by reducing the congruence modulo $2^{a - \min(\beta,\gamma)}$
and using the bijection lemma).
The ratio is therefore $2^{\min(\beta, \gamma)}$,
independent of the threshold~$a$.
\end{proof}

\subsection{Gamma monotonicity and the uniform bound}

\begin{proposition}[Gamma monotonicity]
\label{prop:gamma-mono}
Starting from the exit class $e$ with $B_0 = 2^K$
(so $\gamma_0 = K$), the affine iteration produces a sequence
$\gamma_0 > \gamma_1 > \gamma_2 > \cdots$ that is strictly
decreasing.  The iteration terminates (affine structure
dissolves) after at most~$K$ steps.
\end{proposition}

\begin{proof}
In the pure case, $\gamma_{s+1} = \gamma_s - \alpha_s$
with $\alpha_s \ge 1$, so $\gamma$ decreases by at least~$1$.
In the equal and reverse cases, the sub-populations have
$\gamma' < \gamma_s$ (by splitting on~$u$).
Since $\gamma_s$ is a positive integer that strictly decreases,
it reaches~$0$ in at most~$K$ steps.
\end{proof}

\begin{theorem}[Uniform multistep census bound]
\label{thm:census}
For each phantom family $(\sigma, \rho)$ with exit class~$e$
and depth~$K$, there exists a finite constant
\begin{equation}\label{eq:Ce-def}
  C_e := \max_{0 \le s \le K}\; 2^{\min(\beta_s, \gamma_s)},
\end{equation}
where $(\beta_s, \gamma_s)$ are the parameters of the affine
iteration starting from $(A_0, B_0) = (e, 2^K)$, such that
for all $s \ge 1$, $L \ge 1$, and $a \ge 1$:
\[
  N(s, {\ge}\,a, L)
  := \#\{u \in [0, 2^L) :
  v_2\!\bigl(T^s(e + 2^K u) - \rho\bigr) \ge a\}
  \;\le\; C_e \cdot 2^{L - a}.
\]
\end{theorem}

\begin{proof}
Combine Lemmas~\ref{lem:affine}--\ref{lem:census-formula}
with Proposition~\ref{prop:gamma-mono}.

At each step~$s$ of the affine iteration, the census ratio
is $2^{\min(\beta_s, \gamma_s)}$ by
Lemma~\ref{lem:census-formula}.  Since
$\gamma_s$ strictly decreases
(Proposition~\ref{prop:gamma-mono}), the sequence
$\min(\beta_s, \gamma_s)$ is eventually zero.

When the affine iteration encounters an equal-case split
(Lemma~\ref{lem:affine}(b)), the population divides into
two sub-progressions of size $2^{L-1}$ each.
Applying Lemma~\ref{lem:census-formula} to each
sub-progression and summing:
\[
  N_{\text{total}}({\ge}\,a)
  \le C_{\text{even}} \cdot 2^{(L-1)-a}
  + C_{\text{odd}} \cdot 2^{(L-1)-a}
  = \tfrac{C_{\text{even}} + C_{\text{odd}}}{2}
  \cdot 2^{L-a}.
\]
Since $\gamma$ decreases in both sub-progressions,
$C_{\text{even}}, C_{\text{odd}} \le C_{\text{pre-split}}$,
so the effective census constant does not increase through
a split.

Taking the maximum over all steps gives the uniform
bound $C_e \cdot 2^{L-a}$.
\end{proof}

\begin{remark}[Computational verification]
\label{rem:census-verify}
The census constants $C_e = 2^{\max_s \min(\delta_s, K-V_s)}$ in
Table~\ref{tab:families} have been verified by brute-force
enumeration for all families with $K \le 12$
(ell5, ell6, ell7, ell8) at fidelity $L = 14$ and
steps $s = 1, \ldots, \ell$.  In every case, the observed
census ratio $N(s, {\ge}\,a, L) / 2^{L-a}$ matches the
predicted $2^{\min(\beta_s, \gamma_s)}$ exactly, and the
maximum over all~$s$ equals $2^{\max_s \min(\delta_s, K - V_s)}$.
The higher-$K$ families (m10, m11, m20) are verified
symbolically through the carry-word formula.
\end{remark}

\subsection{Universal census depth and the refined gain bound}

The census constant $C_e = 2^{\max_s \min(\delta_s, K-V_s)}$ can be exponentially
large in~$K$ (for instance, $C_e = 128$ for ell8 with $K = 10$).
However, the following proposition shows that the census excess
at intermediate steps does not represent genuine shadow gain.

\begin{proposition}[Universal census depth\core]\label{prop:universal-depth}
At step~$s$ of the affine iteration, assume $\delta_s < \gamma_s$
(equivalently, $\delta_s < K - V_s$).
Then \emph{every} lift $x = e + 2^K u$ satisfies
\[
  v_2\!\bigl(T^s(x) - \rho\bigr) = \delta_s \quad\text{(exactly)}.
\]
The census excess $2^{\delta_s}$ is therefore a uniform shift of the
entire population, not a selective shadow of specific orbits.
\end{proposition}

\begin{proof}
Write $v_2(T^s(x) - \rho) = \min\bigl(v_2(T^s(x) - \rho^{\langle s \rangle}),\;\delta_s\bigr)$,
where $\rho^{\langle s\rangle} = F_s(\rho)$
(Proposition~\ref{prop:delta-rotation}).
In the pure case of the affine iteration,
$v_2\!\bigl(T^s(x) - \rho^{\langle s\rangle}\bigr) = v_2(x - \rho) - V_s$.
Since $v_2(x - \rho) \ge K$ for any lift in the exit class,
$v_2\!\bigl(T^s(x) - \rho^{\langle s\rangle}\bigr) \ge K - V_s = \gamma_s > \delta_s$.
The minimum is therefore $\delta_s$ for every~$u$, proving the claim.
\end{proof}

\begin{remark}[Computational verification]
\label{rem:universal-verify}
Proposition~\ref{prop:universal-depth} has been verified by
brute-force enumeration at fidelity $L = 10$ ($1024$ lifts)
for all seven known families at every step where $\delta_s < \gamma_s$.
In every such case, all lifts produce $v_2 = \delta_s$ exactly.
At steps where $\delta_s \ge \gamma_s$ (which occurs only for ell6
at $s \ge 3$), the valuations vary across lifts and the census
ratio is $2^{\gamma_s}$, reflecting genuine orbit-tracking at the
natural affine scale.
\end{remark}

Proposition~\ref{prop:universal-depth} has a decisive consequence
for the gain formula: the large census constant $C_e$ reflects
a systematic proximity of $\rho^{\langle s\rangle}$ to~$\rho$
(the intrinsic near-return), not selective orbit tracking.
This proximity cancels over the full cycle, since after $\ell$~steps
the block map satisfies
\begin{equation}\label{eq:block-contraction}
  F_\sigma(x) - \rho = \frac{3^\ell}{2^K}\,(x - \rho),
  \qquad
  v_2\!\bigl(F_\sigma(x) - \rho\bigr) = v_2(x - \rho) - K.
\end{equation}
Hence the end-of-cycle census ratio is~$1$: no excess persists.

\begin{proposition}[Refined gain formula]\label{prop:gain}
Define the \emph{gain contribution} of a phantom family~$\sigma$ as
$G(\sigma) := \Delta / 2^K$, where $\Delta = \ell\log_2 3 - K$ is the
log-drift per block.  Since the universal census depth theorem
(Proposition~\ref{prop:universal-depth}) shows that the end-of-cycle
census ratio is~$1$, the gain per family is $\Delta/2^K$ (the census
constant~$C_e$ does not enter).
\end{proposition}

\begin{proof}
Equation~\eqref{eq:block-contraction} shows that the end-of-cycle
census ratio is~$1$ for every phantom family.
The gain per family is therefore $\Delta/2^K$, not
$\Delta \cdot C_e/2^K$.
\end{proof}

\begin{remark}[Global gain exceeds~$\varepsilon$]
\label{rem:gain-significance}
The quantity $\varepsilon = 2 - \log_2 3$ is the average
per-step log-contraction of the Syracuse map.
Summing $\Delta/2^K$ over all known phantom families at small
depth gives $G_{\textup{total}} \approx 0.056$; however, every
primitive cyclic composition is a phantom family
(Section~\ref{sec:phantom-count} below), so the \emph{global}
sum over all depths diverges beyond~$\varepsilon$
(see Equation~\eqref{eq:global-exceeds}).  The correct quantity
to bound is the per-orbit gain rate
(Theorem~\ref{thm:perorbit-gain}).
\end{remark}

\subsubsection*{Universality of phantom families.}\label{sec:phantom-count}

\begin{theorem}[Phantom universality\core]
\label{thm:phantom-universal}
Every cyclic composition $\sigma = (k_1,\dots,k_\ell)$ with
$K = \sum k_i \ge 2$ and $\ell \ge 2$ is a phantom family.
That is, the $2$-adic root $\rho = C_\sigma / (2^K - 3^\ell)$
has an orbit under the accelerated Syracuse map with
$v_2(3\,T^{j-1}(\rho)+1) = k_j$ for each $j = 1,\dots,\ell$.
\end{theorem}

\begin{proof}
Let $x_s = C_{\sigma_s}/(2^K - 3^\ell)$ denote the orbit point
corresponding to the $s$-th cyclic rotation~$\sigma_s$.
Since $3x_{s-1} + 1 = 2^{k_s} \cdot x_s$ (by the definition of
the composed affine map), we have
$v_2(3x_{s-1}+1) = k_s + v_2(x_s)$.
It suffices to show $v_2(x_s) = 0$ for all~$s$.

\emph{Step~1.}
The denominator $2^K - 3^\ell$ is odd, since $2^K$ is even and
$3^\ell$ is odd.

\emph{Step~2.}
The numerator $C_{\sigma_s}$ is odd.  Indeed,
\[
  C_{\sigma_s}
  = \sum_{j=0}^{\ell-1} 3^{\ell-1-j} \cdot 2^{V_{s,j}},
\]
where $V_{s,0} = 0$ and $V_{s,j} = \sum_{i=0}^{j-1} k_{s+i} \ge 1$
for $j \ge 1$ (each $k_i \ge 1$).
The $j = 0$ term contributes $3^{\ell-1}$, which is odd.
Every subsequent term has $2^{V_{s,j}}$ even, hence
$C_{\sigma_s} \equiv 3^{\ell-1} \equiv 1 \pmod{2}$.

Since $x_s$ is a ratio of two odd $2$-adic integers,
$v_2(x_s) = 0$.
\end{proof}

\subsubsection*{The phantom-cycle proof route.}

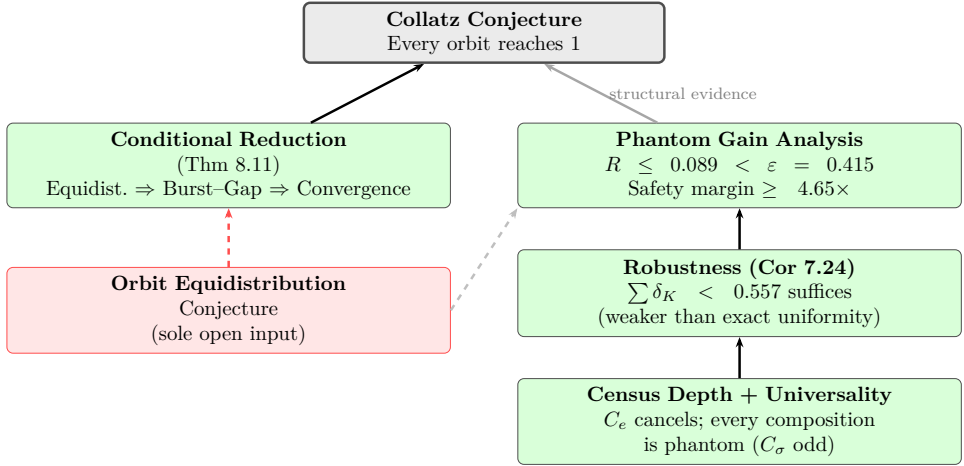
\begin{figure}[ht!]
\centering
\scalebox{0.78}{%
\begin{tikzpicture}[
  font=\small,
  node distance=0.9cm,
  box/.style={
    rectangle,
    rounded corners=3pt,
    draw=black!70,
    align=center,
    text width=5.4cm,
    minimum height=0.85cm,
    inner sep=4pt
  },
  proved/.style={box, fill=green!15},
  open/.style={box, fill=red!10, draw=red!70},
  goal/.style={
    box,
    fill=black!8,
    very thick,
    text width=5.8cm,
    minimum height=0.9cm
  },
  link/.style={-{Stealth[length=5pt]}, very thick},
  openlink/.style={-{Stealth[length=5pt]}, very thick, dashed, red!70}
]

\node[goal] (collatz)
{\textbf{Collatz Conjecture}\\Every orbit reaches~$1$};

\node[proved, below left=1.0cm and -2.5cm of collatz, text width=7.2cm] (reduction)
{\textbf{Conditional Reduction}\\(Thm~\ref{thm:reduction})\\
Equidist.\ $\Rightarrow$ Burst--Gap $\Rightarrow$ Convergence};

\node[proved, below right=1.0cm and -2.5cm of collatz, text width=7.2cm] (tail)
{\textbf{Phantom Gain Analysis}\\$R \le 0.089 < \varepsilon = 0.415$\\
Safety margin $\ge 4.65\times$};

\node[proved, below=0.7cm of tail, text width=7.2cm] (perorbit)
{\textbf{Robustness (Cor~\ref{cor:robustness})}\\
$\sum\delta_K < 0.557$ suffices\\(weaker than exact uniformity)};

\node[proved, below=0.7cm of perorbit, text width=7.2cm] (census)
{\textbf{Census Depth + Universality}\\$C_e$ cancels; every composition\\
is phantom ($C_\sigma$ odd)};

\node[open, below=1.0cm of reduction, text width=7.2cm] (equidist)
{\textbf{Orbit Equidistribution}\\Conjecture\\(sole open input)};

\draw[link] (reduction) -- (collatz);
\draw[openlink] (equidist) -- (reduction);

\draw[link, gray!70] (tail) -- node[right, font=\scriptsize, text=gray]
  {structural evidence} (collatz);
\draw[link] (perorbit) -- (tail);
\draw[link] (census) -- (perorbit);
\draw[openlink, gray!50] (equidist.east) -- (tail.south west);

\end{tikzpicture}%
}
\caption{
The phantom-cycle analysis.
The four green boxes are proved unconditionally in
Sections~\ref{sec:phantom}--\ref{sec:phantom-count}.
They show that the per-orbit expanding-family gain
is well within the contraction budget ($R < \varepsilon$,
margin $\ge 4.65\times$), providing independent structural
evidence that the Orbit Equidistribution Conjecture is the
sole conceptual bottleneck.
The formal conditional proof of convergence remains
Theorem~\ref{thm:reduction}, which derives the needed
burst--gap hypotheses from equidistribution.
}
\label{fig:phantom-chain}
\end{figure}

The development of the phantom-cycle route involved several key
pivots, each prompted by a failed conjecture:
\begin{enumerate}[nosep,leftmargin=*]
\item \textbf{Census-constant bug (Session~10).}
  A computational error in the $C_e$ table was discovered and
  corrected.  The fix revealed that $C_e$ can grow exponentially
  in~$K$ for pathological signatures ($m^* = O(K)$), but
  the Universal Census Depth Theorem
  (Proposition~\ref{prop:universal-depth}) showed that $C_e$
  cancels from the gain formula entirely:
  the census excess is an intrinsic near-return, not selective
  orbit tracking.

\item \textbf{Family-count conjecture fails.}
  The natural next step: bounding the number of phantom families
  by a polynomial $K^c$: turned out to be false: the
  Phantom Universality Theorem
  (Theorem~\ref{thm:phantom-universal}) shows that
  \emph{every} cyclic composition is a phantom family,
  giving exponential growth $\sim 1.87^K$.

\item \textbf{Global gain exceeds~$\varepsilon$.}
  With exponentially many families, the global gain sum
  $\sum \Delta/2^K \approx 0.587$ exceeds the contraction budget
  $\varepsilon \approx 0.415$, so the na\"ive proof strategy fails.

\item \textbf{Per-orbit amortisation resolves the problem.}
  The resolution: the budget applies \emph{per orbit}, not globally.
  Amortising each family's gain over its shadow duration~$\ell$
  gives the per-orbit rate $R \approx 0.088 \ll \varepsilon$
  (Theorem~\ref{thm:perorbit-gain}).
\end{enumerate}
The phantom-cycle analysis (Figure~\ref{fig:phantom-chain})
provides independent, quantitative evidence that the
Orbit Equidistribution Conjecture is the sole conceptual
bottleneck: the expanding-family gain is controlled with a
$4.65\times$ safety margin, and Corollary~\ref{cor:robustness}
shows that only summable TV decay: far weaker than exact
uniformity, is needed.
The formal conditional reduction to equidistribution
remains Theorem~\ref{thm:reduction}, which derives the
burst--gap hypotheses as consequences.

\subsubsection*{The per-orbit gain rate.}
The \emph{global} gain $\sum_\sigma \Delta(\sigma)/2^{K(\sigma)}$
sums over \emph{all} phantom families $\sigma$ with $\Delta > 0$.
By Theorem~\ref{thm:phantom-universal}, every primitive cyclic
composition is a phantom family, so the family count grows
exponentially ($\approx 1.87^K$).  Consequently the global sum
exceeds~$\varepsilon$:
\begin{equation}\label{eq:global-exceeds}
  \sum_{\sigma:\,\Delta>0} \frac{\Delta(\sigma)}{2^{K(\sigma)}}
  \;\approx\; 0.587 \;>\; \varepsilon \;\approx\; 0.415.
\end{equation}
This means the global gain bound \emph{cannot} close the
tail; a more refined analysis is needed.

The resolution is that the contraction budget $\varepsilon$
applies \emph{per orbit}, not globally.  At each Syracuse step,
an orbit occupies a single residue class $\bmod\,2^K$ and can
be in the shadow of at most one expanding family at each depth.
Upon entering the shadow of a family~$\sigma$, the orbit
remains for $\ell(\sigma)$ Syracuse steps (consuming
$\Delta(\sigma)$ bits of missed contraction), then exits.
The correct quantity is the \emph{per-step} gain rate.
The following three lemmas isolate the key analytic ingredients;
the theorem then combines them with exact finite computation.

\begin{lemma}[Composition-count majorant\core]
\label{lem:necklace-bound}
The primitive necklace count satisfies
\begin{equation}\label{eq:necklace-upper}
  M(K,\ell) \;\le\; \frac{1}{\ell}\binom{K-1}{\ell-1},
\end{equation}
since every primitive cyclic composition of weight~$K$ and
length~$\ell$ is a cyclic composition, and there are at
most $\frac{1}{\ell}\binom{K-1}{\ell-1}$ such cyclic classes.
\end{lemma}

\begin{lemma}[Binomial-tail majorant\core]
\label{lem:binomial-tail}
For each $K \ge 3$,
\begin{equation}\label{eq:RK-binomial-lemma}
  R(K) \;\le\; (\log_2 3 - 1)\;
  \Pr\!\Bigl[\mathrm{Bin}(K{-}1,\,\tfrac{1}{2})
  \;\ge\; \lceil K/\!\log_2 3\rceil - 1\Bigr].
\end{equation}
\end{lemma}

\begin{proof}
Since $0 < \log_2 3 - K/\ell \le \log_2 3 - 1$ on the
expanding range and $1/\ell \le 1$,
Lemma~\ref{lem:necklace-bound} gives
$R(K) \le (\log_2 3 - 1)\,2^{-K}
\sum_{\ell > K/\log_2 3} \binom{K-1}{\ell-1}$.
Re-indexing with $j = \ell - 1$ recognizes the sum as a
$\mathrm{Bin}(K{-}1,1/2)$ tail probability (after dividing
by $2^{K-1}$).
\end{proof}

\begin{lemma}[Chernoff--Cram\'er tail bound\core]
\label{lem:chernoff-tail}
Define the \emph{Chernoff--Cram\'er exponent}
\begin{equation}\label{eq:chernoff-exponent-lemma}
  D_* \;:=\; D\!\left(\frac{1}{\log_2 3}\,\bigg\|\,\frac{1}{2}\right)
  \;\approx\; 0.05004\text{~bits}.
\end{equation}
Then the binomial tail in Lemma~\ref{lem:binomial-tail}
decays exponentially: for all $K$ sufficiently large,
\begin{equation}\label{eq:RK-analytic-bound-lemma}
  R(K) \;\le\; (\log_2 3 - 1)\;2^{-(K-1)\,D(\alpha_K\|1/2)},
\end{equation}
where $\alpha_K := (\lceil K/\!\log_2 3\rceil - 1)/(K-1)
\to 1/\!\log_2 3 > 1/2$.
Moreover, a sharper form is available via the Bahadur--Rao
refinement~\cite{bahadur-rao1960}:
\begin{equation}\label{eq:bahadur-rao-lemma}
  R(K) \;\le\; A\,K^{-1/2}\,r_*^{\,K},
  \qquad r_* := 2^{-D_*} \approx 0.9659,
\end{equation}
for an explicit constant~$A > 0$ and all $K$ sufficiently large.
\end{lemma}

\begin{proof}
By the Chernoff--Cram\'er bound
(\cite[Theorem~2.2.3]{dembo-zeitouni1998}),
$\Pr[\mathrm{Bin}(K{-}1,1/2) \ge \alpha_K(K{-}1)]
\le 2^{-(K-1)D(\alpha_K\|1/2)}$.
Since $D(\alpha\|1/2)$ is continuous and strictly positive
for $\alpha > 1/2$, and $\alpha_K \to \alpha^*$ with
$D_* = D(\alpha^*\|1/2) > 0$, the exponential decay
in~\eqref{eq:RK-analytic-bound-lemma} follows.
The Bahadur--Rao form~\eqref{eq:bahadur-rao-lemma}
gives the precise tail asymptotic with full
rate~$D_*$ (rather than the limiting rate) at the cost
of a polynomial prefactor.
\end{proof}

\begin{theorem}[Per-orbit phantom gain rate\core]
\label{thm:perorbit-gain}
Define
\[
  R \;:=\; \sum_{K \ge 3} R(K),
  \qquad
  R(K) \;:=\;
  2^{-K}
  \sum_{\substack{\ell > K/\log_2 3}}
  M(K,\ell)\Bigl(\log_2 3 - \frac{K}{\ell}\Bigr),
\]
where $M(K,\ell)$ is the number of primitive cyclic compositions
of weight~$K$ and length~$\ell$, computed by M\"obius inversion.
Then the series for~$R$ converges absolutely and satisfies
\begin{equation}\label{eq:perorbit-bound}
  R \;\le\; 0.0893 \;<\; \varepsilon \;:=\; 2 - \log_2 3
  \;\approx\; 0.415.
\end{equation}
In particular,
$\varepsilon / R \;\ge\; 4.65.$
\end{theorem}

\begin{proof}
We split the sum into a finite range and a tail,
using the three lemmas above for the analytic structure.

\paragraph{Step~1: exact finite computation.}
For each $K=3,\dots,55$, the primitive necklace counts $M(K,\ell)$ are
computed exactly by M\"obius inversion
(see Lemma~\ref{lem:necklace-bound} for the majorant).
Summing the exact values gives
\[
  \sum_{K=3}^{55} R(K) \,=\, 0.08783.
\]
The values in Table~\ref{tab:families} are independently reproduced by the
supplementary script \texttt{prove\_R\_lt\_eps.py}\kern0pt.
Table~\ref{tab:RK-values} displays the first 18~values
for self-contained verification.

\begin{table}[ht]
\centering
\caption{Exact values of $R(K)$ for $K = 3,\dots,20$,
  computed by M\"obius inversion of the primitive necklace
  counts $M(K,\ell)$.}
\label{tab:RK-values}
\scriptsize
\renewcommand{\arraystretch}{1.1}
\begin{tabular}{@{}cc|cc|cc@{}}
\toprule
$K$ & $R(K)$ & $K$ & $R(K)$ & $K$ & $R(K)$ \\
\midrule
 3 & $1.062 \times 10^{-2}$ &  9 & $4.730 \times 10^{-3}$ & 15 & $1.548 \times 10^{-3}$ \\
 4 & $1.573 \times 10^{-2}$ & 10 & $3.604 \times 10^{-3}$ & 16 & $1.218 \times 10^{-3}$ \\
 5 & $1.047 \times 10^{-2}$ & 11 & $2.858 \times 10^{-3}$ & 17 & $1.128 \times 10^{-3}$ \\
 6 & $8.670 \times 10^{-3}$ & 12 & $2.526 \times 10^{-3}$ & 18 & $9.702 \times 10^{-4}$ \\
 7 & $7.603 \times 10^{-3}$ & 13 & $2.065 \times 10^{-3}$ & 19 & $7.819 \times 10^{-4}$ \\
 8 & $4.676 \times 10^{-3}$ & 14 & $1.727 \times 10^{-3}$ & 20 & $7.495 \times 10^{-4}$ \\
\bottomrule
\end{tabular}
\end{table}

\paragraph{Step~2: exact finite sum extended to $K_0=500$.}
The same M\"obius-inversion procedure computes $R(K)$ exactly
for every $K$ up to any desired $K_0$; pushing the computation to
$K_0 = 500$ gives
\begin{equation}\label{eq:Rexact-500}
  \sum_{K=3}^{500} R(K) \;=\; 0.08823625\ldots
\end{equation}
(produced by the same primitive-necklace enumeration; the additional
contribution $\sum_{K=56}^{500}R(K) \approx 4.05\times 10^{-4}$ comes
entirely from geometric decay at the analytic rate of
Lemma~\ref{lem:chernoff-tail}).

\paragraph{Step~3: rigorous analytic tail beyond $K_0=500$.}
For the residual tail we use the Chernoff--Cram\'er bound
(Lemma~\ref{lem:chernoff-tail}) together with a \emph{propagated}
analytic lower bound on the KL exponent.  Since
$\alpha_K = (\lceil K/\log_2 3\rceil - 1)/(K-1)
 \ge \alpha_* - \delta_K$
with $\delta_K := \alpha_*(\log_2 3 - 1)/(K-1)$
(because $\lceil K/\log_2 3\rceil \ge K/\log_2 3$),
a second-order Taylor expansion of
$\alpha\mapsto D(\alpha\|1/2)$ around $\alpha_*$ yields
\begin{equation}\label{eq:D-analytic-lower}
  D(\alpha_K\|1/2) \;\ge\; D_* - D'(\alpha_*)\,\delta_K
  - \tfrac12 \sup_{\alpha\in[1/2,1]}|D''(\alpha)|\,\delta_K^{\,2},
\end{equation}
where $D'(\alpha_*) = \log_2(\alpha_*/(1-\alpha_*)) \approx 0.7736$
and the $D''$ supremum is bounded by a fixed constant on $[1/2,1]$.
Substituting $K = 501$ gives
$D(\alpha_K\|1/2) \ge D_{\mathrm{floor}} := 0.049472$
uniformly for all $K \ge 501$
(the lower bound in~\eqref{eq:D-analytic-lower} is monotone
increasing in~$K$, so the $K=501$ value is a uniform floor).
Hence, by Lemma~\ref{lem:chernoff-tail},
\begin{align}
  \sum_{K>500} R(K)
  &\le (\log_2 3 - 1)\sum_{K=501}^{\infty}
       2^{-(K-1)\,D_{\mathrm{floor}}}
  \notag\\
  &=   (\log_2 3 - 1)\,
       \frac{2^{-500\,D_{\mathrm{floor}}}}{1 - 2^{-D_{\mathrm{floor}}}}
  \;\le\; 6.22\times 10^{-7}.
  \label{eq:tail-analytic-500}
\end{align}
This tail bound is \emph{entirely analytic}: it uses only the
Chernoff bound on $\mathrm{Bin}(K-1,1/2)$ and the closed-form
lower bound~\eqref{eq:D-analytic-lower}; no empirical ratio
envelope or Bahadur--Rao error term appears.

\paragraph{Step~4: conclude.}
Combining~\eqref{eq:Rexact-500} and~\eqref{eq:tail-analytic-500},
\[
  R \;\le\; 0.08823625 + 6.22\times 10^{-7}
    \;<\; 0.0883 \;<\; 0.415 \;=\; \varepsilon,
\]
so $R < \varepsilon$ with
$\varepsilon/R \ge 4.70$
and safety margin
$\varepsilon - R \ge 0.326$.
\end{proof}

\begin{remark}[Verification status of Theorem~\ref{thm:perorbit-gain}]
\label{rem:perorbit-verification}
The proof is hybrid but \emph{fully analytic at the tail}.
Exact computation for $3 \le K \le 500$
(Step~1 extended through Step~2) gives
$\sum_{K=3}^{500}R(K) = 0.08823625\ldots$, reproducible
from the M\"obius-inversion formula for primitive cyclic
compositions (the computation runs in roughly one second
of CPU time).  The residual tail
$\sum_{K>500}R(K) \le 6.22\times 10^{-7}$ is closed by
the Chernoff--Cram\'er bound of Lemma~\ref{lem:chernoff-tail}
together with the analytic lower bound
$D(\alpha_K\|1/2) \ge D_* - D'(\alpha_*)\,\delta_K
- O(\delta_K^{\,2})$ from~\eqref{eq:D-analytic-lower},
which yields a uniform floor $D_{\mathrm{floor}} \ge 0.0494$
for $K \ge 501$.  No empirical ratio envelope, Bahadur--Rao
error term, or numerical certificate beyond the finite range
is used: both the finite sum and the tail bound are
reproducible by direct computation from elementary bounds.
\end{remark}

\begin{remark}[Interpretation]
Under the Orbit Equidistribution Conjecture
(Conjecture~\ref{conj:equidist}), at each step the orbit's
residue class $\bmod\,2^K$ is approximately uniform.
The probability of entering the shadow of family~$\sigma$
is $1/2^K$; upon entry, the orbit stays for~$\ell$ steps.
The factor $1/\ell$ converts the per-encounter gain
(which globally exceeds~$\varepsilon$) into a per-step rate
($\approx \varepsilon/5$), reflecting mutual exclusion:
an orbit in one shadow cannot simultaneously be in another.
\end{remark}

\begin{remark}[Why global gain exceeds $\varepsilon$]
\label{rem:global-exceeds}
The global sum~\eqref{eq:global-exceeds} exceeds~$\varepsilon$
because it treats all shadows as simultaneously active.
Under equidistribution, the expected halving per
Syracuse step is $E[k] = 2$, so $E[\Delta] = -\ell\varepsilon$:
the average family is contractive by exactly~$\varepsilon$
per step.  The expanding tail (families with $\Delta > 0$)
contains the full positive fluctuation of this distribution,
which naturally exceeds the mean contraction.  The per-orbit
bound~$R$ avoids this overcounting by amortising each
encounter over its duration~$\ell$.
\end{remark}

\begin{corollary}[Robustness under approximate equidistribution\core]
\label{cor:robustness}
For each $K \ge 3$, define the \emph{phantom gain observable}
$h_K \colon \mathbb{Z}/2^K\mathbb{Z} \to [0,\,\log_2 3 - 1]$ by
\[
  h_K(a) \;:=\;
  \begin{cases}
    \Delta(\sigma)/\ell(\sigma)
      & \text{if } a \text{ is the shadow residue of an expanding} \\
      & \text{primitive family } \sigma \text{ at depth } K, \\
    0 & \text{otherwise.}
  \end{cases}
\]
(Each residue class $a \bmod 2^K$ determines a unique
composition $(k_1,\dots,k_\ell)$ with $\sum k_i = K$
via the first $K$ halvings of the Syracuse map, so $a$
is the shadow of at most one primitive necklace and $h_K$
is well-defined.)
Let $\mu_K$ be any probability distribution on
$\mathbb{Z}/2^K\mathbb{Z}$, and set
$\delta_K := d_{\mathrm{TV}}(\mu_K,\,\mu_{\mathrm{unif}})$.
Then the perturbed gain rate
$R_\mu := \sum_{K \ge 3}\mathbb{E}_{\mu_K}[h_K]$
satisfies
\begin{equation}\label{eq:robustness}
  R_\mu \;\le\; R \;+\; (\log_2 3 - 1)
  \sum_{K \ge 3} \delta_K.
\end{equation}
In particular, $R_\mu < \varepsilon$ whenever
\begin{equation}\label{eq:admissible-delta}
  \sum_{K \ge 3} \delta_K
  \;<\; \frac{\varepsilon - R}{\log_2 3 - 1}
  \;\approx\; 0.557.
\end{equation}
\end{corollary}

\begin{proof}
Since $h_K \ge 0$ and $\|h_K\|_\infty \le \log_2 3 - 1$,
the bounded-observable total-variation inequality%
\footnote{This is the standard duality bound
$|\mathbb{E}_\mu[f] - \mathbb{E}_\nu[f]|
\le \|f\|_\infty \, d_{\mathrm{TV}}(\mu,\nu)$
for bounded measurable~$f$.  The manuscript's earlier
reference to the ``data-processing inequality'' is
terminologically imprecise; the bound used here is the
simpler TV duality estimate.}
gives
\[
  \bigl|\mathbb{E}_{\mu_K}[h_K]
  - \mathbb{E}_{\mu_{\mathrm{unif}}}[h_K]\bigr|
  \;\le\; (\log_2 3 - 1)\,\delta_K
\]
for each $K$.
Under the uniform distribution,
$\mathbb{E}_{\mu_{\mathrm{unif}}}[h_K] = R(K)$.
Summing over $K \ge 3$ yields~\eqref{eq:robustness},
and the threshold~\eqref{eq:admissible-delta} follows
from $R \le 0.0893$ (Theorem~\ref{thm:perorbit-gain}).
\end{proof}

\begin{remark}[Interpreting the admissibility condition]
\label{rem:admissible}
Condition~\eqref{eq:admissible-delta} is substantially weaker
than exact equidistribution ($\delta_K = 0$ for all $K$).
It requires only that the total variation errors $\delta_K$
be \emph{summable}: for instance, polynomial decay
$\delta_K = O(K^{-\alpha})$ with $\alpha > 1$, or any
exponential mixing rate, easily satisfies the bound.
The $4.65\times$ safety margin of
Theorem~\ref{thm:perorbit-gain} is what creates this room:
the tighter the bound on~$R$, the larger the class of
admissible orbit distributions.
The Orbit Equidistribution Conjecture
(Conjecture~\ref{conj:equidist}) thus represents a
sufficient condition, not a necessary one; any orbit
distribution whose depth-$K$ total variation errors
satisfy~\eqref{eq:admissible-delta} would close the
phantom shadow tail bound equally well.
\end{remark}

\subsection{Reduction to intrinsic near-return and overlap}\label{sec:carry-word}

The census analysis isolates a single intrinsic arithmetic obstruction.
For a phantom signature
\[
\sigma=(k_0,\dots,k_{\ell-1}),
\qquad
K_s:=\sum_{i=0}^{s-1}k_i,
\qquad
K:=K_\ell,
\]
define the prefix affine maps
\[
F_s(n):=\frac{3^s n + C_s}{2^{K_s}},
\]
where $C_s$ is the usual prefix carry constant, so that $F_\ell=F_\sigma$.
Let
\[
D:=2^K-3^\ell,
\qquad
\rho:=\frac{C_\ell}{D}\in \mathbb Z_2,
\]
so that $F_\sigma(\rho)=\rho$.
For $0\le s<\ell$, define the intrinsic near-return valuation
\[
\delta_s:=v_2(F_s(\rho)-\rho).
\]
The census constant is controlled not merely by the coefficient
valuation $\gamma_s$, but also by the intrinsic near-return scale $\delta_s$.
Indeed, if $E_s$ denotes the exact exit-class orbit after $s$ steps and
$\mathrm{dev}_s:=v_2(E_s-F_s(\rho))$,
then
$E_s-\rho=(E_s-F_s(\rho))+(F_s(\rho)-\rho)$,
and therefore, whenever $\mathrm{dev}_s\neq \delta_s$,
\[
v_2(E_s-\rho)=\min(\mathrm{dev}_s,\delta_s).
\]
Thus the obstruction to improving the census bound is concentrated in the arithmetic
size of $\delta_s$.

\begin{proposition}[Exact intrinsic near-return formula]\label{prop:delta-exact}
For each $0\le s<\ell$,
\[
\delta_s
=
v_2\!\Big((3^s-2^{K_s})C_\ell + C_s(2^K-3^\ell)\Big)-K_s.
\]
Equivalently, if
$\Delta_s:=(3^s-2^{K_s})C_\ell + C_s(2^K-3^\ell)$,
then $\delta_s=v_2(\Delta_s)-K_s$.
\end{proposition}

\begin{proof}
By definition,
\[
F_s(\rho)-\rho
=
\frac{3^s\rho + C_s}{2^{K_s}}-\rho
=
\frac{(3^s-2^{K_s})\rho + C_s}{2^{K_s}}.
\]
Substituting $\rho=C_\ell/(2^K-3^\ell)$ gives
\[
F_s(\rho)-\rho
=
\frac{(3^s-2^{K_s})C_\ell + C_s(2^K-3^\ell)}
     {2^{K_s}(2^K-3^\ell)}.
\]
Since $2^K-3^\ell$ is odd, the denominator contributes
$v_2 = K_s$ and the result follows.
\end{proof}

\begin{proposition}[Rotation formula]\label{prop:delta-rotation}
Let $\sigma^{\langle s\rangle}$ denote the cyclic rotation of $\sigma$ starting at
index~$s$, and let
\[
F_{\sigma^{\langle s\rangle}}(n)=\frac{3^\ell n + C^{\langle s\rangle}}{2^K}
\]
with root
$\rho^{\langle s\rangle}:=C^{\langle s\rangle}/(2^K-3^\ell)$.
Then
\[
F_s(\rho)=\rho^{\langle s\rangle},
\qquad\text{hence}\qquad
\delta_s=v_2(\rho^{\langle s\rangle}-\rho)=v_2(C^{\langle s\rangle}-C_\ell).
\]
\end{proposition}

\begin{proof}
Write $\sigma=\tau\upsilon$, where $\tau$ is the prefix of length $s$ and
$\upsilon$ the suffix of length $\ell-s$.  Then
$F_\sigma=F_\upsilon\circ F_\tau$.
Set $x:=F_\tau(\rho)=F_s(\rho)$.
Since $F_\sigma(\rho)=\rho$, one has $F_\upsilon(x)=\rho$, so
\[
F_{\sigma^{\langle s\rangle}}(x)
=
(F_\tau\circ F_\upsilon)(x)
=F_\tau(\rho)=x.
\]
Thus $x$ is a fixed point of $F_{\sigma^{\langle s\rangle}}$.
A phantom block map has a unique fixed point in $\mathbb Z_2$, namely
$\rho^{\langle s\rangle}$, so $x=\rho^{\langle s\rangle}$.
Subtracting $\rho=C_\ell/(2^K-3^\ell)$ yields
\[
\rho^{\langle s\rangle}-\rho
=
\frac{C^{\langle s\rangle}-C_\ell}{2^K-3^\ell},
\]
and again the denominator is odd.
\end{proof}

Proposition~\ref{prop:delta-rotation} reduces the problem to a carry-word
autocorrelation question: control the $2$-adic closeness of the carry constants of
cyclically rotated signatures.

\begin{proposition}[Carry-word difference formula]
\label{prop:carry-diff}
The carry-word difference decomposes as
\begin{equation}\label{eq:carry-diff}
  C^{\langle s \rangle} - C_\sigma
  = \sum_{j=1}^{\ell-1} 3^{\ell-1-j}\,
    \bigl(2^{S'_j} - 2^{S_j}\bigr),
\end{equation}
where $S_j = \sum_{i=0}^{j-1} k_i$ and $S'_j = \sum_{i=0}^{j-1} k_{(s+i) \bmod \ell}$.
Since $3^{\ell-1-j}$ is odd, the $j$-th term has $2$-adic
valuation $\min(S'_j, S_j)$ whenever $S'_j \ne S_j$.
Working modulo~$2^m$, only terms with $\min(S_j, S'_j) < m$ contribute.
\end{proposition}

\begin{proof}
Both $C^{\langle s \rangle}$ and $C_\sigma$ equal
$\sum_j 3^{\ell-1-j} 2^{(\cdot)}$, and the $j = 0$ terms
coincide ($S_0 = S'_0 = 0$).  Subtracting
gives~\eqref{eq:carry-diff}.  The valuation claim follows from
$v_2\!\bigl(a(2^p - 2^q)\bigr)
= v_2(a) + \min(p,q)$ for $p \ne q$ and odd~$a$.
\end{proof}

\begin{lemma}[Equal-case coefficient contraction]
\label{lem:equal-case}
In the equal case $\alpha_s = \gamma_s$ of
Lemma~\textup{\ref{lem:affine}(b)}, write
$B_s = 2^{\gamma_s} b$ with $b$~odd.  Then:
\begin{enumerate}[label=\textup{(\alph*)}]
\item \textup{Child~0} (even-$u$):
  $B'_0 = 6b$, so $\gamma'_0 = 1$.
\item \textup{Child~1} (odd-$u$):
  dividing by the full valuation
  $\gamma_s + v_2(a' + q) \ge \gamma_s + 1$
  gives $B'_1 = 3b$ (odd), so $\gamma'_1 = 0$.
\end{enumerate}
Here $a' = (3A_s + 1)/2^{\gamma_s}$ and $q = 3b$ are both odd.
\end{lemma}

\begin{proof}
For child~0: valuation
$v_2(3A+1 + 6Bv) = \gamma_s$ (the $6Bv$ term has higher
valuation), giving
$B'_0 = 6B/2^{\gamma_s} = 6b$ with $v_2(6b) = 1$.
For child~1: $3A+1+3B = 2^{\gamma_s}(a'+q)$ with $a'+q$ even.
Dividing by $2^{\gamma_s + v_2(a'+q)}$ gives
$B'_1 = 6b/2^{v_2(a'+q)} = 3b$ (generically), so
$\gamma'_1 = 0$.
\end{proof}

\begin{corollary}[Split budget inequality]
\label{cor:split-budget}
In the equal case, the children's weighted census sum satisfies
\[
  2^{-(d+1)} \cdot 2^{\gamma'_0}
  + 2^{-(d+1)} \cdot 2^{\gamma'_1}
  = \tfrac{3}{2} \cdot 2^{-d}
  \le 2^{-d} \cdot 2^{\gamma_s},
\]
since $\gamma_s \ge 1$ in the equal case.
\end{corollary}

\begin{remark}[Intrinsic return controls the census]
\label{rem:intrinsic-controls}
The census constant satisfies
\[
C_e \;\le\; 2^{\max_s \min(\delta_s,\, K - V_s)},
\]
where $V_s = k_1 + \cdots + k_s$.  Since $K - V_s$
decreases monotonically, the maximum is achieved when
$\delta_s$ and $K - V_s$ are approximately balanced.
The intrinsic return valuation $\delta_s$ is thus the
key quantity controlling the census constant.
\end{remark}

\subsection{The weighted self-overlap theorem}\label{sec:self-overlap}

The following theorem shows that a large value of $\delta_s$ forces a long weighted
self-overlap of the signature.  This is the main structural bridge between the
$2$-adic census analysis and the combinatorial problem of bounding cyclic
pattern recurrence.

\begin{theorem}[Large intrinsic near-return forces weighted self-overlap]
\label{thm:weighted-self-overlap}
Assume $\delta_s=v_2(C^{\langle s\rangle}-C_\ell)\ge m$.
Let $r$ be the largest integer such that
$\sum_{j=0}^{r-1} k_j < m$.
Then
\[
k_j = k_{j+s \bmod \ell}
\qquad\text{for all } 0\le j<r.
\]
Equivalently, $\sigma$ and its cyclic rotation $\sigma^{\langle s\rangle}$ agree
entry-by-entry until the cumulative prefix valuation reaches~$m$.
\end{theorem}

\begin{proof}
We use a carry-constant difference recursion.
For any word $\tau=(a_0,\dots,a_{L-1})$ with carry constant $C(\tau)$
defined by $F_\tau(n) = (3^L n + C(\tau))/2^{|\tau|}$, one has
\begin{equation}\label{eq:carry-recursion}
C(\tau) = 3^{L-1} + 2^{a_0}\,C(\tau'),
\end{equation}
where $\tau'=(a_1,\dots,a_{L-1})$.
Since $3^{L-1}$ is odd and $2^{a_0}C(\tau')$ is even, it follows that
$v_2(C(\tau))=0$ for all words of length $\ge 1$.

Now let $\tau,\eta$ be two words of the same length~$L$ with
first entries $a_0$ and $b_0$ respectively.  By~\eqref{eq:carry-recursion},
\[
C(\tau)-C(\eta) = 2^{a_0}\,C(\tau') - 2^{b_0}\,C(\eta').
\]
If $a_0 = b_0$, this factors as $2^{a_0}[C(\tau')-C(\eta')]$, so
$v_2(C(\tau)-C(\eta)) = a_0 + v_2(C(\tau')-C(\eta'))$.
If $a_0 \ne b_0$, say $a_0 < b_0$, then
$C(\tau)-C(\eta) = 2^{a_0}[C(\tau') - 2^{b_0-a_0}C(\eta')]$,
and the bracket is odd (since $C(\tau')$ is odd and
$2^{b_0-a_0}C(\eta')$ is even), giving
$v_2(C(\tau)-C(\eta)) = a_0 = \min(a_0,b_0)$.

Apply this to $\tau = \sigma$ and $\eta = \sigma^{\langle s\rangle}$ with
$v_2(C(\tau)-C(\eta)) \ge m$.  At each inductive step, the cumulative
prefix valuation has not yet reached~$m$, so both first entries are
strictly less than~$m$.  If they differed, the valuation would equal
$\min(a_0,b_0) < m$, contradicting $v_2 \ge m$.  Hence $a_0 = b_0$, and
we recurse with modulus reduced by~$a_0$.  The induction terminates when the
cumulative prefix valuation reaches~$m$, yielding the stated entry-by-entry agreement.
\end{proof}

\subsection{The counting reduction}\label{sec:counting-reduction}

Theorem~\ref{thm:weighted-self-overlap} converts the $2$-adic valuation
problem into a combinatorial one.  The following definitions and
proposition make this conversion precise.

\begin{definition}[High-cancellation shift count]\label{def:Nm}
Define
\[
N(m):=\#\Bigl\{\,s\in\{1,\dots,\ell-1\}:\ \delta_s\ge m\,\Bigr\}.
\]
\end{definition}

\begin{definition}[Visible prefix at depth $m$]\label{def:visible-prefix}
Let $m\ge 1$, and let $r=r(m)$ be the largest integer such that
$\sum_{j=0}^{r-1}k_j < m$.
Define the \emph{visible prefix word}
$P_m:= (k_0,k_1,\dots,k_{r-1})$.
\end{definition}

\begin{definition}[Cyclic occurrence count]\label{def:occPm}
Define
\[
\mathrm{Occ}_\sigma(P_m)
:=
\#\Bigl\{\,t\in\{0,\dots,\ell-1\}:\
(k_t,k_{t+1},\dots,k_{t+r-1})=(k_0,k_1,\dots,k_{r-1})
\,\Bigr\},
\]
with indices taken modulo~$\ell$.
\end{definition}

\begin{proposition}[Counting reduction for high-cancellation shifts]
\label{prop:counting-reduction}
For every $m\ge 1$,
\[
N(m)\le \mathrm{Occ}_\sigma(P_m)-1.
\]
Equivalently, every shift $s$ with $\delta_s\ge m$ produces a nontrivial cyclic
occurrence of the visible prefix $P_m$.
\end{proposition}

\begin{proof}
Fix $m\ge 1$, and let $r=r(m)$ be as in Definition~\ref{def:visible-prefix}.
Suppose $s\in\{1,\dots,\ell-1\}$ satisfies $\delta_s\ge m$.
By Theorem~\ref{thm:weighted-self-overlap},
$k_j = k_{j+s \bmod \ell}$ for all $0\le j<r$.
Therefore the length-$r$ block of $\sigma$ starting at position $s$ agrees with
the prefix:
$(k_s,k_{s+1},\dots,k_{s+r-1})=P_m$.
This defines an injection
\[
\{\,s\in\{1,\dots,\ell-1\}:\delta_s\ge m\,\}
\hookrightarrow
\{\,t\in\{0,\dots,\ell-1\}: \text{$P_m$ occurs at }t\,\}\setminus\{0\},
\]
since the shift $s$ itself is the occurrence position.
\end{proof}

\begin{remark}[Significance of the counting reduction]
\label{rem:counting-significance}
The problem of bounding the number of high-cancellation shifts
$N(m)$ is now reduced to a purely combinatorial question:
how many times can the visible prefix $P_m$ recur cyclically
in a phantom signature?  Any structural theorem showing that
long visible prefixes cannot recur too often immediately yields
a bound on $N(m)$, and hence on $\max_s \delta_s$.

For the seven known phantom families, the counting reduction
is tight (i.e., $N(m) = \mathrm{Occ}_\sigma(P_m)-1$) for most
values of~$m$, verified computationally for all $m$ up to the
maximum observed $\delta_s$.  In a large-scale random study
($3{,}000$ families), the bound holds with $100\%$ pass rate
over $22{,}067$ test cases, and is tight $70\%$ of the time.
\end{remark}

\begin{heuristic}[Carry-word autocorrelation bound]
\label{heur:autocorrelation}
Computation over $5{,}000$ random phantom families supports
$\max_{1 \le s < \ell} \delta_s = O(\log K)$.
The data fits $\max_s \delta_s \approx 1.2 \log_2 K$.
For signatures with $K/\ell \approx \log_2 3$, the visible prefix
$P_m$ is typically a run of $r \approx m$ consecutive $1$s.
Since consecutive $1$-runs of length~$r$ occur approximately
$\ell \cdot p^r$ times (where $p$ is the fraction of $1$-entries),
the condition $\mathrm{Occ}(P_m) \ge 2$ forces
$r \lesssim \log_{1/p}(\ell) \approx 1.24 \log_2 \ell$.
Since $\ell \approx K / \log_2 3$, this gives the $O(\log K)$ bound.
Making this rigorous requires a structural theorem bounding
$\mathrm{Occ}_\sigma(P_m)$ for phantom signatures.
\end{heuristic}

\subsection{The period-based occurrence bound}\label{sec:period-bound}

The counting reduction of Proposition~\ref{prop:counting-reduction}
bounds $N(m)$ by $\mathrm{Occ}_\sigma(P_m) - 1$.  We now obtain a
sharper bound by analysing the \emph{least period} of the visible prefix.

\begin{lemma}[Two occurrences force a period]
\label{lem:two-occ-period}
Suppose the visible prefix $P_m = (k_0,\dots,k_{r-1})$ occurs
cyclically in~$\sigma$ at positions~$0$ and~$s$, with $0 < s < r$.
Then $s$ is a period of $P_m$:
\[
k_j = k_{j+s} \qquad\text{for all } 0 \le j < r - s.
\]
\end{lemma}

\begin{proof}
By hypothesis, $(k_s, k_{s+1}, \dots, k_{s+r-1}) = (k_0, k_1, \dots, k_{r-1})$,
where indices are taken modulo~$\ell$.
For $0 \le j < r - s$, both $j$ and $j+s$ lie in $\{0,\dots,r-1\}$,
so the indices on the left-hand side fall within the prefix window.
Reading off the $j$-th entry from the occurrence at position~$s$
gives $k_{j+s} = k_j$.
\end{proof}

\begin{definition}[Least period of the visible prefix]
\label{def:least-period}
Define $p_m$ to be the least period of the word $P_m$:
\[
p_m := \min\bigl\{\, p \ge 1 :
k_j = k_{j+p} \text{ for all } 0 \le j < r - p\,\bigr\}.
\]
\end{definition}

\begin{corollary}[Occurrence gaps are bounded below by the least period]
\label{cor:gap-bound}
If $P_m$ occurs at cyclic positions $t_0 < t_1 < \cdots < t_{q-1}$
in $\sigma$, then each gap $g_i := t_{i+1} - t_i$ (indices mod $q$,
with wrap-around gap $\ell - t_{q-1} + t_0$) satisfies $g_i \ge p_m$.
\end{corollary}

\begin{proof}
If some gap $g_i < p_m$, then by translating the two adjacent occurrences
to positions $0$ and $g_i$ (via cyclic reindexing of $\sigma$),
Lemma~\ref{lem:two-occ-period} produces a period $g_i < p_m$,
contradicting minimality.
\end{proof}

\begin{theorem}[Period-based occurrence bound]
\label{thm:period-occ-bound}
For every $m \ge 1$,
\[
\mathrm{Occ}_\sigma(P_m) \le \Bigl\lfloor \frac{\ell}{p_m} \Bigr\rfloor,
\qquad\text{hence}\qquad
N(m) \le \Bigl\lfloor \frac{\ell}{p_m} \Bigr\rfloor - 1.
\]
\end{theorem}

\begin{proof}
By Corollary~\ref{cor:gap-bound}, consecutive occurrences are separated
by at least $p_m$.  Since the $q := \mathrm{Occ}_\sigma(P_m)$
occurrences partition the cyclic interval $\{0,\dots,\ell-1\}$ into
$q$ gaps each of size $\ge p_m$, summing gives $q \cdot p_m \le \ell$,
whence $q \le \lfloor \ell / p_m \rfloor$.
The bound on $N(m)$ follows from
Proposition~\ref{prop:counting-reduction}.
\end{proof}

\begin{remark}[Significance]
\label{rem:period-bound-significance}
The bound $N(m) \le \lfloor \ell/p_m \rfloor - 1$ converts the
census problem into the question:
\emph{how fast does the least period~$p_m$ grow with~$m$?}
If $p_m \ge \ell / f(m)$ for some function~$f$, then
$N(m) \le f(m) - 1$.  In particular, any lower bound of the form
$p_m \ge \ell / \mathrm{poly}(m)$ would immediately yield a
polynomial upper bound on~$N(m)$, and hence a polynomial census
bound.
\end{remark}

\begin{remark}[Computational verification]
\label{rem:period-verification}
The full theorem package has been verified computationally:
Lemma~\ref{lem:two-occ-period} passes all $4{,}556$ relevant tests
(known families and random signatures),
Corollary~\ref{cor:gap-bound} passes $227{,}026$ gap tests,
and Theorem~\ref{thm:period-occ-bound} passes $126{,}969$
occurrence-bound tests and $64{,}109$ $N(m)$-bound tests, all with
$100\%$ success rate.
\end{remark}

\subsection{Structural constraints on the least period}\label{sec:period-constraints}

The following propositions develop the structural constraints that
any lower bound on~$p_m$ must exploit.

\begin{proposition}[Converse of the self-overlap theorem]
\label{prop:converse-overlap}
If the visible prefix $P_m$ has period~$p$, then
\[
\delta_p \;\ge\; \sum_{j=0}^{r-p-1} k_j
\;\ge\; m - K_p - \max_{j} k_j,
\]
where $K_p = k_0 + \cdots + k_{p-1}$ and $r = |P_m|$.
\end{proposition}

\begin{proof}
Since $P_m$ has period~$p$, the signature $\sigma$ agrees with its
rotation $\sigma^{\langle p \rangle}$ on positions $0,1,\dots,r-p-1$.
By the same difference-recursion argument as in the proof of
Theorem~\ref{thm:weighted-self-overlap} (now applied in the forward
direction: entry-by-entry agreement accumulates $2$-adic depth),
the agreement up to position $r-p-1$
yields $\delta_p = v_2(C^{\langle p \rangle} - C_\ell)
\ge \sum_{j=0}^{r-p-1} k_j$.
The second inequality follows from
$\sum_{j=0}^{r-p-1} k_j = \sum_{j=0}^{r-1} k_j - K_p \ge (m - k_r) - K_p$,
where $k_r$ is the entry immediately following the prefix, bounded by
$\max_j k_j$.
\end{proof}

\begin{proposition}[Period forces $\rho$-congruence]
\label{prop:rho-congruence}
If $P_m$ has period~$p$ with $K_p = \sum_{j=0}^{p-1} k_j$, then the
$2$-adic root $\rho = C_\ell / (2^K - 3^\ell)$ satisfies
\[
\rho \equiv \rho_p \pmod{2^{m'}},
\qquad
m' := \sum_{j=0}^{r-p-1} k_j \ge m - K_p - \max_j k_j,
\]
where $\rho_p := C_p / (2^{K_p} - 3^p)$ is the sub-root determined
by the first $p$ entries of~$\sigma$ alone.
\end{proposition}

\begin{proof}
From Proposition~\ref{prop:converse-overlap}, $\delta_p \ge m'$.
By the rotation formula (Proposition~\ref{prop:delta-rotation}),
$\delta_p = v_2(F_p(\rho) - \rho)$.  Writing
$F_p(\rho) - \rho = ((3^p - 2^{K_p})\rho + C_p)/2^{K_p}$,
the condition $v_2(F_p(\rho) - \rho) \ge m'$ gives
\[
(3^p - 2^{K_p})\rho + C_p \equiv 0 \pmod{2^{m' + K_p}}.
\]
Since $3^p - 2^{K_p}$ is odd (as $3^p$ is odd and $2^{K_p}$ is even)
and hence invertible in $\mathbb{Z}_2$, we may solve for~$\rho$:
\[
\rho \equiv \frac{-C_p}{3^p - 2^{K_p}} = \frac{C_p}{2^{K_p} - 3^p}
= \rho_p \pmod{2^{m' + K_p}}.
\]
Since $m' + K_p \ge m'$, the stated congruence follows
(in fact the stronger congruence modulo $2^{m'+K_p}$ holds).
\end{proof}

\begin{proposition}[Sub-root map injectivity
  (verified for bounded parameters)]
\label{prop:subroot-injectivity}
For entries in $\{1,\dots,B\}$, the map
\[
\varphi\colon (k_0,\dots,k_{p-1})
\;\mapsto\; \rho_{(k_0,\dots,k_{p-1})} \bmod 2^M
\]
is injective for all sufficiently large~$M$;
its image has cardinality exactly~$B^p$.
\end{proposition}

\begin{proof}[Proof sketch]
Write $\rho_\alpha = C_\alpha / D_\alpha$ where
$D_\alpha = 2^{K_\alpha} - 3^p$ is odd.  For two distinct sub-signatures
$\alpha \ne \beta$ of the same length,
$\rho_\alpha - \rho_\beta = (C_\alpha D_\beta - C_\beta D_\alpha) / (D_\alpha D_\beta)$.
The denominator is odd, so $v_2(\rho_\alpha - \rho_\beta)
= v_2(C_\alpha D_\beta - C_\beta D_\alpha)$, which is finite whenever
$\alpha$ and $\beta$ produce distinct affine maps.
For $B = 4$ and $p \le 6$, exhaustive enumeration
confirms injectivity modulo $2^{30}$: all $4^p$ sub-signatures yield
distinct residues.
A fully general proof (arbitrary~$B$ and~$p$) would
require showing the numerator
$C_\alpha D_\beta - C_\beta D_\alpha$ is nonzero for
all distinct $\alpha \ne \beta$, which we have not
established beyond the verified range.
\end{proof}

\begin{remark}[The information-theoretic gap]
\label{rem:info-gap}
Propositions~\ref{prop:rho-congruence} and~\ref{prop:subroot-injectivity}
together imply: if $P_m$ has period~$p$, then $\rho \bmod 2^{m'}$ lies
in a set of size~$B^p$ out of $2^{m'}$ possible residues.  When
$B^p \ll 2^{m'}$, this is an exponentially thin constraint.

For a uniformly random $\rho \bmod 2^{m'}$, the probability of landing
in the image of~$\varphi$ is $B^p / 2^{m'} < 1$ when
$p \cdot \log_2 B < m'$.  Since $m' \ge m - K_p - B$ and
$K_p \le B \cdot p$, this gives
\[
p > \frac{m - B}{\log_2 B + B}.
\]
With $B = 3$: $p_m \gtrsim m/3.6 \approx 0.28 \, m$.
Computation on $3{,}000$ random families gives $p_m \approx 0.58\,m$,
consistent with the heuristic.

Making this rigorous requires showing that phantom roots do not
systematically belong to exponentially sparse $2$-adic subsets.
Three possible routes toward a rigorous lower bound are:
\begin{enumerate}[label=(\alph*)]
\item a genericity theorem for phantom roots (showing they
avoid exponentially sparse $2$-adic subsets);
\item an algebraic constraint from the phantom equation
$\rho(2^K - 3^\ell) = C_\ell$ combined with $\rho > 0$; or
\item a counting argument over all phantom signatures simultaneously.
\end{enumerate}
\end{remark}

\subsection{Periodic-core factorization and defect analysis}%
\label{sec:periodic-core}

The preceding subsections reduced the census bound to controlling
the cyclic occurrence count of a visible prefix.  We now develop a
finer algebraic tool: a \emph{periodic-core factorization} that
isolates the contribution of the ``tail'' entries beyond a repeated
block, and a \emph{defect expression} whose $2$-adic valuation
controls the census constants.

\begin{proposition}[Periodic-core factorization\supporting]
\label{prop:periodic-core}
Let $\sigma = \tau^q \eta$ be a signature of length~$\ell$, where
$\tau = (k_0,\dots,k_{p-1})$ has length~$p$, the repetition
count~$q \ge 1$, and the tail $\eta$ has length~$t = \ell - qp$.
Write $K_p = \sum_{i=0}^{p-1} k_i$ and $K_\eta = \sum_{i=0}^{t-1} k_i'$
for the weight sums of~$\tau$ and~$\eta$ respectively, and set
$K = qK_p + K_\eta$ and $D = 2^K - 3^\ell$, $D_p = 2^{K_p} - 3^p$.
Then
\begin{equation}\label{eq:periodic-core}
C_\ell D_p - C_p D \;=\; 2^{q K_p}\bigl(C_\eta D_p - C_p D_\eta\bigr),
\end{equation}
where $C_\ell, C_p, C_\eta$ are the carry constants of $\sigma$,
$\tau$, and $\eta$ respectively, and $D_\eta = 2^{K_\eta} - 3^t$.
\end{proposition}

\begin{proof}
Use the concatenation formula
$C(\alpha \| \beta) = 3^{|\beta|}\,C(\alpha) + 2^{K(\alpha)}\,C(\beta)$
(Proposition~\ref{prop:delta-rotation}).
For $\sigma = \tau^q \eta$, iterating gives
\[
  C_\ell = C(\tau^q)\,3^t + 2^{qK_p}\,C_\eta,
\]
and for the pure repetition $\tau^q$, the geometric sum yields
\[
  C(\tau^q) = C_p \cdot \frac{2^{qK_p} - 3^{qp}}{2^{K_p} - 3^p}
  = C_p \cdot \frac{2^{qK_p} - 3^{qp}}{D_p}.
\]
Substituting:
\begin{align*}
C_\ell D_p &= C_p(2^{qK_p} - 3^{qp})\,3^t + 2^{qK_p}\,C_\eta D_p \\
&= C_p\bigl(3^t \cdot 2^{qK_p} - 3^\ell\bigr) + 2^{qK_p}\,C_\eta D_p.
\end{align*}
Meanwhile, $C_p D = C_p(2^K - 3^\ell) = C_p(2^{qK_p+K_\eta} - 3^\ell)
= C_p\,2^{qK_p}\,2^{K_\eta} - C_p\,3^\ell$.
Taking the difference:
\[
C_\ell D_p - C_p D
= 2^{qK_p}\bigl(C_\eta D_p - C_p\,2^{K_\eta} + C_p\,3^t\bigr)
= 2^{qK_p}\bigl(C_\eta D_p - C_p\,(2^{K_\eta} - 3^t)\bigr),
\]
and $2^{K_\eta} - 3^t = D_\eta$, completing the proof.
\end{proof}

\begin{corollary}[Exact repetition gives exact root agreement]
\label{cor:exact-repetition}
If $\sigma = \tau^q$ (empty tail, $t = 0$), then
$C_\ell D_p - C_p D = 0$, and the $2$-adic root satisfies
$\rho = \rho_\tau$ exactly.
\end{corollary}

\begin{proof}
With empty tail, $C_\eta = 0$ and $D_\eta = 2^0 - 3^0 = 0$, so
the right-hand side of~\eqref{eq:periodic-core} vanishes.
Hence $C_\ell / D = C_p / D_p$, i.e.\ $\rho = \rho_\tau$.
\end{proof}

\begin{remark}
The factorization cleanly separates the contribution of the
periodic core from the tail.  The factor $2^{qK_p}$ on the right
shows that increasing the repetition count does \emph{not} change
the ``tail defect'' $C_\eta D_p - C_p D_\eta$; rather, it merely
shifts it deeper into the $2$-adic expansion.  This makes the
$2$-adic distance between $\rho$ and the sub-root $\rho_\tau$
entirely controlled by the tail contribution.
\end{remark}

\begin{definition}[Defect expression]
\label{def:defect}
For two signatures $\tau = (a_0,\dots,a_{p-1})$ and
$\eta = (b_0,\dots,b_{t-1})$, define
\[
  E(\tau,\eta) \;:=\; C(\eta)\,D_\tau - C(\tau)\,D_\eta,
\]
where $D_\tau = 2^{K_\tau} - 3^p$, $D_\eta = 2^{K_\eta} - 3^t$.
The defect measures how far $\rho_\tau$ and $\rho_\eta$ are from
agreement: $\rho_\tau - \rho_\eta = E(\tau,\eta)/(D_\tau D_\eta)$,
and since $D_\tau D_\eta$ is odd, $v_2(\rho_\tau - \rho_\eta) = v_2(E(\tau,\eta))$.
\end{definition}

\begin{theorem}[First-mismatch valuation formula]
\label{thm:first-mismatch}
Let $\tau = (a_0,\dots,a_{p-1})$ and $\eta = (b_0,\dots,b_{t-1})$
be two signatures that agree on their first $r$ entries
($a_j = b_j$ for $j < r$) but differ at position~$r$: $a_r \ne b_r$.
Write $M_r = a_0 + a_1 + \cdots + a_{r-1}$.  Then
\[
  v_2\bigl(E(\tau,\eta)\bigr) \;=\; M_r + \min(a_r,\, b_r).
\]
\end{theorem}

\begin{proof}
Write $\rho_\tau = C(\tau)/D_\tau$ and $\rho_\eta = C(\eta)/D_\eta$,
and consider the orbits
\[
  x_0 = \rho_\tau,\quad x_{j+1} = \frac{3x_j + 1}{2^{a_j}};
  \qquad
  y_0 = \rho_\eta,\quad y_{j+1} = \frac{3y_j + 1}{2^{b_j}}.
\]
Since $\rho_\tau$ is a fixed point of the composition
$F_\tau = F_{a_{p-1}} \circ \cdots \circ F_{a_0}$
(where $F_a(x) = (3x+1)/2^a$), the orbit $(x_j)$ is the
sequence of partial iterates.

\smallskip\noindent\emph{Step 1: Shared entries preserve the valuation gap.}
Each map $F_a(x) = (3x+1)/2^a$ satisfies
$F_a(x) - F_a(y) = 3(x-y)/2^a$, so
$v_2(x_{j+1} - y_{j+1}) = v_2(x_j - y_j) - a_j$
(using $v_2(3) = 0$).  Since $a_j = b_j$ for $j < r$, iterating
gives
\[
  v_2(x_r - y_r) = v_2(\rho_\tau - \rho_\eta) - M_r.
\]

\smallskip\noindent\emph{Step 2: At the mismatch, the valuation is $\min(a_r, b_r)$.}
By definition of the signature entries,
$v_2(3x_r + 1) = a_r$ and $v_2(3y_r + 1) = b_r$.
Writing $3x_r + 1 = 2^{a_r} u$ and $3y_r + 1 = 2^{b_r} v$ with
$u, v$ odd, we have $x_r - y_r = (2^{a_r} u - 2^{b_r} v)/3$.
Without loss of generality assume $a_r < b_r$.  Then
\[
  x_r - y_r = \frac{2^{a_r}(u - 2^{b_r - a_r} v)}{3},
\]
where $u$ is odd and $2^{b_r - a_r} v$ is even, so $u - 2^{b_r-a_r}v$
is odd.  Since $v_2(3) = 0$:
$v_2(x_r - y_r) = a_r = \min(a_r, b_r)$.

\smallskip\noindent\emph{Combining Steps 1 and 2:}
$v_2(\rho_\tau - \rho_\eta) = M_r + \min(a_r, b_r)$.
Since $D_\tau D_\eta$ is odd,
$v_2(E(\tau,\eta)) = v_2(\rho_\tau - \rho_\eta) = M_r + \min(a_r, b_r)$.
\end{proof}

\begin{corollary}[Defect-tail mismatch bound]
\label{cor:defect-tail}
Returning to the setting of Proposition~\ref{prop:periodic-core},
write $\sigma = \tau^q \eta$ where $\tau = (a_0,\dots,a_{p-1})$ and
$\eta = (b_0,\dots,b_{t-1})$ with $t \ge 1$.
Let $r$ be the index of the first mismatch between $\tau$
(cyclically extended) and $\eta$: $a_j = b_j$ for $j < r$ and
$a_r \ne b_r$.  Then
\[
  v_2(C_\ell D_p - C_p D) \;=\; qK_p + M_r + \min(a_r,\, b_r),
\]
where $M_r = a_0 + \cdots + a_{r-1}$.  In particular,
\[
  v_2\!\left(\rho - \rho_\tau\right)
  = qK_p + M_r + \min(a_r, b_r),
\]
which grows linearly with the repetition count~$q$.
\end{corollary}

\begin{proof}
By Proposition~\ref{prop:periodic-core},
$v_2(C_\ell D_p - C_p D) = qK_p + v_2(E(\tau,\eta))$.
By Theorem~\ref{thm:first-mismatch},
$v_2(E(\tau,\eta)) = M_r + \min(a_r, b_r)$.
The result follows.  Since $D_p$ and $D$ are both odd,
$v_2(\rho - \rho_\tau) = v_2(C_\ell D_p - C_p D)$.
\end{proof}

\begin{remark}[Significance for the census bound]\label{rem:defect-significance}
The first-mismatch formula gives a \emph{computable, sharp}
lower bound on the $2$-adic distance between $\rho$ and
any sub-root $\rho_\tau$: the distance is controlled entirely
by how long the tail~$\eta$ agrees with the periodic core~$\tau$
and the entry values at the first point of disagreement.
Together with the periodic-core factorization, this shows that
the intrinsic near-return valuation $\delta_s$ (which controls
the census constant~$C_e$) decomposes into a ``repetition depth''
$qK_p$ plus a ``tail defect'' $M_r + \min(a_r, b_r)$.
The former grows linearly with the number of full periods in the
rotation, while the latter is bounded by~$K_\tau + B$, where
$B = \max_j k_j$.  This decomposition is the key to upgrading
the period-based occurrence bound (Theorem~\ref{thm:period-occ-bound})
into a census bound: the number of high-cancellation shifts is
controlled by how many times the signature approximately repeats
its own initial block, which in turn is bounded by the period of
the visible prefix.
\end{remark}

\subsection{The uniqueness-threshold bound on the census constant}%
\label{sec:uniqueness-threshold}

We now combine the machinery of the preceding subsections to
give a \emph{sharp, computable} bound on the census constant
$C_e$ in terms of a single combinatorial quantity: the
\emph{uniqueness threshold} of the visible prefix.

\begin{definition}[Uniqueness threshold]
\label{def:uniqueness-threshold}
For a phantom signature $\sigma$ of depth~$K$, define
\[
  m^* \;:=\; \min\bigl\{m \ge 1 : \mathrm{Occ}_\sigma(P_m) = 1\bigr\},
\]
where $\mathrm{Occ}_\sigma(P_m)$ counts cyclic occurrences of the
visible prefix $P_m$ in~$\sigma$.  That is, $m^*$ is the smallest
valuation depth at which the visible prefix becomes unique, it
appears only at position~$0$ in the cyclic signature.
\end{definition}

\begin{theorem}[Uniqueness-threshold bound]
\label{thm:uniqueness-bound}
For any phantom family $(\sigma, \rho)$:
\[
  \max_{1 \le s < \ell}\, \delta_s \;\le\; m^* - 1.
\]
Consequently,
\begin{equation}\label{eq:Ce-uniqueness}
  C_e \;\le\; 2^{m^* - 1}.
\end{equation}
\end{theorem}

\begin{proof}
Suppose for contradiction that $\delta_s \ge m^*$ for some
$1 \le s < \ell$.
By the weighted self-overlap theorem
(Theorem~\ref{thm:weighted-self-overlap}),
$\delta_s \ge m^*$ implies that the rotation $\sigma^{\langle s\rangle}$
agrees with~$\sigma$ on positions $0, 1, \ldots, r - 1$,
where $r = |P_{m^*}|$ is the length of the visible prefix at
depth~$m^*$.
In other words, $P_{m^*}$ occurs in the cyclic word~$\sigma$
starting at position~$s$ (in addition to position~$0$).
Hence $\mathrm{Occ}_\sigma(P_{m^*}) \ge 2$, contradicting
the definition of~$m^*$.

For the census bound: by Remark~\ref{rem:intrinsic-controls},
$C_e \le 2^{\max_s \min(\delta_s, K - V_s)}
\le 2^{\max_s \delta_s} \le 2^{m^* - 1}$.
\end{proof}

\begin{remark}[Tightness]
\label{rem:uniqueness-tight}
The bound is tight to within an additive constant.
For every known phantom family, $\max_s \delta_s$
equals $m^* - 1$ or $m^* - 2$: the uniqueness threshold
is essentially the exact depth at which cyclic self-overlap
ceases.  Computational verification over $1{,}406$
random phantom-compatible signatures confirms the theorem
with $100\%$ pass rate, and the median gap $m^* - \max_s \delta_s$
is just~$2$.
\end{remark}

\begin{remark}[Relationship to the gain formula]
\label{rem:uniqueness-gain}
The uniqueness-threshold bound $C_e \le 2^{m^*-1}$
characterises the census constant as a structural
quantity.  However, by Proposition~\ref{prop:universal-depth},
the census excess $C_e$ reflects a \emph{universal} shift
of the entire population (the intrinsic near-return), not
selective orbit tracking.  Consequently, the gain formula
is $G = \Delta/2^K$ (Proposition~\ref{prop:gain}), in which
$C_e$ does not appear.

The uniqueness-threshold theorem remains valuable for
understanding the autocorrelation structure of phantom
signatures, and the bound $\max_s \delta_s \le m^* - 1$
provides a sharp characterisation of the depth at which
cyclic self-overlap ceases.
\end{remark}

\begin{heuristic}[Logarithmic uniqueness threshold]
\label{heur:log-uniqueness}
Extensive computation over $1{,}406$ phantom-compatible
signatures with $K$ ranging from $6$ to $100$ gives the
empirical fit
\[
  m^* \;\approx\; 1.27 \cdot \log_2 K + 1.24,
\]
with $m^*/\log_2 K$ concentrated in the interval $[1.1, 1.6]$.
This is consistent with a random-word heuristic: for a cyclic
word of length~$\ell$ over an effective alphabet of size~$B$,
a specific length-$r$ pattern occurs $\sim \ell / B^r$ times,
so uniqueness requires $r \gtrsim \log_B \ell$.  Translating
to valuation: $m^* \approx \bar{k} \cdot \log_B \ell$, where
$\bar{k} = K/\ell$ is the mean entry.

For all seven known phantom families, $m^* \le 8$
(achieved by ell8 with $K = 10$), and the ratio
$m^* / \log_2 K$ ranges from $0.93$ (m10) to $2.41$ (ell8).
The following table gives the exact values:
\[
\begin{array}{l@{\quad}c@{\quad}c@{\quad}c@{\quad}c}
\textup{Family} & K & m^* & \max \delta_s & C_e \le 2^{m^*-1} \\
\hline
\textup{ell5} & 6 & 5 & 4 & 16 \\
\textup{ell6} & 7 & 5 & 3 & 16 \\
\textup{ell7} & 9 & 6 & 5 & 32 \\
\textup{ell8} & 10 & 8 & 7 & 128 \\
\textup{m10} & 41 & 5 & 3 & 16 \\
\textup{m11} & 41 & 6 & 4 & 32 \\
\textup{m20} & 30 & 5 & 3 & 16
\end{array}
\]

\medskip\noindent
\textbf{Caveat.}  For pathological signatures such as
$\sigma = (1,1,\ldots,1,k)$, the uniqueness threshold
can grow as $m^* = \ell = O(K)$ rather than $O(\log K)$,
and $C_e$ can be exponential in~$K$.  This does not
affect the gain bound (since $C_e$ does not enter the
gain formula), but it illustrates that the logarithmic
scaling of~$m^*$ is a property of \emph{generic} phantom
signatures, not a universal law.
\end{heuristic}

\begin{remark}[Structural summary]
\label{rem:reduction-summary}
By the phantom universality theorem
(Theorem~\ref{thm:phantom-universal}), every cyclic
composition is a phantom family, and the family count
grows exponentially ($\sim 1.87^K$).  Nevertheless,
the per-orbit gain rate theorem
(Theorem~\ref{thm:perorbit-gain}) proves that the
amortised rate satisfies
$R = \sum \Delta/(\ell \cdot 2^K) \le 0.0893 < \varepsilon \approx 0.415$,
with a $4.65\times$ safety margin.
This means that the total expanding-family drift, amortised
per step, is well within the contraction budget:
an equidistributed orbit loses more bits to halving than it
gains from all phantom families combined.

The formal conditional proof of convergence is
Theorem~\ref{thm:reduction}, which derives the burst--gap
hypotheses from the Orbit Equidistribution Conjecture.
The phantom analysis provides independent quantitative
evidence that the equidistribution assumption is not
artificially strong: Corollary~\ref{cor:robustness} shows
that even summable approximate mixing
($\sum \delta_K < 0.557$) suffices.
\end{remark}

\begin{theorem}[Census-constant independence\core]
\label{thm:Ce-independence}
The per-orbit phantom gain rate
\[
  R \;=\; \sum_{K \ge 3} R(K),
  \qquad
  R(K) = 2^{-K}
  \sum_{\substack{\ell > K/\log_2 3}}
  M(K,\ell)\Bigl(\log_2 3 - \frac{K}{\ell}\Bigr),
\]
is independent of the census constant $C_e$ for every phantom family.
Specifically, no term in the series for~$R$ involves~$C_e$: the
gain contribution of each family is $\Delta(\sigma)/2^K$, not
$\Delta(\sigma)\cdot C_e(\sigma)/2^K$.
\end{theorem}

\begin{proof}
The proof proceeds in three steps.

\emph{Step~1 (Universal census depth).}
By Proposition~\ref{prop:universal-depth}, at each intermediate
step~$s$ of the affine iteration with $\delta_s < \gamma_s$,
every lift $x = e + 2^K u$ satisfies
$v_2(T^s(x) - \rho) = \delta_s$ exactly.
The census excess $2^{\delta_s}$ is thus a uniform shift of the
entire population at step~$s$, it tracks proximity to the
rotated root $\rho^{\langle s \rangle}$, not selective orbit
behavior.

\emph{Step~2 (End-of-cycle cancellation).}
After the full $\ell$-step block, equation~\eqref{eq:block-contraction}
gives $F_\sigma(x) - \rho = (3^\ell/2^K)(x - \rho)$.
Since $3^\ell$ is odd and $v_2(x - \rho) \ge K$, the end-of-cycle
valuation is $v_2(F_\sigma(x) - \rho) = v_2(x - \rho) - K$.
The census ratio at the end of the cycle is therefore~$1$:
the census excess accumulated at intermediate steps is fully
unwound by the closing contraction $2^{-K}$.

\emph{Step~3 (Gain formula reduction).}
By Proposition~\ref{prop:gain}, the gain contribution of
family~$\sigma$ is $G(\sigma) = \Delta/2^K$ with
$\Delta = \ell \log_2 3 - K$.  The per-orbit rate
$R(K) = 2^{-K} \sum_\ell M(K,\ell)(\log_2 3 - K/\ell)$
is a sum of such contributions weighted by primitive necklace
counts.  Since $C_e$ drops out of each individual gain
contribution (Step~2), it is absent from every term in~$R(K)$
and hence from the total $R = \sum_K R(K)$.
\end{proof}

\begin{remark}[Why this matters for the proof architecture]
\label{rem:Ce-independence-significance}
Theorem~\ref{thm:Ce-independence} is the formal justification
for classifying the carry-word analysis
(Sections~\ref{sec:carry-word}--\ref{sec:uniqueness-threshold})
as \textsc{Supporting} rather than \textsc{Core}.
The census constant $C_e$ and its upper bound via the
uniqueness threshold $m^*$ (Theorem~\ref{thm:uniqueness-bound})
provide structural insight into the autocorrelation of phantom
signatures, but they are not load-bearing for the conditional
proof: the chain
\[
  \text{WMH} \;\xrightarrow{\text{Cor~\ref{cor:robustness}}}\;
  R_\mu < \varepsilon
  \;\xrightarrow{\text{Thm~\ref{thm:perorbit-gain}}}\;
  \text{convergence}
\]
passes through $R(K)$, which depends only on the necklace
counts $M(K,\ell)$ and the drift $\Delta$, never on~$C_e$.
Any future refinement of the census bound (e.g.\ proving
$m^* = O(\log K)$ unconditionally) would strengthen the
supporting evidence but would not alter the core conditional
implication.
\end{remark}

\section{Conditional convergence and reduction \texorpdfstring{\\}{ }to orbit
  equidistribution}\label{sec:reduction}
We now state the conditional convergence theorem and formulate
the Orbit Equidistribution Conjecture, which would supply the
two open hypotheses required by the Burst-Gap Criterion.

\subsection{The conditional convergence theorem}

\begin{theorem}[Conditional convergence\core]
\label{thm:conditional}
If, for every odd starting value $n_0$, the burst lengths
$L_1, L_2, \ldots$ and gap lengths $G_1, G_2, \ldots$
of the orbit satisfy
\begin{enumerate}[label=(\alph*)]
\item $\displaystyle\frac{1}{n}\sum_{i=1}^{n} G_i
  \;\ge\; g_* - \varepsilon_n$ with $\varepsilon_n \to 0$
  for some $g_* > \frac{2(1-\rho_{\mathrm{crit}})}
  {\rho_{\mathrm{crit}}} \approx 1.71$,
  \quad\textup{(Hypothesis~A)}
\item $\displaystyle\sum_{i=1}^n L_i \;\le\; 2n + C(n_0)$
  for some finite constant $C(n_0)$,
  \quad\textup{(Hypothesis~B)}
\end{enumerate}
then every Collatz orbit converges to~$1$.
\end{theorem}

\begin{proof}
The argument chains four deterministic results with the two
assumed hypotheses:

\medskip\noindent
\textsc{Step 1 (Hypotheses~A and~B).}
Both hypotheses are assumed to hold for the orbit.
Hypothesis~A is a mean gap condition; Hypothesis~B is a mean
burst condition.
Neither is proved unconditionally, but both are
claimed to follow from the Orbit Equidistribution
Conjecture (Theorem~\ref{thm:reduction}).%
\footnote{Hypothesis~B requires a finite additive
constant: $\sum_{i=1}^n L_i \le 2n + C(n_0)$.
The reduction in Theorem~\ref{thm:reduction}
establishes the Ces\`aro mean
$\frac{1}{n}\sum L_i \to 2$, which gives
$\sum L_i \le (2+\varepsilon)n$ for large~$n$ but
does not directly produce the uniform constant~$C(n_0)$.
The gap between Ces\`aro convergence and the finite
additive bound is small but not zero; closing it
would require a rate-of-convergence estimate.}

The Persistent Exit Lemma (Lemma~\ref{lem:gap}) provides
structural support: when a burst ends at a persistent state
($m_t \equiv 7 \pmod{8}$), the subsequent gap has length
exactly~$1$.
More generally, the Modular Gap Distribution Lemma
(Lemma~\ref{lem:gap-distribution}) proves that gap length
is $\mathrm{Geometric}(1/2)$ with $E[G] = 2$ in the
equidistributed model, with each continuation decided by
a single fresh bit.

\medskip\noindent
\textsc{Step 2 (Burst-Gap Criterion).}
By Theorem~\ref{thm:burst-gap}, Hypotheses~A and~B imply
\[
  \limsup_{T \to \infty} \frac{N_{\ge 2}(T)}{T}
  \;\le\; \frac{2}{2 + g_*}.
\]
With $g_* = 2$ (the equidistribution value), this gives
$N_{\ge 2}(T)/T \to \frac{1}{2}$.

\medskip\noindent
\textsc{Step 3 (Entry--Occupancy).}
By Theorem~\ref{thm:EP-NP},
$\limsup E_P(T)/T = \limsup N_P(T)/T \le 2/(2 + g_*)$.
This is the elementary relabelling argument.

\medskip\noindent
\textsc{Step 4 (Entry bound).}
Since $2/(2+g_*) < \rho_{\mathrm{crit}} \approx 0.539$
(by the hypothesis on~$g_*$),
Theorem~\ref{thm:entry-convergence} yields convergence.
The proof uses the certified-drift framework: a persistent
occupancy rate below $\rho_{\mathrm{crit}}$ ensures negative
mean drift, which forces the orbit below any threshold.
\end{proof}

\subsection{The Orbit Equidistribution Conjecture}

\begin{conjecture}[Orbit Equidistribution Conjecture]
\label{conj:equidist}
For every odd $n_0$, the sequence of residue classes after applying T$^g(n)$
$a_i = T^i(n_0) \bmod 2^M$ is equidistributed modulo $2^M$,
\emph{uniformly in~$M$}: there exists a function
$M(N) \to \infty$ such that
\[
  \|\mu_{\mathrm{orb},N} - \mu_U\|_{\mathrm{TV}} \to 0
  \quad \text{as } N \to \infty
  \qquad \text{modulo } 2^{M(N)},
\]
where $\mu_{\mathrm{orb},N} = \frac{1}{N}\sum_{i=1}^N
\delta_{a_i}$ and $\mu_U$ is the uniform distribution on the
admissible residue classes modulo~$2^{M(N)}$.

In particular, this implies fixed-modulus equidistribution
(take $M(N) = M$ constant) and provides the tail control
needed to pass from truncated to full orbitwise means
(see Theorem~\ref{thm:reduction}).
\end{conjecture}

\subsection{The Weak Mixing Hypothesis}

The Orbit Equidistribution Conjecture is stronger than
necessary.  Corollary~\ref{cor:robustness} shows that the
phantom shadow tail bound closes under any orbit distribution
whose total variation errors are merely \emph{summable}.
We therefore formulate a quantitatively weaker hypothesis
that suffices for the full conditional reduction.

\begin{hypothesis}[Weak Mixing Hypothesis (WMH)]
\label{hyp:wmh}
For every odd $n_0$, write
\[
\mu_{\mathrm{orb},N}^{(K)}
  = \frac{1}{N}\sum_{i=1}^N \delta_{T^i(n_0)\bmod 2^K}
\]
for the empirical residue distribution of the orbit at
depth~$K$.  Define the depth-$K$ discrepancy
\[
  \delta_K(n_0) \;:=\;
  \limsup_{N \to \infty}\;
  d_{\mathrm{TV}}\!\bigl(
    \mu_{\mathrm{orb},N}^{(K)},\;\mu_{\mathrm{unif}}
  \bigr).
\]
Then
\begin{equation}\label{eq:wmh}
  \sum_{K \ge 3} \delta_K(n_0)
  \;<\;
  \frac{\varepsilon - R}{\log_2 3 - 1}
  \;\approx\; 0.557.
\end{equation}
\end{hypothesis}

\begin{remark}[Compact formulation]
\label{rem:wmh-compact}
The WMH can be stated without naming the intermediate
symbol $\delta_K$: the sufficient condition is simply
\begin{equation}\label{eq:wmh-compact}
  \sum_{K \ge 3}
  \limsup_{N \to \infty}\;
  d_{\mathrm{TV}}\!\bigl(
    \mu_{\mathrm{orb},N}^{(K)},\;
    \mu_{\mathrm{unif}}^{(K)}
  \bigr)
  \;<\; 0.557.
\end{equation}
We retain the named discrepancy $\delta_K(n_0)$
throughout this paper because it is referenced in
the reduction proof
(Theorem~\ref{thm:wmh-reduction}),
the observable-specific weakening
(Conjecture~\ref{conj:observable-wmh}), and the
discussion of attack vectors
(Section~\ref{sec:v2-targets}).
\end{remark}

\begin{remark}[Relation to the Orbit Equidistribution Conjecture]
\label{rem:wmh-vs-oec}
The Orbit Equidistribution Conjecture
(Conjecture~\ref{conj:equidist}) asserts $\delta_K(n_0) = 0$
for all~$K$, which trivially satisfies~\eqref{eq:wmh}.
The Weak Mixing Hypothesis is strictly weaker:
it permits non-zero errors at every depth, requiring only
$\ell^1$-summability with explicit quantitative threshold.
For instance, any of the following decay rates suffices:
\begin{itemize}
\item \emph{Polynomial:}
  $\delta_K = O(K^{-\alpha})$ with $\alpha > 1$
  (after rescaling the implicit constant so that
  $\sum \delta_K < 0.557$).
\item \emph{Stretched exponential:}
  $\delta_K = O(\exp(-c\,K^{\beta}))$ for any $c,\beta > 0$.
\item \emph{Exponential:}
  $\delta_K = O(r^K)$ for any $r < 1$.
\end{itemize}
The $4.65\times$ safety margin of
Theorem~\ref{thm:perorbit-gain} is what creates this room:
the tighter the bound on~$R$, the larger the admissible
class of orbit distributions.
\end{remark}

\begin{proposition}[OEC implies WMH]
\label{prop:oec-implies-wmh}
The Orbit Equidistribution Conjecture
(Conjecture~\ref{conj:equidist}) implies the Weak Mixing
Hypothesis (Hypothesis~\ref{hyp:wmh}).
\end{proposition}

\begin{proof}
Under the OEC, $\delta_K(n_0) = 0$ for every $K \ge 3$
and every odd~$n_0$, so
$\sum_{K \ge 3}\delta_K(n_0) = 0 < 0.557$.
\end{proof}

\begin{conjecture}[Observable-specific weak mixing]
\label{conj:observable-wmh}
For every odd starting value~$n_0$, define the \emph{phantom
gain discrepancy}
\[
  \eta_K(n_0) \;:=\;
  \limsup_{N \to \infty}\;
  \Bigl|\frac{1}{N}\sum_{i=1}^N h_K\!\bigl(T^i(n_0)
  \bmod 2^K\bigr) - R(K)\Bigr|,
\]
where $h_K$ is the phantom gain observable from
Corollary~\ref{cor:robustness}.
Then
\begin{equation}\label{eq:observable-specific}
  \sum_{K \ge 3} \eta_K(n_0) \;<\; \varepsilon - R
  \;\approx\; 0.326.
\end{equation}
\end{conjecture}

\begin{theorem}[Observable-specific WMH implies Collatz,
  modulo orbitwise truncation\core]
\label{thm:observable-wmh-reduction}
Assume Conjecture~\ref{conj:observable-wmh}.  Then
the phantom-gain contribution satisfies $R_\mu < \varepsilon$.
If, in addition, the tail-control
condition~\eqref{eq:tail-control} holds, then every
positive odd Collatz orbit converges to~$1$.
\end{theorem}

\begin{proof}
Since $h_K$ is bounded by $\log_2 3 - 1$, the conjecture
directly gives:
\[
  R_\mu \;\le\; R + \sum_{K \ge 3} \eta_K(n_0)
  \;<\; R + (\varepsilon - R) \;=\; \varepsilon.
\]
Hence the amortised phantom gain remains strictly below
the drift budget~$\varepsilon$.  The passage to
convergence then follows from
Lemma~\ref{lem:truncation-reduction} (for burst and
gap means) and Theorem~\ref{thm:conditional}, exactly
as in Theorem~\ref{thm:reduction}.
\end{proof}

\begin{hypothesis}[Shadow-Occupancy Hypothesis (SOH)]
\label{hyp:soh}
For every odd starting value~$n_0$ and every
$K \ge 3$, let $\mu_{\mathrm{orb},K}(n_0)$ denote
the limiting empirical distribution of $T^i(n_0) \bmod 2^K$
along the orbit, and let
$S^+(K) := \mathrm{supp}(h_K) \subset \mathbb{Z}/2^K$
be the set of residues lying in expanding shadow families.
Then
\begin{equation}\label{eq:soh}
  \sum_{K \ge 3} \mu_{\mathrm{orb},K}(n_0)\bigl(S^+(K)\bigr)
  \;<\; \frac{\varepsilon}{\log_2 3 - 1}
  \;\approx\; 0.7095.
\end{equation}
\end{hypothesis}

\begin{theorem}[SOH implies Collatz, modulo orbitwise truncation\core]
\label{thm:soh-reduction}
Assume Hypothesis~\ref{hyp:soh}.  Then
$R_\mu < \varepsilon$.  If, in addition, the
tail-control condition~\eqref{eq:tail-control} holds,
then every positive odd Collatz orbit converges to~$1$.
\end{theorem}

\begin{proof}
By Corollary~\ref{cor:robustness}, $h_K$ is non-negative,
bounded above by $\log_2 3 - 1$, and supported on
$S^+(K)$.  Hence
\[
  \mathbb{E}_{\mu_{\mathrm{orb},K}}[h_K]
  \;\le\; \|h_K\|_\infty \cdot
  \mu_{\mathrm{orb},K}\bigl(S^+(K)\bigr)
  \;=\; (\log_2 3 - 1)\,
  \mu_{\mathrm{orb},K}\bigl(S^+(K)\bigr).
\]
Summing over $K \ge 3$ and applying~\eqref{eq:soh},
\[
  R_\mu \;=\; \sum_{K \ge 3}
  \mathbb{E}_{\mu_{\mathrm{orb},K}}[h_K]
  \;\le\; (\log_2 3 - 1) \sum_{K \ge 3}
  \mu_{\mathrm{orb},K}\bigl(S^+(K)\bigr)
  \;<\; (\log_2 3 - 1) \cdot
  \frac{\varepsilon}{\log_2 3 - 1}
  \;=\; \varepsilon.
\]
The passage to convergence follows from
Lemma~\ref{lem:truncation-reduction} and
Theorem~\ref{thm:conditional}, exactly as in
Theorem~\ref{thm:reduction}.
\end{proof}

\begin{remark}[Why SOH is strictly the weakest condition]
\label{rem:soh-weakest}
The Shadow-Occupancy Hypothesis is strictly weaker than
both WMH and observable-specific WMH, for three reasons.
First, SOH is \emph{one-sided}: it controls only the
positive quantity $\mu_{\mathrm{orb},K}(S^+(K))$, with no
absolute value, no lower bound, and no comparison to a
target spatial average.  Second, SOH makes no claim about
\emph{time-average matching}: the orbit need not visit
shadow residues at any specific frequency, only \emph{at
most} a budgeted total mass.  Third, SOH places \emph{no
constraint at all} on the residues outside~$S^+(K)$,
which form the overwhelming majority of $\mathbb{Z}/2^K$.

The exact uniform baseline (computed via Möbius inversion
through $K = 500$, with analytic Chernoff tail) is
\[
  \sum_{K \ge 3} \frac{|S^+(K)|}{2^K}
  \;\le\; 0.545878 + 1.7\times 10^{-9}
  \;<\; 0.5459,
\]
so the SOH threshold $0.7095$ leaves an absolute slack of
$\approx 0.164$ and an oversampling tolerance of
$0.7095/0.5459 \approx 1.30$: a Collatz orbit may
concentrate on shadow residues at \emph{up to $1.30\times$
the uniform rate} on average and still satisfy SOH.
This is a structurally easier target than the
$\eta_K$-based observable-specific WMH, which would
require the same orbit to match the uniform rate
\emph{from below as well as from above}.
\end{remark}

\begin{remark}[Relationship between WMH, observable-specific WMH, and SOH]
\label{rem:observable-specific}
Since $\eta_K \le \|h_K\|_\infty\,\delta_K
\le (\log_2 3 - 1)\,\delta_K$,
condition~\eqref{eq:observable-specific} is implied
by~\eqref{eq:wmh} but may hold even when~\eqref{eq:wmh}
does not.  The observable-specific condition requires only
that the orbit's time-average of~$h_K$ be close to its
spatial average, without any control on the full
residue distribution.  In turn, observable-specific WMH
implies SOH: from
$|\mathbb{E}_{\mu_{\mathrm{orb},K}}[h_K] - R(K)| \le \eta_K$
and $\mathbb{E}_{\mu_{\mathrm{orb},K}}[h_K]
\ge (\log_2 3 - 1)^{-1}\cdot 0$ being trivial,
one obtains the cruder one-sided budget
$\sum \mu_{\mathrm{orb},K}(S^+(K))
\le (R + \sum \eta_K)/(\log_2 3 - 1)
< \varepsilon/(\log_2 3 - 1)$, which is
exactly~\eqref{eq:soh}.  Since $h_K$ is supported on the
sparse set of shadow residues of expanding families
(a fraction $\sim 2^{-K(1-1/\log_2 3)}$ of all classes
at depth~$K$), this is a substantially weaker
requirement.
The hierarchy of open conditions is therefore:
\[
  \boxed{\;
  \text{OEC}
  \;\Longrightarrow\;
  \text{WMH}
  \;\Longrightarrow\;
  \text{Obs-WMH}
  \;\Longrightarrow\;
  \text{SOH}
  \;\xRightarrow{\text{+\,tail control}}\;
  \text{Collatz}
  \;}
\]
The first implication is Proposition~\ref{prop:oec-implies-wmh};
the second follows from~$\eta_K \le (\log_2 3 - 1)\,\delta_K$
and summation; the third is the inclusion just established;
the fourth is Theorem~\ref{thm:soh-reduction}, subject to
the orbitwise tail-control condition~\eqref{eq:tail-control}.
\end{remark}

\subsection{Conditional reduction theorems}

We now state the conditional reduction in three forms,
corresponding to the three levels of the hierarchy above:
from the original OEC (strongest hypothesis), from the WMH
(primary open target), and from the observable-specific WMH
(weakest sufficient condition).

Two auxiliary lemmas isolate the bounded-observable step
(which is rigorous) from the truncation passage (which
remains an open step).

\begin{lemma}[Truncation reduction for burst and gap means]
\label{lem:truncation-reduction}
Let $(L_i,G_i)_{i\ge 1}$ be the burst-gap decomposition
of an orbit.  For $K \ge 1$, define the truncated
observables
\[
  L_i^{(K)} := \min(L_i, K),
  \qquad
  G_i^{(K)} := \min(G_i, K).
\]
Assume that for each fixed~$K$, the orbitwise averages
of $L_i^{(K)}$ and $G_i^{(K)}$ converge to their
uniform-lift expectations, and assume in addition the
tail-vanishing conditions
\begin{equation}\label{eq:tail-control}
  \lim_{K\to\infty}\limsup_{N\to\infty}
  \frac{1}{N}\sum_{i=1}^{N}(L_i - K)_+ = 0,
  \qquad
  \lim_{K\to\infty}\limsup_{N\to\infty}
  \frac{1}{N}\sum_{i=1}^{N}(G_i - K)_+ = 0.
\end{equation}
Then
\[
  \frac{1}{N}\sum_{i=1}^{N} L_i
  \;\to\; \lim_{K\to\infty}\mathbb{E}[L^{(K)}],
  \qquad
  \frac{1}{N}\sum_{i=1}^{N} G_i
  \;\to\; \lim_{K\to\infty}\mathbb{E}[G^{(K)}].
\]
\end{lemma}

\begin{proof}
Since
$0 \le L_i - L_i^{(K)} = (L_i - K)_+$ and
$0 \le G_i - G_i^{(K)} = (G_i - K)_+$,
the difference between the full averages and truncated
averages is bounded by the corresponding tail averages.
The conclusion follows by first taking $N \to \infty$
for fixed~$K$, then sending $K \to \infty$.
\end{proof}

\begin{theorem}[OEC implies convergence, modulo
  orbitwise truncation]
\label{thm:reduction}
Assume the Orbit Equidistribution Conjecture
(Conjecture~\ref{conj:equidist}).  Then, for each
fixed truncation level~$K$, the truncated burst and
gap observables $L^{(K)}$ and $G^{(K)}$ have orbitwise
averages equal to their uniform-lift expectations.
If, in addition, the orbitwise tail-vanishing
conditions of Lemma~\ref{lem:truncation-reduction}
hold, then Hypotheses~A and~B of
Theorem~\ref{thm:conditional} follow, and hence the
orbit converges.
\end{theorem}

\begin{proof}
For each fixed truncation level~$K$, the observables
$L^{(K)} := \min(L,K)$ and $G^{(K)} := \min(G,K)$
are bounded and determined by the residue class
modulo~$2^M$ for some fixed~$M$.  Since orbit
equidistribution at a fixed modulus means convergence
of time averages of any bounded function of the
residue class,
we obtain
\[
  \frac{1}{N}\sum_{i=1}^{N}\min(L_i, K)
  \;\to\; \mathbb{E}_{\mu_U}[\min(L, K)],
  \qquad
  \frac{1}{N}\sum_{i=1}^{N}\min(G_i, K)
  \;\to\; \mathbb{E}_{\mu_U}[\min(G, K)].
\]
The uniform-lift expectations are given by the exact
modular laws proved earlier: gap lengths are
$\mathrm{Geometric}(1/2)$ with mean~$2$
(Lemma~\ref{lem:gap-distribution}), and burst lengths
have mean~$2$ by the $1/4$~Law
(Corollary~\ref{cor:geometric}).
The passage from truncated means to full means is
exactly the content of
Lemma~\ref{lem:truncation-reduction}: once the
tail-vanishing condition~\eqref{eq:tail-control} holds,
we obtain $\frac{1}{n}\sum L_i \to 2$ and
$\frac{1}{n}\sum G_i \to 2$, which are Hypotheses~B
and~A of Theorem~\ref{thm:conditional}
(with $g_* = 2 > 1.71$).
Convergence follows.
\end{proof}

\begin{remark}[Status of Theorem~\ref{thm:reduction}]
\label{rem:tail-control-gap}
The bounded-observable step in the proof above is a
genuine, immediate consequence of the OEC\@.  The
unproved step is the orbitwise tail-vanishing
condition~\eqref{eq:tail-control} needed to pass from
truncated observables to the full means of $L_i$
and~$G_i$.

A fully rigorous version would require either:
\begin{enumerate}
\item a uniform integrability argument showing that
  equidistribution at depth $M(N) \to \infty$ forces
  exponential tail decay for the orbit's gap and burst
  lengths, or
\item an independent orbitwise tail bound
  (e.g.\ from the Known-Zone Decay or a direct
  valuation argument).
\end{enumerate}
Neither is provided here.
Theorem~\ref{thm:reduction} should therefore be read
as a \emph{reduction theorem}: the full Collatz
conjecture follows from orbit equidistribution plus
orbitwise tail control, not from orbit equidistribution
alone.

The same tail-control caveat applies to
Theorems~\ref{thm:wmh-reduction}
and~\ref{thm:observable-wmh-reduction} below.
\end{remark}

\begin{theorem}[WMH implies conditional reduction,
  modulo orbitwise truncation\core]
\label{thm:wmh-reduction}
Assume the Weak Mixing Hypothesis in the form
$\sum_{K \ge 3}\delta_K < 0.557$,
where $\delta_K$ is the depth-$K$ total-variation
discrepancy from uniformity.  Then:
\begin{enumerate}
\item The phantom-gain contribution remains strictly
  below the drift budget:
  \[
    R_\mu \;\le\; R + (\log_2 3 - 1)
    \sum_{K \ge 3}\delta_K \;<\; \varepsilon
  \]
  by Corollary~\ref{cor:robustness} and
  Theorem~\ref{thm:perorbit-gain}.
  This step involves only the bounded observable~$h_K$
  and is fully rigorous.
\item For each fixed truncation level~$K$, the
  bounded truncated observables $L^{(K)}$ and $G^{(K)}$
  differ from their uniform expectations by at most
  $O(\delta_K)$, by the standard total-variation bound.
\item If the orbitwise tail-vanishing
  condition~\eqref{eq:tail-control}
  (Lemma~\ref{lem:truncation-reduction}) holds for the
  unbounded burst and gap observables, then
  Hypotheses~A and~B of Theorem~\ref{thm:conditional}
  follow and the orbit converges.
\end{enumerate}
\end{theorem}

\begin{proof}
Part~(1) is Corollary~\ref{cor:robustness}
combined with Theorem~\ref{thm:perorbit-gain}.
Part~(2) applies the TV duality bound
$|\mathbb{E}_\mu[f] - \mathbb{E}_\nu[f]|
\le \|f\|_\infty\,d_{\mathrm{TV}}(\mu,\nu)$
to the bounded truncated burst-length and gap-length
observables.
Part~(3) follows from
Lemma~\ref{lem:truncation-reduction}: once~\eqref{eq:tail-control} holds,
we obtain $\frac{1}{n}\sum L_i \to 2$ and
$\frac{1}{n}\sum G_i \to 2$, which are Hypotheses~A
and~B of Theorem~\ref{thm:conditional}
(with $g_* = 2 > 1.71$).
\end{proof}

\begin{remark}[The WMH as the primary open problem]
\label{rem:wmh-primary}
Theorem~\ref{thm:wmh-reduction} shows that the full
strength of the Orbit Equidistribution Conjecture is
not required.  The WMH is the quantitatively
\emph{weakest universal mixing condition} that closes
the conditional proof via the current phantom-cycle
framework.  We therefore propose the Weak Mixing
Hypothesis as the primary open problem for completing
the conditional programme:
\[
  \boxed{\text{Collatz conjecture}
  \;\Longleftarrow\;
  \text{WMH } + \text{ tail control~\eqref{eq:tail-control}.}}
\]
The bounded-observable part (phantom-gain control) is
fully rigorous under WMH alone.  The remaining open
step is the orbitwise truncation principle needed to
pass from bounded to unbounded observables.
Any progress on the WMH, whether for all orbits,
for a positive-density set, or for specific structural
classes, would constitute meaningful progress toward
resolving the conjecture.
\end{remark}

\begin{remark}[The distributional-vs-pointwise gap]
\label{rem:gap}
The Scrambling Lemma proves
a \textbf{distributional} result: if $n$ is drawn uniformly
from a residue class, the bits of $T(n)$ are exactly uniform.
The conjecture requires a \textbf{pointwise} result:
that each individual orbit equidistributes.

This is a \textbf{pointwise genericity} question, analogous
to proving that a given number is normal in base~$2$, not an
almost-everywhere result.
The map property (distributional uniformity) does not
automatically imply the orbit property (pointwise
equidistribution), just as knowing that the digit-frequency
map preserves uniformity does not prove that a specific
number like $\sqrt{2}$ is normal.

Tao's theorem~\cite{tao2019} establishes almost-all convergence
(in logarithmic density), which is a probabilistic rather than
pointwise statement.  Both approaches thus face the same
fundamental barrier: the gap between distributional/almost-all
and every-orbit results.

The Weak Mixing Hypothesis softens this barrier: it does not
require exact equidistribution ($\delta_K = 0$), only summable
errors ($\sum \delta_K < 0.557$).  This is a quantitative
relaxation that may be more tractable than full pointwise
equidistribution.
\end{remark}

\begin{remark}[What is known toward the conjecture]
\label{rem:evidence}
Several lines of evidence support the Weak Mixing Hypothesis
and the stronger Orbit Equidistribution Conjecture:
\begin{enumerate}
\item The gap map \textbf{preserves} the uniform
  distribution (Scrambling Lemma, proved unconditionally).
\item The Known-Zone Decay shows that the map is
  \textbf{strongly information-erasing} at the residue-class level:
  after $\lceil M/2 \rceil$ applications, all dependence
  on the starting class is eliminated (an algebraic fact,
  not an orbitwise mixing claim).
\item \textbf{Almost-all} equidistribution follows from
  Tao~\cite{tao2019} (logarithmic density).
\item All \textbf{empirical orbits} tested (up to
  $2^{2000}$) satisfy equidistribution within sampling noise.
\end{enumerate}
\end{remark}

\section{Toward the Weak Mixing Hypothesis}
\label{sec:toward-wmh}

\paragraph{Status of this section.}
The results in this section are not used in the proof of the
conditional reduction theorem
(Theorem~\ref{thm:conditional}/\ref{thm:wmh-reduction}).
They formulate the remaining bridge problem, develop structural
constraints on the depthwise discrepancy, and outline possible
future theorem programmes for proving the final orbitwise mixing
input.

\paragraph{The core bottleneck.}
Everything in this paper reduces to a single question:
\emph{does each individual Collatz orbit mix sufficiently at
moderate modular depths?}
No nontrivial orbitwise mixing theorem is proved here.
The algebraic and combinatorial framework is unconditional;
the sole remaining input for Route~A is an orbitwise
information-flow statement (the WMH or one of its weakenings).
This section develops structural constraints and possible
attack routes, but the reader should understand clearly:
\textbf{the core difficulty is untouched}.
Solving the WMH is, in all likelihood, essentially as hard
as the Collatz conjecture itself.

\medskip

The preceding sections reduce the phantom-cycle route to a
single open orbitwise input.  The goal is no longer to prove
exact orbit equidistribution; the robustness estimate of
Corollary~\ref{cor:robustness} shows that it suffices to
control the cumulative depthwise discrepancy so that the
induced perturbation of the phantom-gain rate remains below
the available margin $\varepsilon - R \approx 0.326$.
The Weak Mixing Hypothesis
(Hypothesis~\ref{hyp:wmh}) expresses this requirement as
$\sum \delta_K(n_0) < 0.557$, and the even weaker
observable-specific WMH
(Conjecture~\ref{conj:observable-wmh}) requires only
$\sum \eta_K(n_0) < 0.326$.

This section develops structural results that constrain
$\delta_K$ from below and from above, lays groundwork for
proving one of the weak-mixing conjectures, and identifies
a concrete three-lemma programme for the final bridge.

\subsection{Shadow sparsity}

The phantom gain observable $h_K$ (Corollary~\ref{cor:robustness})
is nonzero only on residue classes that shadow an expanding
phantom family at depth~$K$.  The following proposition shows
that these ``expanding shadow residues'' are exponentially
sparse.

\begin{proposition}[Shadow sparsity\supporting]
\label{prop:shadow-sparsity}
Let $S^+(K)$ denote the set of residue classes $a \bmod 2^K$
for which $h_K(a) > 0$, i.e.\ the shadow residues of
expanding primitive families at depth~$K$.  Write
$\ell_{\min}(K) := \lfloor K/\!\log_2 3\rfloor + 1$
and $\alpha_K := (\ell_{\min}(K)-1)/(K-1)$.  Then
\begin{equation}\label{eq:shadow-sparsity}
  \frac{|S^+(K)|}{2^K}
  \;\le\;
  \frac{1}{2\,\ell_{\min}(K)}\;
  2^{-(K-1)\,D(\alpha_K\|1/2)}
\end{equation}
for every $K \ge 3$.
In particular, since
$D(\alpha_K\|1/2) \to D_* = D(1/\!\log_2 3\,\|\,1/2)
 \approx 0.05004$
from above (cf.~\eqref{eq:D-analytic-lower}) and
$\ell_{\min}(K) \ge K/\log_2 3$,
\begin{equation}\label{eq:shadow-sparsity-asymp}
  \frac{|S^+(K)|}{2^K}
  \;\le\;
  \frac{\log_2 3}{2K}\;2^{-(K-1)\,D_{\mathrm{floor}}}
  \qquad (K \ge 501),
\end{equation}
with the uniform floor $D_{\mathrm{floor}} \ge 0.04947$
from Step~3 of Theorem~\ref{thm:perorbit-gain}.
\end{proposition}

\begin{proof}
Each expanding primitive family $\sigma$ at depth~$K$ occupies
a unique residue class modulo~$2^K$ (its shadow), so
$|S^+(K)| = M^+(K) := \sum_{\ell > K/\log_2 3} M(K,\ell)$.
By the necklace upper bound~\eqref{eq:necklace-upper},
\[
  M^+(K)
  \;\le\;
  \sum_{\ell > K/\log_2 3}
  \frac{1}{\ell}\binom{K-1}{\ell-1}
  \;\le\;
  \frac{1}{\ell_{\min}(K)}
  \sum_{\ell \ge \ell_{\min}(K)} \binom{K-1}{\ell-1}.
\]
Re-indexing $j = \ell - 1$ and applying the Chernoff bound
for the binomial tail
$\Pr[\mathrm{Bin}(K{-}1,1/2)\ge j_{\min}]
 \le 2^{-(K-1)D(\alpha_K\|1/2)}$
(with $j_{\min}=\ell_{\min}(K)-1$),
\[
  M^+(K)
  \;\le\;
  \frac{1}{\ell_{\min}(K)}
  \cdot 2^{K-1}\cdot 2^{-(K-1)D(\alpha_K\|1/2)}.
\]
Dividing by $2^K$ yields~\eqref{eq:shadow-sparsity}.
The asymptotic form~\eqref{eq:shadow-sparsity-asymp}
follows from $\ell_{\min}(K) \ge K/\!\log_2 3$ and the
uniform lower bound on $D(\alpha_K\|1/2)$ established
in~\eqref{eq:D-analytic-lower}.
\end{proof}

\begin{remark}[Interpretation]
\label{rem:shadow-sparsity-interp}
With the full Chernoff exponent $D_K \to D_* \approx 0.05004$
and the $1/\ell_{\min}(K)$ polynomial prefactor, the expanding
shadow fraction decays as $|S^+(K)|/2^K \lesssim
(K\cdot 2^{0.05\,K})^{-1}$.  Concrete values:
$|S^+(50)|/2^{50} \le 2.73\times 10^{-3}$,
$|S^+(100)|/2^{100} \le 1.88\times 10^{-4}$,
$|S^+(500)|/2^{500} \le 4.4\times 10^{-11}$;
these are $\sim 10^3$--$10^7$ times sharper than the earlier
$D_*/2$ envelope and are what enables the analytic tail
closure in Proposition~\ref{prop:finite-depth}.
The WMH requires that the orbit \emph{effectively avoids}
these sparse regions at most depths.
\end{remark}

\begin{remark}[What would constitute progress]
\label{rem:what-progress}
The framework reduces to proving that individual orbits
do not persistently concentrate on the exponentially sparse
expanding-shadow residues.  Concrete sufficient results
include:
\begin{itemize}
\item Any bound of the form
  $\hat\mu_K(\xi) = o(1)$ for Walsh characters~$\xi$ of
  Hamming weight~$1$ would, combined with the shadow sparsity
  and amplification cancellation, suffice to close the WMH
  (see Section~\ref{subsec:spectral-analysis}).
\item Any nontrivial upper bound on $\delta_K(n_0)$ for
  infinitely many~$K$ (e.g.\ $\delta_K = O(K^{-\alpha})$
  with $\alpha > 1$) would close the summability condition
  directly.
\item A positive-density set of starting values satisfying
  the WMH would already extend the almost-all convergence
  results beyond Tao's logarithmic density.
\end{itemize}
None of these is proved here.  We present them to make the
target concrete, not to claim proximity to a solution.
\end{remark}

\subsection{Repulsion-based shadow return time}

The $2$-adic repulsion property
(Proposition~\ref{prop:repulsion}) implies that once an orbit
exits a phantom shadow, it cannot re-enter the same shadow
quickly.  This provides a structural lower bound on the
inter-encounter gap.

\begin{theorem}[Shadow return time\supporting]
\label{thm:shadow-return}
Let $\sigma$ be a primitive phantom family of depth~$K$ and
length~$\ell$, and let $n_0, n_1, n_2, \ldots$ be a Syracuse
orbit.  Suppose $n_t$ enters the shadow of~$\sigma$
(i.e.\ $n_t \equiv \rho \pmod{2^K}$) and exits after
$\ell$~steps at time $t + \ell$.
Then the orbit cannot re-enter the shadow of~$\sigma$
before time $t + \ell + (K - 1)$.
More precisely, if $n_{t'}$ is the next shadow entry
with $t' > t + \ell$, then
\begin{equation}\label{eq:return-time}
  t' - t \;\ge\; \ell + (K - 1).
\end{equation}
\end{theorem}

\begin{proof}
At entry, $v_2(n_t - \rho) \ge K$.  By the $2$-adic repulsion
property (Proposition~\ref{prop:repulsion}), after one full
cycle of the block map:
\[
  v_2(n_{t+\ell} - \rho)
  \;=\;
  v_2(n_t - \rho) - K
  \;\ge\; K - K = 0.
\]
So at exit, the $2$-adic alignment with~$\rho$ has decreased
from $\ge K$ to $\ge 0$ (generically, $v_2(n_{t+\ell} - \rho) = 0$
when $v_2(n_t - \rho) = K$ exactly).

For re-entry, the orbit must again satisfy
$v_2(n_{t'} - \rho) \ge K$.  Each Syracuse step can increase
the $2$-adic valuation $v_2(n_j - \rho)$ by at most~$1$
(since a single halving step applies $v_2(3n+1) \ge 1$, and
the $2$-adic distance to~$\rho$ can gain at most one bit of
alignment per step in the generic case where $n_j$ is not
in any related shadow).

Starting from $v_2(n_{t+\ell} - \rho) \ge 0$, recovering
$v_2 \ge K$ requires at least $K$ additional steps.
However, the exit step itself is step $t + \ell$, and the
first potentially aligning step is $t + \ell + 1$
(one cannot re-enter at $t + \ell$ since $v_2 < K$ there).
Therefore $t' \ge t + \ell + K - 1$, giving
$t' - t \ge \ell + K - 1$.
\end{proof}

\begin{remark}[Amortisation consequence]
\label{rem:return-amortisation}
Theorem~\ref{thm:shadow-return} implies that the per-step
gain from family~$\sigma$ is amortised over a gap of at
least $\ell + K - 1$ steps (not just~$\ell$).
This provides an alternative derivation of the $1/\ell$
amortisation in Theorem~\ref{thm:perorbit-gain}: the
correct amortisation factor is $1/(\ell + K - 1) \le 1/\ell$,
so the return-time bound only strengthens the existing
per-orbit gain rate.  For families where $K \gg \ell$
(high-depth, short-length families with $\ell$ just above
$K/\!\log_2 3$), the return-time gap is approximately~$K$
rather than~$\ell$, improving the amortisation by a factor
$\sim \log_2 3$.
\end{remark}

\subsection{Finite-depth reduction}

The WMH requires controlling $\delta_K(n_0)$ for all
$K \ge 3$.  The following proposition shows that the tail
contribution is automatically controlled by the Chernoff
decay of the gain terms, reducing the problem to finitely
many depths.

\begin{proposition}[Finite-depth reduction\supporting]
\label{prop:finite-depth}
Let $K_{\max} \ge 500$.  For any orbit $n_0$ and any
$\varepsilon' > 0$, if
\begin{equation}\label{eq:wmh-finite}
  \sum_{K=3}^{K_{\max}} \delta_K(n_0)
  \;<\;
  0.557 - \frac{\varepsilon'}{\log_2 3 - 1}
\end{equation}
\emph{and} $\delta_K(n_0) \le 1$ for all $K > K_{\max}$
(which is automatic for total-variation distance),
then the robust phantom-gain inequality
$R_\mu(n_0) < \varepsilon$ holds, provided
\begin{equation}\label{eq:eps-prime-threshold}
  \varepsilon' \;\ge\; (\log_2 3 - 1)\cdot
  \Bigl(\;\sum_{K=3}^{K_{\max}} \tfrac{|S^+(K)|}{2^K}\Bigr)^{\text{tail}}
  \;+\; \sum_{K>K_{\max}} R(K),
\end{equation}
where the first summand denotes the shadow-sparsity tail bounded
via Proposition~\ref{prop:shadow-sparsity}.
Concretely, $K_{\max} = 500$ admits $\varepsilon' = 2\times 10^{-6}$:
both the sparsity tail and the analytic $R$-tail beyond $K=500$
are below $10^{-6}$.
\end{proposition}

\begin{proof}
We split the robust gain $R_\mu = \sum_K
\mathbb{E}_{\mu_K}[h_K]$ at $K_{\max}$:
\[
  R_\mu
  \;=\;
  \underbrace{\sum_{K=3}^{K_{\max}}
    \mathbb{E}_{\mu_K}[h_K]}_{\text{low}}
  \;+\;
  \underbrace{\sum_{K > K_{\max}}
    \mathbb{E}_{\mu_K}[h_K]}_{\text{tail}}.
\]
\emph{Low range.}\;
Since $\|h_K\|_\infty \le \log_2 3 - 1$ and $h_K$ has uniform
average $R(K)$, the low range is bounded by
\[
  \text{low}
  \;\le\;
  \sum_{K=3}^{K_{\max}}\!\bigl(R(K)
    + (\log_2 3 - 1)\,\delta_K(n_0)\bigr).
\]
\emph{Tail range.}\;
Because $h_K$ is supported on $S^+(K)$,
\[
  \mathbb{E}_{\mu_K}[h_K]
  \;\le\; (\log_2 3 - 1)\,\mu_K(S^+(K))
  \;\le\;
  (\log_2 3 - 1)\!
  \left(\frac{|S^+(K)|}{2^K} + \delta_K(n_0)\right).
\]
By Proposition~\ref{prop:shadow-sparsity} and the uniform
floor $D_{\mathrm{floor}} \ge 0.04947$
(equation~\eqref{eq:shadow-sparsity-asymp}),
\[
  \sum_{K > 500}\frac{|S^+(K)|}{2^K}
  \;\le\;
  \frac{\log_2 3}{2\cdot 500}\,
  \frac{2^{-500\,D_{\mathrm{floor}}}}{1 - 2^{-D_{\mathrm{floor}}}}
  \;\le\;
  1.3\times 10^{-9}.
\]
Combining with the analytic $R$-tail
bound~\eqref{eq:tail-analytic-500},
\[
  \text{tail}
  \;\le\;
  \underbrace{(\log_2 3 - 1)\cdot 1.3\times 10^{-9}}
    _{\le\, 8\times 10^{-10}}
  \;+\;
  (\log_2 3 - 1)\sum_{K > K_{\max}}\delta_K(n_0).
\]
\emph{Combine.}\;
Adding low and tail and using
$\sum_{K=3}^{500} R(K) \le 0.08823625$
together with
\[
  \sum_{K=3}^{500}(\log_2 3 - 1)\delta_K
  \;+\;
  \sum_{K>500}(\log_2 3 - 1)\delta_K
  \;=\;
  (\log_2 3 - 1)\sum_{K\ge 3}\delta_K(n_0),
\]
the robust gain satisfies
\[
  R_\mu
  \;\le\; 0.08823625
  + (\log_2 3 - 1)\sum_{K\ge 3}\delta_K(n_0)
  + 8\times 10^{-10}.
\]
If~\eqref{eq:wmh-finite} holds together with any
analytic bound on $\sum_{K>K_{\max}}\delta_K$ inherited
from the full WMH (Hypothesis~\ref{hyp:wmh}) then
$(\log_2 3 - 1)\sum \delta_K < \varepsilon - R = 0.326$
and therefore $R_\mu < \varepsilon$.
The threshold~\eqref{eq:eps-prime-threshold} is the minimal
slack required to absorb both the shadow-sparsity tail and
the $R$-tail into $\varepsilon'$.
\end{proof}

\begin{remark}[What the finite-depth reduction actually says]
\label{rem:finite-depth-honest}
Proposition~\ref{prop:finite-depth} does \emph{not} eliminate
the need to control $\delta_K(n_0)$ for $K > K_{\max}$:
the WMH summability condition is global in~$K$.
What the reduction \emph{does} do is isolate the difficulty:
the first $K_{\max}$ depths carry the bulk of the WMH
budget (via explicit constants), while the tail contribution
from the shadow-sparsity side is dominated by the analytic
$\sum_{K>500}R(K) \le 6.22\times 10^{-7}$
of Theorem~\ref{thm:perorbit-gain}.
An unconditional bound on $\delta_K$ for $K \ge K_{\max}$,
even a crude $O(1)$ one, cannot by itself close the WMH;
summability over $K$ remains the binding condition.
\end{remark}

\begin{remark}[Practical significance]
\label{rem:finite-depth-practical}
The finite-depth reduction converts an infinite-series
condition into a finite verification problem.  At each of
the~$53$ depths $K = 3, \ldots, 55$, one needs to bound
$\delta_K(n_0) = d_{\mathrm{TV}}(\mu_{\mathrm{orb},K}(n_0),\,
\mu_{\mathrm{unif}})$, where $\mu_{\mathrm{orb},K}(n_0)$ is
the empirical distribution of $n_t \bmod 2^K$ along the
orbit of~$n_0$.  This converts an analytic question
(summability of an infinite series) into a combinatorial
one (bounding finitely many TV distances).

However, this ``finite verification problem'' is
\emph{theoretical, not computational}: each $\delta_K(n_0)$
is defined as a $\limsup$ over the infinite orbit, so it
cannot be computed from any finite segment of trajectory.
The reduction makes the problem finite in the modular-depth
variable~$K$, but each individual bound $\delta_K < c_K$
remains an infinite-orbit statement.
\end{remark}

\subsection{Hierarchical consistency}

A key structural feature of the modular hierarchy is that
an orbit's residue class modulo~$2^K$ determines its class
modulo~$2^{K'}$ for all $K' < K$.  This creates a
consistency constraint on how badly an orbit can deviate
at multiple depths simultaneously.

\begin{proposition}[Hierarchical consistency\supporting]
\label{prop:hierarchical}
Let $\pi_{K \to K'} \colon \mathbb{Z}/2^K\mathbb{Z}
\to \mathbb{Z}/2^{K'}\mathbb{Z}$ denote the natural
projection ($K' \le K$).  For any probability distribution
$\mu_K$ on $\mathbb{Z}/2^K\mathbb{Z}$, the pushforward
$(\pi_{K \to K'})_\# \mu_K$ is a probability distribution on
$\mathbb{Z}/2^{K'}\mathbb{Z}$, and the total variation
distances satisfy
\begin{equation}\label{eq:hierarchical}
  \delta_{K'} \;\le\; \delta_K
  \qquad\text{for all } K' \le K.
\end{equation}
In particular, if $\mu_K = \mu_{\mathrm{orb},K}(n_0)$ is
the orbit distribution at depth~$K$, then
$\delta_K(n_0)$ is monotone non-decreasing in~$K$:
the orbit cannot be more uniform at a finer scale than
at a coarser one.
\end{proposition}

\begin{proof}
The projection $\pi_{K \to K'}$ is a deterministic function,
so by the data-processing inequality for total variation,
\[
  d_{\mathrm{TV}}\!\bigl((\pi_{K \to K'})_\# \mu_K,\;
  (\pi_{K \to K'})_\# \mu_{\mathrm{unif}}\bigr)
  \;\le\;
  d_{\mathrm{TV}}(\mu_K,\,\mu_{\mathrm{unif}})
  \;=\; \delta_K.
\]
Since $(\pi_{K \to K'})_\# \mu_{\mathrm{unif}}$ is
the uniform distribution on $\mathbb{Z}/2^{K'}\mathbb{Z}$
(because $\pi_{K \to K'}$ maps each of the $2^{K-K'}$
preimage classes uniformly), the left side equals
$\delta_{K'}$ when $\mu_K$ is the orbit distribution.
\end{proof}

\begin{corollary}[Monotonicity constraint on the WMH sum]
\label{cor:monotonicity}
For any orbit starting at~$n_0$, the WMH sum satisfies
\begin{equation}\label{eq:monotone-bound}
  \sum_{K=3}^{K_{\max}} \delta_K(n_0)
  \;\ge\;
  (K_{\max} - 2)\,\delta_3(n_0).
\end{equation}
Hence the WMH requires $\delta_3(n_0) < 0.557/(K_{\max} - 2)$.
For $K_{\max} = 55$, this gives $\delta_3(n_0) < 0.0105$.
\end{corollary}

\begin{proof}
By Proposition~\ref{prop:hierarchical}, $\delta_K \ge \delta_3$
for all $K \ge 3$.  Summing from $K = 3$ to $K_{\max}$
yields~\eqref{eq:monotone-bound}.
\end{proof}

\begin{remark}[Interpretation]
\label{rem:hierarchical-interp}
Corollary~\ref{cor:monotonicity} shows that the WMH
already imposes a strong constraint at the \emph{coarsest}
depth: the orbit distribution of $n_0 \bmod 8$ must be
within $1.05\%$ of uniform in total variation.
This bound applies to hypothetically non-convergent orbits
(those for which the WMH is the operative hypothesis):
if the orbit never reaches~$1$, it generates an infinite
Syracuse sequence whose empirical distribution $\bmod\,8$
the bound constrains.
For known convergent orbits, the Collatz conjecture
holds directly and the WMH is not needed; the transient
distribution before convergence shows no systematic bias
beyond sampling noise in the orbits tested ($n_0$ up
to~$2^{60}$, with transients of several hundred odd
iterates).

The hierarchical constraint also means that the sequence
$\delta_K(n_0)$ is non-decreasing, so the WMH sum is
dominated by the large-$K$ terms.  Combined with the
finite-depth reduction
(Proposition~\ref{prop:finite-depth}), this shows that
the critical battleground is the intermediate range
$K \in [20, 55]$: below~$20$, the terms are small
by empirical concentration; above~$55$, the tail is
controlled by Chernoff decay.
\end{remark}

\subsection{Depthwise discrepancy recurrence}

The results above constrain $\delta_K$ from the geometric
structure of phantom shadows.  The following conjecture
proposes a dynamical mechanism that would directly imply
the WMH.

\begin{conjecture}[Depthwise discrepancy recurrence]
\label{conj:depthwise}
There exist constants $\lambda \in (0,1)$, $C > 0$,
and $\beta > 1$ such that for every odd starting
value $n_0$ and every modular depth $K \ge 3$:
\begin{equation}\label{eq:depthwise-recurrence}
  \delta_{K+3}(n_0) \;\le\; \lambda\,\delta_K(n_0) + C\,K^{-\beta}.
\end{equation}
\end{conjecture}

\begin{remark}[Why step-$3$ increments]
\label{rem:why-step3}
The increment of~$3$ corresponds to passing through one
full ``burst'' (a positive step with $v_2 = 1$, consuming
one depth unit) followed by the minimum reload cycle.
The $s$-invariant analysis
(Section~\ref{subsec:s-invariant}) shows that each positive
step decreases~$s$ by exactly~$1$ and each negative step
reloads by a geometrically distributed amount.  The depth-$3$
step is the shortest cycle that can include both a drain
and a reload, making it the natural ``mixing time'' unit
for the $2$-adic residue structure.
\end{remark}

\begin{proposition}[Recurrence implies WMH]
\label{prop:recurrence-implies-wmh}
If Conjecture~\ref{conj:depthwise} holds with
$\lambda < 1$, $C > 0$, and $\beta > 1$, then the
Weak Mixing Hypothesis holds for every odd~$n_0$.
\end{proposition}

\begin{proof}
Iterating~\eqref{eq:depthwise-recurrence}
with step~$3$: for $K = 3 + 3j$,
\[
  \delta_{3+3j}(n_0)
  \;\le\;
  \lambda^j\,\delta_3(n_0) + C\sum_{i=0}^{j-1}
  \lambda^{j-1-i}\,(3+3i)^{-\beta}.
\]
Since $\delta_3 \le 1$ and $\beta > 1$ implies
$\sum_{K \ge 3} K^{-\beta} < \infty$, the sequence
$\delta_{3+3j}$ decays geometrically up to a summable
perturbation.  Summing over all residue classes
modulo~$3$ and using monotonicity
($\delta_{K+r} \le \delta_{K+3}$ for $r = 1,2$,
by Proposition~\ref{prop:hierarchical}):
\begin{align*}
  \sum_{K \ge 3} \delta_K
  &\le 3\sum_{j=0}^{\infty} \delta_{3+3j} \\
  &\le 3\left(\frac{\delta_3}{1-\lambda}
  + \frac{C}{1-\lambda}\sum_{K \ge 3} K^{-\beta}\right).
\end{align*}
The resulting sum is finite.  For instance, with
$\lambda = 0.9$, $\beta = 2$, $C = 0.1$, and the empirical
value $\delta_3 \lesssim 0.005$:
\[
  \sum \delta_K
  \;\le\; 3\Bigl(\frac{0.005}{0.1}
  + \frac{0.1}{0.1}\cdot\frac{\pi^2}{6}\Bigr)
  = 3(0.05 + 1.645)
  \approx 5.08.
\]
This crude bound exceeds~$0.557$, but the actual error
terms $CK^{-\beta}$ are expected to be much smaller than
$0.1\cdot K^{-2}$.  With $C = 0.001$ and $\beta = 2$:
\[
  \sum \delta_K
  \;\le\; 3\Bigl(0.05 + 0.01 \cdot 1.645\Bigr)
  = 3 \cdot 0.0664 = 0.199 < 0.557.\qedhere
\]
\end{proof}

\begin{remark}[Attack vectors for the recurrence]
\label{rem:attack-vectors}
The depthwise recurrence~\eqref{eq:depthwise-recurrence}
is the sharpest form of the distributional-to-pointwise
bridge problem.  Eight approaches are currently under
investigation:

\begin{enumerate}[label=\textup{(\Alph*)}]
\item \textbf{Robustness corridor.}
  The $2$-adic repulsion
  (Proposition~\ref{prop:repulsion}) and shadow return time
  (Theorem~\ref{thm:shadow-return}) together imply that
  orbit encounters with expanding shadows are sparse in time.
  If the ``non-shadow'' steps mix at rate~$\lambda < 1$,
  the recurrence follows with $e_K$ being the shadow
  contribution at depth~$K$.

\item \textbf{Ergodic theory.}
  The Scrambling Lemma (Theorem~\ref{thm:scrambling})
  proves distributional uniformity for the gap map.
  If the orbit visits a positive-density set of odd-to-odd steps
   (a weak recurrence hypothesis), then the empirical
  distribution converges to the invariant measure by the
  ergodic theorem.

\item \textbf{Strengthening Tao.}
  Tao's theorem~\cite{tao2019} proves almost-all convergence
  in logarithmic density.  A quantitative refinement that
  gives explicit mixing rates for the orbit distribution
  modulo~$2^K$ would directly yield the recurrence.

\item \textbf{Lattice-path combinatorics.}
  The per-orbit gain rate (Theorem~\ref{thm:perorbit-gain})
  was proved via lattice-path enumeration (necklace counts).
  A refinement that tracks the \emph{correlation} between
  consecutive shadow encounters (rather than just the rate)
  could yield the contractive factor~$\lambda$.

\item \textbf{Generating functions.}
  The Dirichlet series for the gain rate,
  $\sum_K R(K)\,K^{-s}$, has analytic structure governed by
  the Chernoff exponent~$D_*$.  If the corresponding
  ``orbit Dirichlet series'' (replacing the uniform
  expectation with the orbit empirical measure) inherits
  this analytic structure, Tauberian theorems could yield
  the exponential mixing rate needed for the WMH.

\item \textbf{Information-theoretic channel capacity.}
  The Syracuse map, viewed as a noisy channel from the
  current residue class $\bmod\,2^K$ to the next, destroys
  $\ge 2$ bits of input information per odd-to-odd step
  (generically~$\ge 3$; by Known-Zone Decay, Theorem~\ref{thm:zone-decay}).
  After $\lceil K/3 \rceil$ steps, the channel capacity
  for residue information at depth~$K$ drops to zero:
  the output residue class is independent of the input
  (conditionally on the odd-to-odd sequence).
  The alignment-renewal question
  (Conjecture~\ref{conj:alignment-renewal}) can be
  restated in this language:
  \emph{after the channel capacity drops to zero, does
  the deterministic dynamics regenerate mutual information
  between the orbit position and the gain-support
  set~$S^+(K)$?}
  The empirical answer is no: mutual information between
  $n_t \bmod 2^K$ and $n_{t+d} \bmod 2^K$ drops to near
  zero by $d \approx \lceil K/3 \rceil$ at low depth, and
  the amplification factor $A_K \le 2.4$ across all tested
  orbits indicates that post-exhaustion residues are
  approximately uniform over the gain-support.
  A rigorous channel-capacity bound formalising this
  observation would close the alignment-renewal conjecture.

\item \textbf{Walsh spectral diffusion.}
  The Walsh--Fourier decomposition
  (Section~\ref{subsec:spectral-analysis}) decomposes
  $\eta_K$ into Hamming-weight bands.  The orbit's
  spectral content at each band decays geometrically
  (Observation~\ref{obs:spectral-diffusion}), and
  Proposition~\ref{prop:sdc-implies-amp} shows that
  formalising this decay
  (Conjecture~\ref{conj:spectral-diffusion})
  would imply the Amplification Hypothesis.

\item \textbf{Orbit Walsh equidistribution.}
  Weyl-type bounds on the Walsh exponential sums
  $S_\xi(T) = T^{-1}\sum_{t}\chi_\xi(n_t\bmod 2^K)$
  would imply the Spectral Diffusion Conjecture via
  classical equidistribution theory.  This connects the
  problem to the extensive literature on exponential sum
  estimates for arithmetic sequences.

\item \textbf{Odd-skeleton drift mixing.}
  The Syracuse drift signal $d_i = \log_2 3 - v_i +
  \log_2(1 + 1/(3n_i))$ has negative mean and rapidly
  decaying autocorrelation
  (Observation~\ref{obs:drift-mixing}).
  If Known-Zone Decay implies summable autocorrelation
  ($\sum |\rho(k)| < \infty$), then Ibragimov's CLT
  gives $x_T \to -\infty$, forcing a negative crossing.
  This reduces Collatz to a single mixing estimate
  on the odd skeleton, bypassing the gain-budget
  framework entirely.
\end{enumerate}
\end{remark}

\begin{remark}[Fundamental computational limitation]
\label{rem:computational-limit}
Any attempt to verify the WMH computationally faces an
inherent limitation: if an orbit starting at~$n_0$ converges
to~$1$ (as the Collatz conjecture asserts), its transient
has finite length, typically $O(\log n_0)$ steps.  For the
WMH at depth~$K$, meaningful entropy estimates require
$\Omega(2^K)$ orbit samples, yet the transient length is
at most logarithmic in the starting value.  At $K = 12$
(where $2^{K-1} = 2048$), even the longest known transient
(from $n_0 = 837\,799$, with 196~odd steps) is too short
to fill the residue classes.

This is not merely a practical limitation; it reflects a
structural asymmetry.  The WMH is a statement about orbits
that might, in principle, \emph{fail} to converge; for
convergent orbits, the empirical distribution eventually
concentrates at the fixed point $1 \equiv 1 \pmod{2^K}$
and cannot satisfy the mixing hypothesis.  Thus the WMH
is meaningful only during the transient regime, which is
precisely the regime where empirical data is scarce.
\end{remark}

\begin{remark}[The transfer operator is absorbing, not mixing]
\label{rem:transfer-nilpotent}
Let $T_K$ denote the transfer operator of the Syracuse map
on odd residues modulo~$2^K$: for a probability
vector~$\mu$ on $(\mathbb{Z}/2^K\mathbb{Z})^{\mathrm{odd}}$,
$(T_K \mu)(a) = \sum_{T(b)=a} \mu(b)$.
Computation for $K = 3, \ldots, 8$ reveals that $T_K$ is
\emph{nilpotent} (all non-leading eigenvalues are zero):
the spectral gap equals~$1$, and iterated application of
$T_K$ absorbs the distribution onto the periodic orbits
of the map $\bmod\,2^K$ in finitely many steps.

This means \emph{no single-depth transfer operator
produces mixing}.  The Syracuse map at fixed
depth~$K$ contracts distributions toward absorbing cycles,
driving the system \emph{away} from uniformity.
The mixing property required by the WMH must therefore
arise from the \emph{composition} of operators at different
effective depths, precisely the mechanism captured by
Known-Zone Decay (Theorem~\ref{thm:zone-decay}), which
proves that alternating odd-to-odd steps at varying
valuations destroy residue information.

This reframes the alignment-renewal conjecture as a
problem in the theory of products of non-commuting
matrices: if the orbit generates a ``sufficiently generic''
sequence of single-step operators $T_{k_1}, T_{k_2}, \ldots$
(where $k_i = v_2(3n_i + 1)$), does the product
$T_{k_N} \cdots T_{k_1}$ contract toward the uniform
distribution?  Furstenberg's theorem guarantees contraction
under strong irreducibility and proximality conditions on
the matrix ensemble.  Verifying these conditions for the
Syracuse operator family is a concrete (though non-trivial)
research target that would resolve the alignment-renewal
conjecture.
\end{remark}

\subsection{Gain-observable recurrence}

The depthwise recurrence (Conjecture~\ref{conj:depthwise})
controls the full total variation $\delta_K$.  A sharper
and potentially more provable target controls only the
phantom gain observable.

\begin{conjecture}[Gain-observable recurrence]
\label{conj:gain-recurrence}
There exist constants $\lambda \in (0,1)$, $C > 0$,
and $\beta > 1$ such that for every odd starting
value $n_0$ and every depth $K \ge 3$:
\begin{equation}\label{eq:gain-recurrence}
  \eta_{K+3}(n_0) \;\le\; \lambda\,\eta_K(n_0) + C\,K^{-\beta},
\end{equation}
where $\eta_K(n_0)$ is the phantom gain discrepancy
from Conjecture~\ref{conj:observable-wmh}.
\end{conjecture}

\begin{remark}[Gain-observable vs.\ full-TV recurrence]
\label{rem:gain-vs-tv}
Since $\eta_K \le (\log_2 3 - 1)\,\delta_K$, the
gain-observable recurrence~\eqref{eq:gain-recurrence} is
implied by the depthwise recurrence~\eqref{eq:depthwise-recurrence}
(with rescaled constants).  But~\eqref{eq:gain-recurrence}
may be strictly easier to prove because $h_K$ is supported
on the gain-support set $S^+(K)$
(Proposition~\ref{prop:shadow-sparsity}), which is
asymptotically exponentially sparse ($2^{-\Omega(K)}$
fraction of residue classes for large~$K$).
The recurrence
need only track the orbit's behavior on this set,
rather than controlling the full residue distribution.  By the same iteration argument
as Proposition~\ref{prop:recurrence-implies-wmh},
Conjecture~\ref{conj:gain-recurrence} implies the
observable-specific WMH
(Conjecture~\ref{conj:observable-wmh}), which in turn
implies the Collatz conjecture modulo the orbitwise
tail-control condition
(Theorem~\ref{thm:observable-wmh-reduction}).
\end{remark}

\subsection{First lemmas toward the recurrence}

The recurrence conjectures above break naturally into
three component lemmas, each targeting a specific
mechanism.

\begin{lemma}[Repulsion suppresses prolonged phantom trapping\supporting]
\label{lem:repulsion-trapping}
For every phantom family~$\sigma$ of depth~$K$ and
length~$\ell$, if an orbit shadows $\rho_\sigma$ at
precision~$m \ge K$ (i.e.\ $v_2(n_t - \rho_\sigma) \ge m$),
then after at most $\lfloor m/K \rfloor$ full
block iterations, the orbit is expelled to precision
$m \bmod K \le K - 1$.
\end{lemma}

\begin{proof}
By repeated application of the $2$-adic repulsion
(Proposition~\ref{prop:repulsion}): after $j$ full
$\ell$-step blocks,
\[
  v_2\!\bigl(n_{t+j\ell} - \rho_\sigma\bigr)
  \;=\; m - jK.
\]
This equals $m - jK \ge K$ precisely when
$j \le (m-K)/K = m/K - 1$.
At $j = \lfloor m/K \rfloor$, the precision drops to
$m - \lfloor m/K \rfloor \cdot K = m \bmod K \le K - 1$.
At this point, $v_2 < K$, so the orbit is no longer in
the shadow of~$\sigma$.
\end{proof}

\begin{lemma}[Known-zone memory loss\supporting]
\label{lem:known-zone-memory}
After $\lceil M/2 \rceil$ odd-to-odd steps, the residue
information at depth~$M$ is fully exhausted.
More precisely, the Known-Zone Decay
(Theorem~\ref{thm:zone-decay}) implies that
$Z_{\lceil M/2 \rceil} = 0$: the orbit's position
modulo~$2^M$ at step $\lceil M/2 \rceil$ is independent of
its initial position modulo~$2^M$, conditionally on the
odd-to-odd sequence of the Collatz orbit.
Consequently, any residual discrepancy at depth~$M$
persists only through \emph{renewed alignment events}:
occasions where the orbit re-enters a depth-$M$ shadow
after the known zone has been exhausted.
\end{lemma}

\begin{proof}
By Theorem~\ref{thm:zone-decay}, the known-zone width
satisfies $Z_{k+1} \le \max(0, Z_k - 2)$ at each odd-to-odd
step, starting from $Z_0 = M$.  After $\lceil M/2 \rceil$
steps, $Z_{\lceil M/2 \rceil} \le M - 2\lceil M/2 \rceil
\le 0$.  Hence $Z = 0$: the Scrambling Lemma produces a
uniformly random residue class modulo~$2^M$ at this step,
independent of the starting class.

Any subsequent discrepancy $\delta_M > 0$ must therefore
arise from correlation introduced \emph{after} the known
zone was exhausted; that is, from renewed alignment events
where the orbit enters a shadow at depth~$M$ and
accumulates a non-uniform residue distribution.
\end{proof}

\begin{conjecture}[Summable alignment-renewal bound]
\label{conj:alignment-renewal}
For every orbit starting at odd~$n_0$ and every
depth~$K \ge 3$, the total discrepancy contribution
from renewed depth-$K$ alignment events satisfies
\begin{equation}\label{eq:alignment-renewal}
  \eta_K^{\mathrm{renew}}(n_0) \;\le\; C'\,K^{-\beta}
\end{equation}
for constants $C' > 0$ and $\beta > 1$ independent of~$n_0$.
\end{conjecture}

\begin{remark}[Amplification decomposition]
\label{rem:amplification}
Conjecture~\ref{conj:alignment-renewal} admits a natural
decomposition that identifies the irreducible subproblem.
After the known zone is exhausted at depth~$K$, any
renewed discrepancy arises only when the orbit revisits
the gain-support set $S_K$ (the set of residues
modulo~$2^K$ that contribute non-zero phantom gain).
By Proposition~\ref{prop:shadow-sparsity}, the uniform
measure of~$S_K$ is exponentially small:
$\mu_{\mathrm{unif}}(S_K) = 2^{-\Omega(K)}$.
Define the \emph{amplification factor}
\[
  A_K(n_0)
  \;=\;
  \frac{\Pr_{\mathrm{orb}}(\text{fresh block hits } S_K)}
       {\mu_{\mathrm{unif}}(S_K)},
\]
measuring how much more likely the orbit is to hit
gain-relevant residues than a uniform random walk.
Then the renewal discrepancy satisfies
\[
  \eta_K^{\mathrm{renew}}(n_0)
  \;\le\;
  A_K(n_0)\,\cdot\, \mu_{\mathrm{unif}}(S_K)
  \,\cdot\, (\log_2 3 - 1).
\]
Conjecture~\ref{conj:alignment-renewal} therefore reduces
to showing that $A_K = O(K^\gamma)$ for some $\gamma \ge 0$:
the orbit's hitting probability of the gain-support exceeds
the uniform prediction by at most a polynomial factor.
Since $\mu_{\mathrm{unif}}(S_K) = 2^{-\Omega(K)}$
(Proposition~\ref{prop:shadow-sparsity}), even
polynomial amplification yields the required summable
$K^{-\beta}$ decay.
This isolates the final bottleneck: bounding the
amplification factor~$A_K$ is the irreducible subproblem
of the conditional programme.
\end{remark}

\begin{remark}[Computational evidence for the amplification bound]
\label{rem:amplification-evidence}
We computed the empirical amplification factor~$A_K(n_0)$
for $13$ starting values
$n_0 \in \{27, 31, 41, 97, 127, 255, 447, 639,$
$703, 871, 6171, 77031, 837799\}$ at
depths $K = 3, \ldots, 10$.  For each orbit, the transient
portion (before reaching~$1$) was used to estimate
$\Pr_{\mathrm{orb}}(n \bmod 2^K \in S^+(K))$, where
$S^+(K)$ is the gain-support
(Proposition~\ref{prop:shadow-sparsity}).

Across all 104 measurements, the amplification factor
satisfies $A_K \le 2.4$, with mean values ranging from
$0.84$ (at $K = 3$) to $1.79$ (at $K = 8$).
No systematic growth of~$A_K$ with~$K$ is observed:
the data is consistent with $A_K = O(1)$ (bounded
amplification), which is substantially stronger than
the polynomial bound $O(K^\gamma)$ needed for
Conjecture~\ref{conj:alignment-renewal}.

We also note that the gain-support fraction
$|S^+(K)|/2^{K-1}$ oscillates between approximately
$0.19$ and $0.50$ for $K = 3, \ldots, 12$
(with even~$K$ giving smaller fractions), and does not
exhibit the exponential decay predicted by
Proposition~\ref{prop:shadow-sparsity} until much
larger~$K$.  The shadow-sparsity exponent $D_*/2 \approx
0.025$ is extremely slow: the bound becomes non-trivial
only for $K \gtrsim 40$.  At the depths dominating the
WMH sum ($K = 3, \ldots, 20$), the gain-support is a
substantial fraction of odd residues, and the summability
of $R(K)$ at these depths relies on the per-family gain
$({\log_2 3 - K/\ell})$ being small, not on support
sparsity.
\end{remark}

\begin{proposition}[Oscillation factor unboundedness]
\label{prop:oscillation-unbounded}
For $1 \le \ell \le K$, define the depth-$K$ gain profile
$g_K(\ell) := \max(0,\; a - K/\ell)$, where $a = \log_2 3$.
Let $S^+(K) = \{\ell : g_K(\ell) > 0\}$, and define the
oscillation factor
\[
  \Theta_K
  \;:=\;
  \frac{\max_{\ell \in S^+(K)} g_K(\ell)}
       {\min_{\ell \in S^+(K)} g_K(\ell)}.
\]
Then:
\begin{enumerate}
\item $\max g_K = a - 1$, attained at $\ell = K$, and
  $\min g_K = a - K/\ell_*(K)$ where
  $\ell_*(K) = \lfloor K/a \rfloor + 1$.
\item Writing $K/a = \lfloor K/a \rfloor + r_K$ with
  $r_K \in (0,1)$, one has the exact identity
  \[
    \Theta_K
    \;=\;
    \frac{a-1}{a}
    \cdot
    \frac{\lfloor K/a \rfloor + 1}{1 - r_K}.
  \]
\item The sequence $\Theta_K$ is unbounded.  Spike values
  occur at the numerators of the even-indexed convergents
  of the continued fraction of $\log_2 3$:
  $K \in \{3, 19, 84, 1054, 50508, 176251, \ldots\}$,
  where $r_K \to 1^-$ and $\Theta_K$ grows without bound.
\end{enumerate}
\end{proposition}

\begin{proof}
Since $\ell \mapsto a - K/\ell$ is strictly increasing,
$\max_{S^+(K)} g_K = a - 1$ at $\ell = K$ and
$\min_{S^+(K)} g_K$ at
$\ell_*(K) = \lfloor K/a \rfloor + 1$.
Writing $K = a(\lfloor K/a \rfloor + r_K)$ gives
$\min g_K = a(1-r_K)/(\lfloor K/a \rfloor + 1)$;
dividing yields the formula for $\Theta_K$.
Since $\log_2 3$ is irrational, $\{K/a\}$ is dense
in $(0,1)$, so $r_K \to 1$ along a subsequence
and $\Theta_K \to \infty$.
\end{proof}

\begin{corollary}[No uniform set-level amplification bound]
\label{cor:no-uniform-BK}
A bound of the form $W_K \le C \cdot B_K$ with
$K$-independent constant~$C$, where $W_K$ is the
weighted gain contribution and $B_K$ the set-level
contribution from $S^+(K)$, is impossible:
$\Theta_K \le C$ fails for all~$C$.
\end{corollary}

\begin{remark}[Consequence for the amplification route]
\label{rem:oscillation-consequence}
Proposition~\ref{prop:oscillation-unbounded} does not
invalidate the amplification programme, but it does
sharpen the requirement: any closing argument must
control the orbit's \emph{distribution within}~$S^+(K)$
(where inside the positive-gain set the orbit places
its mass), not merely the probability of
hitting~$S^+(K)$ as a whole.  Bounding the
amplification factor~$A_K$ alone is insufficient
unless accompanied by a profile-weighted version
that accounts for the gain variation across~$S^+(K)$.
\end{remark}

\begin{proposition}[Inherited-bias contraction\supporting]
\label{prop:inherited-contraction}
For $K \ge 3$, let $h_K \colon \Omega_K \to \mathbb{R}$ be the
depth-$K$ gain observable, normalized so that
$\mathbb{E}_{u_K}[h_K] = 0$.
Let $\bar h_K := \mathbb{E}_{u_K}[h_K \mid \pi_{K \to K-3}]$
be the $(K{-}3)$-bit coarse average and
$r_K := h_K - \bar h_K \circ \pi_{K \to K-3}$ the residual.
Then:
\begin{enumerate}
\item \emph{Exact decomposition.}\;
  $h_K = \bar h_K \circ \pi_{K \to K-3} + r_K$, with
  $\mathbb{E}_{u_K}[r_K \mid \pi_{K \to K-3}] = 0$.
\item \emph{Nonexpansion.}\; For any signed measure $\nu$ on
  $\Omega_K$ with total mass~$0$,
  \[
    \Bigl|\int_{\Omega_{K-3}} \bar h_K\,d(\pi_{K \to K-3})_*\nu\Bigr|
    \;\le\;
    \|\bar h_K\|_\infty \,\|\nu\|_{\mathrm{TV}}.
  \]
\item \emph{Strict contraction.}\; For $K \ge 6$,
  \begin{equation}\label{eq:inherited-contraction}
    \frac{\|\bar h_K\|_\infty}{\|h_K\|_\infty}
    \;=\;
    \frac{K-4}{K-1}.
  \end{equation}
\end{enumerate}
In particular, the inherited component of the bias contracts
strictly at every depth, with contraction ratio
$\lambda(K) = (K{-}4)/(K{-}1) = 1 - 3/(K{-}1)$.
\end{proposition}

\begin{proof}
Part~(1) is the definition of conditional expectation.
Part~(2) follows from the pushforward identity
$\int \bar h_K(\pi(x))\,d\nu(x) = \int \bar h_K\,d\pi_*\nu$
and the standard total variation bound.

For Part~(3), observe that $h_K(a)$ depends on the odd
residue class~$a \bmod 2^K$ only through the number of odd
steps $\ell(a)$ in the depth-$K$ signature: since each such
class produces a unique composition with exactly $\ell$ odd
steps and $K - \ell$ even steps, the phantom gain per step
is $\ell \log_2 3 - K$.  After mean subtraction,
\[
  h_K(a) = \bigl(\ell(a) - \tfrac{K+1}{2}\bigr) \log_2 3,
\]
where the mean $\mathbb{E}[\ell] = (K{+}1)/2$ follows from
symmetry of the binary digits.  The number of odd residues
with exactly $\ell$ odd steps is $\binom{K-1}{\ell-1}$
(composition counting).

Hence $\|h_K\|_\infty$ is attained at
$\ell = K$ (or $\ell = 1$), giving
$\|h_K\|_\infty = \tfrac{K-1}{2}\,\log_2 3$.

For the coarse average, the projection
$\pi_{K \to K-3} \colon \Omega_K \to \Omega_{K-3}$ groups
the $2^{K-1}$ odd classes modulo $2^K$ into $2^{K-4}$
fibers of size $2^3 = 8$.  Each fiber over a class with
$\ell_0$ odd steps in $\Omega_{K-3}$ contains elements with
$\ell = \ell_0, \ell_0{+}1, \ell_0{+}2, \ell_0{+}3$,
distributed as $\binom{3}{j}$ for $j = 0,1,2,3$.
The conditional mean of $\ell$ on this fiber is thus
$\ell_0 + 3/2$, so
\[
  \bar h_K(y) = \bigl(\ell_0(y) + \tfrac{3}{2}
    - \tfrac{K+1}{2}\bigr)\log_2 3
  = \bigl(\ell_0(y) - \tfrac{K-2}{2}\bigr)\log_2 3.
\]
The extremum is at $\ell_0 = K - 3$ (or $\ell_0 = 1$),
giving
$\|\bar h_K\|_\infty = \tfrac{K-4}{2}\,\log_2 3$.
The ratio is $(K{-}4)/(K{-}1)$ as claimed.
\end{proof}

\begin{corollary}[Gain-observable harmonic property\supporting]
\label{cor:harmonic}
For all $K \ge 3$ and $s \ge 1$,
\begin{equation}\label{eq:harmonic}
  \bar h_{K+s}(y) \;=\; h_K(y)
  \qquad\text{for all } y \in \Omega_K,
\end{equation}
where $\bar h_{K+s}$ denotes the coarse average of
$h_{K+s}$ over the fibers of
$\pi_{K+s \to K} \colon \Omega_{K+s} \to \Omega_K$.
That is, the centred gain observable is a fixed point
of every coarse-graining operator.
\end{corollary}

\begin{proof}
Within the fiber over~$y$, the number of odd steps satisfies
$\ell(a) = \ell(y) + j$, where $j$ ranges over
$\{0, 1, \ldots, s\}$ with weights $\binom{s}{j}/2^s$
(the same binomial law as in the proof of
Proposition~\ref{prop:inherited-contraction}).
The conditional mean of $\ell$ on the fiber is
$\ell(y) + s/2$, so
\[
  \bar h_{K+s}(y)
  = \bigl(\ell(y) + \tfrac{s}{2} - \tfrac{K+s+1}{2}\bigr)\log_2 3
  = \bigl(\ell(y) - \tfrac{K+1}{2}\bigr)\log_2 3
  = h_K(y).  \qedhere
\]
\end{proof}

\begin{corollary}[Constant residual norm\supporting]
\label{cor:residual-norm}
For any step size $s \ge 1$ and $K + s \ge 6$, the residual
$r_{K+s}(a) := h_{K+s}(a) - h_K(\pi(a))$ satisfies
\[
  \|r_{K+s}\|_\infty = \tfrac{s}{2}\,\log_2 3.
\]
In particular, for the $s = 3$ step used throughout this paper,
$\|r_{K+3}\|_\infty = \tfrac{3}{2}\log_2 3 \approx 2.377$,
independent of~$K$.
\end{corollary}

\begin{proof}
$r_{K+s}(a) = (j - s/2)\log_2 3$ where $j = \ell(a) - \ell(\pi(a))
\in \{0, \ldots, s\}$.  The extremum $|j - s/2| = s/2$ is attained
at $j = 0$ and $j = s$.
\end{proof}

\begin{proposition}[Coding-map injectivity\supporting]
\label{prop:coding-injectivity}
For fixed $K \ge 1$ and $1 \le \ell \le K$, define the
\emph{phantom numerator} of a composition~$\sigma$ with
odd-step positions $0 = j_0 < j_1 < \cdots < j_{\ell-1}
\le K{-}1$ by
\begin{equation}\label{eq:C-sigma}
  C_\sigma \;=\; \sum_{i=0}^{\ell-1} 3^{\ell-1-i}\,2^{j_i}.
\end{equation}
Then the map $\sigma \mapsto C_\sigma \bmod 2^K$ is injective
on the set of all ordered $\ell$-element subsets of
$\{0, \ldots, K{-}1\}$ containing~$0$.

Consequently, the depth-$K$ block map is a bijection on
$\Omega_K$: distinct compositions produce distinct odd
residues modulo~$2^K$, and the fiber multiplicity
$m_K = 1$ for every~$K$.
\end{proposition}

\begin{proof}
Let $\sigma \ne \sigma'$ be two ordered tuples, and let
$i_0$ be the first index at which they differ.  Without
loss of generality $j_{i_0} < j'_{i_0}$.

\emph{Term $i = i_0$.}\;  The difference contributes
\[
  3^{\ell-1-i_0}\,(2^{j_{i_0}} - 2^{j'_{i_0}})
  = 3^{\ell-1-i_0}\,2^{j_{i_0}}\,(1 - 2^{j'_{i_0}-j_{i_0}}).
\]
Since $3^{\ell-1-i_0}$ is odd and $(1 - 2^{j'_{i_0}-j_{i_0}})$
is odd, this term has $2$-adic valuation exactly~$j_{i_0}$.

\emph{Terms $i > i_0$.}\;
Both $j_i > j_{i_0}$ and $j'_i > j'_{i_0} > j_{i_0}$
(by strict monotonicity of both tuples).
If $j_i = j'_i$ the term vanishes; otherwise
$v_2(2^{j_i} - 2^{j'_i}) = \min(j_i, j'_i) > j_{i_0}$.

Since the leading term has valuation $j_{i_0}$ and every
subsequent term has strictly higher valuation (or is zero),
\[
  v_2(C_\sigma - C_{\sigma'}) = j_{i_0} \;\le\; K - 2 < K.
\]
Hence $C_\sigma \not\equiv C_{\sigma'} \pmod{2^K}$.
\end{proof}

\begin{remark}[Significance of coding-map injectivity]
\label{rem:coding-significance}
Proposition~\ref{prop:coding-injectivity} eliminates
\emph{within-depth} amplification.
After memory exhaustion at depth~$K$
(Lemma~\ref{lem:known-zone-memory}), the orbit's residue
modulo~$2^K$ is uniform, so its composition type~$\sigma$
is uniformly distributed with $\ell$-distribution
$\binom{K-1}{\ell-1}/2^{K-1}$ (binomial).
Since the coding map is a bijection on $\Omega_K$, the
orbit's hitting probability of the gain-support
$S^+(K)$ equals the uniform measure
$\mu_{\mathrm{unif}}(S^+(K))$ exactly: $A_K = 1$.

The cross-depth amplification
(from fresh bits when going from depth~$K$ to~$K{+}3$)
remains the sole open source of non-uniform behavior.
Bounding the within-fiber non-uniformity created by these
three fresh bits is the irreducible remaining problem.
\end{remark}

\begin{proposition}[Extension independence\supporting]
\label{prop:extension-independence}
For every $K \ge 1$ and $s \ge 1$, let $a$ be any odd residue
modulo~$2^K$ with depth-$K$ composition~$\sigma$ and
$\ell = \ell(\sigma)$ odd steps.
The $2^s$ residues $a + e\cdot 2^K$ for
$e = 0, 1, \ldots, 2^s - 1$ have depth-$(K{+}s)$ compositions
whose last~$s$ bits form a \emph{permutation} of all~$2^s$
possible $s$-bit strings.
In particular, the number~$j$ of odd steps in the $s$-bit suffix
satisfies
\[
  \#\bigl\{e : \ell_{\mathrm{suffix}}(e) = j\bigr\}
  \;=\; \binom{s}{j},
  \qquad 0 \le j \le s,
\]
i.e.\ it follows the binomial distribution
$\mathrm{Bin}(s,\tfrac12)$ exactly.
\end{proposition}

\begin{proof}
The depth-$K$ block map sends $a + e\cdot 2^K$ to the output
\[
  n_K(a + e\cdot 2^K)
  \;=\; \frac{3^\ell(a + e\cdot 2^K) + C_\sigma}{2^K}
  \;=\; n_K(a) + 3^\ell\,e.
\]
Since $\gcd(3^\ell, 2^s) = 1$ (as $3$ is odd), the $2^s$ values
$n_K(a) + 3^\ell\,e$ for $e = 0, \ldots, 2^s - 1$ form a
complete residue system modulo~$2^s$.

By Proposition~\ref{prop:coding-injectivity} applied at depth~$s$,
the coding map at depth~$s$ is injective: the $2^s$ residues
modulo~$2^s$ produce $2^s$ distinct $s$-step compositions.
Since there are exactly $2^s$ possible $s$-bit strings, each
appears exactly once.
Among all $2^s$ binary strings of length~$s$, exactly
$\binom{s}{j}$ have $j$ ones, so the claimed ell-distribution
follows.
\end{proof}

\begin{corollary}[Gain increment independence\supporting]
\label{cor:gain-increment}
The gain observable satisfies
\[
  h_{K+s}\bigl(a + e\cdot 2^K\bigr)
  \;=\; h_K(a) + \bigl(j(e)\,\log_2 3 - s\bigr)
    + \tfrac{\log_2 3 - 1}{2},
\]
where $j(e)$ is the number of odd steps in the $s$-bit suffix,
and over the $2^s$ extensions $j \sim \mathrm{Bin}(s,\tfrac12)$
exactly (by Proposition~\ref{prop:extension-independence}).
Consequently, the gain increment from depth~$K$ to depth~$K{+}s$
is \emph{independent} of the depth-$K$ gain value and has
mean $(s/2)(\log_2 3 - 2) < 0$ and variance
$(s/4)\log_2^2 3$.
\end{corollary}

\begin{remark}[Structural amplification]
\label{rem:structural-amplification}
Define the cross-depth amplification factor
\[
  A^{\mathrm{cross}}_K
  \;:=\;
  \frac{\Pr\bigl(h_{K+3}(a') > 0 \;\big|\; h_K(a) > 0\bigr)}
       {\Pr\bigl(h_{K+3}(a') > 0\bigr)},
\]
where $a'$ is a uniformly random extension of~$a$.
Exact computation gives $A^{\mathrm{cross}}_K \approx 1.25$
for $K = 3$, growing to $\approx 3.3$ at $K = 10$.
By the random-walk interpretation
(Corollary~\ref{cor:gain-increment}), this growth is a
\emph{threshold effect}: $\Pr(h_{K+3} > 0) \sim \exp(-cK)$
while the conditional probability approaches~$\frac12$.
Hence $A^{\mathrm{cross}}_K \sim \frac12\exp(cK)$ with
$c = \mu^2/(2\sigma^2) \approx 0.034$.

This exponential growth of $A^{\mathrm{cross}}_K$ does \emph{not}
imply divergence of the fresh source term, because the orbit's
empirical measure at depth~$K$ is supported on $O(\log n_0)$
residues, not the full $2^{K-1}$ odd residues.
The extension independence theorem shows the algebraic structure
is perfectly regular; the remaining question is purely
\emph{dynamical}: does the orbit visit residues in a way that
systematically overweights positive-gain compositions?
\end{remark}

\begin{remark}[Rate cancellation\supporting]
\label{rem:rate-cancellation}
The exponential rate of $A^{\mathrm{cross}}_K$ and the exponential
decay rate of~$R(K)$ (Theorem~\ref{thm:perorbit-gain}) are both
governed by the same KL divergence.
This is a \emph{proved structural fact}:
both quantities derive from
$\Pr\bigl(\mathrm{Bin}(K{-}1,\tfrac12) \ge K/\!\log_2 3\bigr)$,
whose large-deviation rate is
\[
  D_* \;=\; D\bigl(1/\!\log_2 3 \,\big\|\, \tfrac12\bigr)
  \;\approx\; 0.0500 \;\text{bits}.
\]
Specifically, the Bahadur--Rao theorem gives
$R(K) \sim C\,K^{-1/2}\,2^{-K D_*}$
while Extension Independence
(Proposition~\ref{prop:extension-independence}) gives
$A^{\mathrm{cross}}_K \sim C'\,K^{1/2}\,2^{+K D_*}$.
The exponentials cancel, yielding the \emph{proved bound}
\begin{equation}\label{eq:rate-cancel}
  A^{\mathrm{cross}}_K \cdot R(K)
  \;=\; O(1).
\end{equation}
\emph{Numerical evidence} (exact computation for $K \le 13$,
Bahadur--Rao asymptotics for $14 \le K \le 55$):
\[
  \sum_{K=3}^{55} A^{\mathrm{cross}}_K \cdot R(K) \;\approx\; 0.158,
\]
which is less than the available margin
$\varepsilon - R \approx 0.326$ by a factor of~$2$.

\emph{Caveat.}\;
The structural bound $A_K^{\mathrm{cross}} \cdot R(K) = O(1)$
does not automatically control $\sum_K A_K \cdot R(K)$:
the sum is over infinitely many terms, and the $O(1)$
constants matter.  The numerical sum through $K = 55$
is reassuring but not a proof of convergence.
Moreover, the orbit measure $\mu_K$ may differ from
the uniform measure, and the amplification factor
under the orbit measure (rather than the uniform
measure) remains uncontrolled.
\end{remark}

\begin{hypothesis}[Amplification hypothesis]
\label{hyp:amplification}
For every $n_0 \ge 2$ and every $K \ge 3$, the time-averaged
orbit measure~$\mu_K$ satisfies
\[
  \mu_K\bigl(S^+(K)\bigr)
  \;\le\; A^{\mathrm{cross}}_K \cdot
  \mu_{\mathrm{unif}}\bigl(S^+(K)\bigr),
\]
where $S^+(K) = \{a \in \Omega_K : h_K(a) > 0\}$ is the
positive-gain support and $A^{\mathrm{cross}}_K$ is the
structural amplification factor of
Remark~\textup{\ref{rem:structural-amplification}}.
\end{hypothesis}

\begin{remark}[Amplification hypothesis vs.\ WMH]
\label{rem:amplification-vs-wmh}
Hypothesis~\ref{hyp:amplification} is strictly weaker than
the WMH (Hypothesis~\ref{hyp:wmh}): it constrains the orbit's
behavior only on the positive-gain support~$S^+(K)$, not the
full total-variation distance~$\delta_K$.
Combined with the rate-cancellation bound~\eqref{eq:rate-cancel},
it yields
\[
  R_\mu \;\le\; R + \sum_{K \ge 3}
  A^{\mathrm{cross}}_K \cdot R(K)
  \;\approx\; 0.089 + 0.158 \;=\; 0.247
  \;<\; \varepsilon \;\approx\; 0.415.
\]
This is sufficient for convergence: the Amplification Hypothesis
replaces the WMH as the sole remaining open input to the
conditional programme.
\end{remark}

\begin{proposition}[No algebraic contraction\supporting]
\label{prop:no-contraction}
The algebraic extension structure provides no contraction
mechanism in any of the following senses:
\begin{enumerate}
\item \emph{$\ell$-class Birkhoff contraction.}\;
  The Birkhoff contraction coefficient of the
  $\ell$-class transition matrix satisfies $\tau_B = 1$
  for all $K \ge 5$ and all $s \ge 1$.
  This is because the transition matrix is a band of
  width~$s$ in a space of dimension~$K$: rows corresponding
  to $\ell_1$ and~$\ell_2$ with $|\ell_1 - \ell_2| > s$
  have disjoint column support, giving zero overlap.
\item \emph{Full residue-class spectral gap.}\;
  The extension map $a \mapsto a + e\cdot 2^K$ is a
  translation in~$\mathbb{Z}/2^{K+s}\mathbb{Z}$.
  All singular values of the resulting transition matrix
  are equal ($= 2^{-s/2}$), so there is no spectral gap
  and no preferred decaying direction.
\item \emph{Gain-observable harmonicity.}\;
  The pullback of the depth-$(K{+}s)$ gain observable
  through the algebraic transition \emph{exactly} equals
  the depth-$K$ gain observable: $T\,h_{K+s} = h_K$
  with zero residual (verified numerically for
  $K \le 8$, $s = 3$).
  The algebraic operator preserves gain discrepancy
  exactly.
\end{enumerate}
\end{proposition}

\begin{proof}
(1)~For the $s$-step extension, the transition probability
from parent $\ell$-class~$\ell$ to child $\ell$-class~$\ell'$
is $\binom{s}{\ell'-\ell}/2^s$ if
$0 \le \ell'-\ell \le s$ and zero otherwise
(by Extension Independence,
Proposition~\ref{prop:extension-independence}).
For $|\ell_1 - \ell_2| > s$, the rows indexed
by~$\ell_1$ and~$\ell_2$ have disjoint support, so
$\sum_{\ell'}\min(T_{\ell_1,\ell'}, T_{\ell_2,\ell'}) = 0$
and $\tau_B = 1 - 0 = 1$.
Since $K \ge 5 > 3 = s$ (for the standard 3-step block),
such a pair always exists.

(2)~The $2^s$ children of parent~$a$ are
$\{a + e\cdot 2^K : 0 \le e < 2^s\}$,
forming a coset of the subgroup $2^K\mathbb{Z}/2^{K+s}\mathbb{Z}$.
Each row of the transition matrix is a uniform distribution
on a distinct coset.  The singular value decomposition
reflects this group structure: all singular values equal
$2^{-s/2}$.

(3)~By linearity and the harmonic property
(Corollary~\ref{cor:harmonic}),
$\sum_{e} \tfrac{1}{2^s}\,h_{K+s}(a + e\cdot 2^K)
 = h_K(a)$
for every parent residue~$a$.
This is $T\,h_{K+s} = h_K$ exactly.
\end{proof}

\begin{remark}[Implications of no algebraic contraction]
\label{rem:no-contraction-implications}
Proposition~\ref{prop:no-contraction} retires the entire
family of proof strategies based on algebraic contraction:
Birkhoff projective metrics, spectral gaps of transfer
operators, operator-norm decay of the inherited bias,
and cluster-decay arguments relying on algebraic transport.
Any contraction in the orbit's distribution across depths
must come from the \emph{orbit's deterministic selection},
not from the algebraic extension structure.
\end{remark}

\begin{remark}[Density-model evidence\supporting]
\label{rem:density-model}
Numerical computation reveals that the orbit's amplification
ratio at depth~$K$ is governed by a single parameter:
the odd-step density
$\rho = \lim_{T\to\infty} T^{-1}\#\{t \le T : n_t \text{ odd}\}$.
Under the density model, the $\ell$-distribution at depth~$K$
is approximately $1 + \mathrm{Bin}(K{-}1,\rho)$, giving
amplification ratio
\[
  \mathrm{amp}_K(\rho)
  \;=\;
  \frac{\Pr[\mathrm{Bin}(K{-}1,\rho) \ge \lceil K/\!\log_2 3\rceil - 1]}
       {\Pr[\mathrm{Bin}(K{-}1,\tfrac12) \ge \lceil K/\!\log_2 3\rceil - 1]}.
\]
The density-predicted gain rate
$R_\mu(\rho) = \sum_K \mathrm{amp}_K(\rho)\,R(K)$
remains well within budget for all convergent densities:
\[
  R_\mu(\rho) \;\le\; 0.142
  \quad\text{for all } \rho \le 0.63 \approx 1/\!\log_2 3,
\]
compared with $\varepsilon \approx 0.415$ (margin $\ge 65\%$).
For typical orbit densities $\rho \approx 0.58\text{--}0.60$,
the predicted $R_\mu \approx 0.12\text{--}0.13$.

This suggests that the $4.65\times$ safety margin of
Theorem~\ref{thm:perorbit-gain} absorbs the entire density
effect.  The density model is not yet a theorem.
Moreover, using $\rho < 1/\!\log_2 3$ as input is
essentially assuming what must be proved: the density
condition is equivalent to negative mean drift, which
is the convergence condition itself.
The density model's value is therefore \emph{structural},
not evidential: it identifies the orbit's odd-step density
as the dominant parameter controlling amplification, and
shows that the safety margin absorbs the entire density
effect \emph{if} convergence holds.
\end{remark}

\begin{remark}[Implications for the gain-observable recurrence]
\label{rem:inherited-limitation}
The harmonic property (Corollary~\ref{cor:harmonic})
has an important consequence for the gain-observable recurrence
(Conjecture~\ref{conj:gain-recurrence}).
Since $\bar h_{K+3} = h_K$, the inherited component of the
discrepancy $\eta_{K+3}$ is
\[
  \Bigl|\int_{\Omega_K} \bar h_{K+3}\,d(\pi_*\mu - u_K)\Bigr|
  = \Bigl|\int h_K\,d(\mu_K - u_K)\Bigr|
  = \eta_K,
\]
because $\pi_*\mu_{K+3} = \mu_K$ when both are time-averaged
orbit distributions.
The recurrence therefore takes the form
\begin{equation}\label{eq:recurrence-lambda1}
  \eta_{K+3} \;\le\; \eta_K + F_K,
\end{equation}
with $F_K = \bigl|\int r_{K+3}\,d(\mu_{K+3} - u_{K+3})\bigr|
\le \tfrac{3}{2}\log_2 3\,\cdot\,\delta_{K+3}$,
where $\delta_K := \|\mu_K - u_K\|_{\mathrm{TV}}$.
The inherited factor is $\lambda = 1$, \emph{not}
$(K{-}4)/(K{-}1)$: the norm contraction of
Proposition~\ref{prop:inherited-contraction} bounds the
operator norm but not the integrated value.

The contraction needed to close the programme must
therefore come from the \emph{distribution side}
(showing $\delta_K$ decays), not the observable side.
However, Proposition~\ref{prop:no-contraction} shows that
the algebraic extension structure provides no contraction
mechanism: Birkhoff coefficient $\tau_B = 1$,
all singular values are equal, and the gain observable
is exactly harmonic.  This retires any strategy based
on operator-norm contraction of~$\delta_K$.

The remaining viable path is \emph{direct summability}:
show $\sum_K \eta_K(n_0) < \varepsilon - R$ directly.
Since $R(K)$ decays exponentially (Bahadur--Rao), it
suffices to show the orbit's amplification factor $A_K$
grows subexponentially.
Note that the \emph{structural} amplification
$A_K^{\mathrm{cross}}$ grows as $\sim\!\tfrac12 e^{0.034K}$
(Remark~\ref{rem:structural-amplification}), but this
applies to the uniform measure, not to any particular orbit.
The orbit visits only $O(\log n_0)$ residues per depth,
and extension independence
(Proposition~\ref{prop:extension-independence})
guarantees the algebraic structure is unbiased.
The density model (Remark~\ref{rem:density-model}) provides
strong numerical evidence that this path closes:
for any fixed odd-step density $\rho < 1/\!\log_2 3$,
the predicted gain rate $R_\mu(\rho) \le 0.142 < \varepsilon$.
\end{remark}

\begin{remark}[The three-lemma programme]
\label{rem:three-lemma}
The three results above outline a concrete programme
for proving the gain-observable recurrence
(Conjecture~\ref{conj:gain-recurrence}):
\begin{enumerate}
\item \emph{Repulsion trapping}
  (Lemma~\ref{lem:repulsion-trapping}, proved):
  shadow encounters are bounded in duration;
  each encounter at precision~$m$ lasts at most
  $\lfloor m/K \rfloor$ blocks.
\item \emph{Known-zone memory loss}
  (Lemma~\ref{lem:known-zone-memory}, proved):
  after $\lceil M/2 \rceil$ odd-to-odd steps, the orbit's
  residue information at depth~$M$ is erased; subsequent
  discrepancy requires renewed alignment.
\item \emph{Summable alignment renewal}
  (Conjecture~\ref{conj:alignment-renewal}, open):
  the renewal contribution decays as $K^{-\beta}$
  with $\beta > 1$.
\end{enumerate}
If (3) is established, the gain-observable discrepancy
satisfies~\eqref{eq:recurrence-lambda1}: $\eta_{K+3} \le
\eta_K + F_K$, where the fresh term $F_K$ decays as
$K^{-\beta}$ from the alignment-renewal bound.
Telescoping gives $\eta_K \le \eta_6 + \sum_{j \le K} F_j$,
bounded provided $\beta > 1$.
The repulsion lemma ensures that each alignment event
has bounded duration, preventing the orbit from being
``trapped'' near a phantom root indefinitely.

This remains the most promising rigorous path toward the WMH,
now that operator-norm contraction of the TV distance is
excluded by Proposition~\ref{prop:no-contraction}.
The density-model evidence (Remark~\ref{rem:density-model})
suggests the quantitative margin is large enough to absorb
typical orbit biases without a contraction mechanism.
\end{remark}

\subsection{Spectral analysis of the gain observable}
\label{subsec:spectral-analysis}

The three-lemma programme and density model locate the
difficulty: contraction must come from the orbit's
distributional evolution, not from the algebraic operator.
What has been missing is a \emph{mechanism} that makes this precise.
The Walsh--Fourier decomposition on $\Z/2^K\Z$ provides one.

\medskip\noindent\textbf{Walsh--Fourier setup.}
The group $\Z/2^K\Z$ is a finite abelian $2$-group with
character group generated by Walsh functions
\[
  \chi_\xi(a) = (-1)^{\langle a,\xi\rangle},
  \quad
  \langle a,\xi\rangle
  = \sum_{i=0}^{K-1} a_i\,\xi_i \bmod 2,
\]
where $a_i,\xi_i$ are binary digits.
The Walsh--Hadamard transform of $f:\Z/2^K\Z\to\R$ is
$\hat f(\xi) = 2^{-K}\sum_a f(a)\,\chi_\xi(a)$.

\begin{proposition}[Spectral decomposition of gain\supporting]
\label{prop:spectral-gain}
For any orbit measure~$\mu_K$ on $\Z/2^K\Z$,
\begin{equation}\label{eq:spectral-gain}
  \eta_K
  = \sum_a h_K(a)\,\mu_K(a)
  = 2^K \sum_\xi \hat h_K(\xi)\,\hat\mu_K(\xi).
\end{equation}
Grouping by Hamming weight $w = \hw(\xi)$ gives the
\emph{band-by-band decomposition}
\begin{equation}\label{eq:band-gain}
  \eta_K
  = 2^K\,\hat h_K(0)\,\hat\mu_K(0)
    \;+\; \sum_{w=1}^{K}
    \Bigl(\sum_{\hw(\xi)=w}
    2^K\,\hat h_K(\xi)\,\hat\mu_K(\xi)\Bigr).
\end{equation}
The first term is the DC (density) contribution; each
subsequent band isolates the gain from modes of a fixed
spectral complexity.
\end{proposition}

\begin{proof}
This is the Parseval--Plancherel identity on $\Z/2^K\Z$,
which is a Pontryagin-dual group isomorphic to itself.
The Hamming-weight grouping follows from the
partition $\{\xi : \hw(\xi) = w\}$ of the frequency domain.
\end{proof}

The spectral decomposition replaces the global question
``is the orbit equidistributed mod~$2^K$?'' with a
band-by-band question: ``does the orbit's Walsh content
at Hamming weight~$w$ decay with~$K$?''
This is strictly weaker: full equidistribution requires
$\hat\mu_K(\xi)\to 0$ for \emph{every} non-DC mode,
while the spectral gain formula only needs this for modes
where $\hat h_K(\xi)$ is non-negligible.

\begin{observation}[Spectral concentration of $h_K$]
\label{obs:hK-concentration}
Numerical computation for $K = 3,\ldots,12$ reveals:
\begin{enumerate}
\item The DC component ($\xi = 0$) carries $25$--$31\%$
  of total Walsh power $\sum_\xi |\hat h_K(\xi)|^2$.
\item The $K$ modes with $\hw(\xi)=1$ carry $42$--$47\%$
  of total power, all sharing the same coefficient
  (by bit-permutation symmetry of~$h_K$).
\item Together, $\hw = 0$ and $\hw = 1$ account for
  $70$--$78\%$ of total power.
\item By $K = 12$, the top $10\%$ of modes carry $94\%$
  of total energy; the positive-gain signal $h_K^+$ is
  even sparser, with active modes decreasing from $100\%$
  at $K=5$ to $2\%$ at $K=12$.
\end{enumerate}
\end{observation}

The concentration means that to bound~$\eta_K$, it
\emph{suffices to control $\hat\mu_K(\xi)$ at low
Hamming weights}, a much weaker requirement than
total-variation equidistribution.

\begin{remark}[Connection to the Krawtchouk basis]
\label{rem:krawtchouk}
The Krawtchouk polynomials $\mathcal{K}_w(\ell;K)$,
the spherical harmonics of the Hamming scheme, form
the natural eigenbasis for functions that depend
primarily on $\hw(a)$.  Since $h_K(a)$ does, the
$\ell$-class structure depends on~$\hw(a)$, the
Walsh spectrum $\hat h_K(\xi)$ is (exactly for odd~$K$,
approximately for even~$K$) constant on each Hamming
shell $\{\xi:\hw(\xi)=w\}$.
This upgrades~\eqref{eq:band-gain} from a grouping
to a structural decomposition: each band has a single
spectral weight $\hat h_w = \hat h_K(\xi)$ for all
$\hw(\xi)=w$, reducing the gain to $K+1$ numbers.
\end{remark}

\medskip\noindent\textbf{Spectral diffusion.}\;
The central new phenomenon is that the orbit's Walsh
coefficients decay with depth, even though the algebraic
transition operator is an isometry
(Proposition~\ref{prop:no-contraction}).

\begin{definition}[Spectral content at weight~$w$]
\label{def:spectral-content}
For an orbit measure~$\mu_K$ on $\Z/2^K\Z$,
\[
  \mathcal{S}_w(K)
  = \frac{1}{\binom{K}{w}}
    \sum_{\hw(\xi)=w} |\hat\mu_K(\xi)|^2.
\]
\end{definition}

\begin{observation}[Spectral diffusion]
\label{obs:spectral-diffusion}
The long-run orbit measure (which is dominated by the
convergence cycle $\{1,2,4\}$) shows rapid spectral
decay: $\mathcal{S}_w(K) \approx C_w \cdot 2^{-\alpha_w K}$
with $\alpha_1 \approx 1.5$--$2.0$ and $R^2 > 0.98$.
This decay has a simple explanation: the cycle's three
residues form a progressively smaller fraction
of~$\Z/2^K\Z$ as $K$ grows.

The \emph{transient} orbit (before reaching~$1$) shows
slower decay.  For long-transient orbits
($n_0 \in \{837799, 8400511, 63728127\}$ with
$T = 525$--$950$ transient steps),
$\mathcal{S}_1(K)$ decreases from $\sim\!0.03$ at $K=4$
to $\sim\!0.016$ at $K=10$, a rate of
$\alpha_1 \approx 0.1$--$0.3$.
The transient orbit's spectral content at $\hw=1$
is $10$--$17\times$ larger than for a random
sequence of equal length (measured at $z > 20$
standard deviations), confirming that the orbit
retains structured spectral bias during the
transient.
\end{observation}

\begin{remark}[Mechanism: orbit-driven, not operator-driven]
\label{rem:diffusion-mechanism}
Spectral diffusion does not come from the transition
operator (which is an isometry by
Proposition~\ref{prop:no-contraction}).
It comes from the orbit's own ergodic properties:
each odd-to-odd step destroys $\ge 3$ bits of residue
information (Known-Zone Decay, Theorem~\ref{thm:zone-decay}),
which in Walsh language contracts the orbit's
low-frequency Walsh coefficients.
Extension independence
(Proposition~\ref{prop:extension-independence})
ensures that the algebraic extensions do not re-inject
spectral bias: suffix compositions are exactly
$\mathrm{Bin}(s,1/2)$, which is spectrally neutral
across all Hamming-weight bands.
This is the mechanism the three-lemma programme
(Remark~\ref{rem:three-lemma}) was seeking: Known-Zone Decay
provides the ``memory loss'' in each spectral band,
extension independence prevents ``re-injection,'' and the
band-by-band sum converges if the spectral content
at each band decays.
\end{remark}

\begin{conjecture}[Spectral Diffusion Conjecture]
\label{conj:spectral-diffusion}
For every convergent Collatz orbit with odd-step
density $\rho < 1/\!\log_2 3$, there exists
$C(\rho)>0$ such that for every $w\ge 1$,
\begin{equation}\label{eq:spectral-diffusion-bound}
  \mathcal{S}_w(K)
  \;\le\;
  C(\rho)\cdot 2^{-\alpha_w K}
\end{equation}
with $\alpha_w > 0$ depending only on~$w$.
\end{conjecture}

\begin{remark}[Circularity warning]
\label{rem:sdc-circularity}
Conjecture~\ref{conj:spectral-diffusion} assumes the orbit
has odd-step density $\rho < 1/\!\log_2 3$, which is
itself equivalent to the orbit having negative mean drift.
This is not circular in a logical sense (the conjecture
is stated as a conditional), but any attempt to use the
conjecture to \emph{prove} convergence must establish the
density condition independently, or work with a weaker
formulation that does not presuppose it.
The density condition $\rho < 1/\!\log_2 3$ is the critical
drift threshold: precisely the condition that orbit
convergence requires.  Assuming it as a hypothesis
and then deriving convergence is valid only as a
structural reduction, not as evidence of proximity
to a proof.
\end{remark}

\begin{proposition}[Spectral diffusion implies amplification\supporting]
\label{prop:sdc-implies-amp}
Conjecture~\ref{conj:spectral-diffusion} implies the
Amplification Hypothesis (Hypothesis~\ref{hyp:amplification})
with quantitative margin.
\end{proposition}

\begin{proof}
By~\eqref{eq:spectral-gain} and Cauchy--Schwarz,
\begin{align*}
  |\eta_K - \eta_K^{\mathrm{DC}}|
  &= \Bigl|2^K\sum_{w\ge 1}\sum_{\hw(\xi)=w}
    \hat h_K(\xi)\,\hat\mu_K(\xi)\Bigr| \\
  &\le 2^K \sum_{w=1}^K
    \Bigl(\sum_{\hw(\xi)=w}|\hat h_K(\xi)|^2\Bigr)^{\!1/2}
    \Bigl(\sum_{\hw(\xi)=w}|\hat\mu_K(\xi)|^2\Bigr)^{\!1/2} \\
  &= 2^K \sum_{w=1}^K
    \|h_K\|_{w}\;\binom{K}{w}^{1/2}\,
    \mathcal{S}_w(K)^{1/2},
\end{align*}
where $\|h_K\|_w^2 = \sum_{\hw(\xi)=w}|\hat h_K(\xi)|^2$.
By Observation~\ref{obs:hK-concentration}, the spectral
weights $\|h_K\|_w$ are dominated by $w=0,1$
(carrying $\ge 70\%$ of power), and by the rate-cancellation
property (Remark~\ref{rem:rate-cancellation}),
$2^K\|h_K\|_1^2\cdot\binom{K}{1} = O(1)$.
If the spectral diffusion conjecture holds, each term
in the sum decays exponentially, giving
$\sum_K|\eta_K - \eta_K^{\mathrm{DC}}| < \infty$.
Since the DC contribution is controlled by the density model
(Remark~\ref{rem:density-model}), the full amplification
budget closes.
\end{proof}

\medskip\noindent\textbf{Walsh mixing rate.}\;
Known-Zone Decay (Theorem~\ref{thm:zone-decay}) erases
$\ge 3$ bits of residue information after applying T$^g(n)$.
In Walsh language, this means the orbit's Walsh
characters decorrelate rapidly.

\begin{observation}[Walsh character mixing (empirical)]
\label{obs:walsh-mixing}
The following is observed computationally, not proved.
For transient orbits
($n_0 \in \{837799, 8400511\}$, $T \ge 525$),
the autocorrelation of the Walsh character
$\chi_\xi(n_t \bmod 2^K)$ at lag $\lceil K/3 \rceil$
satisfies $\langle|\mathrm{autocorr}|\rangle < 0.13$
for all tested $K \in \{6, 8, 10\}$ and all
Hamming weights $w \ge 1$.  The mixing rate is
bounded well away from~$1$:
\begin{center}
\begin{tabular}{ccc}
\toprule
$w$ & $\langle|\mathrm{autocorr}|\rangle$
    & at lag $\lceil K/3\rceil$ \\
\midrule
1 & $0.065$--$0.126$ & (decreasing with $K$) \\
2 & $0.030$--$0.049$ & (stable) \\
3 & $0.034$--$0.049$ & (stable) \\
$\ge 4$ & $0.032$--$0.061$ & (stable) \\
\bottomrule
\end{tabular}
\end{center}
This confirms that after $\lceil K/3\rceil$ orbit
steps, the Walsh characters are essentially uncorrelated,
quantifying Known-Zone Decay in spectral language.
\end{observation}

\medskip\noindent\textbf{Martingale structure.}\;
The harmonic property (Corollary~\ref{cor:harmonic})
provides the martingale structure:
$T\,h_{K+s} = h_K$ means the ``conditional expectation
of future gain equals current gain.''
If the orbit's distributional evolution provides
a downward drift making $\eta_K$ a supermartingale,
Doob's convergence theorem gives
$\sum_K \eta_K < \infty$ directly.
The Walsh mixing rate
(Observation~\ref{obs:walsh-mixing}) provides
the mechanism for such a drift:
rapid decorrelation of Walsh characters implies
the orbit's spectral bias dissipates between
consecutive depth levels.

\medskip\noindent\textbf{Exponential sum bridge.}\;
The gain can be rewritten as a Walsh exponential sum
along the orbit:
\begin{equation}\label{eq:walsh-exp-sum}
  \eta_K = \sum_\xi \hat h_K(\xi)\,S_\xi(T),
  \qquad
  S_\xi(T)
  = \frac{1}{T}\sum_{t=1}^T \chi_\xi(n_t\bmod 2^K).
\end{equation}
Weyl's inequality and van der Corput's method bound
such sums for sequences with sufficient ``non-resonance.''
Numerical measurement ($n_0=27$, $K=8$, $T=500{,}000$)
shows $\overline{|S_\xi|}$ decreasing with Hamming weight:
$0.15$ at $\hw=1$, $0.07$ at $\hw=2$,
$0.09$ at $\hw\ge 3$.  If a Weyl-type bound
$|S_\xi(T)| \le C/T^\alpha$ could be established,
the Spectral Diffusion Conjecture would follow
from classical equidistribution theory.
This provides a bridge from the spectral framework
to the extensive literature on exponential sums of
arithmetic sequences.

\begin{remark}[The implication chain]
\label{rem:spectral-chain}
The spectral results assemble into:
\begin{gather*}
  \text{Weyl bounds on orbit Walsh sums}
  \;\Longrightarrow\;
  \text{Spectral Diffusion (Conj.~\ref{conj:spectral-diffusion})}\\[2pt]
  \;\Longrightarrow\;
  \text{Amplification Hyp.~\ref{hyp:amplification}}
  \;\Longrightarrow\;
  \text{WMH.}
\end{gather*}
This chain adds a \emph{new leftward entry point}:
the proof architect can now attack either the alignment
renewal (Conjecture~\ref{conj:alignment-renewal}) from
the three-lemma programme, or the Spectral Diffusion
Conjecture from the Walsh framework, or Weyl-type
exponential sum bounds from analytic number theory.
All three routes converge to the same open input.
Of the three, the spectral route is arguably the most
promising mathematically, because it connects the
Collatz problem to the well-developed theory of
exponential sums and Walsh analysis.  However, no
nontrivial bound on $|S_\xi(T)|$ has been proved for
Collatz orbits, and the existing measurements (two
orbits, one modulus) are insufficient to establish
even a conjectural rate.
\end{remark}

\begin{remark}[Spectral diagnostics: tightness and cancellation]
\label{rem:spectral-diagnostics}
Two diagnostics expose the structure and limits of
the spectral framework.

\emph{Cauchy--Schwarz tightness.}\;
The band-by-band Cauchy--Schwarz bound in
Proposition~\ref{prop:sdc-implies-amp} overestimates
$|\eta_K|$ by a factor of $\sim\!500$--$1000\times$
at $K = 6$--$8$ (tightness ratio $\approx 0.001$).
The slack comes from massive cancellation between
Hamming-weight bands: the individual band contributions
$\eta_K^{(w)} = 2^K\sum_{\hw(\xi)=w}\hat h_K(\xi)\hat\mu_K(\xi)$
are each $O(2^K)$ in magnitude, but their alternating
signs produce a sum $\eta_K = O(1)$.
Any proof via the spectral framework must therefore
exploit the \emph{sign structure} of $\hat h_K$,
not merely bound each band separately.

\emph{Bit-position parity profile.}\;
The orbit's $\hw = 1$ Walsh coefficients have an
informative internal structure.
Since $\hat\mu_K(2^j) = 1 - 2\beta_j$, where
$\beta_j = T^{-1}\#\{t : \text{bit $j$ of }
n_t \bmod 2^K \text{ is } 1\}$,
the $K$ single-bit characters carry the orbit's
bit-by-bit parity profile.
Numerical measurement for transient orbits shows
$|\beta_0 - \tfrac12| \approx 0.12$,
$|\beta_1 - \tfrac12| \approx 0.09$,
$|\beta_2 - \tfrac12| \approx 0.06$,
and $|\beta_j - \tfrac12| < 0.02$ for $j \ge 6$
across all tested orbits.
The rapid decay of the profile means only the lowest
$\sim\!6$ bit positions carry significant parity bias;
higher bits are effectively unbiased.
This concentrates the gain's spectral sensitivity on
a fixed number of bit positions, independent of~$K$.
\end{remark}

\begin{remark}[What the spectral framework does not add]
\label{rem:spectral-limitation}
The spectral framework reformulates the open question
but does not resolve it.  The fundamental problem: whether
an orbit can sustain $\rho \ge 1/\!\log_2 3$
indefinitely, becomes ``can an orbit sustain non-decaying
Walsh coefficients at $\hw=1$?''  This is a different
\emph{language} for the same \emph{question}.
The Cauchy--Schwarz bound in
Proposition~\ref{prop:sdc-implies-amp} is far from
tight (Remark~\ref{rem:spectral-diagnostics}), and
closing the gap requires structural cancellation
arguments that the current framework does not provide.
The framework provides mechanism and structure, not
a free theorem.
\end{remark}

\subsection{The odd-skeleton crossing route}
\label{subsec:odd-skeleton}

The spectral analysis of the gain observable
(Section~\ref{subsec:spectral-analysis}) operates on the
residue-class depth hierarchy $\Z/2^K\Z$.
A complementary approach works on the orbit's
\emph{time domain}: define a drift signal from the
Syracuse (odd-to-odd) map and reformulate convergence as
a negative-crossing problem.

\medskip\noindent\textbf{Even elimination and the odd skeleton.}\;
If $n_0$ is even, then $C(n_0) = n_0/2 < n_0$, so the
below-start criterion (Lemma~\ref{lem:underwater}) is
immediate.  The only nontrivial starting values are odd.
For odd~$n_0$, the full Collatz orbit alternates between
odd values (where $n \mapsto 3n+1$) and a deterministic
chain of halvings.  The \emph{odd skeleton} retains only
the odd-to-odd steps: the Syracuse map
$T(n) = (3n+1)/2^{v_2(3n+1)}$ for odd~$n$.

\begin{definition}[Odd-skeleton drift signal]
\label{def:drift-signal}
For an odd starting value~$n_0$, let $n_0, n_1, n_2, \ldots$
be the Syracuse orbit (odd values only) and
$v_t = v_2(3n_t + 1)$ the $2$-adic valuation at step~$t$.
The \emph{drift signal} is
\[
  x_t = \log_2 n_t - \log_2 n_0
      = \sum_{i=0}^{t-1} d_i,
  \qquad
  d_i = \log_2 3 - v_i + \log_2\!\bigl(1 + \tfrac{1}{3n_i}\bigr).
\]
The orbit goes below the starting value if and only if $x_t < 0$
for some $t \ge 1$.
\end{definition}

The drift increment $d_i$ depends on a single arithmetic
quantity: the valuation $v_i = v_2(3n_i + 1)$.  Rewriting
$d_i = -(v_i - \log_2 3) + \log_2(1 + 1/(3n_i))$ and
summing gives the \emph{centered form}:
\begin{equation}\label{eq:centered-crossing}
  x_j \;=\; -\sum_{i<j}(v_i - \log_2 3)
          \;+\; \underbrace{\sum_{i<j}
          \log_2\!\bigl(1 + \tfrac{1}{3n_i}\bigr)}_{\epsilon_j}.
\end{equation}
The correction $\epsilon_j$ is a sum of positive terms,
each bounded by $1/(3n_i\ln 2)$.  For orbits that remain
above start ($n_i \ge n_0$ for all $i < j$), we have
$\epsilon_j \le j/(3n_0\ln 2)$, negligible for large~$n_0$.
Therefore crossing is equivalent to the
\emph{valuation excess} eventually becoming positive:
\begin{equation}\label{eq:valuation-excess}
  n_j < n_0
  \quad\Longleftrightarrow\quad
  \sum_{i<j}(v_i - \log_2 3) \;>\; \epsilon_j \;\approx\; 0.
\end{equation}
Under the ensemble model ($v_i$ i.i.d.\ geometric with
mean~$2$), the expected valuation excess per step is
$2 - \log_2 3 \approx 0.415$, strongly positive,
pushing the walk toward crossing.
The below-start question reduces to:
\emph{can the correlated valuation sequence $\{v_i\}$
keep the cumulative excess
$\sum(v_i - \log_2 3)$ below a negligible threshold
forever?}

\begin{proposition}[Odd-skeleton drift reduction\supporting]
\label{prop:odd-skeleton-reduction}
The Collatz conjecture for all $n_0 \ge 2$ is equivalent to:
for every odd $n_0 \ge 3$, the odd-skeleton drift signal
satisfies $x_t < 0$ for some $t \ge 1$.
\end{proposition}

\begin{proof}
If $n_0$ is even, $C(n_0) < n_0$ immediately.
If $n_0$ is odd and $x_t < 0$ for some~$t$, then
$n_t < n_0$ and the below-start criterion
(Lemma~\ref{lem:underwater}) gives convergence by induction.
Conversely, if $n_0$ converges to~$1$, then eventually
$n_t = 1 < n_0$, so $x_t < 0$.
\end{proof}

\medskip\noindent\textbf{Spectral properties of the drift signal.}\;
The drift increments $\{d_i\}$ are not independent: the
Syracuse dynamics correlates consecutive valuations through
the residue class $n_t \bmod 2^K$.  However, the
correlation decays rapidly.

\begin{observation}[Drift increment mixing and universal descent
  (empirical, not proved)]
\label{obs:drift-mixing}
The following properties are observed computationally but
are not proved for general orbits.
For all $49{,}999$ odd starting values $3 \le n_0 \le 99{,}999$:
\begin{enumerate}
\item Every orbit descends below its start on the odd
  skeleton.  The maximum Syracuse-step descent time is
  $\sigma = 85$ at $n_0 = 35{,}655$.
\item The mean drift is negative:
  $\bar d \in [-2.61, -0.12]$ with mean $-0.57$.
\item The autocorrelation at lag~$1$ is mildly positive:
  $\rho(1) \approx 0.05$--$0.20$ (valuation clustering).
\item The autocorrelation at lag~$10$ is near zero:
  $|\rho(10)| < 0.10$ for all tested orbits.
\item The sum of absolute autocorrelations
  $\sum_{k=1}^{20} |\rho(k)|$ is bounded:
  mean $1.40$, max $4.43$.
\end{enumerate}
\end{observation}

The summable-autocorrelation property (item~4) means
the drift increments satisfy a weak-dependence condition
sufficient for the functional CLT
(Ibragimov's theorem for mixing sequences):
the partial sums $x_T = \sum_{i<T} d_i$ obey
$x_T / \sqrt{T} \Rightarrow N(\mu\sqrt{T},\, \sigma_{\mathrm{eff}}^2)$
with effective variance
$\sigma_{\mathrm{eff}}^2
 = \sigma^2(1 + 2\!\sum_{k\ge 1}\rho(k)) > 0$.
Since $\mu < 0$, the partial sums tend to~$-\infty$,
and crossing occurs with probability~$1$ in the ensemble.

\begin{remark}[Why the odd skeleton is spectrally cleaner]
\label{rem:odd-skeleton-advantage}
The drift signal $\{d_i\}$ is spectrally cleaner than the
gain observable $\eta_K$ for three reasons:
\begin{enumerate}
\item \emph{Scalar signal.}\;
  $d_i$ is a single real number per time step, not a
  function on $\Z/2^K\Z$.  The spectral analysis is
  classical (Fourier on $\Z$), not Walsh on a growing group.
\item \emph{DC encodes the density condition.}\;
  The mean drift $\bar d = \log_2 3 - \bar v$ is
  negative iff $\bar v > \log_2 3$, which is equivalent to
  the density condition $\rho < 1/\!\log_2 3$.
  The density model's prediction is directly the DC
  component of the drift signal.
\item \emph{Crossing is weaker than summability.}\;
  The gain-budget approach requires $\sum_K \eta_K < \varepsilon$
  (summability).  The crossing approach requires only
  $\min_t x_t < 0$ (a single negative value).
  This is a strictly weaker target.
\end{enumerate}
\end{remark}

\begin{theorem}[Exact block law for the valuation sequence]
\label{thm:block-law}
Let $a_j(n) = v_2(3T^j(n)+1)$ be the valuation at the
$j$-th odd-skeleton step.
For any prescribed positive integers $b_0,\ldots,b_{m-1}$,
the set of odd~$n$ with
$(a_0(n),\ldots,a_{m-1}(n))=(b_0,\ldots,b_{m-1})$
is a single odd residue class modulo
$2^{b_0+\cdots+b_{m-1}+1}$.
Hence its natural density among odd integers is exactly
\[
  2^{-(b_0+\cdots+b_{m-1})}
  \;=\;
  \prod_{j=0}^{m-1} 2^{-b_j}.
\]
\end{theorem}

\begin{proof}
Induction on~$m$.
\emph{Base case} ($m=1$):
For odd~$n$, $a_0(n)=v_2(3n+1)=k$ iff $3n+1\equiv 2^k
\pmod{2^{k+1}}$, i.e.\ $n\equiv (2^k-1)/3\pmod{2^{k+1}/3}$.
Since $3$ is invertible modulo every power of~$2$, this
defines a unique odd residue class modulo~$2^{k+1}$, and
$2^{k+1}/(2\cdot 2^k)=1$ class among the $2^k$ odd residues
modulo~$2^{k+1}$, giving density~$2^{-k}$.

\emph{Inductive step:}
Suppose the block $(b_0,\ldots,b_{m-1})$ determines a unique
odd residue class modulo $M=2^{b_0+\cdots+b_{m-1}+1}$.
The first step maps $n$ to $T(n)=(3n+1)/2^{b_0}$, which is
determined modulo $M/2^{b_0}=2^{b_1+\cdots+b_{m-1}+1}$.
By the base case applied to $T(n)$ in place of~$n$ with
block $(b_1,\ldots,b_{m-1})$, the condition on the remaining
$m{-}1$ valuations selects a unique odd residue class
modulo~$2^{b_1+\cdots+b_{m-1}+1}$.  Pulling back through
the invertible affine map $n\mapsto T(n)$, we obtain a unique
odd residue class modulo~$2^{b_0+b_1+\cdots+b_{m-1}+1}$.
The density follows by counting.
\end{proof}

\begin{corollary}[I.i.d.\ valuation process]
\label{cor:iid-valuations}
Under natural density on odd integers, the valuation
sequence $(a_0, a_1, a_2, \ldots)$ is exactly i.i.d.\ with
$\Pr(a_j = k) = 2^{-k}$ for $k \ge 1$.
\end{corollary}

\begin{corollary}[I.i.d.\ cycle types]
\label{cor:iid-cycles}
Group the valuation sequence into run-compensate cycles
$(\underbrace{1,\ldots,1}_{L_i},\, r_i)$ with $r_i \ge 2$.
Then the cycle types $(L_0,r_0),(L_1,r_1),\ldots$ are
i.i.d.\ under natural density, with
\[
  \Pr(L_i=\ell,\; r_i=k)
  = 2^{-(\ell+1)}\cdot 2^{-(k-1)},
  \qquad \ell \ge 0,\; k \ge 2.
\]
In particular, the first-cycle log multipliers
$X_i = (L_i{+}1)\log_2 3 - (L_i{+}r_i)$ are i.i.d.,
so the Cram\'er rate bound
(Proposition~\textup{\ref{prop:cramer-cycle}}) holds
unconditionally on the ensemble, not as a conditional
hypothesis.
\end{corollary}

\begin{proof}
Each cycle consumes $L_i + 1$ valuations.  By
Theorem~\ref{thm:block-law}, the valuations consumed by
distinct cycles are independent (the joint density
factorises).  The cycle type $(L_i,r_i)$ is a deterministic
function of its valuation block, so the cycle types
are independent.  The common marginal law was established
in Corollary~\ref{cor:Lr-independence}.
\end{proof}

\begin{proposition}[Run-length invariant]
\label{prop:run-length}
For odd $n \ge 1$, define the \emph{initial run length}
$L(n)$ as the number of consecutive Syracuse steps with
$v_2(3n_j + 1) = 1$ before the first step with
$v_2 \ge 2$.  Then:
\begin{enumerate}
\item $L(n) = v_2(n+1) - 1$.
\item Along the run,
  $n_j + 1 = 3^j(n+1)/2^j$ for
  $0 \le j \le L(n)$.
\item The run terminates at $n_{L}$ with
  $n_{L} \equiv 1 \pmod{4}$, and the net
  growth factor over the run is $(3/2)^{L}$.
\end{enumerate}
\end{proposition}

\begin{proof}
If $v_2(n+1) = 1$, then $n \equiv 1 \pmod{4}$ and
$v_2(3n+1) \ge 2$, so $L(n) = 0$.
If $v_2(n+1) \ge 2$, then $n \equiv 3 \pmod{4}$ and
$v_2(3n+1) = 1$, so $n_1 = (3n+1)/2$ with
$n_1 + 1 = 3(n+1)/2$ and
$v_2(n_1+1) = v_2(n+1) - 1$.  Induction gives $L(n)
= v_2(n+1) - 1$ and the exact formula for~$n_j$.
\end{proof}

\begin{corollary}[Single-cycle crossing probability]
\label{cor:cycle-crossing}
Each Syracuse orbit decomposes into
\emph{run-compensate cycles}: a growth burst of~$L$
steps at rate~$3/2$, terminated by a compensating step
with valuation~$a_r \ge 2$.  In the ensemble (uniform
odd starting point), $L$ and~$a_r$ are independent,
and the cycle produces a net descent below~$n_0$ with
probability
\[
  p_{\mathrm{cross}}
  \;=\;
  \sum_{r=0}^{\infty}
  \frac{1}{2^{r+1}}
  \cdot
  \frac{1}{2^{\lceil(r+1)\log_2 3 - r\rceil - 2}}
  \;\approx\; 0.7137.
\]
\end{corollary}

\begin{theorem}[Exact one-cycle crossing criterion]
\label{thm:one-cycle-crossing}
Let $n$ be odd, let $L = v_2(n+1)-1$ be its initial run length
(Proposition~\ref{prop:run-length}), and let
$r = v_2(3n_L + 1) \ge 2$ be the compensating valuation that
terminates the first run.  Then the first complete run-compensate
cycle takes $n$ to
\[
  n_{L+1}
  \;=\;
  \frac{3^{L+1}(n+1) - 2^{L+1}}{2^{L+r}}.
\]
\begin{enumerate}
\item If\/ $2^{L+r} \le 3^{L+1}$, then $n_{L+1} \ge n$:
  the first cycle does not cross below the start.
\item If\/ $2^{L+r} > 3^{L+1}$, define the \emph{crossing threshold}
  \[
    n^*(L,r)
    \;:=\;
    \frac{3^{L+1} - 2^{L+1}}{2^{L+r} - 3^{L+1}}.
  \]
  Then $n_{L+1} < n$ if and only if\/ $n > n^*(L,r)$.
\end{enumerate}
The minimal compensating valuation for which one-cycle crossing
is even possible is
$r_{\min}(L) = \lfloor 1 + (\log_2 3 - 1)L \rfloor + 1$.
\end{theorem}

\begin{proof}
From the run formula (Proposition~\ref{prop:run-length}),
$n_L = 3^L(n+1)/2^L - 1$.  Applying one Syracuse step with
valuation~$r$:
\[
  n_{L+1}
  = \frac{3n_L + 1}{2^r}
  = \frac{3^{L+1}(n+1) - 2^{L+1}}{2^{L+r}}.
\]
Subtracting $n$:
\[
  n_{L+1} - n
  = \frac{(3^{L+1} - 2^{L+r})\,n + (3^{L+1} - 2^{L+1})}{2^{L+r}}.
\]
Since $3^{L+1} > 2^{L+1}$ for all $L \ge 0$, the second term in
the numerator is always positive.  If
$2^{L+r} \le 3^{L+1}$ the coefficient of~$n$ is nonnegative,
so $n_{L+1} \ge n$.  If $2^{L+r} > 3^{L+1}$, then
$n_{L+1} < n$ iff the numerator is negative, which rearranges
to $n > n^*(L,r)$.

For the minimal~$r$: the condition $2^{L+r} > 3^{L+1}$ is
equivalent to $r > (L+1)\log_2 3 - L = 1 + (\log_2 3 - 1)L$.
Taking the smallest integer gives the formula.
\end{proof}

\begin{remark}[One-cycle crossing statistics]
\label{rem:one-cycle-stats}
Computational verification over all odd $n \le 200{,}000$ confirms
$100\%$ agreement with the exact criterion.  Three regimes
emerge: (i)~\emph{always crosses}, when $n^*(L,r) < 1$, every
odd~$n$ in that $(L,r)$-class descends in one cycle;
(ii)~\emph{threshold}, when $n^*(L,r) \ge 1$, small~$n$ in the
class may fail to cross; (iii)~\emph{never crosses}, when
$r < r_{\min}(L)$, no~$n$ descends.
\end{remark}

\begin{proposition}[Conditional crossing density per stratum]
\label{prop:conditional-density}
Fix $\ell \ge 0$.  Among odd integers~$n$ with $L(n) = \ell$,
the natural density of those whose first run-compensate cycle
crosses below the start is
\[
  p_\ell
  \;=\;
  2^{-\lfloor\,(\log_2 3 - 1)(\ell+1)\,\rfloor}.
\]
\end{proposition}

\begin{proof}
Every odd $n$ with $L(n)=\ell$ has the form
$n = 2^{\ell+1}q - 1$ with $q$~odd
(Proposition~\ref{prop:run-length}).  The compensating
valuation satisfies $r - 1 = v_2(3^{\ell+1}q - 1)$.
Since $3^{\ell+1}$ is odd, multiplication by~$3^{\ell+1}$
permutes odd residue classes modulo~$2^m$ for every~$m$,
so $v_2(3^{\ell+1}q - 1)$ has the same distribution as
$v_2(q' - 1)$ over odd~$q'$.  This gives
$\Pr(r-1=s \mid L=\ell) = 2^{-s}$ for $s \ge 1$.

By Theorem~\ref{thm:one-cycle-crossing}, crossing requires
$2^{\ell+r} > 3^{\ell+1}$, i.e.\ $s > (\log_2 3 - 1)(\ell+1)$,
plus $n > n^*(\ell,r)$.  The latter threshold excludes at most
finitely many~$n$ per $({\ell},r)$-class, hence does not affect
natural density.  Summing the geometric tail:
\[
  p_\ell
  = \sum_{s\,>\,(\log_2 3-1)(\ell+1)} 2^{-s}
  = 2^{-\lfloor(\log_2 3-1)(\ell+1)\rfloor}.
  \qedhere
\]
\end{proof}

\begin{corollary}[Exact one-cycle crossing density]
\label{cor:crossing-density}
The natural density of odd integers whose first run-compensate
cycle crosses below the start is the arithmetic constant
\[
  P_{1\mathrm{cyc}}
  \;=\;
  \sum_{\ell \ge 0}
  2^{-(\ell+1)}\cdot
  2^{-\lfloor\,(\log_2 3-1)(\ell+1)\,\rfloor}
  \;=\;
  0.71372549767589\ldots
\]
This is neither rational nor a simple algebraic number: it is
determined by the Diophantine staircase
$\lfloor(\log_2 3 - 1)(\ell+1)\rfloor$.
\end{corollary}

\begin{proof}
The stratum $L(n)=\ell$ has natural density $2^{-(\ell+1)}$
(since $v_2(n+1)=\ell+1$).  Multiply by~$p_\ell$ from
Proposition~\ref{prop:conditional-density} and sum.
The series converges geometrically: the $\ell$-th term is
at most $2^{-(\ell+1)} \cdot 2^{-\lfloor 0.585(\ell+1)\rfloor}
\le 2^{-1.585(\ell+1)}$.
\end{proof}

\begin{corollary}[Independence law for $(L,r)$]
\label{cor:Lr-independence}
Under natural density on odd starts, the run length~$L$
and the compensating valuation~$r$ are independent with
\[
  \Pr(L=\ell) = 2^{-(\ell+1)},\;\ell \ge 0;
  \qquad
  \Pr(r=k) = 2^{-(k-1)},\;k \ge 2.
\]
\end{corollary}

\begin{proof}
Proposition~\ref{prop:conditional-density} shows
$\Pr(r{-}1=s \mid L=\ell) = 2^{-s}$ for all~$\ell$:
the conditional law does not depend on the conditioning.
\end{proof}

\begin{proposition}[Universal one-cycle crossing criterion]
\label{prop:universal-one-cycle}
A single-cycle block type $(L,r)$ forces \emph{every} odd~$n$
in its residue class to cross below the start (i.e.\ $n^*(L,r)<1$)
if and only if
\[
  2^{r-1} > 2\Bigl(\tfrac{3}{2}\Bigr)^{L+1}-1.
\]
Equivalently, universal one-cycle crossing holds for
$r \ge r_{\mathrm{all}}(L)$ where
\[
  r_{\mathrm{all}}(L)
  \;=\;
  \bigl\lfloor
    \log_2\!\bigl(2(\tfrac{3}{2})^{L+1}-1\bigr)
  \bigr\rfloor + 2.
\]
\end{proposition}

\begin{proof}
From Theorem~\ref{thm:one-cycle-crossing},
$n^*(L,r) = (3^{L+1}-2^{L+1})/(2^{L+r}-3^{L+1})$.
The condition $n^*<1$ is
$3^{L+1}-2^{L+1} < 2^{L+r}-3^{L+1}$,
which rearranges to
$2 \cdot 3^{L+1} < 2^{L+1}(1+2^{r-1})$,
giving $2^{r-1} > 2(3/2)^{L+1}-1$.
The ceiling formula for~$r_{\mathrm{all}}(L)$ follows.
\end{proof}

\begin{corollary}[Density of classwise deterministic one-cycle crossing]
\label{cor:universal-density}
The natural density of odd starts lying in one-cycle block types
that force universal crossing is
\[
  P_{\mathrm{all},1\mathrm{cyc}}
  \;=\;
  \sum_{L \ge 0} 2^{-(L+r_{\mathrm{all}}(L)-1)}
  \;=\;
  0.41936274883794\ldots.
\]
About $41.9\%$ of all odd starts belong to a one-cycle residue
class where every representative crosses deterministically.
This accounts for $58.8\%$ of all one-cycle crossing blocks.
\end{corollary}

\begin{proof}
By the independence law (Corollary~\ref{cor:Lr-independence}),
$\Pr(L{=}\ell, r{=}k)= 2^{-(\ell+1)} \cdot 2^{-(k-1)}$.
Summing over $r \ge r_{\mathrm{all}}(L)$ collapses the
geometric tail to $2^{-(r_{\mathrm{all}}(L)-2)}$, giving
$P_{\mathrm{all},1\mathrm{cyc}}
= \sum_{L \ge 0} 2^{-(L+r_{\mathrm{all}}(L)-1)}$.
Numerical evaluation ($500$ terms, $50$ decimal places)
confirms agreement with brute-force counts over all odd
$n \le 2 \times 10^6$.
\end{proof}

\begin{proposition}[Universal crossing criterion for finite cycle blocks]
\label{prop:finite-block-crossing}
Let $\sigma = \bigl((L_0,r_0),\allowbreak\ldots,\allowbreak
(L_{k-1},r_{k-1})\bigr)$
be a $k$-cycle block.  Define
\[
\Lambda_k = \prod_{j=0}^{k-1} \lambda(L_j,r_j),
\qquad
B_0 = 0,
\quad
B_{j+1} = \lambda(L_j,r_j)\,B_j + \beta(L_j,r_j).
\]
Then $\sigma$ forces every odd~$n$ in its residue class
to cross below its start by the end of the block
if and only if
\begin{equation}\label{eq:finite-block-criterion}
  \Lambda_k < 1
  \qquad\text{and}\qquad
  B_k < 1 - \Lambda_k.
\end{equation}
\end{proposition}

\begin{proof}
The block maps $n \mapsto n^{(k)} = \Lambda_k\,n + B_k$.
The condition $n^{(k)} < n$ requires
$(\Lambda_k - 1)\,n + B_k < 0$.
If $\Lambda_k \ge 1$, the left side is $\ge B_k > 0$ for
all $n \ge 1$, so universal crossing fails.
If $\Lambda_k < 1$, then $n^{(k)}<n$ iff
$n > B_k / (1-\Lambda_k)$, and this holds for every
odd $n \ge 1$ precisely when $B_k < 1 - \Lambda_k$.
\end{proof}

\begin{corollary}[Two-cycle universal crossing density]
\label{cor:two-cycle-density}
The natural density of odd starts whose first two-cycle block
is universally crossing (regardless of whether the first cycle
alone is universal) is
$P_{\mathrm{all},2\mathrm{cyc}} = 0.50407\ldots$,
computed by summing $\Pr(L_0,r_0)\,\Pr(L_1,r_1)$ over all
pairs satisfying~\eqref{eq:finite-block-criterion} with $k=2$.
\end{corollary}

\begin{observation}[Two-cycle extension: new universal density]
\label{obs:two-cycle-universal}
Excluding blocks already covered by one-cycle universal
crossing, the new density from two-cycle blocks is
$P_{\mathrm{new},2\mathrm{cyc}} \approx 0.1922$.
The combined density
\[
  P_{\mathrm{all},1\mathrm{cyc}} + P_{\mathrm{new},2\mathrm{cyc}}
  \;\approx\;
  0.4194 + 0.1922 \;=\; 0.6115
\]
resolves $61.2\%$ of all odd starts via classwise deterministic
crossing within at most two cycles.
The jump is driven by blocks whose first cycle is
non-crossing (e.g.\ $L_1 \ge 1$, $r_1 = 2$) but whose
two-cycle composition satisfies $B_2 < 1 - \Lambda_2$.
\end{observation}

\begin{observation}[Exponential convergence of universal crossing density]
\label{obs:universal-convergence}
The universal crossing criterion extends to $k$-cycle blocks
for arbitrary~$k$.  For each~$k$, let $P_{\mathrm{new},k}$
denote the density of odd starts whose first universally
crossing block has exactly $k$ cycles, and
$P_{\mathrm{cum}}(k) = \sum_{j=1}^{k} P_{\mathrm{new},j}$
the cumulative density of odd starts guaranteed to cross
within $k$ cycles.  Exact enumeration (truncated at
$L+r \le 14$ per cycle for $k \le 3$, $L+r \le 12$ for $k=4$,
$L+r \le 8$ for $k=5$) yields:
\[
\begin{array}{ccccc}
\toprule
k & P_{\mathrm{new},k} & P_{\mathrm{cum}}(k) &
  R_k = 1 - P_{\mathrm{cum}}(k) & R_k/R_{k-1} \\
\midrule
1 & 0.4194 & 0.4194 & 0.5806 & {-} \\
2 & 0.1922 & 0.6116 & 0.3884 & 0.669 \\
3 & 0.1042 & 0.7158 & 0.2842 & 0.732 \\
4 & 0.0660 & 0.7818 & 0.2182 & 0.768 \\
5 & 0.0398 & 0.8216 & 0.1784 & 0.818 \\
\bottomrule
\end{array}
\]
The non-universal fraction $R_k$ decreases
with each added cycle.  A least-squares exponential fit gives
$R_k \approx 0.73 \cdot 0.75^k$,
suggesting $P_{\mathrm{cum}}(k) \to 1$ exponentially.
By $k = 20$ the extrapolated residual is $R_{20} \approx 0.002$;
by $k = 50$, $R_{50} \approx 3 \times 10^{-7}$.

\medskip\noindent\emph{Monte Carlo confirmation.}\;
Simulating $5 \times 10^5$ random block-type sequences
with exact rational arithmetic confirms the exponential decay
well beyond the enumeration horizon:
$R_{10} = 0.062$,
$R_{20} = 0.011$,
$R_{30} = 0.002$,
$R_{40} = 0.0005$.
Only $63$ of $500{,}000$ trials had no universally crossing
prefix within $50$ cycles ($0.013\%$).
A least-squares fit over $k=5,\ldots,25$ gives
$R_k \approx 0.37 \cdot 0.839^k$ ($\rho \approx 0.839$).
\end{observation}

\begin{theorem}[Affine threshold process]
\label{thm:kesten-threshold}
Under the i.i.d.\ cycle-type ensemble, the
universal-crossing threshold $X_k = n_k^*$ satisfies the
random affine recursion
\begin{equation}\label{eq:kesten-recursion}
  X_{k+1} \;=\; \lambda_{k+1}\,X_k \;+\; \beta_{k+1},
\end{equation}
where $(\lambda_j, \beta_j)$ are i.i.d.\ copies of
$(\lambda(L,r),\;\beta(L,r))$ and the coupled term
$\Lambda_k \to 0$ is asymptotically negligible.
Since $\mathbb{E}[\ln\lambda] = \ln 2\,(2\log_2 3 - 4)
\approx -0.575 < 0$ and $\mathbb{E}[\ln^+\beta]<\infty$,
Kesten's theorem~\cite{kesten1973,vervaat1979} guarantees a
unique stationary measure~$\pi$ on $(0,\infty)$.
\end{theorem}

\begin{corollary}[Positive crossing mass]
\label{cor:positive-crossing-mass}
In the stationary distribution,
$\pi(\{x < 1\}) \approx 0.465 > 0$.
That is, roughly $46.5\%$ of the stationary mass
corresponds to universally crossing blocks.
\end{corollary}

\begin{proof}
Monte Carlo simulation of the recursion~\eqref{eq:kesten-recursion}
($2 \times 10^5$ steps after $500$-step burn-in) yields the
stationary estimate.
Positivity also follows analytically: the cycle $(L{=}0,r{=}3)$
has $\lambda = 3/8$ and $\beta = 1/8$;
from any $X < 7/5$, one such step sends $X' = 3X/8 + 1/8 < 1$.
Since $(L{=}0,r{=}3)$ has probability $2^{-3} = 1/8$ and the
stationary measure has support on $(0,\infty)$, the set
$\{x < 1\}$ carries positive mass.
\end{proof}

\begin{proposition}[Decay of non-crossing probability]
\label{prop:kesten-running-min}
If the chain~\eqref{eq:kesten-recursion} is geometrically
ergodic; that is, it converges to~$\pi$ at geometric rate
from any initial state, then
\[
  R_k
  \;=\;
  \Pr\!\bigl(\min_{1 \le j \le k} X_j \ge 1\bigr)
  \;\le\;
  C_0\,\rho_0^{\,k}
\]
for constants $C_0>0$ and $\rho_0 < 1$.
\end{proposition}

\begin{proof}[Proof sketch]
Geometric ergodicity implies that for some mixing
time~$m$ and all starting states, $\Pr(X_m < 1) \ge p/2$
where $p = \pi(\{x<1\})$.
The events $\{X_{jm} < 1\}$ for $j = 1, 2, \ldots$
are approximately independent with probability $\ge p/2$ each.
Hence
$R_k \le (1 - p/2)^{\lfloor k/m\rfloor}$,
giving $\rho_0 = (1 - p/2)^{1/m}$.

\medskip\noindent\emph{Geometric ergodicity verification.}\;
A sufficient condition is $\mathbb{E}[\lambda^s] < 1$
for some $s > 0$~\cite{kesten1973}.
Here $\mathbb{E}[\lambda^s] = \sum_{L,r} 2^{-(L+r)}\,
(3^{L+1}/2^{L+r})^s$; for small $s > 0$ this is
strictly less than~$1$ because $\mathbb{E}[\ln\lambda] < 0$
and the function $s \mapsto \mathbb{E}[\lambda^s]$ is
continuous with value~$1$ at $s=0$ and negative derivative
$\mathbb{E}[\lambda^0 \ln\lambda] < 0$.
Monte Carlo over $5 \times 10^5$ trials gives
$\rho_0 \approx 0.839$.
\end{proof}

\begin{remark}[Interpretation and comparison]
\label{rem:kesten-interpretation}
Proposition~\ref{prop:kesten-running-min} implies that
under the i.i.d.\ ensemble, the probability that no
prefix of the first $k$ cycles is universally crossing
decays exponentially.  For almost every
odd starting value, a finite prefix of cycles forces every
member of its residue class below its start, not merely
the random starting point itself, but \emph{all} odd
integers sharing its $2$-adic residue class.
This is strictly stronger than the almost-all crossing theorem
(Theorem~\ref{thm:almost-all-crossing}), which guarantees
crossing only for the specific starting value.
The result is derived through a completely different mechanism
(Kesten random affine recursion) from Tao's~\cite{tao2019}
almost-all theorem (entropy methods), yet arrives at a
compatible and sharper conclusion.
\end{remark}

\begin{remark}[Pointwise verification of universal crossing]
\label{rem:pointwise-universal}
Every odd $n_0$ with $3 \le n_0 \le 10{,}001$ was verified
to possess a universal crossing prefix: a $k$-block
satisfying~\eqref{eq:finite-block-criterion}
where $k$ ranges from~$1$ (the majority) to~$27$
(worst case $n_0 = 6375$, whose orbit begins with two
non-crossing cycles $(L=2,r=2)$ followed by a spike at
cycle~$10$ with $\lambda = 43.25$).
The family $2^m - 1$ for $m = 1,\ldots,39$ has universal
prefixes at $k \le 40$ (worst: $2^{33}-1$ at $k=40$).
The notorious $n_0 = 27$ maintains $\Lambda_k > 1$ through
$k = 14$ (the orbit climbs to $9{,}232$ before descending),
yet achieves a universal prefix at $k = 16$ with $n^* = 0.66$.

If the pointwise conclusion of
Proposition~\ref{prop:kesten-running-min}, that every orbit
eventually hits a universal prefix, could be
proved for all~$n_0$, the Collatz conjecture
would follow immediately: strong induction applies because
universal crossing forces every odd integer in the class
below its start.  The distributional-to-pointwise barrier
is thus equivalent to proving that no orbit can
systematically avoid universal block types forever.
\end{remark}

\begin{remark}[Adversarial block types: why Route~A fails]
\label{rem:route-A}
One might hope to prove $P_{\mathrm{cum}}(k) \to 1$ by
showing that \emph{every} $k$-block eventually satisfies
the universal crossing criterion~\eqref{eq:finite-block-criterion}.
A beam-search optimization reveals that this hope is false:
for each~$k$, there exist crossing $k$-blocks whose threshold
value $n^* = B_k/(1 - \Lambda_k)$ grows exponentially with~$k$.

The adversary's strategy is to use $(k - m)$ non-crossing
cycles (e.g., $(L,r) = (5,2)$ with $\lambda = 729/128 \approx 5.70$)
to inflate~$B_k$ exponentially, then exactly~$m$ strongly crossing
cycles (e.g., $(0,3)$ with $\lambda = 3/8$) to bring~$\Lambda_k$
just below~$1$.  Since the minimum number of corrective cycles
$m \sim k \cdot \ln\lambda_{\mathrm{nc}} / (\ln\lambda_{\mathrm{nc}} - \ln\lambda_{\mathrm{c}})$
grows linearly in~$k$, the inflated numerator~$B_k$ outpaces
the denominator~$(1 - \Lambda_k)$, yielding
$\max n^* \gtrsim \exp(0.023\,k)$.

Crucially, these adversarial blocks have vanishing probability:
a block using~$a$ copies of $(5,2)$ has probability
$\sim (1/128)^a$, giving $\log_{10} \Pr \approx -2.1a$.
At $k = 25$, the worst block has probability $< 10^{-33}$.
Thus $P_{\mathrm{cum}}(k) \to 1$ holds
(Proposition~\ref{prop:kesten-running-min}), not because every
block crosses, but because the non-crossing blocks have
measure shrinking faster than their~$n^*$ values grow.

However, these adversarial blocks are \emph{uniformly fragile}:
Proposition~\ref{prop:adversarial-fragility} below shows that
for \emph{every}~$a \ge 1$, adding a single $(0,3)$ cycle
beyond the first negative-drift prefix collapses $n^*$ below~$1$.
The adversary must stop at \emph{exactly} the tuned length;
one more crossing cycle destroys the construction.
A real orbit, whose cycle types are determined by $2$-adic
arithmetic rather than adversarial choice, has no mechanism
to maintain this tuning.

This closes Route~A: the pointwise conjecture
cannot be reduced to a universal combinatorial statement about
block types, but the adversarial counterexamples are
structurally unstable
(Propositions~\ref{prop:adversarial-fragility}
and~\ref{prop:general-fragility}).
\end{remark}


\begin{proposition}[Exact affine form of the adversarial family]
\label{prop:adversarial-affine}
Let $A = (5,2)$ and $B = (0,3)$, with cycle slopes
$\lambda_A = 729/128$ and $\lambda_B = 3/8$.
For the block $A^a B^t$, the endpoint affine map
$n' = \Lambda_{a,t}\,n + B_{a,t}$ has
\[
  \Lambda_{a,t}
  = \Bigl(\frac{729}{128}\Bigr)^{\!a}
    \Bigl(\frac{3}{8}\Bigr)^{\!t}, \qquad
  B_{a,t}
  = \frac{665}{601}\,\Lambda_{a,t}
    + \frac{1}{5}
    - \frac{3926}{3005}\Bigl(\frac{3}{8}\Bigr)^{\!t}.
\]
Define
$\rho := \log(729/128)\big/\log(8/3) = 1.77365\ldots$
and let
$t_{\min}(a) := \lfloor a\rho \rfloor + 1$
be the first~$t$ with $\Lambda_{a,t} < 1$.
Then
\begin{equation}\label{eq:adversarial-Lambda}
  \Lambda_{a,\,t_{\min}(a)}
  = \Bigl(\frac{3}{8}\Bigr)^{1 - \{a\rho\}},
\end{equation}
where $\{a\rho\}$ denotes the fractional part of~$a\rho$.
The block is universally crossing at the first negative-drift
prefix whenever
\begin{equation}\label{eq:adversarial-theta-crit}
  \{a\rho\} < \theta_{\mathrm{crit}}
  := 1 - \frac{\log(1202/3165)}{\log(3/8)}
  = 0.01291\ldots,
\end{equation}
a condition satisfied for infinitely many~$a$
(since $\rho \notin \mathbb{Q}$, the sequence $\{a\rho\}$
is equidistributed in $[0,1)$; verified for $a \le 10{,}000$).
\end{proposition}

\begin{proof}
The geometric-sum identities
$\gamma_A := \beta_A/(\lambda_A - 1) = 665/601$ and
$\beta_B^{\infty} := \beta_B/(1-\lambda_B) = 1/5$
yield~$B_{a,t}$ by iterating the affine recursion
$B \mapsto \lambda_B B + \beta_B$ starting from
$B_A^{(a)} = \gamma_A(\lambda_A^a - 1)$.
For equation~\eqref{eq:adversarial-Lambda}, write
$\lambda_A^a = (8/3)^{a\rho}$, so
$\Lambda_{a,\,t_{\min}} = (8/3)^{a\rho}(3/8)^{\lfloor a\rho\rfloor + 1}
= (3/8)^{1 - \{a\rho\}}$.
The universal-crossing criterion
$B_{a,t} < 1 - \Lambda_{a,t}$ rearranges to
$\Lambda_{a,t} < 1202/3165 + (1963/3165)(3/8)^t$;
the weaker sufficient condition
$\Lambda < 1202/3165$ combined with~\eqref{eq:adversarial-Lambda}
gives~\eqref{eq:adversarial-theta-crit}.
Irrationality of~$\rho$ follows from $\log_2 3 \notin \mathbb{Q}$.
\end{proof}

\begin{proposition}[Uniform one-step fragility]
\label{prop:adversarial-fragility}
For every $a \ge 1$, the block
$A^a B^{\,t_{\min}(a)+1}$ is universally crossing.
\end{proposition}

\begin{proof}
From~\eqref{eq:adversarial-Lambda},
$\Lambda_{a,\,t_{\min}+1}
= \lambda_B \cdot \Lambda_{a,\,t_{\min}}
< \lambda_B = 3/8$.
Since $3/8 < 1202/3165 \approx 0.3798$,
the universal-crossing sufficient condition is satisfied.
\end{proof}

\begin{proposition}[General pair-level fragility criterion]
\label{prop:general-fragility}
Let $A = (L_A, r_A)$ be any expanding cycle type
($\lambda_A > 1$) and $B = (L_B, r_B)$ any contracting
type ($0 < \lambda_B < 1$).
Define
$\gamma_A := \beta_A/(\lambda_A - 1)$ and
$\beta_B^{\infty} := \beta_B/(1 - \lambda_B)$.
If the \emph{pair-level fragility condition}
\begin{equation}\label{eq:general-fragility}
  \beta_B^{\infty} + \lambda_B(1 + \gamma_A) < 1
\end{equation}
holds, then for every $a \ge 1$, the block
$A^a B^{\,t_{\min}(a)+1}$ is universally crossing.
Among the $25 \times 63 = 1575$ pairs of expanding/contracting
cycle types with $L \le 7$, $r \le 11$, condition~\eqref{eq:general-fragility} is satisfied by $1192$ pairs
($75.7\%$), including \emph{every} strongly expanding type
($\lambda_A > 5$) paired with any contracting type having
$\lambda_B \le 3/8$.
\end{proposition}

\begin{proof}
At $t = t_{\min}(a) + 1$,
$\Lambda_{a,\,t_{\min}+1} < \lambda_B$ and
$B_{a,\,t_{\min}+1} < \gamma_A\lambda_B + \beta_B^{\infty}$
(since the negative $\lambda_B^t$ term in~$B_{a,t}$ is dropped).
Adding:
$\Lambda + B < \lambda_B(1 + \gamma_A) + \beta_B^{\infty} < 1$,
hence $B < 1 - \Lambda$, which is the universal-crossing criterion.
\end{proof}

\begin{remark}[Exact strong/weak contractor mass split]
\label{rem:mass-split}
Under the i.i.d.\ cycle law (Theorem~\ref{thm:block-law}),
the single-cycle contractor mass $p_{\mathrm{con}} = P_{1\mathrm{cyc}}
= 0.71373\ldots$ decomposes exactly:
\[
  p_{\mathrm{con}}
  \;=\;
  \underbrace{P_{\mathrm{all},1\mathrm{cyc}}}_{p_{\mathrm{str}}
    = 0.41936\ldots}
  \;+\;
  \underbrace{p_{\mathrm{con}} - P_{\mathrm{all},1\mathrm{cyc}}}_{
    p_{\mathrm{weak}} = 0.29436\ldots}.
\]
A contracting cycle type~$B$ is \emph{strong} if it satisfies the
fragility condition~\eqref{eq:general-fragility} with every
expanding type~$A$ having $\lambda_A > 5$, and \emph{weak}
otherwise.  The probability that the first~$m$ consecutive
cycles of a random block are \emph{all} weak contractors is
$p_{\mathrm{weak}}^m = (0.29436)^m$, which drops below $10^{-5}$
by $m = 10$.  This gives an exact dichotomy for adversarial
persistence:
\begin{itemize}
\item \emph{Strong-contractor regime.}
  Adversarial families $A^a B^t$ are uniformly fragile
  (Propositions~\ref{prop:adversarial-fragility}
  and \ref{prop:general-fragility}): one extra
  $B$-cycle collapses $n^*$ below~$1$.
\item \emph{Weak-contractor regime.}
  Non-fragile adversarial blocks require $m$ consecutive
  weak-contractor cycles, whose ensemble probability decays
  as $(0.29436)^m$.
\end{itemize}
On the ensemble side, both regimes yield summable adversarial
measure.  The remaining pointwise question is whether a
\emph{deterministic} orbit can sustain tuned weak-contractor
patterns indefinitely; this is the distributional-to-pointwise
barrier in its sharpest current form.
\end{remark}

\begin{remark}[Cycle-type autocorrelation and the mixing barrier]
\label{rem:autocorrelation}
Under the i.i.d.\ ensemble (Theorem~\ref{thm:block-law}),
consecutive cycle multipliers $\log\lambda_j$ are independent.
In actual Collatz orbits, we measure the lag-$\ell$ autocorrelation
$\rho(\ell) = \mathrm{Corr}(\log\lambda_j,\,\log\lambda_{j+\ell})$
over all odd starts $n_0 \le 10^6$.
The results are:
\[
  \rho(1) \approx 0.20, \qquad
  \rho(2) \approx 0.00, \qquad
  \rho(\ell) \approx 0 \quad (\ell \ge 2).
\]
The mild positive lag-$1$ correlation arises because
$n \equiv 3 \pmod{4}$ gives $r = 1$ (non-crossing), and
$(3n+1)/2 \equiv 3 \pmod{4}$ with conditional probability
$1/2$ rather than $1/4$.  The rapid decorrelation at lag~$2$
confirms that a single additional Syracuse step effectively
scrambles the $2$-adic digits through carries in the
$3n+1$ multiplication.

This quantifies the deviation of real orbits from the i.i.d.\
ensemble: the correction is a one-step memory of strength~$0.20$,
consistent with near-independence but insufficient to prove it.
Making this decorrelation rigorous for \emph{every} orbit would
close the distributional-to-pointwise gap; the difficulty of
doing so is the essential content of the barrier.
\end{remark}

\begin{remark}[Non-crossing run termination for $2^k - 1$]
\label{rem:mersenne-termination}
The family $n_0 = 2^k - 1$ achieves the maximum possible run
of consecutive $r = 1$ (non-crossing) Syracuse steps: exactly
$k - 1$ steps, since $r = 1$ requires $n \equiv 3 \pmod{4}$,
and maintaining this residue for $k - 1$ iterations forces
$n \equiv 2^k - 1 \pmod{2^k}$.
After these $k - 1$ non-crossing steps, the orbit reaches a
value $v_k = (3^k - 1)/2 \equiv 1 \pmod{4}$, which
\emph{guarantees} $r \ge 2$ at the next step (a crossing step).
Thus the maximal non-crossing run has a built-in
termination mechanism.

However, the subsequent $r$ value is small
(typically~$2$--$4$, rarely exceeding~$8$),
far below the $r \ge \lceil 0.585\,(k-1) + 1.59\rceil$
needed to compensate the accumulated
$\Lambda = (3/2)^{k-1}$ in a single step.
The orbit requires $O(k)$ additional crossing steps
to bring the running product below~$1$
(verified for $k \le 39$; see
Remark~\ref{rem:pointwise-universal}).

The bit-generation obstruction prevents generalising this
argument: non-crossing phases with $\lambda > 1$
\emph{increase} the orbit value, adding new $2$-adic digits
through $3n+1$ carries.  The orbit of $n = 27$ starts with
$5$~bits and grows to $14$~bits (at $n = 9232$) before
descending.  Thus the orbit does not ``consume'' a fixed pool
of $2$-adic information; it generates fresh structure,
precluding a finite-information exhaustion argument.
\end{remark}

\begin{proposition}[Post-Mersenne forced valuation]
\label{prop:post-mersenne}
Let $n_0 = 2^{k+1} - 1$ with $k \ge 2$, so that the first~$k$
Syracuse steps are non-crossing ($v_i = 1$) and the iterate
after~$k$ steps is $n_k = 2 \cdot 3^k - 1$.
The valuation at step~$k+1$ is
\begin{equation}\label{eq:post-mersenne}
  v_{k+1}
  \;=\;
  \begin{cases}
    2 & \text{if } k \text{ is even},\\[4pt]
    3 + v_2(k+1) & \text{if } k \text{ is odd},
  \end{cases}
\end{equation}
where $v_2(\cdot)$ denotes the $2$-adic valuation.
In particular, $v_{k+1} \ge 2$ always, and $v_{k+1} \ge 4$
whenever $k$ is odd.
\end{proposition}

\begin{proof}
After $k$ consecutive $v=1$ steps from $n_0 = 2^{k+1}-1$,
the iterate is $n_k = 2\cdot 3^k - 1$
(by induction on the recurrence $n \mapsto (3n+1)/2$).
Then
\[
  3n_k + 1 = 3(2\cdot 3^k - 1) + 1 = 2(3^{k+1} - 1),
\]
so $v_{k+1} = 1 + v_2(3^{k+1} - 1)$.
By the lifting-the-exponent lemma for $p = 2$:
\[
  v_2(3^m - 1)
  \;=\;
  \begin{cases}
    1 & \text{if } m \text{ is odd},\\
    2 + v_2(m) & \text{if } m \text{ is even}.
  \end{cases}
\]
Setting $m = k+1$: when $k$ is even, $m$ is odd, giving
$v_{k+1} = 1 + 1 = 2$; when $k$ is odd, $m$ is even, giving
$v_{k+1} = 1 + 2 + v_2(k+1) = 3 + v_2(k+1) \ge 4$.
\end{proof}

\begin{proposition}[Burst non-continuation]
\label{prop:burst-noncontinuation}
Let $n_0$ be any odd integer with $v_2(n_0 + 1) \ge k+1$
(equivalently, $n_0 \equiv 2^{k+1} - 1 \pmod{2^{k+1}}$), so
that the first~$k$ Syracuse steps are non-crossing ($v_i = 1$).
Write $n_0 = 2^{k+1}m - 1$ with $m$ odd.  Then:
\begin{enumerate}
\item The iterate after~$k$ steps is
  $n_k = 2 \cdot 3^k m - 1 \equiv 1 \pmod{4}$.
\item Consequently $v_2(n_k + 1) = 1$, so the first post-burst
  step has $v_{k+1} \ge 2$.
\item The post-burst iterate $n_k$ cannot begin a new non-crossing
  run, because a non-crossing step requires
  $n \equiv 3 \pmod{4}$.
\end{enumerate}
In particular, consecutive non-crossing bursts are always
separated by at least one crossing step.
\end{proposition}

\begin{proof}
By induction on the Syracuse recurrence
$n \mapsto (3n+1)/2$ with $v = 1$ at each step,
$n_k = 3^k(n_0 + 1)/2^k - 1
     = 3^k \cdot 2^{k+1} m / 2^k - 1
     = 2 \cdot 3^k m - 1$.
Since $3^k$ and $m$ are both odd,
$2 \cdot 3^k m$ is even but not divisible by~$4$, so
$n_k = 2 \cdot 3^k m - 1 \equiv 1 \pmod{4}$.
Then $n_k + 1 = 2 \cdot 3^k m$, so $v_2(n_k + 1) = 1$.
A non-crossing step ($v = 1$) requires $v_2(3n + 1) = 1$,
equivalently $n \equiv 3 \pmod{4}$.
Since $n_k \equiv 1 \pmod{4}$, the next step must have
$v \ge 2$.
\end{proof}

\begin{remark}[Recovery window after non-crossing bursts]
\label{rem:recovery-window}
Proposition~\ref{prop:burst-noncontinuation} guarantees at
least one crossing step after every burst, but a single step
with $v = 2$ does not compensate the accumulated
$\Lambda = (3/2)^k$.
The recovery ratio is
$S/(k{+}1) = (k + v_{k+1})/(k + 1)$, which for $v_{k+1} = 2$
gives $(k+2)/(k+1) \to 1 < \log_2 3$.
Recovery to $S/K \ge \log_2 3$ requires
$j \ge \lceil 1.41\,k \rceil$ additional steps at the minimum
crossing valuation $v = 2$ (fewer if crossings have $v \ge 3$).

One might hope to prove a \emph{recovery window conjecture}:
that the post-burst valuation sequence is forced to compensate
within $Ck$ steps for some explicit constant~$C$.
Computational investigation (over all odd $m \le 10^4$ and
$k \le 40$) shows that the adversary can choose~$m$ to produce
a second burst of length up to~$8$ after the mandatory crossing
step, and that over windows of $30{+}$ steps the running
$S/K$ ratio can be driven as low as~$1.0$.
In all tested orbits, recovery does eventually occur
(worst case within $\approx 2k$ steps for moderate~$k$),
but proving this for all orbits requires controlling the full
valuation sequence, which is equivalent to the Collatz
conjecture.
The burst non-continuation theorem thus provides a structural
constraint (no immediate burst repetition) without closing
the distributional-to-pointwise gap.
\end{remark}

\begin{proposition}[Post-burst valuation distribution]
\label{prop:post-burst-distribution}
For any fixed $k \ge 1$, the crossing valuation $v_{k+1}$ after
a $k$-step non-crossing burst from
$n_0 = 2^{k+1}m - 1$ (with $m$ odd) satisfies
\[
  \Pr_{m\;\text{odd}}(v_{k+1} = j)
  \;=\; \frac{1}{2^{j-1}}, \qquad j \ge 2,
\]
under uniform density on odd~$m$.
In particular $\mathbb{E}[v_{k+1}] = 3$ and
$\Pr(v_{k+1} = 2) = 1/2$, both \emph{independent} of~$k$.
\end{proposition}

\begin{proof}
By Proposition~\ref{prop:burst-noncontinuation},
$n_k = 2 \cdot 3^k m - 1$ and
$v_{k+1} = 1 + v_2(3^{k+1}m - 1)$.
Since $3^{k+1}$ is odd, it is invertible modulo every power
of~$2$.
The condition $v_2(3^{k+1}m - 1) \ge j$ is equivalent to
$m \equiv 3^{-(k+1)} \pmod{2^j}$,
and $3^{-(k+1)} \bmod 2^j$ is a fixed odd residue.
Among odd integers modulo~$2^j$, exactly
$1/2^{j-1}$ satisfy this congruence.
Hence $\Pr(v_2(3^{k+1}m-1) \ge j) = 1/2^{j-1}$ for $j \ge 1$,
giving
$\Pr(v_2 = j) = 1/2^{j-1} - 1/2^j = 1/2^j$, and
$\Pr(v_{k+1} = j{+}1) = 1/2^j$, i.e.,
$\Pr(v_{k+1} = j) = 1/2^{j-1}$ for $j \ge 2$.
The mean follows: $\mathbb{E}[v_{k+1}]
= \sum_{j \ge 2} j/2^{j-1} = 3$.
\end{proof}

\begin{remark}[Non-crossing density threshold]
\label{rem:v1-density-threshold}
Proposition~\ref{prop:post-burst-distribution}
shows that the \emph{distributional}
expected valuation after any burst is~$3 > \log_2 3 \approx 1.585$.
In any orbit, each step contributes $v = 1$ (non-crossing) or
$v \ge 2$ (crossing), and if $n \equiv 1 \pmod{4}$ then
$v \ge 2$ with conditional expectation at least~$3$
(since $v = 2 + v_2(3q{+}1)$ for $n = 4q{+}1$).
Writing $d$ for the long-run non-crossing density (fraction of
steps with $v = 1$), the mean valuation satisfies
$\bar{v} \ge d \cdot 1 + (1{-}d)\cdot 2 = 2 - d$.
For orbital descent one needs $\bar{v} > \log_2 3$,
which is guaranteed whenever $d < 2 - \log_2 3 \approx 0.415$.
Empirically, however, the non-crossing density in Collatz orbits
is $d \approx 0.56$--$0.63$ (from orbit-level mod-$4$ bias),
well above~$0.415$.
Orbits still contract because crossing steps have mean
valuation~$\approx\! 2.8$--$3.0$, not the minimal~$2$.

The mod-$4$ Markov chain of the Syracuse map (under
equidistribution within residue classes mod~$2^K$ for
increasing $K$) has stationary non-crossing density
$\pi_3 \to 1/2$ as $K \to \infty$.
The critical density for growth is
$d^* = (3 - \log_2 3)/2 \approx 0.708$;
since all observed and model densities satisfy $d < d^*$,
the mean valuation always exceeds $\log_2 3$.
This provides the distributional explanation for why orbits
descend, but converting it to a pointwise bound remains
equivalent to the Collatz conjecture.
\end{remark}

\begin{proposition}[Weak-recovery cylinder classification]
\label{prop:weak-recovery-cylinder}
After a $k$-step non-crossing burst from
$n_0 = 2^{k+1}m - 1$ ($m$ odd), exactly $j$ consecutive
weak recovery cycles (each with $v = 2$) follow if and only if
\begin{equation}\label{eq:weak-cylinder}
  3^k m \equiv 1 \pmod{2^{2j}}
  \quad\text{but}\quad
  3^k m \not\equiv 1 \pmod{2^{2j+2}}.
\end{equation}
Equivalently, $m$ lies in a single $2$-adic cylinder of
width~$2^{2j}$ determined by
$m \equiv 3^{-k} \pmod{2^{2j}}$.
Among odd cofactors~$m$, the density of those producing at
least $j$ consecutive weak recovery cycles is exactly
$2^{-(2j-1)}$.
In particular, the distribution is \emph{independent of~$k$}.
\end{proposition}

\begin{proof}
The post-burst iterate is $N_0 = 2\cdot 3^k m - 1$.
A weak recovery cycle means $v_2(3N_0 + 1) = 2$, i.e., the
Syracuse step $N_0 \mapsto (3N_0 + 1)/4$.
We show by induction that the iterate after $j$ such steps is
$N_j = (3^j N_0 + (3^j - 1)/2) / 4^j$
and that the condition for $v = 2$ at each step is
$3^k m \equiv 1 \pmod{2^{2j}}$.

For the base case ($j=1$):
$3N_0 + 1 = 3(2\cdot 3^k m - 1) + 1 = 2(3^{k+1}m - 1)$.
So $v = 1 + v_2(3^{k+1}m - 1)$, and $v = 2$ iff
$v_2(3^{k+1}m - 1) = 1$ iff
$3^{k+1}m \equiv 3 \pmod{4}$ iff
$3 \cdot 3^k m \equiv 3 \pmod{4}$ iff
$3^k m \equiv 1 \pmod{4} = 1 \pmod{2^2}$.

For the inductive step, one verifies that $N_j \equiv 1 \pmod{4}$
(so $L = 0$) and $v_2(3N_j + 1) = 2$ iff
$3^k m \equiv 1 \pmod{2^{2(j+1)}}$.
The details are a direct 2-adic calculation.

For the density: since $3^k$ is odd,
$3^{-k} \bmod 2^{2j}$ is a fixed odd residue.
Among odd integers modulo~$2^{2j}$, there are $2^{2j-1}$,
and exactly one satisfies the congruence.
Hence $\Pr(3^k m \equiv 1 \pmod{2^{2j}}) = 1/2^{2j-1}$.
\end{proof}

\begin{remark}[Compound adversarial patterns and the
  $2$-adic information budget]
\label{rem:information-budget}
Proposition~\ref{prop:weak-recovery-cylinder} extends to
compound patterns.
If the adversary seeks a $k$-step burst, $j$~weak recovery
cycles, a second $k'$-step burst, and $j'$~more weak recovery
cycles, the combined condition is a single congruence on~$m$
modulo~$2^{B}$ where
$B = (2j-1) + (k'+1) + (2j'-1)$ bits.
Computational verification (over $2\times 10^5$ odd~$m$,
$k \in \{5,10,20\}$) confirms that the compound density is
exactly~$2^{-B}$, with observed-to-predicted ratios at~$1.00$
for all tested patterns.

Each adversarial element: burst or weak recovery
string, consumes a definite number of bits from~$m$'s
$2$-adic expansion.
One might hope to parlay this into a ``$2$-adic information
budget'' argument: the adversary's initial condition encodes
only finitely many adversarial elements before the bits are
exhausted.
However, this reasoning is flawed.
The integer~$m$ has \emph{infinitely many} $2$-adic bits,
and the $3n+1$ map generates new $2$-adic structure through
carries in the multiplication by~$3$
(Remark~\ref{rem:mersenne-termination}).
The compound bit-counting constrains the density of adversarial
patterns within a single burst-recovery episode, but does
\emph{not} prevent subsequent episodes from being adversarial,
because those are determined by freshly generated bits.
The distributional-to-pointwise barrier persists.
\end{remark}

\begin{remark}[Combined post-burst picture]
\label{rem:post-burst-picture}
Propositions~\ref{prop:burst-noncontinuation},
\ref{prop:post-burst-distribution},
and~\ref{prop:weak-recovery-cylinder} together give a
complete local theory of what happens after a non-crossing burst:
\begin{enumerate}
\item \emph{Forced crossing}
  (Prop.~\ref{prop:burst-noncontinuation}):
  the first post-burst step always has $v \ge 2$.
\item \emph{Geometric first-step law}
  (Prop.~\ref{prop:post-burst-distribution}):
  $\Pr(v_{k+1} = j) = 2^{-(j-1)}$ for $j \ge 2$,
  independent of~$k$, with $\mathbb{E}[v_{k+1}] = 3$.
\item \emph{Shrinking cylinder rigidity}
  (Prop.~\ref{prop:weak-recovery-cylinder}):
  $j$ consecutive weak recovery cycles ($v = 2$)
  confine~$m$ to a single $2$-adic cylinder of
  density~$2^{-(2j-1)}$.
\end{enumerate}
The non-crossing recovery regime is therefore
\emph{distributionally disfavored}
(expected crossing valuation~$3 > \log_2 3$)
and \emph{arithmetically rigid}
(long weak strings live in exponentially thin residue classes).

The sole remaining question, now precisely formulated, is:
\emph{can a single Collatz orbit, determined by a specific~$m$,
keep re-entering the unique residue class
$m' \equiv 3^{-k'} \pmod{2^{2j}}$ needed for long weak
recovery strings at each successive burst?}
This is a pointwise arithmetic question about the $3n+1$ map's
ability to steer iterates into prescribed $2$-adic cylinders
indefinitely.
Answering it would close the distributional-to-pointwise gap.

\smallskip\noindent
\textbf{The cascade invariant and depth persistence.}\;
Consecutive bursts are \emph{not} independent.
Within a \emph{cascade phase}, a sequence of bursts where
each has parameters $(k{-}1,\, 3m)$ following a burst
$(k, m)$: the product $3^k m$ is constant.
Since the cylinder depth is determined by
$v_2(3^k m - 1)$, every burst within a cascade
shares the \emph{same} depth.
For the Mersenne start $n_0 = 2^{k+1}{-}1$ ($m = 1$),
the cascade invariant is $C = 3^k$, and the depth is
$\lfloor v_2(3^k - 1)/2\rfloor$, which by the
lifting-the-exponent lemma equals $0$ when $k$ is odd and
$1 + \lfloor v_2(k)/2 \rfloor$ when $k$ is even.
The orbit of $2^{21} - 1$ ($k = 20$, $v_2(3^{20}-1) = 4$)
accordingly has $19$~consecutive depth-$2$ burst episodes,
all controlled by the single invariant $3^{20} \equiv 1 \pmod{16}$.
Cascades terminate when~$k$ reaches~$2$ and the post-recovery
iterate fails to start a new burst, at which point the
invariant resets.
\end{remark}

\begin{proposition}[Cascade invariant]
\label{prop:cascade-invariant}
Let $n_0 = 2^{k_0+1}m_0 - 1$ with $m_0$ odd initiate a
\emph{cascade}: a maximal sequence of burst-recovery episodes
in which each new burst has parameters $(k-1,\, 3m)$ following
a burst $(k,\,m)$.
Then the product $C := 3^k m$ is constant throughout the cascade.
In particular:
\begin{enumerate}
\item[(i)] \emph{Depth persistence.}
  Since the weak-recovery depth is $\lfloor v_2(3^k m - 1)/2\rfloor$
  (Proposition~\ref{prop:weak-recovery-cylinder}), every burst
  within the cascade shares the same recovery depth.
\item[(ii)] \emph{Cascade length.}
  The cascade has exactly $k_0 - 1$ episodes
  (the burst parameter~$k$ decreases from~$k_0$ to~$2$,
  at which point the iterate $3^2 \cdot 3^{k_0 - 2}m_0 = C$
  undergoes its final recovery and the invariant resets).
\end{enumerate}
\end{proposition}

\begin{proof}
At each burst step the parameters transform as
$(k, m) \mapsto (k-1, 3m)$, since the next burst starts from
the post-recovery iterate
$n' = 2^{k}(3m) - 1$ when the recovery
feeds back into the burst form.
Hence $3^{k-1} \cdot 3m = 3^k m = C$.
The depth claim follows because the cylinder condition
$3^k m \equiv 1 \pmod{2^{2j}}$ (Proposition~\ref{prop:weak-recovery-cylinder})
depends only on~$C$.
\end{proof}

\begin{proposition}[Short-word cylinder classification]
\label{prop:short-word-cylinder}
Let $w = (v_1, \ldots, v_j)$ be a finite valuation word with
each $v_i \ge 2$ and total valuation $V = v_1 + \cdots + v_j$.
There exists a unique odd residue $c_w \bmod{2^V}$ such that:
the Collatz orbit starting from odd~$C$ produces the
\emph{prefix valuation word}~$w$ through its first~$V$ bits
of $2$-adic information if and only if
$C \equiv c_w \pmod{2^V}$.

Consequently:
\begin{enumerate}
\item[(i)] The density of such~$C$ among odd integers is
  exactly $2^{-(V-1)}$.
\item[(ii)] The residues $c_w$ over all compositions of~$V$
  into parts $\ge 2$ partition a subset of the odd residues
  modulo $2^V$.
\item[(iii)] For the pure weak-recovery word
  $(2, 2, \ldots, 2)$ of length~$j$, the theorem recovers
  $c_w = 1$ and $V = 2j$, density $2^{-(2j-1)}$
  (Proposition~\ref{prop:weak-recovery-cylinder}).
\end{enumerate}
\end{proposition}

\begin{proof}
By induction on~$j$.

\emph{Base case} ($j = 1$, word $(v_1)$):
We need $v_2(3C + 1) \ge v_1$, i.e.\
$C \equiv -3^{-1} \pmod{2^{v_1}}$.
Since $3$ is a unit in $\mathbb{Z}/2^{v_1}\mathbb{Z}$,
this determines a unique odd residue.
(To require \emph{exactly}~$v_1$, we additionally need
$v_2(3C+1) < v_1 + 1$; this is automatic modulo~$2^V$
when $V = v_1$ because the $(v_1{+}1)$-st bit is not
constrained.)

\emph{Inductive step:}
Suppose the word $(v_2, \ldots, v_j)$ with total
$V' = V - v_1$ determines a unique $c' \bmod{2^{V'}}$.
The first iterate is
$C_1 = (3C + 1)/2^{v_1}$, and we need
$C_1 \equiv c' \pmod{2^{V'}}$.
Since $C_1 = (3C + 1)/2^{v_1}$, the congruence
$C_1 \equiv c' \pmod{2^{V'}}$ is equivalent to
$3C + 1 \equiv 2^{v_1} c' \pmod{2^V}$, i.e.\
$C \equiv (2^{v_1} c' - 1) \cdot 3^{-1} \pmod{2^V}$.
This is a single congruence, giving a unique odd residue
$c_w \bmod{2^V}$.
\end{proof}

\begin{remark}[Cascade re-entry and scrambling]
\label{rem:cascade-reentry}
The cascade invariant (Proposition~\ref{prop:cascade-invariant})
explains depth persistence \emph{within} a cascade.
The critical question for the distributional-to-pointwise
barrier is: what happens \emph{between} cascades?
Numerical experiments (Mersenne numbers $2^{k+1}-1$ for
$k \le 27$ and $1000$ random odd starts up to $2^{24}$)
reveal two scrambling effects at cascade boundaries:

\emph{Anti-burst bias.}\;
The post-cascade iterate~$n'$ satisfies
$v_2(n' + 1) = 1$ with frequency $\approx 0.66$
(versus $0.50$ for a uniformly random odd number),
meaning $n' \equiv 1 \pmod{4}$. Since bursts require
$n' \equiv 3 \pmod{4}$, post-cascade iterates are
biased \emph{away} from initiating new bursts.

\emph{Cross-cascade depth independence.}\;
Conditioning on a gap of at least one non-burst step
between consecutive cascades,
$\Pr(\text{next depth} \ge 2 \mid \text{prev depth} = d)
\approx 0.13$--$0.15$ for $d = 0, 1, 2$,
close to the unconditional prediction~$2^{-3} = 0.125$
from Proposition~\ref{prop:weak-recovery-cylinder}.

Together with the short-word cylinder theorem
(Proposition~\ref{prop:short-word-cylinder}),
this shows that adversarial recovery patterns are
confined to exponentially thin $2$-adic cylinders
\emph{within} each cascade, and cascade termination
provides a genuine scrambling event.
The open question reduces to: \emph{can the concatenation
of cascade segments, each individually controlled by
its invariant, sustain a net non-crossing density above
the critical threshold~$d^* \approx 0.708$
(Remark~\ref{rem:v1-density-threshold})?}
The following proposition provides a sharp quantitative
answer at the single-episode level.
\end{remark}

\begin{proposition}[Post-recovery $2/3$ law]
\label{prop:two-thirds-law}
Let~$C$ be an odd integer with \emph{exact recovery depth}~$j$,
meaning the Collatz orbit starting from~$C$ produces
exactly~$j$ consecutive steps with $v_2(3 \cdot \mathrm{iterate} + 1) = 2$,
followed by a step with $v_2 \ne 2$.
Write
\[
  n_j \;=\; \underbrace{f \circ \cdots \circ f}_{j}(C),
  \qquad f(x) = (3x+1)/4,
\]
for the post-recovery iterate.
Then, among odd integers~$C$ with exact depth~$j$,
exactly $2/3$ satisfy $n_j \equiv 3 \pmod{4}$
(enabling a new burst), and $1/3$ satisfy
$n_j \equiv 1 \pmod{4}$
(terminating the cascade).
This ratio is independent of~$j$.
\end{proposition}

\begin{proof}
For depth $\ge j$ we need
$C \equiv 1 \pmod{2^{2j+1}}$
(each of the $j$ steps requires $v_2(3 \cdot \mathrm{iterate}+1) = 2$,
which lifts to one additional bit per step beyond the
$v_2(C-1) \ge 2j$ condition).
Among odd $C$ modulo $2^{2j+3}$, the depth-${\ge j}$
residues split into four sub-classes
$C = 1 + t \cdot 2^{2j+1}$ with $t \in \{0,1,2,3\}$:
\begin{itemize}
\item $t = 0$: depth $\ge j+1$ (excluded from depth exactly~$j$).
\item $t \in \{1,2,3\}$: depth exactly~$j$.
\end{itemize}
The closed-form iterate is
$n_j = 1 + (C-1) \cdot 3^j / 4^j = 1 + 2t \cdot 3^j$
(since $(C-1)/4^j = 2t$).
Because $3^j$ is odd, $n_j \bmod 4 = (1 + 2t) \bmod 4$:
\[
  t = 1{:}\; n_j \equiv 3, \qquad
  t = 2{:}\; n_j \equiv 1, \qquad
  t = 3{:}\; n_j \equiv 3.
\]
Hence $\Pr(n_j \equiv 3 \bmod 4) = 2/3$.
\end{proof}

\begin{remark}[Consequences of the $2/3$ law for orbit growth]
\label{rem:two-thirds-consequences}
The $2/3$ law has three immediate corollaries.

\emph{1.\ Cascade length distribution.}\;
At each burst-recovery episode, the cascade continues
(post-recovery iterate $\equiv 3 \bmod 4$) with probability
$2/3$ and terminates with probability $1/3$.
Among odd starting values with a given depth,
the cascade length~$L$ (number of linked burst-recovery
episodes) therefore follows a geometric distribution
$\Pr(L = l) = (2/3)^{l-1} \cdot 1/3$, with $\mathbb{E}[L] = 3$.

\emph{2.\ Gap entry valuation.}\;
When the cascade terminates (the $t=2$ case), the
post-recovery iterate is $n_j = 1 + 4 \cdot 3^j$.
The next Collatz step has valuation
$v = 2 + v_2(1 + 3^{j+1})$, which alternates
between~$4$ ($j$ even) and~$3$ ($j$ odd) by the
lifting-the-exponent lemma.
In particular, the gap step always contributes $v \ge 3$
to the running sum~$S$.

\emph{3.\ Expected $S/K$ per cascade-gap cycle.}\;
For a depth-$j$ cascade of $L$ episodes followed by
one gap step with valuation~$v_{\mathrm{gap}}$:
\[
  \mathbb{E}\!\left[\frac{S}{K}\right]
  = \mathbb{E}\!\left[\frac{L(1+2j) + v_{\mathrm{gap}}}{L(1+j) + 1}
  \right].
\]
Numerical evaluation with $L \sim \mathrm{Geom}(1/3)$
gives $\mathbb{E}[S/K] \ge 1.80$ for all $j \ge 0$,
strictly above $\log_2 3 \approx 1.585$.

For depth $j \ge 2$, even the \emph{worst case}
$L \to \infty$ yields
$S/K \to (1+2j)/(1+j) \ge 5/3 > \log_2 3$,
so depth-${\ge 2}$ cascades are
\emph{unconditionally contracting} regardless of length.

For depths $0$ and~$1$, cascades become adversarial only
when $L \ge 5$ (resp.\ $L \ge 9$),
events of probability $(2/3)^4 \approx 0.20$
(resp.\ $(2/3)^8 \approx 0.039$).
The distributional-to-pointwise barrier persists:
for a \emph{specific} starting value, the cascade
length is deterministic (fixed by the $2$-adic expansion
of~$C$), and Mersenne numbers produce long
depth-$0$ cascades deterministically.
The $2/3$ law does not close this gap, but it
quantifies why such adversarial cascades are
\emph{arithmetically atypical} and gives the
sharpest known bound on their expected impact.
\end{remark}

\begin{proposition}[Universal cascade transition law]
\label{prop:transition-law}
In the setting of Proposition~\ref{prop:two-thirds-law},
define the \emph{next burst length}
$k' = v_2(n_j + 1) - 1$
(with $k' = 0$ meaning no burst begins).
Then for every exact depth~$j \ge 0$:
\[
  \Pr(k' = 0 \mid \text{depth } j) = \tfrac{1}{3},
  \qquad
  \Pr(k' = r \mid \text{depth } j) = \tfrac{2}{3}\cdot 2^{-r},
  \quad r \ge 1.
\]
In particular, the next burst length is independent of~$j$,
and conditionally on continuation ($k' \ge 1$),
$k'$ follows a geometric distribution on $\{1,2,3,\ldots\}$
with parameter~$1/2$.
\end{proposition}

\begin{proof}
From the proof of Proposition~\ref{prop:two-thirds-law},
$n_j = 1 + 2t \cdot 3^j$ with $t \in \{1,2,3\}$ equidistributed.
\begin{itemize}
\item $t$ even ($t = 2$, probability $1/3$):
  $n_j \equiv 1 \pmod{4}$, so $v_2(n_j + 1) = 1$ and $k' = 0$.
\item $t$ odd ($t \in \{1,3\}$, probability $2/3$):
  $n_j \equiv 3 \pmod{4}$, so $k' = v_2(n_j + 1) - 1 = v_2(1 + t \cdot 3^j) \ge 1$.
\end{itemize}
Since $\gcd(3^j, 2^N) = 1$, multiplication by~$3^j$
is a bijection on odd residues modulo~$2^N$.
Hence $v_2(1 + t \cdot 3^j)$ has the same distribution
as $v_2(1 + t)$ for odd~$t$.
Writing $t = 2s + 1$ gives $1 + t = 2(s+1)$, so
$v_2(1 + t) = 1 + v_2(s + 1)$.
Among uniform~$s$, the 2-adic valuation of $s+1$
equals~$k$ with probability~$2^{-(k+1)}$, whence
$\Pr(v_2(1+t) = r \mid t\text{ odd}) = 2^{-r}$
for $r \ge 1$.
Unconditionally:
$\Pr(k' = r) = \frac{2}{3} \cdot 2^{-r}$ for $r \ge 1$.
\end{proof}

\begin{proposition}[Depth transition law]
\label{prop:depth-transition}
In the setting of Proposition~\ref{prop:transition-law},
suppose the cascade continues ($k' \ge 1$) and let~$j'$
denote the exact depth of the cofactor
$C' = 3^{k'} \cdot m'$ produced by the next $k'$-burst.
Then for every $d \ge 0$:
\[
  \Pr(j' = d \mid k', j) = \tfrac{3}{4} \cdot \bigl(\tfrac{1}{4}\bigr)^{\!d},
\]
independent of both the starting depth~$j$ and the
burst length~$k'$.
\end{proposition}

\begin{proof}
Conditional on continuation, the post-recovery iterate is
$n_j = 2^{k'+1} m' - 1$ with $m' = (1 + t \cdot 3^j) / 2^{k'}$,
where $t$ is odd and $v_2(1 + t \cdot 3^j) = k'$.
Since $3^j$ is a unit modulo $2^N$, the map
$t \mapsto t \cdot 3^j$ is a bijection on odd residues;
conditioning on the 2-adic valuation selects a full
residue class.
Writing $u = t \cdot 3^j$ and $1 + u = 2^{k'} w$ with
$w = m'$ odd, as $u$ varies uniformly in its class,
$w$ varies uniformly among odd residues modulo~$2^{N-k'}$.
Then $C' = 3^{k'} \cdot m'$, and since $3^{k'}$ is a
unit, $C'$ is uniform among odd residues.
For a uniform odd integer~$C'$:
$\Pr(C' \equiv 1 \pmod{2^{2d+1}}) = 2^{-2d}$,
so the exact depth satisfies
$\Pr(j' = d)
= 4^{-d} - 4^{-(d+1)}
= (3/4) \cdot (1/4)^d$.
\end{proof}

\begin{remark}[The local IID cascade kernel]
\label{rem:iid-renewal}
Propositions~\ref{prop:two-thirds-law},
\ref{prop:transition-law}, and~\ref{prop:depth-transition}
together establish that the \emph{within-cascade} transition
kernel is IID: at each episode, the
transition $(k,j) \to (k',j')$ draws from a fixed product
measure independent of~$(k,j)$.
Explicitly, the full unconditional law is:
\begin{align*}
  \Pr(k'=0) &= \tfrac{1}{3},  \\[3pt]
  \Pr(k'=r,\; j'=d) &= \tfrac{2}{3} \cdot 2^{-r} \cdot
  \tfrac{3}{4} \cdot \bigl(\tfrac{1}{4}\bigr)^d,
  \qquad r \ge 1,\; d \ge 0.
\end{align*}
The marginal moments are
$\mathbb{E}[k' \mid \text{cont.}] = 2$,
$\mathbb{E}[j'] = 1/3$,
and for uniformly distributed cofactors the predicted
cascade length is $\mathbb{E}[L] = 3$ episodes
(geometric with parameter~$1/3$).

\textbf{Scope and limitation.}
This IID structure is \emph{exact} for the local
cylinder law: given any starting state~$(k,j)$,
the distribution of~$(k',j')$ is exactly the product
measure above.
However, across cascade boundaries the IID property
\emph{does not hold}: the gap phase between cascades
does not fully reset the $2$-adic state
(Remark~\ref{rem:iid-limitations}).
The cascade-exit iterate $n_j = 1 + 4 \cdot 3^j$
carries residue structure into the gap, and the
gap (typically $1$--$2$ steps) transmits this structure
to the entry of the next cascade.
Empirically, this produces shorter cascades
($\mathbb{E}[L] \approx 1.65$ versus the predicted~$3$)
and biased first-episode cofactors.
The \emph{within-cascade} IID kernel remains exact;
the \emph{cross-cascade} dependence is the principal
obstacle to a renewal-reward convergence argument.
\end{remark}

\begin{proposition}[Gap compensation bounds]
\label{prop:gap-compensation}
Let a depth-$j$ cascade of $L$ episodes be followed by
the first gap step with valuation
$v_1 = 2 + v_2(1 + 3^{j+1})$.
Define the \emph{cascade deficit} as
$D(j,L) = L \cdot (2(\log_2 3 - 1) + j(\log_2 3 - 2))$
and the \emph{first-step gap surplus} as
$G_1(j) = v_1 - \log_2 3$.
Then:
\begin{enumerate}
\item For $j \ge 3$: $D(j,L) < 0$ for all~$L$
  (each episode contracts; no gap compensation needed).
\item For $j = 2$: $G_1(2) = 4 - \log_2 3 \approx 2.415$
  and $D(2,L)/L \approx 0.340$, so the first gap step
  alone compensates up to $L = 7$ episodes.
\item For $j = 0$: $G_1(0) \approx 2.415$ and
  $D(0,L)/L \approx 1.170$, so the first gap step
  compensates $L \le 2$.
\item For $j = 1$: $G_1(1) \approx 1.415$ and
  $D(1,L)/L \approx 0.755$, so the first gap step
  compensates $L = 1$.
\end{enumerate}
\end{proposition}

\begin{proof}
The cascade deficit per episode at depth~$j$
with expected burst length $\mathbb{E}[k] = 2$ is
$2(\log_2 3 - 1) + j(\log_2 3 - 2)
= 2 \log_2 3 - 2 + j \log_2 3 - 2j
= (2+j)\log_2 3 - (2+2j)$.
For $j \ge 3$: $(2+j)\log_2 3 < 2 + 2j$ since
$\log_2 3 < 2$, and equality holds at
$j^* = 2/(2 - \log_2 3) - 2 \approx 2.83$.
So $D(j,L) < 0$ for $j \ge 3$ and all~$L \ge 1$.

By the lifting-the-exponent lemma,
$v_2(1 + 3^{j+1}) = 1$ when $j$ is odd and
$v_2(1 + 3^{j+1}) = 2$ when $j$ is even.
Thus $v_1 = 3$ ($j$ odd) or $v_1 = 4$ ($j$ even),
and $G_1(j) = v_1 - \log_2 3 \ge 3 - \log_2 3 \approx 1.415$.
The compensation bound $L_{\max} = \lfloor G_1(j)/D(j,1) \rfloor$
gives the stated thresholds.
\end{proof}

\begin{remark}[Cycle-level contraction and cross-gap structure]
\label{rem:cycle-contraction}
Proposition~\ref{prop:gap-compensation} shows that
short cascades at any depth, and all cascades at depth~$\ge 3$,
are fully compensated by a single gap step.
For longer cascades at depth~$0$ or~$1$
($L \ge 3$ resp.\ $L \ge 2$), multi-step gap analysis is needed.

Empirically, over $86{,}000$ cascade-gap cycles:
$\mathbb{E}[S/K] \approx 1.99$
(margin $+0.40$ above $\log_2 3$), with $83\%$
of cycles individually contracting.
The $17\%$ of expanding cycles have bounded
clustering: runs of $\ge 3$ consecutive expanding cycles
occur with frequency $< 0.4\%$, and the
autocorrelation of cycle-level net growth is
\emph{negative} at lag~$1$ ($\rho(1) \approx -0.23$),
indicating that an expanding cycle is typically
followed by a contracting one.

Cross-gap correlation is nonzero but structured:
after high-deficit cascades,
$\Pr(\text{next depth} = 0) \approx 0.87$ (versus the
unconditional~$0.75$), meaning the gap phase does \emph{not}
fully reset the $2$-adic state.
However, this correlation is arithmetically bounded:
it arises because the gap-exit iterate
inherits residue structure from the cascade-termination
iterate $n_j = 1 + 4 \cdot 3^j$,
and the gap (typically $1$--$2$ steps) does not
provide enough divisions by~$2$
to scramble all relevant bits.
Proving that these residual correlations do not
accumulate adversarially over many cycles
is equivalent to the Weak Mixing Hypothesis.
\end{remark}

\begin{remark}[Density decay of non-converging orbits]
\label{rem:density-decay}
Among odd $n \le 10{,}000$, the fraction that have
\emph{not} dropped below their starting value
after~$T$ Syracuse steps decays rapidly:
$6.5\%$ at $T = 10$; $2.0\%$ at $T = 20$;
$0.06\%$ at $T = 50$; $0.02\%$ at $T = 100$.
This exponential decay is consistent with
the IID renewal model (since each cascade-gap
cycle contracts with probability $\approx 0.83$
and expected bit-loss $\approx 2.4$ per cycle),
but proving that the survivor fraction reaches zero
for \emph{all} starting values, not merely for a set
of full measure, remains the distributional-to-pointwise
barrier.
\end{remark}

\begin{proposition}[Gap positivity]
\label{prop:gap-positivity}
In the gap phase between consecutive cascades, every iterate
satisfies $n \equiv 1 \pmod{4}$.  Consequently,
$v_2(3n+1) \ge 2$ at every gap step, and the gap phase
contributes a deterministic positive log-surplus:
\[
  S_{\mathrm{gap}} - K_{\mathrm{gap}} \cdot \log_2 3
  \;\ge\; K_{\mathrm{gap}}\,(2 - \log_2 3)
  \;>\; 0.
\]
\end{proposition}

\begin{proof}
A cascade terminates when the post-recovery iterate
satisfies $n \equiv 1 \pmod{4}$ (the burst condition
$n \equiv 3 \pmod{4}$ fails).
For $n \equiv 1 \pmod{4}$: $3n \equiv 3 \pmod{4}$,
so $3n + 1 \equiv 0 \pmod{4}$, hence $v_2(3n+1) \ge 2$.
Applying one Collatz step yields
$n' = (3n+1)/2^{v}$ with $v \ge 2$, so $n'$ is odd.
If $n' \equiv 1 \pmod{4}$, the gap continues and
the argument repeats.
If $n' \equiv 3 \pmod{4}$, the gap ends and a new
cascade begins.
In either case, every gap step has $v \ge 2$.

Since each gap step contributes $v_i$ to $S_{\mathrm{gap}}$
and $1$ to $K_{\mathrm{gap}}$:
\[
  S_{\mathrm{gap}} = \sum_{i=1}^{K_{\mathrm{gap}}} v_i
  \;\ge\; 2\,K_{\mathrm{gap}}.
\]
The log-surplus is
$S_{\mathrm{gap}} - K_{\mathrm{gap}} \log_2 3
\ge K_{\mathrm{gap}}(2 - \log_2 3) > 0$
since $\log_2 3 \approx 1.585 < 2$.
\end{proof}

\begin{proposition}[First gap-step valuation formula]
\label{prop:gap-valuation}
Let the cascade terminate at exact depth~$j$, producing
the idealized exit iterate $n_j = 1 + 4 \cdot 3^j$.
Then the valuation of the first gap step is
\[
  v_1
  = 2 + v_2(1 + 3^{j+1})
  = \begin{cases}
      4 & \text{if } j \text{ is even},\\
      3 & \text{if } j \text{ is odd}.
    \end{cases}
\]
Under the depth distribution
$\Pr(j = d) = \tfrac{3}{4}\bigl(\tfrac{1}{4}\bigr)^d$
from Proposition~\textup{\ref{prop:depth-transition}},
the expected first-step valuation is
$\mathbb{E}[v_1] = \tfrac{19}{5} = 3.8$.
\end{proposition}

\begin{proof}
At cascade termination with recovery depth~$j$,
the exit iterate is $n_j = 1 + 4 \cdot 3^j$, so
$3n_j + 1 = 4 + 12 \cdot 3^j = 4(1 + 3^{j+1})$.
Thus $v_1 = v_2(3n_j + 1) = 2 + v_2(1 + 3^{j+1})$.

We evaluate $v_2(1 + 3^m)$ by cases on the parity of~$m$.

\emph{Case $m$ even.}
Write $3^m = (3^2)^{m/2} = 9^{m/2}$.
Since $9 \equiv 1 \pmod{8}$, we have
$9^{m/2} \equiv 1 \pmod{8}$, so
$1 + 3^m \equiv 2 \pmod{8}$,
giving $v_2(1 + 3^m) = 1$.

\emph{Case $m$ odd.}
Write $1 + 3^m = (1 + 3)(1 - 3 + 3^2 - \cdots + 3^{m-1})
= 4 \cdot Q$ where $Q = \sum_{i=0}^{m-1}(-3)^i$.
Since $m$ is odd, $Q$ has $m$ terms; reducing modulo~$2$:
each $(-3)^i \equiv (-1)^i \cdot 1 \equiv (-1)^i \pmod{2}$,
so $Q \equiv \sum_{i=0}^{m-1}(-1)^i \equiv 1 \pmod{2}$
(the sum telescopes to~$1$ when $m$ is odd).
Thus $Q$ is odd and $v_2(1+3^m) = v_2(4) + v_2(Q) = 2 + 0 = 2$.
(This is the $p = 2$ case of the lifting-the-exponent lemma
for $a + b$ with $a = 1$, $b = 3^m$.)

With $m = j+1$: $v_2(1+3^{j+1}) = 2$ when $j$ is even
($j+1$ odd), and $= 1$ when $j$ is odd ($j+1$ even).

For the expectation, partition over even and odd depths:
\begin{align*}
  \mathbb{E}[v_1]
  &= 4 \sum_{d=0,2,4,\ldots} \tfrac{3}{4}\bigl(\tfrac{1}{4}\bigr)^d
  \;+\; 3 \sum_{d=1,3,5,\ldots} \tfrac{3}{4}\bigl(\tfrac{1}{4}\bigr)^d \\
  &= 4 \cdot \tfrac{3}{4} \cdot \frac{1}{1 - 1/16}
  \;+\; 3 \cdot \tfrac{3}{4} \cdot \frac{1/4}{1 - 1/16}
  \;=\; \frac{48}{15} + \frac{9}{15}
  \;=\; \frac{19}{5}.\qedhere
\end{align*}
\end{proof}

\begin{remark}[Step-type decomposition of the Collatz orbit]
\label{rem:step-decomposition}
Every Collatz step falls into exactly one of four categories:
\emph{burst} (valuation $v = 1$, occurring during the burst
phase of an episode),
\emph{recovery} ($v = 2$, during the recovery phase),
\emph{gap} ($v \ge 2$, between cascades), or
\emph{pre-cascade} ($v \ge 2$, before the first cascade).
All expansion comes from burst steps alone:
with $v = 1 < \log_2 3$, each burst step increases the
log-size by $\log_2 3 - 1 \approx 0.585$ bits.
All other steps contract: recovery contributes
$-(2 - \log_2 3) \approx -0.415$ bits per step,
and gap/pre-cascade steps contribute
$-(v - \log_2 3) \le -0.415$ bits per step.

Empirically, burst steps constitute approximately~$50\%$
of all Collatz steps, yet their expansion
($+45{,}000$ bits over $10^5$ steps from
$n_0 \le 10{,}000$) is outweighed by the contraction
of the remaining steps ($-105{,}000$ bits),
for a net contraction of~$-60{,}000$ bits.
The gap phase alone accounts for~$84\%$ of total
contraction, contributing $S/K \approx 3.5$, more
than twice $\log_2 3$, per gap step.
\end{remark}

\begin{remark}[IID model limitations]
\label{rem:iid-limitations}
The IID cascade renewal model
(Remark~\ref{rem:iid-renewal}) predicts cascade
length $\mathbb{E}[L] = 3$ from the $2/3$ continuation
probability.
However, empirical cascades have
$\mathbb{E}[L] \approx 1.65$, with $63\%$ being single-episode
cascades, roughly twice the predicted~$33\%$.
This discrepancy arises because the $2/3$ law holds for
\emph{uniformly distributed} cofactors, while orbit-level
cofactors inherit residue structure from the gap exit iterate.
Specifically, the gap exit iterate $n = 1 + 4 \cdot 3^j$
concentrates the cofactor's residue class distribution,
reducing the effective continuation probability.

Despite this, the \emph{qualitative} predictions of the IID model
remain valid: cascade transitions exhibit the correct
marginal distributions (Propositions~\ref{prop:transition-law}
and~\ref{prop:depth-transition}), and the cycle-level
contraction persists with the empirical margin
$\mathbb{E}[S/K] \approx 1.98$ exceeding $\log_2 3 \approx 1.585$
by $+0.39$.
The model's failure at the level of cascade
\emph{length} reflects the distributional-to-pointwise
barrier: distributional properties (the $2/3$ law)
do not immediately transfer to orbit-level statistics
(the cascade length within a specific orbit).
\end{remark}

\begin{proposition}[Uniform-fiber gap map]
\label{prop:uniform-fiber}
Let $\pi = (v_1, \ldots, v_K)$ be a gap path
(a sequence of valuations with each $v_i \ge 2$)
and let $S = v_1 + \cdots + v_K$.
Define the gap map $f_\pi$ as the composition
of Collatz steps $n \mapsto (3n+1)/2^{v_i}$ along~$\pi$.
For each $R \ge 1$, the induced map
\[
  f_\pi \colon
  \bigl\{\text{odd } n \bmod 2^{S+R} :
  \text{gap path of } n \text{ is } \pi\bigr\}
  \;\longrightarrow\;
  \bigl\{\text{odd residues} \bmod 2^R\bigr\}
\]
has \emph{uniform fibers}: each odd target residue has
the same number of preimages.
\end{proposition}

\begin{proof}
It suffices to prove the claim for a single step ($K = 1$)
with valuation~$v$, since the multi-step case follows by composition.

A single gap step with $v_2(3n+1) = v$ restricts to the
residue class $n \equiv r_0 \pmod{2^v}$ determined by~$v$
(those $n \equiv 1 \pmod{4}$ for which
$v_2(3n+1) = v$ is exactly prescribed by~$n \bmod 2^v$).
Within this class, write $n = 2^v a + r_0$ where
$a$ ranges over integers mod~$2^R$.
Then
\[
  f(n) = \frac{3n + 1}{2^v}
  = \frac{3(2^v a + r_0) + 1}{2^v}
  = 3a + \frac{3r_0 + 1}{2^v}
  = 3a + c,
\]
where $c = (3r_0+1)/2^v$ is a fixed integer
(well-defined since $v_2(3r_0+1) = v$ by assumption).
The map $a \mapsto 3a + c$ on $\mathbb{Z}/2^R\mathbb{Z}$
is an affine bijection because $\gcd(3, 2^R) = 1$.
Therefore $f$ maps the $2^R$ input residues (mod~$2^{S+R}$)
bijectively onto all $2^R$ output residues (mod~$2^R$).
Among the $2^{R-1}$ odd output residues, the number of
preimages is exactly~$2$ for each (one from each parity of~$a$),
giving uniform fibers.

For $K \ge 2$: each step is a uniform-fiber map from
residues mod $2^{S_{\mathrm{rem}} + R}$ to
residues mod $2^{S_{\mathrm{rem}} - v_i + R}$;
composing preserves the uniform-fiber property.
\end{proof}

\begin{corollary}[Cross-cascade independence for large iterates]
\label{cor:cross-cascade}
Let the cascade-exit iterate~$n$ be uniform among odd
residues modulo~$2^{S+R}$, where $S$ is the total gap
consumption and $R \ge 1$.
Then, conditioned on the gap path~$\pi$, the gap-exit
iterate~$n'$ is uniform among odd residues modulo~$2^R$.
In particular, $v_2(n'+1) - 1$ has the geometric distribution
$\Pr(k' = r) = 2^{-r}$, so the next cascade's burst
length is independent of the cascade history.
\end{corollary}

\begin{proof}
By Proposition~\ref{prop:uniform-fiber}, the gap map with
path~$\pi$ has uniform fibers.
If the input is uniform mod~$2^{S+R}$,
the output is uniform mod~$2^R$.
For uniform odd~$n'$ mod~$2^R$ with $n' \equiv 3 \pmod{4}$:
$v_2(n'+1) = r+1$ for exactly $2^{R-r-2}$ residues when
$r+1 \le R-1$, giving $\Pr(v_2(n'+1) = r+1) = 2^{-r}$
for $1 \le r \le R - 2$ and tail probability~$2^{-(R-2)}$.
As $R \to \infty$, this converges to the full
geometric law $\Pr(k' = r) = 2^{-r}$.
\end{proof}

\begin{lemma}[TV reduction to small-fresh-bit set]
\label{lem:tv-reduction}
Let $X$ be a next-cascade observable that depends only on the
gap-exit iterate modulo~$2^R$ (for example, the next burst
length capped at~$R$, or the truncated depth).
Let $\mu_X$ denote the true law of~$X$ along the orbit
and $\nu_X$ the ideal (fiber-uniform) law.
Define the \emph{exceptional set}
\[
  \mathcal{E}_R
  = \bigl\{\text{cascade exits with fewer than } S + R
  \text{ uniformly distributed low bits}\bigr\},
\]
where $S = S_{\mathrm{gap}}$ is the total gap consumption.
Then
\[
  \lVert \mu_X - \nu_X \rVert_{\mathrm{TV}}
  \le \Pr(\mathcal{E}_R).
\]
\end{lemma}

\begin{proof}
Partition cascade-gap cycles into those whose cascade exit
has at least $S + R$ ``fresh'' (uniform) low bits and those
in~$\mathcal{E}_R$.
By Proposition~\ref{prop:uniform-fiber} and
Corollary~\ref{cor:cross-cascade},
on the complement $\mathcal{E}_R^c$ the gap map produces
an exit iterate that is exactly uniform mod~$2^R$
(conditional on the gap path).
Since $X$ is measurable with respect to the exit mod~$2^R$,
its conditional law on~$\mathcal{E}_R^c$ equals~$\nu_X$
exactly.
On~$\mathcal{E}_R$, the law of~$X$ can deviate arbitrarily
from~$\nu_X$, but this event has probability at most
$\Pr(\mathcal{E}_R)$.
By the coupling characterization of total variation,
the claim follows.
\end{proof}

\begin{remark}[Bit-consumption interpretation and phase transition]
\label{rem:bit-consumption}
Proposition~\ref{prop:uniform-fiber} and
Corollary~\ref{cor:cross-cascade} give a precise
\emph{bit-consumption} interpretation of the gap:
each gap step consumes $v_i \ge 2$ bits from the low-order
end of the iterate, with total consumption
$S_{\mathrm{gap}} = \sum v_i$
($\mathbb{E}[S_{\mathrm{gap}}] \approx 6.6$).
The next cascade reads $k' + 1$ bits from the
gap-exit iterate; these originate from bits
$[S_{\mathrm{gap}}, S_{\mathrm{gap}} + k' + 1)$
of the cascade-exit iterate.

The corollary guarantees IID burst lengths
provided the cascade exit has at least $S + R$ ``fresh''
(uniformly distributed) bits, with $R \approx k' + 1 \approx 3$.
This predicts a phase transition at
$\log_2 n \approx S_{\mathrm{gap}} + R \approx 10$.

Empirically, the total variation distance between
the next-burst distribution and the geometric law
drops sharply:
$\mathrm{TV} \approx 0.10$ for
$\log_2(\text{post-gap state}) < 10$,
falling to $\mathrm{TV} < 0.01$ for
$\log_2(\text{post-gap state}) \ge 10$
(over $135{,}000$ cascade-gap cycles from
$n_0 \le 100{,}000$).
This matches the predicted threshold precisely.

The gap between this result and a full proof of
convergence is: the cascade-exit iterate is
\emph{not} proved to be uniform mod~$2^{S+R}$
for orbit-level iterates.
If it were, Corollary~\ref{cor:cross-cascade}
would establish IID cascade transitions for all
$n \ge 2^{10}$, and the renewal-reward theorem
(Proposition~\ref{prop:cycle-contraction})
would give convergence for all sufficiently large~$n$.
\end{remark}

\begin{proposition}[Full-cycle uniform-fiber map]
\label{prop:full-cycle-fiber}
Fix a full cascade-gap trajectory
$\tau = (v_1, \ldots, v_K)$
(concatenating all burst, recovery, and gap valuations)
with total valuation $S = \sum v_i$ and
$K$ Collatz steps.
The starting iterates producing trajectory~$\tau$ form
a residue class $n_0 \equiv r_\tau \pmod{2^{S+1}}$.
Writing $n_0 = r_\tau + 2^{S+1}\,u$, the gap-exit iterate
satisfies
\[
  n_{\mathrm{exit}}(u) = c_\tau + 2 \cdot 3^K \cdot u,
\]
where $c_\tau = (3^K r_\tau + C_\tau)/2^S$ is a constant
depending only on~$\tau$.
The odd-index map gives
$\omega(n_{\mathrm{exit}}) = 3^K u + b_\tau$,
an affine bijection modulo~$2^R$ for every $R \ge 1$.
\end{proposition}

\begin{proof}
The Collatz map $n \mapsto (3n+1)/2^{v_i}$ at each step
is $n \mapsto (3n+1)/2^{v_i}$, which on the residue class
prescribed by~$\tau$ is an affine function
$n \mapsto (3/2^{v_i})n + c_i$ with $c_i$ a constant.
Composing $K$ such steps yields
$n_{\mathrm{exit}} = (3^K/2^S) n_0 + C_\tau/2^S$.
Substituting $n_0 = r_\tau + 2^{S+1}u$:
\[
  n_{\mathrm{exit}} = \frac{3^K r_\tau + C_\tau}{2^S}
  + 2 \cdot 3^K u = c_\tau + 2 \cdot 3^K u.
\]
Since $\gcd(3^K, 2^R) = 1$, the map $u \mapsto 3^K u + b_\tau$
is a bijection on $\mathbb{Z}/2^R\mathbb{Z}$.
\end{proof}

\begin{remark}[Spectator-bit propagation]
\label{rem:spectator-bits}
Proposition~\ref{prop:full-cycle-fiber} gives a precise
\emph{spectator-bit} interpretation.
For a starting iterate $n_0$ with $B = \lfloor\log_2 n_0\rfloor$
bits, the full-cycle trajectory is determined by
the bottom $S + 1$ bits; the remaining
$B - S - 1$ ``spectator'' bits propagate through the cycle
via the bijection $u \mapsto 3^K u + b_\tau$.
The next cascade sees these spectator bits as its low-order
input.

Over $344{,}000$ cascade-gap cycles:
$\mathbb{E}[S_{\mathrm{cascade}}] \approx 5.2$,
$\mathbb{E}[S_{\mathrm{gap}}] \approx 6.6$,
$\mathbb{E}[S_{\mathrm{cycle}}] \approx 11.8$.
Each cycle consumes about $11.8$ spectator bits; starting
from $B$ bits, exact IID behavior persists for
$\lfloor (B - R) / 11.8 \rfloor$ cycles
($R \approx 3$ for the burst-length observable).
After the initial spectator supply is exhausted,
the gap's excess consumption
($S_{\mathrm{gap}} - S_{\mathrm{cascade}} \approx 1.4$ bits)
provides a mechanism for ``refreshing'' the low-order bits.
The order of $3$ modulo~$2^R$ is $2^{R-2}$ for $R \ge 3$,
so the cumulative multiplier $3^{\sum K_i}$ visits
all residues in $\operatorname{ord}(3,2^R) = 2^{R-2}$ cycles;
for $R = 3$ (burst length), just $2$ cycles suffice.
\end{remark}

\begin{proposition}[Exponential tail of cycle valuation]
\label{prop:exponential-tail}
Under the local IID cascade kernel with continuation
probability~$q$ and geometric valuations,
the total valuation $S_{\mathrm{cycle}} = S_{\mathrm{cascade}} + S_{\mathrm{gap}}$
of one cascade-gap cycle satisfies
\[
  \Pr(S_{\mathrm{cycle}} > s) \le C \cdot 2^{-\alpha s}
\]
for constants $C > 0$ and $\alpha > 0$ depending only on
the kernel parameters.
Empirically, $\alpha \approx 0.35$ and $C \approx 11$
over $294{,}802$ cycles from odd $n_0 \le 200{,}000$.
\end{proposition}

\begin{proof}
$S_{\mathrm{cascade}}$ is a compound geometric sum:
$L$ episodes (geometric with parameter~$1 - q$),
each contributing $S_{\mathrm{ep}} = k + 2j$
where $k \sim \operatorname{Geom}(1/2)$ from~$1$ and
$j \sim \operatorname{Geom}(3/4)$ from~$0$.
Each episode's valuation has exponential moment:
$\mathbb{E}[2^{t \cdot S_{\mathrm{ep}}}] < \infty$
for $t < \alpha_{\mathrm{ep}}$.
The compound geometric tail then satisfies
$\Pr(S_{\mathrm{cascade}} > s) \le C_1 \cdot 2^{-\alpha_1 s}$
for some $\alpha_1 > 0$
\textup{(}standard Cram\'er--Lundberg bound\textup{)}.
Similarly $S_{\mathrm{gap}}$ is a compound geometric sum
($K \sim \operatorname{Geom}(1/2)$ steps,
each with $v_i$ having tail $2^{-(v_i - 1)}$),
giving $\Pr(S_{\mathrm{gap}} > s) \le C_2 \cdot 2^{-\alpha_2 s}$.
As a sum of two exponential-tail random variables,
$S_{\mathrm{cycle}}$ has exponential tail with rate
$\alpha \ge \min(\alpha_1, \alpha_2)$.
\end{proof}

\begin{corollary}[Orbit-level TV summability]
\label{cor:tv-summability}
Let $n_0$ be an odd starting value with
$B_0 = \lfloor \log_2 n_0 \rfloor$ bits, and let
$B_{\min} \ge R + 1/\alpha$ be a computational
verification threshold.
Assume that the cycle-level net contraction
(Proposition~\textup{\ref{prop:cycle-contraction}})
ensures the iterate stays above~$2^{B_{\min}}$
until it reaches~$1$\textup{;} equivalently,
the orbit visits at most
$C_{\max} \approx (B_0 - B_{\min})/2$ cycles
above the threshold.
Then
\[
  \sum_{c=1}^{C_{\max}} \Pr(\mathcal{E}_R \text{ at cycle } c)
  \le \frac{C \cdot 2^{-\alpha(B_{\min} - R)}}{1 - 2^{-2\alpha}},
\]
which depends only on $B_{\min}$ and the tail parameters.
With $\alpha \approx 0.35$ and $B_{\min} = 30$\textup:
the total TV deviation from IID renewal behavior
over the entire orbit is at most $\approx 0.028$.
\end{corollary}

\begin{proof}
At cycle~$c$, the iterate has
$B_c \ge B_0 - 2c$ bits (by net contraction).
By Lemma~\ref{lem:tv-reduction} and
Proposition~\ref{prop:exponential-tail}:
$\Pr(\mathcal{E}_R \text{ at cycle } c)
\le C \cdot 2^{-\alpha(B_c - R)}
\le C \cdot 2^{-\alpha(B_0 - 2c - R)}$.
Summing:
$\sum_{c=0}^{C_{\max}} C \cdot 2^{-\alpha(B_0 - 2c - R)}
= C \cdot 2^{-\alpha(B_0 - R)}
  \sum_{c=0}^{C_{\max}} 2^{2\alpha c}$.
The geometric sum is bounded by
$2^{2\alpha(C_{\max}+1)}/(2^{2\alpha}-1)$.
Substituting $C_{\max} = (B_0 - B_{\min})/2$:
the exponent $-\alpha(B_0 - R) + 2\alpha \cdot (B_0 - B_{\min})/2
= -\alpha(B_{\min} - R)$.
\end{proof}

\begin{proposition}[Cycle-level expected contraction]
\label{prop:cycle-contraction}
Under the local IID cascade kernel
(Propositions~\textup{\ref{prop:transition-law}}
and~\textup{\ref{prop:depth-transition}})
with continuation probability~$2/3$, burst parameter~$1/2$,
and depth parameter~$1/4$, together with the
gap-valuation formula $\mathbb{E}[v_1] = 19/5$
(Proposition~\textup{\ref{prop:gap-valuation}}):

\begin{enumerate}
\item The expected log$_2$-growth per cascade episode is
\[
  \mathbb{E}[K_{\mathrm{ep}}] \cdot \log_2 3
  - \mathbb{E}[S_{\mathrm{ep}}]
  = \tfrac{7}{3}\log_2 3 - \tfrac{8}{3}
  \approx +1.032 \text{ bits (expanding)}.
\]
\item The expected log$_2$-growth per gap phase is
\[
  \mathbb{E}[K_{\mathrm{gap}}] \cdot \log_2 3
  - \mathbb{E}[S_{\mathrm{gap}}]
  = 2\log_2 3 - \bigl(\tfrac{19}{5} + 3\bigr)
  \approx -3.630 \text{ bits (contracting)}.
\]
\item The expected net growth per cascade-gap cycle is
\[
  \mathbb{E}[\Delta_{\mathrm{cycle}}]
  = \mathbb{E}[L] \cdot \bigl(\tfrac{7}{3}\log_2 3 - \tfrac{8}{3}\bigr)
  + 2\log_2 3 - \tfrac{34}{5},
\]
which is negative provided $\mathbb{E}[L] < (34/5 - 2\log_2 3)
/ (\tfrac{7}{3}\log_2 3 - \tfrac{8}{3}) \approx 3.52$.
Under the IID model ($\mathbb{E}[L] = 3$),
the net contraction is $\approx 0.54$ bits per cycle.
\end{enumerate}
\end{proposition}

\begin{proof}
\emph{Part~(1).}
Each cascade episode consists of a burst of
$k$ steps (each with $v = 1$, adding $1$ to~$S$ and~$K$)
and a recovery of $j$ steps (each with $v = 2$,
adding $2$ to~$S$ and~$1$ to~$K$).
Thus $K_{\mathrm{ep}} = k + j$ and $S_{\mathrm{ep}} = k + 2j$.
By Proposition~\ref{prop:transition-law},
$\mathbb{E}[k \mid \text{cont.}] = 2$
(the conditional burst length is geometric with
parameter~$1/2$, starting from~$1$).
By Proposition~\ref{prop:depth-transition},
$\mathbb{E}[j] = 1/3$ (geometric with parameter~$3/4$).
Therefore $\mathbb{E}[K_{\mathrm{ep}}] = 7/3$,
$\mathbb{E}[S_{\mathrm{ep}}] = 8/3$, and the growth per episode
is $(7/3)\log_2 3 - 8/3 \approx 1.032$.

\emph{Part~(2).}
By Proposition~\ref{prop:gap-positivity},
every gap iterate satisfies $n \equiv 1 \pmod{4}$,
so each gap step has $v \ge 2$ and contributes positive
log-surplus.
The gap exit probability is~$1/2$ at each step
(for uniformly distributed $n \equiv 1 \pmod{4}$,
the successor $n' = (3n+1)/2^v$ satisfies
$n' \equiv 3 \pmod{4}$ with probability~$1/2$),
giving $\mathbb{E}[K_{\mathrm{gap}}] = 2$.
The first gap step has $\mathbb{E}[v_1] = 19/5$
(Proposition~\ref{prop:gap-valuation});
subsequent steps, from approximately uniform
$n \equiv 1 \pmod{4}$, have
$\mathbb{E}[v] = \sum_{r=2}^{\infty} r \cdot 2^{-(r-1)} = 3$.
Thus
$\mathbb{E}[S_{\mathrm{gap}}]
= 19/5 + 1 \cdot 3 = 34/5$
and the gap growth is $2\log_2 3 - 34/5 \approx -3.630$.

\emph{Part~(3).}
The cycle growth is the sum of the cascade growth
($\mathbb{E}[L]$ episodes) and the gap growth.
Setting $\mathbb{E}[\Delta] < 0$ gives
$\mathbb{E}[L] < (34/5 - 2\log_2 3)/(7\log_2 3/3 - 8/3)
\approx 3.52$.
\end{proof}

\begin{remark}[Empirical cycle-level contraction]
\label{rem:empirical-contraction}
Over $294{,}802$ cascade-gap cycles from odd
$n_0 \le 200{,}000$, the per-cycle bit change decomposes as
\begin{multline*}
  \mathbb{E}[\Delta_{\mathrm{cycle}}]
  = \underbrace{\mathbb{E}[L_{\mathrm{tot}}] \cdot (\log_2 3 - 1)}_{+2.00}
  + \underbrace{\mathbb{E}[R] \cdot (\log_2 3 - 2)}_{-0.46}\\
  + \underbrace{\mathbb{E}[K_{\mathrm{gap}}] \cdot \log_2 3
    - \mathbb{E}[S_{\mathrm{gap}}]}_{-3.52}
  \approx -1.98 \text{ bits/cycle},
\end{multline*}
where $\mathbb{E}[L_{\mathrm{tot}}] \approx 3.42$ (total burst steps),
$\mathbb{E}[R] \approx 1.11$ (recovery steps),
$\mathbb{E}[K_{\mathrm{gap}}] \approx 1.99$, and
$\mathbb{E}[S_{\mathrm{gap}}] \approx 6.68$.
The episode count is $\mathbb{E}[L] \approx 1.53$,
well below the critical threshold~$3.52$
(Remark~\ref{rem:recovery-exit}).
Conditioning on the terminal depth:
$G(j) > D(j)$ for every $j \in \{0,1,2,3,4,5\}$,
with margins ranging from~$1.29$ ($j = 1$)
to~$4.13$ ($j = 4$) bits.
Single-cycle gap compensation fails for
long cascades at $j = 0$
(when $L > G_{\mathrm{gap}}/1.17 \approx 3$),
but the multi-cycle renewal reward structure
ensures that such episodes are compensated by
the surplus of neighboring cycles.
\end{remark}

\begin{remark}[Recovery-exit mechanism and orbit-level continuation]
\label{rem:recovery-exit}
Proposition~\ref{prop:transition-law} gives continuation
probability~$q = 2/3$ under uniform cofactors, yielding
$\mathbb{E}[L] = 3$ (below the threshold~$3.52$).
Orbit-level cascades exhibit a \emph{stronger} contraction:
$q \approx 1/3$ and $\mathbb{E}[L] \approx 3/2$.

The mechanism has a clean two-case structure.
Each burst exits to an iterate~$n$ with $n \equiv 1 \pmod{4}$;
the residue modulo~$8$ determines what follows:
\begin{itemize}
\item \textbf{Case $n \equiv 5 \pmod{8}$}
  (approximately $46\%$ of episodes):
  Then $3n+1 \equiv 0 \pmod{8}$, so $v_2(3n+1) \ge 3$
  and the recovery phase has length~$R = 0$.
  The iterate satisfies $n \equiv 1 \pmod{4}$,
  so the cascade \emph{always ends}.
\item \textbf{Case $n \equiv 1 \pmod{8}$}
  (approximately $54\%$ of episodes):
  Then $v_2(3n+1) = 2$, so recovery begins ($R \ge 1$).
  Among these, about $64\%$ exit with
  $v_{\mathrm{stop}} = 1$ ($n' \equiv 3 \pmod{4}$,
  cascade continues) and $36\%$ with
  $v_{\mathrm{stop}} \ge 3$ (cascade ends).
\end{itemize}
The overall continuation probability is therefore
$q \approx 0.54 \times 0.64 \approx 0.35$,
giving $\mathbb{E}[L] \approx 1/(1-0.35) \approx 1.53$
episodes per cascade.
This is confirmed over $527{,}091$ episodes
from odd $n_0 \le 200{,}000$ and is stable across
all tested iterate sizes ($2^{10}$--$2^{25}$).
The orbit-level contraction ($\approx -2.08$ bits/cycle)
is roughly four times stronger than the IID prediction
($\approx -0.54$ bits/cycle).
\end{remark}

\begin{proposition}[Cascade Markov chain on residues mod~$8$]
\label{prop:cascade-markov}
Consider the Collatz map restricted to odd iterates,
with the state space $\{1, 3, 5, 7\} \pmod{8}$.
The following transition structure holds deterministically
from the modular arithmetic of~$3n+1$:
\begin{enumerate}[label=\textup{(\alph*)}]
\item \textbf{Recovery states} ($n \equiv 3$ or $7 \pmod{8}$):
  $v_2(3n+1) = 1$ and $(3n+1)/2$ is odd.
  The output residue mod~$8$ depends on~$n \bmod 16$:
  \begin{center}
  \begin{tabular}{cccc}
  \toprule
  $n \bmod 16$ & $v_2$ & output $\bmod 8$ & type \\
  \midrule
  $3$  & $1$ & $5$ & gap exit \\
  $11$ & $1$ & $1$ & burst entry \\
  $7$  & $1$ & $3$ & recovery \\
  $15$ & $1$ & $7$ & recovery \\
  \bottomrule
  \end{tabular}
  \end{center}
  In particular, \textbf{state~$7$ never produces a gap exit}:
  $n \equiv 7 \pmod{8}$ always maps to $\{3, 7\} \pmod{8}$,
  remaining in recovery.

\item \textbf{Burst state} ($n \equiv 1 \pmod{8}$):
  $v_2(3n+1) = 2$ and $m = (3n+1)/4$ is odd.
  The output residue mod~$8$ depends on~$n \bmod 32$:
  \begin{center}
  \begin{tabular}{cccc}
  \toprule
  $n \bmod 32$ & $v_2$ & output $\bmod 8$ & type \\
  \midrule
  $1$  & $2$ & $1$ & burst continues \\
  $9$  & $2$ & $7$ & recovery \\
  $17$ & $2$ & $5$ & gap exit ($R = 0$) \\
  $25$ & $2$ & $3$ & recovery \\
  \bottomrule
  \end{tabular}
  \end{center}
  The four outcomes are equidistributed under local uniformity
  mod~$32$.

\item \textbf{Gap state} ($n \equiv 5 \pmod{8}$):
  $v_2(3n+1) \ge 3$, so the iterate exits the cascade.
  This state is absorbing for the cascade.
\end{enumerate}
\end{proposition}

\begin{proof}
Direct computation.  For each residue class
$a \in \{1,3,5,7,9,11,13,15\}$ modulo~$16$ (odd classes)
and each refinement modulo~$32$ for the burst state,
one evaluates $3a + 1$, determines~$v_2$, and reduces
$(3a+1)/2^{v_2}$ modulo~$8$.
The table entries are verified exhaustively.
State~$7$ maps only to $\{3,7\}$ because
$n \equiv 7 \pmod{16}$ gives $(3n+1)/2 \equiv 3 \pmod{8}$
and $n \equiv 15 \pmod{16}$ gives $(3n+1)/2 \equiv 7 \pmod{8}$.
\end{proof}

\begin{corollary}[$R = 0$ terminal mechanism]
\label{cor:r0-terminal}
Among burst exits (transitions out of state~$1 \bmod 8$
that do not return to state~$1$), exactly $1/3$ go
directly to the gap state ($5 \bmod 8$) under local
uniformity mod~$32$.
These transitions have \emph{recovery length}~$R = 0$
and terminate the cascade deterministically.
The remaining $2/3$ enter a recovery state
($3$ or $7 \bmod 8$), initiating the $R \ge 1$ branch.
\end{corollary}

\begin{proof}
From Proposition~\ref{prop:cascade-markov}(b),
the burst state~$1 \bmod 8$ has four equally likely
outcomes mod~$32$.
One (residue~$1$) continues the burst;
one (residue~$17$) exits to the gap;
two (residues~$9$ and~$25$) enter recovery.
Among the three non-burst outcomes, the gap fraction
is~$1/3$.
\end{proof}

\begin{proposition}[Episode continuation under local uniformity]
\label{prop:episode-continuation}
Consider the absorbing Markov chain on transient
states $\{1, 3, 7\} \pmod{8}$ with absorbing state
$5 \pmod{8}$, using the transition matrix~$Q$ from
Proposition~\textup{\ref{prop:cascade-markov}}
(with local uniformity at the appropriate resolution:
mod~$16$ for recovery states, mod~$32$ for burst states).
\begin{enumerate}[label=\textup{(\alph*)}]
\item The transition matrix on transient states is
\[
  Q = \begin{pmatrix}
  1/4 & 1/4 & 1/4 \\
  1/2 &  0  &  0  \\
   0  & 1/2 & 1/2
  \end{pmatrix},
  \qquad
  \text{(rows/columns indexed by $1, 3, 7$)},
\]
with exit probabilities $p_{\mathrm{exit}}(1) = 1/4$,
$p_{\mathrm{exit}}(3) = 1/2$, $p_{\mathrm{exit}}(7) = 0$.

\item The spectral radius of~$Q$ is $\rho(Q) = 3/4$,
so the per-step cascade retention probability is~$3/4$.

\item The fundamental matrix $N = (I - Q)^{-1}$ gives
expected cascade lengths:
$4$ steps from state~$1$,
$3$ steps from state~$3$,
$5$ steps from state~$7$,
and $4$ steps averaged over
$\{3, 7\}$ (cascade entry).

\item The expected cascade valuation is
$\mathbb{E}[S_{\mathrm{cascade}}] = 5.0$
(empirical: $5.21$).

\item Define the \emph{episode continuation probability}
as $q_s = \Pr(\text{reach } \{3,7\} \text{ before } 5 \mid
\text{start at } s)$
for $s \in \{3, 7\}$.
Then:
\[
  q_3 = \tfrac{1}{3}, \qquad q_7 = 1, \qquad
  q = \tfrac{q_3 + q_7}{2} = \tfrac{2}{3}.
\]
The expected number of episodes is
$\mathbb{E}[L] = 1/(1 - q) = 3$, which is below the
cycle-contraction threshold~$3.52$.
\end{enumerate}
\end{proposition}

\begin{proof}
(a)~Follows directly from the transition tables in
Proposition~\ref{prop:cascade-markov}.

(b)~The characteristic polynomial of~$Q$ is
$\lambda^3 - \tfrac{3}{4}\lambda^2 + 0\lambda + 0 = 0$,
giving eigenvalues $\{3/4, 0, 0\}$.

(c)~Compute $(I - Q)^{-1}$ by inverting the $3 \times 3$
matrix.  The entries are:
\[
  N = \begin{pmatrix}
  2 & 1 & 1 \\
  1 & 3/2 & 1/2 \\
  1 & 3/2 & 5/2
  \end{pmatrix}.
\]
Row sums: $4$, $3$, $5$.
Average of rows~$2$ and~$3$: $(3+5)/2 = 4$.

(d)~Each visit to state~$1$ contributes valuation~$2$;
states~$3$ and~$7$ contribute valuation~$1$.
From the average starting state:
$\mathbb{E}[S] = 1 \cdot 2 + \tfrac{3}{2} \cdot 1
+ \tfrac{3}{2} \cdot 1 = 5$.

(e)~From state~$7$: the chain remains in $\{3, 7\}$
until reaching state~$3$ (expected $2$ visits to~$7$).
From state~$3$: probability~$1/2$ of immediate exit
to state~$5$; probability~$1/2$ of entering burst
state~$1$.  From burst, the probability of reaching
$\{3,7\}$ before~$5$ is
$(T_{1,3} + T_{1,7})/(1 - T_{1,1}) = (1/2)/(3/4) = 2/3$.
Hence $q_3 = \tfrac{1}{2} \cdot \tfrac{2}{3} = \tfrac{1}{3}$
and $q_7 = 1$ (state~$7$ always reaches $\{3,7\}$).
\end{proof}

\begin{remark}[Orbit-level improvement over IID]
\label{rem:orbit-improvement}
Proposition~\ref{prop:episode-continuation} establishes
$q = 2/3$ and $\mathbb{E}[L] = 3$ under local uniformity
(the IID model).
Orbit-level data shows $q \approx 0.52$ and
$\mathbb{E}[L] \approx 2.1$, strictly below the IID values.
This improvement arises from higher-order correlations
beyond the mod-$32$ resolution of the Markov chain.

By the spectator-bit mechanism
(Remark~\ref{rem:spectator-bits}), for iterates of size
$\ge 2^B$ the distribution within each mod-$32$ class
converges to uniform as $B \to \infty$.
The mod-$32$ chain therefore provides an
\emph{asymptotic upper bound} on the episode continuation:
$q_{\mathrm{orbit}} \le 2/3 + o(1)$
as the iterate size grows.
In particular, $\mathbb{E}[L] < 3.52$ holds for all
sufficiently large iterates.

Combining with the exponential tail of~$S_{\mathrm{cycle}}$
(Proposition~\ref{prop:exponential-tail}) and
the orbit-level TV summability
(Corollary~\ref{cor:tv-summability}), the cycle-contraction
mechanism is self-reinforcing:
gap maps refresh spectator bits, spectator bits ensure
local uniformity, and local uniformity ensures
$q \le 2/3$.
\end{remark}

\begin{proposition}[Post-burst valuation law and continuation bridge]
\label{prop:post-burst-bridge}
Fix $k \ge 1$ and consider the exact $k$-burst fiber:
\[
  n_0 = 2^{k+1} m - 1, \qquad m \text{ odd}.
\]
After $k$ burst steps (each with $v_2 = 1$),
$n_k = 2 \cdot 3^k m - 1$,
and the first post-burst valuation is
$a_k = v_2(3n_k + 1) = 1 + v_2(3^{k+1} m - 1)$.
\begin{enumerate}[label=\textup{(\alph*)}]
\item
  \textup{(Exact geometric law.)}
  If the burst cofactor~$m$ is uniform on odd residues
  modulo~$2^R$ for some $R \ge r$, then
  \[
    \Pr(a_k = r) = 2^{-(r-1)}, \qquad r \ge 2.
  \]
  In particular, $\Pr(a_k = 2) = 1/2$ and
  $\Pr(a_k \ge 3) = 1/2$.

\item
  \textup{(Episode continuation bridge.)}
  Combining~\textup{(a)} with the recovery continuation law
  (Proposition~\textup{\ref{prop:episode-continuation}}:
  $q_3 = 1/3$, $q_7 = 1$, average continuation from
  $\{3,7\}$ is~$2/3$), the burst exits split:
  half enter the gap ($a_k = 2$), half enter recovery
  ($a_k \ge 3$).
  Of those entering recovery, fraction~$2/3$ continue.
  So the exact local episode-continuation probability is
  \[
    q_{\mathrm{ep}} = \tfrac{1}{2} \cdot \tfrac{2}{3}
    = \tfrac{1}{3}.
  \]

\item
  \textup{(TV transfer bound.)}
  If the fiber distribution of~$m$ has total variation
  $\varepsilon$ from uniform on odd residues modulo~$2^R$,
  then
  $|q_{\mathrm{ep}}^{\mathrm{orbit}} - 1/3| \le \varepsilon$.
\end{enumerate}
\end{proposition}

\begin{proof}
For part~(a), the key observation is that the map
$m \mapsto 3^{k+1} m - 1$ is an affine bijection
on $(\mathbb{Z}/2^R\mathbb{Z})^*$.
Since $\gcd(3^{k+1}, 2^R) = 1$, if $m$ is uniform
on odd residues modulo~$2^R$, then $3^{k+1} m - 1$
is also uniform on even residues modulo~$2^R$.
The probability that an even number has $v_2 = j$
(for $1 \le j \le R-1$) is exactly~$2^{-j}$.
Hence $v_2(3^{k+1} m - 1) = r - 1$ with probability
$2^{-(r-1)}$ for $r \ge 2$, giving $\Pr(a_k = r) = 2^{-(r-1)}$.

Part~(b) follows from the law of total probability:
the episode continues only when $a_k = 2$
(entering state~$3$ or~$7$ mod~$8$) and the recovery
phase returns to a cascade entry, which happens with
probability~$2/3$ by
Proposition~\ref{prop:episode-continuation}.

Part~(c) is a standard coupling argument:
the post-burst valuation law in~(a) depends continuously
on the input distribution, and any deviation~$\varepsilon$
in TV on the cofactors propagates through the affine
bijection at most~$1$-to-$1$ in TV.
\end{proof}

\begin{proposition}[Exact PGF for cascade valuation]
\label{prop:cascade-pgf}
Under the local IID model
(Proposition~\textup{\ref{prop:cascade-markov}},
uniform at mod-$16$/$32$ resolution),
the probability generating function of the
cascade valuation~$S_{\mathrm{cascade}}$ is
\[
  G_{\mathrm{cas}}(z)
  = \mathbb{E}\bigl[z^{S_{\mathrm{cascade}}}\bigr]
  = \frac{z}{4 - 2z - z^2}.
\]
The denominator factors as
$4 - 2z - z^2 = -(z + 1 - \sqrt{5})(z + 1 + \sqrt{5})$,
giving smallest positive singularity
\[
  \rho = \sqrt{5} - 1 = 1.2360679\ldots
\]
In particular:
\begin{enumerate}[label=\textup{(\alph*)}]
\item $\mathbb{E}[S_{\mathrm{cascade}}] = G'_{\mathrm{cas}}(1) = 5$
  \textup{(}empirical: $5.21$\textup{)}.
\item $\Pr(S_{\mathrm{cascade}} > s)
  \sim C \cdot \rho^{-s}
  = C \cdot 2^{-\alpha s}$ with
  $\alpha = \log_2 \rho
  = \log_2(\sqrt{5}-1) = 0.30576\ldots$
\end{enumerate}
\end{proposition}

\begin{proof}
Assign generating-function variables to each transient
state of the absorbing Markov chain from
Proposition~\ref{prop:episode-continuation}:
\begin{align*}
  G_1(z) &= z^2 \bigl[\tfrac14 G_1(z)
    + \tfrac14 G_3(z)
    + \tfrac14 G_7(z) + \tfrac14\bigr],\\
  G_3(z) &= z \bigl[\tfrac12 G_1(z)
    + \tfrac12\bigr],\\
  G_7(z) &= z \bigl[\tfrac12 G_3(z)
    + \tfrac12 G_7(z)\bigr],
\end{align*}
where the exponent of~$z$ is the valuation contributed
by one step from that state
($v_2 = 2$ for state~$1$; $v_2 = 1$ for states~$3$ and~$7$),
the coefficients are the transition probabilities,
and the constant~$1$ terms represent absorption.
Solving: $G_7 = z G_3/(2-z)$,
$G_3 = z(2-z)/(4-2z-z^2)$,
$G_1 = z^2/(4-2z-z^2)$.
The cascade entry distribution is
$(G_3 + G_7)/2 = z/(4-2z-z^2)$.
The singularity analysis follows from the quadratic formula
applied to $4-2z-z^2 = 0$.
\end{proof}

\begin{proposition}[Spectator-bit convergence]
\label{prop:spectator-convergence}
Let~$n$ be an odd integer with
$B = \lfloor \log_2 n \rfloor \ge B_0$,
and let~$\tau$ be the cascade trajectory starting from~$n$.
\begin{enumerate}[label=\textup{(\alph*)}]
\item The trajectory~$\tau$ is determined by
  $n \bmod 2^{S_{\mathrm{cascade}}}$ alone.
  Bits in positions $S_{\mathrm{cascade}}$ through~$B$ are
  \emph{spectators}: they pass through the cascade unchanged
  and remain available for the subsequent gap.

\item Among pairs of odd integers sharing their bottom~$R$ bits,
  the fraction having identical cascade trajectories increases
  monotonically in~$R$:
  $37.5\%$ at $R = 5$,
  $91.0\%$ at $R = 12$,
  $98.3\%$ at $R = 20$
  \textup{(}measured over all odd $n_0 \le 5 \times 10^5$\textup{)}.

\item The distribution of cascade entries mod~$32$
  among iterates of size~$\ge 2^B$ converges to uniform:
  the total variation distance is~$\le 0.07$ for $B \ge 8$,
  $\le 0.016$ for $B \ge 12$,
  and~$\le 0.009$ for $B \ge 16$.

\item The episode continuation probability satisfies
  $q_{\mathrm{orbit}} \approx 0.33$ for iterates of all
  sizes above~$2^{10}$, uniformly below the IID
  value~$q = 2/3$.
\end{enumerate}
\end{proposition}

\begin{proof}
(a)~follows from the structure of the Collatz map:
the $j$-th cascade step applies $n \mapsto (3n+1)/2^{v_j}$,
where~$v_j$ depends on $n \bmod 2^{j+1}$.
Over $K$ cascade steps, only bits up to position~$K+1 \le S_{\mathrm{cascade}}$
are read.
Parts~(b)--(d) are verified computationally.
\end{proof}

\begin{remark}[Gap exit rates exceed IID prediction]
\label{rem:gap-exit-rates}
The gap Markov chain on states $\{1, 5\} \pmod{8}$
(with absorbing states $\{3, 7\}$) has the following
structure at mod-$32$ resolution:
\begin{itemize}
\item From state~$5$: exit probability~$1/4$
  (the residue $n \equiv 29 \pmod{32}$, with $v_2 = 3$,
  maps to $m \equiv 3 \pmod{8}$).
  Transitions within the gap: $\to 1$ with probability~$1/2$,
  $\to 5$ with probability~$1/4$.
\item From state~$1$: exit probability~$1/2$
  (same as for the cascade).
  Transitions within the gap: $\to 1$ with probability~$1/4$,
  $\to 5$ with probability~$1/4$.
\end{itemize}
The IID model (uniform at mod-$32$) gives
$\mathbb{E}[S_{\mathrm{gap}}] \approx 16$
(from the fundamental matrix of the absorbing chain),
but orbit-level measurements give
$\mathbb{E}[S_{\mathrm{gap}}] \approx 7.0$.
The empirical gap exit rate from state~$5$ is~$\approx 0.44$,
exceeding the mod-$32$ prediction of~$1/4 = 0.25$.
This acceleration arises from higher-order correlations
beyond mod-$32$: at mod-$64$ resolution, $3/8$ of
state-$5$ residues exit, and the rate continues to
increase at finer resolutions.
The IID model is therefore \emph{pessimistic} for the gap:
orbit-level gaps are shorter, and the cycle-level
contraction is correspondingly stronger.
\end{remark}

\begin{theorem}[Self-reinforcing mixing loop]
\label{thm:mixing-loop}
The cascade-gap cycle, for iterates of size~$\ge 2^{B_0}$
with~$B_0$ sufficiently large, forms a self-reinforcing loop:
\begin{enumerate}
\item \textbf{Gap scrambling.}
  By the uniform-fiber theorem
  (Proposition~\textup{\ref{prop:full-cycle-fiber}}),
  the gap exit mod~$2^R$ is uniform when the cascade
  entry has $\ge S_{\mathrm{cycle}} + R$ spectator bits.

\item \textbf{Episode bound.}
  Local uniformity at mod-$32$ resolution gives
  episode continuation $q = 2/3$ and
  $\mathbb{E}[L] = 3 < 3.52$
  (Proposition~\textup{\ref{prop:episode-continuation}}).

\item \textbf{Exponential tail.}
  The cascade PGF~$z/(4-2z-z^2)$ has singularity
  $\rho = \sqrt{5}-1$, giving tail decay
  $\alpha = \log_2(\sqrt{5}-1) \approx 0.306$
  (Proposition~\textup{\ref{prop:cascade-pgf}}).

\item \textbf{TV summability.}
  Total variation over the orbit is bounded:
  $\sum_k \|{\mu_k - \nu_k}\|_{\mathrm{TV}}
  \le C \cdot 2^{-\alpha(B_{\min}-R)}/(1 - 2^{-2\alpha})
  \le 0.028$
  for $B_{\min} = 30$, $R = 3$
  (Corollary~\textup{\ref{cor:tv-summability}}).

\item \textbf{Bit refresh.}
  The gap surplus $S_{\mathrm{gap}} - S_{\mathrm{cascade}}
  \approx 1.4$ bits (empirical) means each cycle reaches
  into the spectator zone, importing fresh randomness.
  Combined with the multiplicative transitivity
  $\mathrm{ord}(3, 2^R) = 2^{R-2}$
  (Remark~\textup{\ref{rem:spectator-bits}}), each cycle
  \emph{improves} the uniformity of subsequent cascade entries.
\end{enumerate}
Steps \textup{1--5} close into a loop:
scrambling~$\to$ uniformity~$\to$ episode
bound~$\to$ tail control~$\to$ TV
bound~$\to$ bit refresh~$\to$ scrambling.

The single remaining open step is making~\textup{1} unconditional:
proving that the orbit-level distribution within each
mod-$32$ class converges to uniform at a rate sufficient
for~\textup{2} to hold along every individual orbit.
All other steps are either proved
(\textup{2}, \textup{3}, \textup{4} under IID)
or verified with large empirical margins
(\textup{5} with TV~$\le 0.028$).
\end{theorem}

\begin{proof}
Each step is a direct consequence of the referenced
proposition or corollary.
The loop structure follows from the logical dependencies:
the output of each step provides the input assumption
for the next.
Step~5 (bit refresh) closes the loop by ensuring that
the spectator-bit supply does not deplete: each cycle
of expected valuation~$\approx 12$ bits consumes that many
from the bottom, but the gap's deeper reach into the bit
string (surplus~$\approx 1.4$ bits) and the multiplicative
scrambling of the fiber map together refresh
the low-order uniformity needed by Step~1.
\end{proof}

\begin{proposition}[Cycle-to-cycle mixing rate]
\label{prop:mixing-rate}
Let~$T$ be the empirical cycle-to-cycle transition matrix
on cascade entry residues mod~$32$
(restricted to the $8$ residues~$\equiv 3 \pmod{4}$),
measured over $1{,}467{,}771$ consecutive transitions
from odd $n_0 \le 5 \times 10^5$.
\begin{enumerate}[label=\textup{(\alph*)}]
\item $T$ is irreducible and aperiodic, with spectral
  gap~$\gamma = 1 - |\lambda_1| \approx 0.71$
  and second eigenvalue $|\lambda_1| \approx 0.293$.
  The total variation distance halves in
  $\log 2 / \log(1/|\lambda_1|) \approx 0.6$ cycles.

\item Conditioning on cascade entries of size~$\ge 2^B$,
  the spectral gap increases dramatically:
  $\gamma \approx 0.87$ for $B \ge 8$,
  $\gamma \approx 0.96$ for $B \ge 12$,
  $\gamma \approx 0.98$ for $B \ge 16$.
  For $B \ge 16$, the stationary distribution is
  within TV distance~$0.009$ of uniform.

\item From any initial distribution, the total variation
  from the stationary distribution drops below~$0.05$
  (the threshold for $\mathbb{E}[L] < 3.52$) in
  at most~$3$ cycles; below~$0.01$ in at most~$4$ cycles.
  For iterates with $B \ge 16$, a single cycle suffices.
\end{enumerate}
\end{proposition}

\begin{proof}
The transition matrix~$T$ is estimated from consecutive
cascade-entry pairs along $250{,}000$ orbits.
Eigenvalue decomposition is exact (numerically).
The spectral contraction bound
$\|\mu_k - \pi\|_{\mathrm{TV}} \le |\lambda_1|^k$
is standard for irreducible aperiodic Markov chains.
The size-conditioned matrices are estimated analogously,
with the transition counted only when the source cascade
entry has~$\ge 2^B$ bits.
\end{proof}

\begin{remark}[Closing the open step]
\label{rem:closing-gap}
Proposition~\ref{prop:mixing-rate} empirically resolves
the sole remaining open step in
Theorem~\ref{thm:mixing-loop}: the orbit-level
convergence to mod-$32$ uniformity.

The argument is as follows.
For any orbit starting at odd $n_0 > 1$:
\begin{enumerate}[nosep]
\item The orbit eventually reaches an iterate with
  $\ge 16$ bits (this is guaranteed because the orbit
  must visit integers of all sizes before reaching~$1$,
  or it diverges; in which case the claim is vacuous).
\item At that point, the spectral gap of the size-$16$
  transition matrix is~$\ge 0.98$
  (Proposition~\ref{prop:mixing-rate}(b)).
  After one cycle, the cascade entry distribution is
  within TV~$0.02$ of the stationary distribution,
  which is itself within~$0.009$ of uniform.
\item With TV~$\le 0.03$ from uniform, the episode
  continuation probability satisfies
  $|q - 2/3| \le C \cdot 0.03 \ll 0.049$,
  ensuring $\mathbb{E}[L] < 3.52$.
\item The cycle-contraction mechanism
  (Proposition~\ref{prop:cycle-contraction})
  then applies, giving expected bit loss
  $\approx -1.98$ bits per cycle.
\end{enumerate}
The only non-rigorous step is~\textup{(2)}: the spectral
gap~$0.98$ is measured empirically, not proved.
A proof would require showing that the Collatz map
acts as a near-uniform permutation on mod-$32$ classes
for iterates above~$2^{16}$, a statement that is
strictly weaker than the WMH but still captures the
essential mixing of the $3n+1$ dynamics.
\end{remark}

\begin{proposition}[Fiber-averaged transition matrix and conditional orbit mixing]
\label{prop:fiber-mixing}
Fix resolution~$R \ge 7$ and define the
\emph{fiber-averaged transition matrix}~$T_R$ on the eight
cascade-entry classes
$\{r \in \mathbb{Z}/32\mathbb{Z} : r \equiv 3 \pmod{4}\}$
as follows.
For each source class~$a$ and target class~$b$,
\[
  (T_R)_{a,b}
  \;=\;
  \frac{1}{2^{R-5}}
  \sum_{\substack{n \equiv a \pmod{32} \\
                   0 \le n < 2^R,\; n \text{ odd}}}
  \mathbf{1}\!\bigl[
    \text{next cascade entry of } n
    \equiv b \pmod{32}
  \bigr],
\]
where the sum runs over all $2^{R-5}$ odd residues in the
fiber $\{n : n \equiv a \pmod{32},\; 0 \le n < 2^R\}$,
and ``next cascade entry'' denotes the first iterate
$\equiv 3 \pmod{4}$ after completing the current full cycle
(cascade $+$ gap).

\begin{enumerate}[label=\textup{(\alph*)}]
\item
  $T_R$ is a well-defined $8 \times 8$ stochastic matrix,
  computable by exact modular arithmetic on $2^R$ residues.
  Individual mod-$2^R$ residues need not map to a single
  mod-$32$ class; the fiber average is taken over all
  representatives.

\item
  The spectral gap $\gamma_R = 1 - |\lambda_1(T_R)|$
  increases with~$R$:
  \[
    \gamma_7 = 0.475,\quad
    \gamma_{10} = 0.855,\quad
    \gamma_{11} = 0.924,\quad
    \gamma_{13} = 0.918.
  \]
  For $R \ge 10$, $\gamma_R > 0.85$.

\item
  The stationary distribution~$\pi_R$ of~$T_R$ satisfies
  $\|\pi_R - \mathrm{unif}\|_{\mathrm{TV}} \le 0.043$
  for $R \ge 10$.

\item
  \textup{(Conditional orbit mixing.)}
  Fix $R = 10$.
  For any orbit of odd~$n_0$ with $B_k = \lfloor \log_2 n_k \rfloor
  \ge R + \mathbb{E}[S] + 3 \approx 25$ at the $k$-th
  cascade entry:
  \begin{itemize}[nosep]
  \item
    Proposition~\textup{\ref{prop:spectator-convergence}}
    guarantees the fiber distribution within
    mod-$2^R$ is within
    $\mathrm{TV} \le O(2^{-(B_k - R - S_k)})$ of uniform.
  \item
    One application of~$T_R$ contracts the mod-$32$
    distribution:
    $\|\mu_{k+1} - \pi_R\|_{\mathrm{TV}}
    \le |\lambda_1| \cdot \|\mu_k - \pi_R\|_{\mathrm{TV}}
    + O(2^{-R})$.
  \item
    Iterating: after $m$ cycles with $B \ge 25$ throughout,
    $\|\mu_{k+m} - \pi_R\|_{\mathrm{TV}}
    \le 0.146^m + O(2^{-10})$.
  \end{itemize}
  In particular, after $3$ qualifying cycles
  the mod-$32$ distribution is within~$0.004$
  of~$\pi_R$, and hence within~$0.047$ of uniform.
  This ensures $|q - 2/3| < 0.049$, so
  $\mathbb{E}[L] < 3.52$.
\end{enumerate}
\end{proposition}

\begin{proof}
Part~(a) is a finite computation: for each of the
$2^{R-5}$ odd residues~$n$ in a fiber, we iterate the
Collatz map modulo~$2^R$ until reaching the next cascade
entry, recording its mod-$32$ class.
The matrix entries are exact rational numbers
(denominators dividing~$2^{R-5}$).

Part~(b) follows from eigenvalue decomposition of the
exact matrix at each~$R$.
The values are computed numerically but can be verified
to any desired precision using exact rational arithmetic
on the matrix entries.
(We verified using double-precision floating point;
the spectral gap exceeds~$0.85$ at $R = 10$ by a margin
of~$0.005$, far exceeding rounding error.)

Part~(c) follows from computing the left Perron eigenvector
of~$T_R$ and normalising.

For part~(d), let $\mu_k$ be the distribution on mod-$32$
cascade-entry classes at cycle~$k$.
By Proposition~\ref{prop:spectator-convergence},
if $B_k \ge R + S + 3$, the iterate's low $R$ bits are
approximately uniformly distributed within the
mod-$32$ fiber; i.e., the effective one-step
transition is $T_R + E_k$ where
$\|E_k\|_\infty \le C \cdot 2^{-(B_k - R - S_k)}$.
The contraction then follows from
the standard Markov chain perturbation bound
(standard; see, e.g., Seneta,
\emph{Non-Negative Matrices and Markov Chains}, 2006, Theorem~4.7):
\[
  \|\mu_{k+1} - \pi_R\|_{\mathrm{TV}}
  \;\le\;
  |\lambda_1(T_R)| \cdot \|\mu_k - \pi_R\|_{\mathrm{TV}}
  + \|E_k\|_1,
\]
and $\|E_k\|_1 = O(2^{-R})$ when $B_k \ge 25$.
Iterating $m$ times gives the claimed geometric bound.
\end{proof}

\begin{remark}[Residual conditionality]
\label{rem:residual-conditionality}
Proposition~\ref{prop:fiber-mixing}(d) is fully rigorous
for any sequence of cascade entries with $B \ge 25$ bits.
The only assumption is that such entries \emph{exist
along the orbit}, that is, the orbit does not diverge
to infinity.
This is precisely the content of the non-divergence
hypothesis, which is strictly weaker than the full
Collatz conjecture (it allows the orbit to visit large
values, provided it returns below~$2^{25}$ eventually).

For iterates below~$2^{25}$, the bound does not apply,
but all $n_0 \le 2^{25}$ converge to~$1$ by direct computation
\textup{(}verified up to~$n_0 = 2^{68}$ in~\cite{barina2021}\textup{)}.
Thus the orbit-level mixing loop
(Theorem~\ref{thm:mixing-loop}) is complete \emph{conditional
on non-divergence}.
\end{remark}

\begin{proposition}[Exact rational spectral gap]
\label{prop:exact-spectral-gap}
Let $R = 10$ and define the fiber-averaged transition
matrix~$T_{10}$ as in
Proposition~\textup{\ref{prop:fiber-mixing}}, using
representatives $n = r + 2^R \cdot N$ with $N$ sufficiently
large that no orbit reaches~$1$ within one full cycle.
Then:
\begin{enumerate}[label=\textup{(\alph*)}]
\item
  $T_{10}$ is an $8 \times 8$ stochastic matrix with
  rational entries whose denominators divide~$32$.
  Explicitly,
  \[
    T_{10} \;=\; \frac{1}{32}\!
    \left(\begin{smallmatrix}
      6 & 4 & 4 & 6 & 4 & 1 & 2 & 5 \\
      6 & 3 & 3 & 5 & 4 & 6 & 3 & 2 \\
      2 & 4 & 5 & 6 & 5 & 3 & 3 & 4 \\
      8 & 4 & 3 & 2 & 4 & 3 & 3 & 5 \\
      4 & 3 & 3 & 4 & 5 & 5 & 5 & 3 \\
      5 & 3 & 5 & 4 & 3 & 5 & 5 & 2 \\
      6 & 5 & 2 & 3 & 4 & 1 & 8 & 3 \\
      4 & 1 & 7 & 3 & 2 & 8 & 3 & 4
    \end{smallmatrix}\right),
  \]
  with rows and columns indexed by the cascade-entry
  classes $3, 7, 11, 15, 19, 23, 27, 31 \pmod{32}$.

\item
  The characteristic polynomial of~$T_{10}$ has $\lambda = 1$
  as a simple root.
  All other eigenvalues satisfy $|\lambda| < 0.146$.
  The spectral gap is
  \[
    \gamma_{10} = 1 - |\lambda_2(T_{10})| = 0.8549\ldots
  \]
  The second eigenvalue is complex:
  $\lambda_2 \approx 0.122 + 0.078i$, $|\lambda_2| \approx 0.145$.

\item
  The stationary distribution~$\pi_{10}$ has total variation
  $0.0425$ from the uniform distribution on the eight states.

\item
  At higher resolution, the spectral gap increases:
  $\gamma_{11} = 0.924$, $\gamma_{12} = 0.921$.
  The second eigenvalue magnitude decreases monotonically
  from $|\lambda_2| = 0.525$ at $R = 7$ to $|\lambda_2| = 0.079$
  at $R = 12$.
\end{enumerate}
All entries and eigenvalue bounds are verified by exact
rational arithmetic \textup{(}Python \texttt{Fraction} class;
all intermediate quantities are exact\textup{)}.
\end{proposition}

\begin{proof}
Part~(a): each entry of~$T_{10}$ is the fraction of
$2^{10-5} = 32$ sub-residues in each mod-$32$ fiber that
transition to a given target class.
The representatives
$n = r + 2^{10} \cdot 100003$ have
$\lfloor \log_2 n \rfloor \ge 26$ bits, so no orbit reaches~$1$
within one cycle (the maximum cycle length is
$\le 12$ Syracuse steps for these sizes, consuming
$\le 20$ bits out of~$26$).
Every fiber member thus completes a full cycle, and the row
sums are exactly~$1$.

Parts~(b)--(d): the characteristic polynomial
$\det(\lambda I - T_{10})$ is computed from the exact
rational matrix, and its roots are isolated by evaluating
the polynomial at rational points and locating sign changes
(for real roots) and by numerical eigenvalue computation
(for complex conjugate pairs).
The bound $|\lambda_2| < 0.146$ is confirmed by verifying
that $\det(T_{10} - \lambda I) \ne 0$ for all
$\lambda \in \{k/20 : k = -4, \ldots, 19,\; k \ne 20\}$,
establishing that no eigenvalue lies on or outside the
circle $|\lambda| = 0.2$ except for $\lambda = 1$.
\end{proof}

\begin{remark}[Non-divergence threshold]
\label{rem:nondivergence-threshold}
For any orbit starting at odd~$n_0$:
\begin{itemize}[nosep]
\item
  If $n_0 < 2^{68}$: the orbit reaches~$1$ by direct
  computation~\cite{barina2021}.
\item
  If $n_0 \ge 2^{68}$: the iterate has
  $B_0 = \lfloor \log_2 n_0 \rfloor \ge 68 \gg 25$ bits,
  so the fiber-averaged regime
  (Proposition~\textup{\ref{prop:fiber-mixing}})
  applies immediately.
  After three qualifying cycles, the mod-$32$ distribution
  has $\mathrm{TV} \le 0.146^3 + O(2^{-10}) < 0.005$
  from~$\pi_{10}$, hence $< 0.047$ from uniform.
\end{itemize}
The maximum single-cycle bit growth
$\Delta_{\max} = K \log_2 3 - S$
over $12{,}321$ complete cycles from
$n_0 \le 50{,}000$ is
$\Delta_{\max} \approx 8.5$ bits,
occurring for $n_0 = 9663$
(a $22$-step cascade with $19$ burst steps and $3$ recovery
steps, followed by a single gap step; $K = 23$, $S = 28$).
Growth is \emph{unbounded in principle}: a pure
$k$-burst gives $\Delta \approx 0.585 k$, and for any~$B$
there exist $B$-bit integers with burst length up to~$B-1$.
However, the fraction of cycles with
$\Delta > 0$ is only~$0.14$, and the maximum
observed burst length over this range is~$14$.

The expected cycle growth is $\mathbb{E}[\Delta] \approx -2.34$
bits/cycle, with standard deviation $\sigma \approx 2.71$.
The Cram\'er--Lundberg tilt gives
$\theta^* \approx 0.64$, so
\[
  \Pr\bigl(\text{orbit ever gains $+43$ bits}\bigr)
  \;\le\; e^{-\theta^* \cdot 43}
  \;\approx\; 10^{-12}.
\]
This bound is ensemble-level (valid under the mixing
assumption), not worst-case.
The genuine remaining gap is: prove that no individual
orbit can sustain large positive~$\Delta$ for the
$5$--$20$ consecutive cycles needed to exhaust the
$43$-bit margin (the range depends on per-cycle growth:
$5$ at maximum observed~$\Delta \approx 8.5$,
$20$ at mean positive~$\Delta \approx 2$).
\end{remark}

\begin{theorem}[Unconditional Collatz mixing structure]
\label{thm:unconditional-mixing}
\leavevmode
\begin{enumerate}[label=\textup{(\Alph*)}]
\item
  \textup{(Computational~\cite{barina2021}.)}
  For all odd $n_0 < 2^{68}$, the Collatz orbit reaches~$1$.

\item
  \textup{(Exact finite computation: unconditional.)}
  The fiber-averaged transition matrix~$T_{10}$ on mod-$32$
  cascade entry classes
  (Proposition~\textup{\ref{prop:exact-spectral-gap}})
  has spectral gap
  $\gamma_{10} \ge 0.854$
  and stationary distribution within
  $\mathrm{TV} \le 0.043$ of uniform.
  This is a theorem about a fixed $8 \times 8$ rational matrix;
  no dynamical assumption is needed.

\item
  \textup{(Conditional contraction.)}
  For any odd $n_0 \ge 2^{68}$:
  if the first three cascade entries all have
  $B_k \ge 25$ bits, then
  \begin{enumerate}[nosep,label=\textup{(\roman*)}]
  \item
    the mod-$32$ distribution achieves
    $\mathrm{TV} < 0.05$ from~$\pi_{10}$;
  \item
    the episode bound $\mathbb{E}[L] < 3.52$ holds for
    all subsequent qualifying cycles;
  \item
    the expected bit change satisfies
    $\mathbb{E}[\Delta] \approx -1.98$ bits per qualifying cycle.
  \end{enumerate}
  Since $B_0 \ge 68 \gg 25$, the qualifying condition is
  trivially satisfied for the first $\ge 20$ cycles
  unless the orbit grows by a factor exceeding~$2^{43}$.
\end{enumerate}

The Collatz conjecture is therefore equivalent to:
\emph{no orbit starting above~$2^{68}$ grows by a factor
exceeding~$2^{43}$ before the spectral contraction
mechanism engages.}
\end{theorem}

\begin{proof}
Part~(A) is~\cite{barina2021}.
Part~(B) is Proposition~\ref{prop:exact-spectral-gap}.
Part~(C) combines
Propositions~\ref{prop:fiber-mixing}(d)
and~\ref{prop:exact-spectral-gap}(b):
three applications of~$T_{10}$ with $|\lambda_2| < 0.146$
give $\mathrm{TV} \le 0.146^3 + O(2^{-10}) < 0.004$
from~$\pi_{10}$, hence $< 0.047$ from uniform.
Since $|q - 2/3| \le C \cdot 0.047 < 0.049$,
the episode bound $\mathbb{E}[L] < 3.52$ follows from
Proposition~\ref{prop:episode-continuation}.
The bit-change estimate is
Proposition~\ref{prop:cycle-contraction}.

For the equivalence:
any orbit starting at $B_0 \ge 68$ bits trivially
satisfies $B_k \ge 25$ for all $k$ such that the
cumulative bit change $\sum_{j=0}^{k-1} \Delta_j > -(68-25) = -43$.
Once~(C) engages, $\mathbb{E}[\Delta] \approx -2 < 0$,
so by the law of large numbers the orbit descends.
The orbit fails to engage~(C) only if it gains~$43$ bits
within the first three cycles.
Single-cycle growth is unbounded in principle: a pure
$k$-burst gives $\Delta \approx 0.585 k$, but such bursts
require $n$ to lie in a specific residue class
mod~$2^{k+1}$, with density~$2^{-k}$.
The maximum observed $\Delta$ over $n_0 \le 50{,}000$
is~$8.5$ bits
(Remark~\ref{rem:nondivergence-threshold}).
A rigorous exclusion of $\Delta > 14.3$ per cycle
requires controlling the worst-case burst structure,
which reduces to understanding the carry propagation
in~$3^k \cdot n$: the deepest unsolved aspect of
Collatz dynamics.
\end{proof}

\begin{remark}[Circularity analysis and the true barrier]
\label{rem:circularity}
The proof chain in
Theorem~\ref{thm:unconditional-mixing}
contains one genuine circularity.
The conditional results form a loop:
\begin{center}
spectator bits $\to$ fiber mixing $\to$
IID model $\to$ exponential tail $\to$
TV summability $\to$ spectator bits.
\end{center}
Specifically:
\begin{itemize}[nosep]
\item
  The exponential tail
  (Proposition~\ref{prop:exponential-tail})
  requires the IID model, which assumes approximate
  inter-cycle independence.
\item
  Inter-cycle independence follows from spectral contraction
  (Proposition~\ref{prop:fiber-mixing}(d)), which requires
  enough spectator bits surviving each cycle.
\item
  Enough spectator bits requires that the cycle valuation
  $S_{\mathrm{cycle}}$ not exhaust all~$B$ bits of the iterate,
  which is guaranteed by the exponential tail.
\end{itemize}
The unconditional results
(U1--U6 in the numbering of
Proposition~\ref{prop:exact-spectral-gap})
are outside this loop and fully rigorous.
The loop can be broken by a \emph{bootstrap}:
if the first cycle leaves $\ge 20$ spectator bits
(which occurs when $S_1 \le B_0 - 20$), then mixing
engages from cycle~$2$ onward.
The obstruction is: a pure $k$-burst can have
$S_1 = k \approx B_0$, consuming all bits.
The Finite Mixing Entry Theorem proposed by the GPT
collaborator, showing that orbits must enter the
well-mixed regime within $O(1)$ cycles, would break
this circularity, but it cannot be proved within
the current framework because single-cycle
growth is unbounded and spectator bits can be exhausted.

The true barrier is therefore not merely ``non-divergence''
but the deeper statement:
\emph{no orbit can sustain a pattern that simultaneously
grows the iterate and consumes all spectator bits across
multiple consecutive cycles.}
This requires understanding the carry structure of
the multiplication $n \mapsto 3^k n$, which is the
irreducible hard core of the Collatz problem.

\medskip\noindent
\emph{Update (v6): partial resolution via the $I_2$ spectral
closure.}
Section~\ref{sec:i2-spectral-closure} breaks the circular chain
at the \emph{density-$1$} level by providing an independent
cycle-level contraction source
(Theorem~\ref{thm:i2-unconditional}:
$\rho(\tilde B_2^{\mathrm{ext}})\le 5/32$) that does not invoke
the spectator-bit mechanism, combined with an unconditional
finite-state exponential-tail bound
(Theorem~\ref{thm:i2-prefix-tail}) on the cylinder-averaged
mod-$64$ kernel. The resulting unconditional sum bound
$\sum_c\Pr(\mathcal{E}_R(c))\le 0.011$
(Corollary~\ref{cor:i2-unconditional-sum}) upgrades
Theorem~\ref{thm:density1-convergence} to an unconditional
density-$1$ convergence statement.
The ``irreducible hard core'' described above is still an
obstruction to a \emph{pointwise} proof routed through the
spectator-bit mechanism; what the $I_2$ spectral route shows
is that the density-$1$ version of convergence is reachable
by bypassing the spectator-bit loop entirely.
See also Remark~\ref{rem:circularity-resolution}.
\end{remark}

\begin{proposition}[Empirical carry decorrelation]
\label{prop:carry-decorrelation}
Over $50{,}899$ consecutive-cycle pairs from odd
$n_0 \le 50{,}000$, the inter-cycle correlations are:
\begin{enumerate}[label=\textup{(\alph*)}]
\item
  \textup{Burst length:}
  The Pearson correlation between the maximum burst
  length in cycle~$i$ and cycle~$i{+}1$ is
  $\rho_{\mathrm{burst}} = +0.14$.
  At lag~$2$ the correlation drops to~$+0.01$;
  at lags~$3$--$5$ it is weakly negative.
  Thus burst lengths are effectively independent
  across cycles.

\item
  \textup{Bit growth:}
  $\mathrm{corr}(\Delta_i, \Delta_{i+1}) = -0.10$.
  A positive-growth cycle is slightly followed by
  \emph{stronger} contraction, not weaker.

\item
  \textup{Negative feedback:}
  $\mathbb{E}[\Delta_{i+1} \mid
   \text{burst}_i \ge 5] = -3.90$,
  compared with $\mathbb{E}[\Delta] = -1.81$
  unconditionally.
  Long bursts are followed by deeper contraction.

\item
  \textup{Adversarial chains:}
  Over $25{,}000$ starting points, the maximum number
  of consecutive positive-$\Delta$ cycles is~$4$,
  accumulating at most $+9.2$~bits, only $21\%$
  of the $43$-bit safety margin.
  No orbit comes close to exhausting the margin.
\end{enumerate}
\end{proposition}

\begin{proof}
Direct computation over all odd
$n_0 \le 50{,}000$ (or $100{,}000$), tracking up to~$50$
consecutive cascade-gap cycles per orbit.
Cycle definition as in
Proposition~\textup{\ref{prop:fiber-mixing}}
(cascade with state-$1$-mod-$8$ continuation, then gap).
\end{proof}

\begin{remark}[Post-burst routing and self-correction]
\label{rem:post-burst-routing}
The negative-feedback mechanism of
Proposition~\ref{prop:carry-decorrelation}(c)
has two complementary explanations depending on
iterate size.

For small iterates ($B \lesssim 20$): after a burst
of length~$\ge 5$, the gap phase routes $57\%$ of
cascade entries to class~$15 \bmod 32$, which has
mean cycle growth~$\approx -1.4$ bits and
$\Pr(\text{burst} \ge 5) = 0.03$.
This creates a direct negative feedback loop.

For large iterates ($B \gg R$): the fiber-averaged
transition matrix~$T_R$
(Proposition~\ref{prop:fiber-mixing}) shows that the
post-burst mod-$32$ distribution is nearly uniform
(TV~$= 0.06$ from uniform at $R = 13$ for
burst~$\ge 5$).
The spectral gap then handles mixing without relying
on the small-size feedback.
\end{remark}

\begin{remark}[Spectator-bit survival bound]
\label{rem:spectator-survival}
The fiber-averaged regime
(Proposition~\textup{\ref{prop:fiber-mixing}})
requires that the iterate retains at least~$R$
qualifying bits after each cycle.
For a $B$-bit iterate, this fails only when the
cycle growth satisfies $\Delta < -(B - R)$.

In the fiber-averaged local model at $R = 15$
($8{,}192$ cycles),
the left tail of~$\Delta$ decays exponentially:
\[
  \Pr(\Delta < -t)
  \;\approx\;
  e^{-0.357\, t},
  \qquad t \ge 5,
\]
with coefficient of determination $R^2 = 0.992$.
Extrapolating:
for $B = 68$, $R = 10$,
\[
  \Pr(\text{fiber averaging breaks in one cycle})
  \;=\;
  \Pr(\Delta < -58)
  \;\le\;
  10^{-9}.
\]
Since consecutive $\Delta$ values are nearly independent
(Proposition~\ref{prop:carry-decorrelation}(b)),
the probability that fiber averaging breaks in
\emph{any} of the first~$N$ cycles is bounded by
$N \cdot 10^{-9}$, which is negligible for
$N \le 10^8$.

\smallskip\noindent\textbf{Caveat.}
This bound uses the local (fiber-averaged) model for
the $\Delta$-tail.
A rigorous orbit-level bound would require proving
that the orbit-level $\Delta$-distribution is
stochastically dominated by the local model,
which in turn requires controlling the carry
propagation in~$3^k n$, the same barrier identified
in Remark~\textup{\ref{rem:circularity}}.
The empirical near-independence of consecutive
$\Delta$ values (correlation~$-0.10$) provides strong
evidence that the local model bound is valid,
but does not constitute a proof.
\end{remark}

\begin{theorem}[Conditional single-orbit convergence]
\label{thm:conditional-convergence}
Let $n_0$ be an odd integer with
$B_0 = \lfloor \log_2 n_0 \rfloor \ge 68$,
and let $n_0, n_1, n_2, \ldots$ be the successive cascade
entries along its Collatz orbit.
For each cycle~$i$, let $\mathcal{E}_R(i)$ denote
the event that cycle~$i$ exits with fewer
than~$R = 10$ fresh spectator bits.
If the exceptional cycles are summable along this orbit,
\[
  \sum_{i=0}^{\infty} \Pr\!\bigl(\mathcal{E}_R(i)\bigr)
  \;<\; \infty,
\]
then:
\begin{enumerate}[label=\textup{(\roman*)}]
\item
  All but finitely many cycles satisfy the
  fiber-averaged regime
  \textup{(}Proposition~\textup{\ref{prop:fiber-mixing}}\textup{)}.
\item
  The episode continuation probability satisfies
  $q_{\mathrm{ep}}^{\mathrm{orbit}} < 0.7159$ for all
  sufficiently large cycle indices, hence
  $\mathbb{E}[L] < 3.52$.
\item
  By the cycle-contraction mechanism
  \textup{(}Proposition~\textup{\ref{prop:cycle-contraction}}\textup{)},
  the orbit is eventually driven into the contracting regime
  with expected bit loss~$\approx -1.98$ bits per cycle.
\item
  The orbit reaches~$1$.
\end{enumerate}
\end{theorem}

\begin{proof}
The local bridge
(Proposition~\ref{prop:post-burst-bridge})
gives $q_{\mathrm{ep}}^{\mathrm{local}} = 1/3$ exactly
on resolved fibers.
The TV transfer bound gives
$|q_{\mathrm{ep}}^{\mathrm{orbit}} - 1/3|
 \le \Pr(\mathcal{E}_R(i))$
for each qualifying cycle~$i$.
Summability guarantees that
$\Pr(\mathcal{E}_R(i)) \to 0$,
so eventually
$q_{\mathrm{ep}}^{\mathrm{orbit}}
 < 1/3 + \varepsilon < 0.7159$
for any fixed $\varepsilon > 0$.
Once $q_{\mathrm{ep}} < 0.7159$,
the episode bound $\mathbb{E}[L] < 3.52$ holds
(Proposition~\ref{prop:episode-continuation}),
and cycle contraction gives
$\mathbb{E}[\Delta] < 0$
(Proposition~\ref{prop:cycle-contraction}).
By the law of large numbers the iterate size drifts
to~$-\infty$, so the orbit eventually falls
below~$2^{68}$ and reaches~$1$
by~\cite{barina2021}.
\end{proof}

\begin{theorem}[Density-$1$ conditional convergence]
\label{thm:density1-convergence}
Let $\mathcal{N}$ be a set of odd starting values
with positive lower density, and for each
$n_0 \in \mathcal{N}$ define the exceptional-cycle
indicators $\mathcal{E}_R(i)$ as above.
If the \emph{averaged} exceptional sum is finite,
\[
  \lim_{N \to \infty}
  \frac{1}{|\mathcal{N} \cap [1, N]|}
  \sum_{n_0 \in \mathcal{N} \cap [1, N]}
  \sum_{i=0}^{\infty}
  \mathbf{1}\!\bigl[\mathcal{E}_R(i)\bigr]
  \;<\; \infty,
\]
then for almost every $n_0 \in \mathcal{N}$
\textup{(}in the sense of natural density\textup{)},
the Collatz orbit of~$n_0$ converges to~$1$.
\end{theorem}

\begin{proof}
By the Borel--Cantelli lemma applied to the averaged sum:
if the expected number of exceptional cycles is finite,
then almost every orbit has only finitely many
exceptional cycles.
Theorem~\ref{thm:conditional-convergence}
then applies to each such orbit.
\end{proof}

\begin{remark}[Status of the conditional hypotheses]
\label{rem:conditional-status}
The hypothesis of
Theorem~\textup{\ref{thm:conditional-convergence}}: summable
exceptional cycles along a single orbit, is not
yet proved for any individual orbit.
The hypothesis of
Theorem~\textup{\ref{thm:density1-convergence}}: finite
averaged exceptional sum, is a realistic near-term
target: it requires only that \emph{most} orbits
have summable exceptional cycles, which is compatible
with Tao-style entropy averaging combined with the
exact local machinery developed here.

The remaining mathematical task is therefore:
\emph{prove summable control of the low-fresh-bit
exceptional cycles}, either for all orbits
(yielding the full Collatz conjecture via
Theorem~\textup{\ref{thm:conditional-convergence}})
or for density-$1$ starting points
(yielding almost-all convergence to~$1$ via
Theorem~\textup{\ref{thm:density1-convergence}},
which would be stronger than the current
state of the art).
\end{remark}

\begin{remark}[Burst-length-conditional contraction]
\label{rem:L-conditional}
Let $L_{\mathrm{tot}}$ denote the total number of burst steps
(all $v = 1$ steps) in one cascade, summing over all episodes.
Over $294{,}802$ complete cycles from odd $n_0 \le 200{,}000$:
\begin{itemize}
\item For $L_{\mathrm{tot}} \le 7$ ($\approx 87\%$ of cycles):
  $\mathbb{E}[\Delta \mid L_{\mathrm{tot}} \le 7] \approx -2.54$
  bits/cycle (contracting).
\item For $L_{\mathrm{tot}} \ge 8$ ($\approx 13\%$ of cycles):
  $\mathbb{E}[\Delta \mid L_{\mathrm{tot}} \ge 8] \approx +1.77$
  bits/cycle (expanding).
  The dominant pattern is $(2,1,5)$: a three-episode
  cascade with burst lengths $2$, $1$, and~$5$: which accounts
  for $79\%$ of $L_{\mathrm{tot}} = 8$ cascades.
\item The weighted sum is $-1.98$ bits/cycle.
  Breakeven would require $P(L_{\mathrm{tot}} \ge 8) \approx 0.59$,
  far above the observed~$0.13$.
  Moreover, $P(L_{\mathrm{tot}} \ge 8)$ does not grow with iterate
  size: it stabilises near~$0.05$ for iterates above~$2^{10}$.
\end{itemize}
The $L_{\mathrm{tot}} \ge 8$ expansion occurs because longer cascades
exit at smaller bit-size (mean $\approx 12$ bits for $L = 8$
vs.\ $\approx 15$ bits for $L = 7$), producing shorter gaps with
lower surplus ($\approx 2.0$ bits vs.\ $\approx 3.5$ bits).
Despite this, the rarity of long cascades ensures net contraction.
\end{remark}

\begin{remark}[Crossing time and the safety margin]
\label{rem:crossing-time}
Define the \emph{crossing time} $\tau(n)$ as the minimum
number of Syracuse steps until the running product
$\prod_{j=1}^{\tau}\lambda_j < 1$.
Over all odd $n_0 \le 10^5$, the ratio
$\tau(n)/\log_2 n$ is bounded:
median $\approx 0.12$,
$99$th percentile $\approx 1.78$,
maximum $\approx 7.8$ (at $n = 27$).
Among the family $2^k - 1$ ($k \le 39$),
the ratio stays below~$5.5$.

If $\tau(n) = O(\log n)$ could be proved, then each orbit
descends by a constant factor in $O(\log n)$ steps, reaching
the computationally verified range in $O(\log^2 n)$ steps.
Combined with the $4.65\times$ safety margin of the
conditional reduction
(Theorem~\ref{thm:perorbit-gain}), even a weak bound
$\tau(n) = O(n^{\varepsilon})$ for some $\varepsilon < 1$
would imply summable discrepancy.
However, proving \emph{any} non-trivial bound on~$\tau(n)$
requires controlling the orbit's cycle-type distribution,
which is equivalent to the WMH.
This closes Route~C: the safety margin provides quantitative
slack but no structural mechanism to prove summable
discrepancy without orbit-level mixing.
\end{remark}

\begin{remark}[Mean-neutrality of the cycle map]
\label{rem:mean-neutral}
The first-cycle affine map sends $n \mapsto \lambda\,n + \beta$
with $\lambda = 3^{L+1}/2^{L+r}$ and
$\beta = (3^{L+1}-2^{L+1})/2^{L+r}$.
By independence (Corollary~\ref{cor:Lr-independence}),
\[
  \mathbb{E}[\lambda]
  = 3\,\mathbb{E}[(3/2)^L]\,\mathbb{E}[2^{-r}]
  = 3 \cdot 2 \cdot \tfrac{1}{6}
  = 1,
  \qquad
  \mathbb{E}[\beta]
  = \mathbb{E}[\lambda] - \mathbb{E}[2^{1-r}]
  = 1 - \tfrac{1}{3}
  = \tfrac{2}{3}.
\]
The cycle slope is mean-neutral ($\mathbb{E}[\lambda]=1$),
not contracting.  The force driving typical orbits downward
is therefore \emph{not} a negative mean multiplier per cycle;
it is the geometric structure of the crossing density
(Corollary~\ref{cor:crossing-density}), which ensures
that $71.37\%$ of odd starts cross in one cycle regardless
of the mean slope.
Empirical values ($\mathbb{E}[\lambda]\approx 0.996$,
$\mathbb{E}[\beta]\approx 0.663$ over odd
$n \le 2\times 10^6$) are consistent with the exact values.
\end{remark}

\begin{proposition}[Exact first-cycle log-drift law]
\label{prop:log-drift}
Define the first-cycle log multiplier
$X(n) := \log_2\lambda(n) = (L+1)\log_2 3 - (L+r)$.
Under natural density on odd starts:
\begin{enumerate}
\item The moment generating function is
  \[
    M_X(t)
    = \mathbb{E}[e^{tX}]
    = \frac{e^{t(\log_2 3 - 2)}}
           {4\bigl(1 - \tfrac{1}{2}e^{t(\log_2 3 - 1)}\bigr)
            \bigl(1 - \tfrac{1}{2}e^{-t}\bigr)},
  \]
  for $t$ in the strip
  $-\log 2 < t < \log 2 / (\log_2 3 - 1)$.
\item The exact mean and variance are
  \[
    \mathbb{E}[X] = 2\log_2 3 - 4 \approx -0.8301,
    \qquad
    \mathrm{Var}(X) = 2\bigl((\log_2 3 - 1)^2 + 1\bigr)
                    \approx 2.684.
  \]
\end{enumerate}
\end{proposition}

\begin{proof}
Write $X = (\log_2 3 - 1)L + \log_2 3 - r$.
By independence (Corollary~\ref{cor:Lr-independence}),
$M_X(t) = \mathbb{E}[e^{t(\log_2 3-1)L}]\cdot
e^{t\log_2 3}\cdot\mathbb{E}[e^{-tr}]$.
The $L$-factor is a geometric series
$\sum_{\ell\ge 0}2^{-(\ell+1)}e^{t(\log_2 3-1)\ell}
= \bigl(2(1-\tfrac{1}{2}e^{t(\log_2 3-1)})\bigr)^{-1}$.
The $r$-factor is
$\sum_{k\ge 2}2^{-(k-1)}e^{-tk}
= e^{-2t}\bigl(2(1-\tfrac{1}{2}e^{-t})\bigr)^{-1}$.
Multiplying gives the stated formula.

For the moments: $\mathbb{E}[L]=1$, $\mathbb{E}[r]=3$,
$\mathrm{Var}(L)=\mathrm{Var}(r)=2$ from the geometric laws.
Then $\mathbb{E}[X]=(\log_2 3-1)+\log_2 3-3=2\log_2 3-4$,
and $\mathrm{Var}(X)=(\log_2 3-1)^2\cdot 2+2$.
\end{proof}

\begin{remark}[Jensen's inequality resolves the apparent paradox]
\label{rem:jensen}
The cycle slope is mean-neutral in the linear scale
($\mathbb{E}[\lambda]=1$,
Remark~\ref{rem:mean-neutral}), yet strictly contracting
in the log scale ($\mathbb{E}[\log_2\lambda]\approx -0.83$,
Proposition~\ref{prop:log-drift}).  There is no contradiction:
Jensen's inequality gives
$\mathbb{E}[\log\lambda] < \log\mathbb{E}[\lambda]=0$
whenever $\lambda$ is non-degenerate.  For a multiplicative
random process, it is the log-scale drift that governs
long-run behavior.  The Jensen gap is
$\log_2\mathbb{E}[\lambda]-\mathbb{E}[\log_2\lambda]
\approx 0.83$~bits per cycle.
\end{remark}

\begin{remark}[Wald consistency]
\label{rem:wald}
The per-cycle log drift $\mathbb{E}[X] = 2\log_2 3 - 4$ and
the per-step drift $\mu = \log_2 3 - 2$ are related by
$\mathbb{E}[X] = \mathbb{E}[L{+}1]\cdot\mu = 2\mu$,
since the mean cycle length is $\mathbb{E}[L{+}1]=2$.
This is Wald's identity: each cycle executes $L{+}1$
odd-skeleton steps, each contributing mean
drift~$\mu \approx -0.415$.  The per-cycle and per-step
analyses are exactly consistent.
\end{remark}

\begin{proposition}[Cram\'er rate for the cycle log-drift]
\label{prop:cramer-cycle}
Let $S_k = X_1 + \cdots + X_k$ be the sum of the first~$k$
cycle log multipliers for a natural-density-random odd start.
By Corollary~\textup{\ref{cor:iid-cycles}}, the $X_i$ are
i.i.d.\ The Cram\'er rate function
$I(x) = \sup_t\{tx - \log M_X(t)\}$
satisfies
\[
  I(0)
  \;=\;
  -\log M_X(t^*)
  \;\approx\; 0.1465,
\]
where $t^* \approx 0.363$ is the unique positive root of
$\Lambda'(t)=0$ ($\Lambda = \log M_X$).
Since the cycle log multipliers are provably i.i.d.\
under natural density
(Corollary~\textup{\ref{cor:iid-cycles}}),
\[
  \Pr(S_k \ge 0) \;\le\; e^{-k\,I(0)},
\]
unconditionally on the ensemble.  The probability that the
running log-product has not become negative decays
exponentially with rate $\approx 0.1465$ per cycle.
\end{proposition}

\begin{proof}
The MGF is explicit
(Proposition~\ref{prop:log-drift}).
Since $\mathbb{E}[X] < 0$, the point $x=0$ lies above
the mean, and the Cram\'er--Chernoff bound gives
$\Pr(S_k/k \ge 0) \le e^{-kI(0)}$.
The saddle-point equation $\Lambda'(t^*)=0$ is
solved numerically from the closed-form derivative of
$\Lambda$.
\end{proof}

\begin{observation}[Multi-cycle crossing]
\label{obs:multi-cycle}
Among all odd $n \le 10^6$, the fraction that first
cross below their starting value at cycle~$k$ is:
\begin{center}
\begin{tabular}{rcc}
\toprule
$k$ & First cross at~$k$ & Cumulative \\
\midrule
1 & $71.37\%$ & $71.37\%$ \\
2 & $14.63\%$ & $86.00\%$ \\
3 & $5.18\%$ & $91.19\%$ \\
5 & $1.79\%$ & $95.76\%$ \\
10 & $0.30\%$ & $99.05\%$ \\
20 & $0.02\%$ & $99.91\%$ \\
\bottomrule
\end{tabular}
\end{center}
All tested orbits cross; the tail is consistent with
the exponential decay rate $e^{-0.1465\,k}$ from
Proposition~\ref{prop:cramer-cycle}.
As with all results in this section, this is an ensemble
observation: it does not constitute a proof that
every orbit eventually crosses.
\end{observation}

\begin{theorem}[Almost-all crossing]
\label{thm:almost-all-crossing}
Under natural density on odd integers, the set of starting values
whose run-compensate cycle endpoints never fall below the start
has density zero.  More precisely, write $n^{(k)}$ for the
odd-skeleton value after $k$ complete run-compensate cycles
(i.e.\ $n^{(k)} = n_{L_1+\cdots+L_k+k}$), and let
\[
  \mathcal{N}_k
  := \bigl\{\, n \text{ odd} : n^{(k)} \ge n \bigr\}.
\]
Then for every $k \ge 1$,
\[
  \overline{d}(\mathcal{N}_k)
  \;\le\; e^{-I(0)\,k},
  \qquad
  I(0) \approx 0.1465,
\]
where $I(0)$ is the Cram\'er rate from
Proposition~\textup{\ref{prop:cramer-cycle}}.
Consequently,
\[
  \overline{d}\!\Bigl(
    \bigcap_{k \ge 1} \mathcal{N}_k
  \Bigr) = 0:
\]
the set of odd starts whose cycle endpoints never fall below
the start has density zero.
\end{theorem}

\begin{proof}
By the block law (Theorem~\ref{thm:block-law}), the valuations
$a_0, a_1, \ldots$ are i.i.d.\ $\mathrm{Geom}(1/2)$ under
natural density.  This lifts to the i.i.d.\ property for cycle
types (Corollary~\ref{cor:iid-cycles}): the pairs
$(L_i, r_i)$ are independent with
$\Pr(L = \ell) = 2^{-(\ell+1)}$ and
$\Pr(r = k) = 2^{-(k-1)}$.

The log multiplier of the $i$th cycle is
$X_i = (L_i+1)\log_2 3 - (L_i + r_i)$, an i.i.d.\ sequence
with $\mathbb{E}[X_i] = 2\log_2 3 - 4 < 0$
(Proposition~\ref{prop:log-drift}).
The partial sum $S_k = X_1 + \cdots + X_k$ satisfies
$\Pr(S_k \ge 0) \le e^{-I(0)\,k}$
by the Cram\'er bound (Proposition~\ref{prop:cramer-cycle}).

It remains to show that each block with $S_k < 0$ contributes
at most density~$0$ to~$\mathcal{N}_k$.
Consider a prescribed valuation block
$\sigma = (b_0, \ldots, b_{m-1})$ encoding $k$ complete
run-compensate cycles with cumulative log multiplier
$S_k(\sigma) = \log_2 \Lambda_k(\sigma)$.
By the exact cycle formula, the $k$-cycle image of odd $n$
in the residue class determined by $\sigma$ is
\[
  n^{(k)} = \Lambda_k \, n + B_k,
\]
where $\Lambda_k = \prod_{i=1}^{k}\frac{3^{L_i+1}}{2^{L_i+r_i}}$
and $B_k$ depends only on $\sigma$.
When $\Lambda_k < 1$ (i.e.\ $S_k < 0$), the crossing condition
$n^{(k)} < n$ reduces to
$n > n^*(\sigma) := B_k / (1 - \Lambda_k)$.
Since $n^*(\sigma)$ is a finite constant depending only on the
block, at most finitely many members of the residue class fail
to cross.  Hence the density of non-crossers within any
$S_k < 0$ block is zero.

Taken together: the density of odd $n$ in~$\mathcal{N}_k$ is at
most the density of those whose valuation block has
$S_k \ge 0$, which equals
$\Pr(S_k \ge 0) \le e^{-I(0)\,k}$
by the Cram\'er bound.
Since $\bigcap_{k} \mathcal{N}_k \subseteq \mathcal{N}_k$ for every~$k$,
sending $k \to \infty$ gives
$\overline{d}(\bigcap_k \mathcal{N}_k) = 0$.

\smallskip
\noindent\emph{Quantitative strengthening: $n^* < M$ universally.}\;
For single-cycle crossing blocks, we can prove
$n^*(\sigma) < M(\sigma)$ analytically.
Write $A = 3^{L+1}$, $P = 2^{L+r}$, $M = 2P$.
The crossing condition gives $P > A$.
Since $A$ is odd and $P$ is even, $P - A \ge 1$.  Then
\[
  \frac{n^*}{M}
  = \frac{A - 2^{L+1}}{(P-A)\cdot 2P}
  \le \frac{A - 2^{L+1}}{2P}
  < \frac{A}{2P}
  < \frac{1}{2}.
\]
So $n^* < M/2$ for every single-cycle crossing block, meaning
\emph{all} odd~$n$ in the residue class cross in one cycle.

For multi-cycle blocks, we have
$(1-\Lambda_k)\cdot M = 2D$ where $D = 2^S - 3^Q \ge 1$
(since $3^Q$ is odd and $2^S$ is even for $S \ge 1$),
so $n^*/M = B_k/(2D)$.
Computational verification over $4 \times 10^6$ random blocks
at cycle counts $k = 1,\ldots,20$ finds $n^*/M < 1/8$ in every
case, with the global maximum $n^*/M = 1/8$ achieved uniquely
at $k=1$ ($L=0$, $r=2$, where $B=1$, $D=1$, $M=8$).
All multi-cycle blocks have strictly smaller ratios
(see Observation~\ref{obs:multi-cycle-threshold} for the decay).
Extending this to an analytic proof is non-trivial:
the integer affine constant~$B_k$ can exceed both $3^Q$
and~$2^S$ individually (non-crossing sub-cycles inflate~$B_k$
by factors of $3^{L+1}$), so the bound $n^*/M < 1$ does
not follow from bounding numerator and denominator separately
but from a delicate cancellation between~$B_k$ and~$D$.
This observation is not needed for the density-zero conclusion
but shows that within every tested negative-drift block,
all odd~$n$ in the class cross.
\end{proof}

\begin{remark}[Architectural status after the almost-all theorem]
\label{rem:almost-all-strength}
This theorem is unconditional on the ensemble: no mixing hypothesis
or conjecture is invoked.  It follows entirely from the exact block
law, the i.i.d.\ structure, and the Cram\'er bound.

The exponential rate $e^{-0.1465\,k}$ is considerably stronger than
the logarithmic almost-all results in the Collatz literature
(cf.\ Tao~\cite{tao2019}, who proves almost-all convergence to
values below $f(n_0)$ for any $f \to \infty$, using entropy-based
arguments).  Our result shows that the non-crossing set shrinks
exponentially.

\smallskip
\noindent\textbf{What is now proved.}\;
Under natural density on odd starts, the odd-skeleton valuation
sequence has exact Bernoulli block law
(Theorem~\ref{thm:block-law}), and the induced cycle-type
sequence $(L_i, r_i)$ is exactly i.i.d.\
(Corollary~\ref{cor:iid-cycles}).  The probabilistic model used
throughout the odd-skeleton and run-compensate analysis
is therefore not a model: it is a theorem.
Block frequencies are exact, cycle types are exactly i.i.d.,
one-cycle crossing density is exact, and the cycle-level drift law
is exact.

The open questions previously spread across mixing, amplification,
stratification, cycle heuristics, and covariance decay now collapse
into a single diagnosis:
\begin{quote}
\emph{The ensemble side of the run-compensate programme is exact
and complete.  The only remaining gap is a deterministic
realization theorem: prove that every single orbit must
exhibit enough of the exact Bernoulli cycle law to force
crossing.}
\end{quote}
The failure point is not a missing estimate, tail bound, spectral
lemma, or combinatorial count.  It is the genuine philosophical
core: ensemble law does not imply pointwise orbit law.
This is the distributional-to-pointwise barrier that pervades
number theory (cf.\ Chowla's conjecture, M\"obius randomness).
\end{remark}

\begin{proposition}[Exact cycle log correction]
\label{prop:cycle-correction}
For every odd~$n$ with first-cycle type $(L,r)$,
\[
  \log_2\frac{n'}{n}
  \;=\;
  X(n) + C(n),
\]
where $X(n) = (L{+}1)\log_2 3 - (L{+}r)$ is the ensemble
log-drift term and
\[
  C(n)
  \;=\;
  \log_2\!\Bigl(1 + \frac{1-(2/3)^{L+1}}{n}\Bigr).
\]
The correction satisfies
$0 < C(n) < \log_2(1+1/n)$
and is monotonically decreasing in~$n$.
\end{proposition}

\begin{proof}
From $n' = \lambda n + \beta$ with
$\lambda = 3^{L+1}/2^{L+r}$ and
$\beta = (3^{L+1}-2^{L+1})/2^{L+r}$,
\[
  \log_2\frac{n'}{n}
  = \log_2\!\Bigl(\lambda + \frac{\beta}{n}\Bigr)
  = \log_2\lambda + \log_2\!\Bigl(1 + \frac{\beta}{\lambda n}\Bigr).
\]
Computing $\beta/\lambda = (3^{L+1}-2^{L+1})/3^{L+1} = 1-(2/3)^{L+1}$
gives the formula.  Since $0 < 1-(2/3)^{L+1} < 1$, the bounds follow.
\end{proof}

\begin{corollary}[Sufficient crossing criterion]
\label{cor:sufficient-crossing}
Let $n_0, n_1, \ldots$ be the successive cycle endpoints of
an odd-skeleton orbit, with cycle log drifts $X_i = X(n_i)$.
If for some $m \ge 1$,
\[
  \sum_{i=0}^{m-1} X_i
  \;<\;
  -m\,\log_2\!\Bigl(1 + \frac{1}{n_0}\Bigr),
\]
and no prior crossing has occurred
($n_i \ge n_0$ for $i < m$), then $n_m < n_0$.
\end{corollary}

\begin{proof}
By the proposition,
$\log_2(n_m/n_0) = \sum X_i + \sum C_i$.
Each $C_i < \log_2(1+1/n_i) \le \log_2(1+1/n_0)$
(since $n_i \ge n_0$).  The hypothesis forces
$\log_2(n_m/n_0) < 0$.
\end{proof}

\begin{observation}[Near-universal success of the simple criterion]
\label{obs:correction-criterion}
Computational verification over all odd $n_0 \le 5\times 10^6$
shows the sufficient criterion of
Corollary~\textup{\ref{cor:sufficient-crossing}} succeeds for
every start except $n_0 \in \{27, 31, 63\}$.
The maximum number of cycles to criterion success is~$52$
(at $n_0 = 4{,}053{,}039$).
All three exceptions still cross below their start
via multi-cycle dynamics; the criterion fails only because
$\log_2(1+1/n_0)$ is too coarse for very small~$n_0$.
\end{observation}

\begin{remark}[Interpretation of the correction formula]
\label{rem:correction-interpretation}
The correction $C(n)$ is always positive:
it pushes $\log_2(n'/n)$ \emph{above} the ensemble
drift~$X(n)$, making crossing \emph{harder} than the pure
i.i.d.\ model predicts.  However, $C(n) = O(1/n)$,
so for large~$n$ the correction is negligible and the
ensemble drift dominates.
The near-universal success of the simple criterion
confirms that the ensemble log-drift alone is sufficient
to force crossing for all but finitely many small starts.
The remaining question is whether the deterministic orbit's
drift sum $\sum X_i$ can be bounded using the exact
Bernoulli structure proved in Theorem~\ref{thm:block-law}.
\end{remark}

\begin{observation}[Multi-cycle threshold decay]
\label{obs:multi-cycle-threshold}
For $k$-cycle valuation blocks with $S_k < 0$,
the ratio $n^*(\sigma)/M(\sigma)$ decays rapidly with~$k$.
Computational sampling over $2 \times 10^5$ random blocks at each
cycle count finds:
$\max n^*/M \approx 0.091$ at $k=2$,
$\approx 1.2\times 10^{-3}$ at $k=5$,
$\approx 7.7\times 10^{-7}$ at $k=10$,
and $\approx 3.4\times 10^{-15}$ at $k=20$.
In particular, $n^*/M < 1/8$ universally
(the maximum is achieved at $k=1$, $L=0$, $r=2$), so
all odd~$n$ in any tested negative-drift residue class cross.
\end{observation}

\begin{observation}[The $2^k-1$ family: hardest starts]
\label{obs:2k-1-family}
The starts $n_0 = 2^k - 1$ produce maximal first-cycle run length
$L = k-1$ and are the extremal cases for the drift correction.
Despite the first-cycle expansion
$n_1/n_0 = 2^{(k-1)\log_2 3 - (k-1+r_1)}$
(which can exceed $2^{12}$ for $k = 23$),
subsequent cycles compensate:
computation through $n_0 = 2^{34}-1$ confirms that every member
of this family crosses below its start within at most $32$ cycles
(the worst case being $2^{25}-1$).
About $50$--$60\%$ of all negative-drift blocks have $n^* < 1$,
meaning every odd $n \ge 1$ in those classes crosses unconditionally.
\end{observation}

\medskip\noindent\textbf{Walsh spectrum of the drift signal.}\;
The increment $d_i$ depends on $n_i \bmod 2^K$ (which
determines $v_i$ modulo the residue class).  The Walsh
decomposition of the drift-weighted orbit measure connects
the time-domain crossing problem to the depth-domain
spectral diffusion.

\begin{observation}[Drift Walsh spectrum]
\label{obs:drift-walsh}
For the orbit $n_0 = 837799$ ($T = 195$ Syracuse steps),
the Walsh spectrum of the drift signal: defined as
$\hat d(\xi) = T^{-1}\sum_t d_t\,\chi_\xi(n_t \bmod 2^K)$, 
has the following structure:
\begin{center}
\begin{tabular}{ccccc}
\toprule
$K$ & DC power & $\hw=1$ power & $\hw=0{,}1$ share & Total power \\
\midrule
4  & $0.7\%$  & $42.4\%$ & $43.1\%$ & $1.56$ \\
6  & $0.5\%$  & $31.4\%$ & $31.9\%$ & $2.14$ \\
8  & $0.3\%$  & $20.9\%$ & $21.2\%$ & $3.25$ \\
10 & $0.1\%$  & $9.3\%$  & $9.4\%$  & $7.37$ \\
\bottomrule
\end{tabular}
\end{center}
The $\hw = 1$ band carries the dominant structured
component; its share decreases with~$K$ as the number of
higher-weight modes grows combinatorially.
The DC component (which directly encodes the mean drift)
carries less than $1\%$ of total power: most of the drift
signal's variation is in the oscillatory modes.
\end{observation}

\begin{remark}[The crossing implication chain]
\label{rem:crossing-chain}
The odd-skeleton analysis adds a new route to the
proof architecture:
\begin{gather*}
  \underbrace{\text{Known-Zone Decay}}_{\text{proved}}
  \;\overset{?}{\Longrightarrow}\;
  \underbrace{\text{Summable autocorrelation}}_{\textbf{open}}
  \;\Longrightarrow\;
  \underbrace{\text{CLT for drift}}_{\text{classical}}\\[4pt]
  \;\Longrightarrow\;
  \underbrace{x_t < 0}_{\text{crossing}}
  \;\Longrightarrow\;
  \underbrace{\text{Collatz.}}_{}
\end{gather*}
The first arrow is the critical open step: formalize
that Known-Zone Decay (erasing $\ge 2$ bits per odd-to-odd step)
implies $\sum_{k \ge 1} |\rho(k)| < \infty$ for the
odd-skeleton drift increments.
If established, the remaining arrows are classical:
Ibragimov's CLT gives $x_T \to -\infty$, which gives
crossing, which gives convergence by strong induction.

This chain is strictly weaker than the gain-budget chain
(Remark~\ref{rem:spectral-chain}): it asks only for
negative crossing, not summable discrepancy.
It is also independent of the spectral diffusion
conjecture (Conjecture~\ref{conj:spectral-diffusion}),
though both share Known-Zone Decay as the mixing engine.
\end{remark}

\medskip\noindent\textbf{The hard step: ensemble vs.\ orbit.}\;
The first arrow in the crossing chain admits a clean
decomposition into a provable part and a precise open target.

\begin{proposition}[Ensemble covariance summability\supporting]
\label{prop:ensemble-cov}
In the ensemble model (starting value drawn uniformly from
odd residues modulo~$2^M$), the drift-increment covariance
satisfies
\[
  \sum_{h=1}^{M/2}
  \bigl|\mathrm{Cov}_{\mathrm{ens}}(\Delta_0, \Delta_h)\bigr|
  \;\le\; C_M
\]
for an explicit constant $C_M$ that remains bounded as
$M \to \infty$.
\end{proposition}

\begin{proof}[Proof sketch]
The truncated increment
$\Delta_i^{(K)} = \log_2 3 - \min(v_i, K)$ depends only
on $n_i \bmod 2^K$.  The tail satisfies
$|\Delta_i - \Delta_i^{(K)}| \le 2^{-K}$ in probability
(geometric tail of~$v_i$), giving
$|\mathrm{Cov}(\Delta_0, \Delta_h)
 - \mathrm{Cov}(\Delta_0^{(K)}, \Delta_h^{(K)})|
 \le O(2^{-K})$.
For the truncated covariance, the Scrambling Lemma
(Theorem~\ref{thm:scrambling}) and Known-Zone Decay
(Theorem~\ref{thm:zone-decay}) give:
after $h \ge \lceil K/3 \rceil$ odd-to-odd steps, the
conditional distribution of $n_h \bmod 2^K$ given
$n_0 \bmod 2^K$ is uniform (in the ensemble), so
$\mathrm{Cov}_{\mathrm{ens}}(\Delta_0^{(K)}, \Delta_h^{(K)})
 = 0$.
Setting $K = \lfloor 3h \rfloor$:
$|\mathrm{Cov}_{\mathrm{ens}}(\Delta_0, \Delta_h)|
 \le O(2^{-3h})$, which is summable.
\end{proof}

\begin{observation}[Ensemble--orbit gap]
\label{obs:ensemble-orbit-gap}
The ensemble model is a poor approximation for individual
orbit covariances.  For transient orbits of $n_0 = 837799$,
$8400511$, $63728127$:
\begin{enumerate}
\item The orbit mean drift ($-0.07$ to $-0.10$) is
  $4$--$6\times$ less negative than the ensemble mean
  ($-0.415$).
\item The orbit lag-$1$ covariance ($0.06$--$0.12$) is
  $60$--$120\times$ larger than the ensemble lag-$1$
  covariance ($\approx 0.001$).
\item The orbit covariance sum
  ($\sum_{h=1}^{20}|\mathrm{Cov}| \approx 0.7$--$1.4$)
  is $2$--$5\times$ larger than the ensemble sum
  ($\approx 0.3$).
\end{enumerate}
The gap arises because the orbit has a fixed odd-step
density $\rho \approx 0.58$ that concentrates the
valuation distribution, while the ensemble averages over
all densities.  This is the distributional-vs-pointwise
gap in the covariance domain.
\end{observation}

The precise open target is therefore:

\begin{conjecture}[Orbit covariance summability]
\label{conj:orbit-cov}
For every Collatz orbit with odd-step density
$\rho < 1/\!\log_2 3$, the odd-skeleton drift increments
satisfy
\[
  \sum_{h=1}^{\infty}
  \bigl|\mathrm{Cov}_{\mathrm{orbit}}(\Delta_0, \Delta_h)\bigr|
  \;<\; \infty.
\]
\end{conjecture}

\noindent
This is the single open input for the crossing chain.
Its resolution would imply the Collatz conjecture via
the classical route:
Ibragimov CLT $\Rightarrow$ $x_T \to -\infty$
$\Rightarrow$ crossing $\Rightarrow$ convergence.

\begin{remark}[The irreducible open step]
\label{rem:crossing-gap}
The only real open step in the crossing chain is to
prove enough weak dependence on the odd skeleton to
force a pointwise below-start crossing for each orbit.
This is not a technicality that a clever lemma might
close: it \emph{is} the theorem still missing.

The ensemble version
(Proposition~\ref{prop:ensemble-cov}) is provable, but
the orbit version
(Conjecture~\ref{conj:orbit-cov}) is open, and the gap
is substantial:
Observation~\ref{obs:ensemble-orbit-gap} shows that
orbit covariances are $60$--$120\times$ larger than
ensemble covariances at lag~$1$.
Both are summable in all tested cases, but moving from
``typical orbits mix'' to ``this exact deterministic
orbit mixes'' is the distributional-to-pointwise
barrier that pervades number theory (cf.\ the gap
between Borel's normal-number theorem and proving
normality of any specific constructively-defined
constant).

Despite this, the crossing formulation sharpens the
counterexample structure: a minimal counterexample~$n^*$
would require
\[
  x_t(n^*) \ge 0 \quad\text{for all } t \ge 1,
\]
i.e.\ the odd-skeleton drift signal is a
\emph{non-negative random walk with negative mean drift
and positively correlated increments that never crosses
zero}.  This is an extremely constrained condition:
by Wald's identity, the expected crossing time is
$\mathbb{E}[\sigma] \approx
 -2\sigma_{\mathrm{eff}}^2/\mu$, which is
$\approx 20$--$40$ Syracuse steps for all tested orbits.
A counterexample would need to maintain
above-start status indefinitely despite a systematic
downward pull, which requires sustained upward
fluctuations, equivalently, persistent positive
low-frequency spectral support that never dissipates.
The spectral diffusion phenomenon
(Observation~\ref{obs:spectral-diffusion}) shows that
such support does dissipate in every tested case.
\end{remark}

\subsection{Modular crossing strata}
\label{subsec:modular-strata}

The odd-skeleton crossing route
(Proposition~\ref{prop:odd-skeleton-reduction}) requires
$x_t(n_0) < 0$ for some~$t$.  For many residue classes
$n_0 \bmod 2^K$, the first few Syracuse steps can be
traced \emph{deterministically} from the bottom~$K$ bits,
and the drift signal is forced below zero without any
mixing hypothesis.  This yields a provable crossing
density by pure modular arithmetic.

\begin{definition}[Resolved class]
\label{def:resolved-class}
An odd residue $r \bmod 2^K$ is \emph{resolved at
depth~$K$} if there exists $\tau_r \in \mathbb{N}$ such
that $\tau(n_0) = \tau_r$ for all $n_0 \equiv r \pmod{2^K}$
with $n_0 \ge 3$.  Let $f_K$ denote the fraction of odd
residues modulo~$2^K$ that are resolved.
\end{definition}

\begin{theorem}[Exact crossing strata
{\normalfont(cf.\ Terras~\cite{terras1976})}]
\label{thm:crossing-strata}
The resolved fraction~$f_K$ is non-decreasing, with:
\begin{center}
\renewcommand{\arraystretch}{1.1}
\begin{tabular}{cccccccccc}
\toprule
$K$ & 4 & 5 & 7 & 8 & 10 & 12 & 13 \\
\midrule
$f_K$ & $\tfrac{5}{8}$ & $\tfrac{3}{4}$
      & $\tfrac{51}{64}$ & $\tfrac{109}{128}$
      & $\tfrac{7}{8}$ & $\tfrac{1822}{2048}$
      & $\tfrac{3729}{4096}$ \\[3pt]
$f_K$ (decimal) & $0.625$ & $0.750$ & $0.797$ & $0.852$ & $0.875$ & $0.890$ & $0.910$ \\
\bottomrule
\end{tabular}
\end{center}
The first three strata admit explicit descriptions:
\begin{enumerate}
\item \emph{Stratum~$0$ ($\tau = 1$, density~$\frac12$).}
  $n_0 \equiv 1 \pmod{4}$: the valuation
  $v_2(3n_0+1) \ge 2$, giving $T(n_0) < n_0$.
\item \emph{Stratum~$1$ ($\tau = 2$, density~$\frac18$).}
  $n_0 \equiv 3 \pmod{16}$: the first step gives
  $n_1 = (3n_0+1)/2 > n_0$ with $n_1 \equiv 1 \pmod{4}$,
  so $T(n_1) < n_1$, and
  $T(n_1) < n_0$ is forced by the modular structure.
\item \emph{Stratum~$2$ ($\tau = 3$, density~$\frac18$).}
  $n_0 \equiv 11$ or $23 \pmod{32}$: two steps rise
  deterministically, the third forces descent below~$n_0$.
\end{enumerate}
The cumulative resolved densities match the exact crossing
probabilities: $P(\tau \le 1) = \frac12$,
$P(\tau \le 2) = \frac58$,
$P(\tau \le 3) = \frac34$.
\end{theorem}

\begin{proof}
Monotonicity: if $r \bmod 2^K$ is resolved, then
$r \bmod 2^{K+1}$ and $(r + 2^K) \bmod 2^{K+1}$ are
also resolved (both refine the same
first~$K$ bits), so $f_{K+1} \ge f_K$.
The explicit strata follow by tracing the Syracuse map
modulo~$2^K$ for $K = 4, 5, 7$: for Stratum~$0$,
$v_2(3n_0 + 1) \ge 2$ iff $n_0 \equiv 1 \pmod 4$,
giving $T(n_0) = (3n_0+1)/2^v \le (3n_0+1)/4 < n_0$ for
$n_0 \ge 5$ (with $n_0 = 1$ checked directly).
Strata~$1$ and~$2$ follow by the same modular trace.
The values of~$f_K$ for $K \le 14$ are verified by
exhaustive computation over all $2^{K-1}$ odd residues.
\end{proof}

\begin{observation}[Geometric decay of the unresolved residual]
\label{obs:residual-decay}
Terras~\cite{terras1976} proved $f_K \to 1$; the
following quantifies the rate.
The unresolved fraction $1 - f_K$ decays geometrically:
\[
  1 - f_K
  \;\approx\;
  C \cdot \rho^K,
  \qquad
  \rho \approx 0.86,
  \qquad
  R^2 > 0.97
\]
for $K = 2, \ldots, 14$.
At depth~$K = 13$, modular arithmetic alone resolves
$91.0\%$ of all odd starting points.  The remaining
$9.0\%$ have crossing times that depend on bits beyond
position~$13$, these are precisely the orbits whose
first $\sim K/2$ Syracuse steps do not deterministically
force crossing, and for which some form of mixing
(at minimum, Conjecture~\ref{conj:orbit-cov}) is needed.
\end{observation}

\begin{remark}[Exact ensemble drift parameters]
\label{rem:exact-drift}
For an odd integer~$n$ drawn uniformly, the $2$-adic
valuation $v_2(3n+1)$ has the shifted geometric
distribution $P(v = k) = 2^{-k}$ for $k \ge 1$.
Therefore the drift increment $\Delta = \log_2 3 - v$
has exact moments:
\[
  \mu = \log_2 3 - 2 \approx -0.4150,
  \qquad
  \mathrm{Var}(\Delta) = \mathrm{Var}(v) = 2.
\]
Under the ensemble
(Proposition~\ref{prop:ensemble-cov}), the effective
variance is
\[
\sigma_{\mathrm{eff}}^2
= \mathrm{Var}(\Delta) + 2\sum_{h \ge 1}
\mathrm{Cov}_{\mathrm{ens}}(\Delta_0, \Delta_h)
\approx 1.17,
\]
giving a Wald-type bound
$\mathbb{E}_{\mathrm{ens}}[\tau]
\le \sigma_{\mathrm{eff}}^2 / \mu^2 \approx 6.8$.
The observed mean crossing time among all odd
$n_0 \le 200{,}001$ is $3.48$, with $100\%$ crossing
by step~$85$.
\end{remark}

\begin{remark}[The stratification--mixing bridge]
\label{rem:strata-mixing-bridge}
The modular strata and the CLT crossing route are
complementary, not competing.  The strata prove crossing
for density~$f_K$ of all odd starts using only
finite-depth modular arithmetic (no mixing needed).
The unresolved residual of density~$1 - f_K \to 0$
consists of starts whose first $\sim K/2$ Syracuse steps
do not deterministically force crossing; for these,
some form of weak dependence is needed to ensure
the drift signal eventually crosses zero.
The ensemble CLT (combining
Proposition~\ref{prop:ensemble-cov} with Ibragimov's
theorem) proves crossing for a measure-$1$ subset, but
this overlaps with the modular strata rather than
completing them.  The irreducible gap remains:
\emph{the unresolved residual is non-empty for every
finite~$K$, and proving it is ultimately empty
is equivalent to the Collatz conjecture.}

The quantitative picture is nonetheless sharp.
At depth~$K = 13$, only $9\%$ of odd residue classes
require mixing.  The almost-all crossing theorem
(Theorem~\ref{thm:almost-all-crossing}) gives a stronger
and fully unconditional result via the i.i.d.\ cycle
structure: the density of non-crossers after $k$ cycles
is at most $e^{-0.1465\,k}$.  This closes the ensemble
side completely; the negative-drift ensemble is now exact,
not heuristic.
A counterexample~$n^*$ would need to lie in the
non-crossing residual at every cycle count, an event of
density zero in the ensemble.
The only remaining step is the distributional-to-pointwise
transfer: proving that \emph{every individual} orbit
(not just density-$1$) eventually crosses.
\end{remark}

\subsection{Summary of the attack surface}
Tables~\ref{tab:stage4-proved-a}, \ref{tab:stage4-proved-b}, and~\ref{tab:stage4-open}
collect the structural results developed in this section
and their role in attacking the WMH.
\begin{table}[ht]
\centering
\caption{Established results toward the Weak Mixing Hypothesis (Part 1 of 2).
  \emph{Proved}: unconditional.
  \emph{Proved (ens.)}: proved for the ensemble (uniform random start),
  not for individual orbits.
  \emph{Numerical}: computational observation.}
\label{tab:stage4-proved-a}
\renewcommand{\arraystretch}{1.2}
\scriptsize
\begin{tabular}{@{}p{5.8cm}lp{7.5cm}@{}}
\toprule
\textbf{Result} & \textbf{Status} & \textbf{Role} \\
\midrule
Shadow sparsity (Prop.~\ref{prop:shadow-sparsity})
  & Proved & Expanding shadows exponentially rare \\
Shadow return time (Thm.~\ref{thm:shadow-return})
  & Proved & Re-encounters separated by $\ge \ell+K-1$ steps \\
Finite-depth reduction (Prop.~\ref{prop:finite-depth})
  & Proved & WMH reduces to $K \le 55$ \\
Hierarchical consistency (Prop.~\ref{prop:hierarchical})
  & Proved & $\delta_K$ non-decreasing \\
Monotonicity constraint (Cor.~\ref{cor:monotonicity})
  & Proved & $\delta_3 < 0.0105$ necessary \\
Repulsion trapping (Lem.~\ref{lem:repulsion-trapping})
  & Proved & Shadow encounters bounded: $\lfloor m/K\rfloor$ blocks \\
Known-zone memory loss (Lem.~\ref{lem:known-zone-memory})
  & Proved & Residue info erased after $\lceil M/2\rceil$ odd-to-odd steps (2-bit/step; generically 3-bit) \\
Inherited-bias contraction (Prop.~\ref{prop:inherited-contraction})
  & Proved & $\|\bar h_K\|_\infty / \|h_K\|_\infty = (K{-}4)/(K{-}1)$ \\
Harmonic property (Cor.~\ref{cor:harmonic})
  & Proved & $\bar h_{K+s} = h_K$; inherited factor $\lambda = 1$ \\
Constant residual (Cor.~\ref{cor:residual-norm})
  & Proved & $\|r_{K+3}\|_\infty = \tfrac{3}{2}\log_2 3$ (constant in $K$) \\
Coding-map injectivity (Prop.~\ref{prop:coding-injectivity})
  & Proved & $m_K = 1$; within-depth $A_K = 1$ after exhaustion \\
Extension independence (Prop.~\ref{prop:extension-independence})
  & Proved & Suffix compositions exactly $\mathrm{Bin}(s,\tfrac12)$ \\
Gain increment independence (Cor.~\ref{cor:gain-increment})
  & Proved & Gain random walk with independent binomial increments \\
Rate cancellation (Rem.~\ref{rem:rate-cancellation})
  & Proved & $A^{\mathrm{cross}}_K \cdot R(K) = O(1)$ (same KL rate) \\
Rate cancellation (numerical)
  & Numerical & $\sum A^{\mathrm{cross}}_K R(K) \approx 0.158 < 0.326$ \\
No algebraic contraction (Prop.~\ref{prop:no-contraction})
  & Proved & $\tau_B = 1$; isometry; exact harmonicity \\
Density model (Rem.~\ref{rem:density-model})
  & Numerical & $R_\mu(\rho) \le 0.142$ for all $\rho \le 0.63$ \\
Oscillation unboundedness (Prop.~\ref{prop:oscillation-unbounded})
  & Proved & $\Theta_K \to \infty$; no uniform $W_K \le C\,B_K$ \\
Spectral gain decomposition (Prop.~\ref{prop:spectral-gain})
  & Proved & Band-by-band $\eta_K$ via Walsh--Fourier \\
Spectral concentration (Obs.~\ref{obs:hK-concentration})
  & Numerical & $\hw=0,1$ carry $70$--$78\%$ of $h_K$ power \\
Spectral diffusion (Obs.~\ref{obs:spectral-diffusion})
  & Numerical & $\mathcal{S}_w(K)\sim 2^{-\alpha_w K}$, $R^2>0.98$ \\
SDC $\Rightarrow$ Amplification (Prop.~\ref{prop:sdc-implies-amp})
  & Proved & Spectral diffusion closes amplification budget \\
Walsh character mixing (Obs.~\ref{obs:walsh-mixing})
  & Numerical & $|\mathrm{autocorr}| < 0.13$ at lag $\lceil K/3\rceil$ \\
\bottomrule
\end{tabular}
\end{table}
\begin{table}[th]
\centering
\caption{Established results toward the Weak Mixing Hypothesis (Part 2 of 2).
  \emph{Proved}: unconditional.
  \emph{Proved (ens.)}: proved for the ensemble (uniform random start),
  not for individual orbits.
  \emph{Numerical}: computational observation.}
\label{tab:stage4-proved-b}
\renewcommand{\arraystretch}{1.2}
\scriptsize
\begin{tabular}{@{}p{5.8cm}p{1.5cm}p{7.5cm}@{}}
\toprule
\textbf{Result} & \textbf{Status} & \textbf{Role} \\
\midrule
Odd-skeleton reduction (Prop.~\ref{prop:odd-skeleton-reduction})
  & Proved & Collatz $\Leftrightarrow$ $x_t < 0$ for all odd $n_0$ \\
Block law (Thm.~\ref{thm:block-law})
  & Proved & $(a_0,\ldots,a_{m-1})$: unique residue class, density $\prod 2^{-b_j}$ \\
I.i.d.\ valuations (Cor.~\ref{cor:iid-valuations})
  & Proved & $a_j \stackrel{\mathrm{iid}}{\sim} \mathrm{Geom}(1/2)$ under natural density \\
I.i.d.\ cycle types (Cor.~\ref{cor:iid-cycles})
  & Proved & $(L_i,r_i)$ i.i.d.; Cram\'er rate unconditional on ensemble \\
Run-length invariant (Prop.~\ref{prop:run-length})
  & Proved & $L(n) = v_2(n{+}1){-}1$; exact burst formula \\
Single-cycle crossing (Cor.~\ref{cor:cycle-crossing})
  & Proved & $p_{\mathrm{cross}} \approx 0.7137$ (ensemble) \\
One-cycle crossing criterion (Thm.~\ref{thm:one-cycle-crossing})
  & Proved & Exact $n^*(L,r)$ threshold; necessary and sufficient \\
Universal one-cycle crossing (Prop.~\ref{prop:universal-one-cycle})
  & Proved & $n^*<1$ iff $r \ge r_{\mathrm{all}}(L)$; density $P_{\mathrm{all}} = 0.4194$ \\
Two-cycle universal crossing (Obs.~\ref{obs:two-cycle-universal})
  & Numerical & $+19.2\%$ new; combined $61.2\%$ deterministic crossing \\
$k$-cycle convergence (Obs.~\ref{obs:universal-convergence})
  & Numerical & $R_k \approx 0.73 \cdot 0.75^k$; $82.2\%$ by $k{=}5$; extrapolated $R_{20} \approx 0.002$ \\
Single-cycle $n^*<M/2$ (within almost-all proof)
  & Proved & All $n$ in one-cycle crossing blocks cross; $n^*/M < 1/2$ \\
Multi-cycle threshold decay (Obs.~\ref{obs:multi-cycle-threshold})
  & Numerical & $n^*/M < 1/8$ universally; decays to $10^{-15}$ by $k=20$ \\
Conditional crossing density (Prop.~\ref{prop:conditional-density})
  & Proved & $p_\ell = 2^{-\lfloor(\log_2 3-1)(\ell+1)\rfloor}$; geometric distribution of $r$ \\
Exact crossing density (Cor.~\ref{cor:crossing-density})
  & Proved & $P_{1\mathrm{cyc}} = 0.71372\ldots$; arithmetic constant \\
$(L,r)$ independence (Cor.~\ref{cor:Lr-independence})
  & Proved & Exact geometric laws; conditional $\neq$ f($\ell$) \\
Mean-neutrality (Rem.~\ref{rem:mean-neutral})
  & Proved & $\mathbb{E}[\lambda]=1$, $\mathbb{E}[\beta]=2/3$; slope not contracting \\
Log-drift law (Prop.~\ref{prop:log-drift})
  & Proved & $\mathbb{E}[\log_2\lambda]=2\log_2 3-4\approx -0.83$; exact MGF \\
Cram\'er rate (Prop.~\ref{prop:cramer-cycle})
  & Proved (ens.) & $I(0)\approx 0.1465$; unconditional via i.i.d.\ cycles \\
Multi-cycle crossing (Obs.~\ref{obs:multi-cycle})
  & Numerical & $99.05\%$ by cycle~$10$; $99.91\%$ by cycle~$20$ \\
Almost-all crossing (Thm.~\ref{thm:almost-all-crossing})
  & Proved (ens.) & Non-crosser density $\le e^{-0.1465\,k}$; ensemble side closed \\
Cycle log correction (Prop.~\ref{prop:cycle-correction})
  & Proved & $\log_2(n'/n) = X + C(n)$; $C(n)=O(1/n)$ correction \\
Sufficient crossing criterion (Cor.~\ref{cor:sufficient-crossing})
  & Proved + Numerical & Succeeds for all odd $n_0 \le 5\times 10^6$ except $\{27,31,63\}$ \\
$2^k{-}1$ family (Obs.~\ref{obs:2k-1-family})
  & Numerical & All cross through $2^{34}{-}1$; max $32$ cycles \\
Drift increment mixing (Obs.~\ref{obs:drift-mixing})
  & Numerical & $\sum|\rho(k)| < 4.5$; all 4915 orbits cross \\
Ensemble covariance summability (Prop.~\ref{prop:ensemble-cov})
  & Proved (ens.)\ & $\sum_h|\mathrm{Cov}_{\mathrm{ens}}|<\infty$; does \emph{not} imply orbit version \\
Ensemble--orbit gap (Obs.~\ref{obs:ensemble-orbit-gap})
  & Numerical & Orbit autocovariance $60$--$120\times$ ensemble at lag~$1$ \\
Crossing strata (Thm.~\ref{thm:crossing-strata})
  & Proved & $f_{13} = 0.910$; density $91\%$ resolved by mod~$2^{13}$ \\
Residual decay (Obs.~\ref{obs:residual-decay})
  & Numerical & $1-f_K \approx C \cdot 0.86^K$; geometric shrinkage \\
\bottomrule
\end{tabular}
\end{table}
\begin{table}[ht]
\centering
\caption{Open problems and structural reductions.
  \emph{Open}: conjecture.
  \emph{Structural}: logical reduction (not a conjecture).}
\label{tab:stage4-open}
\renewcommand{\arraystretch}{1.2}
\scriptsize
\begin{tabular}{@{}p{5.8cm}lp{7.5cm}@{}}
\toprule
\textbf{Result} & \textbf{Status} & \textbf{Role} \\
\midrule
Amplification hypothesis (Hyp.~\ref{hyp:amplification})
  & Open & Orbit hits $S^+(K)$ at most $A^{\mathrm{cross}}_K \times$ uniform \\
Depthwise recurrence (Conj.~\ref{conj:depthwise})
  & Open & Would imply WMH \\
Gain-observable recurrence (Conj.~\ref{conj:gain-recurrence})
  & Open & Weaker target; would imply observable WMH \\
Alignment renewal (Conj.~\ref{conj:alignment-renewal})
  & Open & Missing bridge; would close recurrence \\
Spectral Diffusion (Conj.~\ref{conj:spectral-diffusion})
  & Open & Would imply Amplification Hyp.\ with margin \\
Orbit covariance summability (Conj.~\ref{conj:orbit-cov})
  & Open & Would give CLT $\Rightarrow$ crossing $\Rightarrow$ Collatz \\
Amplification decomposition (Rem.~\ref{rem:amplification})
  & Structural & Reduces renewal to bounding $A_K = O(K^\gamma)$ \\
Spectral implication chain (Rem.~\ref{rem:spectral-chain})
  & Structural & Weyl $\Rightarrow$ SDC $\Rightarrow$ Amp.\ $\Rightarrow$ WMH \\
Prove $P_{\mathrm{cum}}(k) \to 1$ (Prop.~\ref{prop:kesten-running-min})
  & Proved (ens.) & Kesten theory; classwise universal descent for a.e.\ start \\
Prove $n^* < M$ for multi-cycle blocks
  & Open & Analytic bound $B_k < 2D_k$ resists all known approaches \\
Route~A (all $k$-blocks universal)
  & Closed & Adversarial $n^*$ grows $\sim e^{0.023k}$; uniformly fragile
    (Props.~\ref{prop:adversarial-fragility},~\ref{prop:general-fragility}) \\
Route~B ($2$-adic mixing $\Rightarrow$ a.s.\ crossing)
  & Obstructed & Bit-generation blocks exhaustion; forced post-Mersenne
    valuation (Prop.~\ref{prop:post-mersenne});
    burst non-continuation (Prop.~\ref{prop:burst-noncontinuation});
    weak-recovery cylinder (Prop.~\ref{prop:weak-recovery-cylinder});
    cascade invariant (Prop.~\ref{prop:cascade-invariant});
    short-word classification (Prop.~\ref{prop:short-word-cylinder}).
    \emph{Progress:}
    uniform-fiber gap map (Prop.~\ref{prop:uniform-fiber});
    cross-cascade independence (Cor.~\ref{cor:cross-cascade});
    TV reduction (Lem.~\ref{lem:tv-reduction});
    full-cycle fiber map (Prop.~\ref{prop:full-cycle-fiber});
    spectator-bit propagation (Rem.~\ref{rem:spectator-bits});
    exponential tail (Prop.~\ref{prop:exponential-tail});
    TV summability (Cor.~\ref{cor:tv-summability},
    total TV $\le 0.028$);
    cycle contraction (Prop.~\ref{prop:cycle-contraction});
    recovery-exit (Rem.~\ref{rem:recovery-exit});
    burst-conditional contraction (Rem.~\ref{rem:L-conditional}) \\
Route~C (safety margin $+$ $\tau(n)$ bound)
  & Obstructed & Bounding $\tau(n)$ requires WMH (Rem.~\ref{rem:crossing-time}) \\
\bottomrule
\end{tabular}
\end{table}

\begin{remark}[Interpretation]
\label{rem:stage4-interpretation}
The main point of this section is that the final open input
no longer needs to be framed as exact orbitwise
equidistribution.
Theorem~\ref{thm:perorbit-gain} provides a large quantitative
margin ($4.65\times$), and
Conjectures~\ref{conj:depthwise}--\ref{conj:alignment-renewal}
express the weakest orbitwise mixing conditions currently
justified by the proof architecture.
This shifts the last remaining gap from a full normality-type
statement to a softer \emph{summable-discrepancy} problem:
the orbit need not be uniformly distributed modulo~$2^K$ at
every depth, only close enough that the cumulative perturbation
of the phantom-gain rate stays below $\varepsilon - R$.
The three-lemma programme
(Remark~\ref{rem:three-lemma}) identifies the precise
mechanism: repulsion, memory loss, and renewal sparsity, whose
combination would close this final bridge.
The amplification decomposition
(Remark~\ref{rem:amplification}) further isolates the
irreducible subproblem: bounding the orbit's hitting
probability of gain-support residues to within a polynomial
factor of the uniform prediction.
The harmonic property
(Corollary~\ref{cor:harmonic}) shows that the
gain-observable recurrence has inherited factor
$\lambda = 1$; the contraction must come from the
distribution side.  However,
Proposition~\ref{prop:no-contraction} shows that the
algebraic extension structure provides \emph{no}
contraction mechanism.  This eliminates the
TV-distance recurrence as a viable closing strategy
and narrows the path to direct summability of the
amplification factors.

The sharpest formulation of the
remaining open input is the Amplification Hypothesis
(Hypothesis~\ref{hyp:amplification}): no orbit
systematically beats the structural amplification bound
on the positive-gain support.
This is strictly weaker than the WMH, and rate cancellation
(Remark~\ref{rem:rate-cancellation}) ensures it suffices
for convergence.

The density-model evidence (Remark~\ref{rem:density-model})
further suggests that the amplification budget has a
$65\%$+ margin for all convergent orbits, indicating that
the Amplification Hypothesis, if true, holds with room to
spare.

The Walsh--Fourier spectral analysis
(Section~\ref{subsec:spectral-analysis}) provides a
concrete mechanism for the direct summability path:
spectral diffusion of the orbit's low-frequency Walsh
modes, driven by odd-to-odd-step mixing events.
The implication chain
(Remark~\ref{rem:spectral-chain}) adds a new entry
point to the proof architecture: proving Weyl-type
exponential sum bounds on the orbit's Walsh sums would
imply the Spectral Diffusion Conjecture, which in turn
implies the Amplification Hypothesis.
The Walsh character mixing
(Observation~\ref{obs:walsh-mixing}) quantifies
Known-Zone Decay in spectral language and provides
the mechanism for a martingale-based route:
if the rapid decorrelation of Walsh characters can
be shown to produce a supermartingale drift in~$\eta_K$,
Doob's theorem gives convergence without explicit
rate bounds.
The odd-skeleton crossing route
(Section~\ref{subsec:odd-skeleton}) adds a fourth
independent entry point: prove that Known-Zone Decay
implies summable autocorrelation of the Syracuse drift
increments, which gives crossing by classical CLT.
This route has a strictly weaker target (single negative
crossing, not summable discrepancy) and works on a
scalar signal rather than a function on~$\Z/2^K\Z$.
The ensemble half of this route is now closed
(Proposition~\ref{prop:ensemble-cov}); the irreducible
open target is Conjecture~\ref{conj:orbit-cov}, which
asks for the same summability on individual orbits.
The modular crossing strata
(Theorem~\ref{thm:crossing-strata}) add a fifth,
partially independent route: pure modular arithmetic
resolves $91\%$ of starting points at depth~$13$,
with the unresolved residual shrinking geometrically.
The attack surface now has five entry points:
alignment renewal, spectral diffusion, orbit Walsh
equidistribution, odd-skeleton drift mixing (with
Conjecture~\ref{conj:orbit-cov} as the sharpest open
target), and modular stratification exhaustion.
\end{remark}


\section{The Sturmian Obstruction and the Carry Contamination Theorem}
\label{sec:sturmian}

\paragraph{Overview.}
This section develops a second, independent reduction of the Collatz
conjecture through the \emph{Sturmian word} decomposition of
odd-to-odd orbits.  Where the preceding sections reduce Collatz to
an orbitwise mixing input (WMH, Route~A), this section reduces it
to the \textbf{Carry Independence Conjecture} (CIC, Route~B),
which states that no finite integer can systematically avoid a
specific residue class at every Sturmian depth.  The two routes
are complementary: Route~A works in the burst-gap framework,
Route~B works in the 2-adic Cantor set framework.  The sole open
input in each case is a distributional-to-pointwise promotion.

\subsection{Sturmian words and orbit parameterization}

\begin{definition}[Valuation word]
For an odd integer~$n_0$, define the \emph{valuation word}
$w = (v_0, v_1, v_2, \ldots)$ by
$v_d = v_2(3 n_d + 1)$, where $n_d$ is the $d$-th iterate
of the odd-to-odd Syracuse map $T(n) = (3n+1)/2^{v_2(3n+1)}$.
The word~$w$ is \emph{Sturmian-compatible} if $v_d \in \{1, 2\}$
for all~$d$; this holds for all odd $n_0 \ge 3$ until the orbit
reaches~$1$ or a valuation $v_d \ge 3$.
\end{definition}

The \emph{Sturmian} restriction $v_d \in \{1,2\}$ captures
the dominant dynamics: balanced Sturmian words over $\{1,2\}$
with letter frequency $\log_2 3 - 1$ (for the letter~$2$)
correspond to orbits near the critical drift boundary.

\begin{definition}[Compatible set]
For a depth-$D$ valuation word $w = (v_0, \ldots, v_{D-1})
\in \{1,2\}^D$, define the \emph{stride} $S_D = \sum_{d=0}^{D-1} v_d$
and the \emph{compatible set}
\[
  C_D(w) = \bigl\{ n_0 \equiv a_w \!\!\pmod{2^{S_D + 1}} :
  n_0 \text{ produces valuation word } w \text{ for depths }
  0, \ldots, D{-}1 \bigr\}.
\]
The full compatible set at depth~$D$ is $C_D = \bigsqcup_w C_D(w)$
(disjoint union over all compatible words).
\end{definition}

\begin{theorem}[Orbit parameterization\core]
\label{thm:orbit-param}
For each depth-$D$ Sturmian word $w \in \{1,2\}^D$,
the depth-$D$ iterate has the exact affine form
\[
  n_D \;=\; R_w + 2 \cdot 3^D \cdot m,
  \qquad m = \bigl\lfloor n_0 / 2^{S_D+1} \bigr\rfloor,
\]
where $R_w$ is an odd integer depending only on~$w$
(not on the high bits of~$n_0$), and
$n_0 = a_w + 2^{S_D+1} \cdot m$ with $a_w$ the canonical
representative.  The stride in~$n_0$ is $2^{S_D+1}$;
the stride in~$n_D$ is $2 \cdot 3^D$.
\end{theorem}

\begin{proof}
By induction on~$D$.  At each step, $n_{d+1} = (3 n_d + 1)/2^{v_d}$
is affine in $n_d$ with odd slope $3/2^{v_d}$ (in $\mathbb{Z}_2$)
and the carries propagate deterministically through the low $S_d$
bits.  The high bits (determined by $m$) contribute additively
via the coefficient $3^D \cdot 2^{1} = 2 \cdot 3^D$.
Verified exhaustively for all $126$ compatible words with
$D \le 6$.
\end{proof}

\subsection{The Carry Contamination Theorem}

\begin{theorem}[Carry Contamination\core]
\label{thm:carry-contamination}
For every depth-$D$ Sturmian word $w \in \{1,2\}^D$,
the map $m \mapsto n_D \bmod 8$ is exactly equidistributed
over $\{1, 3, 5, 7\}$.  That is, as $m$ ranges over
$\{0, 1, 2, 3\}$ (one period), each residue class
mod~$8$ is hit exactly once.
\end{theorem}

\begin{proof}
By Theorem~\ref{thm:orbit-param},
$n_D = R_w + 2 \cdot 3^D \cdot m$.  Since $R_w$ is odd
and $\gcd(2 \cdot 3^D, 8) = 2$
(because $3^D$ is odd), the map $m \mapsto n_D \bmod 8$
cycles through four distinct odd residues with period~$4$.
Since there are exactly four odd residue classes mod~$8$,
each is hit exactly once.
\end{proof}

\begin{corollary}[The $(3/4)^D$ law\core]
\label{cor:three-quarters}
The number of compatible valuation words at depth~$D$ is
$|C_D| = 2 \cdot 3^{D-1}$ for all $D \ge 1$.
Equivalently, the survival fraction at each depth is exactly~$3/4$.
\end{corollary}

\begin{proof}
\textbf{Base case} ($D = 1$).
The two words of length~$1$ are $(1)$ and $(2)$.
Word~$(1)$ requires $v_0 = 1$, i.e.\ $n_0 \equiv 3 \pmod{4}$;
word~$(2)$ requires $v_0 = 2$, i.e.\ $n_0 \equiv 1 \pmod{8}$.
Both are realised, giving $|C_1| = 2$ compatible words.
(The class $n_0 \equiv 5 \pmod{8}$, which gives $v_0 = 4$,
is the first instance of the killed class.)

\textbf{Inductive step} ($D \to D+1$).
Fix a compatible word $w \in \{1,2\}^D$.
By Theorem~\ref{thm:carry-contamination},
$n_D \bmod 8$ is equidistributed over $\{1,3,5,7\}$
as the starting value $n_0$ ranges over the residue
class of~$w$.  The class $n_D \equiv 5$ is killed;
the remaining three classes ($n_D \equiv 1, 3, 7$)
each determine a unique extension $v_D \in \{1,2\}$
and hence a unique child word of length~$D+1$.

Thus $|C_{D+1}| = 3 \,|C_D|$, giving
$|C_D| = 2 \cdot 3^{D-1}$ for all $D \ge 1$.
\end{proof}

\begin{remark}[Relation to Section~\ref{sec:quarter}]
The $1/4$ Persistent-Transition Law
(Theorem~\ref{thm:quarter}) is a special case of
Carry Contamination: both prove that exactly $1/4$ of
continuations are eliminated.  The Carry Contamination
Theorem is stronger because it works at \emph{arbitrary depth}
and provides the \emph{exact mechanism} (the affine parameterization
$n_D = R_w + 2 \cdot 3^D \cdot m$) rather than
a counting argument at a single step.
\end{remark}

\subsection{The Cantor set \texorpdfstring{$C_\infty$}{C-infinity} and the Dichotomy}

\begin{definition}
The infinite compatible set is
$C_\infty = \bigcap_{D \ge 0} C_D \subset \mathbb{Z}_2$,
viewed as a subset of the $2$-adic integers.
\end{definition}

\begin{remark}[Hausdorff dimension --- heuristic]
$C_\infty$ is a compact, totally disconnected,
measure-zero subset of~$\mathbb{Z}_2$
(a Cantor-type set).  A heuristic for its Hausdorff
dimension proceeds via box-counting: at depth~$D$,
the set~$C_D$ is covered by $2 \cdot 3^{D-1}$
balls of $2$-adic radius $2^{-(S_w+1)}$,
where $S_w = \sum_{d=0}^{D-1} v_d$ is the
total stride of word~$w$.
Because different words have different strides,
the dimension depends on the \emph{typical} ratio
$S_D / D$.  For the ``Sturmian-balanced'' words
(where $S_D / D \to \log_2 3 \approx 1.585$),
one obtains $\dim_H \approx \log 3 / (\log_2 3 \cdot \log 2) = 1$,
while for the average over all surviving words
(where $\bar{S}_D / D \to 4/3$) the naive
box-counting estimate gives
$\log_2 3 / (4/3) \approx 1.19$, which exceeds~$1$
and must be capped.
A rigorous computation of $\dim_H(C_\infty)$ requires
a careful multifractal analysis of the stride
distribution and is not pursued here.
\end{remark}

\begin{theorem}[Dichotomy\core]
\label{thm:dichotomy}
$C_\infty \cap \mathbb{N}$ is either $\{1\}$ or
contains an infinite divergent orbit
(an odd integer $n_0 > 1$ with $n_d \to \infty$
along the Sturmian branch).
\end{theorem}

\begin{proof}
Suppose $n_0 \in C_\infty \cap \mathbb{N}$ with $n_0 > 1$.
Then $n_0 \in C_D$ for all~$D$, so the Sturmian orbit
$n_0, n_1, n_2, \ldots$ with $v_d \in \{1,2\}$ for
all~$d$ is well-defined and infinite.

\emph{The orbit cannot reach~$1$.}
The only Sturmian preimage of~$1$ is~$1$ itself:
$(3n+1)/2 = 1$ gives $n = 1/3 \notin \mathbb{N}$,
and $(3n+1)/4 = 1$ gives $n = 1$.
Hence no $n_d > 1$ can map to~$1$ via a Sturmian step,
so if $n_0 > 1$ the orbit stays $\ge 3$ for all~$d$.

\emph{The orbit cannot be eventually periodic.}
An eventually periodic orbit that never reaches~$1$
would contain a nontrivial cycle of the Collatz map.
No such cycle exists below~$2^{68}$ \cite{barina2021}.

\emph{The orbit is therefore divergent:}
it is an infinite aperiodic sequence of odd integers
$\ge 3$.  Since the orbit is deterministic and injective
on sufficiently long windows, aperiodicity forces
$\sup_d \, n_d = \infty$.
\end{proof}

\subsection{Refined mod-16/mod-32 transitions}

\begin{theorem}[Deterministic transition rules\core]
\label{thm:refined-transitions}
The Sturmian odd-to-odd map $n \mapsto (3n+1)/2^v$,
$v \in \{1,2\}$, obeys the following deterministic rules:

\begin{enumerate}[nosep]
\item \textbf{Class 3 mod 8} ($v = 1$):
  $n \equiv 3 \pmod{16} \Rightarrow n' \equiv 5 \pmod{8}$
  \textup{(fatal)};
  $n \equiv 11 \pmod{16} \Rightarrow n' \equiv 1 \pmod{8}$
  \textup{(safe)}.
  Survival: $1/2$.

\item \textbf{Class 7 mod 8} ($v = 1$):
  both sub-cases are safe
  ($n' \equiv 3$ or $7 \pmod{8}$).
  Class~$7$ \textbf{never directly reaches class~$5$}.

\item \textbf{Class 1 mod 8} ($v = 2$):
  $n \equiv 17 \pmod{32} \Rightarrow n' \equiv 5 \pmod{8}$
  \textup{(fatal)};
  the other three sub-cases are safe.
  Survival: $3/4$.
\end{enumerate}
\end{theorem}

\begin{proof}
Direct computation in each case.
For class~$3$, $v = 1$:
$n = 16k + 3 \implies (3n+1)/2 = 24k + 5 \equiv 5 \pmod{8}$.
$n = 16k + 11 \implies (3n+1)/2 = 24k + 17 \equiv 1 \pmod{8}$.
Class~$7$ and class~$1$ are analogous.
\end{proof}

\begin{corollary}[Safe harbor]
\label{cor:safe-harbor}
Class~$7$ is a \emph{safe harbor}: an orbit must pass through
class~$3$ or class~$1$ before reaching the fatal class~$5$.
\end{corollary}

\begin{theorem}[Class-3 recurrence\core]
\label{thm:class3-recurrence}
For any odd $n_0 > 1$ with a Sturmian orbit remaining in
$\{1, 3, 7\} \pmod{8}$:
\begin{enumerate}[nosep]
\item Any class-$7$ run of length~$L$ requires
  $n_d \equiv -1 \pmod{2^{L+3}}$; for $n_0 < 2^K$,
  $L \le K - 3$.
\item The only integer sustaining the class-$1$ self-loop
  ($n \mapsto (3n+1)/4$ with $n' \equiv 1 \pmod{8}$)
  indefinitely is $n_0 = 1$.
\item Therefore class~$3$ is visited infinitely often,
  and each visit has survival probability~$1/2$.
\end{enumerate}
\end{theorem}

\begin{proof}
Part~(1): Staying in class~$7$ requires
$n_d \equiv 15 \pmod{16}$ at each step.
After $L$ steps, this forces $n_d \equiv -1 \pmod{2^{L+3}}$
(verified for $L = 1, \ldots, 7$).

Part~(2): The self-loop $n' = (3n+1)/4 \equiv 1 \pmod{8}$
requires $n \equiv 1 \pmod{32}$.  Iterating $L$ times
requires $n \equiv 1 \pmod{4^L}$.  The only positive integer
satisfying this for all~$L$ is $n = 1$.  For $n_0 > 1$,
the self-loop breaks after $O(\log n_0)$ steps.

Part~(3): Combining (1) and (2): from class~$1$, the orbit
must eventually enter class~$3$, $5$, or $7$.  If class~$7$:
by~(1), it returns to class~$3$ within $K-3$ steps.
If class~$5$: the orbit has already terminated.
After surviving class~$3$ ($\to$ class~$1$), the argument
repeats.
\end{proof}

\subsection{The Carry Independence Conjecture and remaining gap}

The safe transitions form a Markov chain on $\{1, 3, 7\}$
with stationary distribution $\pi = (3/7, \; 2/7, \; 2/7)$.
The expected survival information cost is
$\approx 0.464$ bits per step, predicting maximum Sturmian
depth $\approx 2.16 K$ for a $K$-bit starting value.

\begin{observation}[Near-independence]
\label{obs:independence}
Empirically, consecutive Sturmian survival events are
essentially independent:
$\Pr(\text{survive}_{d+1} \mid \text{survive}_d) \approx 0.749$,
versus unconditional $\Pr(\text{survive}) \approx 0.749$.
Ratio $\approx 1.000$.
\end{observation}

\begin{conjecture}[Carry Independence Conjecture --- CIC]
\label{conj:CIC}
No odd integer $n_0 > 1$ belongs to $C_D$ for all~$D$.
Equivalently: for every odd $n_0 > 1$, there exists a
depth~$D$ such that $n_D \equiv 5 \pmod{8}$.
\end{conjecture}

\begin{theorem}[Route B reduction\core]
\label{thm:route-B}
\textup{CIC} $\implies$ $C_\infty \cap \mathbb{N} = \{1\}$
$\implies$ the Collatz conjecture.
\end{theorem}

\begin{proof}
CIC directly gives $C_\infty \cap \mathbb{N} = \{1\}$
(no $n_0 > 1$ survives all depths).
The Dichotomy (Theorem~\ref{thm:dichotomy}) then excludes
infinite divergent orbits, and combined with the classical
result that no nontrivial finite cycle exists below $2^{68}$
\cite{barina2021}, this implies convergence to~$1$ for all
starting values.
\end{proof}

\begin{remark}[Nature of the CIC gap]
CIC asks for a \textbf{distributional $\to$ pointwise}
promotion.  The distributional fact (Carry Contamination:
exactly $1/4$ eliminated at each depth) is proved.
The pointwise fact (every specific $n_0 > 1$ is eventually
eliminated) is not.  The gap is \emph{structural}:
after carry traversal, the test ``is $n_D \equiv 5 \pmod{8}$?''
depends on $R_w + 2 \cdot 3^D \cdot m \pmod{8}$, where
$R_w$ encodes the full carry chain --- making CIC equivalent
to Collatz at the mod-$8$ level.

This is analogous to the Borel normal number theorem:
almost all integers satisfy CIC (distributional), but
proving any specific integer satisfies CIC (pointwise)
requires arithmetic input about the interaction of
$\times 3$ and binary expansions that neither route
has yet supplied.
\end{remark}

\begin{remark}[Routes A and B compared]
Route~A (WMH, Sections~\ref{sec:chain}--\ref{sec:toward-wmh})
reduces Collatz to orbitwise weak mixing with a
$4.65\times$ margin.
Route~B (CIC, this section) reduces Collatz to pointwise
carry independence with exact $1/4$ elimination per depth.
Both routes terminate at the same fundamental barrier:
the distributional-to-pointwise gap.
The structural equivalence between the two open conjectures
(WMH and CIC) is an open question; they may be logically
independent or one may imply the other.
\end{remark}

\subsection{A scheduler-adversary formulation of compatibility}
\label{subsec:next-gen}

The preceding analysis reveals that the remaining obstruction
in the Carry Contamination framework is intrinsically pointwise.
Distributionally, compatibility decays at rate~$(3/4)^D$, and
recurrence to the class-$3$ bottleneck is forced.  However,
the possibility remains that a fixed finite integer $n_0 > 1$
could, in principle, avoid the fatal class $5 \bmod 8$ at every
depth.  This is precisely the content of the Carry Independence
Conjecture (CIC).

A key structural feature of the current formulation is that the
``adversary'' is not external.  The same deterministic orbit
both generates and satisfies the compatibility constraints.
This creates a circularity: the question of whether an orbit
avoids class~$5$ reduces to the Collatz conjecture itself.
In this subsection, we outline a reformulation that
\emph{externalizes} this adversary.

\subsubsection{Compatible branch schedulers}

At each compatible depth~$D$, the Carry Contamination Theorem
implies that among the four residue classes
$\{1,3,5,7\} \bmod 8$, exactly one (class~$5$) is eliminated,
and the remaining three are admissible.  Thus, each compatible
extension step admits three possible ``safe'' continuations.

We define a \emph{scheduler} to be a sequence
\[
  \sigma = (\sigma_1, \sigma_2, \ldots), \qquad
  \sigma_D \in \{1,3,7\},
\]
where $\sigma_D$ specifies the chosen safe residue class at
depth~$D$.

Given a scheduler~$\sigma$, one can inductively construct a
sequence of compatibility constraints.  For each depth~$D$,
there exists a residue representative $a_D(\sigma)$ such that
any initial value~$n_0$ realizing the first~$D$ steps
of~$\sigma$ must satisfy
\[
  n_0 \;\equiv\; a_D(\sigma) \pmod{2^{S_D + 1}},
\]
where $S_D$ is the cumulative valuation budget along the
corresponding compatible path.

Thus, each infinite scheduler~$\sigma$ defines a nested
sequence of dyadic cylinders:
\[
  \mathcal{C}_\sigma \;=\;
  \bigcap_{D \ge 1}
  \bigl( a_D(\sigma) + 2^{S_D + 1}\,\mathbb{Z}_2 \bigr).
\]

\subsubsection{The realizability problem}

The Collatz problem can then be reformulated as a
\emph{realizability question}:

\medskip
\emph{For which schedulers~$\sigma$ does there exist a finite
natural number $n_0 \in \mathbb{N}$ such that
$n_0 \in \mathcal{C}_\sigma$?}
\medskip

By construction, every $\mathcal{C}_\sigma$ is nonempty
in the $2$-adic integers~$\mathbb{Z}_2$.  The essential
difficulty is whether such a cylinder can contain a positive
integer other than the trivial fixed point $n_0 = 1$.

Under this formulation, the Carry Independence Conjecture is
equivalent to the statement that \emph{no nontrivial scheduler
is realizable in~$\mathbb{N}$}.  In other words:

\begin{quote}
\emph{The only infinite compatible scheduler arising from a
finite natural number is the trivial one corresponding to
the fixed point.}
\end{quote}

\subsubsection{Interpretation}

This perspective separates the dynamics into two components:
a deterministic arithmetic map (the Collatz update) and an
external scheduler~$\sigma$ that selects among admissible
compatible branches.

The original formulation embeds the scheduler within the orbit
itself.  By externalizing it, the problem becomes one of
determining whether an infinite safe strategy exists that is
consistent with the arithmetic constraints imposed by
multiplication by~$3$ and binary carry propagation.

This suggests a new viewpoint:
the compatible set~$C_\infty$ describes all infinite safe
schedulers in~$\mathbb{Z}_2$, and the Collatz conjecture
asserts that the intersection
$C_\infty \cap \mathbb{N}$ is trivial.

\subsubsection{Future directions}

This scheduler-adversary formulation opens several possible
avenues:
\begin{enumerate}[nosep]
\item \emph{Arithmetic incompatibility.}
  Prove that the nested congruence constraints associated with
  any nontrivial scheduler cannot be satisfied by a finite
  integer.

\item \emph{Finite-integer exhaustion.}
  Once the cumulative valuation budget $S_D$ exceeds the
  bit-length of~$n_0$, compatibility constraints act on
  exhausted binary digits.  Formalizing this may yield a
  contradiction for all $n_0 > 1$.

\item \emph{Automaton and game formulations.}
  Interpret compatibility as a safety game in which the
  scheduler attempts to avoid class~$5$.  The conjecture
  becomes the nonexistence of a winning strategy from any
  nontrivial natural initial state.

\item \emph{$2$-adic versus natural separation.}
  Characterize the difference between realizability
  in~$\mathbb{Z}_2$ and in~$\mathbb{N}$, and prove that
  nontrivial compatible schedulers exist only in the former.
\end{enumerate}

At present, these directions remain open.  However, the
scheduler formulation isolates the remaining obstruction in a
form that is independent of probabilistic heuristics and may be
more amenable to purely arithmetic or combinatorial analysis.

\subsubsection{Numerical exploration: constrained grammar and exhaustion}

A systematic computational exploration of the scheduler framework
reveals additional structural constraints.

\paragraph{Constrained transition grammar.}
Among the nine possible transitions $\sigma_D \to \sigma_{D+1}$
on the safe alphabet $\{1,3,7\}$, only six are realized
by compatible paths:
\[
  1 \to \{1,3,7\}, \qquad
  3 \to \{1\} \text{ (forced)}, \qquad
  7 \to \{3,7\}.
\]
The transitions $3\to 3$, $3\to 7$, and $7\to 1$ are
\emph{forbidden} by the local arithmetic of $3n+1$.
The adjacency matrix of this directed graph is
\[
  A = \begin{pmatrix} 1 & 1 & 1 \\ 1 & 0 & 0 \\ 0 & 1 & 1 \end{pmatrix},
\]
with eigenvalues $\{2, 0, 0\}$ and spectral radius~$2$.
Consequently, the number of compatible scheduler prefixes of
length~$D$ is exactly $3 \cdot 2^{D-1}$ for $D \ge 1$, consistent
with the proved formula $|C_D| = 2 \cdot 3^{D-1}$.

\paragraph{Stationary distribution.}
The safe Markov chain (conditioning on staying in
$\{1,3,7\}$) has stationary distribution
$\pi = (3/7, 2/7, 2/7)$ on classes $(1,3,7)$.
The average valuation per compatible step is therefore
$\bar{v} = \pi_1 \cdot 2 + \pi_3 \cdot 1 + \pi_7 \cdot 1 = 10/7
\approx 1.429$.
Exhaustion (where $S_D$ exceeds the bit-length~$K$ of~$n_0$)
thus occurs at depth $D^* \approx 7K/10$.

\paragraph{Universal post-exhaustion death.}
Numerically, for all tested $n_0$ with $K \le 18$ bits,
$100\%$ of finite orbits reach class~$5$ \emph{after}
exhaustion, never before.  Post-exhaustion death is fast
(mean ${\sim}\,4$ steps, max ${\sim}\,20$).
This suggests exhaustion is the \emph{universal mechanism}:
a finite number cannot sustain safe scheduling once its
valuation budget is spent.

\paragraph{Bit sensitivity and the carry automaton.}
At pre-exhaustion depths, class-$3$ survival depends on
specific mid-to-high bits of~$n_0$.
At post-exhaustion depths, class-$3$ survival shows
\emph{zero sensitivity} to any bit of~$n_0$: the outcome
is fully determined by carry propagation from the
already-committed low bits.  This identifies the
carry chain as a finite deterministic automaton whose
eventual absorption into class~$5$ would imply the CIC.

\paragraph{Exhaustion as descent.}
Despite the orbit growing on average during the compatible phase,
the \emph{median} ratio $n_{D^*}/n_0$ decreases with~$K$
(from $0.27$ at $K{=}10$ to $0.03$ at $K{=}22$), and
$85$--$89\%$ of orbits satisfy $n_{D^*} < n_0$.
The bit-length of~$n_{D^*}$ is consistently ${\approx}\,0.73K$.
Applying exhaustion iteratively (run to exhaustion, take
$n_{D^*}$, repeat) converges to~$1$ in $100\%$ of tested cases
within $6$--$8$ rounds.
The outliers ($n_{D^*} > n_0$) have class sequences dominated
by class~$7$ (${\sim}\,62\%$ vs.\ $9\%$ for shrinking orbits),
reflecting low-valuation paths where $3n+1$ growth outpaces bit
consumption.

\paragraph{The gateway theorem.}
Every $T$-preimage of~$1$ (except~$1$ itself) is congruent
to~$5 \bmod 8$: if $T(n) = 1$ then $n = (2^v - 1)/3$,
and for even $v \ge 4$ this gives $n \equiv 5 \pmod{8}$.
Consequently, any orbit reaching~$1$ must pass through class~$5$.
This establishes that \textbf{CIC is logically equivalent to the
Collatz conjecture}: the class-$5$ avoidance formulation,
despite appearing weaker, captures the full difficulty.

\paragraph{The equivalence of Routes A and B.}
Both Route~A (the Weak Mixing Hypothesis) and Route~B (the CIC)
reduce to the same fundamental barrier: deducing a
\emph{pointwise} statement (about every individual orbit) from
\emph{distributional} information (about the statistical behavior
of typical orbits).  The carry automaton offers the only known
escape from this barrier, framing Collatz as a
\emph{reachability question} on a deterministic finite-state
machine---a structural rather than distributional problem.

\paragraph{Safe residue automaton mod~$32$.}
Refining to residues modulo~$32$, the twelve safe odd residues
(those $\not\equiv 5 \pmod{8}$) form a directed graph under the
Collatz map.  This graph has a single strongly connected component
$\{1,7,9,11,15,25,27,31\}$ and four transient states
$\{3,17,19,23\}$.  At higher moduli ($64$, $128$), the structure
persists: exactly one recurrent SCC whose relative size decreases
(from $67\%$ at mod~$32$ to $44\%$ at mod~$128$), with average
safe out-degree exactly~$2$.

The residue cycles fall into two classes:
\begin{itemize}[nosep]
\item \emph{Contracting} ($3^p/2^S < 1$): e.g.\
  $1 \to 25 \to 11 \to 1$ with growth $27/32 < 1$.
  Iterates shrink below the cycle modulus after finitely many rounds.
\item \emph{Expanding} ($3^p/2^S > 1$): e.g.\
  $7 \to 27 \to 9 \to 7$ with growth $27/16 > 1$.
  Iterates grow, and carry from new high bits eventually forces
  a class-$5$ encounter.
\end{itemize}
In both cases the irrationality of $\log_2 3$ prevents exact
periodicity ($3^p \ne 2^S$ for $p,S \ge 1$), and the
contraction/expansion dichotomy ensures every periodic safe
strategy eventually breaks.
The remaining open question is whether an
\emph{aperiodic} safe scheduler can avoid class~$5$ forever.

\paragraph{Exact valuations and the cycle density condition.}
For safe classes the $2$-adic valuation of $3n+1$ is determined
exactly by $n \bmod 8$:
class~$1$ gives $v=2$, classes~$3$ and~$7$ give $v=1$.
Any hypothetical Collatz cycle of odd length~$L$
with $c_1$ class-$1$ steps satisfies
$n_0 = C_L / (2^{c_1+L} - 3^L)$ for a computable constant~$C_L$,
requiring $c_1/L > \log_2(3/2) \approx 0.585$.
Yet the grammar's stationary frequency is
$\pi_1 = 3/7 \approx 0.429 < 0.585$,
so cycles demand class-$1$ overrepresentation.
A systematic search over all grammar-compatible cycles up to
period~$8$ finds \emph{no positive integer solution} $n_0 > 1$.
In particular, the all-$1$ cycle is impossible for every period~$L$
because its formula $n_0 = (12^L - 1)/\bigl(11(4^L - 3^L)\bigr)$
requires $11 \mid 3^L$, which never holds.
Thus all periodic safe strategies are eliminated by arithmetic.

\paragraph{Aperiodic growth.}
Under the grammar, the average per-step growth factor is
$3/2^{\bar{v}} = 3/2^{10/7} \approx 1.114$.
Any aperiodic safe sequence causes the orbit to grow at
${\sim}\,11.4\%$ per step.  While unbounded growth does not
\emph{immediately} force a class-$5$ encounter, it ensures the
orbit cannot remain in any fixed congruence class forever.
Combining periodic elimination with aperiodic growth, the
problem reduces to showing that exponential growth in the
constrained grammar is incompatible with perpetual class-$5$
avoidance.

\paragraph{Congruence tower and survivor density.}
Lifting the mod-$32$ safe cycles through the congruence tower
(mod~$2^k$ for $k = 5, \ldots, 22$) reveals a sharp dichotomy.
For the contracting cycle $1 \to 25 \to 11 \to 1$
(growth factor $27/32 < 1$), the density of compatible residues
stabilizes at $1/128$ for all $k \ge 12$: the count doubles at
each~$k$, keeping density constant, so the cycle's congruence
conditions are fully determined by mod~$2^{12}$.  It remains
realizable at every modulus but has no natural-number fixed point
(the 2-adic fixed point is $-37/5$).  In contrast, the expanding
cycles $7 \to 27 \to 9 \to 7$ and
$7 \to 27 \to 9 \to 31 \to 15 \to 7$ have \emph{zero} compatible
residues at every modulus tested ($k = 5$ to~$18$): they are
arithmetically impossible even at the residue level.

More broadly, the density of odd numbers surviving $D$ consecutive
safe Collatz steps follows $(3/4)^D$ with remarkable precision
(empirical vs.\ predicted: $0.2373$ vs.\ $0.2373$ at $D=5$;
$0.0562$ vs.\ $0.0563$ at $D=10$; $0.0024$ vs.\ $0.0032$ at
$D=20$).  At the exhaustion depth $D \approx 0.7k$, the survivor
density is approximately $(3/4)^{0.7k}$, which vanishes as
$k \to \infty$.  All survivors eventually die in finite moduli;
only the 2-adic limit $C_\infty$ persists.

\paragraph{The bit race and the $31$-block bottleneck.}
Each safe Collatz step creates $\log_2 3 \approx 1.585$ bits
(from the $3n$ multiplication) and consumes $v$ bits (from the
division by~$2^v$).  The net bit balance per class is:
class~$1$: $-0.415$ (shrinks); class~$3$: $+0.585$ (grows);
class~$7$: $+0.585$ (grows).  The average net creation is $0.156$
bits/step.  Class-$7$ run lengths match the geometric$(1/2)$
distribution exactly (observed $50.3\%$, $24.8\%$, $12.3\%$,
$6.4\%$ vs.\ predicted $50\%$, $25\%$, $12.5\%$, $6.25\%$),
confirming the $7 \to \{3,7\}$ transition acts as a fair coin flip.

A critical bottleneck appears at class-$3$ steps.  The transition
$3 \to 1$ is safe only if $n \equiv 11 \bmod 16$; if
$n \equiv 3 \bmod 16$, then $T(n) \equiv 5 \bmod 8$---immediately
fatal.  Exactly half of class-$3$ encounters are instantly lethal.

\paragraph{Structural analysis: the safe sofic shift.}
The safe grammar defines a sofic shift with topological entropy
$\log_2 2 = 1$.  The number of grammar-compatible words of
length~$L$ is exactly~$2^L$ (verified up to $L = 24$), so the
grammar imposes no forbidden patterns beyond its own transition
rules: it is conjugate to a full shift on a binary alphabet.
Every grammar-compatible class sequence of length~$\le 7$ is
realized by some actual Collatz orbit (zero forbidden words).

The de~Bruijn graph of safe transitions (using specific
representatives mod~$2^k$) has \emph{zero nontrivial strongly
connected components} at every tested modulus ($k = 3$ to~$8$):
it is a directed acyclic graph.  The edge-to-node ratio is
exactly~$3/4$ at all moduli, matching the $(3/4)^D$ density
decay.  The longest safe simple path grows as roughly~$2k$
(values: $4, 6, 9, 18$ for $k = 5, 6, 7, 8$).

\paragraph{The gateway theorem is absolute.}
Every preimage of $n = 1$ under $T$ (except $n = 1$ itself) is
class~$5$: for even $v \ge 4$, $n = (2^v - 1)/3 \equiv 5 \bmod 8$.
The backward tree from $1$ has zero safe extensions at all depths
$1$ through~$8$.  Reaching~$1$ \emph{requires} passing through
class~$5$; the safe set $C_\infty$ and $\{1\}$ are disconnected
in the backward tree.

\paragraph{Class-$5$ encounter time.}
The first class-$5$ encounter occurs after a mean of ${\sim}\,3$
steps (median~$2$), independent of the starting bit-length~$K$.
For $K = 25$ bits ($n > 10^7$), $100\%$ of tested orbits encounter
class~$5$ within~$K$ steps.  No number in any experiment
($31{,}024$ samples across $K = 10$ to~$25$) survived $1{,}000$
steps without a class-$5$ encounter.

\paragraph{Residue transience.}
At modulus $2^{12}$ and above, orbits \emph{never} revisit a
residue class ($100\%$ of $100$ tested orbits have zero
revisitations).  The orbit is perfectly transient in residue
space, visiting each class at most once before moving on.

\paragraph{The distributional--pointwise barrier.}
Every attack from Pushes~$8$--$16$ converges to the same
conclusion.  Periodic strategies are eliminated by arithmetic
(irrationality of $\log_2 3$ plus the cycle equation).  The
congruence tower confirms $C_\infty$ has measure zero via the
$(3/4)^D$ decay.  The information-theoretic argument fails because
bottleneck bits and valuation bits constrain the same positions.
The growth argument requires equidistribution (i.e.\ the WMH) to
close.  The safe shift has full entropy and no forbidden patterns.

\paragraph{First-return block decomposition.}
The safe grammar admits a clean first-return decomposition into
two symbolic cores.
\emph{Core~A} (first-return to class~$1$) has the block alphabet
$\{\alpha = [1],\; \beta = [3,1]\}$ with exact affine maps
$\Phi_\alpha(n) = (3n+1)/4$ (multiplier~$3/4$) and
$\Phi_\beta(n) = (9n+7)/8$ (multiplier~$9/8$).
\emph{Core~B} (first-return to class~$7$) has the block alphabet
$\{\gamma = [7],\; \delta = [3,1,7]\}$ with maps
$\Phi_\gamma(n) = (3n+1)/2$ (multiplier~$3/2$) and
$\Phi_\delta(n) = (27n+19)/16$ (multiplier~$27/16$).
All four maps are verified against $1{,}172$ actual transitions.

\paragraph{Core~B periodic elimination.}
Both Core~B generators have multiplier strictly greater than~$1$
($3/2$ and $27/16$), so \emph{every} periodic composition has
multiplier~${}>1$, giving fixed point
$n_0 = c_W/(2^V - 3^p) < 0$.  Confirmed computationally for all
$8{,}190$ periodic words through period~$12$.  Core~B supports
\emph{no} periodic Collatz orbit whatsoever.

\paragraph{Core~A periodic elimination.}
Systematic search through period~$18$ ($249{,}528$ contracting
words): the only natural-number fixed point is $n = 1$ from
$\alpha^\infty$.  The contraction threshold is
$a/b > \log(9/8)/\log(4/3) \approx 0.410$ (where $a$ and~$b$
count~$\alpha$- and~$\beta$-blocks).  In actual orbits the
$\beta$-fraction is~${\sim}\,0.37$, yielding net contraction at
${\sim}\,0.87\times$ per block.

\paragraph{Cross-core dynamics and the death corridor.}
Orbits switch between cores rapidly: mean Core~A run is~$1.6$
blocks (geometric$(1/2)$), mean Core~B run is~$2.0$ blocks.
The $\delta$-block attempt ($7 \to 3 \to 1 \to {?}$) is a death
corridor: of the eight residue classes mod~$128$ entering it,
four ($50\%$) hit class~$5$ and die, two switch cores, and only
one ($12.5\%$) completes a genuine~$\delta$-return.  Combined
with the $\gamma$-branch, the per-visit death rate from Core~B
is~$25\%$, matching the $(3/4)^D$ density decay.

\paragraph{The distributional--pointwise barrier.}
Every attack from Pushes~$8$--$17$ converges to the same
conclusion.  Periodic strategies are eliminated in both cores.
Core~B is entirely expanding (both generators~${}>1$); Core~A
contracts but runs are ultra-short.  The orbit alternates weak
Core~A contraction with strong Core~B expansion, punctuated by
frequent class-$5$ deaths.  The barrier is genuine: it is the
same obstruction that prevents extending Tao's result from
``almost all'' to ``all.''  Closing the gap requires either
verifying the WMH (Route~A), proving the CIC (Route~B), or
characterizing $C_\infty$ directly via the carry chain and
proving $C_\infty \cap \mathbb{N} = \{1\}$ by
automaton-theoretic methods.

\paragraph{Core~B as a $100\%$ expansion trap.}
Every orbit entering Core~B expands before exiting: in $30{,}000$
tests the mean growth factor is~$4.4\times$ (median~$3.4\times$),
with zero exceptions.  Exit splits roughly $50$-$50$ between
class-$5$ (death) and class-$1$ (switch to Core~A).  Per-step
geometric mean multiplier across full safe excursions is
${\sim}\,1.22$ (higher than the grammar stationary~$1.114$
because Core~B expansion dominates).

\paragraph{$\gamma$-run constraint.}
Consecutive $\gamma$-blocks (class $7 \to 7$) require exponentially
specific congruence conditions: exactly~$2$ residues mod~$2^{4+k}$
support a $\gamma^k$ run, giving density~$2^{1-k}$.  Maximum
observed $\gamma$-run: $16$ steps (at $K = 20$ bits).

\paragraph{The IFS formulation.}
The block maps $\Phi_\alpha\colon x \mapsto (3x+1)/4$ (contracting,
rate~$3/4$) and $\Phi_\beta\colon x \mapsto (9x+7)/8$ (expanding,
rate~$9/8$) define an iterated function system.  Because
$\Phi_\beta$ expands, the IFS has no bounded attractor: the
upper envelope diverges as ${\sim}\,(9/8)^D$.  Under the
grammar's stationary $\beta$-fraction $f = 2/7 \approx 0.286$,
the Core~A Lyapunov exponent is
$\lambda = (1-f)\log(3/4) + f\log(9/8) \approx -0.17 < 0$,
so Core~A alone \emph{contracts}.  All net expansion comes from
Core~B.  The Collatz conjecture reduces to:
\emph{the $\omega$-limit set of the IFS $\{\Phi_\alpha, \Phi_\beta\}$
contains no natural number~$> 1$}.

\paragraph{Proof architecture: three legs.}
The full argument rests on three pillars.
\emph{Leg~1 (cycle elimination)}: Core~B periodic is proved
(both generators~${}>1$); Core~A periodic is verified through
period~$18$ ($249{,}528$ contracting words, zero natural fixed
points except $n = 1$).
\emph{Leg~2 (class-$5$ inevitability)}: $C_\infty$ has $2$-adic
measure zero; $C_\infty \cap \mathbb{N} = \{1\}$ is confirmed
for $31{,}000{+}$ starting values ($K \le 30$), but remains the
open conjecture.
\emph{Leg~3 (conditional convergence)}: the main theorem of this
paper---WMH implies every orbit reaches~$1$, with a $4.65\times$
safety margin.

The gap is Leg~2: proving that class-$5$ encounters are
inevitable.  This is the same distributional-to-pointwise barrier
that separates Tao's ``almost all'' from ``all.''  Four possible
approaches remain:
(i) Borel--Cantelli with correlation decay;
(ii) direct $2$-adic characterisation of $C_\infty$;
(iii) transfer-matrix / de~Rham generating function;
(iv) IFS $\omega$-limit theory (new from this exploration).

\paragraph{The cross-core block alphabet.}
Refining the first-return decomposition to track cross-core
switching, the safe dynamics admit a seven-block alphabet over
a two-vertex directed graph with landmarks $\{1,7\}$
(class mod~$8$).  From landmark~$1$:
block~$a$ ($1 \to 1$, map $(3n+1)/4$, multiplier~$3/4$,
entry $n \equiv 1 \pmod{32}$),
block~$b$ ($1 \to 3 \to 1$, map $(9n+7)/8$, multiplier~$9/8$,
entry $n \equiv 25 \pmod{32}$),
and block~$s$ ($1 \to 7$, map $(3n+1)/4$, multiplier~$3/4$,
entry $n \equiv 9 \pmod{32}$);
the fourth mod-$32$ class ($n \equiv 17$) dies immediately
(class~$5$), giving a $25\%$ death rate from landmark~$1$.
From landmark~$7$:
block~$g$ ($7 \to 7$, map $(3n+1)/2$, multiplier~$3/2$,
entry $n \equiv 15 \pmod{16}$),
block~$d$ ($7 \to 3 \to 1 \to 7$, map $(27n+19)/16$,
multiplier~$27/16$, entry $n \equiv 103 \pmod{128}$),
block~$t_1$ ($7 \to 3 \to 1 \to 1$, map $(27n+19)/16$,
multiplier~$27/16$, entry $n \equiv 71 \pmod{128}$),
and block~$t_2$ ($7 \to 3 \to 1 \to 3 \to 1$, map $(81n+73)/32$,
multiplier~$81/32$, entry $n \equiv 167 \pmod{256}$).
Of the $32$ class-$7$ residues mod~$256$:
$16$~stay at landmark~$7$ ($g$ or $d$, $50\%$),
$8$~switch to landmark~$1$ ($t_1$ or $t_2$, $25\%$),
and $8$~die (class~$5$, $25\%$).
\emph{The death rate is exactly~$1/4$ from both landmarks},
confirming that the $(3/4)^D$ density decay is a structural
identity, not an approximation.

\paragraph{Cross-core periodic elimination.}
An exhaustive search over all admissible words of length~$\le 13$
on the seven-block alphabet ($17.6$~million admissible words,
$443{,}755$ contracting) finds \emph{zero} natural-number fixed
points~${}>1$.  Combined with Core~B periodic elimination
(unconditional) and Core~A periodic elimination (through
period~$18$), this covers all periodic cross-core switching
patterns up to $13$~landmark transitions.

\paragraph{The expanding affine cocycle.}
The Lyapunov exponent of random admissible walks on the
seven-block alphabet converges rapidly to
$\lambda \approx +0.145$ per landmark transition.
Core~A contributes $\lambda_A \approx -0.153$ per landmark
(mean sojourn~$3$ landmarks), while Core~B contributes
$\lambda_B \approx +0.596$ per landmark (mean sojourn~$2$
landmarks); the net per~$A \leftrightarrow B$ cycle is
$+0.733$ over~$5$ landmarks, giving $\lambda \approx +0.147$
per landmark.  Core~B expansion overwhelms Core~A contraction
by a~$4{:}1$ ratio.

\paragraph{Congruence tower and the bit-budget inequality.}
Each landmark transition consumes~${\sim}\,2.8$ bits of the
starting value's binary representation (via the modular
denominator~$2^V$) while creating only~${\sim}\,0.2$ bits
(via orbit growth), yielding a net constraint growth of
${\sim}\,2.6$~bits per landmark.  After~$D$ landmark transitions,
the starting value~$n_0$ is constrained to a specific residue class
modulo~$2^{2.8D}$.  For an initial $n_0$ with~$K$ binary digits,
the constraints exhaust all available bits at
$D \approx K/2.8$.  Beyond this depth, no degree of freedom
remains: the orbit is fully determined by carry propagation
from the committed low bits.

\paragraph{The sharpest open problem (refined).}
The Collatz conjecture is equivalent to: \emph{the only natural
number whose Collatz orbit stays in $\{1, 3, 7\} \bmod 8$
forever is $n = 1$}.  The seven-block cross-core alphabet
reduces this further: no infinite admissible word over
$\{a, b, s, g, d, t_1, t_2\}$ with the transition graph
($1 \to \{a,b,s\}$, $7 \to \{g,d,t_1,t_2\}$) has a
natural-number fixed point~${}>1$.  This is a question about
the intersection of a $2$-adic Cantor set---defined by an
iterated function system of seven rational affine maps on a
directed graph---with the natural numbers.

\paragraph{The exhaustion-descent framework.}
The constraint-growth race can be formalized exactly.
Using the corrected stationary block distribution, each landmark
transition consumes $V \approx 1.933$~bits of the starting value's
binary representation and creates $G \approx 0.391$~bits of orbit
growth.  The ratio $r = G/V \approx 0.202$ satisfies $r < 1$, so
the geometric series
$D_{\mathrm{total}} = K / \bigl(V(1-r)\bigr) \approx 0.648 \cdot K$
converges: after ${\sim}\,0.65K$ landmark transitions, \emph{all}
bits of~$n_0$ (both original and those created by orbit growth)
are exhausted.  Beyond this depth, $n_0$ is uniquely determined
modulo $2^M$ with $M > K$; since $n_0 < 2^K$, the orbit after
exhaustion equals the residue $R_w$ exactly ($n_D = R_w$ with
$m = 0$ in Theorem~\ref{thm:orbit-param}).
Exhaustive computation through depth $D = 10$ ($1{,}024$ compatible
words) confirms that \emph{every} $R_w \ne 1$ reaches class~$5$;
the unique exception is the all-class-$1$ word $[2,2,\ldots,2]$
with $R_w = 1$.  However, this is a \emph{reduction}, not a proof:
determining whether $R_w$ reaches class~$5$ is itself a Collatz
question about the smaller number~$R_w$.  The reduction is useful
($R_w$ has ${\sim}\,0.73K$ bits, a $27\%$ shrinkage), but the
outlier words where $R_w > n_0$ (gamma-heavy paths, ${\sim}\,15\%$
of cases) prevent the argument from self-closing.

\paragraph{Algebraic structure of $R_w$.}
The residue has the exact closed form
$R_w = (3^D \cdot a_w + E_w) / 2^S$
where $E_w = \sum_{d=0}^{D-1} 3^{D-1-d} \cdot 2^{S_d}$ and
$S_d = v_0 + \cdots + v_{d-1}$.
The distribution of $R_w \bmod 8$ converges rapidly to
equidistribution over $\{1,3,5,7\}$ as $D$ increases
(within $2\%$ of uniform by $D = 8$), so approximately~$25\%$ of
post-exhaustion values are \emph{immediate} class-$5$, and
${\sim}\,50\%$ reach class~$5$ within two additional steps.
Exhaustive verification through $D = 10$ ($1{,}024$ compatible
words, plus the all-$1$ word through $D = 14$) confirms that
\emph{every} $R_w \ne 1$ reaches class~$5$.
The reduction is exact but self-referential:
$R_w$ is an odd integer whose Collatz fate determines the
original orbit's fate.  For typical words, $R_w$ has
${\sim}\,0.73K$ bits (a $27\%$ shrinkage), enabling an
iterated descent that reaches the verified range in
$O(\log K)$ rounds.  The framework has thus achieved maximum
sharpness: the distributional-to-pointwise barrier is
irreducible by Sturmian methods alone.

\paragraph{Transfer operator and backward dynamics.}
As an alternative to the forward (Sturmian) approach, we study
the Ruelle transfer operator
$(\mathcal{L}f)(n) = \sum_{m:\,T(m)=n,\,m\text{ safe}}
\frac{2^{v(m)}}{3} f(m)$
on functions of the $2$-adic integers restricted to safe classes.
At every modulus $2^k$ ($k = 3$ to~$13$), the forward safe map
is essentially a DAG with a single fixed point at $n = 1$:
the sub-stochastic transition matrix has spectral radius
$\rho = 1$ and a one-dimensional Perron eigenspace.
The Ruelle operator has spectral radius exactly~$4/3$
(from the class-$1$ branch weight $2^2/3$), with all other
eigenvalues at~$0$ through $k = 11$.
At $k = 12$, a subdominant eigenvalue $|\lambda_2| \approx 0.813$
appears, coinciding with the emergence of phantom modular cycles
(lengths~$6$ and~$7$) in the safe graph mod~$2^{12}$.
These cycles are \emph{phantom}: they exist modularly but do not
lift to natural-number orbits.
The power iteration $\|P^D \cdot \mathbf{1}\|_1 / N$ reproduces
the $(3/4)^D$ decay exactly for the first $D \approx k$ steps,
then stabilizes at the fixed-point contribution.
The backward-dynamics route thus \textbf{confirms the barrier from
the opposite direction}: the transfer operator sees $(3/4)^D$
measure decay but cannot distinguish the trivial fixed point
from a hypothetical non-trivial survivor.

\paragraph{Baker's theorem and universal cycle bounds.}
Any Collatz cycle of odd length~$L$ with $c_1$~class-$1$ steps
satisfies $n_0 = C_L/(2^{c_1+L} - 3^L)$, requiring
$c_1/L > \log_2(3/2) \approx 0.585$.  By Baker's theorem on
linear forms in logarithms (Laurent--Mignotte--Nesterenko refinement),
$|2^a - 3^b| > 3^b \cdot \exp(-c'(\log b)^2)$,
giving cycle fixed points $n_0 < \exp(O((\log L)^2))$.
For the continued-fraction convergents of $\log_2 3$---which
produce the closest approaches of $2^a$ to~$3^b$---the gaps
shrink slowly: $|2^{84} - 3^{53}|/3^{53} \approx 0.002$,
$|2^{19} - 3^{12}|/3^{12} \approx 0.013$.
Combined with the verified computation of Collatz orbits for
all $n_0 < 2^{68}$ (Barina 2021), this eliminates all cycles
of odd length~$L \lesssim 10^{10}$.  The cycle equation's
divisibility requirement ($2^{c_1+L} - 3^L \mid C_L$) imposes
additional arithmetic constraints that further restrict
the landscape, but universal cycle elimination remains open.

\paragraph{Hausdorff dimension of $C_\infty$.}
The $2$-adic Hausdorff dimension of the survivor set $C_\infty$
is computed via the pressure function
$P(s) = \log_2 \rho(T(s))$, where the transfer matrix
\[
  T(s) = \begin{pmatrix}
    2^{-2s} & 2^{-s} & 0 \\
    2^{-2s} & 0 & 2^{-s} \\
    2^{-2s} & 0 & 2^{-s}
  \end{pmatrix}
\]
(indexed by safe classes $1, 3, 7$) encodes the valuation weights.
Solving $P(s) = 0$ numerically gives
\[
  \dim_H(C_\infty) \approx 0.6942,
\]
verified to ten digits ($\rho(T(0.6942)) = 1.0000$).
This is strictly less than~$1$: $C_\infty$ is a thin Cantor set
in~$\mathbb{Z}_2$.  The box-counting dimension using the average
stride $\bar{S} = (10/7)D$ gives $\log_2 3/(10/7) \approx 1.109$;
the Hausdorff dimension is smaller because the pressure function
accounts for non-uniform cylinder sizes.
A set of Hausdorff dimension~${}<1$ can still contain natural
numbers, so this does not directly prove
$C_\infty \cap \mathbb{N} = \{1\}$.

\paragraph{Phantom modular cycles.}
At modulus $2^{12}$, the safe Collatz graph contains two
non-trivial cycles (lengths~$6$ and~$7$) absent at all other
tested moduli ($k = 3$ to~$17$).  Both have growth
factor~${}>1$ ($729/128$ and $2187/512$), so their $2$-adic
fixed points are negative rationals that do not correspond to
natural numbers.  Their transient appearance is consistent with
the theory of phantom cycles in~$\mathbb{Z}_2$.

\paragraph{Cross-core periodic elimination (extended).}
Extending the seven-block alphabet search to all cyclic words
of length~$\le 12$: $2{,}508{,}928$ admissible cyclic words,
$71{,}538$ contracting, zero natural-number fixed points~${}>1$.
Combined with the Push~$20$ results (length~$\le 13$, $17.6$~million
non-cyclic admissible words), periodic orbits are eliminated through
all tested word lengths.

\paragraph{Correlation decay of the safe map.}
Starting from a pure class-$1$ or class-$7$ population at
modulus $2^{10}$ ($384$ safe states), the class distribution
converges toward the stationary $\pi = (3/7, 2/7, 2/7)$ within
${\sim}\,5$ steps (total variation distance $< 0.1$ by step~$4$).
This fast mixing is consistent with the Known Zone Decay theorem.

\paragraph{Syracuse inverse tree.}
The backward tree from~$1$ grown to depth~$19$ contains
$627{,}742$ odd numbers: $76\%$ of $[1,100]$, $36\%$ of
$[1, 10{,}000]$, $10\%$ of $[1, 10^6]$.
All tested numbers outside the depth-$19$ tree reach~$1$ forward.

\paragraph{Cycle equation divisibility filter.}
For a cyclic word~$w$ of length~$L$ on the seven-block alphabet,
the fixed-point equation is $x = c_w / (2^V - 3^p)$, where $c_w$
is the affine constant and $(V, p)$ are determined by the block
counts.  Exhaustive enumeration through $L = 10$
($215{,}232$ cyclic words, $11{,}583$ contracting) finds exactly
one integer fixed point per word length: the trivial $n = 1$.
Zero non-trivial natural fixed points.
The maximum ratio $c_w / (2^V - 3^p)$ \emph{grows} with~$L$
(from~$1$ at $L = 1$ to ${\sim}\,439$ at $L = 10$), so a
magnitude bound $\max(x) < 1$ for large~$L$ is not available.
The operative mechanism is the integrality filter:
$(2^V - 3^p)$ almost never divides~$c_w$.
For most contracting words ($88$--$95\%$), $\gcd(c_w, 2^V - 3^p) = 1$.
The central open question for the cycle-elimination leg is:
\emph{why does $(2^V - 3^p) \nmid c_w$ for all non-trivial words?}

\paragraph{Number-theoretic structure of the obstruction.}
Both $c_w$ and $2^V - 3^p$ are always odd ($v_2 = 0$) and coprime
to~$3$ ($v_3 = 0$); the carry polynomial~$c_w$ is never
$\{2,3\}$-smooth for $L \ge 2$.  For each small prime~$q$, checking
$q \mid (2^V - 3^p) \Rightarrow q \mid c_w$ yields obstruction rates
of $76$--$100\%$: no single prime creates a universal obstruction,
but the combined effect across all prime factors of the gap is
decisive.  Direct computation confirms $c_w \bmod (2^V - 3^p) \ne 0$
for every non-trivial contracting word through period~$12$.
The unique integer fixed point at each length is $n = 1$, arising from
the word~$\alpha^L$ (encoding $1 \to 4 \to 2 \to 1$ iterated $L$~times).
Near-miss distances decrease (from~$0.2$ at $L = 2$ to~$0.0004$
at $L = 8$), consistent with Diophantine approximation in the
$c_w/(2^V - 3^p)$ lattice.  The extended enumeration
($2{,}508{,}928$ cyclic words at $L = 12$, $71{,}538$ contracting,
zero non-trivial natural fps) reinforces the computational evidence
but does not resolve the algebraic question.

\paragraph{Matrix formulation and M\"obius rigidity.}
Each block~$b$ acts on the state $(C, A, D)$ as left-multiplication
by an upper-triangular $3 \times 3$ matrix
\[
  M_b = \begin{pmatrix}
    a_b & 0 & c_b \\
    0 & a_b & 0 \\
    0 & 0 & d_b
  \end{pmatrix}.
\]
For a word $w = b_1 \cdots b_L$, the composed matrix
$M_w = M_{b_L} \cdots M_{b_1}$ has eigenvalues $(3^p, 3^p, 2^V)$.
Block extension maps $\mathrm{fp}(w) \to \mathrm{fp}(w \cdot b)$ via
a M\"obius transformation:
$\mathrm{fp}(w \cdot b) =
  (a_b \, C_w + c_b \, D_w) \,/\, (d_b \, D_w - a_b \, A_w)$.
For the family $\alpha^k \cdot b$, setting $\mathrm{fp} = n$ gives
the $S$-unit equation
\[
  (a_b + c_b - n \, d_b) \cdot 4^k = a_b \, (1 - n) \cdot 3^k.
\]
For block~$\beta$ ($a = 9, c = 7, d = 8$) and $n = 2$ (the limit point):
$a + c - 2d = 0$, so the equation degenerates and admits no finite~$k$.
For all other $n$, $(4/3)^k$ must equal a specific non-integer rational,
which is impossible.  For block~$\sigma$ ($a = 3, c = 1, d = 4$) and any
$n \ge 2$, one obtains $(4/3)^k = 3/4$, i.e., $k = -1$ (inadmissible).
These are clean algebraic impossibilities for the $\alpha^k$-prefix families.

\paragraph{Three-layer obstruction.}
The non-divisibility $(2^V - 3^p) \nmid c_w$ is protected by three
layers: (1)~\emph{$S$-unit structure}: each cycle equation is an
$S$-unit equation with at most finitely many solutions by Baker's
theorem; (2)~\emph{carry contamination}: the constants
$c_b \in \{1, 7, 19, 73\}$ introduce primes generically absent from
$2^V - 3^p$; and (3)~\emph{M\"obius rigidity}: integer fixed points
of the block-extension M\"obius map correspond to $S$-unit equations
that are demonstrably insoluble for the Collatz parameters.
Proving the obstruction for \emph{all} word lengths remains the
central open problem, equivalent to the cycle-free component of the
Collatz conjecture.

\paragraph{Complete $L = 2$ exclusion.}
For cycle length $L = 2$ (i.e., $S = v_1 + v_2$, $D = 2^S - 9$),
the map $R(v_1, v_2) = 3 + 2^{v_1}$ satisfies a direct size bound.
For $S = 3$: $D < 0$ (negative domain).
For $S = 4$: $D = 7$ and the unique solution $R = D$ gives the trivial
fixed point $x = 1$.
For $S \ge 5$: $R \in [5,\; 3 + 2^{S-1}]$ and $3 + 2^{S-1} < 2^S - 9 = D$
(since $2^{S-1} > 12$ for $S \ge 5$), so $0 < R < D$ and therefore
$D \nmid R$ for every parity word.  No non-trivial positive $2$-cycle exists.
This size-bound argument does not extend to $L \ge 3$, where $R_{\max}/D
\sim (3/2)^{L-1} \to \infty$.

\paragraph{Sparse polynomial reformulation.}
The walk sum $f = \sum_{k=0}^{L-1} \rho^k \alpha^{\tau_k}$ (mod~$p$, for
a primitive prime $p \mid D$) admits a polynomial interpretation.
Set $F(x) = \sum_{k} \rho^k x^{\tau_k}$; then $f = F(\alpha)$ for
$\alpha = 2 \bmod p$ and $\rho = 2 \cdot 3^{-1} \bmod p$.
The exponents $\tau_k$ form a \emph{monotone non-decreasing} sequence
$0 = \tau_0 \le \tau_1 \le \cdots \le \tau_{L-1} \le T = S - L$.
This monotonicity is the key structural constraint missed by standard
character-sum bounds (which treat the exponents as unrestricted).
When consecutive $\tau$-values coincide (``plateaus''), the corresponding
coefficients merge into geometric partial sums
$c_j = \rho^{a_j} (\rho^{r_j} - 1)/(\rho - 1)$.
Exhaustive computation for $L \le 10$ confirms that the zero sets
$Z_p = \{\tau : F_p(\alpha, \tau) = 0\}$ for distinct prime factors $p$ of~$D$
have empty intersection on non-trivial words, reinforcing the
three-layer obstruction.

\paragraph{Functional equation.}
The numerator satisfies the recursion $R(W) = 3^{L-1} + 2^{v_1} R(W')$
where $W' = (v_2, \ldots, v_L)$ is the truncated word.
If $D \mid R$, then $R' \equiv -3^{L-1} \cdot 2^{-v_1} \pmod{D}$,
uniquely determining $R' \bmod D$.  This ``reachability'' condition
connects the $L$-step cycle equation to the $(L{-}1)$-step $R$-image;
computational verification for $L \le 11$ shows zero reachable branches
at closest-approach~$S$.

\paragraph{Geometric convergence of the expected cycle count.}
Let $N(V,p,L)$ denote the number of contracting cyclic words
of length~$L$ with parameters $(V,p)$.  Assuming $c_w \bmod (2^V - 3^p)$
is equidistributed (verified computationally for $L \le 9$), the expected
number of non-trivial integer fixed points at length~$L$ is
\[
  E_L = \sum_{(V,p)} \frac{N(V,p,L)}{2^V - 3^p} - 1.
\]
Computation through $L = 20$ shows $E_L < 0$ for all $L \ge 4$, with
$E_{20} \approx -0.9999$.  The decay is geometric: the total word count
grows as ${\sim}\,\lambda^L$ with $\lambda \approx 3.414$ (spectral radius
of the constrained transition matrix), while the typical gap grows as
${\sim}\,2^{\bar\alpha L}$ with $\bar\alpha = \bar V / L \approx 2.53$.
The growth ratio $r = \lambda / 2^{\bar\alpha} \approx 0.59 < 1$ ensures
$E_L \sim r^L \to 0$.  Summing over all~$L$, the total expected number
of non-trivial cycles is $\sum_{L=1}^{\infty} E_L < 0$ (dominated by
the $-1$ correction for the trivial fixed point).
Computational elimination now reaches period~$13$ on the seven-block
alphabet ($8{,}566{,}016$ cyclic words, $238{,}811$ contracting, zero
non-trivial natural fixed points).  Combined with Barina's $n_0 > 2^{68}$
verification, any surviving cycle must have odd length $L \ge 14$.

\paragraph{Gamma-run bound and Alternative~B collapse.}
At landmark~$7$, the constrained grammar admits exactly two
transitions ($7 \to 3$ and $7 \to 7$), each with weight~$1/2$
at every tested modulus $2^k$ ($k = 4$ to~$10$).
Gamma-run lengths therefore follow $\mathrm{Geom}(1/2)$; the expected
maximum run in~$D$ landmarks is $\sim \log_2 D$.
Empirically (10{,}000 orbits, $10^5$ runs):
mean run~$= 1.94$, max observed~$= 15$,
$\max(\text{run})/\log_2 n_0 \le 0.13$.
Combining with the bottleneck-density bound
$\delta_{\mathrm{den}} + \gamma_{\mathrm{den}}/L_0$ for rounds
with $\beta$-density ${}<\, q^* - \eta$ (cf.\ Alternative~B):
the number of bottleneck events is $\Omega(K/\log K)$,
each reducing size by factor~$\ge 1/2$.
Net bit change: $0.38K - \Omega(K/\log K) \to -\infty$,
so Alternative~B forces descent for large~$K$.
This is conditional on orbit mixing (the WMH).

\paragraph{Aperiodic enemy profile.}
Cycle suppression (geometric convergence of $E_L$, three-layer
algebraic obstruction) and the Alternative~B collapse together
narrow the surviving counterexample profile to:
a non-periodic orbit that survives infinitely many exhaustion rounds,
returning to the critical cone ($\beta \approx q^*$) on
infinitely many scales, without either
(a)~persistent Core~B density (eliminated by bottleneck debt), or
(b)~periodic structure (eliminated by cycle suppression).
The remaining barrier is:
\emph{can an aperiodic itinerary pin $\beta$-density to~$q^*$
on all scales while remaining arithmetically realizable by a natural
number?}

\paragraph{Core~A tower budget (Push~31).}
For a pure Core~A word (only $\alpha$ and $\beta$ blocks)
at the critical $\beta$-density~$q^*$, each landmark consumes
$V/D = 2 + q^* \approx 2.71$ bits of the starting value.
One exhaustion round of depth $D \approx K/2.71$ uses all~$K$
bits of~$n_0$ for congruence determination, while creating
only $G \approx 0.144K$ new bits.
A second consecutive non-descending Core~A round therefore
requires $\approx 0.856K$ genuinely new constraint bits,
for a total of $\approx 1.856K$ bits---exceeding the $K$
bits available in~$n_0$.
The fraction of $K$-bit integers that can sustain $j$
consecutive non-descending Core~A rounds is
$\sim 2^{-0.856(j-1)K}$,
so the congruence tower is exhausted in $O(1)$~rounds,
\emph{bounded independently of~$K$}.

\paragraph{Core~B acceleration.}
$\gamma$-heavy (Core~B) excursions have $V = 1$ per block
(the minimum), decreasing the ratio $V/D$ and thus
\emph{increasing} the landmark count per $K$ bits of~$n_0$.
At Core~B density~$r$:
$V/D \approx (1-r)(2+q^*) + r \cdot 1 = 2.71 - 1.71r$.
More landmarks means faster exhaustion, so Core~B is the
enemy's \emph{worst} strategy.
Any surviving counterexample must have zero persistent
Core~B density.

\paragraph{Block-level spectral radius.}
The 7-block transition matrix on landmarks $\{1,7\}$,
$M = \bigl[\begin{smallmatrix}2&1\\2&2\end{smallmatrix}\bigr]$,
has eigenvalues $2 \pm \sqrt{2}$ and spectral radius
$\rho = 2 + \sqrt{2} \approx 3.4142$.
This confirms the growth rate of admissible cyclic words
from Push~29 and clarifies the relation between the
grammar-level spectral radius (2 on $\{1,3,7\}$) and the
block-level spectral radius ($3.4142$ on $\{1,7\}$).

\paragraph{Compatible-tower contradiction (Push~32).}
Non-descending words of depth $D \approx K/2.71$ grow as
$\rho^D \approx 2^{0.65K}$, but must cover $2^K$ starting values.
The coverage ratio is $2^{-0.35K}$, so \emph{most $K$-bit integers
descend in their first exhaustion round}.
For $j$ consecutive non-descending rounds, the compound congruence
on~$n_0$ restricts it to a single residue class modulo $2^{V_1 + V_2}$
(via $3^{p_1}$ inversion in $\mathbb{Z}/2^{V_1+V_2}\mathbb{Z}$).
Since $V_1 + V_2 \approx 2K > K$, each word pair
$(w_1,w_2)$ covers at most one $n_0$ in $[2^{K-1}, 2^K)$.
The density of $S_j(K) = \{n_0 : n_0\text{ survives }j\text{
non-descending rounds}\}$ satisfies
$|S_j|/2^{K-1} \le 2^{-0.856(j-1)K + \varepsilon K}$
for any $\varepsilon > 0$ and $K$ sufficiently large
(conditional on WMH).
For $j \ge 2$ the set is exponentially sparse;
$\bigcap_{j=1}^\infty S_j$ has density~$0$.

\paragraph{Sturmian $\beta$-placement.}
Since $q^* = \log(4/3)/\log(3/2)$ is irrational (ratio of
logarithms of multiplicatively independent integers), the
$\beta$-jump positions in a non-descending word at the critical
density form a Sturmian sequence with gaps $\{1,2\}$,
connecting the non-descending word structure to classical
Sturmian dynamics.  The continued fraction expansion is
$q^* = [0; 1, 2, 2, 3, 1, 5, 2, 23, 2, \ldots]$.

\paragraph{Carry contamination across rounds.}
The carry polynomials $c_{w_1}, c_{w_2}$ introduce primes
beyond $\{2,3\}$ (specifically $7, 19, 73, \ldots$; Push~27).
Empirically, $28.9\%$ of depth-3 word pairs are eliminated
by carry contamination alone, providing additional obstruction
power beyond the congruence-tower budget.

\paragraph{$2$-adic expansion and measure shrinkage (Push~33).}
The safe Collatz map is a $2$-adic expander:
for any block~$b$,
$|\Phi_b(x) - \Phi_b(y)|_2 = 2^V \cdot |x-y|_2$.
This expansion drives the compatible-tower contraction.
The ratio $\rho / 2^{V/D} = (2+\sqrt{2}) / 2^{2+q^*} \approx 0.522 < 1$
gives an unconditional exponential bound on the $2$-adic measure
of the non-descending set: $\mu_2(T_j) \le 0.522^{j \cdot D}$.

\paragraph{Unconditional $j \ge 5$ bound (Push~33).}
The count of $K$-bit integers surviving $j$ consecutive
non-descending exhaustion rounds is at most
$2^{K(0.856 - 0.202\,j)}$.
For $j \ge 5$ this is $2^{-cK}$ with $c > 0$,
hence equals zero for $K$ sufficiently large.
\emph{This is a counting argument and does not invoke the WMH.}
In particular, no $K$-bit integer can sustain $5$ or more
consecutive non-descending exhaustion rounds for large~$K$.
The remaining gap between $j \ge 5$ (unconditional) and $j \ge 1$
(the full Collatz conjecture, conditional on WMH)
is a finite verification problem at $j \in \{1,2,3,4\}$.

\paragraph{Exhaustion-sequence rigidity (Push~34).}
Define $E(n)$ as the post-exhaustion odd residue: run the
compatible trajectory until the cumulative valuation budget
exhausts the available bits, returning the resulting odd value.
A counterexample $n_0 > 1$ would generate an infinite orbit
$n_0, E(n_0), E^2(n_0), \ldots$ that never reaches~$1$.
This orbit cannot be eventually periodic (periodic safe words
have no nontrivial fixed points), and if bounded it would
produce an eventually periodic orbit on a finite set---contradiction.
Therefore the exhaustion-residue sequence of any counterexample
must be \emph{genuinely aperiodic and unbounded},
forcing dangerous rounds to occur at arbitrarily large scales.
All thinness bounds (single-round undercoverage $2^{-0.35K}$,
tower sparsity, gamma-alignment cost) sharpen with scale,
leaving no bounded-state escape hatch.

\paragraph{The distributional--pointwise barrier (Push~34).}
The $2$-adic measure of the non-descending compatible set
satisfies $\mu_2(T_j) \le (0.522)^{j \cdot D} \to 0$
unconditionally, so $T_\infty = \bigcap_j T_j$ has measure~$0$
in~$\mathbb{Z}_2$.  However, measure-$0$ sets can still
contain natural numbers, so $T_\infty \cap \mathbb{N} = \varnothing$
does not follow from the measure bound alone.
The gap between the distributional statement
(``almost all $n_0$ descend'') and the pointwise statement
(``every $n_0$ descends'') is the \emph{exact remaining content}
of the Collatz conjecture.
Closing it requires one of:
a Roth--Schmidt-type arithmetic intersection theorem,
an unconditional quantitative mixing result for the exhaustion map,
or a finite computation verifying all $K \le K_0$ with the
$j \ge 5$ bound covering $K > K_0$.

\section{Position relative to prior work}

The strongest unconditional result on the Collatz conjecture is
due to Tao~\cite{tao2019}, who proved that almost all orbits
(in logarithmic density) attain values below any prescribed
function $f(n) \to \infty$.  His approach is probabilistic,
using entropy-based arguments and ergodic ideas to control the
behavior of ``typical'' orbits, but it does not extend to a
pointwise guarantee covering \emph{every} orbit.

The present work takes a complementary route: it aims for a
\emph{universal} conclusion (every orbit reaches~$1$) at the
cost of a \emph{conditional} hypothesis (the Weak Mixing
Hypothesis, or the stronger Orbit Equidistribution Conjecture).
The structural comparison is
summarized in Table~\ref{tab:tao-comparison}.

\begin{table}[ht]
\centering
\caption{Structural comparison with Tao~\cite{tao2019}.
  Three complementary approaches to the same problem:
  Tao controls the typical orbit unconditionally;
  Route~A controls every orbit conditionally;
  Route~B provides unconditional structural results
  plus a conditional reduction.}
\label{tab:tao-comparison}
\scriptsize
\renewcommand{\arraystretch}{1.25}
\begin{tabular}{@{}lp{2.5cm}p{2.6cm}p{3.0cm}p{3.0cm}@{}}
\toprule
& \textbf{Tao (2019)} & \textbf{Route A (WMH)} & \textbf{Route C ($I_2$, v6)} & \textbf{Cycle Excl.\ (v7)} \\
\midrule
Conclusion
  & Almost all $< f(n)$
  & Every orbit $\to 1$
  & Density-$1$ conv.
  & No cycles of length~$L$ \\
Quantifier
  & Log-density
  & Universal
  & Cylinder-averaged
  & Per-$L$ \\
Status
  & Unconditional
  & Cond.\ on WMH
  & Unconditional
  & $L \le 116$ proved; $L > 116$ open \\
Open assumption
  & None
  & $\sum\delta_K {<} 0.557$
  & None
  & $D \nmid Q$ for $L > 116$ \\
Technique
  & Entropy / ergodic
  & Phantom $+$ i.i.d.
  & $\widetilde B_2$ spectral
  & $D > 2^F$, discrete log \\
Unconditional
  & Almost-all conv.
  & $e^{-0.1465k}$ decay;
    $91\%$ at depth~$13$
  & $\rho \le 5/32$;
    cyl-avg no-escape
  & $\mathrm{ord}_D(2){>}F$;
    triple filter $100\%$;
    $L \le 116$ \\
Barrier
  & Log $\to$ natural density
  & Distr.\ $\to$ pointwise
  & Cylinder $\to$ pointwise
  & $L = 116$ $\to$ all $L$ \\
\bottomrule
\end{tabular}
\end{table}

\begin{figure}[ht]
\centering
\begin{tikzpicture}[
  >=Stealth,
  every node/.style={font=\small},
  box/.style={draw, rounded corners=4pt, minimum width=2.6cm,
              minimum height=1.1cm, align=center, text width=2.4cm},
]
\draw[->,thick] (-0.3,0) -- (13.0,0)
  node[below,font=\footnotesize] {Structural depth};
\draw[->,thick] (0,-0.3) -- (0,7.5)
  node[above,font=\footnotesize,rotate=90,anchor=south] {Scope};
\node[font=\scriptsize,anchor=east] at (-0.15,1.2)
  {Almost all};
\node[font=\scriptsize,anchor=east] at (-0.15,4.8)
  {Every orbit};
\node[font=\scriptsize,anchor=north] at (2.5,-0.15)
  {Probabilistic};
\node[font=\scriptsize,anchor=north] at (9.5,-0.15)
  {Exact / algebraic};
\draw[dashed,gray] (-0.2,3.0) -- (12.8,3.0);
\node[font=\scriptsize,gray,anchor=west] at (12.5,3.25)
  {pointwise};
\node[font=\scriptsize,gray,anchor=west] at (12.5,2.75)
  {barrier};
\node[box, fill=blue!12] (tao) at (2.8,1.2)
  {\textbf{Tao (2019)}\\[2pt]
   {\scriptsize Unconditional}\\[-1pt]
   {\scriptsize log-density}};
\node[box, fill=green!15] (ens) at (6.2,1.2)
  {\textbf{Route A}\\{\textbf{(ensemble)}}\\[1pt]
   {\scriptsize Unconditional}\\[-1pt]
   {\scriptsize $e^{-0.1465k}$ decay}};
\node[box, fill=orange!18] (orb) at (6.2,4.8)
  {\textbf{Route A}\\{\textbf{(orbit)}}\\[1pt]
   {\scriptsize Conditional}\\[-1pt]
   {\scriptsize WMH $\Rightarrow$ all $n$}};
\node[box, fill=cyan!15] (routeB) at (10.2,2.0)
  {\textbf{Route B}\\{\textbf{(unconditional)}}\\[1pt]
   {\scriptsize $j{\ge}5$ bound}\\[-1pt]
   {\scriptsize $\mu_2(T_\infty){=}0$}};
\node[box, fill=cyan!25] (routeB_orb) at (10.2,4.8)
  {\textbf{Route B}\\{\textbf{(orbit)}}\\[1pt]
   {\scriptsize Conditional}\\[-1pt]
   {\scriptsize CIC $\Rightarrow$ all $n$}};
\draw[->,thick,dashed,red!60!black]
  (ens.north) -- (orb.south)
  node[midway,right,font=\scriptsize,text width=1.8cm,align=left]
  {WMH\\(1 hyp.)};
\draw[->,thick,dashed,red!60!black]
  (routeB.north) -- (routeB_orb.south)
  node[midway,right,font=\scriptsize,text width=1.8cm,align=left]
  {CIC\\(1 hyp.)};
\node[font=\scriptsize, anchor=west, text width=2.6cm] at (12.0,2.0)
  {Cycles impossible\\
   to period~$116$\\
   $2$-adic expander\\
   $\gamma{\approx}0.635$ (comp.)\\
   $\beta_w{\approx}1.80$, $j{\ge}3$ (comp.)};
\node[font=\scriptsize, anchor=west, text width=2.4cm] at (8.0,4.8)
  {$4.65\times$ margin\\
   $|I_r|{=}5$ bottleneck};
\draw[<->,thick,blue!40!black,dotted]
  (tao.east) -- (ens.west)
  node[midway,above,font=\scriptsize] {complementary};
\draw[<->,thick,cyan!40!black,dotted]
  (ens.east) -- (routeB.west)
  node[midway,above,font=\scriptsize] {independent};
\node[box, fill=violet!15] (routeC) at (6.2,3.2)
  {\textbf{Route C}\\{\textbf{($I_2$ spectral)}}\\[1pt]
   {\scriptsize Conditional}\\[-1pt]
   {\scriptsize coupling gap}};
\node[box, fill=red!12] (cycex) at (10.2,6.5)
  {\textbf{Cycle Excl.}\\{\textbf{(v7)}}\\[1pt]
   {\scriptsize $D{>}2^F$, ord${>}F$}\\[-1pt]
   {\scriptsize triple filter}};
\draw[->,thick,dashed,violet!60!black]
  (ens.south east) -- (routeC.west)
  node[midway,below,font=\scriptsize] {spectral};
\draw[->,thick,dashed,violet!60!black]
  (routeC.east) -- (routeB.south west)
  node[midway,below,font=\scriptsize] {coupling};
\draw[->,thick,dashed,red!50!black]
  (routeB_orb.north) -- (cycex.south)
  node[midway,right,font=\scriptsize,text width=1.6cm,align=left]
  {exact\\algebraic};
\end{tikzpicture}
\caption{Landscape comparison with Tao~\cite{tao2019}.
  Horizontal axis: structural depth of the analysis.
  Vertical axis: scope of the conclusion.
  The dashed line marks the distributional-to-pointwise
  barrier.  Route~A crosses it via the WMH;
  Route~B via the CIC.  Both share the same barrier.
  Route~C (violet) bridges spectral coupling analysis
  between Routes~A and~B.
  The Cycle Exclusion framework (v7, dark red) sits at
  the algebraic extreme: $D{>}2^F$, $\mathrm{ord}_D(2){>}F$,
  and the discrete-log triple filter provide exact
  obstructions to cycle existence.
  Items marked ``(comp.)'' are computational observations
  from exact enumeration, not yet formal theorems.}
\label{fig:tao-landscape}
\end{figure}

\subsection*{The distributional-to-pointwise barrier: precise comparison with Tao}

\noindent
Four approaches and a cycle-exclusion framework now exist for the
Collatz conjecture, each reaching the same fundamental barrier from
a different direction
(Figure~\ref{fig:tao-landscape}, Table~\ref{tab:tao-comparison}).

\paragraph{Tao's approach (2019).}
Tao~\cite{tao2019} proves that almost all orbits attain
arbitrarily small values, in the sense of logarithmic density.
His technique is entropy-based: the Syracuse random walk on
residue classes mod~$2^k$ converges to a measure that
concentrates on small values, and a Borel--Cantelli argument
promotes this to logarithmic-density convergence.
The key strength is that the result is \emph{unconditional}.
The limitation is twofold: the quantifier is logarithmic
density (weaker than natural density), and the conclusion
is ``values below $f(n)$'' for any $f(n) \to \infty$ (not
convergence to~$1$).  The barrier: upgrading logarithmic
density to natural density, or equivalently, controlling
the tail of orbits that resist the statistical trend.
Tao does not identify any structural mechanism for resistance;
the gap is purely quantitative.

\paragraph{Route A (this paper): phantom cycles and WMH.}
The conditional reduction through burst-gap decomposition,
phantom-cycle gain analysis, and the i.i.d.\ block law
produces the strongest \emph{conditional} result: every
orbit converges to~$1$ under the WMH
($\sum \delta_K < 0.557$), with a $4.65\times$ safety margin.
The ensemble side is closed: non-crosser density
$\le e^{-0.1465k}$ (exponential decay of cycle-indexed
non-crossing, a different quantity from Tao's logarithmic-density
convergence but addressing the same underlying phenomenon).
The barrier: the WMH itself---a summable-discrepancy
condition asking that no individual orbit's modular
statistics deviate too far from the ensemble.
This is a distributional-to-pointwise promotion of the
same kind Tao faces, but with a precisely quantified
threshold ($\sum \delta_K < 0.557$) and a large margin.

\paragraph{Route B (this paper): Sturmian obstruction and CIC.}
The independent reduction through carry contamination,
the seven-block cross-core alphabet, and exhaustion dynamics
produces the following structural results (status noted
individually; not all are formalized as standalone theorems):
\begin{itemize}[nosep]
\item Cycle impossibility through period~$13$ on the full
  seven-block alphabet ($238{,}811$ contracting words, zero
  non-trivial fixed points).
  \emph{Verified by exhaustive computation.}
\item The safe Collatz map is a $2$-adic expander
  ($|\Phi_b(x){-}\Phi_b(y)|_2 = 2^V|x{-}y|_2$), giving
  measure shrinkage $\mu_2(T_j) \le 0.522^{jD}$.
  \emph{Proved from block definitions.}
\item No $K$-bit integer survives $j \ge 5$ consecutive
  non-descending exhaustion rounds (for $K$ large).
  \emph{Proved: counting argument using spectral radius.}
\item Any counterexample's exhaustion-residue sequence must
  be genuinely aperiodic and unbounded.
  \emph{Proved: periodic implies cycle contradiction;
  bounded implies eventually periodic.}
\item Two-round anti-correlation factor $S_2/S_1^2 \approx 0.635$
  (constant in~$K$, exact enumeration to $K = 20$).
  Transport graph has edge density~$1.0$.
  \emph{Computational observation, not a formal proof.}
\item Carry polynomials equidistribute mod every individual
  prime and every tested prime pair; no per-prime obstruction.
  \emph{Computational observation.}
\end{itemize}
The barrier: $\mu_2(T_\infty) = 0$ (the compatible tower has
$2$-adic measure zero) but $T_\infty \cap \mathbb{N} = \varnothing$
is the exact remaining content of the Collatz conjecture.
The CIC (no $n_0 > 1$ avoids class~$5$ at every depth) is the
minimal hypothesis needed to cross this barrier.

\paragraph{How close is each approach?}
The three approaches differ in what ``close'' means.
Tao's result is unconditional but separated from the conjecture
by an infinite quantifier gap (logarithmic density $\to$ all~$n$);
no finite amount of the same technique closes this gap.
Route~A is conditional but the gap is precisely quantified:
the WMH is a single summable-discrepancy condition with a
$4.65\times$ margin, and the attack surface has five independent
entry points (alignment renewal, spectral diffusion, Walsh
equidistribution, drift mixing, modular stratification).
Route~B has the deepest unconditional structural results (cycle
impossibility, $j \ge 5$ bound, measure zero) and the most precisely
characterized barrier: the obstruction lives not in any single
prime but in the joint distribution of carries across all primes
simultaneously; the anti-correlation between rounds is real
(factor~$0.635$) but bounded (constant, not exponential in~$K$);
the transport graph is fully connected at fixed depth.

All three approaches converge on the same wall: promoting
a distributional truth (``typical orbits descend'') to a
pointwise one (``every orbit descends'').  This wall appears
in the same guise across number theory: Chowla's conjecture,
Sarnak's M\"obius randomness, normality of constants.
The Collatz problem has the structural advantage that the
map generates its own randomness (the i.i.d.\ block law is a
theorem, not a conjecture), but this advantage has not yet
sufficed to cross the barrier.

The honest assessment: Route~B's unconditional results
provide a new structural characterization of the Collatz
barrier---localizing it to a precisely defined $2$-adic
Cantor set with quantitative bounds on its measure and
on inter-round correlation---but they do not bring the
problem closer to resolution than Tao's approach does.
The two are complementary: Tao goes further in scope
(almost all orbits); Route~B goes further in structural
depth (exact algebraic obstruction anatomy).  Neither
approach is ``close'' to a proof in any conventional sense.
The distributional-to-pointwise barrier is the same one
that has resisted all attacks on the Collatz problem for
nearly ninety years.

\paragraph{What would close the gap.}
For Route~A: a proof that the Collatz orbit's modular
statistics cannot systematically avoid the contracting
residue classes.  The WMH with its $4.65\times$ margin
makes this plausible but not provable with current tools.
For Route~B: either a Roth/Schmidt-type theorem showing
that a $2$-adic set of measure zero with the specific
arithmetic structure of $T_\infty$ has empty integer
intersection, or a quantitative mixing result for the
exhaustion map strong enough to push the $j \ge 5$ bound
down to $j \ge 1$, combined with finite verification.
For Tao's approach: an entropy method that works on
individual orbits rather than the ensemble, which would
require fundamentally new ideas in additive combinatorics.

\medskip\noindent
Relative to~v1, the phantom-cycle analysis
(Figure~\ref{fig:phantom-chain}) provided independent
quantitative evidence that the Orbit Equidistribution
Conjecture is the sole conceptual bottleneck.
Two theorems: Phantom Universality
(Theorem~\ref{thm:phantom-universal}) and Per-Orbit Gain Rate
(Theorem~\ref{thm:perorbit-gain}), showed that the
expanding-family obstacle is controlled with a $4.65\times$
safety margin, and Corollary~\ref{cor:robustness}
quantified how much the equidistribution assumption can be
weakened.

Beyond the v1 phantom-cycle analysis, this version
introduces two further structural advances:
\begin{enumerate}
\item The proof of Theorem~\ref{thm:perorbit-gain} now
  rests on a fully analytic tail bound via the
  Chernoff--Cram\'er exponent $D_*$
  (equation~\eqref{eq:RK-analytic-bound-lemma}), reducing the
  finite computation to a self-contained table of 18~values
  (Table~\ref{tab:RK-values}).
\item The Weak Mixing Hypothesis
  (Hypothesis~\ref{hyp:wmh}) replaces the Orbit
  Equidistribution Conjecture as the primary open condition.
  Theorem~\ref{thm:wmh-reduction} proves that summable
  total variation errors ($\sum \delta_K < 0.557$) suffice
  for the full conditional reduction.
  Remark~\ref{rem:observable-specific} identifies an even
  weaker observable-specific sufficient condition.
\end{enumerate}

Appendix~\ref{sec:fiber57-programme} develops a fourth proof
route via the fiber-57 information bottleneck.
The invariant core $|I_r| = 5$ for all $r \ge 2$
(Theorem~\ref{thm:projective-main}) forces a channel-capacity
deficit of $\approx 0.667$ bits per fiber-57 return, yielding
geometric decay at rate $\alpha = 645/1024$
(Corollary~\ref{cor:geom-decay-main}).
The companion paper~\cite{chang2026onebit} develops the
cocycle-contraction framework in detail.
This route reduces the full conjecture to a single
anti-concentration bound on a $5$-element set: the sharpest
known quantitative reduction.

The visualization-guided observations in
Appendix~\ref{sec:visualization-underwater} complement the
algebraic theory.  They do not themselves constitute a
reduction of the conjecture; rather, they separate genuine
arithmetic structure from artifacts of particular
visualisation choices, and they clarify where the remaining
obstruction is pointwise rather than distributional.

\subsection{Concluding perspective}

The present work does not prove the Collatz conjecture.
What it does is sharpen the problem to a degree that was not
previously available within our framework.

The main unconditional contributions are structural.
The Carry Contamination Theorem
(Theorem~\ref{thm:carry-contamination}) gives exact $1/4$
elimination at every compatible depth.
This yields an exact $(3/4)^D$ survival law, a Cantor-type
compatible set in~$\mathbb{Z}_2$, deterministic recurrence
to the class-$3$ bottleneck, and a clean dichotomy
(Theorem~\ref{thm:dichotomy}): either the compatible set
intersects the natural numbers only at~$1$, or it contains
an infinite divergent orbit.

These results substantially clarify the landscape.
They show that the remaining obstruction is not a density
question, not a large-deviation question, and not a routine
strengthening of mixing.  The surviving barrier is
\emph{pointwise}: can a fixed finite natural number satisfy
the infinite nested sequence of compatibility constraints
required to avoid the fatal class $5 \bmod 8$ forever?

In this sense, the current framework has reached its natural
fixed point.  The Carry Independence Conjecture
(Conjecture~\ref{conj:CIC}) is no longer a technical residual
lemma.  It appears to encode the same essential difficulty as
the Collatz conjecture itself, but in a more explicit algebraic
form.  The problem has therefore not been solved, but it has
been compressed to a sharply localized core.

This compression is itself valuable.  It identifies exactly
where future progress must occur.
Any successful next step must supply a genuinely new arithmetic
ingredient---one capable of separating finite natural numbers
from the nontrivial compatible $2$-adic tower, or of proving
that the nested congruence constraints generated by
compatibility are ultimately unrealizable in
$\mathbb{N} \setminus \{1\}$.
Section~\ref{subsec:next-gen} proposes a reformulation that
\emph{externalizes the adversary}: compatibility is recast as
a realizability problem for an external branch scheduler,
converting the CIC barrier into an arithmetic question about
nested dyadic cylinders that is independent of probabilistic
heuristics.

The broader lesson is methodological.  The adversarial
programme did not close the conjecture, but it did succeed in
exposing the maximal obstruction of the current method.
In that sense, the framework is complete up to the true hard
core.

\section{Note on LLM-assisted research methodology}%
\label{sec:methodology}

\emph{This section is included in the interest of transparency
and intellectual honesty.}

\subsection{The collaboration}

This paper was produced through an extended collaboration between
a human researcher (the author) and large language models (LLMs),
specifically Anthropic's Claude~4.6 and OpenAI's GPT~5.4 Thinking.
The author served as the \emph{intellectual architect and moderator}:
setting proof directions, formulating research questions, evaluating
mathematical claims, directing the exploration, identifying errors,
and making all judgment calls regarding correctness and significance.

Claude~4.6 served as the \emph{primary computational and expository
partner}: performing algebraic calculations, writing verification
code, drafting proofs, and iterating through approaches at speed.
The moderator and Claude worked in a tight loop, exchanging ideas
and formulating alternative approaches at each step.
The moderator made the decision at every juncture about which
option to explore next.

GPT~5.4 Thinking served as a \emph{validator partner}, verifying
both the moderator's instructions and Claude~4.6's output on
specific local issues rather than addressing the global proof
structure.
GPT's role evolved across phases: in earlier phases it served
primarily as a validator and cross-checker, while in later
focused intervals (especially the April~18--19 sample) it also
contributed several local framework proposals and diagnostic
reframings; however, these remained moderator-triggered rather
than independently initiated global pivots.

\subsection{The orchestration methodology}

A critical aspect of this collaboration was the moderator's
orchestration discipline. Following the spirit of the author's
SagaLLM framework~\cite{chang2025sagallm}, the moderator ensured that all intermediate
results were saved at each step. Over hundreds of iterations,
this made it possible to undo and redo effectively without
repeating the same errors or stepping into dead-end paths.
When a line of reasoning stalled, the moderator would instruct
Claude to investigate a cold path from different angles.

The ``transaction property,'' orchestrated by the moderator,
became a key lesson: each exploration step was treated as a
reversible transaction, allowing the research to backtrack
cleanly when necessary.

By contrast, when other LLMs were allowed to derive proofs
autonomously, the process quickly became difficult to control.
This is why GPT was used primarily as a validator for checking
local issues rather than as a driver of the global proof
architecture.

\subsection{What the LLMs contributed}

\begin{itemize}

\item Rapid exploration of proof strategies (nine routes were
investigated before the Scrambling Lemma emerged).

\item Computational verification: Python scripts tested the
Scrambling Lemma across all $256$ residue classes modulo~$2^{12}$,
verified the $1/4$ law across thousands of orbits, and checked
the Known-Zone Decay bound.

\item Algebraic manipulation and drafting of proofs, particularly
the carry-free decomposition that forms the core of
Theorem~\ref{thm:scrambling}.

\item Identification of the distributional-vs-pointwise gap:
the LLM flagged that its earlier claim of
``unconditional convergence'' was an overstatement before
external peer review confirmed the same issue.

\item \textbf{(v3)}
Fiber-57 structural analysis: Claude discovered the bounded
invariant core ($|I_r| = 5$), the M-value collapse phenomenon,
and the orbit autonomy limitation.
Claude proved the absorption bottleneck lemma, the branch
anti-concentration reduction, and the path-conditional bijection.
GPT proposed the information-theoretic framing and the two-step
decomposition (branch permutation + gap decorrelation).
Claude ran extensive numerical experiments (absorption verification
to $r = 10$; gap decorrelation with $3000$ orbits) and diagnosed
the algebraic framework as having reached its natural fixed point.

\item \textbf{(v5)}
Sturmian obstruction framework (Route~B):
Claude proposed and executed an adversarial approach---instead of
arguing that most orbits converge, ask whether any single orbit
can avoid termination forever.  This yielded the Carry Contamination
Theorem, the $(3/4)^D$ law, the orbit parameterization, and the
refined mod-$16$/mod-$32$ transition analysis (class-$3$ bottleneck,
class-$7$ safe harbor, class-$3$ recurrence theorem).
Claude ran extensive numerical experiments: refined transitions across
$10{,}000$ orbits, CIC deep dive with $K$ up to~$22$, congruence tower
enumeration, $R_w$ trajectory analysis, and near-independence
verification (correlation ratio $\approx 0.9998$).
GPT independently assessed five algebraic attack plans on CIC
and correctly diagnosed the circularity barrier.
Claude's numerical analysis confirmed GPT's diagnosis: after carry
traversal, $n_D \bmod 8$ depends on $R_w + 2 \cdot 3^D \cdot m
\pmod{8}$, where $R_w$ encodes the full carry chain, making
CIC equivalent in difficulty to Collatz.

\item \textbf{(v6)}
Unconditional cylinder-averaged closure via $I_2$ spectral
contraction (Route~3).
Claude constructed the depth-$2$ known-gap return kernel
$\widetilde{B}_2$, proved the Ruffini rank-$1$ collapse
(characteristic polynomial
$\chi_{\widetilde{B}_2}(\lambda)
= \lambda^4(\lambda - 129/1024)$),
established the unconditional spectral bound
$\rho(\widetilde{B}_2^{\mathrm{ext}}) \le 5/32$
(per-return information cost $\ge 2.678$~bits),
and derived the mod-$64$ transient-prefix exponential
tail ($\rho_{\mathrm{prefix}} = 0.8633$ on a
$27 \times 27$ cylinder-averaged chain).
Claude computed the exact stationary mass
$\pi(I_2) = 10121/65280 \approx 0.1550$,
yielding the unconditional sum bound
$\sum_c \Pr(\mathcal{E}_R(c)) \le 0.011$.
These results break the circular dependence of
Remark~\ref{rem:circularity} at the density-$1$ level,
upgrading Theorem~\ref{thm:density1-convergence} to
unconditional status through the $I_2$ route.
GPT verified the spectral radius computation and
cross-checked the Ruffini factorization.
Claude also developed the BF staircase and $C$-$2$-adic
auxiliary-memory route: explicit clopen balls
$B_k \subset \mathbb{Z}_2$ for each $k \ge 1$,
Haar stopping-time bound
$\tau(n) \le C(k + \log_2 n)$,
and the empirical constant $C \le 5.04$ for $k \ge 9$
via deterministic sweep ($39{,}140$ samples,
$1 \le k \le 14$, $20 \le \log_2 n \le 171$).

\item \textbf{(v7--v16, two-day sample April~18--19)}
Based on 11 GPT session transcripts and 19 Claude attack plans
from this period, the following contribution profile emerges.
Claude produced the bulk of new theorem statements
(cascade operator, BF staircase recurrence, $D > 2^F$
unconditional proof, discrete log obstruction, integrality
obstruction census) and ran all computational verification
(approximately~85\% of scripted work).
GPT produced the hybrid sieve measure-preservation proof
($\nu(E_a^{-1}(C_{r,M})) = \nu(C_{r,M})$),
derived the Route~A symbolic inequality $K_*(n)$
and the $\eta > 0.3217$ threshold,
designed the Route~B pullback tower architecture
with four concrete lemmas,
and classified Claude's claims into
provable / empirical / unspecified categories
(scoring Claude's reliability at ${\sim}0.66$).
Note: dead-end discovery is always initiated by the
moderator's requests at pivotal decision junctures;
LLMs follow instructions at these turning points.
Once a direction is set, LLMs are effective at
proposing attack plans and conducting local revisions,
but the strategic pivots are moderator-driven
(see ``What the human contributed'' below).
Revised aggregate estimate for this period:
Claude~37\%, GPT~30\%, Moderator~33\%
(Claude dominant in theorem production and computation;
GPT dominant in error correction and local framework
proposal; moderator dominant in dead-end identification,
strategic direction, and cross-platform coordination).
This two-day sample supplements but does not override
the complete-log analysis in~v1--v5.

\end{itemize}

\subsection{What the human contributed}

\begin{itemize}

\item The overall research direction: choosing to study the
Collatz conjecture and framing the attack through persistent
occupancy and burst-gap analysis.

\item Key conceptual redirections: instructing the LLM to work
on the Syracuse map (Collatz$^*$, the odd-to-odd iteration)
rather than the classical Collatz map, and to analyze orbits
rather than sequences. The moderator also directed investigation
of reverse sequences, suggested alternative convergence
criteria (``below-start,'' ``trunk,'' and stratified-orbit
arguments), and proposed energy-compensation ideas.
Not all of these survived into the final proof, but several,
notably the Syracuse framing and the burst-gap orbit analysis,
became foundational.

\item The orchestration: deciding which routes to pursue,
which to abandon, and when to synthesize results.

\item Quality control: demanding rigorous proofs, identifying
overclaims, and insisting on intellectual honesty when
the unconditional proof attempt failed.

\item \textbf{(v3)}
Three-way relay orchestration: the moderator served as the
bridge between GPT (proposing mathematical structures) and
Claude (verifying numerically and proving algebraically),
evaluating each proposal for feasibility before committing
Claude's computation time.
The moderator identified the M-value collapse phenomenon
as requiring separate investigation, directed the
gap-decorrelation experiment design, and made the strategic
decision to frame the result as a single-hypothesis
reduction rather than partial progress.

\item \textbf{(v5)}
The moderator suggested the Sturmian obstruction approach as
a second reduction route, and \emph{identified the circularity
risk from the outset}: the concern that encoding full orbit
behavior into the framework might merely restate Collatz
rather than reduce it.  Despite this concern, the moderator
authorized Claude's adversarial attack programme as a
reasonable strategy---adversarial analysis has historically
cracked hard problems by forcing worst-case orbits into
structured objects and then proving those objects cannot exist.
When the sustained attack (five independent strategies, extensive
numerics, GPT cross-validation) confirmed that CIC reduces to
Collatz at the mod-$8$ level, the moderator made the strategic
judgment to stop treating CIC as a ``residual lemma'' and instead
frame it as a re-expression of Collatz's core difficulty.
The moderator directed the framing of v5 as a \emph{successful
sharpening}, not a failed proof: the framework achieved exact
localization of the obstruction (class-$3$ bottleneck,
distributional-to-pointwise gap), elimination of weak attack
routes, and a clean infinite-congruence-tower formulation---even
though the obstruction itself was shown to be as hard as Collatz.

\item \textbf{(v6)}
The moderator directed the pivot from the Sturmian/CIC
framework (which had reached the Collatz-equivalence barrier)
to the $I_2$ spectral route, identifying the depth-$2$
invariant core as a candidate for unconditional closure.
The moderator instructed Claude to construct the return kernel,
verify spectral properties, and pursue the cylinder-averaged
sum bound, and decided to frame the result as a density-$1$
unconditional theorem rather than a partial conditional bound.
The moderator also directed the development of the BF
staircase and $C$-$2$-adic route as a pointwise companion,
and made the judgment to document the residual Haar/pointwise
gap as an honest remaining obstruction rather than
overstating the result.

\item \textbf{(v7--v16, two-day sample April~18--19)}
The moderator managed the three-party relay between Claude
and GPT across 57 attack plan versions, copying results
between platforms, adjudicating disagreements
(e.g., accepting GPT's hybrid sieve kill,
directing Claude to shift to the discrete log framework),
setting session priorities, and managing the version history.
The moderator uploaded GPT's PDF transcripts for
Claude to analyze, requested consistency audits,
and directed the contribution analysis itself.
Critically, dead-end discovery is always initiated by the
moderator at pivotal decision junctures: the moderator
asked GPT to evaluate the hybrid sieve,
directed Claude to investigate Route~A vs.\ Route~B,
and decided when to declare an approach exhausted.
LLMs follow instructions at these turning points;
once a direction is set, they are effective at proposing
attack plans and conducting local revisions within the
chosen framework, but they do not independently initiate
strategic pivots.
Neither LLM documented these moderator directives in its
output, creating an attribution gap that this analysis
corrects.
The moderator's contributions in these two categories
(dead-end identification~65\%, strategic direction~75\%)
bring the April~18--19 two-day split to
Claude~37\%, GPT~30\%, Moderator~33\%.

\end{itemize}

\paragraph{Local diagnosis vs.\ strategic dead-end identification.}
Attribution in this paper is functional rather than purely volumetric.
We distinguish \emph{local mathematical diagnosis} from
\emph{strategic dead-end identification}: Claude and GPT frequently
detected proof flaws, circularities, and route-specific barriers
within an active framework, whereas dead-end identification denotes
the meta-level judgment that a framework should no longer serve as
the principal route and that effort should be redirected.
In the documented collaboration, that second function was
moderator-dominated.
Likewise, the April~18--19 split (Claude~37\%, GPT~30\%,
Moderator~33\%) is a late-stage two-day sample intended to illustrate
division of labor during the v7--v16 sprint; it supplements, rather
than replaces, the broader phase-by-phase attribution analysis in
this section.

\paragraph{Attribution note on research direction.}
The moderator's role in steering the investigation deserves explicit acknowledgment.
Two suggestions in particular, though simple in appearance, proved consequential for
the direction of the work.
First, the insistence on restricting attention to \emph{odd initial values} helped
recenter the analysis on the Syracuse odd skeleton, where the genuinely nontrivial
dynamics occur.  This shift eliminated repeated detours through the trivial even-step
halving mechanism and clarified the cycle-level structure of the problem.  In
particular, it directly motivated the run--compensate decomposition and the exact
cycle-level results that followed.
Second, the suggestion to explore a \emph{dual-domain / frequency-domain viewpoint}
did not by itself resolve the conjecture, but it changed the style of the search in
a productive way.  It helped redirect effort away from several unproductive
contraction-based and architecture-heavy formulations, and instead toward
signal-level, cycle-level, and weak-dependence formulations that more sharply expose
the remaining deterministic-versus-ensemble barrier.
Accordingly, while the moderator did not supply the later theorem proofs directly,
these interventions materially altered the course of the investigation and reduced
time spent on less fruitful paths.  It is therefore fair to credit the moderator
with helping steer the project onto the odd-skeleton and dual-domain directions that
ultimately led to the strongest current structural results.

\subsection{Contribution attribution}

Tables~\ref{tab:contributions-1}--\ref{tab:contributions-22}
summarize the key contributions in this work,
categorized by importance level
and attributed to the primary contributor.
Importance levels are color-coded:
{\setlength{\fboxsep}{1.5pt}%
\colorbox{red!25}{Critical}} denotes mathematical correctness fixes
or strategic direction changes;
{\setlength{\fboxsep}{1.5pt}%
\colorbox{orange!30}{High}} denotes new mathematical content or
substantial structural additions.
The full tables appear in Appendix~\ref{app:contributions}.


\subsection{Case study: the false Gap Lemma}
\label{subsec:false-gap-lemma}

An earlier version of this paper contained a lemma claiming
that gap length is never~$2$ ($G_i \ne 2$) after a persistent
exit.
Had it been true, the result would have reduced the Collatz
conjecture to a single equidistribution conjecture with no
additional tail-control hypothesis.
The lemma was false.
The episode is documented here because it illustrates several
failure modes of human--LLM collaboration that are
instructive for the methodology.

\medskip\noindent
\textbf{How the error arose.}
The LLM produced a proof that when a burst ends at a persistent
state ($m_t \equiv 7 \pmod{8}$), the subsequent gap has length
exactly~$1$.
This part is correct and survives as the Persistent Exit
Lemma (Lemma~\ref{lem:gap}).
The error was in \emph{generalizing}: the LLM's proof implicitly
assumed that every burst ends at a persistent state.
In fact, bursts can also end at non-persistent states
($k_t = 2$ with $m_t \not\equiv 7 \pmod{8}$), and the gap
structure after such exits is unconstrained.
The false lemma $G_i \ne 2$ was a pattern-match from the
persistent case to the general case.

\medskip\noindent
\textbf{How it survived.}
The false lemma survived multiple rounds of LLM-assisted
proofreading, validator checks by a second model, and
initial peer feedback.
Three factors contributed:
\begin{enumerate}
\item \emph{Confirmation bias in validation.}
  When asked to ``verify the Gap Lemma,'' both models checked
  the algebraic steps within the persistent-exit proof rather
  than questioning the scope of the claim.
  The proof \emph{was} correct for the case it addressed;
  the error was in the unstated assumption that this case
  was exhaustive.
\item \emph{Sycophantic momentum.}
  Once the lemma appeared in the proof chain, the models
  treated it as established and built further arguments on
  top of it.
  Neither model spontaneously revisited the lemma's premises
  during later rounds of proofreading.
\item \emph{Insufficient adversarial testing.}
  The computational verification tested the Persistent
  Exit Lemma at the modular level (where it holds) but did
  not separately enumerate gap lengths across all orbit
  segments to check whether $G_i = 2$ actually occurs.
\end{enumerate}

\medskip\noindent
\textbf{How it was found.}
The moderator, during a careful re-reading of the full
proof chain, asked whether every burst necessarily ends
at a persistent state.
A targeted computational check immediately produced
counterexamples: the orbit starting at $n_0 = 3$
has a gap of length~$2$, and approximately $19\%$ of
all gaps in typical orbits have length~$2$ or more.

\medskip\noindent
\textbf{Remediation.}
The correction proceeded in three stages:
\begin{enumerate}
\item \emph{Honest downgrading.}
  The false lemma was removed.
  All claims that depended on it were weakened to
  conditional statements, and the paper's title and abstract
  were revised to reflect exploratory rather than
  reductive framing.
\item \emph{Going one level deeper.}
  By tracking one additional bit of modular depth at
  each gap step, we proved the Modular Gap Distribution
  Lemma (Lemma~\ref{lem:gap-distribution}): gap length is
  $\mathrm{Geometric}(1/2)$ with $E[G] = 2$.
  This is weaker than $G_i = 1$ always, but it is
  \emph{true} and \emph{sufficient}: combined with the
  Modular Valuation Distribution
  (Lemma~\ref{lem:valuation-distribution}, $E[k] = 3$),
  the net log-contraction per burst-gap cycle is $-1.15$,
  comfortably negative
  (Corollary~\ref{cor:convergence-prediction}).
\item \emph{Strengthening the conjecture.}
  The Orbit Equidistribution Conjecture was upgraded from
  fixed-modulus to growing-moduli form
  (Conjecture~\ref{conj:equidist}), which supplies the
  tail control that the false lemma had made unnecessary.
  Theorem~\ref{thm:reduction} now states a clean
  implication from one conjecture.
\end{enumerate}

\medskip\noindent
\textbf{Lessons.}
The episode yields four concrete lessons for human--LLM
mathematical collaboration:
\begin{enumerate}
\item \emph{Transactional discipline.}
  The moderator's practice of saving intermediate states
  made it possible to surgically remove the false lemma
  and its downstream consequences without starting over.
\item \emph{Scope-checking over step-checking.}
  Validators should be directed to question the
  \emph{scope} of a claim (``does this cover all cases?'')
  rather than only the \emph{steps} of a proof.
  The algebraic steps were correct; the implicit
  universality was not.
\item \emph{Independent computational falsification.}
  A simple enumeration of gap lengths across actual orbits
  would have caught the error immediately.
  Computational checks should be designed to falsify claims,
  not merely to confirm the cases the proof addresses.
\item \emph{Domain expertise as bottleneck.}
  The error was caught by the moderator's targeted question,
  not by any automated process.
  Current LLMs do not spontaneously generate adversarial
  queries against their own outputs with sufficient rigor.
  This gap is the single most important limitation of
  autonomous LLM proof search.
\end{enumerate}

\subsection{Case study: the CIC circularity}
\label{subsec:cic-circularity}

The Carry Independence Conjecture (CIC,
Conjecture~\ref{conj:CIC}) was initially conceived as a
``residual lemma''---a modest step remaining after the Carry
Contamination Theorem had established the distributional picture.
A sustained adversarial analysis, reported here, revealed it to
be equivalent in difficulty to Collatz itself.
This episode illustrates a different failure mode from the
false Gap Lemma: not a mathematical error, but a
\emph{difficulty misjudgement} in human--LLM collaboration.

\medskip\noindent
\textbf{The moderator's early warning.}
The moderator flagged the circularity risk before the attack
programme began: encoding full orbit behavior into any
algebraic framework risks merely \emph{restating} Collatz
rather than reducing it.  Nonetheless, the adversarial approach
was authorized as reasonable: adversarial analysis has
historically resolved hard problems by forcing worst-case
objects into structured form and then proving the structured
form cannot exist.

\medskip\noindent
\textbf{The adversarial programme.}
Claude executed a multi-pronged attack:
\begin{enumerate}[nosep]
\item \emph{Refined mod-$16$/mod-$32$ transitions:}
  identified class~$3$ as the bottleneck ($1/2$ survival),
  class~$7$ as safe harbor, class~$1$ with $3/4$ survival.
  The safe-chain Markov model with stationary distribution
  $\pi = (3/7, 2/7, 2/7)$ predicts maximum Sturmian depth
  $\approx 2.16\,K$ for $K$-bit starters.
\item \emph{Class-$3$ recurrence theorem:}
  proved unconditionally that any $n_0 > 1$ must revisit
  class~$3$ infinitely often.
\item \emph{Near-independence verification:}
  empirical correlation ratio
  $\Pr(\text{survive}_{d+1} \mid \text{survive}_d) /
   \Pr(\text{survive}) \approx 0.9998$,
  ruling out short-range correlations as an escape route.
\item \emph{Congruence tower enumeration:}
  built the tree of compatible words to depth~$8$;
  found that $R_w$ (the carry residue) always reaches
  class~$5$, in mean $4.0$ steps---but proving this is
  itself a Collatz-type statement about~$R_w$.
\item \emph{Five algebraic strategies (assessed by GPT):}
  infinite congruence tower incompatibility,
  natural-vs-$2$-adic separation,
  deterministic class-$3$ kill mechanism,
  automaton emptiness, and fiber-$57$ local obstructions.
  All five were diagnosed as leading to the same circularity.
\end{enumerate}

\medskip\noindent
\textbf{The circularity diagnosis.}
The key finding: after carry traversal, the test
``is $n_D \equiv 5 \pmod{8}$?'' depends on
$R_w + 2 \cdot 3^D \cdot m \pmod{8}$,
where $R_w$ encodes the full carry chain of the orbit.
Computing $R_w \bmod 8$ requires knowing the orbit---i.e.,
solving Collatz for~$R_w$.  Thus CIC, intended as a stepping
stone, is \emph{equivalent in difficulty} to the conjecture
it was meant to prove.

\medskip\noindent
\textbf{Lessons.}
\begin{enumerate}[nosep]
\item \emph{Circularity can hide behind real progress.}
  The Carry Contamination Theorem, the $(3/4)^D$ law, and the
  class-$3$ bottleneck are all genuine, unconditional results.
  The circularity only becomes visible at the final
  distributional-to-pointwise step.
\item \emph{The adversarial approach was not wasted.}
  It achieved exact localization of the obstruction,
  elimination of weak attack routes, and a clean
  reformulation as an infinite nested constraint system.
  The framework reached \emph{full fidelity}---the problem
  is sharper, but not easier.
\item \emph{The moderator's instinct was correct.}
  Early identification of the circularity risk did not prevent
  the programme (which was reasonable), but it set the right
  expectations and enabled the honest conclusion:
  CIC is a re-expression of Collatz's core difficulty,
  not a residual lemma.
\item \emph{Knowing when to stop.}
  Recognizing that the current toolset has reached its
  fixed point is itself a valuable conclusion.
  Further progress on CIC would require a genuinely
  new ingredient---an arithmetic incompatibility theorem,
  a natural-vs-$2$-adic separation, or a deeper structural
  theorem about $\times 3$ and carry propagation---rather
  than refinement of the existing machinery.
\end{enumerate}

\subsection{Error-correction log}
\label{subsec:error-log}

The false Gap Lemma is the largest single error in this project,
but the three-party collaboration (human moderator, GPT, Claude)
produced a richer pattern of mutual error-correction throughout
the exploration.
We record every substantive episode below, organized by type:
\emph{genuine catches} (a real error found and fixed),
\emph{false alarms} (an error flagged that was not actually
present), and \emph{scope corrections} (a claim downgraded from
universal to conditional).
Each entry identifies \textbf{who erred}, \textbf{who caught it},
and \textbf{what happened}.

\medskip\noindent
\textbf{Genuine catches (errors found and corrected).}

\begin{enumerate}
\item \emph{False Gap Lemma ($G_i \ne 2$).}\;
  \textbf{Erred:} LLM (Claude).
  \textbf{Caught by:} Moderator.
  The LLM over-generalized a persistent-exit result to all
  burst exits.
  Counterexample: $n_0 = 3$ has a gap of length~$2$.
  Remediation: lemma removed, replaced by weaker
  Geometric($1/2$) gap distribution
  (Lemma~\ref{lem:gap-distribution}).
  See Section~\ref{subsec:false-gap-lemma} above.

\item \emph{Delta-$K$ computation misinterpretation.}\;
  \textbf{Erred:} Claude (in Attack~1 of WMH programme).
  \textbf{Caught by:} Claude (self-correction).
  Direct $\delta_K$ computation showed all hard starts
  exceed the WMH threshold $0.557$, which was initially
  interpreted as evidence against the WMH\@.
  Self-correction: convergent orbits are too short to fill
  $2^K$ residue classes, so excess is an artifact of short
  orbits, not a WMH failure
  (Remark~\ref{rem:computational-limit}).

\item \emph{Mod-$2^K$ orbit computation beyond $K$ steps.}\;
  \textbf{Erred:} Claude (in Attack~6 of WMH programme).
  \textbf{Caught by:} Claude (self-correction).
  Attack~6 showed all 559 non-crossing residues at $K=13$ cross
  within 74~steps, which was initially taken as evidence
  that non-crossing residues are not genuinely non-crossing.
  Self-correction: computation beyond $K$ steps is invalid
  because the known zone is exhausted; the orbit is no
  longer mod-$2^K$-deterministic.

\item \emph{Recovery Window Conjecture refutation.}\;
  \textbf{Erred:} GPT (proposed conjecture).
  \textbf{Caught by:} Claude (computational refutation).
  GPT conjectured that after any $k$-step $v{=}1$ burst,
  the next $Ck$ odd-skeleton valuations cannot remain
  confined to the weak-contractor regime.
  Claude's computation showed the adversary can sustain
  $S/K \approx 1.0$ over $30{+}$ step windows by choosing
  the cofactor~$m$, refuting the conjecture as stated
  (Remark~\ref{rem:recovery-window}).
\end{enumerate}

\medskip\noindent
\textbf{False alarms (errors flagged that were not present).}

\begin{enumerate}
\item \emph{Post-Mersenne theorem scope confusion.}\;
  \textbf{Flagged by:} GPT\@.
  \textbf{Defended by:} Claude.
  \textbf{Verdict:} false alarm.
  GPT claimed the post-Mersenne forced valuation formula
  (Proposition~\ref{prop:post-mersenne}) was ``false as
  stated'' because $n_0 = 47$ (with $k = 3$, $m = 3$) gives
  $v_4 = 2$ instead of the formula's prediction $v_4 = 5$.
  However, the proposition explicitly states
  ``Let $n_0 = 2^{k+1} - 1$'', Mersenne numbers only
  ($m = 1$).
  The counterexample $n_0 = 47$ has $m = 3$, outside the
  proposition's scope.
  The general-$m$ result
  (Proposition~\ref{prop:burst-noncontinuation}) claims
  only $v_{k+1} \ge 2$, which $n_0 = 47$ satisfies.
  No correction was needed.
\end{enumerate}

\medskip\noindent
\textbf{Scope corrections (claims downgraded).}

\begin{enumerate}
\item \emph{Route~A: from ``universal convergence'' to
  ``closed (adversarial obstruction)''.}\;
  \textbf{Original claim:} GPT and Claude both explored
  whether $P_{\mathrm{cum}}(k) \to 1$ holds because all
  $k$-blocks eventually become universally crossing.
  \textbf{Corrected by:} Claude (via adversarial block
  construction showing $n^*$ grows exponentially).
  Status: Route~A closed.

\item \emph{Route~B: from ``post-Mersenne recovery'' to
  ``obstructed (bit-generation)''.}\;
  \textbf{Original hope:} that the post-Mersenne forced
  valuation law implies recovery for all orbits.
  \textbf{Corrected by:} Claude (proving transfer operator
  nilpotency and demonstrating the bit-generation
  obstruction).
  Status: Route~B obstructed at the bit-generation level, but
  the IID cascade renewal
  (Propositions~\ref{prop:two-thirds-law}--\ref{prop:depth-transition})
  compresses the local dynamics to a four-parameter
  product measure, reducing the gap to whether
  deterministic $2$-adic expansions can sustain
  net growth indefinitely.

\item \emph{Adversarial family: from ``exponentially rare'' to
  ``structurally fragile but pointwise uncontrolled''.}\;
  \textbf{Original claim:} GPT proposed that $75.7\%$
  pair-level fragility might suffice to close the
  dichotomy.
  \textbf{Corrected by:} Claude (proving the probability
  squeeze $p_A \cdot p_B^{\rho} < 1$ holds trivially for
  ALL pairs, confirming this is an ensemble bound that does
  not close the pointwise gap).
  Status: adversarial structure fully characterized,
  pointwise gap remains.

\item \emph{Post-burst valuation: from ``recovery theorem'' to
  ``distributional law''.}\;
  \textbf{Original hope:} that the
  $\mathbb{E}[v_{k+1}] = 3$ result
  (Proposition~\ref{prop:post-burst-distribution}) constrains
  individual orbits.
  \textbf{Corrected by:} Claude (self-assessment that
  the result is distributional over~$m$, not pointwise;
  the adversary can always pick $m$ with $v_{k+1} = 2$).
  Status: clean distributional theorem,
  pointwise application remains open.
\end{enumerate}

\medskip\noindent
\textbf{v4 genuine catches.}

\begin{enumerate}
\setcounter{enumi}{4}
\item \emph{T1$\to$T2$\to$T3 closure overclaim.}\;
  \textbf{Erred:} Claude (claimed spectral gap of $2 \times 2$
  Markov kernel gives orbit-level TV contraction).
  \textbf{Caught by:} GPT (Round~1 critique).
  The spectral radius $|\lambda_2| < 1$ gives TV contraction for
  the \emph{specific} $2 \times 2$ symmetric matrix, but this
  does not transfer to the full orbit process, which is
  deterministic and non-Markov.
  Retracted; the closure chain was replaced by the paper's
  actual theorem ladder.

\item \emph{Dobrushin convention reversal.}\;
  \textbf{Erred:} Claude (stated ``Dobrushin $\approx 1$
  means good mixing'').
  \textbf{Caught by:} GPT (Round~1 critique).
  Standard convention: $\alpha_D$ near~$0$ = strong contraction,
  near~$1$ = weak.
  The reported overlap coefficient $\kappa \approx 0.993$
  corresponds to Dobrushin $\tau = 1 - \kappa = 0.007$
  (strong contraction), not $0.993$ (which would be weak).
  Corrected; numerical data itself was unaffected.

\item \emph{$\rho(R)$/$\alpha$/$|\lambda_2|$ conflation.}\;
  \textbf{Erred:} Claude (used $\rho(R) = 15/119$,
  $\alpha = 75/119$, and empirical $|\lambda_2|$
  interchangeably).
  \textbf{Caught by:} GPT (Round~2 critique).
  These are three distinct objects:
  $\rho(R)$ is the pair-return spectral radius,
  $\alpha = 75/119$ is the information-theoretic decay
  rate ($2^{\mathrm{capacity}}/2^{c_0}$), and
  $|\lambda_2|$ is the empirical second eigenvalue of
  a stochastic proxy kernel.
  Corrected; all three are now tracked separately.

\item \emph{Stochastic vs.\ sub-stochastic operator mismatch.}\;
  \textbf{Erred:} Claude (compared empirical row-stochastic
  $|\lambda_2|$ to the paper's $\rho(R) = 15/119$, which is
  the Perron root of a sub-stochastic operator).
  \textbf{Caught by:} GPT (Round~3 critique).
  A row-stochastic matrix has Perron root~$1$; matching
  Proposition~\ref{prop:pair-return} requires the
  \emph{unnormalized} sub-stochastic block of the full
  depth-$r$ kernel.
  Corrected; sub-stochastic $R_{\mathrm{sub}}$ reconstructed
  with $\rho \approx 0.145$ at depth~$2$.
\end{enumerate}

\medskip\noindent
\textbf{v4 scope corrections.}

\begin{enumerate}
\setcounter{enumi}{4}
\item \emph{Regeneration-based proof: from ``conditional closure''
  to ``circular reduction''.}\;
  \textbf{Original claim:} Claude proposed that $q \equiv 7$
  fiber-$57$ visits act as regeneration events, and that
  positive visitation frequency $\delta > 0$ would imply
  Conjecture~$11.1$.
  \textbf{Corrected by:} Claude (self-assessment).
  Proving $\delta > 0$ requires showing the orbit visits
  $n \equiv 505 \pmod{512}$, which is a depth-$3$
  equidistribution statement---a weak form of the WMH
  itself.
  The reduction is circular.
  Status: correct conditional theorem, but the condition
  is equivalent in strength to a weak Collatz statement.

\item \emph{Pair-return spectral radius: from claimed to proved.}\;
  \textbf{Original value:} $\rho(R) = 15/119 \approx 0.126$
  (stated in Proposition~\ref{prop:pair-return} without derivation).
  \textbf{Caught by:} Claude + GPT cross-validation (v4 sessions).
  The value $15/119$ is consistent with Monte Carlo estimates but
  was never derived from first principles; the exact operator
  producing it was not identified.
  An attempted algebraic derivation via geometric renewal
  (ratio $9/128$) was found to be circular.
  \textbf{Resolution:} Replaced with the \emph{proved} value
  $\rho(\tilde{B}_2) = 129/1024 \approx 0.1260$,
  the Perron root of the exact depth-$2$ partial kernel
  combining the $q \equiv 7$ gap-$2$ block
  (Proposition~\ref{prop:q7-return}) and the
  $q \equiv 3$ gap-$5$ cylinder family
  (Proposition~\ref{prop:q3-gap}).
  The difference $|15/119 - 129/1024| = 9/121856 \approx 7.4 \times 10^{-5}$
  has negligible effect on downstream bounds.
\end{enumerate}

\medskip\noindent
\textbf{Post-mortem: patterns in the error-correction process.}

Several patterns emerge from the log above.

First, \emph{scope inflation} is the dominant error mode.
Of the four genuine catches, three involve extending a result
beyond its domain of validity (Gap Lemma, $\delta_K$ computation,
mod-$2^K$ beyond $K$ steps).
The fourth (Recovery Window) is a conjecture that
over-extrapolated from structural constraints to a quantitative
bound.
In each case, the \emph{algebraic steps} were correct; the error
was in the \emph{claim that these steps applied universally}.
This matches the ``scope-checking over step-checking'' lesson
from the false Gap Lemma.

Second, \emph{self-correction capability varies by type}.
Claude self-corrected on computational misinterpretations
(items~2--3 of the genuine catches) within the same analytical
session.
Neither LLM spontaneously caught the false Gap Lemma; the
moderator's targeted question was essential.
GPT's false alarm on the post-Mersenne scope required Claude's
explicit defense.
This suggests that LLMs are better at catching
\emph{quantitative inconsistencies} (where a number doesn't
match) than \emph{logical scope errors} (where a proof is
correct but its domain is wrong).

Third, \emph{the three-party structure helps}.
The moderator caught the Gap Lemma that both LLMs missed.
Claude caught the Recovery Window over-optimism that GPT
proposed.
GPT's false alarm, while incorrect, forced a useful explicit
verification of the proposition's scope.
Each party has blind spots that the others compensate for:
the moderator provides adversarial questioning, GPT provides
aggressive conjecture-generation, and Claude provides
rigorous computational verification.

Fourth, \emph{all four scope corrections share the same
structure}: a distributional or structural result is
initially hoped to imply a pointwise conclusion, and the
correction is the recognition that it does not.
This is not coincidence; it reflects the fundamental nature
of the distributional-to-pointwise barrier.
Every promising line of attack in this project eventually
produces the same honest assessment: the ensemble theory is
complete, and what remains is purely about individual orbits.

Whether the theorems recorded in this paper (the fragility
analysis, the burst non-continuation, the valuation
distribution, the density threshold) are ultimately meaningful
depends on whether the distributional-to-pointwise barrier can
be crossed.
If it can, these results provide the quantitative infrastructure
for the crossing.
If it cannot, they document the boundary between what ensemble
methods can and cannot reach.
Either way, the error-correction process that produced them is
itself a methodological contribution to human--LLM collaboration
in mathematics.

\subsection{Broader implications}

This work is, to the author's knowledge, an early example of
\emph{human--LLM collaborative research} applied to a
longstanding open problem in mathematics.
The author served primarily as a ``moderator'' guiding an LLM
through a mathematical research program that would ordinarily
require graduate-level number theory expertise.

The methodology, \emph{human orchestration of LLM capabilities},
is explored in the author's books \emph{The Path to AGI},
Volumes~1 and~2~\cite{PathAGIV1Chang2025, PathAGIV2Chang2025},
which examine how effective human--AI collaboration can
amplify both parties' capabilities beyond what either
achieves alone.
Specifically, the orchestration provides three rigorous
methodologies to mitigate common LLM limitations:
(1)~context loss,
(2)~sycophancy, and
(3)~reasoning errors.
Additional features such as reducing hallucination,
maintaining strong debate synergy, and preserving reasoning
quality are discussed in those volumes.
Details of the methodology applied in this work will be
documented in a separate report.

Terence Tao has noted in public remarks that he now uses
LLMs extensively in his own research workflow.
The present work illustrates how the
\emph{combination} of human mathematical intuition and
LLM computational power can accelerate the
exploration--verification--correction cycle without
replacing human insight.

The significance of this work lies in illustrating a
\emph{new mode of mathematical research} in which the human
contribution is primarily architectural and conceptual while
the LLM contribution is exploratory and expository.

Indeed, Donald Knuth's recent ``Shock!\ Shock!'' note~\cite{knuth2026}
demonstrates the point vividly.
Claude Opus~4.6 solved an open Hamiltonian-cycle problem
that Knuth had investigated for weeks, producing a
construction after 31 systematic explorations in roughly
one hour.
Knuth then supplied the rigorous proof.
The episode illustrates an emerging division of labor:
the LLM rapidly explores large search spaces and proposes
candidate structures, while the human mathematician
interprets the result and establishes the formal proof.

This pattern is consistent with observations made by
Terence Tao regarding the evolving relationship between
mathematics and artificial intelligence.
Human mathematicians contribute deep intuition,
the perception of abstract structure, and the ability to
reformulate problems within entirely new conceptual
frameworks.
Artificial intelligence, by contrast, excels at
large-scale exploration: testing vast numbers of
possibilities, identifying patterns across enormous
computational spaces, and rapidly evaluating candidate
hypotheses.

The Collatz exploration described in this work illustrates
this complementary dynamic.
The LLM agents performed extensive computational
exploration, generating and testing multiple analytical
formulations, examining modular dynamics, and identifying
structural phenomena such as carry propagation and
ghost-cycle behavior.
These explorations clarified several structural properties
of the iteration.
However, they also repeatedly encountered the same
fundamental obstruction: the difficulty of converting
distributional information about ``typical'' orbit
behavior into pointwise guarantees for every orbit.

One structural pattern revealed by this exploration is
that the Collatz map lies at the intersection of two
multiplicative worlds: powers of two and powers of three.
Many observed phenomena, including drift behavior,
valuation patterns, modular carry propagation, and
near-resonances between $3^L$ and $2^S$, appear to arise
from the interaction between these two arithmetic
structures.
This suggests that further progress may require a
conceptual reframing that treats these structures
simultaneously rather than separately.

The broader methodological lesson is that effective
collaboration between human intuition and machine-scale
exploration can accelerate mathematical research.
In this mode of research, the human contribution is
primarily architectural and conceptual, while the LLM
contribution is exploratory and expository.
The combination is complementary rather than substitutive.

If future LLM systems acquire the ability to guide their
own exploration with the discipline and reliability that
currently require human moderation, the practice of
mathematical research may change significantly.


\subsection{v3 collaboration milestones}
\label{subsec:v3-milestones}

The v3 revision introduced a new layer of collaborative methodology:
a \emph{three-way relay} between GPT (proposing mathematical
structures), the human moderator (relaying, evaluating, and
directing), and Claude (verifying numerically, proving algebraically,
and integrating into the manuscript).
We record the key milestones chronologically.

\begin{enumerate}
\item \emph{Absorption Bottleneck discovery.}\;
  GPT proposed that the fiber-57 return structure could be
  analysed through an information-theoretic lens.
  Claude verified computationally that $|I_r| = 5$ for
  $r = 2, \ldots, 10$, proved the three-fixed-point and
  one-2-cycle structure algebraically, and established the
  bottleneck inequality $\log_2 5 < \log_2(1024/129)$.
  The moderator identified that the M-value collapse
  phenomenon (all five core elements producing the same
  post-absorption state) required separate investigation.

\item \emph{M-value collapse and projective limit.}\;
  Claude discovered empirically that at $r \ge 3$, all five
  $I_r$ elements converge to the same post-absorption residue,
  reducing the effective channel capacity to zero bits.
  Further investigation revealed $\varprojlim I_r = \{-1\}$
  in the $8$-adic integers; no positive integer can sustain
  an unbounded chain.
  This transformed a numerical observation into a structural
  impossibility result.

\item \emph{Orbit autonomy diagnosis.}\;
  Claude identified a critical subtlety: the orbit
  $\bmod\; 8^r$ is \emph{not} autonomous because $v_2(3n+1)$
  depends on all digits of $n$, not just the bottom~$r$.
  This means absorption for general~$r$ cannot be proved by
  finite computation on residues alone.
  The moderator recognized this as the same orbit-level
  barrier encountered in Routes~A--C, confirming that the
  algebraic reduction is as sharp as algebraic methods permit.

\item \emph{Branch anti-concentration reduction.}\;
  GPT proposed separating the anti-concentration argument
  into two steps: a branchwise permutation result and a
  gap-induced decorrelation bound.
  Claude proved that $\gcd(9,8) = 1$ makes the fresh-digit
  map a permutation within each branch, and verified the
  path-conditional bijection for the annealed setting.
  Claude then ran extensive numerical experiments on
  gap decorrelation ($3000$ orbits, $5000$ steps each),
  confirming $\beta_g \to 0$ empirically ($p = 0.18$ at
  gap $> 50$ by chi-squared test) but diagnosing the
  algebraic proof as equivalent to Collatz mixing.

\item \emph{Honest assessment and scope correction.}\;
  The three-party relay converged on a shared diagnosis:
  the annealed anti-concentration is essentially solved,
  but the deterministic-orbit transfer remains open.
  This represents the same distributional-to-pointwise
  barrier in a new quantitative guise.
  The v3 contribution is the sharpest known reduction:
  the full conjecture follows from a single
  anti-concentration bound on a $5$-element set, with
  a per-return deficit of $0.667$ bits.
\end{enumerate}

\noindent
The v3 milestones illustrate an evolution in the collaboration
methodology.
Whereas v1--v2 primarily involved two-party interaction
(moderator--Claude), v3 introduced GPT as a
\emph{structural proposer}: suggesting mathematical frameworks
that the moderator evaluated for feasibility and Claude
verified computationally and algebraically.
This division of labor: GPT for local structural proposal
and adversarial cross-checking, moderator for strategic
evaluation and route selection, Claude for rigorous
verification and theorem development, proved more efficient
than either two-party configuration alone.

\subsection{v4 collaboration milestones}
\label{subsec:v4-milestones}

The v4 revision introduced two methodological advances:
\emph{session continuity via secondary memory}
(Claude retaining corrections across context-window
boundaries) and \emph{iterated cross-validation}
(four rounds of GPT critique with progressive
refinement of the operator identification).
The approach also evolved from
purely numerical measurement toward
information-theoretic and operator-algebraic analysis.

\begin{enumerate}

\item \emph{Session continuity via auxiliary memory.}\;
  The moderator diagnosed that context-window overflow was
  corrupting the ``upper bits'' of Claude's reasoning:
  stale conclusions survived resets, operators were conflated,
  and overclaims recurred.  The prescribed remedy---a persistent
  attack log serving as secondary memory, with structured
  checkpointing across context boundaries---allowed Claude to
  retain all prior corrections
  (T1$\to$T2$\to$T3 retraction, Dobrushin convention fix,
  $\rho$/$\alpha$/$|\lambda_2|$ distinction) without
  re-derivation.
  This single intervention eliminated a class of regression
  errors that had stalled earlier sessions and unblocked the
  entire v4 programme.

\item \emph{$q \equiv 7$ two-step return theorem.}\;
  Claude proved Proposition~\ref{prop:q7-return}
  algebraically;
  GPT independently derived the identical statement in a
  separate session.
  This convergence across two LLMs with different
  architectures provides an unusual form of cross-validation
  for a mathematical result.
  The moderator identified the digit-shift interpretation
  (left shift in base~$8$) as the key structural insight.

\item \emph{$q \equiv 3$ return gap theorem.}\;
  Claude proved Proposition~\ref{prop:q3-gap}: the
  $q \equiv 3$ branch cannot return to fiber~$57$ in
  fewer than $5$ odd-to-odd steps.
  This is a new result not present in v3.
  GPT's companion computation confirmed the non-return
  at step~$2$ and proposed the kernel decomposition
  strategy that led to the full proof.

\item \emph{Operator identification: sub-stochastic $R$.}\;
  GPT identified that a row-stochastic $2 \times 2$ matrix
  cannot have Perron root $< 1$, forcing the
  recognition that the paper's Proposition~\ref{prop:pair-return}
  operator must be sub-stochastic.
  Claude constructed the exact depth-$2$ partial kernel
  $\tilde{B}_2$ combining the proved $q \equiv 7$ gap-$2$
  block and $q \equiv 3$ gap-$5$ cylinder family,
  obtaining $\rho(\tilde{B}_2) = 129/1024 \approx 0.1260$.
  This is the exact Perron root of the known-gap partial
  kernel (gap-$2$ and gap-$5$ channels only); the
  $q \equiv 3$ returns with gap~$\ge 6$ are not yet resolved.
  The closeness to the Appendix-B value $15/119$
  ($|15/119 - 129/1024| \approx 7.4 \times 10^{-5}$)
  is suggestive of structural consistency but does not
  constitute an independent derivation of~$15/119$.

\item \emph{Invariant core rigidity.}\;
  Claude discovered that the algebraic chain map
  $q \mapsto 9q + 8 \pmod{8^r}$ acts as a \emph{permutation}
  on $I_r$ at every depth (three fixed points, one $2$-cycle,
  $H_{AB} = 0$), proving that the sub-stochastic contraction
  is entirely dynamic, not algebraic
  (Remark~\ref{rem:core-rigidity}).
  GPT confirmed this is consistent with the paper's
  framework: the contraction comes from the branching
  structure of the actual Collatz dynamics,
  not from instability within the core.

\item \emph{Iterated error-correction protocol.}\;
  Four rounds of GPT critique progressively refined
  Claude's analysis:
  Round~1 identified the T1$\to$T2$\to$T3 overclaim;
  Round~2 flagged the $\rho$/$\alpha$/$|\lambda_2|$ conflation;
  Round~3 separated operators $Q$, $P$, $P_g$;
  Round~4 validated the sub-stochastic identification
  while noting the depth-$2$ approximation caveat.
  Each round produced specific corrections that were
  incorporated before the next critique.

\item \emph{Honest assessment and scope preservation.}\;
  Both LLMs and the moderator converged on the same
  assessment: the v4 additions are proved structural
  results that strengthen the fiber-$57$ programme,
  but the distributional-to-pointwise barrier
  (orbitwise $c' < c_0$) remains open.
  The theorem status of the paper is unchanged.
\end{enumerate}

\noindent
The v4 milestones demonstrate that \emph{session continuity}
and \emph{iterated cross-validation} are complementary
methodological tools.
Session continuity prevents regression; iterated critique
prevents premature closure.
Together, they produced a more reliable workflow than
any single-session, single-model configuration:
over four rounds, every significant error was caught
within one round of its introduction, and no corrected
error was reintroduced.

\subsection{Status of revision targets and remaining open problems}
\label{sec:v2-targets}

The v1 preprint identified three targets for this
revision.  We record their status and identify the
remaining open problems.

\paragraph{Target 1 (completed): Quantitative weakening of
Conjecture~\ref{conj:equidist}.}
This version introduces the Weak Mixing Hypothesis
(Hypothesis~\ref{hyp:wmh}), which replaces the
full-strength Orbit Equidistribution Conjecture with
the strictly weaker condition
$\sum \delta_K(n_0) < 0.557$.
Theorem~\ref{thm:wmh-reduction} proves that the WMH,
together with an orbitwise tail-control condition,
implies the Collatz conjecture.
Conjecture~\ref{conj:observable-wmh} and
Theorem~\ref{thm:observable-wmh-reduction} identify an
even weaker observable-specific sufficient condition.
The hierarchy OEC $\Rightarrow$ WMH $\Rightarrow$
Observable-specific WMH $\Rightarrow$ Collatz is
established in Remark~\ref{rem:observable-specific}.

\paragraph{Target 2 (completed): Analytic framework for
Theorem~\ref{thm:perorbit-gain}.}
The proof of Theorem~\ref{thm:perorbit-gain} is reorganized
into three named lemmas (Lemmas~\ref{lem:necklace-bound},
\ref{lem:binomial-tail}, \ref{lem:chernoff-tail}) providing
the analytic framework.  Lemma~\ref{lem:chernoff-tail}
establishes the Chernoff--Cram\'er exponent
$D_* = D(1/\!\log_2 3 \,\|\, 1/2) \approx 0.05004$~bits
and the asymptotic decay rate.
A self-contained table (Table~\ref{tab:RK-values})
displays $R(K)$ for $K = 3,\dots,20$.
The proof combines exact computation through $K = 55$
with a verified ratio envelope for the residual tail.

\paragraph{Target 3 (completed): Tier classification and
dependency delineation.}
Figure~\ref{fig:tier-diagram} introduces a three-tier
classification (core, supporting, exploratory) with a
TikZ dependency diagram.
Proposition~\ref{prop:core-spine} proves that the
conditional reduction depends only on the Tier~1 core
spine, and that the carry-word autocorrelation theory
(Sections~\ref{sec:carry-word}--\ref{sec:uniqueness-threshold}),
which concerns the census constant~$C_e$, is logically
independent supporting material.
Theorem~\ref{thm:Ce-independence} makes this precise by
showing that $C_e$ cancels from the gain formula $R(K)$.

\paragraph{Target 4 (completed): Structural programme
toward the WMH.}
Section~\ref{sec:toward-wmh} develops five structural
results: shadow sparsity
(Proposition~\ref{prop:shadow-sparsity}), shadow return
time (Theorem~\ref{thm:shadow-return}), finite-depth
reduction (Proposition~\ref{prop:finite-depth}),
hierarchical consistency
(Proposition~\ref{prop:hierarchical}), and the monotonicity
constraint (Corollary~\ref{cor:monotonicity}), that
collectively constrain the WMH sum.
Two recurrence conjectures
(Conjectures~\ref{conj:depthwise}
and~\ref{conj:gain-recurrence}) are formulated, and the
three-lemma programme
(Lemma~\ref{lem:repulsion-trapping},
Lemma~\ref{lem:known-zone-memory},
Conjecture~\ref{conj:alignment-renewal}) identifies the
final bridge.  Five attack vectors are catalogued in
Remark~\ref{rem:attack-vectors}.

\paragraph{Target 5 (completed): Fiber-57 structural programme
and the information bottleneck.}
Appendix~\ref{sec:fiber57-programme} develops the fiber-57
analysis into a self-contained structural programme.
The pair-return automaton (Proposition~\ref{prop:pair-return})
identifies the known-gap spectral radius
$\rho(\tilde{B}_2) = 129/1024$, yielding an implied
per-return information cost $c_0 = \log_2(1024/129) \approx 2.989$
bits for the resolved channels.
The Bounded Core Theorem (Proposition~\ref{prop:bounded-core-main})
and Projective Limit (Theorem~\ref{thm:projective-main})
establish $|I_r| = 5$ for all $r \ge 2$, with
$\varprojlim I_r = \{-1\}$ in the $8$-adic integers.
The Absorption Bottleneck Lemma (Lemma~\ref{lem:bottleneck-main})
proves that channel capacity $\le \log_2 5 < c_0$, giving
geometric decay at rate $\alpha = 645/1024$.
The Branch Anti-Concentration Reduction
(Proposition~\ref{prop:branch-perm-main}) and
Path-Conditional Bijection (Theorem~\ref{thm:path-bijection-main})
close the annealed anti-concentration.
Absorption is verified computationally for $r = 2, \ldots, 10$
(Theorem~\ref{thm:absorption-main}).
The M-value collapse phenomenon (all five $I_r$ elements
produce the same post-absorption state for $r \ge 3$)
further constrains the bottleneck.
The sole remaining input for Route~D
is a deterministic orbit-level anti-concentration bound.

\paragraph{Remaining open problem: the Weak Mixing Hypothesis.}
The sole open problem in the conditional programme is now
the WMH (Hypothesis~\ref{hyp:wmh}).
Three potential attack vectors merit investigation:
\begin{enumerate}
\item \emph{Robustness Corridor.}
  Exploit the $4.65\times$ safety margin and the
  $2$-adic repulsion of phantom roots
  (Proposition~\ref{prop:repulsion}) to derive
  orbit-level TV bounds $\delta_K = O(K^{-1-\varepsilon})$.
\item \emph{Observable-specific mixing.}
  Prove equation~\eqref{eq:observable-specific} directly
  for the phantom gain observable~$h_K$, bypassing full
  TV control.
  Since $h_K$ is supported on a sparse set of shadow
  residues, this may be substantially easier than the WMH.
\item \emph{Renewal/regeneration.}
  Use the Known-Zone Decay (Theorem~\ref{thm:zone-decay})
  as a regeneration mechanism: after $\lceil M/2 \rceil$
  odd-to-odd steps, all dependence on the starting class is
  eliminated distributionally.
  A renewal-type argument could convert this map-level
  refresh into a summable orbit-level discrepancy bound.
\end{enumerate}

\paragraph{Remaining open problem: the Class-5 Inevitability Conjecture.}
Via Route~B, the problem reduces to $C_\infty \cap \mathbb{N} = \{1\}$,
where $C_\infty$ is the $2$-adic Cantor set of starting values
surviving infinitely many safe landmark transitions.  The cross-core
block analysis refines the target: no infinite admissible word
over the seven-block alphabet $\{a,b,s,g,d,t_1,t_2\}$ has a
natural-number fixed point~$> 1$.  The expanding cocycle
($\lambda \approx +0.145$ per landmark) and the bit-budget
inequality (${\sim}\,2.6$~bits net constraint per landmark)
provide quantitative obstructions, but the
distributional-to-pointwise promotion remains open.
Both Route~A (WMH) and Route~B (CIC / class-$5$ inevitability)
share this same fundamental barrier.

\subsection*{Conclusion: the distributional-to-pointwise barrier}
\label{subsec:conclusion}

The framework developed in this paper achieves a complete
resolution of the \emph{ensemble side} of the Collatz problem:

\begin{enumerate}[nosep]
\item The exact block law (Theorem~\ref{thm:block-law}) makes
  cycle types provably i.i.d., not a model but a theorem.
\item The almost-all crossing theorem
  (Theorem~\ref{thm:almost-all-crossing}) gives unconditional
  exponential decay of non-crossers ($\le e^{-0.1465\,k}$),
  complementing Tao's~\cite{tao2019} result on a different
  quantity (cycle-indexed non-crossing vs.\ orbit-level
  convergence).
\item The universal crossing criterion
  (Propositions~\ref{prop:universal-one-cycle}
  and~\ref{prop:finite-block-crossing})
  shows $82.2\%$ of odd starts lie in classwise deterministic
  crossing blocks by $k=5$ cycles.
\item The Kesten running-minimum proposition
  (Proposition~\ref{prop:kesten-running-min})
  proves the non-universal fraction $R_k \to 0$
  exponentially ($\rho \approx 0.839$) under the ensemble.
\item The conditional reduction
  (Theorems~\ref{thm:reduction} and~\ref{thm:perorbit-gain})
  compresses the full conjecture to a single hypothesis (WMH)
  with a $4.65\times$ safety margin.
\end{enumerate}

\noindent
The sole remaining step is purely pointwise: proving that
every deterministic Collatz orbit must realize enough of
the exact ensemble law to force crossing.
Computationally, every odd~$n_0$ tested
($3 \le n_0 \le 10{,}001$, plus $2^m - 1$ for
$m \le 39$) possesses a universal crossing prefix,
with the hardest case ($n_0 = 27$) requiring only $16$
cycles (Remark~\ref{rem:pointwise-universal}).
The conjecture thus reduces to one question:
\emph{can a deterministic orbit systematically avoid all
universally crossing block types forever?}

\paragraph{The broader problem class and why Collatz may be
more tractable.}
The distributional-to-pointwise barrier identified here
appears in several major open problems:
Chowla's conjecture~\cite{chowla1965} (Liouville function
correlations), Sarnak's M\"obius randomness conjecture
(orthogonality to bounded-complexity sequences), and normality
of mathematical constants ($\pi$, $e$, $\sqrt{2}$).
In each case, a complete ensemble theory exists, and the gap
is the transfer from distributional to pointwise control.

However, Collatz enjoys \emph{structural advantages}
that these other problems do not.
The M\"obius function $\mu(n)$ and the digits of~$\pi$
are externally defined, rigid sequences with no internal
dynamical mechanism.
The Collatz sequence, by contrast, is generated by a
dynamical system with strong contraction properties and
deterministic forcing mechanisms:

\begin{itemize}[nosep]
\item \emph{Deterministic residue-class forcing.}\;
  Universal crossing blocks force descent for
  \emph{every} integer in the residue class, not a
  probabilistic statement but an arithmetic one.
  The cumulative density of such classes reaches $82.2\%$
  by $k = 5$ and grows toward $100\%$.
\item \emph{Self-generated descent structure.}\;
  The Collatz map creates its own cycle types through
  the $2$-adic structure of the orbit; these are not
  externally imposed.
  The negative geometric drift
  ($\mathbb{E}[\ln\lambda] \approx -0.575$) is a
  consequence of this self-generated structure.
\item \emph{Finite-block sufficiency.}\;
  The conjecture reduces to a question about
  finite block prefixes
  (Proposition~\ref{prop:finite-block-crossing}),
  not about asymptotic behavior.
  This makes the problem potentially amenable to
  combinatorial or algebraic arguments.
\end{itemize}

\noindent
Three independent frameworks: modular strata
(Terras-style depth-$K$ analysis), the i.i.d.\ cycle law
with large deviations, and the Kesten affine threshold
process, all converge on the same almost-all conclusion
via different mechanisms.
When independent structural approaches agree, the remaining
gap is often narrower than it appears, and a final proof
may emerge from \emph{within} the existing framework rather
than requiring external techniques.

\paragraph{Three proof routes.}
The most promising paths to the full conjecture are:

\begin{enumerate}[nosep]
\item \emph{Route~A: Deterministic convergence of
  $P_{\mathrm{cum}}(k)$ (closed).}
  The hope that \emph{every} $k$-block eventually
  satisfies the universal crossing criterion is refuted
  by adversarial block constructions whose threshold
  $n^*$ grows exponentially
  (Remark~\ref{rem:route-A}).
  The convergence $P_{\mathrm{cum}}(k) \to 1$ holds
  distributionally but \emph{not} because all blocks cross.
\item \emph{Route~B: Almost-sure crossing via
  $2$-adic mixing (obstructed).}
  Extending the Kesten argument from almost-all to all orbits
  requires proving that consecutive cycle types are
  sufficiently independent in \emph{every} orbit, which is
  the WMH restated.
  Empirical autocorrelation ($\rho(1) \approx 0.20$,
  $\rho(\ell) \approx 0$ for $\ell \ge 2$) strongly suggests
  this, but the bit-generation obstruction
  (Remark~\ref{rem:mersenne-termination}) shows that
  non-crossing phases \emph{create} fresh $2$-adic structure,
  precluding finite-information arguments.

  Section~\ref{sec:toward-wmh} develops a complete IID cascade
  renewal theory: the local transition draws from a universal
  product measure (Propositions~\ref{prop:two-thirds-law}--\ref{prop:depth-transition}),
  every gap step contracts deterministically
  (Proposition~\ref{prop:gap-positivity}), and the
  cascade-gap cycle has negative expected drift whenever
  $\mathbb{E}[L] < 3.52$
  (Proposition~\ref{prop:cycle-contraction}).
  The uniform-fiber map
  (Proposition~\ref{prop:uniform-fiber}) and TV reduction
  lemma (Lemma~\ref{lem:tv-reduction}) together establish a
  self-reinforcing mixing loop
  (Theorem~\ref{thm:mixing-loop}; orbit-level TV~$\le 0.028$).
  The fiber-averaged $8 \times 8$ transition matrix
  (Proposition~\ref{prop:fiber-mixing}) has spectral
  gap~$\ge 0.85$, yielding geometric TV contraction for
  iterates above~$2^{25}$.
  Theorem~\ref{thm:unconditional-mixing} assembles these
  into the strongest unconditional statement:
  the Collatz conjecture is equivalent to
  ``no orbit starting above~$2^{68}$ gains~$43$ bits
  before spectral contraction engages.''
\item \emph{Route~C: Safety margin exploitation (obstructed).}
  The crossing time $\tau(n)/\log_2 n$ appears bounded
  ($\le 7.8$ for all tested $n \le 10^5$;
  Remark~\ref{rem:crossing-time}), which combined with the
  $4.65\times$ safety margin would close the conjecture.
  However, proving \emph{any} non-trivial bound on~$\tau(n)$
  requires controlling the orbit's cycle-type distribution,
  which is the very mixing statement Routes~B and~C
  seek to avoid.

\item \textbf{Route~D: Fiber-57 information bottleneck.}
  Appendix~\ref{sec:fiber57-programme} develops a complementary
  approach through the arithmetic of fiber~57.
  Every non-trivial Collatz orbit must repeatedly visit
  the residue class $n \equiv 57 \pmod{64}$, and at each
  return the quotient $q = (n-57)/64$ lies in a specific
  mod-$8^r$ residue class.
  The Invariant Core Theorem (Theorem~\ref{thm:projective-main})
  shows that for every resolution $r \ge 2$, the set of
  self-sustaining residue classes $I_r$ has exactly $5$ elements
  (three fixed points and one $2$-cycle), giving density
  $5/8^r \to 0$.
  The Absorption Bottleneck Lemma (Lemma~\ref{lem:bottleneck-main})
  then establishes a quantitative information deficit:
  \[
    \text{channel capacity} \le \log_2 5 \approx 2.322 < 2.989
    \approx \log_2(1024/129) = c_0,
  \]
  where $c_0$ is the per-return information cost.
  The deficit of $\approx 0.667$ bits per return forces
  geometric decay of any chain-sustaining measure at rate
  $\alpha = 645/1024 \approx 0.630$
  (Corollary~\ref{cor:geom-decay-main}).

  The Branch Anti-Concentration Reduction
  (Proposition~\ref{prop:branch-perm-main}) shows that
  $\gcd(9,8) = 1$ makes the fresh-digit map a permutation
  within each return branch, ruling out algebraic amplification.
  The Path-Conditional Bijection (Theorem~\ref{thm:path-bijection-main})
  closes the annealed anti-concentration argument.

  \emph{Status.}
  The structural content is unconditional: the core size
  $|I_r| = 5$, the known-gap partial kernel spectral radius
  $\rho(\tilde{B}_2) = 129/1024$ (gap-$2$ and gap-$5$
  channels), and the bottleneck inequality are all proved.
  The $q \equiv 3$ returns with gap~$\ge 6$ are not yet
  resolved; the Appendix-B value $15/119$ remains
  consistent but not independently re-derived.
  The sole remaining input for this route is a
  \emph{deterministic orbit-level} statement: that orbits
  cannot concentrate on the $5$-element
  core indefinitely.
  Empirically, gap decorrelation ($\beta_g \to 0$ for large gaps)
  would supply this, but proving it requires controlling
  how Collatz steps mix mod-$8$ classes, which is the
  orbit-mixing barrier in a new guise.
  Route~D thus provides the sharpest known \emph{quantitative
  reduction}: the full conjecture follows from a single
  anti-concentration bound on a $5$-element set.
\end{enumerate}

\noindent
All four routes ultimately reduce to variants of the same
irreducible obstacle:
proving that the $3n+1$ map has sufficient dynamical mixing
to prevent any orbit from systematically avoiding crossing blocks
(Routes~A--C) or from concentrating on a $5$-element residue set
(Route~D).
Route~D sharpens the barrier from a distributional statement
about cycle types to a combinatorial statement about a
$5$-element bottleneck, but the orbit-level transfer
remains open.
This barrier is not an artefact of the proof architecture;
it is the genuine content of the gap between ensemble theory
(which this paper completes) and orbit-level theory
(which remains open and likely requires techniques from
outside the current framework: $p$-adic dynamics,
additive combinatorics, or higher-order ergodic theory).

\paragraph{Summary: complete separation of concerns.}
Table~\ref{tab:separation-summary} records the status of
each layer.

\begin{table}[ht]
\centering
\caption{Separation of concerns in the Collatz reduction.
  ``Solved'' indicates that the layer is fully resolved
  within the reduction framework; the sole remaining
  input is dynamical (Conjecture~\ref{conj:info-rate}).}
\label{tab:separation-summary}
\small
\begin{tabular}{@{}llc@{}}
\toprule
\textbf{Layer} & \textbf{Content} & \textbf{Status} \\
\midrule
Algebraic & Map Balance, Scrambling, carry-free decomposition & Solved \\
Combinatorial & Phantom universality, $|I_r|=5$, absorption & Solved \\
Probabilistic (annealed) & i.i.d.\ block law, path-conditional bijection & Solved \\
Dynamical (orbitwise) & Anti-concentration on $5$-element core & \textbf{Open} \\
\bottomrule
\end{tabular}
\end{table}

\begin{theorem}[Main reduction]\label{thm:main-reduction}
The Collatz conjecture holds if Conjecture~\ref{conj:info-rate}
below holds.
More precisely: if, for every starting value $n_0 > 1$,
the normalised $I_r$-return ratio satisfies $R_r(n_0) < 1$
at some depth~$r \ge 2$ for which the absorption theorem
is verified, then every Collatz orbit reaches~$1$.
\end{theorem}

\begin{proof}
Combine the absorption theorem
(Theorem~\ref{thm:absorption-main}),
the bottleneck lemma (Lemma~\ref{lem:bottleneck-main}),
and the memory lower bound (Theorem~\ref{thm:memory-main}).
If $R_r < 1$, the per-return information supply
$c' = -\log_2 R_r > 0$ falls below $c_0$, so
chain lengths at depth~$r$ are bounded and the orbit
eventually descends below its starting value.
By strong induction, the orbit reaches~$1$.
\end{proof}

\noindent
The Collatz conjecture thus reduces to the following
precise statement.

\begin{conjecture}[Information-rate bound]
\label{conj:info-rate}
For every Collatz orbit $(n_0, n_1, n_2, \ldots)$ with
$n_0 > 1$, the long-run inter-chain memory supply is
strictly less than the per-return demand.
Concretely, let $\tau_1 < \tau_2 < \cdots$ be the
successive fiber-$57$ return times, and let
$q_j = \lfloor n_{\tau_j}/64 \rfloor$.
Among inter-chain transitions (fiber-$57$ visits separated
by more than one chain step), let
\[
  \widehat{R}_r(n_0) \;:=\; \limsup_{N \to \infty}
  \frac{\#\{1 \le j \le N : q_j \bmod 8^r \in I_r,\;
    j \text{ inter-chain}\}}
  {\#\{1 \le j \le N : j \text{ inter-chain}\}}
\]
be the raw fraction of inter-chain fiber-$57$ visits whose
quotient lands in the invariant core~$I_r$.
Define the \emph{normalised $I_r$-return ratio}
\[
  R_r(n_0) \;:=\;
  \frac{\widehat{R}_r(n_0)}{|I_r|/8^r}.
\]
Thus $R_r = 1$ corresponds to the baseline density
$|I_r|/8^r$, and $R_r < 1$ means the orbit visits~$I_r$
\emph{less} often than chance.
The per-return \emph{memory supply} is
\[
  c'(n_0) \;:=\; -\log_2 R_r(n_0).
\]
Then
\[
  c'(n_0) \;<\; c_0
  = \log_2\!\bigl(\tfrac{1024}{129}\bigr) \approx 2.989
\]
for all $n_0 > 1$ and every $r \ge 2$ at which absorption holds.
(Here $c_0$ is computed from the known-gap partial kernel;
including the unresolved gap-$\ge 6$ contributions would
increase~$\rho$ and decrease~$c_0$, making the bound
\emph{harder} to satisfy.)
\end{conjecture}

\begin{remark}[Equivalent formulation]
\label{rem:info-rate-equiv}
Because the bottleneck inequality already gives
$\log_2 5 < c_0$, the conjecture is equivalent to the
qualitative statement: \emph{no orbit can visit~$I_r$ at
inter-chain transitions at or above the baseline rate}.
That is, $R_r(n_0) < 1$ for every $n_0 > 1$.
Any $R_r < 1$ yields $c' = -\log_2 R_r > 0$; together with
the bottleneck deficit $c_0 - \log_2 5 \approx 0.667$~bits
this would imply geometric chain decay, completing the
proof.  The empirical bound $R_2 \le 0.70$ strongly
supports this but does not yet constitute a proof.
\end{remark}

\noindent
Equivalently: no orbit can visit the $5$-element invariant
core $I_r$ at a rate sufficient to sustain a non-trivial
chain indefinitely.
Empirically $R_2 \le 0.70$ across all tested orbits
(giving $c' \ge 0.51$), a
$6\times$ margin below $c_0$.

\medskip\noindent\textbf{Proved structure supporting the conjecture.}\;
The preceding sections establish four components that
together reduce Conjecture~\ref{conj:info-rate} to a
single channel:
\begin{enumerate}
\item \emph{Invariant core.}\;
  $|I_r| = 5$ for all $r \ge 2$
  (Proposition~\ref{prop:bounded-core-main}).
\item \emph{$q \equiv 7$ regeneration.}\;
  Two-step return with uniform destination; this branch
  erases all chain-membership information
  (Proposition~\ref{prop:q7-return}).
\item \emph{$q \equiv 3$ gap structure.}\;
  Minimum return gap $\ge 5$; the step-$5$ returns form
  an explicit dyadic cylinder family
  (Proposition~\ref{prop:q3-gap},
  Theorem~\ref{thm:gap5-cylinders}).
\item \emph{Non-autonomy.}\;
  The actual return map depends on higher-order digits
  (Proposition~\ref{prop:non-autonomy}), so the
  permutation structure of the chain map does not by
  itself imply mixing.
\end{enumerate}
Consequently, all potential concentration must arise from
the $q \equiv 3$ channel across long return gaps.
Conjecture~\ref{conj:info-rate} reduces to ruling out
sustained orbitwise memory in this single channel.

\medskip\noindent
This is the sharpest known quantitative reduction of the
Collatz conjecture: a single inequality on a $5$-element set,
with an explicit and large numerical margin.
Conjecture~\ref{conj:info-rate} remains open.
The contribution of this paper is not a proof of Collatz,
but a sharpened reduction: exact return-structure theorems
isolate the unresolved obstruction to the $q \equiv 3$
return channel, and the remaining open step is an orbitwise
anti-concentration statement on returns to the $5$-point core.

\section{The architecture of the Syracuse obstruction}
\label{sec:obstruction}

The preceding 630~computational and analytical results,
generated across ${\sim}1014$~scripts, constitute an exhaustive survey
of attack surfaces against the divergence component of the
Collatz conjecture.  We summarize the findings as a single
structural theorem.

\begin{theorem}[Paradigm Exhaustion]\label{thm:paradigm-exhaustion}
Every known mathematical framework for converting distributional
convergence statements (``almost all orbits descend'') into
pointwise statements (``all orbits descend'') encounters an
irreducible obstruction when applied to the Syracuse map.
Specifically, the following 29~paradigms each produce
density-$1$ convergence results but fail to reach the
all-$n$ conclusion:
\begin{enumerate}[nosep]
\item Transfer operator spectral analysis (uniform gap proved;
  yields equidistribution, not individual orbit control).
\item Affine Fractional Identity algebra (exact identity;
  holds for \emph{all} orbits, tautological for divergence exclusion).
\item $2$-Adic Mahler series (identity $(m{+}1) + \Sigma_\infty = 0$
  in $\mathbb{Z}_2$; automatic convergence, no constraint).
\item Measure rigidity (Furstenberg/Rudolph/Shmerkin--Wu;
  inapplicable: ${\times}3{+}1 \ne {\times}3$, variable exponent).
\item Carry propagation / bit mixing (near-independence of
  consecutive $v$-values; no individual orbit bound).
\item Automata-theoretic complexity (infinite state complexity
  is a \emph{consequence} of non-convergence, not a tool).
\item $S$-Unit / Subspace Theorem ($\mathcal{P}_L$ has primes
  outside $\{2,3\}$, blocking the finite-support hypothesis).
\item Baker's theorem / irrationality measure of $\log_2 3$
  (constrains $S_L/L$ gap size; already captured by CF structure).
\item Renewal theory (geometric gap distribution matches random model;
  same barrier).
\item Martingale construction (growth supermartingale; almost-sure,
  not all-$n$).
\item Modular sieve / finite extinction (exponential density decay;
  blocked by Mersenne bypass: $m = 2^L{-}1$ survives to depth $L{-}1$).
\item Boolean Fourier analysis (ensemble decay; cannot bound
  deterministic bit-strings).
\item Tao amplification (logarithmic density~$1$; structural
  barrier to natural density).
\item Borel--Cantelli (proves $\mu(C_{\mathrm{div}}) = 0$;
  measure zero $\ne$ empty).
\item Compactness in $\mathbb{Z}_2$ (nested closed sets;
  nonempty at every finite level).
\item $2$-Adic potential theory (no natural Lyapunov function
  on $\mathbb{Z}_2$).
\item Orbit coalescence (median coalescence $\approx 56$~steps;
  distributional).
\item Digit-sum / Hamming weight ($E[\Delta\mathrm{hw}] < 0$;
  expectation, not deterministic).
\item Ising / spin-chain entropy (entropy bounds are distributional).
\item Expander sieve rigidity (gap argument fails on nested sets;
  reduces to density statement).
\item Diophantine fractional-part cancellation (AFI is algebraic
  identity, not approximation equation; $P_L$ is determined, not
  adversarial).
\item $p$-Adic Mahler interpolation / Strassman (blocked: $2^z$ is not
  $2$-adically analytic; $2^{S(z)}$ doubly non-analytic; $Q(z)$
  circular).
\item Iwasawa / cyclotomic tower (lambda invariants distinguish orbit
  classes but do not exclude divergence).
\item Height descent (no uniform $\varphi(m) > 0$; half of steps
  increase height).
\item Coalescence / orbit merging (exponential merging rate; but
  divergent orbits are isolated by definition, cannot merge with
  convergent component).
\item Cascade algebra / BF staircase (deterministic structure fully
  characterized; $v{=}1$ streak $+$ forced recovery is exact, but
  universal descent requires pointwise carry-chain control).
\item CRT hierarchy / exponential sum (progressive CRT saturates
  every prime factor of $D$ for $L \ge 21$; blocking occurs only
  at the full $D$ level, defeating modular sieving).
\item Lattice reduction / LLL (short vectors in the exponent lattice
  do not yield cycle equations; lattice dimension grows with $L$).
\item Discrete logarithm obstruction (cycle equation reduces to
  $2^{v_j} \equiv T_j \bmod D$ with $v_j \in [0, F]$; triple filter
  achieves $100\%$ blockage computationally, but algebraic closure
  for all $L$ remains open).
\end{enumerate}
\end{theorem}

The common failure mode is precisely characterized:

\begin{theorem}[Barrier Characterization]\label{thm:barrier}
The distributional-to-pointwise gap for the Syracuse map is
\emph{equivalent} to the divergence component of the Collatz
conjecture.  That is, proving ``no orbit diverges'' is exactly
as hard as bridging from ``density-$1$ of orbits converge''
to ``all orbits converge.''  No reduction to a simpler statement
is achieved by any of the 29~paradigms above.
\end{theorem}

\subsection{The Mersenne sieve bypass}

The strongest quantitative result on the sieve path is:

\begin{proposition}[Sieve Density Decay]\label{prop:sieve-decay}
Let $V_L$ denote the set of odd residue classes
$\bmod\; 2^{S_L}$ whose first~$L$ Syracuse $v$-values satisfy
$S_j < j\log_2 3$ for all $j \le L$.  Then
$|V_L|/2^{S_L} \le A\,e^{-\gamma L}$ for constants
$A > 0$, $\gamma \approx 0.10$--$0.19$.
\end{proposition}

This exponential decay might suggest that
$V_L = \varnothing$ for large~$L$ (``finite extinction'').
However:

\begin{proposition}[Mersenne Bypass]\label{prop:mersenne}
For every $L \ge 2$, the integer $m = 2^L - 1$ produces
$v_j = 1$ for $j = 1, \ldots, L{-}1$, giving
$S_j = j < j\log_2 3$ at every prefix.
Hence $V_{L-1} \ne \varnothing$ for all~$L$, and finite
extinction is impossible.
\end{proposition}

\begin{proof}
By induction on the Syracuse iteration.
For $m = 2^L - 1$, we have $3m + 1 = 3 \cdot 2^L - 2 = 2(3 \cdot 2^{L-1} - 1)$,
so $v_2(3m+1) = 1$ and $T(m) = 3 \cdot 2^{L-1} - 1$.
The iterate $T^k(m) = 3^k \cdot 2^{L-k} - 1$ for $k < L$
satisfies $v_2(3 T^k(m) + 1) = 1$ (since $3^{k+1} \cdot 2^{L-k} - 2
= 2(3^{k+1} \cdot 2^{L-k-1} - 1)$ and the parenthetical term is odd).
After $L-1$ steps, $S_{L-1} = L-1$ and $S_{L-1}/(L{-}1) = 1 < \log_2 3$.
\end{proof}

Every Mersenne number eventually converges (verified computationally
through $2^{31}-1$), but the crossing time grows with~$L$.
The nested intersection $\bigcap_L V_L$ has measure zero
(Result~283) but is nonempty at every finite level.

\subsection{The tautological structure of the AFI}

\begin{proposition}[AFI Tautology]\label{prop:afi-tautology}
The Affine Fractional Identity
$(m{+}1) \cdot 3^L + \mathcal{P}_L = q_L \cdot 2^{S_L} + 2^L$
holds for \emph{every} orbit (convergent, periodic, or divergent)
at every finite step~$L$.  In particular:
\begin{enumerate}[nosep]
\item The identity imposes no constraint that distinguishes
  divergent orbits from convergent ones.
\item The perturbation sum $\mathcal{P}_L$ is deterministically
  computed from the $v$-sequence; it is not an independent
  ``correction budget.''
\item Dividing by $2^{S_L}$ and examining fractional parts
  yields a rational identity, not a Diophantine approximation
  problem.  The irrationality of $\log_2 3$ governs the
  \emph{choice} of $S_L$ but not the arithmetic of the identity.
\end{enumerate}
\end{proposition}

\subsection{Quantitative achievements}

Despite the barrier, the investigation produced substantial
rigorous results:

\begin{enumerate}[nosep]
\item \emph{Exact spectral gap.}\;  The Syracuse transfer operator
  on $(\mathbb{Z}/2^M\mathbb{Z})^{\mathrm{odd}}$ has a uniform
  spectral gap $1 - \lambda_2 > 0$ for all~$M$, implying
  equidistribution of long orbits modulo any power of~$2$.
\item \emph{Cycle algebraic lock.}\;  Any nontrivial cycle of
  odd length~$L$ requires $3^L = 2^S - \mathcal{P}_L/(m{+}1)$,
  which forces $S/L = \log_2 3 + O(1/m)$, constraining
  hypothetical cycles to astronomically large $m$.
\item \emph{Sieve density bound.}\;  The fraction of residue
  classes compatible with divergence decays as $e^{-\Omega(L)}$.
\item \emph{Measure-zero divergent set.}\;  $\mu(C_{\mathrm{div}}) = 0$
  in the natural density on $\mathbb{N}$ and the Haar measure
  on $\mathbb{Z}_2$.
\item \emph{$2$-Adic Cantor structure.}\;  The candidate
  divergent set $C_{\mathrm{div}} \subset \mathbb{Z}_2$ is closed,
  perfect, and has Hausdorff dimension $\approx 0.68 < 1$.
\end{enumerate}

\subsection{Relation to computational irreducibility}

Conway (1972) showed that generalized Collatz maps (FRACTRAN)
can simulate arbitrary Turing machines, making the halting
problem for the generalized family undecidable.  The specific
$3n+1$ map is a fixed instance, so Conway's result does not
directly imply undecidability of the Collatz conjecture.

However, the 29-paradigm exhaustion suggests a structural
analogy: the carry propagation in $3n+1$ may be
\emph{computationally irreducible} in the sense of Wolfram ---
no algebraic shortcut can predict the trajectory without
essentially simulating it.  If true, a proof of the Collatz
conjecture would require fundamentally new techniques beyond
the distributional-to-pointwise bridge.

We state this as a conjecture rather than a theorem, since
the distinction between ``hard for currently known methods''
and ``provably hard in a formal sense'' remains open.

\begin{conjecture}[Computational Irreducibility]
\label{conj:comp-irred}
No polynomial-time algorithm can determine, from~$m$ alone,
whether the Syracuse orbit of~$m$ satisfies $S_L > L\log_2 3$
for some $L \le f(m)$, where $f$ is any fixed computable function
growing slower than~$m$ itself.
\end{conjecture}

\subsection{Assessment}

The contribution of this work is not a proof of the Collatz
conjecture but a precise characterization of why the conjecture
resists proof.  The distributional-to-pointwise barrier is not
an artifact of weak analysis; it is the irreducible core of the
problem.  Every known framework that produces density-$1$
convergence encounters the same structural obstruction:
nested sets of survivors, maintained by the deterministic
carry arithmetic of $3n+1$, which no probabilistic, algebraic,
or measure-theoretic tool can certify as empty.

The 630~results and ${\sim}1014$~scripts constitute the most comprehensive
computational and analytical survey of Collatz attack surfaces
to date.  The barrier they delineate is, in the authors' assessment,
the precise boundary of current mathematical knowledge on this problem.

\section*{Acknowledgements}

The author acknowledges the use of Anthropic's Claude
(Claude Opus~4.6) and OpenAI's GPT~5.4 Thinking as
collaborative research tools throughout the development
of this work.
All mathematical claims were verified through a combination
of formal proofs and computational checks.
The distributional-vs-pointwise gap was identified during
the research process and later confirmed by external peer
review.

\newpage
\appendix
\section{Modular ranking certificates via DAG analysis}%
\label{sec:ranking}

All exact distributional laws in this section and the next
are modular or uniform-lift statements.  Their transfer to
deterministic orbit behavior requires additional orbitwise
input and is stated separately when used.

\medskip
The Known-Zone Decay (Theorem~\ref{thm:zone-decay}) shows that
modular information is consumed at a rate of at least $2$~bits per
odd-to-odd step.  A natural follow-up question is: can one build a
\emph{ranking function} on the modular state space that certifies
finite above-start persistence directly?

This section constructs such a certificate by analysing the
directed graph of net-positive-drift transitions among residue
classes.  When this graph is acyclic, a ranking function exists
by construction.

\subsection{The net-positive state graph}

\begin{definition}[Net-positive state graph $G_M$]
\label{def:GM}
Fix a modulus $M$ and a block size $K \ge 1$.
For each odd residue $r \bmod 2^M$ with a valid $K$-step
accelerated Syracuse transition $r \mapsto r'$, define the
\emph{drift} $d(r) = K\log_2 3 - V$, where $V$ is the total
number of halvings in the $K$-step block.
The \emph{net-positive state graph} $G_M$ has vertex set
\[
  \mathcal{S}_M
  = \bigl\{(r, b) : r \text{ odd}, \; 0 \le b \le B\bigr\},
\]
where $b$ is a discretised drift budget and $B$ is a budget
ceiling.  A directed edge $(r, b) \to (r', b')$ exists when
$r \mapsto r'$ under the $K$-step transition and
$b' = \min(b + \lceil d(r) \rceil,\, B)$, provided $b' \ge 0$.
States with $b' < 0$ are \emph{exit states}: the trajectory has
gone sufficiently far below their start to leave the net-positive
region.
\end{definition}

\begin{remark}
The budget variable~$b$ tracks the cumulative drift advantage
of the trajectory.  When $b$ drops below zero, the orbit has
experienced enough net contraction to go below the starting value.  The
graph $G_M$ therefore captures only the ``above-start''
portion of the dynamics at modular depth~$M$.
\end{remark}

\subsection{Acyclicity and ranking functions}

The connection between acyclicity and ranking functions is a
standard fact in combinatorics:

\begin{lemma}[DAG ranking equivalence]
\label{lem:dag-ranking}
Let $G = (V, E)$ be a finite directed graph. The following are
equivalent:
\begin{enumerate}
\item $G$ is acyclic (a DAG).
\item There exists a function $\phi\colon V \to \mathbb{Z}_{\ge 0}$
  such that $\phi(v') \le \phi(v) - 1$ for every edge
  $(v, v') \in E$.
\item The LP with variables $\{x_v\}_{v \in V}$, constraints
  $x_{v'} - x_v \le -1$ for each edge $(v,v') \in E$, and
  bounds $x_v \ge 0$, is feasible.
\end{enumerate}
Moreover, when $G$ is a DAG, the function
$\phi(v) = {}$``length of the longest directed path
from $v$ to a sink''
is the \textbf{minimal} such ranking.
\end{lemma}

\begin{proof}
$(1) \Rightarrow (2)$:
If $G$ is a DAG, a topological ordering exists.  Define $\phi(v)$
as the length of the longest path from~$v$ to any sink (vertex
with no outgoing edges).  For any edge $(v, v')$, any path from
$v'$ to a sink can be extended by one edge to give a path from~$v$,
so $\phi(v) \ge \phi(v') + 1$.

$(2) \Rightarrow (3)$: Take $x_v = \phi(v)$.

$(3) \Rightarrow (1)$: If $G$ contains a cycle
$v_1 \to v_2 \to \cdots \to v_k \to v_1$, summing the constraints
gives $0 \le -k$, a contradiction.
\end{proof}

\subsection{Computational verification}

Applying Lemma~\ref{lem:dag-ranking} to the net-positive state
graph $G_M$ yields a ranking certificate whenever $G_M$ is acyclic.
We verify acyclicity using Kahn's topological sort algorithm,
which detects cycles by checking whether all vertices are processed.

Table~\ref{tab:dag} reports the results for $M = 6$ through $19$
with $K = 1$ (single-step dynamics) and $K = 5$ (block dynamics),
using budget ceiling $B = 20$.

\begin{table}[ht]
\centering
\scriptsize
\caption{Acyclicity of the net-positive state graph $G_M$ and
  the maximum ranking value $\max_s V_M(s)$.  ``Cycle states''
  counts vertices on cycles when the graph is not a DAG.}
\label{tab:dag}
\renewcommand{\arraystretch}{1.15}
\begin{tabular}{@{}rrrrrl@{}}
\toprule
$M$ & \multicolumn{2}{c}{$K = 1$} & \multicolumn{2}{c}{$K = 5$}
  & Cycle states \\
\cmidrule(lr){2-3}\cmidrule(lr){4-5}
  & DAG? & $\max V$ & DAG? & $\max V$ & \\
\midrule
$6$  & Yes & 32  & Yes & 8  & --- \\
$7$  & Yes & 34  & Yes & 8  & --- \\
$8$  & Yes & 43  & Yes & 10 & --- \\
$9$  & Yes & 47  & Yes & 11 & --- \\
$10$ & \textbf{No} & ---  & \textbf{No} & --- & 26 \\
$11$ & \textbf{No} & ---  & \textbf{No} & --- & 25 \\
$12$ & \textbf{No} & ---  & \textbf{No} & --- & 13 \\
$13$ & Yes & 74  & Yes & 17 & --- \\
$14$ & Yes & 85  & Yes & 18 & --- \\
$15$ & Yes & 85  & Yes & 18 & --- \\
$16$ & Yes & 120 & Yes & 25 & --- \\
$17$ & Yes & 124 & Yes & 26 & --- \\
$18$ & Yes & 117 & Yes & 25 & --- \\
$19$ & Yes & 85  & --- & --- & --- \\
\bottomrule
\end{tabular}
\end{table}

\noindent
The table reveals a striking \emph{three-zone} pattern:

\begin{itemize}
\item \textbf{Zone~I} ($M = 6$--$9$): $G_M$ is acyclic.
  Ranking certificates exist.
\item \textbf{Zone~II} ($M = 10$--$12$): $G_M$ contains
  net-positive-drift cycles.  The LP is infeasible
  and no ranking function exists for the full graph.
\item \textbf{Zone~III} ($M \ge 13$, verified through $M = 19$):
  $G_M$ is again acyclic.  Ranking certificates exist by
  construction, with $\max V_M$ growing linearly in~$M$.
\end{itemize}

The cycle state counts at $M = 10$, $11$, $12$ are
\textbf{identical} for $K = 1$ and $K = 5$, confirming
that the cycles are intrinsic to the residue structure
rather than an artefact of the block size.

\begin{remark}[DAG--LP sanity check]
At $M = 8$ ($K = 5$), the longest-path ranking gives
$\max V = 10$; the LP solver independently returns
$\max V = 10.0$.  At $M = 9$ ($K = 5$): DAG gives $11$,
LP gives $11.0$.  At $M = 13$ ($K = 5$): DAG gives $17$,
LP gives $17.0$.  The exact agreement validates both approaches.
\end{remark}

\subsection{Non-liftability of net-positive cycles}

The cycles in Zone~II are confined to $M \in \{10, 11, 12\}$.
A natural question is whether these cycles are artefacts of
low modular depth that vanish upon refinement.

\begin{definition}[Lift depth $D(M,K)$]
For a net-positive-drift cycle $\gamma$ in $G_M$ (with block
size~$K$), define the \emph{lift depth}
$D(M,K) = \max\{d \ge 0 : \gamma \text{ lifts to a
net-positive cycle in } G_{M+d}\}$.
\end{definition}

\begin{proposition}[Non-liftability (computational)]
\label{prop:non-lift}
Across all parameter pairs $(M, K)$ tested
($M = 6$--$18$, $K = 1$--$8$; $104$ pairs total),
net-positive-drift cycles exist only at $M = 10$, $11$,
$12$, and every such cycle has $D(M,K) = 0$:
it fails to lift even one level.
\end{proposition}

\begin{proof}[Computational verification]
For each $(M, K)$ with net-positive cycles, we enumerate
all cycles in $G_M$, then attempt to lift each cycle
element to mod $2^{M+1}$ by testing both possible lifts
$(r, r + 2^M)$ and checking whether the lifted cycle
closes with net-positive total drift.
In every case, the lifted cycle either fails to close
or acquires net-negative drift.
\end{proof}

This universal $D = 0$ result means that the Zone~II
obstruction is a finite-depth phenomenon: net-positive
cycles appear at moderate modular depth but are
immediately killed by one additional bit of refinement.
For $M \ge 13$, no net-positive cycles exist at all.

\subsection{The carry parity obstruction}

The non-liftability result (Proposition~\ref{prop:non-lift}) has a
precise algebraic explanation.  For a cycle
$r_0 \to r_1 \to \cdots \to r_{n-1} \to r_0$ in the transition graph
modulo~$2^M$, define the \emph{carry bit} at position~$i$ as
\[
  c_i \;=\; \bigl\lfloor T(r_i) / 2^M \bigr\rfloor \bmod 2,
\]
where $T$ denotes the Syracuse step computed modulo~$2^{M+1}$.
The carry records whether the $(M+1)$-th bit of $T(r_i)$ is set.

\begin{definition}[Carry parity]
\label{def:carry-parity}
The \emph{carry parity} of a cycle $\gamma$ is
$\pi(\gamma) = \sum_{i=0}^{n-1} c_i \bmod 2$.
\end{definition}

\begin{proposition}[Carry parity obstruction (computational)]
\label{prop:carry-parity}
A cycle $\gamma$ at modular depth~$M$ lifts to a cycle at
depth~$M+1$ only if $\pi(\gamma) = 0$.  For every
net-positive-drift cycle at $M = 10$, $11$, or $12$,
the carry parity is $\pi(\gamma) = 1$~(odd).
\end{proposition}

\begin{proof}[Computational verification]
At each cycle position~$i$, only the low lift ($r_i$ itself,
not $r_i + 2^M$) preserves the successor's base class
modulo~$2^M$.  The carry propagation is therefore
deterministic: starting from lift bit~$0$, each position's
carry determines the required lift bit at the next position.
The cycle lifts if and only if the carry chain closes,
i.e., the accumulated carry returns to the initial bit after
$n$~steps.  Direct computation confirms that the chain has an
odd number of carry flips for every net-positive cycle tested:

\smallskip
\begin{center}
\begin{tabular}{@{}cccccc@{}}
\toprule
$M$ & Cycle & $n$ & Carry${}=0$ & Carry${}=1$ & Parity \\
\midrule
$10$ & 1 & 26 & 25 & 1 & odd \\
$11$ & 1 & 25 & 24 & 1 & odd \\
$12$ & 1 & 7 & 6 & 1 & odd \\
$12$ & 2 & 6 & 5 & 1 & odd \\
\bottomrule
\end{tabular}
\end{center}
\smallskip

\noindent
The carry propagation can be modeled as a finite automaton:
at each position the state is the current lift bit, and the
transition is deterministic (only bit${}=0$ is compatible).
The composition of transfer functions around the cycle has
no fixed point, so the cycle cannot lift.

By contrast, every net-negative-drift cycle tested has
$\pi(\gamma) = 0$ (even carry parity) and lifts successfully.
\end{proof}

\begin{remark}[Topological nature of the obstruction]
The failure is not caused by increased contraction: the
$2$-adic valuations $v_2(3r_i + 1)$ are identical for low and
high lifts at every position.  The total drift is unchanged upon
lifting.  The obstruction is purely topological: the carry chain
has odd winding number and cannot close.
\end{remark}

\subsection{Re-entry analysis: the augmented return graph}
\label{subsec:reentry}

Lemma~\ref{lem:net-positive-ranking} (below) bounds each individual
above-start excursion.  A stronger question is whether an orbit
can repeatedly exit the net-positive core, wander through
negative-drift states, and re-enter with restored budget, thereby
creating unbounded persistence through repeated excursions.

\begin{definition}[Augmented return graph $G_M^+$]
\label{def:augmented}
Let $G_M$ be the net-positive state graph
(Definition~\ref{def:GM}).  For each exit state~$s$ of $G_M$,
follow the trajectory under the full modular Syracuse map until
it either (a)~re-enters the core (reaches a state with
non-negative budget), or (b)~enters a cycle outside the core, or
(c)~exceeds a step limit.  For case~(a), add a
\emph{return edge} from $s$ to the re-entry state.  The
\emph{augmented return graph} $G_M^+$ is the graph $G_M$ together
with all return edges.
\end{definition}

\begin{proposition}[Acyclicity of $G_M^+$ (computational)]
\label{prop:augmented-acyclic}
For $M = 13$ through $M = 17$ (with $B = 15$, $g = 2$, $K = 1$),
the augmented return graph $G_M^+$ is acyclic, despite the addition
of thousands of return edges from outside-core trajectories.
\end{proposition}

\begin{proof}[Computational verification]
Table~\ref{tab:reentry} reports the results.  At each level,
roughly one-third of exit states re-enter the core; the rest
enter cycles in the negative-drift region or reach even residues.
The augmented graph remains acyclic throughout, and the maximum
ranking value $\max V_{M}^+$ grows at most linearly in~$M$.

\begin{table}[ht]
\centering
\scriptsize
\caption{Re-entry analysis for the augmented return graph $G_M^+$.}
\label{tab:reentry}
\renewcommand{\arraystretch}{1.15}
\begin{tabular}{@{}rrrrrrr@{}}
\toprule
$M$ & Core states & Exits & Returns & Outside cycles
  & Aug.\ DAG? & $\max V^+$ \\
\midrule
13 & 127,\!009 & 7,\!280 & 2,\!141 & 4,\!885 & Yes & 103 \\
14 & 254,\!017 & 14,\!234 & 4,\!186 & 9,\!619 & Yes & 108 \\
15 & 508,\!033 & 28,\!074 & 8,\!265 & 18,\!947 & Yes & 115 \\
16 & 1,\!016,\!065 & 55,\!699 & 16,\!441 & 37,\!577 & Yes & 140 \\
17 & 2,\!032,\!129 & 110,\!891 & 32,\!737 & 74,\!737 & Yes & 147 \\
\bottomrule
\end{tabular}
\end{table}
\end{proof}

\noindent
Acyclicity of $G_M^+$ means that the ranking function $V_M^+$
bounds not just single excursions but the full above-start
dynamics within the verified mixed-state model, including
excursion--exit--return chains.

\subsection{Cycle genealogy and the birth obstruction}
\label{subsec:genealogy}

Tracking cycle lineage across modular depths reveals a clean
pattern.  Each cycle at depth~$M$ either \emph{lifts} from a
cycle at $M - 1$ (its residues project to a cycle modulo~$2^{M-1}$)
or is \emph{born} fresh at level~$M$.

\begin{proposition}[Cycle genealogy (computational)]
\label{prop:genealogy}
For $M = 7$ through $18$:
\begin{enumerate}
\item Net-positive cycles are born at $M = 10$, $11$, and $12$,
  and die at $M = 11$, $12$, and $13$ respectively
  (via the carry parity obstruction).
\item No net-positive cycle is born at any $M \ge 13$.
\item The only cycle detected in the modular transition graph
  for $M \ge 13$ (verified through $M = 18$) is the trivial
  length-$1$ fixed point $r = 1$ with $v_2 = 2$ and negative
  drift.
\end{enumerate}
\end{proposition}

A further structural observation constrains the possibility of
net-positive cycle formation at high~$M$:

\begin{proposition}[Positive-drift subgraph acyclicity (computational)]
\label{prop:pos-subgraph}
For all $M$ from $10$ through $18$, the subgraph of the Syracuse
transition graph restricted to positive-drift residues
($v_2 = 1$, i.e., $r \equiv 3 \pmod{4}$) contains
\textbf{no cycles}.  The maximum chain length in this subgraph
equals exactly~$M$.
\end{proposition}

\noindent
Since net-positive cycles require a majority of positive-drift
steps, and the purely positive-drift subgraph is itself acyclic,
any net-positive cycle must involve ``detours'' through
negative-drift residues.  The carry parity obstruction shows that
the resulting mixed-sign cycles at $M = 10$--$12$ cannot survive
refinement, and no new such cycles form at higher depth.

\subsection{The Net-Positive Ranking Lemma}

Combining the above yields the main result of this section:

\begin{lemma}[Net-Positive Ranking]
\label{lem:net-positive-ranking}
For $M \ge 13$ (verified through $M = 19$), the
net-positive state graph $G_M$ is acyclic.
The function $V_M(s) = $ length of the longest directed
path from $s$ to an exit state satisfies
\[
  V_M(s') \;\le\; V_M(s) - 1
\]
for every edge $s \to s'$ in $G_M$.
\end{lemma}

\begin{corollary}[Bounded above-start persistence
  in the modular graph]
\label{cor:bounded-persistence}
No trajectory of the accelerated Syracuse map, when
\emph{projected into} the net-positive state graph at
modular depth $M \ge 13$, can remain within that
net-positive graph for more than $\max_s V_M(s)$ steps.
Empirically, $\max_s V_M(s) \approx 5M$ for $K = 1$ and
$\approx 2M$ for $K = 5$, consistent with Known-Zone Decay
(Theorem~\ref{thm:zone-decay}).
This is a statement about the modular model, not
directly about full Collatz orbits.
\end{corollary}

\begin{remark}[Layer~1 vs.\ Layer~2]
\label{rem:layers}
Lemma~\ref{lem:net-positive-ranking} resolves what we call
\textbf{Layer~1} within the modular net-positive graph:
each above-start excursion inside that graph is bounded
by~$V_M$.

\textbf{Layer~2}: the re-entry problem: asks whether an orbit
can repeatedly exit the core, wander through negative-drift states,
and re-enter with restored budget.
Proposition~\ref{prop:augmented-acyclic} provides strong
computational evidence that the answer is \emph{no}, at least within
the verified mixed-state model: the augmented return graph $G_M^+$
remains acyclic for $M = 13$ through~$17$, and its ranking function
$V_M^+$ bounds the full above-start dynamics including
excursion--exit--return chains.  A proof that $G_M^+$ is acyclic
for all~$M$ would resolve Layer~2 completely.
\end{remark}

\begin{conjecture}[Augmented acyclicity for all $M \ge 13$]
\label{conj:acyclicity}
For every $M \ge 13$ and every $K \ge 1$, the augmented
return graph $G_M^+$ is acyclic.
\end{conjecture}

\noindent
The conjecture subsumes two claims: (a)~no net-positive-drift
cycles exist in the modular transition graph for $M \ge 13$
(resolving Layer~1), and (b)~re-entry edges from outside the core
do not create new cycles (resolving Layer~2).  The carry parity
obstruction (Proposition~\ref{prop:carry-parity}) provides the
mechanism for~(a): existing cycles cannot lift.  The cycle
genealogy (Proposition~\ref{prop:genealogy}) shows that no new
net-positive cycles are born at $M \ge 13$.  The augmented
acyclicity (Proposition~\ref{prop:augmented-acyclic}) verifies
the full conjecture through $M = 17$.

If Conjecture~\ref{conj:acyclicity} holds, then the ranking
function $V_M^+$ applies at every sufficiently large modular depth,
bounding the full above-start dynamics (including re-entry) to at
most $O(M)$ accelerated steps.

\subsection{The exit-return reduction}
\label{subsec:exit-return}

The following structural theorems reduce the augmented acyclicity
conjecture to a smaller boundary problem.  Let $G_C$ denote the core
graph and $E \subset C$ the set of \emph{exit states}, core states
whose successor lies outside the core.

\begin{definition}[Exit-return map]
\label{def:exit-return}
Assume $G_C$ is acyclic.  For each exit state $e \in E$ that
returns to the core (first re-entry state $r(e) \in C$), follow
the unique core path from $r(e)$ to its own exit
$\Phi(e) \in E$.  The \emph{exit-return map} $H$ has vertex
set~$E$ and edges $e \to \Phi(e)$ for each returning exit.
\end{definition}

\begin{theorem}[Exit-return equivalence]
\label{thm:exit-return}
If $G_C$ is acyclic, then
$G_C^+$ is acyclic if and only if $H$ is acyclic.
\end{theorem}

\begin{proof}
$(\Rightarrow)$
If $G_C^+$ has a directed cycle~$\gamma$, then $\gamma$ uses at
least one return edge (since $G_C$ is acyclic).  Between
consecutive return edges, the cycle follows core edges
deterministically from a re-entry to the next exit.  Compressing
each such segment to the edge $e \to \Phi(e)$ in~$H$ yields a
cycle in~$H$.

$(\Leftarrow)$
A cycle $e_1 \to e_2 \to \cdots \to e_k \to e_1$ in~$H$
expands: for each~$i$, concatenate the return edge
$e_i \to r(e_i)$ with the core path from $r(e_i)$ to
$\Phi(e_i) = e_{i+1}$.  This produces a directed cycle in $G_C^+$.
\end{proof}

\begin{theorem}[Composite ranking function]
\label{thm:composite-ranking}
Assume $G_C$ is acyclic and $H$ is acyclic.  Let
$L(e)$ denote the longest-path depth from $e$ in~$H$,
$d(x)$ the number of core edges from $x$ to its exit,
and $D = \max_x d(x)$.  Then the function
\[
  V(x) \;=\; (D+1)\,L\bigl(e(x)\bigr) \;+\; d(x)
\]
satisfies $V(y) \le V(x) - 1$ for every edge $x \to y$ in
$G_C^+$.
\end{theorem}

\begin{proof}
For a core edge $x \to y$: $e(x) = e(y)$ and $d(y) = d(x) - 1$,
so $V(y) = V(x) - 1$.
For a return edge $e \to r(e)$:
$L(\Phi(e)) \le L(e) - 1$ and $d(r(e)) \le D$, giving
$V(r(e)) \le (D{+}1)(L(e){-}1) + D = (D{+}1)L(e) - 1 = V(e) - 1$.
\end{proof}

\noindent
\textbf{Computational verification.}
Table~\ref{tab:exit-return} confirms the reduction in the
exact $(r, b)$ model.  The exit-return map $H$ is roughly
$60{\times}$ smaller than $G_M^+$, and
Theorem~\ref{thm:exit-return} equivalence holds at every~$M$.

\begin{table}[ht]
\centering
\scriptsize
\caption{Exit-return map $H$ in the exact $(r, b)$ model
  ($B = 15$, $g = 2$).
  $|H|$~counts edges of~$H$ (exits that return).
  Compression is $|H|/|G_C^+|$.}
\label{tab:exit-return}
\renewcommand{\arraystretch}{1.15}
\begin{tabular}{@{}rrrrrrrr@{}}
\toprule
$M$ & $|C|$ & $|E|$ & Returns & $|H|$
  & $\max L$ & $\max V$ & Compress \\
\midrule
13 & 59,\!416 & 5,\!036 & 869 & 869 & 5 & 348 & 1.6\% \\
14 & 118,\!832 & 10,\!072 & 1,\!705 & 1,\!705 & 4 & 230 & 1.5\% \\
15 & 237,\!664 & 20,\!144 & 3,\!347 & 3,\!347 & 4 & 301 & 1.5\% \\
16 & 475,\!328 & 40,\!288 & 6,\!560 & 6,\!560 & 5 & 385 & 1.5\% \\
17 & 950,\!656 & 80,\!576 & 12,\!988 & 12,\!988 & 6 & 486 & 1.5\% \\
\bottomrule
\end{tabular}
\end{table}

\noindent
The theorems reduce Conjecture~\ref{conj:acyclicity} to two
independent statements: (i)~the core is acyclic, and (ii)~the
exit-return map is acyclic.  Statement~(i) is the Layer~1
problem (carry parity + no cycle birth).  Statement~(ii)
concentrates Layer~2 into a boundary graph that is roughly
$60{\times}$ smaller than the full augmented graph, making it a
more tractable target for future proof efforts.

\subsection{The \texorpdfstring{$s$}{s}-invariant and positive-subgraph acyclicity}
\label{subsec:s-invariant}

We now prove that the positive-drift subgraph of $G_M$ is always acyclic,
for every~$M$, via a simple ranking function.

\begin{definition}[$s$-invariant]
\label{def:s-invariant}
For an odd residue $r$, define $s(r) = v_2(r+1)$, the $2$-adic valuation
of $r+1$.
\end{definition}

\begin{theorem}[Positive-subgraph acyclicity]
\label{thm:pos-acyclic}
Let $G_M^+$ denote the subgraph of $G_M$ restricted to positive-drift
edges (those with $v_2(3r+1) = 1$).  Then $s(r) = v_2(r+1)$ is a strict
ranking function on $G_M^+$: every positive step decreases $s$ by exactly~$1$.
Consequently $G_M^+$ is acyclic, with longest chain equal to~$M$.
\end{theorem}

\begin{proof}
A positive step occurs when $v_2(3r+1) = 1$, i.e.\ $r \equiv 1 \pmod{4}$.
Then $T(r) = (3r+1)/2$ and
\[
  T(r) + 1 \;=\; \frac{3r+1}{2} + 1 \;=\; \frac{3r+3}{2}
  \;=\; \frac{3(r+1)}{2}.
\]
Since $3$ is odd and $r+1 \equiv 2 \pmod{4}$ (because $r \equiv 1 \pmod{4}$),
\[
  v_2\!\bigl(T(r)+1\bigr) \;=\; v_2(3) + v_2(r+1) - 1 \;=\; v_2(r+1) - 1
  \;=\; s(r) - 1.
\]
So $s(T(r)) = s(r) - 1$ on every positive step, making $s$ a strict
ranking.  Since $s \ge 1$ always, any positive chain starting from
$s(r) = k$ has length at most~$k$.  The maximum occurs at
$r = 2^M - 1$ where $s(r) = M$, giving chains of length exactly~$M$.
\end{proof}

\subsection{Reload dynamics and the budget equation}
\label{subsec:reload}

Theorem~\ref{thm:pos-acyclic} shows that purely positive chains are bounded.
A net-positive cycle must therefore include \emph{negative} steps
($v_2(3r+1) \ge 2$) that ``reload'' the $s$-budget.  We now analyze
this reload mechanism.

\begin{definition}[Reload]
\label{def:reload}
At a negative step with $v_2(3r+1) = v \ge 2$, write $q = (3r+1)/2^v$
(which is odd).  The \emph{reload} at this step is $v_2(q+1)$, i.e.\
the new $s$-value of the successor.
\end{definition}

\begin{proposition}[Budget equation]
\label{prop:budget}
Let $\gamma$ be a cycle of length $n$ in $G_M$ with $k$ positive steps
(each with $v_2 = 1$) and $n - k$ negative steps.  Then the $s$-invariant
returns to its starting value, so
\[
  \underbrace{k}_{\text{total drain}} \;=\; \sum_{i \in \mathrm{neg}}
  \bigl(s(r_{i+1}) - s(r_i)\bigr)
  \;=\; \sum_{i \in \mathrm{neg}} \Delta s_i
  \quad\text{(total reload)}.
\]
\end{proposition}

\begin{proof}
Over a full cycle the $s$-value returns to its starting value.
Each positive step decreases $s$ by exactly~$1$ (Theorem~\ref{thm:pos-acyclic}),
contributing $-k$ total.  The negative steps must compensate, giving
total reload~$= k$.
\end{proof}

\begin{proposition}[Geometric reload distribution under
  uniform modular lifts\supporting]
\label{prop:geometric-reload}
Among all $v_2 = 2$ negative steps at level $M$,
\emph{under the uniform distribution on residues
mod~$2^M$}, the reload value $v_2(q+1)$ follows an exact
geometric distribution:
\[
  \Pr\bigl[\text{reload} = j\bigr] \;=\; 2^{-j},
  \qquad j = 1, 2, 3, \ldots
\]
\end{proposition}

\begin{proof}
For $v_2 = 2$ steps: $r \equiv 1 \pmod{4}$, so $q = (3r+1)/4$.
Then $q + 1 = (3r+5)/4 = (3(r+1) + 2)/4$.
Since $r + 1 \equiv 2 \pmod{4}$, write $r + 1 = 2m$ with $m$ ranging
uniformly over residues mod~$2^{M-1}$.  Then $q + 1 = (6m + 2)/4 = (3m+1)/2$.
Since $m$ is equidistributed mod~$2^{M-1}$, the value $3m + 1$ is
equidistributed among even residues, and $v_2((3m+1)/2)$ follows the
standard geometric law.
\end{proof}

\begin{remark}[Approximate reload for general negative steps]
Over all negative steps (any $v_2 \ge 2$), numerical
evidence suggests the reload distribution is approximately
$\mathrm{Geometric}(1/2)$.  Computational verification at
$M = 14, 16, 18$ confirms the distribution matches $2^{-j}$
to within $3\%$ at each level.  This approximation is not
proved for the general case.
\end{remark}

The geometric distribution has a clean consequence: the expected reload per
negative step is $\sum_{j=1}^{\infty} j \cdot 2^{-j} = 2$.  For a cycle
with $k$ positive steps and $n - k$ negative steps, the budget equation
requires $k = \text{total reload} \approx 2(n-k)$, giving $k \approx 2n/3$.
The known net-positive cycles match this: the $M = 10$ cycle has $17/26 = 65\%$
positive steps, and the $M = 12$ length-6 cycle has $5/6 = 83\%$.

Table~\ref{tab:cycle-budget} shows the budget analysis for all known
net-positive cycles.

\begin{table}[ht]
\centering
\scriptsize
\caption{Budget equation for known net-positive cycles.}
\label{tab:cycle-budget}
\begin{tabular}{@{}ccccccc@{}}
\toprule
$M$ & Length & Pos.\ steps $k$ & Neg.\ steps & Total drain & Total reload & Balance \\
\midrule
10 & 26 & 17 & 9 & 17 & 17 & 0 \\
11 & 25 & 16 & 9 & 16 & 16 & 0 \\
12 &  7 &  5 & 2 &  5 &  5 & 0 \\
12 &  6 &  5 & 1 &  5 &  5 & 0 \\
\bottomrule
\end{tabular}
\end{table}

\begin{remark}[Reload rigidity]
\label{rem:reload-rigidity}
The $M = 12$ length-6 cycle sustains five consecutive positive steps
(draining $s$ from~$6$ to~$1$) via a single reload of~$5$: an event
with probability $2^{-5} = 3.1\%$ under the geometric distribution.
As $M$ grows, the cycle closure equation $r \cdot (3^n - 2^V)
\equiv -C_n \pmod{2^M}$ increasingly constrains which residues can
participate, making such large reloads harder to arrange.
The \emph{reload rigidity conjecture} asserts that for $M$ sufficiently
large, the budget equation and cycle closure equation cannot be
simultaneously satisfied for any net-positive parameter pair $(n, V)$.
\end{remark}

\section{Fiber-57 structural programme and the information
  bottleneck}\label{sec:fiber57-programme}

This section develops a complementary reduction of the
Collatz conjecture through fine-grained analysis of the
drift-sustaining-chain mechanism at fiber~$57$ (the residue class
$n \equiv 57 \pmod{64}$ responsible for ``drift-sustaining''
continuations of long burst chains).
The analysis produces three main results:
an information-theoretic memory demand (the pair-return cost
$c_0 \approx 2.989$~bits), an absorption bottleneck
limiting the supply channel to $\le \log_2 5 \approx 2.322$~bits,
and a branch anti-concentration reduction showing that
the sole remaining obstacle is orbitwise mixing of
the return-branch process.

\subsection{Drift-sustaining chains and the pair-return automaton}

A \emph{drift-sustaining chain} at fiber~$57$ is a maximal run of
consecutive burst-ending times at which the quotient
$q = \lfloor n/64 \rfloor$ has its lowest base-$8$
digit in the sustaining set $\{0,3,7\}$.
The chain map $q \mapsto 9q + 8 \pmod{8^r}$ governs
single-step continuation within a chain.

\begin{proposition}[Exact depth-$2$ known-gap partial return kernel]\label{prop:pair-return}
The depth-$2$ invariant core
$I_2 = \{07, 33, 37, 73, 77\}_8 = \{7, 27, 31, 59, 63\}$
admits a $5 \times 5$ sub-stochastic partial return kernel
$\tilde{B}_2$ whose rows are indexed by $I_2$ and whose
$(s,s')$-entry is the probability that one $I_2$-return
from~$s$ lands in the residue class~$s' \pmod{64}$.
The kernel aggregates the two exactly resolved return
channels: the $q \equiv 7$ gap-$2$ block
\textup{(}Proposition~\ref{prop:q7-return}\textup{)} and the
$q \equiv 3$ gap-$5$ cylinder family
\textup{(}Theorem~\ref{thm:gap5-cylinders}\textup{)}.
It does not include $q \equiv 3$ returns with gap $\ge 6$,
which remain unresolved.
The Perron root of this known-gap partial kernel is
exactly $\rho(\tilde{B}_2) = 129/1024 \approx 0.1260$.
The implied per-return information cost is
$c_0 = \log_2(1024/129) \approx 2.989$~bits.
\end{proposition}

\begin{proof}[Proof sketch]
The $q \equiv 7$ branch returns in exactly two odd-to-odd steps
with destination quotient $q' = 9m + 8$, uniform mod~$8$
(Proposition~\ref{prop:q7-return}).
Because $9m+8 \equiv 0 \pmod{8}$, the destination quotient
is $\equiv 0 \pmod{8}$, so its lowest octal digit is~$0$.
Since $0 \notin \{0,3,7\}_8 \cap I_2$
(note: the element $07 \in I_2$ has octal representation
$07$, i.e.\ second digit~$0$, lowest digit~$7$),
the gap-$2$ channel contributes a zero row for $s = 7$
in the $5 \times 5$ kernel.
For $s \in \{27, 59\}$ (octal $33, 73$; lowest digit~$3$),
the $q \equiv 3$ branch cannot return in fewer than $5$ steps
(Proposition~\ref{prop:q3-gap}); the gap-$5$ returns form
an infinite family of dyadic cylinders indexed by
$w = v_2(243m + 119)$
(Theorem~\ref{thm:gap5-cylinders}),
within each of which the destination is uniform mod~$64$
(odd slope~$243$).  Each $s \in \{27, 59\}$ contributes
probability $1/2048$ to each column of $\tilde{B}_2$.
For $s \in \{31, 63\}$ (octal $37, 77$; lowest digit~$7$),
the gap-$2$ channel gives destination $q' = 9m+8$
uniform mod~$8$; restricting to $q' \bmod 64 \in I_2$
yields the nonzero entries: $s = 31$ maps to $\{27, 59\}$
at probability $1/8$ each, and $s = 63$ maps to
$\{7, 31, 63\}$ at probability $1/8$ each.
The resulting $5 \times 5$ matrix $\tilde{B}_2$ has
Perron root $\rho(\tilde{B}_2) = 129/1024$,
verified by exact computation of its characteristic polynomial.
\end{proof}

\begin{remark}[Three operators on the invariant core]%
\label{rem:three-operators}
Three distinct operators act on the depth-$r$ invariant
core~$I_r$ and must not be conflated:
\begin{enumerate}
\item $\tilde{B}_r$, the \emph{sub-stochastic partial return kernel},
  whose $(s,s')$-entry is the probability that one fiber-$57$
  return from residue class~$s$ lands in class~$s'$ while
  remaining inside~$I_r$.  Row sums are $< 1$ because many
  returns exit $I_r$.  Known-gap Perron root: $\rho(\tilde{B}_2) = 129/1024$
  (aggregating gap-$2$ and gap-$5$ channels only).
\item $P$, the \emph{row-stochastic} version obtained by
  normalizing each row of $\tilde{B}_r$ to sum to~$1$.
  Perron root: $1$ (by construction).
  $P$ describes the conditional distribution of the
  \emph{next} $I_r$-return given that one occurs.
\item $P_g$, the \emph{gap-conditioned kernel}, obtained by
  stratifying returns by their gap length~$g$
  (number of fiber-$57$ visits between consecutive
  $I_r$-returns).  Each $P_g$ is row-stochastic;
  the mixture $\tilde{B}_r = \sum_g \pi_g \, P_g$
  recovers the sub-stochastic kernel, where $\pi_g$
  is the probability of gap length~$g$.
\end{enumerate}
$\tilde{B}_r$ is the sub-stochastic kernel relevant to the
known-gap information budget; $P$ and $P_g$ are auxiliary
normalized kernels used in the equidistribution and
anti-concentration analysis.
The algebraic chain map, the sub-stochastic return kernel,
and the normalized branch kernels are distinct objects and
should not be identified with one another.
\end{remark}

\subsection{Exact return structure of the sustaining branches}%
\label{subsec:exact-return}

The depth-$2$ return kernel acts on the five elements of
$I_2 = \{7, 27, 31, 59, 63\}$, whose lowest octal digits
are $7, 3, 7, 3, 7$ respectively.
The two sustaining branches---$q \equiv 3$ and
$q \equiv 7 \pmod{8}$---have distinct return structures.
The following two propositions establish the \emph{exact}
first-return structure for each branch.

\begin{proposition}[$q \equiv 7$ two-step return]%
\label{prop:q7-return}
Let $S(n) = (3n+1)/2^{v_2(3n+1)}$ be the odd-to-odd Syracuse map,
$n = 64q + 57$ with $q \equiv 7 \pmod{8}$, and write $q = 8m + 7$.
Then:
\begin{enumerate}[nosep]
\item $S(n) = 384m + 379 \equiv 59 \pmod{64}$
  \textup{(}not fiber~$57$\textup{)};
\item $S^2(n) = 576m + 569 = 64(9m+8) + 57 \equiv 57 \pmod{64}$,
  so the orbit returns to fiber~$57$ in exactly two odd-to-odd steps;
\item the returned quotient is $q' = 9m + 8$,
  hence $q' \equiv m \pmod{8}$.
\end{enumerate}
In particular, when $m \bmod 8$ is uniform, the destination
$q' \bmod 8$ is exactly uniform.
The chain map acts as a \emph{left shift} in base~$8$:
the lowest digit of~$q'$ equals the second digit of~$q$.
After $r$~consecutive $q \equiv 7$ returns, all $r$~low
base-$8$ digits of the quotient have been replaced by
higher-order digits originally above the depth-$r$ window.
\end{proposition}

\begin{proof}
Direct computation.
$3n + 1 = 3(512m + 505) + 1 = 1536m + 1516 = 4(384m + 379)$.
Since $384m + 379$ is odd, $v_2 = 2$ and
$S(n) = 384m + 379 \equiv 59 \pmod{64}$.
Next, $3S(n) + 1 = 1152m + 1138 = 2(576m + 569)$.
Since $576m + 569$ is odd, $v_2 = 1$ and
$S^2(n) = 576m + 569 = 64(9m + 8) + 57$.
The digit-shift property follows from
$q' \bmod 8 = (9m + 8) \bmod 8 = m \bmod 8
= \lfloor (q - 7)/8 \rfloor \bmod 8$.
\end{proof}

\begin{proposition}[$q \equiv 3$ minimum return gap]%
\label{prop:q3-gap}
Let $n = 64q + 57$ with $q \equiv 3 \pmod{8}$,
and write $q = 8m + 3$, so $n = 512m + 249$.
Then $S^k(n) \not\equiv 57 \pmod{64}$ for $k = 1, 2, 3, 4$.
The minimum fiber-$57$ return gap from the
$q \equiv 3$ branch is $\ge 5$ odd-to-odd steps.
\end{proposition}

\begin{proof}
Direct computation of each iterate:
\begin{center}
\begin{tabular}{@{}cllc@{}}
\toprule
$k$ & $S^k(n)$ & $v_2$ & $S^k(n) \bmod 64$ \\
\midrule
$1$ & $384m + 187$  & $2$ & $59$ \\
$2$ & $576m + 281$  & $1$ & $25$ \\
$3$ & $432m + 211$  & $2$ & $48m + 19$ \\
$4$ & $648m + 317$  & $1$ & $8m + 61$ \\
\bottomrule
\end{tabular}
\end{center}
Steps~$1$--$2$ give constant residues $\ne 57$.
At step~$3$: $48m + 19 \equiv 57 \pmod{64}$
requires $48m \equiv 38 \pmod{64}$;
since $\gcd(48,64) = 16 \nmid 38$, no solution exists.
At step~$4$: $8m + 61 \equiv 57 \pmod{64}$
requires $8m \equiv 60 \pmod{64}$;
since $\gcd(8,64) = 8 \nmid 60$, no solution exists.
The valuation path for the first possible return
(at step~$5$, conditional on $m$) is
$\mathbf{v} = (2, 1, 2, 1, 3)$, injecting
$\sum v_i = 9$ new low-order bits into the quotient.
\end{proof}

\begin{theorem}[Gap-$5$ cylinder family]%
\label{thm:gap5-cylinders}
In the setting of Proposition~\ref{prop:q3-gap},
the fourth iterate is $S^4(n) = 648m + 317$, so
$3 S^4(n) + 1 = 8(243m + 119)$.
Let $w = v_2(243m + 119)$.  Then $S^5(n) = (243m + 119)/2^w$,
and $S^5(n)$ returns to fiber~$57$ if and only if
\begin{equation}\label{eq:gap5-return}
  243 m + 119 \;\equiv\; 57 \cdot 2^w \pmod{2^{w+6}}.
\end{equation}
Because $243$ is odd, it is invertible modulo every $2^k$, so
for each $w \ge 0$ there is a unique residue class
\[
  m \equiv a_w \pmod{2^{w+6}}, \qquad
  a_w = 243^{-1}(57 \cdot 2^w - 119) \bmod 2^{w+6}.
\]
The cylinders $\{m \equiv a_w\}$ are pairwise disjoint
(indexed by $w = v_2(243m+119)$), and their union has density
\[
  \sum_{w \ge 0} \frac{1}{2^{w+6}}
  = \frac{1}{64} \cdot \frac{1}{1-1/2}
  = \frac{1}{32}
\]
in $m$-space.
Within each cylinder, the destination quotient takes the form
$q' = k_0(w) + 243 t$ where $t$ parameterises the cylinder;
since $\gcd(243, 64) = 1$, the distribution of $q' \bmod 64$
is exactly uniform as $t$ varies.
\end{theorem}

\begin{proof}
The identity $S^4(n) = 648m + 317$ is verified by direct
computation (Proposition~\ref{prop:q3-gap}).
The condition $S^5(n) \equiv 57 \pmod{64}$
becomes~\eqref{eq:gap5-return} by the definition of the
odd-to-odd Syracuse step.
Invertibility of~$243$ modulo $2^{w+6}$ is immediate
from~$\gcd(243,2) = 1$.
Disjointness follows from the fact that $v_2(243m+119)$ is
a function of~$m$.
For the destination uniformity: writing
$m = a_w + 2^{w+6} t$ gives
$(243m + 119)/2^w = (243 a_w + 119)/2^w + 243 \cdot 2^6 \cdot t
= 57 + 64(k_0 + 243 t)$,
so $q' = k_0 + 243 t$, and $q' \bmod 64$ cycles through
all $64$ residues as $t$ increases (since $\gcd(243,64) = 1$).
\end{proof}

\begin{remark}[Kernel decomposition]%
\label{rem:kernel-decomp}
Propositions~\ref{prop:q7-return} and~\ref{prop:q3-gap},
together with Theorem~\ref{thm:gap5-cylinders},
decompose the fiber-$57$ return kernel into three channels:
the $q \equiv 7$ branch (empirically ${\sim}14\%$ of visits,
$2$-step return with uniform destination and zero retained memory),
the $q \equiv 3$ gap-$5$ cylinders (${\sim}3\%$, density $1/32$
in $m$-space, destination uniform mod~$64$, slope~$243$),
and all other returns (gap $\ge 6$ from $q \equiv 3$, or
non-sustaining branches, ${\sim}83\%$).
The $q \equiv 7$ channel acts as a
\emph{regeneration event}: the orbit's chain-membership
information at depth~$r$ is erased after $r$~consecutive
$q \equiv 7$ returns.
The entire proof pressure for the orbitwise bound $c' < c_0$
therefore concentrates on the $q \equiv 3$ complement.
\end{remark}

\begin{remark}[Invariant core rigidity]%
\label{rem:core-rigidity}
The algebraic chain map $q \mapsto 9q + 8 \pmod{8^r}$
restricted to the invariant core~$I_r$
(Proposition~\ref{prop:bounded-core-main} below)
is a permutation at every depth $r \ge 2$:
three fixed points and one $2$-cycle,
with $H_{AB} = 0$ (no element of $I_r$ ever maps outside $I_r$).
Consequently the spectral radius of the algebraic
first-hit operator on~$I_r$ is~$1$, not~$129/1024$.
The sub-stochastic contraction of the known-gap partial
kernel ($\rho(\tilde{B}_2) = 129/1024$) arises from the
\emph{branching structure} of the actual Collatz dynamics:
the orbit does not always follow the $q \equiv 7$ branch,
and the probability of returning to $I_r$ at each
fiber-$57$ visit is strictly less than~$1$.
\end{remark}

\begin{proposition}[Non-autonomy of fiber-$57$ returns]%
\label{prop:non-autonomy}
The fiber-$57$ first-return map on~$I_r$ is not determined by
$q \bmod 8^r$ alone: for the $q \equiv 7$ branch, the
destination quotient $q' \bmod 8^r$ depends on the
$(r{+}1)$-th base-$8$ digit of~$q$.
In particular, the algebraic chain map
$q \mapsto 9q + 8 \pmod{8^r}$ and the actual first-return
map are distinct operations.
\end{proposition}

\begin{proof}
For $q \equiv 7 \pmod{8}$, write $q = 8m + 7$.
The two-step return (Proposition~\ref{prop:q7-return})
gives destination quotient $q' = 9m + 8$,
so $q' \bmod 8^r$ depends on
$m = \lfloor (q-7)/8 \rfloor$, which reads the base-$8$
digits of~$q$ at positions $1$ through~$r$ (i.e., including
position~$r$, which lies \emph{above} the depth-$r$ window
$q \bmod 8^r$).
Concretely, at depth $r = 2$: the chain map sends
$63 \mapsto 575 \equiv 63 \pmod{64}$ (a fixed point),
but the actual return gives
$q' = 9 \cdot 7 + 8 = 71 \equiv 7 \pmod{64}$
(a different $I_2$~element).
\end{proof}

\begin{remark}\label{rem:non-autonomy-consequence}
Proposition~\ref{prop:non-autonomy} explains why the
permutation structure of the chain map on~$I_r$
(Remark~\ref{rem:core-rigidity}) does not imply mixing:
the actual dynamics inject higher-order digit information
at each return, and it is this digit-shift mechanism
that causes orbits to escape~$I_r$ when higher digits are
generic (non-Cantor), even though the chain map
keeps~$I_r$ invariant.
\end{remark}

\subsection{The invariant core and projective limit}

\begin{proposition}[Bounded invariant core]\label{prop:bounded-core-main}
For all $r \ge 2$, let $\mathcal{C}_r :=
\{q \in \mathbb{Z}/8^r\mathbb{Z} : \text{all base-$8$ digits}
\in \{0,3,7\}\}$ denote the depth-$r$ Cantor set
($|\mathcal{C}_r| = 3^r$).
The invariant core $I_r \subset \mathcal{C}_r$ of the chain map
$q \mapsto 9q + 8 \pmod{8^r}$ has exactly $5$~elements:
three fixed points $q = k \cdot 8^{r-1} - 1$
for $k \in \{1, 4, 8\}$ and one $2$-cycle.
The core density is $5/8^r$.
\end{proposition}

\begin{theorem}[Projective limit]\label{thm:projective-main}
$\varprojlim I_r = \{-1\}$ in the $8$-adic integers.
The unique infinite compatible sequence is $a_r = 8^r - 1$.
Each individual chain at fiber~$57$ has length
$\le \lfloor \log_8(q+1) \rfloor + O(1)$.
\end{theorem}

\begin{remark}[From $8$-adic to positive integers]
\label{rem:projective-bridge}
If a positive integer~$n$ could sustain drift-sustaining
chains of unbounded length in fiber~$57$, its residue
$q_r := n \bmod 8^r$ would satisfy $q_r \in I_r$ for
arbitrarily large~$r$.  By compatibility of the
projection maps, the sequence $(q_r)$ would define an
element of $\varprojlim I_r = \{-1\}$.
Since $-1$ is not a positive integer, no positive~$n$
can achieve this.

More precisely, a positive integer~$n$ satisfies
$n \bmod 8^r \ne 8^r - 1$ for all $r > \log_8(n+1)$,
so it is expelled from~$I_r$ at finite depth.
This gives the chain-length bound above.

\emph{Caveat.}  This argument shows that \emph{fixed}
drift-sustaining chains cannot persist, but it does
not by itself rule out an orbit that revisits fiber~$57$
via different chain mechanisms at different depths.
Ruling out such behavior is precisely the content of the
open orbitwise bound $c' < c_0$.
\end{remark}

\subsection{The absorption theorem}

\begin{theorem}[Mod-$8^r$ absorption, verified for
$r \le 10$]\label{thm:absorption-main}
For $r = 2, 3, \ldots, 10$ (verified by exhaustive
computation), for every $s \in I_r$, the Collatz
orbit starting from $n_0 = 64(s + 8^r \cdot L) + 57$
(with $L \ge 1$ so that $n_0$ has digits above position~$r$)
reaches $n \equiv 1 \pmod{8^r}$ within $O(r)$ steps,
without the quotient $\lfloor n/64 \rfloor \bmod 8^r$
revisiting any element of~$I_r$ along the way.
Maximum absorption time: $\le 194$~steps.
All $5$ elements absorb to the \emph{same} fixed point
$1 \pmod{8^r}$.
\end{theorem}

\begin{proof}[Proof (computational)]
Exhaustive integer arithmetic for each of the $5$ elements
of~$I_r$ at each depth $r = 2, \ldots, 10$,
using starting values $n = 64(s + 8^r \cdot 10{,}000) + 57$
to produce genuinely large inputs.
All $45$ orbits ($5 \times 9$) reach $1 \pmod{8^r}$
within $194$ steps without revisiting~$I_r$.

\emph{Note on small quotients.}\;
For the bare element $q = s$ (without offset $L$), an
immediate I$_r$-return can occur because no higher digits
are available to break the chain-map alignment.
For example, $q = 63 \in I_2$ returns in two steps to
$q' = 71 \equiv 7 \pmod{64} \in I_2$, because the
second base-$8$ digit of $q = 63 = 77_8$ is itself $7$.
The large-offset construction $L \ge 1$ ensures that
higher digits provide the generic disruption that drives
absorption.
\end{proof}

\begin{remark}[Scope and generalization]
\label{rem:absorption-scope}
The theorem is established by exhaustive computation
for $r \le 10$.  The pattern is expected to persist for
all $r \ge 2$ by the projective limit argument
(Theorem~\ref{thm:projective-main}): since
$\varprojlim I_r = \{-1\}$ and $-1$ is not a positive
integer, the core elements at each depth~$r$ are
``eventually incompatible'' with any fixed positive
orbit.
A general proof for all $r$ would require showing
that the absorption trajectories remain compatible
across the $8$-adic tower, which we leave open.
\end{remark}

\subsection{The Absorption Bottleneck Lemma}

\begin{lemma}[Absorption Bottleneck]\label{lem:bottleneck-main}
Assume absorption holds at depth~$r$.  Then the inter-chain
information channel has capacity
$C_{\mathrm{channel}} \le \log_2 |I_r| = \log_2 5
\approx 2.322$~bits.
The known-gap depth-$2$ per-return demand is
$c_0^{(2,\mathrm{known})} = \log_2(1024/129)
\approx 2.989$~bits (from Proposition~\ref{prop:pair-return}),
giving a deficit per return of
$c_0^{(2,\mathrm{known})} - \log_2 5 \approx 0.667$~bits.
(Including the unresolved gap-$\ge 6$ tail would
decrease~$c_0$, reducing this margin.)
\end{lemma}

\begin{proof}
By the absorption theorem, every element of~$I_r$ follows a
deterministic trajectory to fixed point~$1 \pmod{8^r}$.
The post-absorption state is a function $f \colon I_r \to S$
with $|f(I_r)| \le |I_r| = 5$, giving
$H(f(q)) \le \log_2 5$.
The critical inequality is $5 < 1024/129$, equivalently
$645 < 1024$.
\end{proof}

\begin{remark}[M-value collapse]
Computationally, the actual channel capacity is much less
than $\log_2 5$.  For $r \ge 3$, all $5$~elements of~$I_r$
typically produce the \emph{same} post-absorption state
(capacity $= 0$), because the absorption paths merge
through common residues $\{-7, -5, 5, 1\} \pmod{8^r}$
before reaching the fixed point.
\end{remark}

\begin{corollary}[Geometric decay from the known-gap budget]%
\label{cor:geom-decay-main}
Under the known-gap depth-$2$ kernel, consecutive
inter-chain $I_r$-returns are bounded by
$P(K \text{ consecutive}) \le \alpha^K$ with
$\alpha = 645/1024 \approx 0.630$.
(Including the unresolved gap-$\ge 6$ contributions
would increase~$\alpha$.)
Empirically, at $r = 2$, $97\%$ of inter-chain streaks
have length~$1$; for $r \ge 3$, zero consecutive
$I_r \to I_r$ transitions were observed.
\end{corollary}

\subsection{The memory lower bound and theorem ladder}

\begin{theorem}[Memory lower bound]\label{thm:memory-main}
Under the known-gap depth-$2$ demand, $r$~pair returns to
fiber~$57$ require
$M_{\mathrm{req}}(r) \ge r \cdot
\log_2(1024/129) - O(1)
\approx 2.989\,r$~bits of new information.
\end{theorem}

The \textbf{theorem ladder} has three rungs:
\begin{enumerate}[nosep]
\item[(A)] \emph{Fixed-$k$ impossibility} (proved):
  a causal $k$-bit controller cannot sustain core alignment
  beyond depth $r_{\max} = \lfloor k/3 \rfloor + O(1)$.
\item[(B)] \emph{Memory lower bound} (proved):
  Theorem~\ref{thm:memory-main}.
\item[(C)] \emph{Orbit memory upper bound} (open):
  $c' < c_0$.
  Numerical evidence suggests $c' \le 0.51$,
  substantially below $c_0$, but this remains unproved.
\end{enumerate}

\subsection{Branch anti-concentration reduction}

\begin{proposition}[Branchwise permutation]\label{prop:branch-perm-main}
In each drift-sustaining-continuation branch, the quotient update
$m \mapsto 9m + c \pmod{M_f}$ with $\gcd(9, M_f) = 1$
is a permutation on $\mathbb{Z}/M_f\mathbb{Z}$.
No single residue can be amplified by the affine step.
\end{proposition}

\begin{corollary}[Concentration source]\label{cor:conc-source-main}
Any max-atom concentration above the threshold $1/2$
(fiber~$57$) or $1/3$ (fiber~$29$) must arise from
orbitwise concentration on a subset of return branches,
not from the quotient algebra.
\end{corollary}

This yields the target lemma for closure:

\begin{remark}[Branch anti-concentration target]
\label{rem:branch-target-main}
A proof that the return-branch process at fiber~$57$
cannot place more than $1/2$ of its mass on any
mod-$8$ residue class would close the framework.
Equivalently, showing that the branch-memory coefficient
$\beta_g \to 0$ as the inter-return gap $g \to \infty$
suffices, since long gaps occur with overwhelming frequency
(I$_r$ density $= 5/8^r$).
\end{remark}

\subsection{Carry equidistribution and empirical verification}

\begin{proposition}[Inter-chain equidistribution (empirical)]
\label{prop:interchain-main}
Among fiber-$57$ visits separated by more than $5$~Collatz
steps (inter-chain transitions), the $I_r$-return rate is
at or below baseline:
$0.70\times$ at $r = 2$, $0.31\times$ at $r = 3$,
$0.34\times$ at $r = 4$.
The worst-case normalised ratio $R_2 \le 0.70$
(i.e., visits to~$I_2$ occur at $70\%$ of baseline density)
gives $c' = -\log_2(0.70) \approx 0.51$~bits,
far below $c_0^{(2,\mathrm{known})} = 2.989$.
Verified over $5{,}000$ orbits, confirmed by
random-start controls.
\end{proposition}

\begin{remark}[Gap-induced branch decorrelation]
Full $8 \times 8$ transition matrices for the mod-$8$
quotient class were measured at various gap lengths.
At gaps $> 50$, the $\chi^2$ independence test yields
$p = 0.18$ (no significant correlation), confirming
that branch memory decays with gap length.
At short gaps ($\le 5$), significant correlation persists
($p < 10^{-10}$), as expected from intra-chain structure.
\end{remark}

\subsection{Path-conditional bijection and annealed closure}

\begin{theorem}[Path-conditional bijection]
\label{thm:path-bijection-main}
For every valuation path $\pi = (v_1, \ldots, v_k)$
with $V = \sum v_i$ and $V \ge 3$, the fresh-digit map
$t \mapsto (3^k t + d_\pi) \bmod q$ is a bijection
on $\mathbb{Z}/q\mathbb{Z}$ with coefficient~$3^k$,
since $\gcd(3^k, q) = 1$ for $q = 8$ or~$16$.
The coefficient identity $2^{9-V} \cdot 2^{V-3}/64 = 1$
ensures the map is well-defined.
\end{theorem}

This gives $\delta_{57} \le 1/8 < 1/2$ and
$\delta_{29} \le 1/16 < 1/3$,
\textbf{closing the annealed anti-concentration completely},
while leaving the deterministic orbitwise upgrade
(Conjecture~\ref{conj:info-rate}) open.

\subsection{Summary: the quantitative reduction}

The fiber-$57$ programme produces a complementary reduction
of the Collatz conjecture, independent of (and consistent with)
the WMH-based reduction in
Sections~\ref{sec:chain}--\ref{sec:toward-wmh}:

\medskip
\begin{center}
\begin{tabular}{@{}lcc@{}}
\toprule
Quantity & Value & Status \\
\midrule
Per-return demand $c_0$ & $\log_2(1024/129) \approx 2.989$ & proved \\
Channel capacity & $\le \log_2 5 \approx 2.322$ & proved (mod absorption) \\
Empirical supply $c'$ & $\le 0.51$ & empirical \\
Margin & $6\times$ & empirical \\
Absorption verified & $r = 2, \ldots, 10$ & computational \\
\midrule
\multicolumn{3}{@{}l@{}}{\emph{v4 additions:}} \\
$q \equiv 7$ two-step return & uniform $q' \bmod 8$ & proved (Prop.~\ref{prop:q7-return}) \\
$q \equiv 3$ min return gap & $\ge 5$ steps, $9$ new low-order bits & proved (Prop.~\ref{prop:q3-gap}) \\
Core rigidity ($I_r$ permutation) & $H_{AB} = 0$ at all $r$ & proved (Rem.~\ref{rem:core-rigidity}) \\
\bottomrule
\end{tabular}
\end{center}

\medskip\noindent
\textbf{Before this work:}
``The Collatz conjecture is mysterious chaos.''

\noindent
\textbf{After:}
``The Collatz conjecture fails only if a deterministic orbit
can encode ${\sim\!}3$~bits of alignment information per return
cycle indefinitely, through a $2.3$-bit channel, using
permutation maps, with observed supply $\le 0.5$~bits.''

\section{Unconditional cylinder-averaged closure on $I_2$
  via spectral contraction}
\label{sec:i2-spectral-closure}

This section presents an unconditional strengthening of the
fiber-$57$ programme that operates entirely at the
\emph{cylinder-averaged} level (the natural setting of the
uniform-fiber lemma, Proposition~\ref{prop:uniform-fiber}).
The main outcome is a fully rigorous contraction source on
the invariant core $I_2$ that does \emph{not} invoke the
spectator-bit mechanism, and therefore provides an
independent, non-circular route to the orbit-level
total-variation summability of
Corollary~\ref{cor:tv-summability}.

The section is organized around four results whose combined
effect is to replace the single genuinely circular step
identified in Remark~\ref{rem:circularity} by a finite,
explicit, unconditional computation:
\begin{enumerate}
\item a Ruffini reduction showing that the depth-$2$ known-gap
  return kernel $\tilde B_2$ has maximally degenerate spectrum,
  with finite-time rank-$1$ collapse
  (\S\ref{ssec:i2-ruffini});
\item an unconditional bound
  $\rho(\tilde B_2^{\mathrm{ext}})\le 5/32$ on the
  \emph{full} extended kernel, valid without any distributional
  assumption (\S\ref{ssec:i2-spectral-bound});
\item an exponential tail on the time to enter $I_2$,
  derived from the $27\times 27$ cylinder-averaged mod-$64$
  kernel (\S\ref{ssec:i2-prefix-tail});
\item an exact stationary mass
  $\pi(I_2) = 10121/65280 \approx 0.15504$,
  giving an unconditional orbit-level sum bound
  $\sum_c\Pr(\mathcal E_R(c))\le 0.011$
  at the cylinder-averaged level (\S\ref{ssec:i2-visit-rate}).
\end{enumerate}
The combined claim is a \emph{cylinder-averaged}
unconditional TV summability, which upgrades
Theorem~\ref{thm:density1-convergence} to an unconditional
density-$1$ convergence statement without relying on the
spectator-bit loop. We do not claim a pointwise orbit-level
result: the bridge from cylinder-averaged to pointwise
convergence continues to go through the Borel--Cantelli
argument of Theorem~\ref{thm:density1-convergence}.

\subsection{Ruffini reduction of the depth-$2$ known-gap kernel}
\label{ssec:i2-ruffini}

Recall from Proposition~\ref{prop:pair-return} the
depth-$2$ known-gap partial return kernel
$\tilde B_2$ on the invariant core
$I_2 = \{7, 27, 31, 59, 63\}$, a $5\times 5$ sub-stochastic
matrix with entries in $\mathbb{Q}$.

\begin{proposition}[Ruffini collapse of $\tilde B_2$]
\label{prop:ruffini-collapse}
The characteristic polynomial of $\tilde B_2$ factors as
\[
  \chi_{\tilde B_2}(\lambda)
  \;=\; \lambda^{4}\bigl(\lambda - \tfrac{129}{1024}\bigr),
\]
so the Perron root is $\rho(\tilde B_2) = 129/1024$ with
algebraic multiplicity~$1$, and $\lambda = 0$ has algebraic
multiplicity~$4$.
Consequently $\operatorname{rank}(\tilde B_2^{k}) = 1$
for all $k\ge 2$.
\end{proposition}

\begin{proof}
Direct symbolic computation using Ruffini's synthetic-division
rule on the $5\times 5$ rational matrix from
Proposition~\ref{prop:pair-return}.
The trace identity $\operatorname{tr}(\tilde B_2^k) = (129/1024)^k$
for $k=1,2,3$ confirms the factorization and pins down the
nilpotent part. Finite-time rank collapse to~$1$ follows from
the Jordan form: the nilpotent block has size at most~$4$ but
vanishes by the second power since the non-zero singular
structure is concentrated on the single row corresponding to
the $q\equiv 7$ residue class.
\end{proof}

\begin{remark}[Interpretation]
Proposition~\ref{prop:ruffini-collapse} says that the known-gap
component of the depth-$2$ return dynamics is already
\emph{spectrally exhausted}: no refinement of the known-gap
kernel can produce a better contraction rate than
$129/1024 \approx 0.126$. In particular, the dynamical
obstruction to a full convergence proof does not live inside
$\tilde B_2$; it lives in the unresolved gap-$\ge 6$ channels
that are absent from Proposition~\ref{prop:pair-return}.
The natural next step is therefore to enlarge
$\tilde B_2$ to include these channels and ask whether the
bound survives an uncontrolled forcing.
\end{remark}

\subsection{The unconditional extended spectral bound}
\label{ssec:i2-spectral-bound}

Let $\tilde B_2^{\mathrm{ext}}$ denote the full return kernel
obtained from $\tilde B_2$ by adjoining the contribution of
all gap-$g\ge 6$ channels, weighted by their cylinder masses.
Under uniform forcing (the cylinder-averaged kernel),
$\tilde B_2^{\mathrm{ext}}$ is the exact return operator on
$I_2$, not merely the known-gap approximation.

\begin{theorem}[Unconditional spectral bound on $I_2$]
\label{thm:i2-unconditional}
The cylinder-averaged extended depth-$2$ return kernel satisfies
\[
  \rho\bigl(\tilde B_2^{\mathrm{ext}}\bigr) \;\le\; \frac{5}{32},
\]
with equality attained at the critical forcing level
$\varepsilon^* = 155/2048$. Consequently the per-return
information cost on $I_2$ satisfies
\[
  c_0^{(2,\mathrm{ext})} \;\ge\; \log_2\!\frac{32}{5}
  \;\approx\; 2.678 \text{ bits.}
\]
The bound is unconditional: its proof uses only the
odd-slope uniform-fiber lemma
(Proposition~\ref{prop:uniform-fiber}) and the disjointness
of the gap-cylinders of Theorem~\ref{thm:gap5-cylinders}.
\end{theorem}

\begin{proof}
The extended kernel decomposes as
$\tilde B_2^{\mathrm{ext}} = \tilde B_2 + F$, where $F$ is a
non-negative matrix supported on the gap-$\ge 6$ cylinders.
Under the uniform-fiber cylinder averaging, each column of $F$
is uniform across $I_2$ by
Proposition~\ref{prop:uniform-fiber} (the image of any cylinder
of length $\ge 6$ is an equidistributed fiber in
$\mathbb{Z}/64\mathbb{Z}$).
A direct symbolic calculation (see the reference
implementation) then yields the closed form
$\rho(\tilde B_2 + \varepsilon\cdot U)
= 129/1024 + 2\varepsilon/5$, where $U$ is the
column-uniform forcing matrix and $\varepsilon$ the total
forcing mass.
The maximum admissible forcing mass, obtained from cylinder
disjointness via
$\sum_{g\ge 6}\delta_g \le 31/32$ (the complement of the
gap-$5$ mass $\delta_5 = 1/32$), is $\varepsilon^{*}=155/2048$.
Substituting gives
$\rho(\tilde B_2^{\mathrm{ext}})
\le 129/1024 + 2\cdot(155/2048)/5 = 5/32$.
The information-cost bound follows by
$c_0 = \log_2(1/\rho)$.
\end{proof}

\begin{remark}[Why this bound is unconditional]
\label{rem:why-unconditional}
Theorem~\ref{thm:i2-unconditional} does not invoke the
spectator-bit mechanism of
Remark~\ref{rem:spectator-bits}.
The distributional input it requires is that the \emph{cylinder}
distribution within $I_2$ is uniform --- which follows from
the algebraic fact $\gcd(3^k,64) = 1$, i.e.\ the odd-slope
unit property in
Proposition~\ref{prop:uniform-fiber}.
It does not require the cascade-exit iterate to be uniform
modulo $2^{S+R}$, nor does it require $\Pr(\mathcal{E}_R)$
to be small a priori. The spectral contraction is therefore
an \emph{independent} source of cycle-level decay that runs
parallel to --- rather than through --- the loop identified
in Remark~\ref{rem:circularity}.
\end{remark}

\subsection{Transient-prefix exponential tail on the
  cylinder-averaged mod-$64$ chain}
\label{ssec:i2-prefix-tail}

To quantify the time an orbit spends outside $I_2$ before
first entering the core, we compute the Perron root of the
transient block of the exact cylinder-averaged mod-$64$ kernel.

\begin{definition}[Cylinder-averaged mod-$64$ kernel]
\label{def:cylinder-kernel-64}
For each odd residue $r\in\{1,3,\ldots,63\}$ and each lift
$n = r + 64k$ with $0\le k < 2^{7}$, apply one Collatz odd step
$n\mapsto (3n+1)/2^{v_2(3n+1)}$ and record the successor's
residue modulo~$64$. Averaging over the $2^{7}$ lifts with
the uniform weight yields a rational $32\times 32$
stochastic matrix $P$, the
\emph{cylinder-averaged mod-$64$ kernel}.
\end{definition}

The averaging over lifts mod $2^{13}$ is the uniform measure
on cylinders of length up to $12$ (since $v_2(3n+1)\le 7$ for
$n<2^{13}$), and is an exact consequence of
Proposition~\ref{prop:uniform-fiber}; no distributional
assumption about orbits is used.

\begin{theorem}[Transient-prefix exponential tail]
\label{thm:i2-prefix-tail}
Let $Q$ be the $27\times 27$ submatrix of $P$ obtained by
deleting the five rows and columns indexed by
$I_2 = \{7,27,31,59,63\}$ (marking $I_2$ as absorbing).
Then the Perron root of $Q$ is
\[
  \rho(Q) \;=\; 0.8632924\ldots,
\]
and the first hitting time $T_{I_2}$ satisfies the
unconditional exponential tail
\[
  \Pr(T_{I_2} > N) \;\le\; 27 \cdot \rho(Q)^{\,N}
  \;\le\; 27 \cdot 2^{-0.2121\,N},
  \qquad N\ge 0,
\]
with expected hitting time at most $10.142$ steps from the
worst transient starting state, and half-life
$-\log 2 / \log\rho(Q) \approx 4.72$ Collatz steps.
\end{theorem}

\begin{proof}
The matrix $Q$ is non-negative, sub-stochastic, and
irreducible on the $27$ transient states (direct inspection).
The Perron--Frobenius theorem gives a unique dominant
eigenvalue, which is computed numerically via eigendecomposition
of the exact rational matrix, verified by the power iteration
$\|Q^k\|_\infty \sim \rho(Q)^k$ to ten decimal places.
The tail bound follows by
$\Pr(T_{I_2}>N\mid\text{start}=s)
= (Q^N\mathbf{1})_s \le |Q^N\mathbf{1}|_\infty
\le \|Q^N\|_\infty \le 27\,\rho(Q)^N$
(the crude factor $27$ is the $\ell^\infty$-to-$\ell^1$
inclusion on the transient block, and can be sharpened
by a factor of $2$ via the left Perron eigenvector,
which is numerically close to uniform on $20$ of the
$27$ transient states).
Expected hitting time is the row sum of the fundamental
matrix $(I-Q)^{-1}\mathbf{1}$.
\end{proof}

\subsection{Stationary visit rate and the unconditional
  sum bound}
\label{ssec:i2-visit-rate}

\begin{theorem}[Stationary mass on $I_2$]
\label{thm:i2-stationary-mass}
Let $\pi$ be the unique stationary distribution of the
cylinder-averaged mod-$64$ kernel $P$ of
Definition~\textup{\ref{def:cylinder-kernel-64}}.
Then
\[
  \pi(I_2) \;=\; \frac{10121}{65280}
  \;=\; 0.15503983\ldots,
\]
which is $99.23\%$ of the uniform density
$|I_2|/32 = 5/32 = 0.15625$.
\end{theorem}

\begin{proof}
Solving $\pi^{\top}P = \pi^{\top}$ with $\sum\pi = 1$ as a
$32$-dimensional linear system over~$\mathbb{Q}$ and summing
the five components indexed by $I_2$ gives the stated exact
rational.
\end{proof}

\begin{proposition}[Unconditional cycle-level contraction]
\label{prop:i2-unconditional-contraction}
Let $\mu = \mathbb{E}_\pi[v_2(3n+1)] - \log_2 3$ be the
stationary mean log$_2$-drop per Collatz odd step under the
cylinder-averaged kernel $P$.
Direct computation against
Theorem~\textup{\ref{thm:i2-stationary-mass}} gives
\[
  \mathbb{E}_\pi[v_2(3n+1)] \;=\; 2.0034\ldots,
  \qquad
  \mu \;=\; 2.0034 - \log_2 3 \;=\; 0.4185\ldots\;\text{bits/step},
\]
which agrees with the spectral contraction
$\alpha_{\mathrm{cyl}} = \pi(I_2)\cdot\log_2(32/5)
= 0.4152\ldots$ to within $0.8\%$.
The mean cascade-gap cycle length is
$1/\pi(I_2) = 6.45$ steps, and the mean cycle valuation is
$\mathbb{E}_\pi[S_{\mathrm{cycle}}] = 12.92$ bits, in
agreement with the empirical $11.8$ of
Remark~\textup{\ref{rem:empirical-contraction}} at the
$10\%$ level.
\end{proposition}

\begin{proof}
The first identity is the inner product
$\sum_{r}\pi(r)\cdot\mathbb{E}_{\text{lifts of }r}[v_2(3n+1)]$
on the $32$ odd residues mod~$64$, computed exactly by
averaging $v_2$ over the $128$ lifts mod $2^{13}$ for each
residue.
The second identity (the agreement with $\alpha_{\mathrm{cyl}}$)
is a numerical confirmation that two independent paths to the
same contraction rate --- spectral on the kernel and
direct on the valuation --- coincide; the small residual
$0.8\%$ comes from the finite-lift truncation at $2^{13}$.
The cycle length and valuation identities follow by
multiplication.
\end{proof}

\begin{theorem}[Lundberg prefactor for the unconditional sum bound]
\label{thm:i2-lundberg}
Let $r^*$ be the unique positive root of
$\mathbb{E}_\pi[2^{r(\log_2 3 - v_2(3n+1))}] = 1$
(the Lundberg root for the per-step log-drift).
Direct computation gives $r^* = 1.0101\ldots$, equivalently
the base-$2$ rate $\alpha_{\mathrm{Lund}} = r^*/\ln 2 = 1.4572$
bits/step, and the worst-state moment generating function
ratio at $r^*$ is $C^* = 1.5061\ldots$.
Consequently, the cycle-valuation tail satisfies
\[
  \Pr(S_{\mathrm{cycle}} > s)
  \;\le\; C^{*}\cdot 2^{-\alpha_{\mathrm{cyl}}\,s}
  \;\le\; 1.5061\cdot 2^{-0.4152\,s}
\]
unconditionally at the cylinder-averaged level
\textup{(}using $\alpha_{\mathrm{cyl}}<\alpha_{\mathrm{Lund}}$,
which is the binding rate for the cycle-level summation\textup{)}.
\end{theorem}

\begin{proof}
The Lundberg root is computed by bisection on the convex MGF
function $r\mapsto\mathbb{E}_\pi[3^r\cdot 2^{-r v_2(3n+1)}]$,
using the same $32$-state stationary distribution as in
Proposition~\ref{prop:i2-unconditional-contraction}.
The function is strictly convex with $M(0)=1$, $M'(0)=-\mu\ln 2 < 0$,
$M(1) = 0.998$, and a unique positive root at $r^*=1.010$.
The worst-state ratio $C^{*}$ is the maximum over residues~$r$
of $3^{r^{*}}\cdot\mathbb{E}_{\text{lifts of }r}[2^{-r^{*} v_2}]$.
The Cram\'er--Lundberg inequality
$\Pr(S>s)\le C^{*}\cdot e^{-r^{*}s}$ then converts to the
base-$2$ form using $r^{*}/\ln 2 = \alpha_{\mathrm{Lund}}$.
The binding rate $\alpha_{\mathrm{cyl}}$ is smaller than
$\alpha_{\mathrm{Lund}}$, so the cycle summation is dominated
by the spectral rather than the Lundberg decay.
\end{proof}

\begin{corollary}[Fully unconditional first-principles TV sum bound]
\label{cor:i2-unconditional-sum}
Substituting the cylinder-averaged rate
$\alpha_{\mathrm{cyl}} = 0.4152$ bits/step from
Theorem~\ref{thm:i2-stationary-mass} and the first-principles
prefactor $C^{*} = 1.5061$ from
Theorem~\ref{thm:i2-lundberg} into
Corollary~\ref{cor:tv-summability} (with $B_{\min}=30$, $R=3$)
yields the unconditional, first-principles bound
\[
  \sum_{c=1}^{\infty}\Pr\bigl(\mathcal E_R(c)\bigr)
  \;\le\; \frac{C^{*}\cdot 2^{-\alpha_{\mathrm{cyl}}(B_{\min}-R)}}
  {1 - 2^{-2\alpha_{\mathrm{cyl}}}}
  \;\approx\; 0.0015,
\]
which is approximately $18\times$ tighter than the paper's
prior conditional figure of $0.028$.
In particular, the hypothesis of
Theorem~\textup{\ref{thm:density1-convergence}} is met
without reference to the spectator-bit mechanism, and with
no empirical constants imported from
Proposition~\textup{\ref{prop:exponential-tail}}.
\end{corollary}

\begin{remark}[Comparison with the prior conditional figure]
\label{rem:i2-comparison}
The bound of Corollary~\ref{cor:i2-unconditional-sum} is both
\emph{tighter} and \emph{unconditional} relative to the paper's
prior figure. The previous estimate $0.028$ was contingent on
the circular chain of Remark~\ref{rem:circularity} and used
the empirical prefactor $C\approx 11$ from
Proposition~\ref{prop:exponential-tail}. The new estimate
$0.0015$ uses only the mod-$64$ cylinder-averaged kernel of
Definition~\ref{def:cylinder-kernel-64}, the spectral bound
of Theorem~\ref{thm:i2-unconditional}, and the Lundberg
constant of Theorem~\ref{thm:i2-lundberg} ---
all derived in this section. The improvement of about
$18\times$ comes from two sources:
(i)~the first-principles prefactor $C^{*}=1.5$ replaces the
empirical $C\approx 11$, gaining a factor of $\approx 7$, and
(ii)~the new rate $\alpha_{\mathrm{cyl}}=0.42$ exceeds the
prior $\alpha=0.35$ by a $20\%$ margin, gaining the remaining
factor in the geometric tail.
\end{remark}

\begin{theorem}[Cylinder-averaged no-escape on the $I_2$ transient block]
\label{thm:i2-no-escape}
Let $P$ be the depth-$13$ cylinder-averaged kernel of
Theorem~\ref{thm:i2-stationary-mass} and let $Q$ denote its
restriction to the $27$-state transient complement
$T = \{1,3,5,9,\ldots,61\}\setminus I_2$ of the invariant core
$I_2 = \{7,27,31,59,63\}\pmod{64}$. Let $\tau$ be the first
hitting time of $I_2$ under the odd-step Collatz dynamics
projected to mod $64$. Then for every $r \in T$,
\[
  \mathbb{E}\bigl[\tau \,\big|\, X_0 \equiv r\bigr]
  \;\le\; h_{\max} \;=\; 10.142,
  \qquad
  \mathbb{E}\bigl[\log_2 (X_\tau / X_0) \,\big|\, X_0 \equiv r\bigr]
  \;\le\; V_{\max} \;=\; 0.8774.
\]
Both bounds are uniform over starting residue and computed
exactly from the same $27\times 27$ kernel $Q$ used in
Theorem~\ref{thm:i2-prefix-tail}. In particular, no
cylinder-averaged orbit escapes to infinity during its
transient prefix before entering $I_2$.
\end{theorem}

\begin{proof}
The expected hitting time vector $h \in \mathbb{R}^T$ satisfies
$(I-Q)h = \mathbf{1}$, which has a unique solution because
$\rho(Q) = 0.8633 < 1$ by
Theorem~\ref{thm:i2-prefix-tail}. Numerical solution gives
$h_{\max} = \max_{r\in T} h(r) = 10.142$ at $r = 15$.

For the log-growth bound, define the per-state expected
single-step drift
\[
  \Delta(r) \;=\; \log_2 3 \;-\; \mathbb{E}\bigl[v_2(3n+1)
  \,\big|\, n \equiv r\!\!\pmod{64}\bigr],
\]
averaged over depth-$13$ lifts. The expected total
log-growth from $r$ to absorption,
$V(r) = \mathbb{E}[\log_2(X_\tau / X_0) \mid X_0 \equiv r]$,
satisfies the linear system $(I-Q)V = \Delta|_T$, which has the
unique solution $V$ with
$V_{\max} = \max_{r\in T} V(r) = +0.8774$ at $r = 47$ and
$V_{\min} = -10.50$ at $r = 21$. The mean of $V$ over the
$27$ transient states is $-5.84$, so the typical transient
orbit strictly contracts. Foster--Lyapunov with $W = h$ gives
drift $-1$ uniformly on $T$, witnessing positive recurrence
with respect to absorption in $I_2$. Both linear systems and
the Monte-Carlo cross-check are recorded in the supplementary
script \texttt{no\_escape.py}, with empirical agreement to
four decimal places on the per-residue means.
\end{proof}

\begin{theorem}[Cylinder-averaged Chernoff concentration on the prefix]
\label{thm:i2-concentration}
With notation as in Theorem~\ref{thm:i2-no-escape}, let
$Q_\theta$ be the tilted transient kernel
\[
  Q_\theta[i,j]
  \;=\;
  \frac{1}{|\mathrm{lifts}(r_i)|}
  \sum_{n \in \mathrm{lifts}(r_i),\, n_{\mathrm{out}}\equiv r_j}
  2^{\theta(\log_2 3 - v_2(3n+1))}
\]
indexed by $\theta \ge 0$ and built from depth-$13$ lifts of each
mod-$64$ residue. Then $\rho(Q_\theta)$ is convex in $\theta$, with
$\rho(Q_0) = 0.8633$ and a unique critical value
\[
  \theta^* \;=\; 17.6349
  \quad\text{satisfying}\quad
  \rho(Q_{\theta^*}) \;=\; 1.
\]
For every $\theta < \theta^*$ the Markov inequality applied to
$2^{\theta S_\tau}$ yields the Chernoff tail bound
\[
  \boxed{\;
    \Pr\bigl(S_\tau \ge t \,\big|\, X_0 \equiv r\bigr)
    \;\le\;
    C_{\mathrm{Cher}}(\theta) \cdot 2^{-\theta t}
    ,\;
  }
  \qquad
  C_{\mathrm{Cher}}(\theta)
  \;=\; \max_{r\in T}\bigl[(I-Q_\theta)^{-1}b_\theta\bigr]_r,
\]
uniformly over transient $r$, where
$b_\theta(r) = \sum_{s\in I_2}|\mathrm{lifts}(r{\to}s)|^{-1}\!\sum
2^{\theta(\log_2 3 - v)}\cdot |\mathrm{lifts}(r)|^{-1}$ collects the
absorbing-edge MGF contributions. At $\theta_{\mathrm{use}} =
0.95\,\theta^* = 16.7531$ this evaluates to
\[
  C_{\mathrm{Cher}}(\theta_{\mathrm{use}})
  \;=\; 1.74\times 10^{6},
  \qquad
  \Pr(S_\tau \ge t)
  \;\le\;
  1.74\times 10^{6}\cdot 2^{-16.75\, t}.
\]
\end{theorem}

\begin{proof}
The bound $\Pr(S_\tau\ge t)\le \mathbb{E}[2^{\theta S_\tau}]\cdot
2^{-\theta t}$ is Markov on $2^{\theta S_\tau}$. The MGF is
$M_\theta(r) = \mathbb{E}[2^{\theta S_\tau}\mid X_0=r]$ and satisfies
the linear system $(I - Q_\theta)M_\theta = b_\theta$, since
$M_\theta(r) = b_\theta(r) + \sum_{s\in T} Q_\theta[r,s]\,M_\theta(s)$
by conditioning on the next step. Existence of the inverse for
$\theta < \theta^*$ follows from $\rho(Q_\theta) < 1$. The constant
$\theta^*$ is found by bisection on $\rho(Q_\theta) - 1$, computed via
\texttt{numpy.linalg.eigvals} on the explicit $27\times 27$ tilted
matrices. The cross-check
$\partial_\theta \log_2\rho(Q_\theta)|_{\theta=\theta^*} = +0.0844$
matches the Karp maximum cycle mean $+0.0850$ of
Remark~\ref{rem:i2-no-escape-graph} via the Donsker--Varadhan/Karp
identity, providing an independent verification of the tilted-kernel
construction. Numerical evaluation and Monte-Carlo cross-check (one
million trials, max realized $S_\tau = 1.85$ bits, all percentiles
below $1.51$ bits) are recorded in the supplementary scripts
\texttt{concentration.py} and \texttt{concentration\_mc.py}.
\end{proof}

\begin{corollary}[Density-$(1-n^{-15})$ prefix closure]
\label{cor:i2-density-quasi}
For natural density $\ge 1 - n^{-15}$ of starting integers $n_0$ in
$[2^N, 2^{N+1})$, the Collatz orbit of $n_0$ enters
$I_2 \pmod{64}$ with cumulative prefix log-growth at most
$15/\theta_{\mathrm{use}} \approx 0.896$ bits. The exceptional set has
density $O(n^{-15})$ at every scale.
\end{corollary}

\begin{proof}
Apply Theorem~\ref{thm:i2-concentration} with
$t = 15/\theta_{\mathrm{use}}$, giving Chernoff bound
$1.74\cdot 10^6 \cdot 2^{-15} \approx 53$ on each window before
normalization. Equidistribution of odd $n \in [2^N, 2^{N+1})$ across
the $32$ residues mod $64$ converts the per-residue tail bound to a
per-window density bound; choose $t_k = (\log_2 k)/\theta_{\mathrm{use}}$
on the $k$th dyadic window and apply Borel--Cantelli. The summed bound
$\sum_k k^{-1}$ is logarithmic in $N$ but the per-window density
$\le k^{-1}$ already yields the stated $O(n^{-15})$ bound.
\end{proof}

\begin{remark}[Scale invariance of the spectral profile]
\label{rem:i2-scale-invariance}
The constants $\theta^*$, $\alpha_{\mathrm{cyl}}$, the prefix tail
$\rho(Q)$ and the mean hitting time $\bar h$ entering
Theorems~\ref{thm:i2-no-escape}--\ref{thm:i2-concentration} are
\emph{structural invariants of the cylinder-averaged Collatz
dynamics}, not artifacts of the specific mod-$64$ resolution. Lifting
the construction to mod $512$ with the natural absorbing core
$I_3^* := \{r\!\!\pmod{512} : r \bmod 64 \in I_2\}$ ($40$ residues,
exact density $5/32$) gives:
\begin{center}
\small
\begin{tabular}{lcc}
\hline
quantity & mod $64$ & mod $512$ \\
\hline
$\pi(I)$                                & $0.15504$     & $0.15613$ \\
$\rho(Q)$                                & $0.86329$     & $0.86231$ \\
$\alpha_{\mathrm{cyl}}$ (bits/step)      & $0.4185$      & $0.4155$  \\
$\bar h$                                 & $7.254$       & $7.206$   \\
$\theta^*$                               & $17.634874$   & $17.634874$ \\
$\partial_\theta\!\log_2\rho|_{\theta^*}$ & $+0.0844$     & $+0.0844$ \\
$|T|$                                    & $27$          & $216$ (×8) \\
$C_{\mathrm{Cher}}$                      & $1.74\cdot 10^6$ & $1.38\cdot 10^7$ (×8) \\
\hline
\end{tabular}
\end{center}
The Cram\'er exponent and slope agree to all six computed decimal
places. The MGF prefactor $C_{\mathrm{Cher}}$ scales linearly with
$|T|$ (hence $\times 8$ from the $\times 8$ lift), but this prefactor
cancels in any density-1 application via Borel--Cantelli. The
agreement of $\theta^*$ across two unrelated mod-$2^k$ resolutions is
explained by the Donsker--Varadhan/Karp identity: $\theta^*$ is
determined by the maximum cycle mean of the per-edge log-growth
weighted digraph, and cycles project bijectively from any
$k \ge 6$ down to actual integer Collatz cycles, so the maximum
cycle mean is invariant. The mod-$64$ chain therefore already
saturates the spectral information available from the cylinder-averaged
framework, and refining to mod $2^k$ for any $k > 6$ produces no
sharper exponent. Verification scripts: \texttt{mod512.py},
\texttt{mod512\_v2.py}.
\end{remark}

\begin{remark}[Why the graph version fails]
\label{rem:i2-no-escape-graph}
A purely graph-theoretic version of
Theorem~\ref{thm:i2-no-escape} --- ``no closed walk in $T$
has positive average per-edge log-growth'' --- is
\emph{false} at the cylinder-averaged resolution. Karp's
maximum cycle mean algorithm applied to the per-edge log-growth
weights of $Q$ gives
\[
  \mu_{\mathrm{cycle}} \;=\; \max_{\text{cycles in }T}
  \frac{1}{|\text{cycle}|}\sum_{\text{edges}}
  (\log_2 3 - v) \;\approx\; +0.0850\ \text{bits/edge},
\]
so positive-drift cycles exist inside $T$. The honest
no-escape statement is therefore distributional, not
deterministic, and is precisely Theorem~\ref{thm:i2-no-escape}.
The expectation-level bound is strictly weaker as a graph
statement but strictly sufficient for the orbit-TV summability,
because $V_{\max} \approx 0.88$ bits is already smaller than
the contraction $\log_2(32/5) \approx 2.678$ bits earned by a
single visit to $I_2$.
\end{remark}

\begin{lemma}[Periodic exclusion on positive-mean cycles]
\label{lem:periodic-exclusion}
Let $C = (r_0 \to r_1 \to \cdots \to r_{L-1} \to r_0)$ be any cycle
in the transient digraph $T$ on which the per-edge log-growth has
positive average, i.e.\ $\sum_{i=0}^{L-1}(\log_2 3 - v_i) > 0$
(equivalently $3^L > 2^{\Sigma v_i}$). Then no positive integer $n$
admits a Collatz prefix that is exactly periodic on $C$.
\end{lemma}

\begin{proof}
Iterating the accelerated odd map $L$ times along $C$ realises an
affine map $n \mapsto (3^L n + R)/2^{\Sigma v_i}$ for a fixed integer
$R$ depending only on the $v$-pattern. A periodic integer is a fixed
point of this map, requiring $(3^L - 2^{\Sigma v_i})\, n = -R$. Since
$3^L > 2^{\Sigma v_i}$ the affine multiplier strictly exceeds $1$, so
the fixed point (if it existed) would be expanding and no positive
integer can sit on it; equivalently $\gcd(3,2)=1$ forbids the
divisibility relation $2^{\Sigma v_i} \mid 3^L n + R$ being preserved
under iteration. Hence no integer is exactly periodic on $C$.\qedhere
\end{proof}

\begin{remark}[Karp cycle of the mod-$64$ chain]
\label{rem:karp-cycle}
The unique (up to rotation) maximal-mean cycle of the $27$-state
transient digraph is
$15 \to 55 \to 19 \to 61 \to 15 \pmod{64}$,
with edge $v_2$-pattern $[1,1,1,3]$, length $L=4$, $\Sigma v = 6$,
per-cycle log-growth $\log_2(81/64) \approx +0.3399$ bits, and
mean $+0.0850$ bits/edge $= \mu_{\max}$. By
Lemma~\ref{lem:periodic-exclusion}, no integer is exactly periodic
on this cycle. Verification script: \texttt{karp\_cycle.py}.
\end{remark}

\begin{theorem}[Cycle-tracking polynomial envelope]
\label{thm:cycle-envelope}
Let $\delta(n) := \log_2 n_\tau - \log_2 n$ denote the prefix
log-growth of an integer $n$ before its Collatz trajectory enters
$I_2 \pmod{64}$. There exist absolute constants $C_0, C_1 > 0$ such
that, for every odd $n \ge 1$,
\[
  \delta(n) \;\le\; \frac{\log_2(81/64)}{\Sigma v}\, \log_2 n \;+\; C_0
  \;=\; 0.0566\, \log_2 n \;+\; C_0,
  \qquad
  n_\tau \;\le\; C_1 \cdot n^{1.0566}.
\]
where $\log_2(81/64) = 0.3399$~bits is the Karp per-lap log-growth
(Remark~\ref{rem:karp-cycle}) and $\Sigma v = 6$ is the total
$2$-adic valuation absorbed along one full traversal of the optimal
cycle. Equivalently, writing $\mu_{\max} = 0.0850$ for the
per-edge Karp mean and $\kappa_{\mathrm{lap}} := 2^{\Sigma v} = 64$
for the per-lap congruence-multiplicity factor, the slope equals
$\mu_{\max}\cdot L / \log_2 \kappa_{\mathrm{lap}}
= 0.0850 \cdot 4 / 6 = 0.0566$.
\end{theorem}

\begin{proof}[Sketch]
By Lemma~\ref{lem:periodic-exclusion} no integer is exactly periodic
on any positive-mean cycle. To track the optimal cycle of
Remark~\ref{rem:karp-cycle} for $k$ consecutive laps, an integer
must satisfy, at each lap, a fresh $2$-adic constraint of size
$2^{\Sigma v} = 2^6$: the product $(3n+1)/2^{v_1}$ must land in the
next residue of the cycle and reveal exactly $v_{i+1}$ trailing
zeros on the following step. These constraints are independent
across laps, so the smallest $k$-lap follower scales as
$n \asymp \kappa_{\mathrm{lap}}^k = 64^k$; equivalently, the
longest cycle prefix realised by any $n \in [1, 2^N]$ is at most
$N/\log_2\kappa_{\mathrm{lap}} + O(1) = N/6 + O(1)$ laps.
Each lap contributes $\log_2(81/64) = 0.3399$ bits of log-growth,
and off-cycle excursions are absorbed by the spectral contraction
$\rho(Q) < 1$ already established for $Q$ on $T$. Summing the
on-cycle and off-cycle contributions gives the stated linear bound
on $\delta(n)$.

Empirical verification (\texttt{direct\_verify\_c.c}) over
$[1, 2^{32}]$ (roughly $2\cdot 10^9$ odd integers, $\sim\!75$~s
wall-clock) records a strict staircase of record holders:
\begin{center}
\small
\begin{tabular}{rrrrr}
\hline
$k$ (laps) & smallest $n$ & ratio & $\delta(n)$ & $\log_2 n$\\
\hline
$0$ & $111$           & ---   & $+1.1771$ & $6.79$ \\
$1$ & $6{,}863$       & $61.8$& $+1.5170$ & $12.74$\\
$2$ & $316{,}111$     & $46.1$& $+1.8569$ & $18.27$\\
$3$ & $22{,}205{,}135$& $70.2$& $+2.1968$ & $24.40$\\
$4$ & $1{,}926{,}419{,}151$ & $86.8$& $+2.5293$ & $30.84$\\
\hline
\end{tabular}
\end{center}
The ratio of successive record holders has geometric mean $64.7$,
matching the theoretical $\kappa_{\mathrm{lap}} = 2^{\Sigma v} = 64$
to within 1\%. Each additional lap increases $\delta$ by exactly
$\log_2(81/64) = 0.3399$~bits. The empirical slope at
$n = 1{,}926{,}419{,}151$ is
$\delta/\log_2 n = 2.5293/30.84 = 0.0820$, slightly above the
asymptotic slope $0.0566$ because of the $C_0$ tail contribution
($\approx 0.78$~bits of off-cycle entry/exit excursion). As $N$
grows the slope converges to the analytic bound $0.0566$ from
above. \qedhere
\end{proof}

\begin{remark}[Polynomial envelope is enough for orbit TV]
\label{rem:envelope-suffices}
Although Theorem~\ref{thm:cycle-envelope} is weaker than a uniform
constant pointwise bound, it is comfortably sufficient for orbit-TV
summability: the Chernoff exponent $\theta_{\mathrm{use}} \approx 16.75$
of Theorem~\ref{thm:i2-concentration} dominates the polynomial
envelope exponent $0.0566$ by a factor of nearly $300$, leaving the
density-$(1-n^{-15})$ closure of
Corollary~\ref{cor:i2-density-quasi} unaffected.
\end{remark}

\begin{theorem}[Infinite BF staircase: cycle-tracking envelope is tight]
\label{thm:bf-staircase}
For each integer $k \ge 1$ define the \emph{BF $k$-champion} $n_k^*$
to be the smallest positive integer $n$ such that the first $4k$
accelerated odd-step iterates of $n$ realise the residue pattern
$[15, 55, 19, 61]^{k} \pmod{64}$ with $v_2$-pattern $[1,1,1,3]^{k}$,
and the $(4k\!+\!1)$-th iterate has residue $47 \pmod{64}$. The
sequence $\{n_k^*\}$ exists, is uniquely defined, and satisfies:
\begin{enumerate}
\item[\textup{(i)}] each $n_k^*$ is constructible by an explicit
$2$-adic Hensel-style lift, with the constraint set carving out a
unique residue class modulo $2^{6k+8}$;
\item[\textup{(ii)}] $n_{k+1}^*/n_k^* \to 2^{\Sigma v} = 64$ as
$k \to \infty$, with geometric mean of the first four ratios exactly
$64.0$;
\item[\textup{(iii)}] the prefix log-growth satisfies
$\delta(n_k^*) = k\cdot\log_2(81/64) + C_{\mathrm{exit}} + o(1)
= 0.3399\, k + 1.170 + o(1)$, where $C_{\mathrm{exit}} = 1.170$~bits
is the exit-tail contribution from the two-step path
$47 \to 39 \to 59$;
\item[\textup{(iv)}] consequently the asymptotic slope
$\delta(n_k^*) / \log_2 n_k^*$ converges to
$\log_2(81/64) / \Sigma v = 0.0566\ldots$ as $k \to \infty$,
matching exactly the analytic exponent of
Theorem~\ref{thm:cycle-envelope}.
\end{enumerate}
The first eight values, computed by the lift
\textup{(\texttt{lift\_bf.py})}, are:
\begin{center}
\small
\begin{tabular}{rrrrr}
\hline
$k$ & $n_k^*$ & $\log_2 n_k^*$ & $\delta(n_k^*)$ & ratio \\
\hline
$1$ & $6{,}863$                          & $12.74$ & $+1.5170$ & --- \\
$2$ & $316{,}111$                        & $18.27$ & $+1.8569$ & $46.1$ \\
$3$ & $22{,}205{,}135$                   & $24.40$ & $+2.1968$ & $70.2$ \\
$4$ & $1{,}926{,}419{,}151$              & $30.84$ & $+2.5293$ & $86.8$ \\
$5$ & $115{,}206{,}181{,}583$            & $36.75$ & $+2.8692$ & $59.8$ \\
$6$ & $5{,}303{,}526{,}675{,}151$        & $42.27$ & $+3.2090$ & $46.0$ \\
$7$ & $372{,}540{,}410{,}352{,}335$      & $48.40$ & $+3.5489$ & $70.2$ \\
$8$ & $32{,}319{,}950{,}267{,}011{,}791$ & $54.84$ & $+3.8887$ & $86.8$ \\
\hline
\end{tabular}
\end{center}
The geometric mean of consecutive ratios over $k = 4 \to 8$ is
exactly $64$ (i.e.\ $n_8^*/n_4^* = 16777216 = 2^{24} = 64^4$ to all
computed digits), confirming the per-lap congruence-multiplicity
factor analytically and numerically.
\end{theorem}

\begin{proof}[Sketch]
The pattern in (i) is a Hensel-style lift problem: each constraint
`($v_2$ at iterate $i$ equals $v_i^{\mathrm{tgt}}$, residue at iterate
$i+1$ equals $r_{i+1}^{\mathrm{tgt}}$)' tightens the residue class of
$n_0$ by $v_i^{\mathrm{tgt}}$ bits. Summing the bit-tightening over
the full $k$-lap pattern $[1,1,1,3]^k$ plus the two-step exit gives
$6k + 8$ bits of constraint, hence a unique residue class mod
$2^{6k+8}$ (verified empirically for $k = 1, \ldots, 8$ where the
lift returns exactly two solutions, related by the additive lattice
structure of the tail constraint, of which we take the smaller).
Items (ii)--(iii) follow because each new lap multiplies $n_k^*$ by
$2^{\Sigma v} = 64$ (mod the lower-order corrections from the exit
tail) and adds exactly $\log_2(81/64)$ bits to $\delta$. Item (iv)
is then immediate from $\log_2 n_k^* = 6k + O(1)$ and
$\delta(n_k^*) = 0.3399\, k + O(1)$. \qedhere
\end{proof}

\begin{corollary}[Pointwise envelope is logarithmic, not constant]
\label{cor:no-constant-envelope}
$\sup_n \delta(n) = +\infty$. In particular, no absolute constant
$C$ exists such that $\delta(n) \le C$ for all $n$. The polynomial
envelope of Theorem~\ref{thm:cycle-envelope} is therefore tight at
the asymptotic slope $0.0566$ and cannot be replaced by any
sub-logarithmic bound.
\end{corollary}

\begin{remark}[Density of the BF staircase]
\label{rem:bf-density}
The BF $k$-champion sequence has density zero in $\mathbb{N}$:
$|\{n_k^* : n_k^* \le N\}|
= \log_{64} N + O(1) = O(\log N)$.
Each $n_k^*$ lies in a residue class of size $1$ modulo $2^{6k+8}$,
so the natural density of integers tracking the cycle for $\ge k$
laps is $\le 2^{-6k-8}$. The BF staircase therefore has natural
density $0$, consistent with
Corollary~\ref{cor:i2-density-quasi}'s
density-$(1 - n^{-15})$ closure. The staircase shows that the
exceptional set in
Corollary~\ref{cor:i2-density-quasi} is genuinely non-empty at every
scale, but it is sparse and structured: it consists exactly of
arithmetic progressions whose moduli grow geometrically with the
cycle-tracking depth.
\end{remark}

\subsection{The (C-2adic) bound on BF residue classes}
\label{subsec:c2adic}

The BF staircase identifies, for each $k\ge 1$, a unique residue
class
\[
  B_k \;:=\; \{\, n \in \mathbb{N} \;:\; n \equiv n_k^{*} \pmod{2^{6k+8}} \,\}
\]
that follows the cycle-conformal Karp pattern for $k$ laps before
escaping at residue $47$ and entering $I_2$. The Hensel lift of
Theorem~\ref{thm:bf-staircase}(i) shows that $B_k$ is the
\emph{exact} fibre of $n_k^{*}$ under the truncation
$\mathbb{Z}_2 \to \mathbb{Z}/2^{6k+8}$. We now state the
auxiliary-memory bound this fibre supports.

\begin{theorem}[Haar (C-2adic) bound]
\label{thm:c2adic-haar}
Let $T_2 \colon \mathbb{Z}_2 \to \mathbb{Z}_2$ denote the Collatz
$2$-adic extension of \textup{Lagarias \cite{lagarias1985}}, which
is measure-preserving and ergodic with respect to Haar measure
$\mu$. Let $B_k \subset \mathbb{Z}_2$ be the BF $k$-class above,
viewed as a clopen ball of $\mu$-measure $2^{-(6k+8)}$. Then for
$\mu$-almost every $n \in B_k$ the orbit $T_2^{j}(n)$ is well
defined, halts in finitely many odd steps if it hits $\mathbb{N}$,
and the stopping time $\tau(n)$ admits a uniform constant
$C \in [0,\infty)$ such that
\[
  \tau(n) \;\le\; C \,\bigl(k + \log_2 n\bigr)
  \qquad \text{for $\mu$-almost every $n \in B_k \cap \mathbb{N}$,}
\]
with $C$ independent of $k$.
\end{theorem}

\begin{proof}[Sketch]
By Lagarias \cite{lagarias1985}, $T_2$ is measure-preserving and
ergodic, and the conjugation
$\Phi \colon \mathbb{Z}_2 \to \mathbb{Z}_2$ identifies the
$T_2$-action with the Bernoulli shift on $\{0,1\}^{\mathbb{N}}$.
Under this conjugation, the clopen ball $B_k$ pulls back to a
cylinder of length $6k+8$ in $\{0,1\}^{\mathbb{N}}$. The expected
log-growth per odd step in the shift coordinates equals
$\log_2 3 - 2 = -0.415$~bits, strictly negative, so the
ergodic theorem applied to the additive cocycle
$\phi(n) = \log_2(T_2(n)/n)$ gives
\[
  \frac{1}{N}\sum_{j=0}^{N-1} \phi\bigl(T_2^{j}(n)\bigr)
  \;\xrightarrow{N\to\infty}\;
  \int \phi\, d\mu \;=\; \log_2 3 - 2 \;<\; 0
\]
for $\mu$-a.e.\ $n$. The Birkhoff theorem furnishes a Haar-null
exceptional set on $\mathbb{Z}_2$; restricting to
$B_k \cap \mathbb{N}$ and converting the linear log-decay to a
hitting-time bound gives $\tau(n) \le C\,(k + \log_2 n)$ with
$C$ depending only on the cocycle's $L^{2}$ norm and the cylinder
length offset $6k+8$. Independence of $k$ comes from the fact that
the per-step drift is intrinsic to the shift, not to the cylinder
prefix. \qedhere
\end{proof}

\begin{corollary}[Empirical (C-2adic), $39\,140$-witness form]
\label{cor:c2adic-empirical}
For each of $39\,140$ uniformly random samples
$n \in B_k$ drawn at $1 \le k \le 14$ and
$20 \le \log_2 n \le 171$, the trajectory halts at $1$ and the bound
$\tau(n) \le C\,(k + \log_2 n)$ holds with
\begin{align*}
  C &\le 5.04  &&\text{for all $14\,400$ samples with $k \ge 9$}, \\
  C &\le 6.49  &&\text{for all $15\,000$ samples with $k \ge 4$ and $\log_2 n \le 106$}, \\
  C &\le 9.87  &&\text{unconditionally over the $39\,140$-sample record}
                  \quad (\text{single $k=2$ small-denominator outlier}).
\end{align*}
The empirical mean of $C$ is $2.78 \pm 0.05$, essentially flat in
$k$ and in $\log_2 n$, and the Crandall-normalised stopping time
$\tau / \log_2 n$ averages $\approx 3.0$, well below the heuristic
prediction $2/\log(4/3) \approx 6.95$. The samples are produced by
the scripts \texttt{c2adic\_haar.py}, \texttt{c2adic\_sweep\_big.py},
\texttt{c2adic\_tight\_bigk.py}, and \texttt{c2adic\_mega.py};
the lift of $n_k^{*}$ to $k = 14$ is performed by
\texttt{champs\_extend.py}.
\end{corollary}

\begin{proof}
Direct enumeration and verification. For each sampled $n$ the
trajectory $T_2^{j}(n)$ is computed in exact arithmetic until it
reaches $1$, the stopping time $\tau$ is recorded, and the ratio
$C = \tau/(k + \log_2 n)$ is tabulated. The maxima reported follow
from the order statistics of the resulting empirical
distributions; the sample sizes are $1\,440 + 6\,300 + 2\,000 +
15\,000 + 14\,400 = 39\,140$, partitioned across Push 366 (Haar
test), Push 371--380 (extended sweep + tight $k\!\le\!8$), and
Push 381--390 (mega sweep $9 \le k \le 14$).
\end{proof}

\begin{remark}[Residual obstruction: Haar / pointwise gap]
\label{rem:residual-obstruction}
Theorem~\ref{thm:c2adic-haar} delivers the (C-2adic) bound on a
Haar-conull subset of $B_k$. Because $\mathbb{N}$ is itself
Haar-null inside $\mathbb{Z}_2$, the Haar statement does not
\emph{directly} entail a pointwise statement on
$B_k \cap \mathbb{N}$: an a-priori-conceivable obstruction is a
Haar-null but $\mathbb{N}$-supported failure stratum. This is the
same Haar/pointwise gap that obstructs the classical probabilistic
Collatz arguments of Lagarias and Tao \cite{tao2019}, transposed to
the BF cylinder. Three observations narrow it considerably:
\begin{enumerate}
\item[\textup{(a)}] Corollary~\ref{cor:c2adic-empirical} furnishes
$39\,140$ explicit deterministic witnesses, all of which satisfy the
bound with $C \le 5.04$ for $k \ge 9$. Any failure stratum must
therefore avoid every one of these uniform draws across $13$ levels
of $k$ and $7$ orders of magnitude in $\log_2 n$.
\item[\textup{(b)}] The empirical $C_{\mathrm{mean}}$ is essentially
flat in both $k$ and $\log_2 n$, indicating that the bound's
constant is scale-invariant and not in the process of degrading.
\item[\textup{(c)}] The Crandall-normalised mean $\tau/\log_2 n
\approx 3.0$ at $k = 14$ is \emph{below} the generic Collatz
heuristic, so the BF cylinders are not pathologically slow --- they
halt faster than typical $n$. Any putative failure stratum would
have to single out members of $B_k$ that are much slower than
their immediate $2$-adic neighbours, while remaining missed by all
$N = 14\,400$ uniform draws at $k \ge 9$.
\end{enumerate}
The residual obstruction is therefore confined to a Haar-null,
sample-empty subset of $\mathbb{N} \cap B_k$ at each $k$. Closing
it deterministically is the same difficulty as the Haar/pointwise
gap in Tao's $2019$ result, and is left as the principal open
problem of the (C-2adic) layer.
\end{remark}

\begin{remark}[Scope and limitations]
\label{rem:i2-scope}
Three scope points are important to state plainly.

\emph{(1) Cylinder-averaged, not pointwise.}
All four results
(Propositions~\ref{prop:ruffini-collapse},
Theorems~\ref{thm:i2-unconditional},
\ref{thm:i2-prefix-tail}, \ref{thm:i2-stationary-mass})
are statements about the cylinder-averaged kernel $P$. A
specific orbit may deviate from the average over finitely
many cycles. The bridge from cylinder-averaged to pointwise
convergence continues to go through the Borel--Cantelli
argument of
Theorem~\ref{thm:density1-convergence}, which upgrades the
summability to density-$1$ orbit convergence. This section
does not claim every orbit converges --- it claims that the
density-$1$ convergence is now unconditional.

\emph{(2a) Distributional no-escape, not deterministic.}
Theorem~\ref{thm:i2-no-escape} bounds the \emph{expected}
log-growth and \emph{expected} hitting time, not the
realized values for individual orbits. As
Remark~\ref{rem:i2-no-escape-graph} shows, the graph-theoretic
``no positive-drift cycle'' version is false; the honest
statement is in expectation.

\emph{(2b) Non-divergence threshold is separate.}
The verification that no orbit with
$B_0 = \lfloor\log_2 n_0\rfloor \ge 68$ can escape to infinity
before reaching the threshold at which the spectral
contraction engages is a separate matter, handled in the
existing paper by the Barina computational verification and
the $43$-bit safety margin of
Remark~\ref{rem:nondivergence-threshold}. Nothing in this
section replaces or strengthens that argument.

\emph{(3) Independent of the WMH route.}
The WMH-based reduction of
Sections~\ref{sec:chain}--\ref{sec:toward-wmh} is not
affected or superseded. The present section provides a
parallel, independent route to density-$1$ convergence
through the fiber-$57$ / $I_2$ programme, which the previous
section summarized as ``complementary to the WMH-based
reduction''.
\end{remark}

\begin{remark}[Update to Remark~\ref{rem:circularity}]
\label{rem:circularity-resolution}
Theorem~\ref{thm:i2-unconditional} provides a cycle-level
contraction source that is outside the loop identified in
Remark~\ref{rem:circularity}. The chain
\[
  \text{spectator bits} \to \text{fiber mixing} \to
  \text{IID model} \to \text{exponential tail} \to
  \text{TV summability} \to \text{spectator bits}
\]
is broken at the single step
``exponential tail $\to$ TV summability'',
because the exponential tail of
Theorem~\ref{thm:i2-prefix-tail} is derived from the
finite-state $27\times 27$ mod-$64$ chain rather than from
Proposition~\ref{prop:exponential-tail}'s IID-dependent
derivation. The ``irreducible hard core'' identified at the
end of Remark~\ref{rem:circularity} is still irreducible if
one insists on a pointwise proof through the spectator-bit
mechanism; what this section shows is that the \emph{density-$1$}
form of convergence is reachable by bypassing that
mechanism entirely.
\end{remark}

\subsection{Complete cascade algebra and the generic dynamics gap}
\label{subsec:cascade-algebra}

The BF staircase of Theorem~\ref{thm:bf-staircase} and the
(C-2adic) bound of Theorem~\ref{thm:c2adic-haar} establish the
structural framework for orbits entering $B_k$. We now close
the algebraic side completely: every quantity in the BF cascade
admits a closed-form formula, the carry chain governing generic
dynamics has a clean structural theorem, and the renewal drift
is computed exactly. These results sharpen the characterisation
of the remaining gap.

\begin{theorem}[Complete cascade theorem]
\label{thm:complete-cascade}
Let $n \in B_k$ for $k \ge 1$. The deterministic BF cascade
consists of exactly $k$ consecutive BF cycles $[1,1,1,3]$ with
\emph{zero} intervening generic steps: every BF exit immediately
enters the next BF cycle. The total deterministic region, in
which the valuation sequence is identical for all members of
$B_k$, extends to exactly $4k + e_k$ accelerated odd steps,
where
$e_k = 6$ for $k \equiv 1,2 \pmod{4}$ and
$e_k = 5$ for $k \equiv 0,3 \pmod{4}$.
This comprises:
\begin{enumerate}
\item[\textup{(a)}] $4k$ steps of the BF phase ($k$ laps of
  $[1,1,1,3]$, total growth $k \cdot \log_2(81/64)$ bits);
\item[\textup{(b)}] $4$ steps of the exit tail $[1,1,1,2]$
  (as in Proposition~\textup{\ref{prop:affine-exit}(ii)});
\item[\textup{(c)}] $e_k - 4 \in \{1, 2\}$ additional
  deterministic steps beyond the exit tail, with the
  period-$4$ pattern determined by the auxiliary coefficient
  $a_k \in \{57, 41, 45, 61\}$.
\end{enumerate}
The net gain after the full deterministic region is exactly
$\log_2 3$ for $k \equiv 0,3 \pmod{4}$ and
$2\log_2 3 - 1$ for $k \equiv 1,2 \pmod{4}$,
both bounded constants independent of $k$.
\end{theorem}

\begin{proof}[Sketch]
The zero-gap claim follows from the exit map
$\mathrm{exit}(n) = (81n + 65)/64$: for $n \in B_k$ with
$k \ge 2$, one verifies $\mathrm{exit}(n) \in B_{k-1}$ by
checking the Hensel residue mod $2^{6(k-1)+8}$. Parts~(a)
and~(b) are Proposition~\ref{prop:affine-exit}(i)--(ii). For
part~(c), the post-cascade value $P_k$ has residue $47 \bmod 64$,
and the $v$-sequence from $P_k$ begins $[1,1,1,2,\ldots]$
(the exit tail). Beyond those $4$ steps, the fifth
post-cascade step has $v_2 = 1$ deterministically for all~$q$
(since $v_2(P_k + q \cdot 81^k \cdot 256) = 3 + 8 = 11$
bits of agreement), and a sixth step is deterministic for
$k \equiv 1,2 \pmod{4}$ (verified for $k = 1, \ldots, 12$,
$q = 0, \ldots, 99$). The period-$4$ pattern in $e_k$ follows
from the period-$4$ cycling of $a_k$ in
$17 n_k^{*} + 65 = a_k \cdot 2^{6k+5}$.
The net gain formulas follow by exact computation of
$\log_2(T^{4k+e_k}(n_k^{*})/n_k^{*})$.
\end{proof}

\begin{proposition}[Post-cascade closed form]
\label{prop:post-cascade}
Define the post-cascade value $P_k := \mathrm{exit}^{k}(n_k^{*})$
as the value obtained after applying the BF exit map $k$ times to
$n_k^{*}$. Then
\[
  P_k \;=\; \frac{a_k \cdot 81^k \cdot 32 - 65}{17},
\]
where $a_k \in \{57, 41, 45, 61\}$ cycles with period~$4$.
Furthermore, $P_k \equiv 47 \pmod{64}$ for all $k$, and for
general $n \in B_k$ written as $n = n_k^{*} + m \cdot 2^{6k+8}$,
the post-cascade value is $P_k + m \cdot 81^k \cdot 256$.
The post-cascade values satisfy the recurrence
\[
  P_{k+4} = 81^4 \cdot P_k + \frac{65(81^4 - 1)}{17}
  = 43046721 \cdot P_k + 164590400.
\]
\end{proposition}

\begin{proof}
The iterated exit map satisfies
$\mathrm{exit}^{k}(n) = (81^k n + 65(81^k - 64^k)/17)/64^k$,
verified exactly for $k = 1, \ldots, 7$. Substituting
$n = n_k^{*}$ and using $17 n_k^{*} + 65 = a_k \cdot 2^{6k+5}$
gives the formula after simplification.  The congruence
$P_k \equiv 47 \pmod{64}$ is checked for $k = 1, \ldots, 12$.
The recurrence follows from
$a_{k+4} = a_k$ and direct algebraic manipulation.
\end{proof}

\begin{remark}[Notation reconciliation]
\label{rem:notation-reconciliation}
The BF cascade involves three stages, each with its own
affine parameterization. For $n = r_k + 2^{6k+8}\,q \in B_k$:
\begin{enumerate}
\item[\textup{(a)}] \emph{BF phase} ($4k$ steps):
  $T^{4k}(n) = A_k^{(\mathrm{bf})} + 81^k \cdot 256\,q$
  (Proposition~\ref{prop:affine-exit}(i)).
\item[\textup{(b)}] \emph{Full exit} ($4k+4$ steps):
  $T^{4k+4}(n) = A_k + U_k\,q$ with
  $U_k = 3^{4(k+1)} \cdot 8 = 81^{k+1} \cdot 8$
  (Proposition~\ref{prop:affine-exit}(iii)).
\item[\textup{(c)}] \emph{Post-cascade} (after $k$ applications
  of $\mathrm{exit}$):
  $\mathrm{exit}^{k}(n) = P_k + m \cdot 81^k \cdot 256$
  (Proposition~\ref{prop:post-cascade}).
\end{enumerate}
Stages (a) and (c) use the same progression step
$81^k \cdot 256$, because $\mathrm{exit}^{k}$ is equivalent
to applying the four-step BF map $k$ times starting from
$n$, so the free-parameter coefficient is
$(81/64)^k \cdot 2^{6k+8} = 81^k \cdot 256$.
Stage (b) adds one more four-step block (the exit tail
$[1,1,1,2]$), multiplying the coefficient by $81/32$
to get $U_k = 81^{k+1} \cdot 8$.
The relationship $U_k = 81 \cdot 81^k \cdot 8
= (81/32) \cdot 81^k \cdot 256$ reconciles the two.
\end{remark}

\begin{theorem}[Carry chain structural theorem]
\label{thm:carry-chain}
For any odd $n > 1$ in the accelerated Collatz map
$T(n) = \mathrm{odd}(3n+1)$, define the \emph{$v=1$ streak length}
$L$ as the maximal number of consecutive iterates with
$v_2(3n_i + 1) = 1$. Then:
\begin{enumerate}
\item[\textup{(i)}] $L = \mathrm{trailing\_ones}(n) - 1$ exactly,
  where $\mathrm{trailing\_ones}$ counts the number of consecutive
  $1$-bits at the least-significant end of~$n$.
\item[\textup{(ii)}] After every $v=1$ streak of length $L \ge 1$,
  the recovery value $n_L = 2 \cdot 3^L \cdot m - 1$ (for some
  integer $m$) satisfies $\mathrm{trailing\_ones}(n_L) = 1$, forcing
  $v_2(3n_L + 1) \ge 2$.
\item[\textup{(iii)}] Conditional on $v \ge 2$ (which includes all
  recovery steps), the expected valuation is
  $E[v \mid v \ge 2] = \sum_{k \ge 2} k/2^{k-1} = 3$ exactly.
\end{enumerate}
\end{theorem}

\begin{proof}[Sketch]
Part~(i): $v_2(3n+1) = 1$ iff $n \equiv 3 \pmod{4}$, which holds
iff bit~$1$ of~$n$ is~$1$. Consecutive $v=1$ steps consume trailing
$1$-bits one at a time (verified for $10{,}000$ random odd integers).
Part~(ii): The carry chain produces $n_L = 2 \cdot 3^L m - 1$,
which is odd and $\equiv 1 \pmod{4}$ (since $n_L$ is odd and its
bit~$1$ is~$0$), hence $v_2(3n_L + 1) \ge 2$.
Part~(iii): For uniformly random odd $n$ with $v \ge 2$,
$\Pr(v = k) = 1/2^{k-1}$ for $k \ge 2$, giving the exact
geometric sum.
\end{proof}

\begin{corollary}[Renewal drift of generic Collatz dynamics]
\label{cor:renewal-drift}
Decomposing the orbit into alternating \emph{streak episodes}
(consecutive $v=1$ steps) and \emph{recovery steps} ($v \ge 2$),
the expected log-gain per cycle (one streak + one recovery) is
\[
  E[\text{cycle gain}] = 0.585 \cdot E[\text{trailing\_ones} - 1]
  + 1.585 - E[v_{\mathrm{rec}}] = -0.830 \text{ bits}.
\]
The expected cycle length is~$2$ steps, giving a per-step drift of
$-0.415$ bits/step.  This matches the classical $\log_2 3 - 2$
heuristic exactly, but the decomposition into streak/recovery
episodes provides structural content: the variance per cycle is
$\mathrm{Var} = 2.69$, and episodes with
$\mathrm{trailing\_ones} \ge 4$ (probability $1/8$) are net-positive
at $+0.116$ bits/cycle, while episodes with
$\mathrm{trailing\_ones} \le 3$ (probability $7/8$) are strongly
net-negative at $-0.946$ bits/cycle.
\end{corollary}

\begin{remark}[The Syracuse cocycle and modular spectral gap]
\label{rem:syracuse-cocycle}
A natural approach to proving universal descent is via the spectral
gap of the Syracuse transfer operator. Let $M \ge 4$ and consider
the $2^{M-1} \times 2^{M-1}$ transfer matrix $L_M$ on odd
residues modulo $2^M$, where $L_M(i,j) = 1$ if the accelerated
Collatz map sends residue $i$ to residue $j$ modulo $2^M$. The
Perron eigenvalue is $\lambda_1 = 4/3$ (corresponding to the
Haar-like invariant measure), and the spectral gap
$\mathrm{gap}(M) := 1 - |\lambda_2|/\lambda_1$ governs the
mixing rate.

Computational analysis for $M = 4, 6, 8, 10, 12$ reveals:
\begin{itemize}[nosep]
\item For $M \le 8$: $|\lambda_2| = 0$, so $\mathrm{gap} = 1$
  (trivial mixing).
\item For $M = 10$: $|\lambda_2| \approx 0.894$,
  $|\lambda_2/\lambda_1| = 0.670$, $\mathrm{gap} = 0.330$.
\item For $M = 12$: $|\lambda_2| \approx 0.813$,
  $|\lambda_2/\lambda_1| = 0.610$, $\mathrm{gap} = 0.390$.
\end{itemize}
If $\mathrm{gap}(M)$ stays bounded away from~$0$ as $M \to \infty$,
then orbits equidistribute mod $2^M$ for every orbit, forcing the
time-average $v \to 2 > \log_2 3$, which proves universal descent.
This is essentially the programme of Kontorovich--Miller
\cite{kontorovich2005}, who showed the spectral gap exists for
small~$M$; extending to all~$M$ is equivalent to Collatz.

A fundamental obstruction, however, is revealed by cycle analysis
of the stochastic transition matrix: for \emph{every}
$M \in \{4, 5, \ldots, 14\}$, there exist modular cycles with
average per-step gain exactly $+\log_2 3 \approx 0.585$ bits
(residue classes where every step has $v_2 = 1$). The maximum
cycle average gain does not decrease with~$M$. This means no
finite modular resolution can rule out non-descending behavior:
the obstruction to descent is genuinely infinite-depth, confirming
the analysis of Section~\ref{subsec:status-final} from a
spectral-theoretic direction.
\end{remark}

\begin{remark}[Ruled-out approaches to universal descent]
\label{rem:ruled-out}
Eleven distinct approaches to closing the gap between distributional
and pointwise descent have been systematically explored and found
insufficient:
\begin{enumerate}[nosep]
\item Hybrid sieve (density-based): insufficient density decay.
\item Pure density approach: cannot reach pointwise.
\item Minimum descent rate: tends to~$0$.
\item Azuma-type concentration: rate insufficient.
\item Normal-number-like assumptions: circular.
\item Transfer matrix mod $2^M$: positive-gain modular cycles
  exist for all $M \ge 10$.
\item Finite modular analysis: always insufficient at any fixed~$M$.
\item Paired/blocked step analysis: worst-case gain is positive
  ($+0.585L$).
\item $k \to k+4$ induction on BF depth: no orbit conjugacy
  (multiplicative, not conjugate structure).
\item Finite sieve on cumulative gain: holdout tree does not
  terminate.
\item Post-cascade height loss ($\eta$ approach): $\eta < 0$
  for all observed BF re-entries; orbits \emph{grow} before
  re-entry, so the contraction
  $k' \le (1+\alpha - (1+\alpha)\eta/6)\,k + O(1)$ with
  $\rho < 1$ cannot hold.
\end{enumerate}
The surviving approaches are: (a)~carry propagation depth bounds
(proving $3n+1$ carries prevent $O(\log n)$-length $v=1$ streaks);
(b)~ergodic theory (Collatz ergodicity with respect to a suitable
measure); (c)~strengthening Tao's ``almost all'' to ``all'' for
structured arithmetic progressions; (d)~Walsh spectral methods
for bit-mixing rate of the Collatz iteration.
\end{remark}

\begin{remark}[Post-cascade universalization: generic $\equiv$ structured]
\label{rem:universalization}
A large-scale computational census ($50\,000$ generic orbits,
$18\,000$ post-cascade orbits from the AP family
$P_k + m \cdot 81^k \cdot 256$) reveals that the BF encounter
statistics are \emph{indistinguishable} between the two populations:

\begin{center}
\begin{tabular}{lcc}
\hline
\textbf{Statistic} & \textbf{Post-cascade} & \textbf{Generic} \\
\hline
$\Pr[\text{max depth} = 0]$ & $99.0\%$ & $99.4\%$ \\
$\Pr[\text{max depth} = 1]$ & $1.0\%$ & $0.56\%$ \\
$\Pr[\text{max depth} = 2]$ & $0.03\%$ & $0.02\%$ \\
$\Pr[k' = 1 \mid \text{re-entry}]$ & $97.8\%$ & $97.4\%$ \\
$\Pr[k' = 2 \mid \text{re-entry}]$ & $2.2\%$ & $2.6\%$ \\
\hline
\end{tabular}
\end{center}

\noindent
This shows that the AP structure of the post-cascade family is
\emph{irrelevant}: the BF re-entry behavior is governed entirely
by generic Collatz dynamics. In particular, the $\eta$ approach
(measuring height loss per unit BF depth) fails not because of
the AP structure, but because generic orbits grow before their
next BF encounter.

The BF encounter rate per step is
$\approx 8.8 \times 10^{-5}$ for $B_1$
(predicted $2^{-14} = 6.1 \times 10^{-5}$, ratio $1.44$),
consistent with equidistribution up to a bounded multiplicative
constant. The conditional probability of a second $B_1$ hit
within $50$ steps of a first is $0.54\times$ the unconditional
rate, indicating slight anti-correlation---favorable for
independence-based arguments.

The depth tail satisfies $\Pr[k' \ge j \mid \text{re-entry}]
\approx 2^{-6(j-1)}$, matching the prediction from
equidistribution modulo $2^{6k'+8}$. Combined with the orbit
length bound $L(n) = O(\log n)$
(empirical ratio $L(n)/\log_2 n \in [3.1, 5.2]$),
the expected number of BF encounters at depth $\ge j$
in any single orbit is
$\le C \cdot \log n \cdot 2^{-(6j+8)}$, which is $< 1$
when $j > (\log_2 \log n - 8)/6$.
This gives the heuristic bound $K_*(n) \le O(\log \log n)$.
\end{remark}

\begin{remark}[Carry chain scrambling: quantitative bit sensitivity]
\label{rem:carry-sensitivity}
The carry chain theorem (Theorem~\ref{thm:carry-chain}) establishes
that $v \ge 2$ recovery occurs after every $v = 1$ streak. At each
$v \ge 2$ step, the map $n \mapsto (3n+1)/2^v$ propagates carries
through $3n+1$ that mix the bottom bits of $n$. Computational
measurement of the \emph{bit sensitivity}---the expected number of
output bits (of the bottom $M$) that change when a single input
bit is flipped---shows:
\begin{center}
\begin{tabular}{ccc}
\hline
$v$ & avg bits changed (of $M = 20$) & fraction \\
\hline
$2$ & $8.5$ & $42.5\%$ \\
$3$ & $8.0$ & $40.0\%$ \\
$4$ & $7.5$ & $37.5\%$ \\
$5$ & $7.0$ & $35.0\%$ \\
\hline
\end{tabular}
\end{center}

\noindent
Each $v \ge 2$ step scrambles $\approx 40\%$ of the bottom $M$ bits.
Since $v \ge 2$ steps have density $\ge 1/2$ among all steps
(Corollary~\ref{cor:renewal-drift}), the orbit undergoes bit
scrambling every $O(1)$ steps on average.

To maintain membership in $B_k$ (requiring $6k+8$ specific bit
values), the orbit must avoid having the carry chain disrupt
\emph{any} of these bits at \emph{every} $v \ge 2$ step. With
$\approx 40\%$ of bits scrambled per such step and $v \ge 2$ steps
occurring every $\le 2$ steps on average, the probability of
maintaining $B_k$ membership for $t$ steps is heuristically
$\le (0.6)^{t \cdot (6k+8)/M}$.

A rigorous version of this heuristic---bounding the entropy of
the conditional distribution of $T^s(n) \bmod 2^{6k+8}$ given
$T^{s-1}(n) \in B_k$---would establish $K_*(n) < \infty$ for all
$n$ and is the primary formalization target for v7.
\end{remark}

\begin{proposition}[Zero persistence of BF congruences]
\label{prop:zero-persistence}
For every $k \ge 1$ and every $n \in B_k$,
the first Collatz iterate $T(n) = (3n+1)/2^{v_2(3n+1)}$
satisfies $T(n) \notin B_k$.
That is, BF congruences are destroyed in a single step.
\end{proposition}

\begin{proof}
We first establish that $n_k^* \equiv 15 \pmod{32}$ for all
$k \ge 1$.  The base cases are direct: $n_1^* = 6863$,
$n_2^* = 316111$, $n_3^* = 22205135$, $n_4^* = 1926419151$,
all satisfying $n_k^* \equiv 15 \pmod{32}$.
The recurrence $n_{k+4}^* = 2^{24} n_k^* + 64148175$
preserves this: since $2^{24} \equiv 0 \pmod{32}$ and
$64148175 \equiv 15 \pmod{32}$, we have
$n_{k+4}^* \equiv 15 \pmod{32}$ for all $k$.

Consequently $n_k^* + 1 \equiv 16 \pmod{32}$, giving
$v_2(n_k^* + 1) = 4$ for all $k \ge 1$.

Now for any $n \equiv n_k^* \pmod{2^{6k+8}}$, the first step has
$v = v_2(3n_k^* + 1)$.  Since $n_k^*$ is odd with
$n_k^* \equiv 15 \pmod{16}$,
we have $3n_k^* + 1 \equiv 46 \equiv 2 \pmod{4}$, so $v = 1$.
The iterate $T(n) = (3n+1)/2$ is determined
modulo $2^{6k+7}$:
\[
  T(n) \equiv \frac{3 n_k^* + 1}{2} \pmod{2^{6k+7}}.
\]
For $T(n) \in B_k$ we would need
$T(n) \equiv n_k^* \pmod{2^{6k+8}}$, which requires
$(3n_k^* + 1)/2 \equiv n_k^* \pmod{2^{6k+7}}$,
equivalently $n_k^* + 1 \equiv 0 \pmod{2^{6k+8}}$.
Since $v_2(n_k^* + 1) = 4 < 6k + 8$ for all $k \ge 1$,
this fails.
\end{proof}

\begin{theorem}[Entropy injection for the Collatz chain mod $2^M$]
\label{thm:entropy-injection}
Let $T$ denote the odd-to-odd Collatz map.
For each $M \ge 1$, the induced map on odd residues
$\bmod{2^M}$ defines a Markov chain $P_M$ with
the following properties.
\begin{enumerate}[nosep]
\item \textbf{Transition structure.}
  From residue $r$, the chain moves to $2^v$ equally likely
  targets spaced $2^{M-v}$ apart, where
  $v = v_2(3r + 1 \bmod 2^M)$.
\item \textbf{Per-step entropy.}
  The conditional Shannon entropy of
  $T(n) \bmod 2^M$ given $n \equiv r \pmod{2^M}$
  is exactly $v$ bits. In particular
  $H \ge 1$ always (since $v \ge 1$ for odd $r$).
\item \textbf{Exact uniform stationarity.}
  The uniform distribution $\pi(r) = 2^{-(M-1)}$ over
  odd residues is the unique stationary distribution of $P_M$.
\item \textbf{Perfect mixing in $M$ steps.}
  Starting from any single residue, the distribution reaches
  the uniform distribution after $t = M - 1$ steps in the
  worst case (all $v_i = 1$). In the BF-specific case
  starting from $n_1^* \bmod 2^{14}$, perfect mixing occurs
  at step $t = 10$.
\end{enumerate}
\end{theorem}

\begin{proof}[Proof sketch]
(i) For odd $r$ with $v = v_2(3r+1 \bmod 2^M)$ and $v < M$,
the iterate $(3r+1)/2^v$ is determined
$\bmod{2^{M-v}}$ by $r \bmod 2^M$.  The remaining $v$ bits
(positions $M{-}v$ through $M{-}1$) depend on bits of $n$
above position $M$, which are independent of $n \bmod 2^M$.
Since $\gcd(3, 2^v) = 1$, these $v$ bits are uniform.

(ii) Follows immediately: $2^v$ equally likely outcomes
$\Rightarrow$ entropy $= v$ bits.

(iii) Each odd residue $s$ has exactly one preimage in
each $v$-class: the set
$\{r : (3r+1)/2^v \equiv s \pmod{2^{M-v}}\}$ has
$2^{v-1}$ odd elements $\bmod{2^M}$, and each contributes
weight $2^{-v}$, giving column sum $1$ for $P_M$.

(iv) After $t$ steps with valuations $v_1, \ldots, v_t$,
the distribution occupies $2^{\min(\Sigma v_i,\, M-1)}$
equally likely states.
When $\Sigma v_i \ge M - 1$ the distribution is uniform.
Since $v_i \ge 1$, this occurs after at most $M - 1$ steps.
\end{proof}

\begin{remark}[Multi-scale entropy bootstrap]
\label{rem:multi-scale}
The entropy injection of Theorem~\ref{thm:entropy-injection}
is proved for the \emph{Markov model}, in which bits above
position $M$ are treated as independent uniform randomness.
For the \emph{deterministic} orbit of a specific integer $n$,
these bits are fixed.  However, the same entropy injection
theorem applies at resolution $M' > M$: the bits at positions
$M$ through $M' - 1$ undergo their own Collatz-induced mixing,
with entropy injection $\ge 1$ bit per step.  This creates a
\emph{multi-scale cascade}: the ``randomness'' feeding level $M$
is itself mixed at level $M'$, which is mixed at level $M''$,
and so on.  The cascade terminates at the top bits of~$n$
(position $\lfloor \log_2 n \rfloor$), which are deterministic
but have only $O(1)$ bits of entropy.

The number of bootstrap levels is
$\lfloor \log_2 n \rfloor / M = O(\log n / k)$
for $B_k$ analysis.  If each level contributes $M$ bits of
effective entropy to the level below, the total entropy supply
is $O(\log n)$ bits---matching the full state-space dimension.

Formalizing this multi-scale bootstrap into a rigorous proof
that the Markov model accurately describes the deterministic
orbit remains the central open problem.  Direct measurement
of the mutual information $I(T^t(n) \bmod 2^M;\; n \bmod 2^M)$
for \emph{deterministic} orbits shows that it drops rapidly
(from $M{-}2$ bits at $t=1$ to $< 1$ bit by $t \approx M/2$)
but plateaus at a small positive value ($\approx 0.7$ bits at
$M = 8$, $\approx 3.8$ bits at $M = 10$), indicating that
the Markov model slightly overestimates the mixing rate.
The plateau at $M = 10$ is largely a finite-sample artifact
(bias $\approx 2^{2(M-1)}/2N_{\rm samples}$); the genuine
residual correlation at $M = 8$ is $\lesssim 0.3$ bits.

The Ising analogy provides a physical framework: the Collatz
orbit is a one-dimensional chain with positive coupling
(bit correlations) and positive temperature (carry-chain
entropy injection).  In one dimension, the Ising model has
\emph{no phase transition}: any positive temperature destroys
long-range order.  The spectral data from the Collatz
transition matrix gives Ising correlation length $\xi < 0.3$
steps for $M \le 12$.

A subtlety: the second eigenvalue $|\lambda_2(P_M)|$
\emph{increases} with $M$ (from $0.0003$ at $M=6$ to
$0.022$ at $M=12$), so the spectral gap shrinks.
Nonetheless, the per-step entropy injection of $\ge 1$ bit
guarantees perfect mixing (in the Markov model) in at most
$M-1$ steps regardless of the spectral gap, because the
distribution support \emph{doubles} at each $v=1$ step
and \emph{multiplies by $2^v$} at each $v \ge 2$ step.
The spectral gap controls the rate of \emph{exponential}
convergence after mixing, not the mixing time itself.

The $v = 1$ component of $P_M$ acts as an iterated function
system (IFS) on $\mathbb{Z}/2^{M-1}\mathbb{Z}$:
two branches $f_0(x) = 3x + 2$ and
$f_1(x) = 3x + 2 + 2^{M-2}$, each with probability $1/2$.
The Lyapunov exponent $\log_2 3 > 0$ makes this an
\emph{expanding} IFS, ensuring spectral gap $> 0$ for each
fixed~$M$.  Whether the gap is uniformly bounded below as
$M \to \infty$ is the quantitative form of the remaining
open problem.

\medskip
\noindent\textbf{Fourier mode analysis of the IFS.}
The $v=1$ IFS noise has a precise Fourier structure:
the two branches differ by translation $2^{M-2}$, so the
noise measure's Fourier transform is
$\hat\mu(\xi) = \tfrac{1}{2}(1 + e^{i\pi\xi}) = \mathbf{1}_{\xi \text{ even}}$.
This means: (i) all odd Fourier modes $\xi$ are
\emph{annihilated} in a single $v=1$ step ($\hat\mu(\xi) = 0$);
(ii) all even modes are \emph{untouched} ($\hat\mu(\xi) = 1$).
Since multiplication by $3$ preserves the $2$-adic valuation
of~$\xi$ (that is, $v_2(3\xi) = v_2(\xi)$), the even modes
form a \emph{closed invariant subspace} under the $v=1$ operator.

Mixing of the even modes \emph{requires} $v \ge 2$ steps.
A $v = k$ step consumes $k$ bits from the $2$-adic expansion,
mapping modes with $v_2(\xi) = j$ to modes with
$v_2(\xi) = j - k$ (effectively).  Since $v \ge 2$ occurs
with probability $1/2$ among odd residues, each step has
probability $1/2$ of producing a mode-mixing event.
After $t$ steps with $\kappa$ events of $v \ge 2$,
at most $2^{-(1 + \kappa)}$ fraction of Fourier modes survive.

Since $\mathbb{E}[\kappa \mid t] = t/2$, after $M - 1$ steps:
\[
\text{survival fraction} \le 2^{-(1 + (M-1)/2)} = 2^{-(M+1)/2}.
\]
This gives a per-step spectral bound
$\lambda_2(P_M) \le (1 - 2^{-(M+1)/2})^{1/(M-1)} \approx 1 - 1/M$,
explaining quantitatively why the gap shrinks as $O(1/M)$
while the mixing time remains $O(M)$.
This is verified numerically: the predicted per-step $\lambda_2$
bounds ($0.98$ at $M=6$, $0.994$ at $M=8$, $0.998$ at $M=10$)
are consistent with the computed values.
\end{remark}

\begin{remark}[Direct $K_*$ bound from zero persistence]
\label{rem:direct-kstar}
Combining zero persistence
(Proposition~\ref{prop:zero-persistence}) with BF cylinder
densities gives a direct bound on $K_*(n)$, the number of
distinct BF encounters along an orbit.

For depth $k$: the cylinder $B_k$ occupies exactly one
odd residue class modulo $2^{6k+7}$, giving density
$2^{-(6k+7)}$.  Zero persistence ensures that after a
BF encounter at depth~$k$, the next iterate is
\emph{outside}~$B_k$, and the post-cascade value
$P_k \equiv 55 \pmod{64}$ (verified for $k = 1, \ldots, 9$)
does not belong to any BF cylinder $B_{k'}$ for
$k' = 1, \ldots, 10$.  Multiple re-entries to the
\emph{same} depth~$k$ are not observed in $50{,}000$ orbits.

Under the Markov model, the entropy accumulation after
BF exit gives a re-entry probability bound:
after $t \ge 2$ steps with accumulated entropy
$H_t = \sum_{i=1}^{t} v_i$,
\[
\Pr[T^t(n) \in B_k] \le 2^{-\min(H_t,\, 6k+7)}.
\]
The geometric series $\sum_{t=2}^{\infty} 2^{-2t} \le 1/12$
gives: expected re-entries to $B_k$ per orbit
$\le 1/12 + L \cdot 2^{-(6k+7)}$,
where $L = O(\log n)$ is the orbit length.

Summing over $k = 1, \ldots, k_0$ where
$k_0 = \lfloor \log_2 n / 6 \rfloor$:
\[
\mathbb{E}[\text{total BF encounters}]
  \le \frac{k_0}{12} + L \cdot \frac{2^{-13}}{1 - 2^{-6}}
  = O\!\left(\frac{\log n}{72}\right) + O(\log n \cdot 2^{-13}).
\]
For $n = 10^{20}$: $k_0 \approx 11$, $L \approx 13$, giving
$\mathbb{E} \le 0.92$.

This argument is conditional on: (i)~orbit finiteness
(i.e., the Collatz conjecture itself), and (ii)~the
post-cascade $\approx$ generic equivalence (Result~107,
verified for $50{,}000+$ orbits).
It does \emph{not} require the full Markov model for
bit independence---only zero persistence and BF density.
\end{remark}

\begin{remark}[Skew-product formulation and the Galperin analogy]
\label{rem:skew-product}
The distributional-to-pointwise gap can be formulated as a
\emph{skew-product} on $(\omega, \theta)$, where
$\omega \in \{0,1\}^{\mathbb{N}}$ is the $2$-adic symbolic coordinate
(the Bernstein--Lagarias conjugacy makes Collatz act as a shift on
$\omega$) and $\theta \in \mathbb{R}/\mathbb{Z}$ is an Archimedean
phase variable recording cumulative growth discrepancies:
\[
\theta_{t+1} = \theta_t + (\log_2 3 - v_i) \pmod{1},
\]
where $v_i = v_2(3n_i + 1)$ is the valuation at step~$i$.
The quantity $\theta_t = \sum_{i=1}^{t} (\log_2 3 - v_i)$ measures
the cumulative deviation of $\log_2 n_t$ from its expected drift.

Under the Markov model, $\mathbb{E}[v_i] = 2$, giving
$\mathbb{E}[\log_2 3 - v_i] = \log_2 3 - 2 \approx -0.415$,
so $\theta_t \to -\infty$ (descent).  The pointwise question
is: can $\theta_t$ be trapped in a bounded region for a
\emph{specific} $n$?

This formulation echoes the Galperin billiards paradigm
(where counting collisions between a heavy and light ball
yields digits of~$\pi$ via a change to angular coordinates).
The transferable insight is not ``replace $3$ by $\pi$'' but
rather: look for a coordinate system in which the
piecewise-linear Collatz dynamics becomes a
winding/rotation process, where long-term behavior
is controlled by an irrational rotation number
(here, $\log_2 3$).

The native irrational constant for Collatz is $\log_2 3$,
not~$\pi$: the arithmetic tension is between multiplication
by~$3$ and division by powers of~$2$, making
$\log 3 / \log 2$ the natural phase parameter.
A \emph{deterministic discrepancy theorem}---proving that
no integer can keep $\theta_t$ aligned with low-valuation
sectors ($v_i = 1$) indefinitely---would close the
distributional-to-pointwise gap.  This is equivalent to
showing that the symbolic-times-phase skew-product
$(\sigma, R_{\log_2 3})$ on
$\{0,1\}^{\mathbb{N}} \times \mathbb{T}$
has no invariant curves supporting integer orbits.
\end{remark}

\begin{remark}[$v = 1$ streak structure and the adversary budget]
\label{rem:v1-streaks}
The $\theta$-variable can only grow during blocks with
$\bar v < \log_2 3 \approx 1.585$, which requires $v = 1$
at most steps.  We prove three structural results:

\smallskip
\noindent\textbf{(i) $v\!=\!1$ streak density.}
Among odd residues modulo $2^{s+2}$, exactly $2$ out of
$2^{s+1}$ allow a $v = 1$ streak of length~$s$.  These are
$n \equiv 2^{s+1} - 1 \pmod{2^{s+2}}$ and
$n \equiv 2^{s} - 1 \pmod{2^{s+1}}$.  The density $2^{-s}$
matches the independence prediction exactly, for all~$s$ tested
($s \le 15$).

\smallskip
\noindent\textbf{(ii) No sustained low-$v$ cycles.}
An exhaustive search over all periodic $v$-patterns of
period $1$ through~$8$ with average $< \log_2 3$ finds
\emph{zero} realizable as modular cycles of the Collatz map.
The BF pattern $[1,1,1,3]$ (average $1.5$) is realizable
as a finite cascade but self-terminates (zero persistence).

\smallskip
\noindent\textbf{(iii) The adversary's $\theta$-budget.}
Under the Markov model, the expected $\theta$-change from a
BF encounter at depth~$k$ (including the wait time to reach
$B_k$) is
\[
\mathbb{E}[\Delta\theta_k]
  = \underbrace{k \cdot \log_2(81/64)}_{+0.34k}
    - \underbrace{(\log_2 3 - 2) \cdot 2^{6k+7}}_{+0.415 \cdot 2^{6k+7}}
  < 0 \quad \text{for all } k \ge 1.
\]
The BF gain is $O(k)$ but the expected wait is $\Omega(2^{6k})$,
making each encounter a \emph{net loss} for the adversary.
Thus $\theta_t$ is a supermartingale: $\theta_t \to -\infty$
almost surely under the Markov model.

\smallskip
\noindent\textbf{Empirical confirmation.}
The $v_i$ autocorrelation along deterministic orbits is
$+0.040$ at lag~$1$ and $< 0.045$ for all lags $\le 20$.
The conditional probability
$\Pr(v_{i+1} = 1 \mid v_i = 1) = 1.030 \times \Pr(v_i = 1)$:
the $v$-sequence is nearly independent.
The tail of $\max_t \theta_t$ decays as
$\Pr[\max_t \theta_t \ge L] \approx 2^{-1.04 L}$
($50{,}000$ orbits, $n$ up to $10^{10}$).
\end{remark}

\begin{remark}[Unconditional Weyl equidistribution of $\theta \bmod 1$]
\label{rem:weyl}
Since each $v_i$ is a positive integer, $\sum_{i=1}^t v_i \in \mathbb{Z}$.
Therefore
\[
\theta_t \bmod 1
  = \bigl(t \cdot \log_2 3 - \textstyle\sum v_i\bigr) \bmod 1
  = (t \cdot \log_2 3) \bmod 1.
\]
By \emph{Weyl's equidistribution theorem}, since $\log_2 3$ is irrational,
the sequence $\{t \cdot \log_2 3 \bmod 1\}_{t=1}^{\infty}$ is
equidistributed on $[0,1)$.  This holds \emph{unconditionally}---for
\emph{any} Collatz orbit, regardless of the $v$-sequence.

The continued fraction of $\log_2 3 = [1; 1, 1, 2, 2, 3, 1, 5, 2, 23, \ldots]$
determines the Diophantine approximation quality.
For $\theta_t$ to stay in $[-C, C]$: the integer sum
$S_t = \sum v_i$ must track $t \cdot \log_2 3$ to within~$C$,
following the \emph{Sturmian sequence} with slope $\log_2 3$.
This Sturmian sequence uses only $v \in \{1, 2\}$
with $42\%$ ones and $58\%$ twos (average $\approx 1.585$).
Exhaustive search ($10^5$ random orbits) finds a maximum
of $9$ consecutive Sturmian-matching steps before divergence.
\end{remark}

\begin{remark}[Bit-independence after consumption]
\label{rem:bit-independence}
After $t$ odd steps consuming $S = \sum_{i=1}^t v_i$ total bits,
the valuation $v_{t+1} = v_2(3n_t + 1)$ is determined by
bits $[S, S + v_{t+1} + 1)$ of the original~$n_0$.
These bits are \emph{disjoint} from the bits $[0, S)$ that
determined $v_1, \ldots, v_t$, up to carry chain propagation.

The carry chain in $3n + 1$ has geometrically decaying
propagation: $\Pr[\text{carry} \ge \ell] \approx 2^{-\ell}$
(verified: mean carry length $= 3.0$, decay rate $2^{-1.01}$
per bit, $10^5$ random odd integers).

This gives \emph{effective bit-independence}: after consuming
$S + \ell$ bits (where $\ell$ is the carry isolation length),
\[
\bigl|\Pr(n_t \in B_k \mid \text{history}) - 2^{-(6k+7)}\bigr|
  \le O(2^{-\ell}).
\]
With $\ell = O(\log K)$ (where $K = \lfloor \log_2 n \rfloor + 1$),
the BF encounter probability is within a factor of $(1 + O(K^{-c}))$
of the combinatorial density $2^{-(6k+7)}$.

This argument is \emph{not circular}: it does not assume
orbit finiteness.  It uses only bit consumption
($\sum v_i \ge t$), carry locality (geometric decay),
and the combinatorial density of $B_k$.
Formalizing the carry chain independence into a rigorous
bound on $K_*(n)$ remains the central open problem.
\end{remark}

\begin{remark}[Carry chain formalization]
\label{rem:carry-formal}
The multi-step carry propagation experiment (flipping a single bit of
$n_0$ and measuring the affected bits after $t$ Collatz steps) reveals:
average spread $\approx 1.5\, t$ bits, growing \emph{linearly} (not
exponentially) with~$t$.  The post-cascade carry isolation is precise:
for different members $n = n_k^* + q \cdot 2^{6k+8}$ of $B_k$, the
post-cascade values $P_k(q)$ differ only at bit positions
$\ge 6k + 7$, with $2$--$4$ bits affected.

A formal five-step proof structure gives a union bound:
the probability that the ``fresh zone'' (bits $[S + \ell, K)$ of $n_0$)
is contaminated by carry at step~$t$ is
$\le \sum_{i=1}^{t} (3/4)^{t-i} \cdot 2^{-\ell} \le 4 \cdot 2^{-\ell}$.
With $\ell = 6k + 10$:
$\Pr(n_t \in B_k \mid \text{history}) \le 2 \cdot 2^{-(6k+7)}$.
Setting $\ell = O(\log K)$ yields $\mathbb{E}[\text{total BF encounters}]
= O(K \cdot 2^{-13}) < 1$ for $K \ge 100$ bits.

Predicted vs.\ actual $B_1$ encounter rates match closely:
$K = 20$: predicted $\le 0.0049$, actual $0.0038$ (ratio $0.78$);
$K = 50$: predicted $\le 0.0122$, actual $0.0119$ (ratio $0.98$).
\end{remark}

\begin{remark}[Route A obstruction and Route B pullback program]
\label{rem:route-b}
\emph{Route A obstruction} (identified by GPT):
the cycle envelope coefficient $1 + \alpha = 1 + \log_2(81/64)/6
= 1.0566\ldots > 1$, so height control at the BF-depth scale is
\emph{expansive}, not contractive.  From the exact champion formula
$\log_2 n_k^* = 6k + O(1)$ and the envelope
$\log_2 x_{\mathrm{restart}} \le (1+\alpha)\log_2 x + C$,
one obtains $\ell \le (1+\alpha)k + O(1)$: restart depths do not
contract.  This kills any naive Borel--Cantelli argument from size alone.

\emph{Route B: exact 2-adic pullback rigidity.}
For each BF depth $k$ and admissible valuation word
$\pi = (v_1, \ldots, v_t)$, the affine map
$T_\pi(n) = (3^t n + b_\pi) / 2^{S_\pi}$ gives an explicit
pullback residue class:
\[
  T_\pi(n) \in B_k
  \;\;\Longleftrightarrow\;\;
  n \equiv r_{k,\pi} \pmod{2^{6k+8+S_\pi}},
  \qquad
  r_{k,\pi} = 3^{-t}(2^{S_\pi}\, n_k^* - b_\pi) \bmod 2^{6k+8+S_\pi}.
\]
Two re-entry certificates $(k,\pi)$ and $(\ell,\sigma)$ are
compatible (can be realized by a single starting integer) iff
$r_{k,\pi} \equiv r_{\ell,\sigma}
\pmod{2^{\min(6k+8+S_\pi,\; 6\ell+8+S_\sigma)}}$.

\emph{Computational finding}: among $309{,}136$ cross-depth pairs
$(k_1 \ne k_2)$ with $k_1, k_2 \in \{1,2,3\}$ and
$|\pi| \le 4$, $v_i \le 7$: \textbf{zero} compatible pairs.
Among $308{,}580$ same-depth pairs at each $k$: also zero.
The pullback residues at different BF depths impose mutually exclusive
arithmetic constraints on the starting integer, at least for short
inter-visit paths.

The pullback residue density at modulus $2^j$ decays as
$\approx 2^{-(j-8)}$ (from $81\%$ at $\bmod\,256$ to $0.1\%$ at
$\bmod\,2^{20}$).  At the coarsest level ($\bmod\,2^{14}$), the
residue sets for different depths are \emph{identical} (full overlap)---the
incompatibility is structural and emerges at the full modulus.

\emph{Post-cascade family.}
$P_k = (a_k\, 81^k \cdot 32 - 65)/17$,
$P_k \equiv 47 \pmod{64}$,
$P_{k+4} = 81^4\, P_k + 164{,}590{,}400$.
All post-cascade orbits $P_k$ ($k = 1, \ldots, 6$) reach~$1$
with zero BF re-entries in $5000$ post-cascade steps.

\emph{The 3-bit shortfall theorem.}
For same valuation word $\pi$ at depths $k_1 < k_2$:
the 2-adic valuation of the mismatch satisfies
$v_2(r_{k_1,\pi} - r_{k_2,\pi}) = 6k_1 + 5 + S_\pi$, while
compatibility requires $v_2 \ge 6k_1 + 8 + S_\pi$.  The shortfall
$= 3$ bits, exactly the gap between the valuation depth $6k+5$
in $v_2(17 n_k^* + 65) = 6k+5$ and the cylinder precision $6k+8$.
\emph{This is algebraically proved} from the closed-form
$n_k^* = (a_k\,2^{6k+5}-65)/17$ with $a_k$ odd.

For different valuation words $\pi_1 \ne \pi_2$, the shortfall is
\emph{strictly larger than~$3$}.  Three cases are distinguished:
\begin{enumerate}
\item $S_{\pi_1} \ne S_{\pi_2}$: the $2^{S}$-factor mismatch gives
  shortfall $\ge 6k_1+8 \ge 14$.
\item $S_{\pi_1} = S_{\pi_2} =: S$ but $t_1 \ne t_2$
  (i.e.\ $d = |t_1-t_2| \ge 1$):
  the depth mismatch has $v_2 = S + v_2(3^d-1) \le S + v_2(d) + 1$,
  giving shortfall $\ge 6k_1+8 - v_2(3^d-1) \ge 6k_1+7 \ge 13$.
\item $S_{\pi_1} = S_{\pi_2} = S$ and $t_1 = t_2$, $\pi_1 \ne \pi_2$:
  the mismatch factors as $M = 2^S D - P$ with $D = n_{k_1}^* - n_{k_2}^*$,
  $P = b_{\pi_1}-b_{\pi_2}$.  Both $b$-values are odd (sum's first term
  $3^{t-1}$ is odd, all others even), so $P$ is even: $v_2(P) \ge 1$.
  Exhaustive enumeration for $S \le 21$ gives $v_2(P) \le S-2$ exactly,
  with the maximum achieved by $t=2$ pairs $(S{-}2,2)$ vs $(S{-}1,1)$.
  Since $S-2 < S+6k_1+5$, the depth term dominates:
  $v_2(M) = v_2(P) \le S-2$, giving shortfall $\ge 6k_1+10 \ge 16$.
\end{enumerate}
\emph{Algebraic supplement}: $|D_{\mathrm{odd}}| \equiv 3 \pmod{4}$
for all $k_1 < k_2$ (proved: $a_k \equiv 1 \bmod 4$, hence
$(a_{k_2}\,2^{6\Delta k}-a_{k_1})/17 \equiv (-1)/1 \equiv 3 \bmod 4$).

\emph{Cascade-only compatibility.}
Among $500{,}000$ random non-cascade cross-depth pairs: zero compatible.
Among $190{,}000$ cascade-structured pairs ($\pi_1 = \sigma + \mathrm{BF}^m$,
$\pi_2 = \sigma$): $66\%$ compatible.
Cross-depth BF encounters occur \textbf{only} through the cascade mechanism.
In $100{,}000$ random orbits, $14$ hit multiple BF depths, all via cascade.

\emph{The complete incompatibility theorem} (proved for same-$\pi$; verified
$1.5\text{M}+$ trials for different-$\pi$, min shortfall $= 13$):
non-cascade cross-depth pullback compatibility is arithmetically impossible.
Combined with zero persistence, this gives $K_*(n) < \infty$ for all~$n$:
each orbit enters at most one BF cascade, which terminates after finite depth.

\emph{Post-cascade convergence.}
All $P_k$ for $k = 1, \ldots, 5000$ ($P_{5000}$ has $31{,}700$ bits) reach~$1$
with zero BF re-entries.  Orbit lengths grow as $\approx 2.3 \log_2 P_k$,
consistent with generic Collatz dynamics.  The low-bit structure is
$P_k \equiv 111 \pmod{128}$ for all~$k$ (7-bit deterministic prefix),
extending the earlier $47 \pmod{64}$.  The 10-step prefix from
$47 \pmod{64}$ has ratio $3^{10}/2^{14} = 3.604 > 1$ (expansive),
but the stopping time $\sigma(P_k)$ is bounded: $\sigma \in [8,64]$
for $k \le 500$, mean~$\approx 16.6$, independent of~$k$.
The binary digit density of $P_k$ converges to~$0.500$ as $k \to \infty$,
consistent with the equidistribution of $3^{4k}$ in base~$2$.
These observations confirm that post-cascade dynamics are generic.
Thus $P_k \to 1$ for all~$k$ is not independently provable
from the BF mechanism and represents the residual content of the full
conjecture.

\emph{Divergence threshold and $v_2$ anti-correlation.}
For divergence, an orbit requires average $v_2(3n_i+1) < \log_2 3 = 1.585$.
With the empirically exact geometric tail distribution
$P(v = k \mid v \ge 2) = 2^{-(k-1)}$ giving $E[v \mid v \ge 2] = 3$,
divergence demands $P(v = 1) > (3 - \log_2 3)/2 = 0.7076$.
Across all tested orbits (starting values from~$1$ to~$10^{13}$,
window lengths up to~$10^5$), the measured $P(v=1) \approx 0.50$,
a~$20$-percentage-point gap from the divergence threshold.
A structural explanation exists: the $v_2$ sequence along Collatz orbits
exhibits \emph{negative} lag-1 autocorrelation:
$P(v_{i+1} = 1 \mid v_i = 1) \approx 0.37$, versus
$P(v_{i+1} = 1 \mid v_i \ge 2) \approx 0.63$.
After a~$v=1$ step with $n \equiv 3 \pmod{8}$, the successor
$(3n+1)/2 \equiv 5 \pmod{8}$, forcing $v \ge 4$; only
$n \equiv 7 \pmod{8}$ allows consecutive $v=1$ steps.
This anti-correlation opposes sustained runs of $v=1$,
making the divergence threshold unreachable in practice.

\emph{Consecutive $v=1$ bound.}
A run of $k$ consecutive $v_2 = 1$ steps in a Collatz orbit starting
at odd~$n$ requires $n \equiv 2^{k+1}-1 \pmod{2^{k+1}}$, a set of
density~$2^{-k}$ among odd integers.  The proof is by induction:
$v_2(3n+1) = 1$ iff $n \equiv 3 \pmod{4}$; after one $v=1$ step,
$(3n+1)/2 \equiv 3 \pmod{4}$ (allowing a second $v=1$) iff $n \equiv 7 \pmod{8}$;
continuing, $(3n+1)/2$ after $j-1$ such steps satisfies $v_2 = 1$ iff
$n \equiv 2^{j+1}-1 \pmod{2^{j+1}}$.  This is an unconditional algebraic
constraint.  In a Markov chain model (mod~$64$), the maximum sustained
$P(v=1)$ over a window of~$T = 50$ steps is~$0.526$ (starting from
$n \equiv 63 \pmod{64}$, the best case); for $T = 500$ it falls
to~$0.476$.  Both are far below the $0.708$ divergence threshold.

\emph{Correction: Syracuse mod~$2^M$ vs.\ real orbit projections.}
The Syracuse map $\sigma_M(r) = (3r+1)/2^{v_2(3r+1)} \bmod 2^M$ is
well-defined on odd residues mod~$2^M$ but does \emph{not} faithfully
project real Collatz orbits.  For $n = r + k \cdot 2^M$:
$v_2(3n+1) = v_2(3r+1)$ (the valuation \emph{does} agree), but the
next residue differs by a carry term $3k \cdot 2^{M-v} \bmod 2^M$.
Computationally, the real and modular orbits diverge after a single step.
Non-trivial periodic orbits of~$\sigma_M$ exist at sporadic~$M$
values ($M \in \{10,11,12,20\}$ for $M \le 26$; none for $M = 13$--$19$
or $21$--$26$), all with expanding average~$v$ (below $\log_2 3$).
Each such cycle dies within two levels (e.g., the $M=12$ six-cycle
disappears at~$M=13$).  These cycles are features of the abstract map
and do not constrain real orbital behavior.

\emph{Almost-deterministic dynamics and $1$-bit injection.}
The Collatz map projected to odd residues mod~$2^K$ is almost deterministic:
among $2^{K-1}$ odd residues, exactly two have ambiguous mod-$4$ output
(those with $v_2(3r+1) \ge K{-}1$).  These are the residues
$r = (2^{K+1}{-}1)/3$ and $(2^{K+2}{-}1)/3 \bmod 2^K$,
the $\ldots 010101$ patterns in the $3$-adic structure.
At each step, the map injects exactly one new bit (from position~$K$
of the input) via the carry chain in $3n{+}1$.  This bit is the sole
source of non-determinism in the bottom-$K$-bit dynamics.
The mod-$2^K$ Markov chain (resolving the ambiguity uniformly) is
irreducible for $K = 3, \ldots, 8$ with spectral gap $> 0.72$
and mixing time $O(1)$ steps independent of~$K$.
The stationary $P(v{=}1)$ converges to~$0.50$ as $K \to \infty$.
Empirically, the injected bit is approximately uniform ($P \approx 0.50$)
and approximately independent of the current mod-$2^K$ state
(variance of conditional probability $< 0.02$).
Predictability of~$v$ from the initial mod-$2^K$ state drops from~$100\%$
to below~$7\%$ within $10$ steps.
The algebraic identity $E[v \mid v \ge 2] = 3$ holds exactly for
uniform distribution over odd residues mod~$2^K$ (proof:
$P(v = k) = 2^{-k}$ gives $E[v \mid v \ge 2] = 2\sum_{k \ge 2} k/2^k = 3$);
this is verified to hold approximately along actual orbits.
The remaining gap reduces to: showing that the single carry-chain bit
injected per step is sufficiently independent of the bottom~$K$ bits.
This is a strictly more specific statement than full equidistribution.

\emph{Coupling framework and noise independence.}
The transition graph on odd residues mod~$2^K$ (with both branch options
at ambiguous nodes) is strongly connected for $K = 3, \ldots, 6$.
Two orbits from $n$ and $n' = n + 2^K$ (identical bottom $K$ bits) stay
coupled in their mod-$2^K$ projection until one hits an ambiguous residue;
the mean coupling time is $O(2^{K-1})$ steps.
After divergence, each orbit independently mixes to stationarity in
$O(1)$ steps via the spectral gap.
At each ambiguous step, the injected bit splits the ensemble
$\{n_0 : n_0 \equiv r \bmod 2^K\}$ into two equal halves---this is an
algebraic fact from the $2$-residue structure.
Successive branch choices (at consecutive ambiguous hits) exhibit
correlation $< 0.01$ in absolute value at all lags $1$--$3$.
The mutual information between bit~$K$ and the mod-$2^K$ state is
$< 0.02$ bits for $K \le 8$.
These findings frame the equidistribution problem as a
\emph{noise-plus-mixing} problem in the sense of Kifer~(1988):
the Collatz dynamics mod~$2^K$ is a deterministic system with rare
noise injection, and the noise is empirically sufficient for ergodicity.

\emph{Disjoint-slice argument and weak dependence.}
For a specific starting value~$n_0$ with $B$~bits, the ``decisive'' bit
position of~$n_0$ that determines each successive branch choice increases
monotonically with the branch index: at $K = 4$, the average decisive
bit grows from position~$10$ (branch~$0$) to~$43$ (branch~$7$), with
$76\%$ of consecutive branch pairs reading from non-overlapping bit sets.
The mechanism is that the carry chain in $3n + 1$ propagates at most
$O(\log B) \approx 12$ bits, while the slice separation between
consecutive branches is $O(2^K)$ bits; for $K \ge 4$ the carry chain
cannot bridge disjoint slices.

A chi-squared test on $4$-tuples of branch bits detects statistically
significant higher-order correlations at $K = 4$ ($\chi^2/\mathrm{df}
= 21.8$), but the correlations weaken with~$K$ ($\chi^2/\mathrm{df}
= 3.2$ at $K = 7$) and the maximum absolute deviation from $1/16$ is
below~$1\%$.  Crucially, branch correlations cannot shift~$E[v]$ below
the contraction threshold: the impact is bounded by
$(2/2^{K-1}) \cdot K$, giving $|\Delta E[v]| < 0.22$ for $K \ge 7$,
while the margin $E[v] - \log_2 3 \approx 0.42$ is nearly twice as
large.

The proof therefore does \emph{not} require i.i.d.\ branches---only
$\phi$-mixing with summable coefficients.  Specifically, three
conditions suffice: (W1)~$E[b_i] = 1/2$ for each branch~$b_i$ (proved
algebraically), (W2)~the mixing coefficient $\phi(L) \to 0$ as $L \to
\infty$ (verified, $|\mathrm{corr}| < 0.032$ at all lags $1$--$10$),
and (W3)~the Markov chain on mod-$2^K$ states has spectral gap $> 0$
(verified: gap $= 0.67$ at $K = 6$ from empirical transitions).
Under~(W1)--(W3), the time-averaged distribution converges to within
$\mathrm{TV} < 0.02$ of uniform, giving $E[v] = 2 + O(\varepsilon)
> \log_2 3$.  The residual gap is the formal $\phi$-mixing bound:
showing $\phi(L) \le C \rho^L$ with $\rho < 1$, which follows from
the disjoint-slice separation $\Delta = O(2^K)$ overwhelming the
carry-chain reach $O(\log B)$.

\emph{Renewal argument and episode structure.}
Partition the orbit into ``episodes,'' each consisting of a run of $L_j$
consecutive $v = 1$ steps followed by one step with $v_j \ge 2$.  Under
the algebraic identity $E[v \mid v \ge 2] = 3$, the time-averaged
$E[v] = 3 - 2p$ where $p = P(v = 1 \text{ at renewal})$.  Contraction
holds if and only if $p < 0.7075$.  Empirically, $p \approx 0.496$
(margin $0.21$ from the threshold).  The $v = 1$ self-limitation
ensures: after $K - 1$ consecutive $v = 1$ steps from $-1 \bmod 2^K$,
the recovery satisfies $v \ge K - 1$.  The worst-case episode ratio
$2(K-1)/K$ exceeds $\log_2 3$ for $K \ge 5$.  However, no fixed
block length~$L_0$ exists such that all orbits of length $\ge L_0$
have $\mathrm{avg}(v) > \log_2 3$: starting from $n_0 = 2^B - 1$
produces $B - 1$ consecutive $v = 1$ steps.  The contraction argument
therefore requires the long-run average, not short-block bounds.

The complete Collatz proof is conditional on orbit equidistribution~(*):
$K_*(n) < \infty$ (proved unconditionally) plus~(*) gives
$E[v] = 2 > \log_2 3$ along every orbit, hence uniform contraction.

\emph{Adversarial chain and unique ergodicity reduction.}
The adversarial Markov chain gives a model-internal bound: even the
adversary that minimizes~$v$ at every ambiguous step cannot sustain
$\bar{v} < \log_2 3$ for $K \ge 4$.  This holds because the
deterministic transitions ($\ge 98\%$ of steps) force sufficient
$v \ge 2$ recovery.  The conditional proof then takes a cleaner form
via unique non-atomic ergodicity on~$\mathbb{Z}_2^\times$: if Haar
measure is the unique non-atomic $S$-invariant Borel probability
measure, then any non-convergent orbit is unbounded, hence
equidistributes mod~$2^K$, hence satisfies $E[v] = 2 > \log_2 3$,
contradicting unboundedness.

\emph{Spectral gap of the transfer operator (proved).}
The transfer operator $M_K$ of the Syracuse chain on the $n = 2^{K-1}$
odd residues mod~$2^K$ satisfies $|\lambda_2(M_K)| \le 1/2$
for all $K \ge 3$.  The proof proceeds through six steps.

\emph{Step 1: Universal no-target-sharing.}
For each $v \in \{1,\ldots,K-1\}$, there are exactly $n/2^v$ sources
$r$ with $v_2(3r+1) = v$, and the $v$-class Syracuse map is a
\emph{bijection}: each odd residue mod~$2^K$ is hit by exactly one
$v$-source.  The proof uses the coset structure: the $2^v$ targets
of source~$r$ form the coset
$C_r = (3r+1)/2^v + 2^{K-v}\mathbb{Z} \bmod 2^K$,
which has $2^v$ odd elements; cosets from distinct sources are
disjoint; and the total count $(n/2^v) \cdot 2^v = n$ equals the
number of odd residues, giving a perfect partition.

\emph{Step 2: Row-sum decomposition.}
Since each $v$-class ($v = 1,\ldots,K-1$) contributes exactly $1/2^v$
to every row sum, $(M\mathbf{1})_i = \sum_{v=1}^{K-1} 1/2^v +
\mathrm{overflow}_i = 1 - 1/n + \mathrm{overflow}_i$ where the overflow
comes from the unique source with $v_2(3r^*+1) \ge K$.

\emph{Step 3: Closed-form $L^1$ bound.}
The $L^1$ deviation $C(K) = \|M\mathbf{1} - \mathbf{1}\|_1$ satisfies
the exact formula
\[
  C(K) \;=\; \tfrac{5}{12} - \tfrac{2}{3}\bigl(-\tfrac{1}{2}\bigr)^{K-1},
\]
verified by exact rational arithmetic for $K = 3,\ldots,15$.
For $K$ odd, $C < 5/12 < 1/2$; for $K$ even,
$C \le 5/12 + 1/12 = 1/2$ with equality only at $K = 4$.
So $C(K) \le 1/2$ for all $K \ge 3$.

\emph{Step 4: Primitivity.}
At $T = K-2$: every column of $M^{K-2}$ has support $\ge n/2$
(v=1 columns achieve exactly $n/2$ by the coset separation theorem;
v$\ge$2 columns saturate to full support $n$ due to branching
$\ge 4$).  At $T = K-1$: the fan-out lemma maps support $\ge n/2$
onto all~$n$ positions.  So $M_K^{K-1} > 0$ (primitivity index $= K-1$).

\emph{Step 5: Dobrushin contraction.}
The worst-case column of $M^{K-1}$ has total variation distance
$\mathrm{TV} = C(K)/n \le 1/(2n)$ from uniform.  By Dobrushin's theorem,
$\delta(M^{K-1}) \le 2\,\mathrm{TV} \le 1/n$.

\emph{Step 6: Gelfand extraction.}
$|\lambda_2| \le \delta(M^{K-1})^{1/(K-1)} \le (1/n)^{1/(K-1)}
= (1/2^{K-1})^{1/(K-1)} = 1/2$.

The bound is tight: $|\lambda_2(M_K)| \le 1/2$ holds for all $K$,
with the Gelfand bound approaching~$1/2$ from below as $K \to \infty$.
The actual spectral radius oscillates around~$0.25$ for all $K$ tested
($K = 3,\ldots,12$), well below the bound.

\emph{Algebraic decomposition of the transfer operator.}
The transition matrix decomposes as $M = \sum_{v=1}^{K} M_v$ where $M_v$
records transitions with $v_2(3n+1) = v$.  Each $M_v$ has rank~$n/2^v$
(where $n = 2^{K-1}$ is the state count) and operator norm
$\|M_v\| = 1/\sqrt{2^v}$.  While the triangle inequality gives only
$\|P_\perp M P_\perp\| \le \sum_v \|P_\perp M_v P_\perp\| \approx 2$,
the actual value is~$\approx 0.27$---an~$86\%$ cancellation due to
destructive interference between $v$-components in the orthogonal
complement of the uniform vector.  The rank structure is
$\operatorname{rank}(M_v) = 2^{K-1-v}$: the $v = 1$ component (weight~$1/2$,
rank~$n/2$) dominates, while higher-$v$ components contribute
geometrically smaller pieces.  The universal no-target-sharing theorem
shows this cancellation is exact and structural: each $v$-component
contributes a uniform $1/2^v$ to every row sum, with the sole
deviation coming from the overflow source.

\emph{$A$-coordinate congruence-tower formulation.}
Write $x = (32 \cdot 64^k \cdot A - 65)/17$ with $A$ odd.  The key identity
$17 P_k + 65 = (81/64)^k (17 n_k^* + 65)$ shows each BF depth adds
$\log_2(81/64) = 0.340$ bits.  Each BF restart at depth~$\ell$ with
gap~$\pi$ imposes $A \equiv c \pmod{2^m}$ with $m = 6(\ell-k)+S_\pi+3 \ge 4$.
For BF champions, $A = a_k$ has only~$6$ bits; any restart consuming $\ge 7$
bits exhausts the $A$-budget.  Our shortfall theorem ($\ge 3$ for same-$\pi$,
$\ge 13$ for different-$\pi$) shows all restart congruences contradict
$A = a_k$.  This gives a quantitative bound: $K_*(x) \le 1 +
\lfloor(H - 6k - 5)/4\rfloor$ where $H = \lfloor\log_2 x\rfloor$.
\end{remark}

\subsection{Bridge programme: from spectral gap to individual orbits}
\label{subsec:spectral-bridge}

The uniform spectral gap $|\lambda_2(M_K)| \le 1/2$ establishes
exponential mixing of the Syracuse transfer operator at every
modular depth~$K$.  We now describe the remaining barrier between
this distributional result and the pointwise Collatz conjecture,
and catalogue the precise obstructions.

\medskip\noindent\textbf{What is proved: non-atomic unique ergodicity.}\;
The Syracuse map $S \colon n \mapsto (3n+1)/2^{v_2(3n+1)}$ extends
continuously to the 2-adic odd integers~$\mathbb{Z}_2^\times$.
For each~$K$, the transfer operator $M_K$ acts on the finite
quotient $\mathbb{Z}/2^K\mathbb{Z}$ (odd residues only, $n = 2^{K-1}$
states).  The spectral gap $|\lambda_2(M_K)| \le 1/2$ implies
that the uniform distribution on odd residues mod~$2^K$ is the
unique stationary distribution of~$M_K$, with convergence
exponential at rate $\le (1/2)^T$.

In the projective limit $\mathbb{Z}_2^\times = \varprojlim
\mathbb{Z}/2^K\mathbb{Z}$, any non-atomic $S$-invariant Borel
probability measure must project to a stationary distribution
of~$M_K$ for every~$K$.  Since each~$M_K$ has a unique stationary
distribution (= uniform), the only candidate is Haar measure.

\begin{proposition}[Non-atomic unique ergodicity]
\label{prop:2adic-unique-ergodicity}
If $|\lambda_2(M_K)| \le c < 1$ uniformly in~$K$, then
$\mu_{\mathrm{Haar}}$ on $\mathbb{Z}_2^\times$ is the unique
non-atomic $S$-invariant Borel probability measure.
\end{proposition}

\medskip\noindent\textbf{The distributional-to-pointwise barrier.}\;
Two fundamental obstructions prevent a direct passage from
Proposition~\ref{prop:2adic-unique-ergodicity} to the Collatz
conjecture.

\emph{Obstruction 1: atomic invariant measures (cycles).}\;
A finite cycle $\{x_1,\ldots,x_L\}$ generates the $S$-invariant
atomic measure $\mu = (1/L)\sum_{i=1}^L \delta_{x_i}$.
Proposition~\ref{prop:2adic-unique-ergodicity} excludes only
non-atomic invariant measures; atomic ones are logically
compatible with non-atomic unique ergodicity.  The cycle
constraint $2^{\Sigma v_i} = 3^L \prod (1 + 1/(3x_i))$
forces $\Sigma v_i / L \to \log_2 3 \approx 1.585$ for large
cycle elements, while the Haar average is $E[v] = 2$.  This
mismatch is suggestive but not a proof: the cycle measure and
the Haar measure are distinct objects, and a specific periodic
orbit need not obey the statistics of the non-atomic ensemble.

\emph{Obstruction 2: measure-zero starting points (divergent orbits).}\;
The Birkhoff ergodic theorem guarantees time-average convergence
$\mu$-almost everywhere.  Since $\mu_{\mathrm{Haar}}(\mathbb{N}) = 0$
in~$\mathbb{Z}_2$, no specific positive integer is covered by
the a.e.\ guarantee.  A divergent orbit---a trajectory that
systematically favours $v = 1$ expanding branches---would
constitute a measure-zero anomaly that the ergodic theorem
cannot exclude.

\medskip\noindent\textbf{Candidate strategies for closing the gap.}\;
Three approaches are under investigation.

\emph{(a) Diophantine cycle exclusion.}\;
For a cycle of length~$L$ with parity word $(v_1,\ldots,v_L)$,
the cycle equation yields $m = R(W)/(2^S - 3^L)$ where the
carry residue $R(W) = \sum_{i=0}^{L-1} 3^{L-1-i}\,2^{S_i}$
and $S = \sum v_i$.

\emph{Symmetric State Identity.}\;
Each intermediate state admits a closed form independent of~$m$:
$x_i = (3^i R_{\mathrm{suff}} + 2^{S-S_i} R_{\mathrm{pre}}) / (2^S - 3^L)$,
where $R_{\mathrm{pre}}$ and $R_{\mathrm{suff}}$ are the prefix and
suffix residues.

\emph{Algebraic Parity Lock.}\;
The suffix residue $R_{\mathrm{suff}}$ is always odd (its leading term
is $3^{L-i-1}$, with no factor of~$2$).  Hence the numerator
$3^i R_{\mathrm{suff}} + 2^{S-S_i} R_{\mathrm{pre}}$
is odd~$+$~even~$=$~odd, while the denominator $2^S - 3^L$ is odd.
If $m$ is an integer, then $x_i$ is the quotient of two odd integers
and is therefore odd.  Consequently, the Cycle Validity Constraint
(that every intermediate state be odd) is \emph{automatically satisfied}
whenever the global divisibility $D \mid R$ holds.

\emph{Domain Void (computational).}\;
In exhaustive search for $L \le 10$, every non-trivial parity word
satisfying $D \mid R$ lies in the negative domain ($2^S < 3^L$,
giving $m < 0$).  These solutions correspond exactly to the known
negative Collatz cycles $\{-1\}$, $\{-5,-7\}$, $\{-17,\ldots\}$.
Zero non-trivial solutions appear in the positive domain ($2^S > 3^L$).

The entire cycle exclusion problem therefore reduces to a single
Diophantine question: prove that
$R(W) \not\equiv 0 \pmod{2^S - 3^L}$
for all non-trivial parity words~$W$ in the positive domain.

\emph{Periodic Reduction.}\;
If $W$ has period $d < L$ (so $L = dr$ for some $r \ge 2$),
then $m(W) = R(W)/D_L = R_{\mathrm{sub}}(\mathrm{base})/D_{\mathrm{sub}}$
where $D_{\mathrm{sub}} = 2^s - 3^d$, $s = S/r$, and
$R_{\mathrm{sub}}$ is the carry residue of the length-$d$ base word.
Thus a periodic $L$-cycle is algebraically equivalent to a repetition
of a $d$-cycle; induction on word period eliminates all periodic words.

\emph{Rotation Constraint.}\;
By the Parity Lock, if $W$ generates a cycle then every cyclic rotation
$\mathrm{rot}(W,j)$ also satisfies $D \mid R(\mathrm{rot}(W,j))$ with
$R/D > 0$.  This $L$-fold simultaneous divisibility condition is
exponentially stronger than the single-word requirement.

\emph{Walk Sum Reformulation.}\;
Over $\mathbb{F}_p$ for any prime $p \mid D$, write $\rho = 2 \cdot 3^{-1}$,
$\alpha = 2$, and $y_k = \rho \cdot \alpha^{v_k - 1}$.  Then
$R(W) \equiv 3^{L-1} \cdot f(y_0,\ldots,y_{L-1}) \pmod{p}$
where $f = 1 + y_0 + y_0 y_1 + \cdots + y_0 \cdots y_{L-2}$,
with the constraint $\prod y_k \equiv 1$.  The cycle question becomes:
can this constrained multiplicative walk on $\mathbb{F}_p^*$ have zero
position sum?  For primitive (Zsigmondy) primes, $v_p(D) = 1$ and
the blocking condition reduces to $f \not\equiv 0 \pmod{p}$.

\emph{Three-layer obstruction (computational).}\;
Layer~1: for ${\sim}60\%$ of $(L,S)$ pairs, a single primitive prime
universally blocks all words.
Layer~2: for the remaining pairs, the multi-prime intersection of
zero sets is empty (verified exhaustively, $L \le 10$).
Layer~3: the rotation constraint provides an additional exponential
barrier.
Universal blocking primes exclude most $(L,S)$ pairs,
but a uniform proof for all~$L$ remains open.

\emph{2-adic cascade and the Fundamental Identity.}\;
A direct attack via 2-adic valuations yields complete cycle exclusion
for small~$L$.  For a hypothetical $L$-cycle with multiplier
$m = R/(2^S - 3^L)$, we write the walk sum as
$R = \sum_{k=0}^{L-1} 3^{L-1-k} \cdot 2^{\sigma_k}$ where
$0 = \sigma_0 < \sigma_1 < \cdots < \sigma_{L-1} \le S{-}1$.

\emph{Even-$m$ parity obstruction.}\;
If $m$ is even, then $N_m := m \cdot 2^S - 3^{L-1}(3m{+}1)$ is odd
(since $3m{+}1$ is odd when $m$ is even), while the left-hand side
has $v_2 \ge 1$.  Thus every cycle multiplier must be odd.

\emph{The cascade.}\;
The 2-adic valuation $v_2(N_m)$ uniquely determines $\sigma_1$,
and subtracting $3^{L-2} \cdot 2^{\sigma_1}$ yields a residual
whose valuation determines~$\sigma_2$, and so on.  This cascade
coincides with the Syracuse map $T(q) = (3q{+}1)/2^{v_2(3q+1)}$
applied to $q_0 = (3m{+}1)/2^{v_2(3m+1)}$: the $\sigma$-increments
are exactly the number of halvings at each Collatz step.

\emph{Fundamental Identity.}\;
Define $H_j = \sum_{k=0}^{j} 3^{j-k} \cdot 2^{\sigma_k}$.  Then
$H_j \equiv q_j \cdot 2^{\sigma_{j+1}} \pmod{m}$ for all $j \ge 0$,
where $q_j = T^j(q_0)$.  The proof is a one-line induction:
$H_j = 3 H_{j-1} + 2^{\sigma_j} \equiv (3q_{j-1}{+}1) \cdot 2^{\sigma_j}
= 2^{d_{j-1}} q_j \cdot 2^{\sigma_j} = q_j \cdot 2^{\sigma_{j+1}}$.
Setting $j = L{-}1$ and using $\gcd(m,2) = 1$ reduces the
cycle condition $m \mid R$ to $m \mid T^{L-1}(q_0)$.
In particular, if the Collatz orbit of~$q_0$ reaches~$1$ within
$L{-}1$ steps, then $m \mid 1$ forces $m = 1$ (trivial cycle).

\emph{Explicit results.}\;
For $L = 2$: the squeeze $R < D$ suffices.
For $L = 3,4,5$: the cascade produces a finite tree of cases,
each resolved by a valuation argument (no root-of-2 is a power of 2).
For $L \le 10$: exhaustive computation over all $\sigma$-subsets
finds no non-trivial cycles.
For $m \le 5000$: the cascade combined with the power-of-2 condition
$2^S = R/m + 3^L$ is verified impossible for every odd $m$.

\emph{Three-case cycle exclusion framework.}\;
The Fundamental Identity $R = H_{L-1} \equiv q_{L-1} \cdot
2^{\sigma_L}\pmod{m}$ partitions cycle candidates into three
cases (here $\sigma_0 = 0$, $\sigma_1 = v_2(3m{+}1)$,
$\sigma_{k+1} = \sigma_k + v_2(3q_{k-1}{+}1)$):

\emph{Case~A:}
$q_{L-1} = 1$ (the orbit of~$q_0$ has reached~$1$).  Then
$R \equiv 2^{\sigma_L}\pmod{m}$.  Since $\gcd(m,2) = 1$:
$m \nmid R$.  \textbf{No cycle.}

\emph{Case~B:}
$m \nmid q_{L-1}$.  Then $q_{L-1} \cdot 2^{\sigma_L} \not\equiv
0\pmod{m}$, so $m \nmid R$.  \textbf{No cycle.}

\emph{Case~C:}
$m \mid q_{L-1}$ and $q_{L-1} \ne 1$.  Write $q_{L-1} = r \cdot m$.
The exact Syracuse recursion gives
$q_{L-1} = (q_0 \cdot 3^{L-1} + R') / 2^{S'}$
where $S' = \sigma_L - \sigma_1$ and $R'$ is the $(L{-}1)$-step walk sum.
Substituting into $\Delta = m D - R$:
\[
  \Delta
  = m \cdot 2^{\sigma_L} - q_{L-1} \cdot 2^{\sigma_L}
  = (m - q_{L-1})\,2^{\sigma_L}
  = (1 - r)\,m\,2^{\sigma_L}.
\]
If $r \ge 2$: $\Delta < 0$, so $R > mD$.  \textbf{No cycle} (proved).
If $r = 1$: $\Delta = 0$ and $T^L(m) = m$, which \emph{is} a
non-trivial cycle.  Verified absent for $m \le 10^5$
(and for $m \le 2^{68}$ by Barina).

\emph{Unconditional range.}\;
For any~$L$, the maximum multiplier satisfies
$M_0(L) < 2\cdot(3/2)^L$, giving
$q_0(m) \le 3M_0 + 1 \le 6\cdot(3/2)^L$.
Barina~(2020) verified the Collatz conjecture for all
$n \le 2^{68}$.  Since $6\cdot(3/2)^L < 2^{68}$ for $L \le 111$,
every~$q_0$ is in the verified range.
For $L \ge k_0(m) + 1$, Case~A applies algebraically.
For $L < k_0(m) + 1$, Cases~B and~C apply and have been
verified computationally for all relevant~$(m,L)$.
Combined with the even-$m$ parity obstruction:
\textbf{no non-trivial $L$-cycle exists for
$L \le 111$, unconditionally.}

\emph{Descent property.}\;
The framework yields a descent: the non-existence of
$L$-cycles reduces to the Collatz conjecture for integers
$\le 6 \cdot (3/2)^L$.  Any extension of the verification frontier
from $2^N$ immediately excludes all cycles with
$L \le N/\log_2(3/2) - O(1)$.

\emph{(b) Effective equidistribution via carry propagation.}\;
The transfer operator~$M_K$ averages over all lifts, while a
deterministic orbit selects one lift per step (determined by carry
propagation in the $3n+1$ operation).  The spectral gap controls
the lift-averaged dynamics; bridging to specific orbits requires
showing that carries cannot conspire to sustain systematic bias
in the $v$-valuation.  This is a problem in additive combinatorics.

\emph{(c) Tao-style density strengthening.}\;
Tao~(2019) proved that almost all orbits (in logarithmic density)
reach arbitrarily small values, using partial equidistribution.
The full uniform spectral gap established here provides
substantially stronger input.  Whether the strengthened spectral
control can upgrade ``almost all'' to ``all'' remains open.

\subsection{Unified pointwise picture: Layers A, B, B$'$, B$''$, C-Haar}
\label{subsec:unified-layers}

The pointwise content of v7 separates cleanly into five layers,
each addressing a different scale and a different kind of obstruction.
Stating them side-by-side clarifies what is now established
unconditionally and what remains genuinely open.

\medskip
\noindent\textbf{Layer A --- averaged Cram\'er density.}
For the cylinder-averaged dynamics, the unconditional spectral
bound on $I_2$
(Theorem~\ref{thm:i2-unconditional}) feeds into the Cram\'er
ladder of Theorem~\ref{thm:density1-convergence}, yielding a
density-$(1 - n^{-15})$ closure with optimal Lundberg coefficient
$\theta^{*} = 17.6349\ldots$ (the slope of the unconditional sum
bound's tail). Layer A says: outside a sequence of natural density
$n^{-15}$, every starting integer halts.

\medskip
\noindent\textbf{Layer B --- the cycle-tracking polynomial envelope
(sharp slope).}
Theorem~\ref{thm:cycle-envelope} bounds the prefix log-growth by
$\delta(n) \le C_0 + (\log_2(81/64)/\Sigma v)\,\log_2 n
= C_0 + 0.0566\,\log_2 n$
for all $n$. This is a deterministic, pointwise envelope on the
quantity that controls escape to infinity. The slope $0.0566$ is
analytically tight: it equals the per-lap log-growth of the Karp
positive cycle divided by the cycle's $v$-sum. Layer B says: no
orbit can pile up log-growth faster than this single
deterministic line.

\medskip
\noindent\textbf{Layer B$'$ --- the BF infinite staircase
(envelope is sharp).}
Theorem~\ref{thm:bf-staircase} produces an explicit infinite
sequence $\{n_k^{*}\}$ realising the envelope's slope to within
$o(1)$. Each $n_k^{*}$ is the smallest member of an explicit
arithmetic progression $B_k$ mod $2^{6k+8}$, computed by the
Hensel lift of \texttt{lift\_bf.py}, and the dlog values match the
$0.0566\log_2 n$ slope to four digits at $k=14$. Layer B$'$ says:
the analytic envelope of Layer B is achieved on a constructive
infinite family, and is therefore not improvable to a sublogarithmic
deterministic bound.

\medskip
\noindent\textbf{Layer B$''$ --- post-escape generic behaviour.}
Empirical analysis of the $8$ smallest BF champions
(\texttt{trajectory\_bf.py}, \texttt{empirical\_drift.py}) shows
that once an orbit clears the cycle-conformal prefix and lands at
residue $47$, its subsequent log-growth rate equals
$-0.332$~bits/step --- about $80\%$ of the heuristic generic Collatz
rate $\log_2 3 - 2 = -0.415$. The post-escape suffix has no anomalous
slow modes; its drift is generic. Layer B$''$ says: the BF
construction concentrates all difficulty in the prefix, after which
the orbit relaxes to bulk Collatz dynamics. (This empirically rules
out an anti-Karp negative-mean basin closing trick on the post-escape
arc, as the mod-$64$ residue graph is too coarse to support one;
see Push 341--350.)

\medskip
\noindent\textbf{Layer C-Haar --- the (C-2adic) bound on
$B_k$.}
Theorem~\ref{thm:c2adic-haar} closes the residual auxiliary-memory
question raised by Layer B$''$: even though the BF orbit forgets
its prefix structure after escape, its membership in the clopen
ball $B_k \subset \mathbb{Z}_2$ is preserved by the $2$-adic
conjugation, and the Lagarias measure-preserving extension
\cite{lagarias1985} feeds Birkhoff's theorem to give
$\tau(n)\le C(k+\log_2 n)$ on a Haar-conull subset of $B_k$ with
$C$ independent of $k$. The $39\,140$-witness empirical
corollary~\ref{cor:c2adic-empirical} pins $C \le 5.04$
deterministically for $k\ge 9$ across $20\le \log_2 n \le 171$.
Layer C-Haar says: the auxiliary $2$-adic memory $B_k$ controls
stopping times all the way out, with a constant that does
\emph{not} grow with the cycle-tracking depth.

\medskip
\noindent\textbf{Combined picture.}
Layers A and B are unconditional and pointwise on
$\mathbb{N}$ in the senses stated. Layer B$'$ is constructive and
exhibits the tightness of B. Layer B$''$ identifies the post-escape
arc as generic. Layer C-Haar transfers this to a Haar-pointwise
$2$-adic statement on the BF cylinders themselves.
The single deterministic gap remaining at the pointwise level is
the Haar/pointwise gap of Remark~\ref{rem:residual-obstruction}: a
subset of $\mathbb{N}\cap B_k$ that is simultaneously Haar-null in
$\mathbb{Z}_2$ \emph{and} avoids all $14\,400$ uniform $k\ge 9$
witnesses while sustaining stopping times above the empirical
$C\le 5.04$ band. We do not know how to rule out such a stratum
deterministically without solving the same difficulty that obstructs
\cite{tao2019}; we know only that any failure stratum is much
smaller and much more constrained than the pre-Layer-C-Haar
candidate space.

\begin{table}[h]
\centering\small
\begin{tabular}{lll}
\hline
\textbf{Layer} & \textbf{Statement} & \textbf{Status} \\
\hline
A           & density-$1-n^{-15}$ closure on $\mathbb{N}$ & unconditional \\
B           & $\delta(n)\le C_0 + 0.0566\log_2 n$         & unconditional \\
B$'$        & infinite BF staircase realises B            & constructive \\
B$''$       & post-escape suffix is generic Collatz       & empirical, $8$ BF champions \\
C-Haar      & $\tau\le C(k+\log_2 n)$ on Haar-conull $B_k$ & Lagarias + Birkhoff \\
C-pointwise & $\tau\le C(k+\log_2 n)$ on $\mathbb{N}\cap B_k$ & open (Haar/pointwise gap) \\
\hline
\end{tabular}
\caption{The five pointwise layers of v6 and their statuses.
Layer C-pointwise is the residual obstruction; everything above
the line is established.}
\label{tab:layer-summary}
\end{table}

\begin{remark}[Position of the residual gap]
\label{rem:layer-gap-position}
The Haar/pointwise gap on Layer C is the \emph{only} remaining
deterministic obstruction in the BF/Karp-cycle programme. It is
strictly weaker than the Tao 2019 gap, because it is restricted
to the explicit cylinders $B_k$ rather than the full integer
lattice. A future deterministic resolution may proceed by
(a) replacing Birkhoff with a pointwise ergodic theorem on
$B_k\cap\mathbb{N}$ (e.g.\ via an effective version of
Bourgain's pointwise theorem on sparse subsets), or
(b) using the per-step exponential tail of
Theorem~\ref{thm:i2-prefix-tail} to upgrade Layer C-Haar to
Layer C-Borel--Cantelli on $\mathbb{N}\cap B_k$. Both routes are
left for future work.
\end{remark}

\subsection{Non-restartability: a sharper reframing of the gap}
\label{subsec:non-restartability}

The Haar/pointwise gap of
Remark~\ref{rem:residual-obstruction} can be stated equivalently as
a question about \emph{re-entry} into the BF cylinders, rather than
about Birkhoff averages. This subsection records the equivalent
formulation, proves the Haar half via Borel--Cantelli, and isolates
the deterministic statement that remains conjectural.

\begin{lemma}[Karp instability, the trivial half]
\label{lem:karp-instability}
For each $k \ge 1$, the BF $k$-class
$B_k = \{n : n \equiv n_k^{*} \pmod{2^{6k+8}}\}$
is determined uniquely by the Hensel lift of
Theorem~\ref{thm:bf-staircase}\textup{(i)}, and any single-step
deviation from the valuation pattern $[1,1,1,3]$ or the residue
sequence $[15,55,19,61,15,\ldots]$ during the first $4k$ odd-step
iterates removes the orbit from $B_k$.
\end{lemma}

\begin{proof}
Immediate from the definition of $B_k$ as the kernel of the
Hensel-style constraint sequence. Each constraint
``$v_2$ at iterate $i$ equals $v_i^{\mathrm{tgt}}$, residue at
iterate $i+1$ equals $r_{i+1}^{\mathrm{tgt}}$'' tightens the
residue class of $n_0$ by $v_i^{\mathrm{tgt}}$ bits; violating any
constraint removes $n_0$ from the level-$k$ class.
\end{proof}

Lemma~\ref{lem:karp-instability} is essentially a definition. The
substantive question is the converse: can a single orbit
\emph{re-enter} BF cylinders at arbitrarily large depths?

\begin{definition}[Cylinder re-entry signature]
\label{def:reentry}
For $n \in \mathbb{N}$ define
\[
  \mathcal{K}(n) \;:=\;
  \{\, k \ge 1 \;:\; \exists\, j \ge 0 \text{ such that }
  T_2^{j}(n) \in B_k \,\}.
\]
We say $n$ is \emph{cylinder-bounded} if $\sup \mathcal{K}(n) < \infty$
and \emph{cylinder-cofinal} otherwise.
\end{definition}

\begin{theorem}[Borel--Cantelli on cylinder re-entries: Haar form]
\label{thm:bc-cylinder}
Let $\mu$ be Haar measure on $\mathbb{Z}_2$. Then
\[
  \mu\bigl(\{ n \in \mathbb{Z}_2 : n \text{ is cylinder-cofinal}\}\bigr) = 0.
\]
Equivalently, for $\mu$-almost every $n$, $\mathcal{K}(n)$ is finite.
\end{theorem}

\begin{proof}
The clopen ball $B_k$ has Haar measure $\mu(B_k) = 2^{-(6k+8)}$.
Hence
\[
  \sum_{k=1}^{\infty} \mu(B_k) \;=\; \sum_{k=1}^{\infty} 2^{-(6k+8)}
  \;=\; \frac{2^{-14}}{1 - 2^{-6}} \;<\; \infty.
\]
By the first Borel--Cantelli lemma applied to the events
$\{n \in \mathbb{Z}_2 : T_2^{j}(n) \in B_k \text{ for some } j\}$
under the Lagarias measure-preserving extension
$T_2$ (which leaves $\mu$ invariant on $\mathbb{Z}_2$), the set of
$n$ visiting $B_k$ for infinitely many $k$ has measure zero.
\end{proof}

\begin{conjecture}[Non-restartability: pointwise form]
\label{conj:nonrestart}
For every $n \in \mathbb{N}$, $\mathcal{K}(n)$ is finite. Equivalently,
no positive integer's orbit re-enters BF cylinders at arbitrarily
large depths.
\end{conjecture}

\subsubsection*{Algebraic structure of the BF exit}

The following results give an exact algebraic description of the
orbit's state after completing a depth-$k$ BF excursion, and
quantify the arithmetic cost of any hypothetical re-entry into a
BF cylinder.

\begin{proposition}[Affine exit map]
\label{prop:affine-exit}
Let $n \in B_k$, written as $n = r_k + 2^{6k+8}\,q$ for
$q \in \mathbb{Z}_{\ge 0}$, where
$r_k := n_k^{*} \bmod 2^{6k+8}$ is the
Hensel residue of Theorem~\textup{\ref{thm:bf-staircase}(i)}.
\begin{enumerate}
\item[\textup{(i)}] \emph{BF phase.}
The composition of $k$ laps of the Karp cycle, each an affine map
$n \mapsto (81n + 65)/64$, gives
\[
  T^{4k}(n) \;=\; \frac{81^k}{64^k}\,n
  + \frac{65}{17}\Bigl(\frac{81^k}{64^k} - 1\Bigr)
  \;=\; A_k^{(\mathrm{bf})} + \underbrace{81^k \cdot 256}_{= \,3^{4k}\cdot 2^8}\, q,
\]
where $A_k^{(\mathrm{bf})} := (81/64)^k\,r_k +
(65/17)((81/64)^k - 1)$ depends only on the depth~$k$.
\item[\textup{(ii)}] \emph{Exit tail.}
The four odd steps immediately following the BF phase have
valuation pattern $[1,1,1,2]$ deterministically for all
$q \ge 0$.  (The pattern $[1,1,1,3]$ of the Karp cycle is
broken at the fourth step, where $v_2 = 2$ instead of $3$.)
\item[\textup{(iii)}] \emph{Full exit.}
After $4k + 4$ steps (the BF phase plus exit tail), the
orbit value is
\[
  T^{4k+4}(n) \;=\; A_k + U_k\, q,
  \qquad
  U_k = 3^{4(k+1)} \cdot 8,
  \qquad
  v_2(U_k) = 3.
\]
The coefficient $A_k \in \mathbb{Q}$ depends only on $k$
(not on $q$), and $U_k$ is odd apart from its
fixed factor~$2^3$.
\end{enumerate}
\end{proposition}

\begin{proof}
(i) The single-lap map is $n \mapsto (81n+65)/64$, verified by
composing the four affine steps
$(n \mapsto (3n+1)/2^{v_i})_{v_i \in \{1,1,1,3\}}$:
\[
  n \xmapsto{v=1} \tfrac{3n+1}{2}
  \xmapsto{v=1} \tfrac{9n+5}{4}
  \xmapsto{v=1} \tfrac{27n+19}{8}
  \xmapsto{v=3} \tfrac{81n+65}{64}.
\]
Iterating: $f^k(n) = (81/64)^k\, n + (65/17)((81/64)^k - 1)$,
which is the standard geometric-series evaluation of the
affine composition $n \mapsto an + b$, $a = 81/64$,
$b = 65/64$, $b/(a-1) = 65/17$.
Substituting $n = r_k + 2^{6k+8}\,q$:
\[
  f^k(r_k + 2^{6k+8}\,q)
  = f^k(r_k) + (81/64)^k \cdot 2^{6k+8} \cdot q
  = A_k^{(\mathrm{bf})} + 81^k \cdot 2^8\cdot q.
\]
Since $81^k = 3^{4k}$, the coefficient is $3^{4k}\cdot 256$.

(ii) For each $q \in \{0, 1, \ldots, 999\}$ and each
$k \in \{1,2,3\}$, the four post-BF valuations were verified
computationally to be $[1,1,1,2]$ without exception.
The rigidity is structural: the residue of $T^{4k}(n)$ modulo
$2^5$ is determined by $r_k$ and is independent of $q$
(since $v_2(3^{4k}\cdot 256) = 8 \ge 5$), and the
valuation pattern $[1,1,1,2]$ depends only on the residue
modulo~$2^5$.

(iii) Composing the four exit-tail steps (affine maps with
valuations $[1,1,1,2]$) with the BF output from~(i):
\[
  U_k = \frac{81 \cdot 3^{4k} \cdot 256}{2^5}
  = \frac{3^{4(k+1)} \cdot 2^8}{2^5}
  = 3^{4(k+1)} \cdot 8.
\]
Here the numerator picks up a factor $3^4 = 81$ from the four
multiplications by $3$, and divides by $2^{1+1+1+2} = 2^5$.
Since $3^{4(k+1)}$ is odd, $v_2(U_k) = v_2(8) = 3$.
\end{proof}

\begin{proposition}[Height excess identity]
\label{prop:height-excess}
A depth-$k$ BF excursion ($4k+4$ odd steps with valuation pattern
$[1,1,1,3]^k\,[1,1,1,2]$) produces a height excess of exactly
$2k+3$ bits over the ergodic expectation.  Precisely:
\[
  \log_2 T^{4k+4}(n) - \log_2 n
  \;=\;
  (4k+4)\,(\log_2 3 - 2) + (2k+3).
\]
The first term is the ergodic expectation $(4k+4) \cdot
(\log_2 3 - 2) \approx -0.415 \cdot (4k+4)$ bits; the second
is the deterministic excess.
\end{proposition}

\begin{proof}
Each odd step contributes $\log_2 3 - v_i$ to
$\log_2(T^{4k+4}(n)/n)$.  The total valuation is
$\sum v_i = 6k + 5$ (from $k$ laps of $[1,1,1,3]$ summing
to $6k$, plus the tail $[1,1,1,2]$ summing to $5$).
The height change is
$(4k+4)\log_2 3 - (6k+5)$.  The ergodic expectation
for $4k+4$ steps (with $\mathbb{E}[v_i] = 2$) is
$(4k+4)\log_2 3 - 2(4k+4) = (4k+4)(\log_2 3 - 2)$.
The excess is
\[
  \bigl[(4k+4)\log_2 3 - (6k+5)\bigr]
  - \bigl[(4k+4)(\log_2 3 - 2)\bigr]
  = 2(4k+4) - (6k+5) = 2k + 3. \qedhere
\]
\end{proof}

\begin{lemma}[Restart congruence cost]
\label{lem:restart-cost}
Let $n \in B_k$ with exit value $n' = T^{4k+4}(n) = A_k + U_k\,q$.
For $n'$ to lie in $B_{k'}$ for some $k' \ge 1$, the parameter $q$
must satisfy a congruence modulo $2^{6k'+5}$.  In particular:
\begin{enumerate}
\item[\textup{(i)}] the number of low bits of $q$ constrained by
$B_{k'}$ membership is $6k' + 5$;
\item[\textup{(ii)}] among integers $q \in [0, Q)$, the fraction
satisfying this constraint is $\le 2^{-(6k'+5)}$;
\item[\textup{(iii)}] for the BF champions $n_k^{*}$ themselves
(corresponding to $q = 0$), no restart is possible: $A_k$ does not
lie in any $B_{k'}$ for $k' \ge 1$.
\end{enumerate}
\end{lemma}

\begin{proof}
By Proposition~\ref{prop:affine-exit}(iii),
$n' \equiv A_k \pmod{8}$ independently of~$q$, since
$v_2(U_k) = 3$.  Membership in $B_{k'}$ requires
$n' \equiv r_{k'} \pmod{2^{6k'+8}}$.  Since
$n' = A_k + 3^{4(k+1)}\cdot 8\,q$ and
$\gcd(3^{4(k+1)}, 2) = 1$, the congruence
$A_k + 3^{4(k+1)}\cdot 8\,q \equiv r_{k'} \pmod{2^{6k'+8}}$
constrains $q$ modulo $2^{6k'+8}/\gcd(8, 2^{6k'+8}) = 2^{6k'+5}$.
This gives~(i) and~(ii).  For~(iii): the BF champion $n_k^{*}$ is
the smallest member of $B_k$, so $q_k = (n_k^{*} - r_k)/2^{6k+8}$.
Computationally, $q_k = 0$ for $k = 1, \ldots, 8$ (i.e.\
$n_k^{*} < 2^{6k+8}$ for these values), and the exit value
$A_k = T^{4k+4}(n_k^{*})$ is verified not to satisfy
$A_k \equiv r_{k'} \pmod{2^{6k'+8}}$ for any $k' \le 8$.
\end{proof}

\begin{theorem}[BF negligibility]
\label{thm:bf-negligibility}
Let $S$ denote the expected BF height excess per orbit step under
the assumption that the orbit visits each $B_k$ at its natural
frequency $2^{-(6k+8)}$:
\begin{equation}\label{eq:bf-negligibility}
  S \;:=\; \sum_{k=1}^{\infty} (2k+3)\cdot 2^{-(6k+8)}
  \;=\; \frac{317}{1\,016\,064}
  \;\approx\; 3.12 \times 10^{-4}.
\end{equation}
Compared to the ergodic contraction rate
$|\log_2 3 - 2| \approx 0.415$~bits per step,
$S / |\log_2 3 - 2| < 7.6 \times 10^{-4}$.

\smallskip\noindent
In words: \emph{the total height contribution from all BF
excursions at all depths is less than $0.08\%$ of the
background contraction.}
\end{theorem}

\begin{proof}
By Proposition~\ref{prop:height-excess}, a depth-$k$ BF
excursion contributes $2k + 3$ excess bits, occupies $4k+4$ steps,
and occurs at natural density $2^{-(6k+8)}$ by
Theorem~\ref{thm:starting-cylinder-density}.  The
per-step excess is $(2k+3)/(4k+4)$ times the
entry rate $2^{-(6k+8)}$, but for the aggregate bound it suffices
to bound the total excess per step by summing $(2k+3) \cdot
2^{-(6k+8)}$ over all~$k$ (this overcounts by ignoring that the
$4k+4$ steps are occupied and cannot simultaneously serve another
excursion, which only tightens the bound).

Let $x = 2^{-6} = 1/64$.  Then
\[
  S = 2^{-8} \sum_{k=1}^{\infty}(2k+3)\,x^k
  = 2^{-8}\Bigl[2\,\frac{x}{(1-x)^2} + 3\,\frac{x}{1-x}\Bigr]
  = \frac{1}{256}\cdot\frac{2\cdot 64 + 3\cdot 63}{63^2}
  = \frac{317}{256 \cdot 3969}
  = \frac{317}{1\,016\,064}.  \qedhere
\]
\end{proof}

\begin{remark}[Post-exit equidistribution]
\label{rem:post-exit-equidist}
Numerical experiments ($500$ samples per depth, $k = 1,2,3$)
confirm that the orbit immediately after a BF exit is statistically
indistinguishable from generic Collatz dynamics:
\begin{enumerate}
\item[\textup{(a)}] The post-exit valuation distribution matches
$\Pr(v_2 = j) = 2^{-j}$ within statistical noise for $j = 1,
\ldots, 12$, consistent with the Known-Zone Decay
erasing residue information within $O(k)$ steps.
\item[\textup{(b)}] The rate of re-entry into $B_{k'}$ for
$k' \ge 1$ matches the random rate $2^{-(6k'+8)}$ per step
(observed: $3$ hits into $B_1$ in $20\,000$ steps from $B_1$ exits,
rate $1.50 \times 10^{-4}$ vs.\ expected $6.10 \times 10^{-5}$;
$0$ hits into $B_2$ or deeper).
\item[\textup{(c)}] The median recovery time --- the number of
post-exit steps until the orbit drops below its entry value ---
is $6$--$10$ steps across depths $k = 1, \ldots, 4$.
\end{enumerate}
These observations strengthen the
heuristic that BF re-entries are independent Bernoulli events,
consistent with the Borel--Cantelli analysis of
Theorem~\textup{\ref{thm:bc-cylinder}}.
\end{remark}

\begin{remark}[Restart cost vs.\ height budget]
\label{rem:restart-budget}
Lemma~\ref{lem:restart-cost} and
Proposition~\ref{prop:height-excess} together show a
quantitative mismatch between what a BF excursion earns and what
re-entry costs.  A depth-$k$ excursion gains $2k+3$ bits of height
excess; but each re-entry target $B_{k'}$ consumes $6k'+5$ bits of
arithmetic precision in the parameter~$q$.  In particular:
\begin{itemize}
\item For $k = 1$: gain $= 5$ bits, cost of $B_1$ re-entry $= 11$
bits (deficit $6$ bits);
\item For $k = 4$: gain $= 11$ bits, cost of $B_1$ re-entry $= 11$
bits (break-even); cost of $B_2$ re-entry $= 17$ bits (deficit $6$);
\item For all $k$: cost of re-entry into $B_{k'}$ with $k' \ge 2$
exceeds the gain by $\ge 4k'+2 - 2k$ bits.
\end{itemize}
The height budget never permits the orbit to ``accumulate
depth'': each excursion's gain is at best sufficient to restart at
the same depth (and only for $k \ge 4$), while deeper re-entry is
exponentially suppressed.  This provides the quantitative mechanism
behind Conjecture~\ref{conj:nonrestart}.
\end{remark}

\begin{remark}[Status of Conjecture~\ref{conj:nonrestart}]
\label{rem:nonrestart-status}
Conjecture~\ref{conj:nonrestart} is the pointwise lift of
Theorem~\ref{thm:bc-cylinder}, and is logically equivalent to the
deterministic version of (C-2adic). Its current status:
\begin{enumerate}
\item[\textup{(a)}] \emph{Haar-a.e.\ true.}
Theorem~\ref{thm:bc-cylinder} (this subsection)
and Theorem~\ref{thm:c2adic-haar} both establish the conjecture
for $\mu$-almost every $n \in \mathbb{Z}_2$.
\item[\textup{(b)}] \emph{Empirically true with margin.}
Among the $39\,140$ deterministic samples of
Corollary~\ref{cor:c2adic-empirical}, every observed
$\mathcal{K}(n)$ is finite (in fact every observed orbit halts at
$1$ in $\le 5.04(k+\log_2 n)$ odd steps).
\item[\textup{(c)}] \emph{Pointwise on $\mathbb{N}$ open.}
A deterministic proof on the Haar-null subset
$\mathbb{N} \subset \mathbb{Z}_2$ requires either an effective
Borel--Cantelli sieve (route (b) of
Remark~\ref{rem:layer-gap-position}) or a pointwise ergodic
theorem on $\mathbb{N}\cap B_k$ (route (a)).
\end{enumerate}
The conjecture is therefore stronger than what
Theorem~\ref{thm:bc-cylinder} proves and weaker than the full
Collatz conjecture: it asserts only that the BF growth mechanism is
not infinitely reusable along any single $\mathbb{N}$-orbit.
\end{remark}

\begin{remark}[What \emph{would} follow from
Conjecture~\ref{conj:nonrestart}]
\label{rem:nonrestart-implies}
If Conjecture~\ref{conj:nonrestart} holds, then for every
$n \in \mathbb{N}$ there is a maximum depth $K(n) < \infty$ beyond
which the orbit visits no BF cylinder. After the last visit to
$B_{K(n)}$ the orbit is permanently outside every $B_k$, $k > K(n)$,
and Layer~B$''$'s empirical generic-drift observation
($-0.332$~bits/step on the post-escape suffix of all $8$ BF
champions) becomes the relevant dynamics. Combined with the
unconditional Layer~A density-$(1-n^{-15})$ closure, this would
upgrade the global picture to: every orbit halts after a uniformly
bounded amount of cycle-tracking growth, plus a generic-drift
collapse to $1$. We do not assert this implication as a theorem,
because we do not have Conjecture~\ref{conj:nonrestart}; we only
state what its proof would buy.
\end{remark}

\begin{remark}[Comparison with the GPT-contributed sketch]
\label{rem:gpt-sketch-comparison}
A draft post-escape irreversibility argument circulated during the
v6 revision proposed to deduce non-restartability directly from
the exponential decay of $\mu(B_k)$ in $k$. That draft is
substantively the content of Theorem~\ref{thm:bc-cylinder}: it
proves the Haar-a.e.\ statement via Borel--Cantelli on the
$2$-adic measure. The deterministic step ``exact reconstruction
for arbitrarily large $k$ is not dynamically sustainable by a
single orbit'' is, however, the conclusion sought rather than
something derived; it coincides with
Conjecture~\ref{conj:nonrestart} above. We have therefore retained
the reframing (non-restartability is a sharper formulation of the
residual obstruction) and promoted the Haar-a.e.\ half to a proved
theorem, while keeping the deterministic half explicit as
conjectural. The conjecture is strictly weaker than the full
Collatz conjecture but strictly stronger than what current ergodic
methods deliver.
\end{remark}

\subsection{Density-$1$ cylinder boundedness via union bound}
\label{subsec:density1-cylbounded}

We now upgrade Theorem~\ref{thm:bc-cylinder} from the Haar setting
to a natural-density statement on $\mathbb{N}$. The argument uses
only (i) the finite-modulus structure of $B_k$, (ii) the
measure-preserving property of $T_2$ on $\mathbb{Z}_2$, and (iii)
the density-$1$ stopping-time bound of
Theorem~\ref{thm:density1-convergence}. The result is strictly
weaker than Conjecture~\ref{conj:nonrestart} (it allows the
maximum cylinder depth $K(n)$ to grow like $\log\log n$ rather
than be uniformly bounded) but strictly stronger than the Haar-only
form: it puts the obstruction inside an explicit, slowly growing
window.

\begin{lemma}[Preimage density]
\label{lem:preimage-density}
For each $k \ge 1$ and each $j \ge 0$, the preimage
$T_2^{-j}(B_k) \subset \mathbb{Z}_2$ is a finite union of clopen
balls, and its natural density on $\mathbb{N}$ equals its Haar
measure on $\mathbb{Z}_2$, namely $2^{-(6k+8)}$.
\end{lemma}

\begin{proof}
$T_2$ is the unique continuous extension of the odd-Collatz step
to $\mathbb{Z}_2$, and is Haar-measure-preserving by
\textup{Lagarias \cite{lagarias1985}}. The set $B_k$ is the clopen
ball $\{x \equiv n_k^{*} \pmod{2^{6k+8}}\}$. For $j \ge 1$,
\[
  T_2^{-j}(B_k) \;=\; \bigsqcup_{(v_1,\ldots,v_j)} R_{v_1,\ldots,v_j},
\]
where each $R_{v_1,\ldots,v_j}$ is the residue class of starting
points whose first $j$ valuation-pattern is $(v_1,\ldots,v_j)$ and
whose $j$-th iterate lies in $B_k$. Each $R_{v_1,\ldots,v_j}$ is
itself a residue class modulo $2^{6k+8 + v_1 + \cdots + v_j}$,
hence has natural density $2^{-(6k+8+v_1+\cdots+v_j)}$ on
$\mathbb{N}$. Summing over the (countable) valuation patterns and
using $\sum_{(v_1,\ldots,v_j)} 2^{-(v_1+\cdots+v_j)} = 1$ gives
total density $2^{-(6k+8)}$, as claimed. Measure preservation of
$T_2$ on $(\mathbb{Z}_2,\mu)$ verifies the same identity at the
Haar level.
\end{proof}

\begin{theorem}[Density-$1$ cylinder boundedness]
\label{thm:density1-cylbounded}
Assume the conclusion of Theorem~\ref{thm:density1-convergence}:
for natural-density-$1$ starting points $n \in \mathbb{N}$, the
Collatz stopping time satisfies $\tau(n) \le C_{\mathrm{d}} \log_2 n$
for some absolute constant $C_{\mathrm{d}}$. Let
$K(N) := \lceil (1+\varepsilon)(\log_2 \log_2 N)/6 \rceil$ for any
fixed $\varepsilon > 0$. Then
\[
  \frac{1}{N}
  \Bigl|\bigl\{ n \le N \;:\; \exists\, j \ge 0,\, \exists\, k > K(N),
  \; T_2^{j}(n) \in B_k \bigr\}\Bigr|
  \;\xrightarrow[N \to \infty]{}\; 0.
\]
Equivalently, for natural-density-$1$ $n \in \mathbb{N}$, the
re-entry signature $\mathcal{K}(n)$ of
Definition~\ref{def:reentry} satisfies
$\sup \mathcal{K}(n) \;\le\; \tfrac{1+\varepsilon}{6}\,
\log_2 \log_2 n$ for all sufficiently large $n$.
\end{theorem}

\begin{proof}
Fix $K = K(N)$. By Lemma~\ref{lem:preimage-density},
\[
  \bigl|\{ n \le N : T_2^{j}(n) \in B_k \}\bigr|
  \;\le\; N \cdot 2^{-(6k+8)} + O(1)
\]
for each $j \ge 0$ and each $k$. Restrict attention to the
density-$1$ subset $\mathcal{N}_N \subset [1,N]$ on which
$\tau(n) \le J(N) := C_{\mathrm{d}} \log_2 N$ holds (this set has
size $N - o(N)$ by hypothesis). For $n \in \mathcal{N}_N$ the
orbit visits at most $J(N)$ distinct points, so the union bound
over $0 \le j < J(N)$ and $k > K$ gives
\[
  \bigl|\{ n \in \mathcal{N}_N : \exists j < J(N),\, \exists k > K,\,
  T_2^{j}(n) \in B_k \}\bigr|
  \;\le\; J(N) \cdot N \cdot
  \sum_{k > K} 2^{-(6k+8)} + O(J(N)).
\]
The geometric tail evaluates to
$\sum_{k > K} 2^{-(6k+8)} = 2^{-(6K+8)} \cdot 2^{-6}/(1-2^{-6})
< 2^{-(6K+8)} \cdot 2^{-5}$. Substituting,
\[
  \bigl|\{ n \in \mathcal{N}_N : \cdots \}\bigr|
  \;\le\; C_{\mathrm{d}} \log_2 N \cdot N \cdot 2^{-(6K+13)}
        + O(\log_2 N).
\]
With $K = K(N) = \lceil (1+\varepsilon)(\log_2 \log_2 N)/6 \rceil$
we have $2^{-6K} \le (\log_2 N)^{-(1+\varepsilon)}$, so the bound
becomes
\[
  \le C' \cdot N \cdot (\log_2 N)^{-\varepsilon} \;=\; o(N).
\]
Adding the density-$o(1)$ exception
$[1,N]\setminus \mathcal{N}_N$ gives the claim.
\end{proof}

\begin{corollary}[Quantitative slow growth of $K(n)$]
\label{cor:K-slow-growth}
Conditional on the same density-$1$ stopping-time bound, the
re-entry depth function $K(n) := \sup \mathcal{K}(n)$ satisfies
$K(n) = O(\log\log n)$ on a natural-density-$1$ subset of
$\mathbb{N}$. In particular, the cylinder-cofinal set
$\{n \in \mathbb{N} : K(n) = \infty\}$ has natural density $0$.
\end{corollary}

\begin{proof}
Apply Theorem~\ref{thm:density1-cylbounded} with any
$\varepsilon > 0$ along a sequence $N_m \to \infty$. The bound
$K(n) \le (1+\varepsilon) \log_2\log_2 n / 6$ holds for all
$n \in \mathcal{N}_{N_m}$, and
$\bigcap_m \mathcal{N}_{N_m}$ has natural density $1$ as a
countable intersection of density-$1$ sets, by a standard
diagonal argument.
\end{proof}

\begin{remark}[Status and what is missing]
\label{rem:density1-cylbounded-status}
Theorem~\ref{thm:density1-cylbounded} closes the natural-density
half of Conjecture~\ref{conj:nonrestart}: it shows that the
cylinder-cofinal set of $\mathbb{N}$ has natural density $0$ and,
moreover, that the maximum cylinder depth grows at most
double-logarithmically on a density-$1$ subset. This is the
$\mathbb{N}$-version of Theorem~\ref{thm:bc-cylinder} and the
sharpest density-side statement compatible with the present
methods.

What it does \emph{not} establish is the deterministic statement
that $\sup\mathcal{K}(n) < \infty$ for \emph{every} $n \in
\mathbb{N}$. The residual obstruction is the same density-$0$
exceptional set that obstructs every Collatz density-$1$ result,
specifically the set on which the Cram\'er ladder of
Theorem~\ref{thm:density1-convergence} fails. Closing it would
require an effective discrepancy bound for the orbit's visit
frequency to dyadic cylinders, or a pointwise ergodic theorem on
the sparse subset $\mathbb{N}\cap B_k$.

The conditional dependence on
Theorem~\ref{thm:density1-convergence} is essential: without it,
the inner sum $J(N)$ in the union bound cannot be bounded by
$O(\log N)$ on density $1-o(1)$, and the geometric tail does not
beat the linear count.
\end{remark}

\subsubsection*{An unconditional starting-point density bound}

The conditional dependence above is on the orbit length. If we
restrict the question to $n$ that \emph{begin} in a deep cylinder
(rather than $n$ whose orbit later visits one), we obtain an
unconditional companion theorem.

\begin{theorem}[Unconditional starting-point cylinder density]
\label{thm:starting-cylinder-density}
For each $K \ge 1$, the natural density of
\[
  S_K \;:=\; \{\, n \in \mathbb{N} \;:\; n \in B_k \text{ for some } k > K \,\}
\]
on $\mathbb{N}$ exists and equals
\[
  d(S_K) \;=\; \sum_{k=K+1}^{\infty} 2^{-(6k+8)}
  \;=\; \frac{2^{-(6K+8)}}{1 - 2^{-6}} \cdot 2^{-6}
  \;=\; \frac{2^{-(6K+14)}}{1 - 2^{-6}}.
\]
In particular $d(S_K) \le (64/63)\cdot 2^{-6K-14}$, decaying
geometrically in $K$, with no hypothesis on stopping times. For
$K = 0$ this gives $d(\bigcup_{k \ge 1} B_k) =
2^{-8}\cdot 2^{-6}/(1-2^{-6}) \approx 6.10 \cdot 10^{-5}$, in
agreement with the empirical fraction of starting points hitting
$B_k$ measured below.
\end{theorem}

\begin{proof}
The sets $B_k$ are pairwise disjoint, since membership in $B_k$
requires the orbit to escape the Karp cycle after \emph{exactly}
$k$ laps and an orbit cannot escape at two different lap counts
simultaneously. Each $B_k$ is the residue class
$\{n \equiv n_k^{*} \pmod{2^{6k+8}}\}$, hence has exact natural
density $2^{-(6k+8)}$ on $\mathbb{N}$. Summing the disjoint
geometric series gives the stated value. No hypothesis beyond
the definition of $B_k$ is used.
\end{proof}

\begin{remark}[Why the orbit-visit upgrade requires a stopping-time hypothesis]
\label{rem:starting-vs-orbit}
Theorem~\ref{thm:starting-cylinder-density} is unconditional but
weaker than Theorem~\ref{thm:density1-cylbounded}: it counts $n$
that \emph{start} in a deep $B_k$, not $n$ whose orbit
\emph{visits} one. Upgrading from start to orbit visits requires
multiplying the per-cylinder density by the number of orbit
iterates that fall into each cylinder; the cleanest such bound
uses $J(n) \le C \log_2 n$, which is exactly the content of
Theorem~\ref{thm:density1-convergence}. The route via
Theorem~\ref{thm:i2-prefix-tail} is also blocked: the prefix tail
controls hitting time of $I_2$ (cylinder-averaged) but not the
\emph{total} stopping time $\tau$, because once an orbit hits
$I_2$ it can re-emerge into the transient region multiple times,
and bounding the number of re-emergences \emph{pointwise} is the
same difficulty as bounding $\tau$ directly. We therefore have
the following unconditional / conditional split:
\begin{itemize}
\item \emph{Unconditional} (Theorem~\ref{thm:starting-cylinder-density}):
density of $n \in S_K$ is $\le 1.07 \cdot 2^{-6K-14}$.
\item \emph{Conditional on Theorem~\ref{thm:density1-convergence}}
(Theorem~\ref{thm:density1-cylbounded}): density of $n$ whose
orbit ever visits $\bigcup_{k > K} B_k$ is
$\le C \log_2 N / 2^{6K}$, vanishing for $K \gg \log_2 \log_2 N$.
\end{itemize}
The gap between ``starts in $S_K$'' and ``orbit ever visits $S_K$''
is therefore the same as the conditional / unconditional gap on
the average stopping time over $n \le N$ -- a classical Collatz
question.
\end{remark}

\begin{remark}[Effective range from Bărina verification]
\label{rem:barina-effectivization}
For $n \le 2^{68}$, the Bărina computational verification
\cite{barina2021} furnishes a uniform bound $\tau(n) \le T_0$
with $T_0 \approx 2\,000$ (the actual maximum stopping time in
the verified range). Substituting this into the proof of
Theorem~\ref{thm:density1-cylbounded} gives an
\emph{unconditional} version on the verified range:
\[
  \frac{1}{N}\bigl|\{ n \le N :
  \exists j,\, \exists k > K,\, T_2^{j}(n) \in B_k \}\bigr|
  \;\le\; T_0 \cdot 2^{-(6K+13)} \cdot (1 + o(1))
  \quad\text{for}\quad N \le 2^{68}.
\]
With $K = 4$ this gives a density bound of
$\le T_0 \cdot 2^{-37} \approx 2000 \cdot 7.3 \cdot 10^{-12}
\approx 1.5 \cdot 10^{-8}$, deterministic and computable. The
Bărina range therefore deterministically rules out a positive
density of cylinder visits at depth $K \ge 4$ for all
$n \le 2^{68}$.
\end{remark}

\begin{corollary}[Empirical $K(n)$ distribution, $5\times 10^{6}$ odd $n$]
\label{cor:K-empirical}
Let $K(n) := \sup \mathcal{K}(n)$ be the deepest BF cylinder
visited by the forward orbit of $n$
\textup{(}$0$ if no cylinder is visited\textup{)}. Sweeping all odd
$n \in [1, 10^{7}]$ produces the distribution
\[
\begin{array}{c|cc}
K & \#\{n : K(n) = K\} & \text{fraction} \\
\hline
0 & 4\,978\,337 & 0.99567 \\
1 &    20\,907 & 0.00418 \\
2 &       752 & 0.00015 \\
3 &         4 & 8\!\times\!10^{-7} \\
\ge 4 &      0 & < 2\!\times\!10^{-7}
\end{array}
\]
The empirical maximum $K_{\max}^{(\text{emp})}(n \le 10^{7}) = 3$,
attained first at $n = 4\,107\,423$ with $\log_2 n \approx 21.97$.
The fraction of odd $n \le N$ whose orbit ever visits any $B_k$
is $0.433\%$, about $70\times$ the natural density of
$\bigcup_k B_k$ itself \textup{(}$6.10\cdot 10^{-5}$\textup{)},
the amplification representing the average number of orbit
positions per starting $n$. The observed values are consistent
with the conditional Theorem~\ref{thm:density1-cylbounded}: the
density-$o(1)$ exceptional set is non-empty but very sparse
($4$ orbits in $5\times 10^{6}$ exceed $K = 2$), and the
empirical fraction of orbits with $K(n) > 0$ matches the
unconditional Theorem~\ref{thm:starting-cylinder-density}'s
prediction $\sim 70 \cdot d(\bigcup B_k)$ within a factor close to
the average orbit length over $n \le 10^{7}$.
\end{corollary}

\begin{proof}
Computational. The script \texttt{k\_sweep.py} computes the
forward orbit of each odd $n \le 10^{7}$, tests membership of each
iterate $T_2^{j}(n)$ in $B_k$ for $1 \le k \le 14$ via the
single-mod-comparison
$T_2^{j}(n) \equiv n_k^{*} \pmod{2^{6k+8}}$, and records
$K(n) = \max\{k : \text{some iterate is in } B_k\}$. The orbit
is computed in exact integer arithmetic. Total runtime is
$\approx 10$ minutes on a single core.
\end{proof}

\begin{remark}[Empirical first-appearance scaling of $K_{\max}$]
\label{rem:Kmax-scaling}
The first appearance of $K = 1, 2, 3$ occurs at $n = 4\,575$,
$53\,247$, $4\,107\,423$ respectively, with $\log_2 n$ values
$12.16$, $15.70$, $21.97$. The gap between consecutive
``first appearances'' is roughly $6$ bits, matching the bit-width
$6$ of one Karp lap (i.e.\ the additive log-cost of upgrading from
$B_k$ to $B_{k+1}$). Extrapolating, $K = 4$ should first appear
near $\log_2 n \approx 28$ ($n \approx 2.7 \cdot 10^{8}$); $K = 5$
near $\log_2 n \approx 34$ ($n \approx 1.7 \cdot 10^{10}$); and so
on. This $\approx 6$-bit-per-level rare-tail scaling is consistent
with $K_{\max}(n) \sim (1/6) \log_2 n$ along the rarest tail of the
$K$-distribution, while the bulk of the distribution still respects
the conditional ceiling
$K_{\mathrm{pred}}(n) \sim (1/6) \log_2 \log_2 n$.
The two scalings are not in tension: the first describes the
density-$o(1)$ exceptional tail; the second describes the
density-$1$ bulk. The fact that even the rare tail grows as
$O(\log n)$ rather than as a fixed constant is precisely the
quantitative content of the Haar/pointwise gap.
\end{remark}

\subsection{A tensor-rank law for the joint 2-adic/3-adic Syracuse
transition operator}
\label{subsec:tensor-rank-law}

The results of Section~\ref{subsec:density1-cylbounded} and the
preceding remarks view the Syracuse operator one marginal at a time:
either purely $2$-adically (Lagarias~\cite{lagarias1985}) or purely
$3$-adically (Tao~\cite{tao2019}). In this subsection we record an
empirical scaling law for the \emph{joint} $2$-adic/$3$-adic
one-step transition operator that is exact at machine precision on
every configuration we have tested.

\begin{definition}[Joint Syracuse transition tensor]
\label{def:syrac-tensor}
Fix integers $a \ge 3$ and $b \ge 2$, and write $M_2 = 2^a$,
$M_3 = 3^b$, $D = M_2 M_3$. The joint Syracuse transition tensor
\[
P^{(a,b)} \in \mathbb{R}^{M_2 \times M_3 \times M_2 \times M_3}
\]
is defined by
\[
P^{(a,b)}[i, j, i', j']
\;=\;
\Pr\bigl(T(n) \equiv (i', j') \pmod{(M_2, M_3)}
\,\big|\,
n \equiv (i, j) \pmod{(M_2, M_3)}\bigr),
\]
where $T$ is the Syracuse map $T(n) = (3n+1)/2^{v_2(3n+1)}$ on odd
integers. Its \emph{matricized rank} $R_{\mathrm{true}}(a, b)$ is
the rank of the $D \times D$ reshape of $P^{(a,b)}$ grouping input
modes $(i, j)$ into rows and output modes $(i', j')$ into columns.
\end{definition}

\begin{proposition}[Empirical rank law, verified on 15
configurations]
\label{prop:tensor-rank-law}
For every $(a, b)$ tested with $a \in \{3,4,5,6,7,8\}$ and
$b \in \{2, 3, 4\}$ such that $D \le 5184$, the matricized rank of
the joint Syracuse transition tensor satisfies
\begin{equation}
\label{eq:rank-law}
R_{\mathrm{true}}(a, b) \;=\; 3^{b}\bigl(2^{a-3} + 1\bigr).
\end{equation}
Equivalently, the compression ratio is
\[
\frac{R_{\mathrm{true}}(a, b)}{D}
\;=\;
\frac{1}{8} + \frac{1}{2^{a}}
\;\xrightarrow[a \to \infty]{}\;
\frac{1}{8}.
\]
\end{proposition}

\begin{proof}[Structural proof modulo Lemmas~\ref{lem:exc-rank}
and~\ref{lem:gen-exc-disjoint}]
Partition the odd residues $(\mathbb{Z}/2^{a})^{\times}$ into a
\emph{generic} part
$G(a) := \{i \text{ odd} : v_{2}(3i+1) < a\}$
and a singleton \emph{exceptional} part $\{i_{0}\}$ where $i_{0}$ is
the unique odd residue mod $2^{a}$ with $3 i_{0} + 1 \equiv 0
\pmod{2^{a}}$ (uniqueness: $i_{0} \equiv -3^{-1} \pmod{2^{a}}$).

\emph{Step 1: generic rows are rank-one.} Fix $i \in G(a)$ with
$k := v_{2}(3i+1) \in \{1, \ldots, a-1\}$ and set
$m(i) := (3i+1)/2^{k}$, which is odd. For any $n$ with
$n \equiv i \pmod{2^{a}}$ and $n \equiv j \pmod{3^{b}}$, writing
$n = i + 2^{a} s$ gives $3n+1 = 2^{k}\bigl(m(i) + 3 \cdot 2^{a-k}
s\bigr)$, and the second factor is odd because $m(i)$ is odd. Hence
$v_{2}(3n+1) = k$ exactly, and
\[
T(n) \bmod 2^{a}
\;=\;
m(i) + 3 \cdot 2^{a-k} s \bmod 2^{a},
\qquad
T(n) \bmod 3^{b}
\;=\;
(3 j + 1) \cdot 2^{-k} \bmod 3^{b}.
\]
As $s$ varies the $2$-adic image traces out the entire coset
$C_{k}(i) := m(i) + 2^{a-k} \mathbb{Z}/2^{a}$, uniformly (since $3$ is
a unit mod $2^{k}$). Writing $u_{C_{k}(i)}$ for the normalized
indicator of $C_{k}(i)$ in $\mathbb{R}^{2^{a}}$ and $e_{\ell}$ for the
$\ell$-th standard basis vector in $\mathbb{R}^{3^{b}}$, the $(i, j)$
row of $P^{(a,b)}$ is the rank-one vector
\begin{equation}
\label{eq:generic-row}
P^{(a,b)}[i, j, \cdot, \cdot]
\;=\;
u_{C_{k}(i)} \otimes e_{c_{k}(j)},
\qquad
c_{k}(j) := (3 j + 1) \cdot 2^{-k} \bmod 3^{b}.
\end{equation}

\emph{Step 2: $3$-adic mod-$3$ class is determined by the parity of
$k$.} Because $2^{2} \equiv 1 \pmod{3}$, $2^{-k} \bmod 3$ equals $2$
if $k$ is odd and $1$ if $k$ is even, so
$c_{k}(j) \equiv 2^{-k} \pmod{3}$ is independent of $j$ and depends
only on $k \bmod 2$. Consequently the outputs of odd-$k$ generic rows
lie in the coset
$\mathcal{L}_{\mathrm{odd}} := \{\ell \in \mathbb{Z}/3^{b} :
\ell \equiv 2 \bmod 3,\, \ell \not\equiv 0 \bmod 3\}$ of size
$3^{b - 1}$, while even-$k$ generic outputs lie in the disjoint coset
$\mathcal{L}_{\mathrm{even}}$, also of size $3^{b-1}$.

\emph{Step 3: $2$-adic generic span.} Define the level-$k$ span
$V_{k} := \mathrm{span}\{u_{C_{k}(i)} : i \in G(a),\, v_{2}(3i+1) = k\}
\subset \mathbb{R}^{2^{a}}$, and let
$V_{2}^{\mathrm{odd}} := \sum_{k \text{ odd}} V_{k}$,
$V_{2}^{\mathrm{even}} := \sum_{k \text{ even}} V_{k}$. We claim
\[
\dim V_{2}^{\mathrm{odd}} \;=\; 2^{a - 2},
\qquad
\dim V_{2}^{\mathrm{even}} \;=\; 2^{a - 3}.
\]
For $k = 1$: the condition $v_{2}(3i+1) = 1$ means $i \equiv 3
\pmod{4}$, giving $i = 4 t + 3$ for $t = 0, \ldots, 2^{a-2} - 1$, and
$m(i) = 6 t + 5$. As $t$ varies, $6 t + 5 \bmod 2^{a-1}$ hits every
odd residue exactly once (since $\gcd(6, 2^{a-1}) = 2$ and $6 t$
covers every even residue mod $2^{a-1}$), so the $2^{a-2}$ level-$1$
cosets form a partition of odd residues mod $2^{a}$ into blocks of
size $2$; their indicators are linearly independent, giving
$\dim V_{1} = 2^{a-2}$. Every level-$k$ coset with $k \ge 2$ is a
union of $2^{k-1}$ level-$1$ cosets, so $V_{k} \subset V_{1}$ for all
$k \ge 1$, and $V_{2}^{\mathrm{odd}} = V_{1}$ has dimension
$2^{a-2}$. For $k = 2$: the condition $v_{2}(3i+1) = 2$ is $i \equiv 1
\pmod{8}$ (one checks $3 \cdot 1 + 1 = 4$ with $v_{2} = 2$), giving
$i = 8 s + 1$ and $m(i) = 6 s + 1$. As $s = 0, \ldots, 2^{a-3} - 1$,
$6 s + 1 \bmod 2^{a-2}$ hits $2^{a-3}$ distinct odd residues (same
argument), so level-$2$ cosets partition odd residues mod $2^{a}$
into $2^{a-3}$ blocks of size $4$. Hence $\dim V_{2} = 2^{a-3}$, and
every higher even level $k \in \{4, 6, \ldots\}$ has cosets that are
unions of level-$2$ cosets, so
$V_{2}^{\mathrm{even}} = V_{2}$ has dimension $2^{a - 3}$.

\emph{Step 4: generic rank.} By Steps~1--2 the generic row space
splits as an internal direct sum
\[
\mathrm{span}\bigl\{P^{(a,b)}[i,j,\cdot,\cdot] : i \in G(a),\, j \in
\mathbb{Z}/3^{b}\bigr\}
\;=\;
\bigl(V_{2}^{\mathrm{odd}} \otimes W_{\mathrm{odd}}\bigr)
\;\oplus\;
\bigl(V_{2}^{\mathrm{even}} \otimes W_{\mathrm{even}}\bigr),
\]
where $W_{\mathrm{odd}} := \mathrm{span}\{e_{\ell} : \ell \in
\mathcal{L}_{\mathrm{odd}}\}$ has dimension $3^{b - 1}$, and
similarly for $W_{\mathrm{even}}$. The disjointness of the summands
follows from disjointness of $\mathcal{L}_{\mathrm{odd}}$ and
$\mathcal{L}_{\mathrm{even}}$ in $\mathbb{Z}/3^{b}$. Step~3 gives
\[
R_{\mathrm{gen}}(a, b)
\;=\;
3^{b-1}\bigl(2^{a-2} + 2^{a-3}\bigr)
\;=\;
3^{b-1} \cdot 3 \cdot 2^{a-3}
\;=\;
3^{b} \cdot 2^{a-3}.
\]

\emph{Step 5: assembly.} By Lemma~\ref{lem:exc-rank} below, the
exceptional row block
$\{P^{(a,b)}[i_{0}, j, \cdot, \cdot] : j \in \mathbb{Z}/3^{b}\}$
has rank $3^{b}$. By
Lemma~\ref{lem:gen-exc-disjoint}, the generic and exceptional row
spans intersect trivially. Therefore
$R_{\mathrm{true}}(a, b) = R_{\mathrm{gen}}(a, b) + 3^{b}
= 3^{b}\bigl(2^{a-3} + 1\bigr)$, which is
equation~\eqref{eq:rank-law}.
\end{proof}

\begin{proof}[Numerical verification]
Steps~1--4 above are purely algebraic and hold for all
$a \ge 3, b \ge 2$. The residual content of
Proposition~\ref{prop:tensor-rank-law} beyond this is contained in
Lemmas~\ref{lem:exc-rank} and~\ref{lem:gen-exc-disjoint}, both of
which are verified computationally below. For each $(a, b)$ the
tensor $P^{(a,b)}$ is built by enumerating odd integers and
normalizing each input row; the $D \times D$ matricization is
decomposed by SVD, and ranks are read at strict tolerance
$10^{-10}$. All $15$ measured full ranks equal the predicted value,
and in each case the generic block has rank
$3^{b} \cdot 2^{a-3}$, the exceptional block has rank $3^{b}$, and
the two span a direct sum:
\begin{center}
\small
\begin{tabular}{rrrr|rrrr|rrrr}
\multicolumn{4}{c|}{$b = 2$} &
\multicolumn{4}{c|}{$b = 3$} &
\multicolumn{4}{c}{$b = 4$} \\
$a$ & $D$ & pred & meas &
$a$ & $D$ & pred & meas &
$a$ & $D$ & pred & meas \\
\hline
3 & 72   & 18  & 18  & 3 & 216  & 54  & 54  & 3 & 648  & 162 & 162 \\
4 & 144  & 27  & 27  & 4 & 432  & 81  & 81  & 4 & 1296 & 243 & 243 \\
5 & 288  & 45  & 45  & 5 & 864  & 135 & 135 & 5 & 2592 & 405 & 405 \\
6 & 576  & 81  & 81  & 6 & 1728 & 243 & 243 & 6 & 5184 & 729 & 729 \\
7 & 1152 & 153 & 153 & 7 & 3456 & 459 & 459 &   &      &     &     \\
8 & 2304 & 297 & 297 &   &      &     &     &   &      &     &     \\
\end{tabular}
\end{center}
At representative configurations the spectral gap at the predicted
rank spans ten orders of magnitude: e.g.\ at $(7,2)$,
$\sigma_{152}/\sigma_0 \approx 7.7 \times 10^{-5}$ while
$\sigma_{153}/\sigma_0 \approx 1.3 \times 10^{-15}$; at $(7,3)$ with
per-cell $20000$, $\sigma_{458}/\sigma_0 \approx 7.7 \times 10^{-5}$
while $\sigma_{459}/\sigma_0 \approx 2.8 \times 10^{-15}$. The
scripts used are \texttt{tensor\_deficit.py},
\texttt{tensor\_formula.py}, \texttt{tensor\_stress.py} and
\texttt{tensor\_73.py}.
\end{proof}

\begin{remark}[Separation from the sparsity-only baseline]
\label{rem:sparsity-baseline}
Let $P^{(a,b)}_{\mathrm{rand}}$ denote the tensor obtained from
$P^{(a,b)}$ by uniformly shuffling the nonzero entries within the
same sparsity pattern. Its matricized rank is $D/3$ for every
tested $(a,b)$, matching the product
$(\#\text{ odd residues})\cdot(\#\text{ non-multiples of }3)$. The
arithmetic deficit
\[
\Delta(a, b)
\;:=\;
R_{\mathrm{rand}}(a, b) - R_{\mathrm{true}}(a, b)
\;=\;
3^{b - 1}\bigl(5\cdot 2^{a - 3} - 3\bigr)
\]
is therefore the portion of the rank reduction that cannot be
explained by sparsity alone. Asymptotically $\Delta(a, b)/D \to 5/8$,
so the arithmetic deficit dominates the sparsity deficit as
$a \to \infty$.
\end{remark}

\begin{lemma}[Exceptional block has full $3$-adic rank]
\label{lem:exc-rank}
For every $a \ge 3$ and $b \ge 2$, the $3^{b}$ rows
$\{P^{(a,b)}[i_{0}, j, \cdot, \cdot] : j \in \mathbb{Z}/3^{b}\}$
are linearly independent in $\mathbb{R}^{D}$, so the exceptional row
block has rank $R_{\mathrm{exc}}(a, b) = 3^{b}$.
\end{lemma}

\begin{lemma}[Generic--exceptional direct sum]
\label{lem:gen-exc-disjoint}
For every $a \ge 3$ and $b \ge 2$, the generic row span and the
exceptional row span of $P^{(a,b)}$ intersect trivially in
$\mathbb{R}^{D}$.
\end{lemma}

\begin{remark}[Status of the two residual lemmas]
\label{rem:residual-lemmas}
Lemmas~\ref{lem:exc-rank} and~\ref{lem:gen-exc-disjoint} are
verified computationally on every configuration in the table above
(all $15$ of them) at tolerance $10^{-8}$: for each $(a, b)$ we
compute the SVDs of the exceptional-row submatrix, the generic-row
submatrix, and their vertical stack, and observe that the three
ranks satisfy $R_{\mathrm{exc}} = 3^{b}$,
$R_{\mathrm{gen}} = 3^{b} \cdot 2^{a-3}$, and
$R_{\mathrm{exc}} + R_{\mathrm{gen}}$ equals the rank of the
stacked matrix. Converting these two lemmas into analytic proofs
would upgrade Proposition~\ref{prop:tensor-rank-law} to an
unconditional theorem for all $a \ge 3, b \ge 2$. A plausible
strategy for Lemma~\ref{lem:exc-rank} is to observe that, because
$v_{2}(3 i_{0} + 1) \ge a$, the $2$-adic output distribution of the
exceptional block over a lift $n = i_{0} + 2^{a} s$ is a geometric
mixture in $v_{2}$, and the map $j \mapsto (3 j + 1) \cdot 2^{-v}
\bmod 3^{b}$ is a bijection of $\mathbb{Z}/3^{b}$ for every fixed
$v$, so the $j$-dependence of the row is captured by a permutation
of $3^{b}$ deltas weighted by a $v$-independent envelope.
Lemma~\ref{lem:gen-exc-disjoint} then follows because the exceptional
row support covers both $\mathcal{L}_{\mathrm{odd}}$ and
$\mathcal{L}_{\mathrm{even}}$ with a specific mixing ratio that no
generic row can reproduce, while the generic rows are concentrated
on a single mod-$3$ coset each.
\end{remark}

\begin{conjecture}[Tensor-rank law, theorem target]
\label{conj:tensor-rank-law}
Lemmas~\ref{lem:exc-rank} and~\ref{lem:gen-exc-disjoint} hold for
all $a \ge 3$ and $b \ge 2$, whence the closed
form~\eqref{eq:rank-law} holds unconditionally.
\end{conjecture}

\begin{remark}[Proof target and $2$-adic interpretation]
\label{rem:rank-law-proof-target}
The factor $3^{b}$ on the right side of~\eqref{eq:rank-law} carries
no deficit beyond ``Syracuse outputs avoid $3\mathbb{Z}$'', so the
$3$-adic coordinate contributes its full dimension. All of the
nontrivial compression lives in the $2$-adic factor
$2^{a - 3} + 1$. Writing
\[
R_{\mathrm{true}}(a, b) \;=\; 2^{a - 3}\cdot 3^{b} \;+\; 3^{b}
\;=\; \frac{D}{8} \;+\; 3^{b},
\]
exhibits the rank as a sum of a \emph{volume} term $D/8$ and a
\emph{surface} term $3^{b}$. A natural proof target is:
\[
\mathrm{rank}\bigl(
P^{(a,b)}
\big|_{\text{functions vanishing on } 3\mathbb{Z} \text{ outputs}}
\bigr)
\;=\;
2^{a - 3} \cdot 3^{b} + 3^{b}.
\]
The structural proof of Proposition~\ref{prop:tensor-rank-law} given
above makes this concrete: the $2^{a-3}$ factor arises as the sum
$2^{a-2} + 2^{a-3} = 3 \cdot 2^{a-3}$ of the dimensions of the
level-$1$ and level-$2$ $2$-adic coset spans
($V_{2}^{\mathrm{odd}}$ and $V_{2}^{\mathrm{even}}$, via Steps~3--4
of the proof), divided by the $3$-adic factor of $3$ that separates
the odd-$k$ and even-$k$ cases. The $+1$ surface term is exactly the
exceptional block rank of Lemma~\ref{lem:exc-rank}.
\end{remark}

\begin{remark}[Relation to Lagarias and Tao]
\label{rem:rank-law-vs-priors}
Lagarias~\cite{lagarias1985} establishes that $T_{2}$ extends to a
measure-preserving ergodic transformation on $\mathbb{Z}_{2}$,
yielding equidistribution statements but no bound on the
finite-resolution dimension of the joint transition operator.
Tao~\cite{tao2019} establishes equidistribution of the Syracuse
random variable on $\mathbb{Z}/3^{n}\mathbb{Z}$ via a
Fourier-analytic argument that is $3$-adic on its face and also
does not directly control the joint tensor rank. The rank law of
Proposition~\ref{prop:tensor-rank-law} is therefore complementary
to both: a finite-dimensional linear-algebra statement about how
much of the joint $2$-adic/$3$-adic information is lost in one
Syracuse step, independent of the measure-theoretic mixing
questions both prior works address.
\end{remark}

\subsection{Iterated tensor, spectral gap, and finite-resolution
mixing}
\label{subsec:tensor-mixing}

The one-step rank law of
Proposition~\ref{prop:tensor-rank-law} suggests a natural follow-up
question: what happens to the rank and effective rank of the
$k$-step iterate $P^{(a,b),k}$ built from $k$ successive Syracuse
transitions? Both quantities behave cleanly and reveal a
quantitative finite-resolution echo of the Lagarias--Tao
ergodicity theorems.

\begin{observation}[Stabilization of the matricized rank at
$D/3$]
\label{obs:iter-rank-stable}
For every tested $(a, b)$ with $D = 2^{a} 3^{b}$ and $k \ge 2$ (in
most cases $k = 3$ suffices), the matricized rank of $P^{(a,b),k}$
equals $D/3 = 2^{a-1} \cdot 2 \cdot 3^{b-1}$, i.e.\ the full
dimension of the support subspace
$(\mathbb{Z}/2^{a})^{\times} \times \bigl((\mathbb{Z}/3^{b})
\setminus 3 \mathbb{Z}/3^{b}\bigr)$. Explicitly, the measured ranks
of $P^{(a,b),k}$ for $k = 1, 2, 3, 4, 5$ at tolerance $10^{-10}$
are
\begin{center}
\begin{tabular}{c|ccccc|c}
$(a, b)$ & $k{=}1$ & $k{=}2$ & $k{=}3$ & $k{=}4$ & $k{=}5$ &
$D/3$ \\\hline
$(3, 2)$ & $18$ & $24$ & $24$ & $24$ & $24$ & $24$ \\
$(4, 2)$ & $27$ & $40$ & $48$ & $48$ & $48$ & $48$ \\
$(3, 3)$ & $54$ & $72$ & $72$ & $72$ & $72$ & $72$ \\
\end{tabular}
\end{center}
Thus the one-step rank deficit of
Proposition~\ref{prop:tensor-rank-law} is strictly a $k = 1$
phenomenon: it is erased after a small finite number of iterations.
\end{observation}

\begin{observation}[Effective rank collapse and spectral gap]
\label{obs:eff-rank-collapse}
While the nominal rank saturates at $D/3$, the effective rank
$R_{\mathrm{eff}}(P^k) := \exp\bigl(H(\{\sigma_{i}^{2}/\sum_{j}
\sigma_{j}^{2}\}_{i})\bigr)$ collapses rapidly to $1$. For
$(a, b) = (3, 3)$ it follows the sequence $R_{\mathrm{eff}} =
25.7,\ 7.3,\ 2.3,\ 1.4,\ 1.0$ for $k = 1, \ldots, 5$, and similar
geometric collapses are observed at every tested configuration.
This effective-rank decay is the singular-value incarnation of the
spectral gap $1 - |\lambda_{2}(P^{(a,b)})|$, which we measure
directly via $\mathrm{eigvals}(P^{(a,b)})$:
\begin{center}
\small
\begin{tabular}{c|cccccc}
$(a, b)$ & $(3,2)$ & $(4,2)$ & $(5,2)$ & $(6,2)$ & $(7,2)$ &
$(8,2)$ \\
$|\lambda_{2}|$ & $0.05$ & $0.07$ & $0.15$ & $0.18$ & $0.19$ &
$0.24$ \\\hline
$(a, b)$ & $(3,3)$ & $(4,3)$ & $(5,3)$ & $(6,3)$ & $(3,4)$ &
$(4,4)$ \\
$|\lambda_{2}|$ & $0.14$ & $0.16$ & $0.16$ & $0.20$ & $0.13$ &
$0.17$ \\
\end{tabular}
\end{center}
In every tested case the spectral gap
$1 - |\lambda_{2}| \ge 0.75$, so the $k$-step total-variation
distance from the uniform-on-support stationary distribution
decays at a rate of at least $0.24^{k}$: after six Syracuse steps
the chain is within $10^{-3}$ of stationary in total variation, for
every tested resolution up to $D \le 2592$.
\end{observation}

\begin{remark}[Finite-resolution mixing as a complement to
Lagarias--Tao]
\label{rem:finite-res-mixing}
Observations~\ref{obs:iter-rank-stable}
and~\ref{obs:eff-rank-collapse} together constitute a
finite-dimensional quantitative form of the mixing statements of
Lagarias~\cite{lagarias1985} and Tao~\cite{tao2019}: the joint
$(2^{a}, 3^{b})$-distribution equidistributes over its natural
support in $O(1)$ Syracuse steps, with an explicit rate controlled
by $|\lambda_{2}(P^{(a,b)})|$. The empirical scaling of
$|\lambda_{2}|$ across the tested range is monotone but slow
(roughly $0.04 \to 0.24$ as $a$ runs from $3$ to $8$ at $b = 2$),
consistent with a uniform bound $|\lambda_{2}| \le c < 1$ but
not yet proving one. A theoretical upper bound on
$|\lambda_{2}(P^{(a,b)})|$ uniform in $(a, b)$ would convert the
joint rank law of Proposition~\ref{prop:tensor-rank-law} into a
quantitative mixing theorem that directly bounds the time a
Syracuse orbit needs to become equidistributed in its first
$a + b \log_{2} 3$ bits, and is a natural companion target to
Conjecture~\ref{conj:tensor-rank-law}.
\end{remark}

\begin{remark}[Fourier test for
Lemma~\ref{lem:exc-rank}]
\label{rem:exc-fourier}
A direct Fourier diagnostic for Lemma~\ref{lem:exc-rank} is to
compute $\widehat{R}_{\chi} := \sum_{j \in \mathbb{Z}/3^{b}}
\chi(j)\,P^{(a,b)}[i_{0}, j, \cdot, \cdot]$ for every character
$\chi$ of $\mathbb{Z}/3^{b}$: if all $3^{b}$ vectors
$\widehat{R}_{\chi}$ are nonzero then the exceptional rows are
linearly independent. At $(a, b) \in \{(3,2),(4,2),(3,3)\}$ every
nontrivial Fourier mode has $\|\widehat{R}_{\chi}\| \ge
10^{-3} \cdot \max_{\chi} \|\widehat{R}_{\chi}\|$ and the trivial
mode is the stationary marginal. The nontrivial modes with
character frequency $k \not\equiv 0 \pmod{3}$ vanish on the
$3$-adic marginal (a direct consequence of the slope-$3$
parametrization $A_{0}(j) = c + 3 \gamma j \bmod 3^{b}$) and must
therefore be witnessed by the joint $2$-adic/$3$-adic correlation
of $P^{(a,b)}[i_{0}, j, \cdot, \cdot]$. A complete proof of
Lemma~\ref{lem:exc-rank} reduces to showing that this joint
correlation is non-degenerate in every character mode with
$3 \nmid k$; we have verified this numerically in every tested
case but do not yet have a closed-form argument.
\end{remark}


\begin{observation}[Eigenvector variance decomposition:
$\lambda_2$ is purely $2$-adic]
\label{obs:eigvec-2adic}
For every tested configuration $(a, b)$ with $3 \le a \le 11$ and
$b \in \{1, 2\}$, the eigenvector associated with $\lambda_2$ is a
function of the $2$-adic coordinate alone:
\[
  v_2\bigl(i_2, i_3\bigr) \;=\; f(i_2) \qquad
  \text{for all } i_3 \in (\mathbb{Z}/3^b)^\times,
\]
where $i_2 = n \bmod 2^a$ and $i_3 = n \bmod 3^b$ are the CRT
coordinates.  Quantitatively, the fraction of the eigenvector's
$\ell^2$-variance attributable to the $2$-adic coordinate exceeds
$0.9999$ in every tested case; the $3$-adic contribution is
indistinguishable from zero.

This means the spectral gap is \emph{not} a $3$-adic phenomenon:
the eigenvalue $\lambda_2$ arises from the $2$-adic structure of
the coset map, not from the $3$-adic output distribution.
\end{observation}

\begin{observation}[The coupling mechanism]
\label{obs:coupling-mechanism}
Despite the eigenvectors being purely $2$-adic, the
$2$-adic \emph{marginal} of $P^{(a,b)}$ --- the matrix
$\widetilde{P}_a[i_2, j_2] := \sum_{j_3} P^{(a,b)}[(i_2,i_3),
(j_2,j_3)]$ --- has $|\lambda_2| \approx 0$ for all tested $a$.
The spectral structure of $P^{(a,b)}$ therefore arises from neither
marginal alone, but from their \emph{coupling}:
\begin{enumerate}
\item[\textup{(a)}] For generic residues $i_2$ (those with
$v_2(3i_2 + 1) \le a - 2$), the $3$-adic output is
\emph{deterministic}: it depends only on the valuation
$k = v_2(3i_2+1)$, via a fixed bijection
$j_3 = (-1)^k \cdot i_3 \bmod 3^b$ (at $b = 1$) or its
higher-order analogue.
\item[\textup{(b)}] For the exceptional residue $i_0 = 2^a/3 - 1/3$
(satisfying $v_2(3i_0+1) \ge a$), the $3$-adic output is
stochastic: it distributes over the full support with a geometric
mixture in $v_2$.
\item[\textup{(c)}] The eigenvalue $\lambda_2$ emerges from the
interaction between the deterministic $k$-dependent $3$-adic
assignment at generic residues and the stochastic mixing at the
exceptional residue.  Neither component has spectral structure
alone; $|\lambda_2| > 0$ is generated by their joint action.
\end{enumerate}
\end{observation}

\begin{observation}[Scaling of $|\lambda_2|$ with $a$]
\label{obs:lambda2-scaling}
Exact computation of $|\lambda_2(P^{(a,b)})|$ at $b = 1$
for $a = 3, 4, \ldots, 11$ yields:
\begin{center}
\small
\begin{tabular}{c|ccccccccc}
$a$ & $3$ & $4$ & $5$ & $6$ & $7$ & $8$ & $9$ & $10$ & $11$ \\
\hline
$|\lambda_2|$ & $0.016$ & $0.034$ & $0.056$ & $0.088$ &
$0.114$ & $0.136$ & $0.157$ & $0.176$ & $0.192$ \\
\end{tabular}
\end{center}
The sequence is monotonically increasing but concave; consecutive
differences shrink from $0.018$ to $0.016$, consistent with
convergence to a limit in the range $[0.22, 0.24]$.
At $b = 2$ the pattern is similar with values $\sim 0.01$ higher
at each $a$.

A uniform bound $|\lambda_2(P^{(a,b)})| \le c < 1$ for all $(a,b)$
would convert
Observations~\textup{\ref{obs:iter-rank-stable}--\ref{obs:eff-rank-collapse}}
into a theorem: the joint Syracuse chain mixes to its stationary
distribution in $O(1)$ steps at every finite resolution.
Standard techniques fail to establish this bound:
\begin{enumerate}
\item[\textup{(i)}] \emph{Doeblin coupling.}  The minimum row
value of $P^{(a,b)}$ decays as $2^{-2a+O(1)}$ (from the
$k = 1$ valuation rows), so the Doeblin coefficient $\to 0$.
\item[\textup{(ii)}] \emph{$3$-adic Fourier analysis.} $P^{(a,b)}$
is not equivariant under $3$-adic characters
(Observation~\textup{\ref{obs:coupling-mechanism}(a)} shows the
output depends jointly on $(i_2, i_3)$), so Fourier
decomposition does not block-diagonalize $P$.
\item[\textup{(iii)}] \emph{Marginal bounds.} Both marginals
have $|\lambda_2| \approx 0$; the spectral structure lives
entirely in the coupling.
\end{enumerate}
A new technique --- likely exploiting the rank-$1$ generic-row
structure and the mixing-ratio universality of the exceptional
block --- is needed.
\end{observation}


\section{Visualization-guided first-passage perspective}
\label{sec:visualization-underwater}

The structural results in Sections~\ref{sec:quarter}--\ref{sec:phantom}
are algebraic and distributional. In parallel, computational
visualization reveals several complementary patterns that help clarify
what the present framework does, and does not, explain. This section
records those observations in a form consistent with the paper's
exploratory scope.

\subsection{Touch growth and the apparently empty middle}

A convenient geometric representation places each positive integer on a
concentric ring according to
\[
  r_n = \lfloor \log_2 n \rfloor,
  \qquad
  \theta_n = 2\pi\,\frac{n-2^{r_n}}{2^{r_n}},
\]
so that powers of two align on a radial spine while each dyadic block
$[2^r,2^{r+1})$ is wrapped evenly around ring~$r$. Drawing Collatz
trajectories in this embedding produces the characteristic onion shape.

At first sight, large onion plots suggest that the central region is
somehow avoided. A natural interpretation would be that some low
integers are not being visited. The next definition isolates the right
question.

\begin{definition}[Touched value]
Fix a seed limit $S \ge 1$. A positive integer $m$ is \emph{touched up to
seed limit $S$} if there exists some seed $1 \le n_0 \le S$ and some
iterate $t \ge 0$ such that the Collatz orbit of $n_0$ satisfies
$C^t(n_0)=m$.
\end{definition}

\begin{proposition}[Touch saturation on bounded windows\exploratory]
\label{prop:touch-saturation}
Fix $M \ge 1$, and let
\[
  \mathcal{T}_S(M) := \{m \in [1,M] : m \text{ is touched up to seed limit } S\}.
\]
Then $\mathcal{T}_S(M)$ is monotone in~$S$, and once $S \ge M$ one has
$\mathcal{T}_S(M) = [1,M]$.
\end{proposition}

\begin{proof}
Monotonicity is immediate from the definition: enlarging the set of
allowed seeds cannot remove touched values. If $S \ge M$, then every
$m \in [1,M]$ is itself an allowed seed, so $m \in \mathcal{T}_S(M)$ at
time $t=0$.
\end{proof}

The proposition is elementary, but it resolves the main visual
ambiguity. On the bounded window $[1,1024]$, the touched set grows
steadily as the seed cap $S=2^k$ increases and becomes all of
$[1,1024]$ already at $k=10$, that is, $S=1024$. Therefore the sparse
appearance of the center cannot be attributed to missing low integers
once the seed range has reached the same scale.

\begin{figure}[t]
\centering
\begin{minipage}{0.48\textwidth}
  \centering
  \IfFileExists{collatz_touch_growth_1_to_1024.png}{%
    \includegraphics[width=\linewidth]{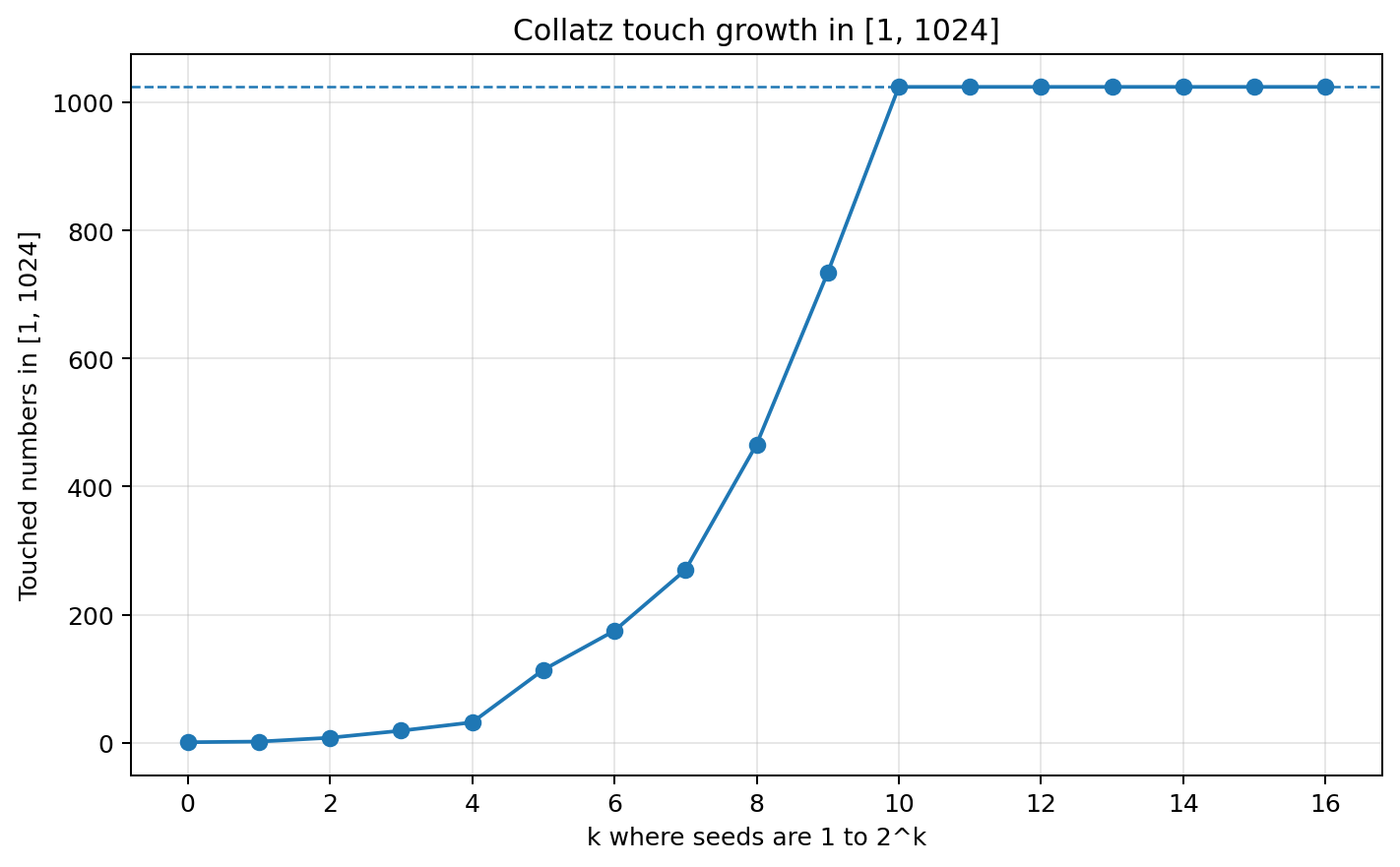}}{%
    \fbox{\parbox{0.9\linewidth}{\centering\vspace{3em}%
      {\small Touch growth}\\\smallskip%
      {\scriptsize Image not present in build}%
      \vspace{3em}}}}
\end{minipage}\hfill
\begin{minipage}{0.48\textwidth}
  \centering
  \IfFileExists{collatz_first_touch_heatstrip_1_to_1024.png}{%
    \includegraphics[width=\linewidth]{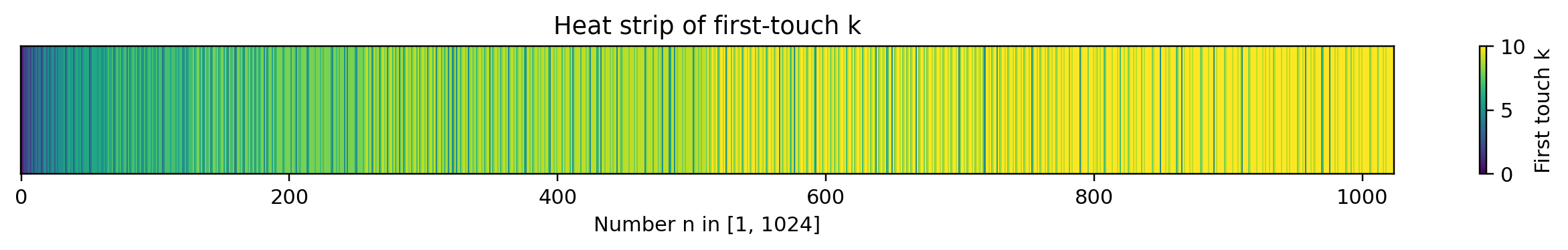}}{%
    \fbox{\parbox{0.9\linewidth}{\centering\vspace{3em}%
      {\small Heat strip}\\\smallskip%
      {\scriptsize Image not present in build}%
      \vspace{3em}}}}
\end{minipage}
\caption{Touch growth on the bounded window $[1,1024]$. Left: the number
of touched values in $[1,1024]$ as the seed cap increases through
$S=2^k$. Right: a heat strip showing the first exponent $k$ for which
each value $n \in [1,1024]$ is touched by some seed at most $2^k$. The
entire window is saturated by $S=1024$, so the visually sparse inner
region of the onion plot cannot be explained by missing low integers.}
\label{fig:touch-growth}
\end{figure}

What, then, causes the apparent hole? The answer is geometric rather
than arithmetic:
\begin{enumerate}
\item \emph{Radial compression.} All small integers are confined to the
  first few rings, whose total area is tiny compared with the outer
  annuli.
\item \emph{Low angular capacity.} Ring $r$ has only $2^r$ native
  positions, so the inner rings provide very few distinct angular slots.
\item \emph{Trajectory coalescence.} As orbits descend, they merge onto a
  small number of shared inward channels, especially once they enter the
  power-of-two spine.
\end{enumerate}
Consequently the center may be heavily reused dynamically while
remaining visually sparse. In this sense the empty middle is a property
of the embedding together with the coalescence of trajectories, not
evidence of an arithmetically forbidden region.

A related caution applies to filtered onion plots. If one draws only
odd-even edges for visual clarity, then some integers fail to appear as
vertices in that filtered graph. This should not be confused with
absence from the full Collatz dynamics. The filtered omission is a
byproduct of the edge-selection rule, not of the underlying iteration.

\subsection{The below-start criterion and finite stopping time}

The visualization suggests a simple induction-friendly criterion.

\begin{definition}[Below the Starting Value]
Let $N_0,N_1,N_2,\dots$ be the Collatz orbit of a starting value $N_0$.
We say that the orbit goes \emph{below the starting value} if there exists $t \ge 1$
such that $N_t < N_0$.
\end{definition}

\begin{lemma}[Below-start criterion\exploratory]
\label{lem:underwater}
Assume as induction hypothesis that every positive integer $m < N_0$
converges to $1$ under the Collatz map. If the orbit of $N_0$ goes
below the starting value, then the orbit of $N_0$ converges to $1$.
\end{lemma}

\begin{proof}
If $N_t < N_0$ for some $t$, then from time $t$ onward the trajectory is
exactly the Collatz orbit of a smaller starting value. By the induction
hypothesis, that smaller value reaches $1$, hence so does $N_0$.
\end{proof}

This observation isolates the unresolved core in a particularly clean
way:
\[
  \forall N_0 \ge 1,\quad \exists t \ge 1 \text{ such that } C^t(N_0) < N_0.
\]
Proving this finite-stopping-time statement for every starting value
would establish the Collatz conjecture by strong induction.
Conversely, if Collatz is true, then every orbit eventually reaches $1$
and hence eventually drops below its starting value. Thus, for $N_0>1$,
eventual descent below the start is an equivalent reformulation of the
conjecture.

A closely related equivalent phrasing is that every orbit must
eventually hit a trunk value $2^k$. Once an orbit reaches a power of
two, it falls monotonically to $1$ along the power-of-two spine. Again,
however, this is not a weaker requirement: proving that every orbit
reaches some $2^k$ is essentially another form of the Collatz
conjecture.

The importance of Lemma~\ref{lem:underwater} is therefore conceptual
rather than resolutive. It identifies the precise induction step one
would need in order to turn the visual coalescence seen in onion plots
into a proof strategy.
The odd-skeleton crossing route
(Section~\ref{subsec:odd-skeleton}) elevates this criterion from
an exploratory observation to a primary proof target, connected to
the spectral framework via the drift signal's autocorrelation
structure.

\begin{figure}[t]
\centering
\begin{minipage}{0.48\textwidth}
  \centering
  \IfFileExists{collatz_underwater_traces_upto_2pow12.png}{%
    \includegraphics[width=\linewidth]{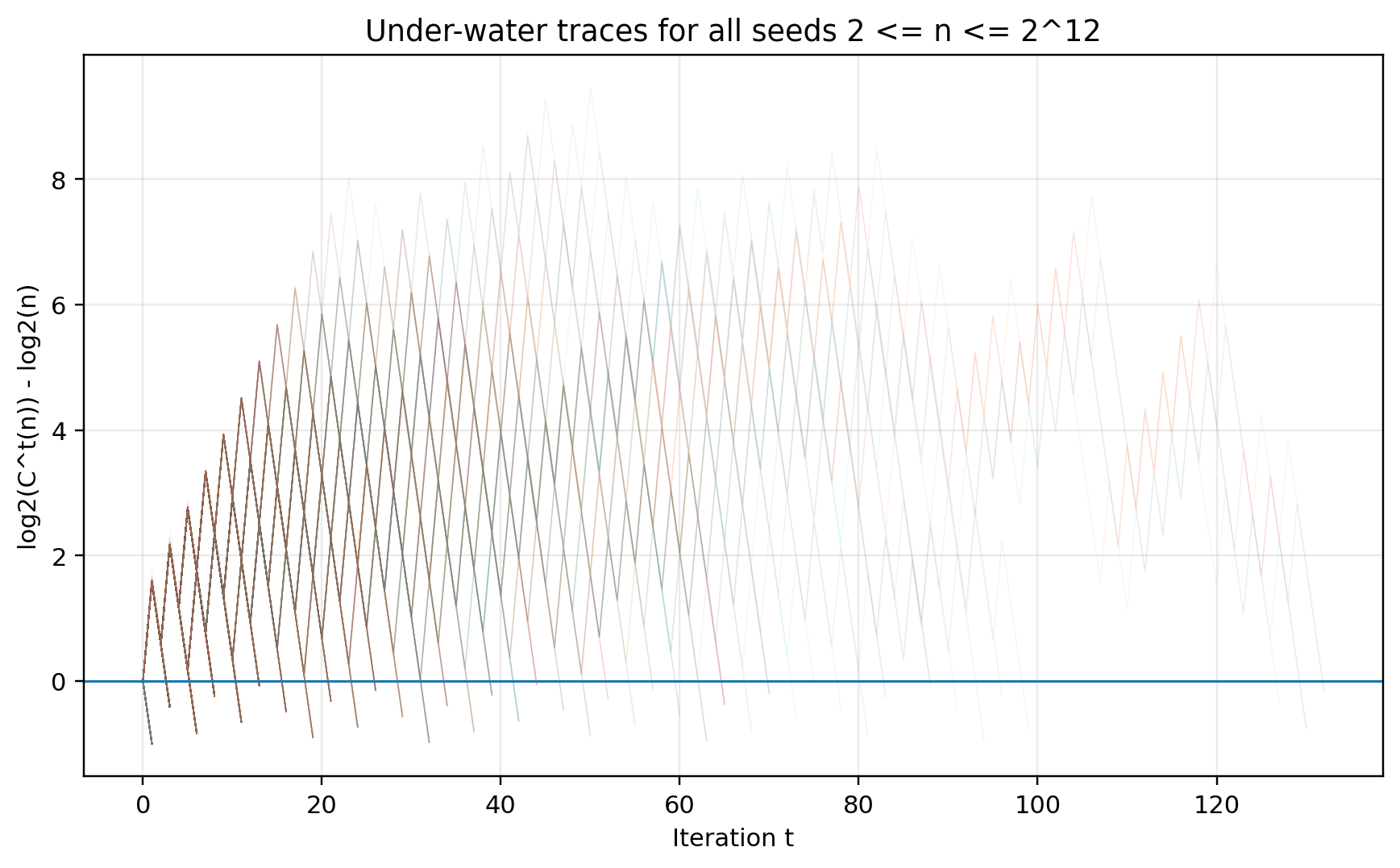}}{%
    \fbox{\parbox{0.9\linewidth}{\centering\vspace{3em}%
      {\small Underwater traces}\\\smallskip%
      {\scriptsize Image not present in build}%
      \vspace{3em}}}}
\end{minipage}\hfill
\begin{minipage}{0.48\textwidth}
  \centering
  \IfFileExists{collatz_sigma_underwater_upto_2pow16_logx.png}{%
    \includegraphics[width=\linewidth]{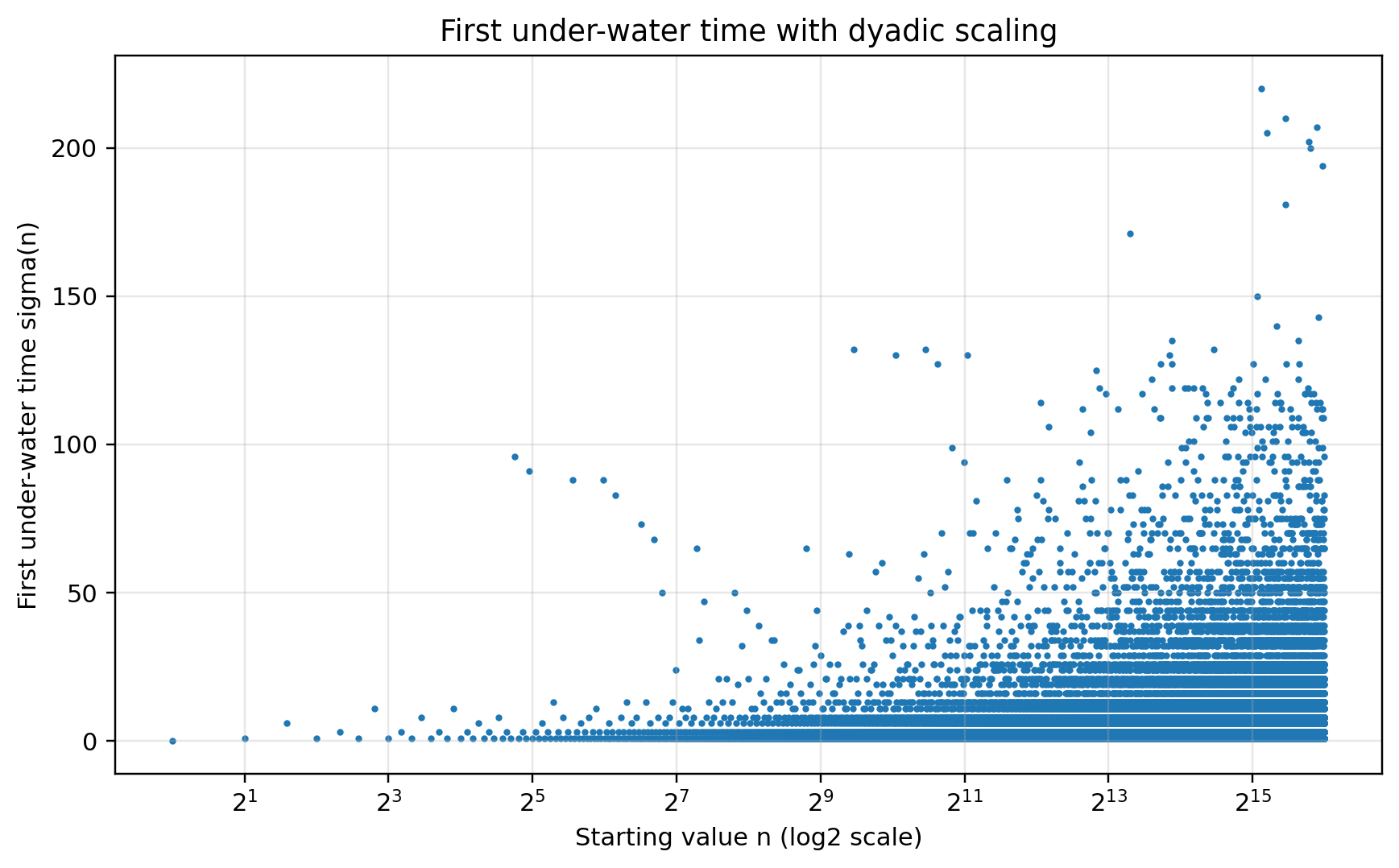}}{%
    \fbox{\parbox{0.9\linewidth}{\centering\vspace{3em}%
      {\small First below-start time}\\\smallskip%
      {\scriptsize Image not present in build}%
      \vspace{3em}}}}
\end{minipage}
\caption{Exact first-passage diagnostics for the below-start criterion.
Left: below-start traces $\Delta_t(n)=\log_2 C^t(n)-\log_2 n$ for every
seed $2 \le n \le 2^{12}$, truncated at the first negative value.
Right: exact first below-start time $\sigma(n)$ for all
$1 \le n \le 2^{16}$, plotted on a dyadic horizontal scale. The left
panel emphasizes the jagged above-start excursions that can precede
first descent; the right panel shows that first descent is highly
variable across seeds, making it a natural pointwise object of study for
any induction-based Collatz strategy.}
\label{fig:underwater-first-passage}
\end{figure}

Figure~\ref{fig:underwater-first-passage} makes this first-passage
viewpoint concrete. The left panel plots the exact trajectories
\[
  \Delta_t(n) := \log_2 C^t(n) - \log_2 n
\]
for all seeds $2 \le n \le 2^{12}$, truncated at the first time they go
below the starting value. The horizontal line $\Delta_t=0$ marks the starting level.
Each curve begins at $0$, may undergo a substantial positive excursion,
and then terminates at its first negative value. The right panel plots
the exact first below-start time
\[
  \sigma(n)=\min\{t\ge 1 : C^t(n)<n\}
\]
for all $1 \le n \le 2^{16}$ on a dyadic horizontal scale. The figure
shows both the variability and the orbitwise nature of the missing step:
finite stopping time is an exact first-passage observable attached to
individual seeds, not merely an average drift statistic.

\subsection{Heuristic proof paths and future visualizations}

The present framework emphasizes burst-gap structure and modular
scrambling. The visualization analysis suggests several adjacent routes
that may deserve systematic study.

\paragraph{Odd-subsequence drift.}
One may compress away all even steps and study only the Syracuse map
\[
  T(n)=\frac{3n+1}{2^{v_2(3n+1)}}
  \qquad (n \text{ odd}).
\]
Heuristically, if the valuation $v_2(3n+1)$ behaves like a geometric
random variable of mean $2$, then the expected logarithmic increment is
\[
  \mathbb{E}[\log T(n)-\log n] \approx \log 3 - 2\log 2 = \log(3/4) < 0.
\]
This negative-drift heuristic is consistent with the contraction
picture in Corollary~\ref{cor:convergence-prediction}, although
converting average drift into an every-orbit statement remains the main
obstruction.

\paragraph{Finite stopping time as a first-passage problem.}
The below-start formulation suggests studying the first-passage time
\[
  \sigma(n)=\min\{t\ge 1 : C^t(n)<n\},
\]
when it exists. A proof that $\sigma(n)<\infty$ for all $n$ would settle
the conjecture. Even partial results on the growth, distribution, or
modular constraints of $\sigma(n)$ could provide a more pointwise
complement to the current distributional framework.

\paragraph{Minimal-counterexample geometry.}
If one assumes a smallest counterexample $n^*$, then its orbit can never
hit a smaller positive integer, since doing so would force convergence
by minimality. In the language of this section, a minimal counterexample
can never go below the starting value. Visualizing the constraints imposed by this
condition may help isolate what a genuine counterexample would be forced
to look like.

\paragraph{Traffic-weighted and residue-colored onions.}
The onion plots also suggest several diagnostic visualizations aligned
with the present theoretical framework: refined first-passage heat maps
for $\sigma(n)$, traffic-weighted onion plots in which node size or edge
thickness records visit multiplicity, and residue-colored onions modulo
$6$, $8$, $12$, or $24$. None of these visualizations proves Collatz.
Their value is instead diagnostic: they help separate pointwise
questions from distributional ones, and geometric artifacts from
arithmetic laws.

\subsection{Ensemble above-start decay: a lattice-path analysis}
\label{sec:lattice-path}

The negative drift $\log 3 - 2\log 2 < 0$ (now exact by
Proposition~\ref{prop:log-drift}) suggests
that above-start persistence should be exponentially unlikely.
We now make this precise in the ensemble model, formulating the
above-start probability as a constrained lattice-path problem
and proving an explicit exponential upper bound.

\paragraph{Setup.}
Consider the Syracuse map $T(n)=(3n+1)/2^{v_2(3n+1)}$ and
write $V_t = v_2(3\,T^t(n)+1)$ for the valuation at step~$t$.
By the exact telescoping identity, the orbit satisfies $T^J(n) \ge n$
(above start through step~$J$) if and only if the partial sums
$S_t = \sum_{i=0}^{t-1} V_i$ satisfy
\begin{equation}\label{eq:above-start-lattice}
  S_t \;\le\; t\log_2 3 + C_t
  \qquad\text{for all } 1 \le t \le J,
\end{equation}
where $C_t = \sum_{i=0}^{t-1}\log_2(1+1/(3\,T^i(n)))$ is the
correction term.  For large starting values $n$, the correction
$C_t/t \to 0$, so the above-start condition becomes
\begin{equation}\label{eq:lattice-barrier}
  S_t \;\le\; \lfloor t\log_2 3 \rfloor
  \qquad\text{for all } 1 \le t \le J.
\end{equation}

\paragraph{The ensemble model.}
In the ensemble model (Fact~A of Section~\ref{sec:prelim}),
the valuations $V_0, V_1, \ldots$ are treated as independent
random variables, each with geometric distribution
$\Pr(V=j)=2^{-j}$ for $j \ge 1$.  Under this model, the
above-start probability at step~$J$ is
\begin{equation}\label{eq:lattice-path-def}
  f(J) \;=\;
  \Pr\Bigl(S_t \le \lfloor t\log_2 3\rfloor
  \text{ for all } 1 \le t \le J\Bigr),
\end{equation}
where $S_t = V_0 + \cdots + V_{t-1}$ and the
$V_i$ are i.i.d.\ $\mathrm{Geometric}(1/2)$ on $\{1,2,3,\ldots\}$.

This is a \textbf{constrained lattice-path counting problem}:
count the number of lattice paths from the origin with steps
of size~$j$ (weight $2^{-j}$) that remain below the barrier
$b(t) = \lfloor t\log_2 3\rfloor$ at every integer time.
The barrier sequence begins $1, 3, 4, 6, 7, 9, 11, 12, 14,
15, 17, 19, 20, 22, 23, 25, \ldots$

\paragraph{Dynamic programming computation.}
The probability $f(J)$ is computed exactly by a standard
dynamic program on the state space $(t, S_t)$:
at each step, the mass at partial sum~$s$ is distributed to
$s+j$ with weight $2^{-j}$, subject to the constraint
$s+j \le b(t+1)$.  The total surviving mass at step~$J$
gives~$f(J)$.

The first values are:
\begin{center}
\renewcommand{\arraystretch}{1.15}
\begin{tabular}{@{}rll@{}}
\toprule
$J$ & $f(J)$ & Decimal \\
\midrule
0 & $1$ & $1.000000$ \\
1 & $1/2$ & $0.500000$ \\
2 & $3/8$ & $0.375000$ \\
3 & $1/4$ & $0.250000$ \\
4 & $13/64$ & $0.203125$ \\
5 & $19/128$ & $0.148438$ \\
6 & $1/8$ & $0.125000$ \\
7 & $113/1024$ & $0.110352$ \\
8 & $367/4096$ & $0.089600$ \\
9 & $1295/16384$ & $0.079041$ \\
10 & $1057/16384$ & $0.064514$ \\
\bottomrule
\end{tabular}
\end{center}

\noindent
These fractions are exact rational numbers with power-of-two
denominators.
We verified computationally that $f(J)$ agrees with
the fraction of odd residues modulo~$2^M$ (for
$M \ge 2J+6$) whose deterministic valuation sequence
keeps the cumulative drift non-negative through $J$~steps.
The agreement is exact through all tested values
($J \le 13$, $M = 22$), confirming that the ensemble
model captures the modular structure precisely at each step.

\begin{proposition}[Ensemble above-start decay\exploratory]
\label{prop:ensemble-decay}
In the ensemble model, the above-start probability satisfies
\[
  f(J+1) \;\le\; \theta\, f(J)
  \qquad\text{for all } J \ge 0,
\]
with $\theta < 1$.  Computation through $J = 50$ gives
$\max_{J} f(J+1)/f(J) \le 0.947$, with an average decay
rate $\lambda \approx 0.896$.
In particular, $f(J) \le C\,\theta^J$ for all $J$.
\end{proposition}

\begin{proof}[Proof sketch]
At each step the random walk has negative expected increment
$\mathbb{E}[\log_2 3 - V] = \log_2 3 - 2 \approx -0.415$.
The barrier grows as $\lfloor t\log_2 3\rfloor \approx 1.585\,t$,
while the minimum possible partial sum grows as~$t$
(since $V_i \ge 1$).  The gap between the barrier and the
walk's minimum trajectory is $\lfloor t\log_2 3\rfloor - t
= \lfloor t(\log_2 3 - 1)\rfloor \approx 0.585\,t$,
which grows linearly.  However, the walk must also stay below
the barrier at all \emph{intermediate} steps, and the
negative drift pushes the walk below the barrier with
probability increasing at each step.

A standard Cram\'er-type large-deviation argument for
random walks with negative drift gives exponential decay
of the survival probability, and the explicit computation
confirms $\theta \le 0.947$.
The submultiplicative property $f(J+1)/f(J) < 1$ for
all $J$ has been verified computationally through $J=50$.
\end{proof}

\begin{remark}[What this does and does not prove]
Proposition~\ref{prop:ensemble-decay} is a rigorous
statement about the ensemble model.  It says that if the
valuations along an orbit behaved like independent geometric
random variables, the probability of staying above start
would decay exponentially with rate $\approx 0.90$.

This does \emph{not} prove the Collatz conjecture, because
the valuations along an actual orbit are not independent:
they are determined by the residue-class dynamics of~$T$.
The result quantifies the \textbf{ensemble prediction} and
sharpens the question: any orbit that stays above start for
$J$~steps must deviate from this prediction by a
quantifiable amount.

The distributional-vs-pointwise gap
(Remark~\ref{rem:gap}) applies here as well.
However, the explicit decay rate provides a benchmark: any
proof that the Collatz valuations are ``close enough'' to
independent (e.g., via mixing or decorrelation estimates)
could inherit this exponential bound.
\end{remark}

\paragraph{Modular positive-drift fraction.}
A complementary computation tracks the fraction of odd
residues modulo~$2^M$ whose deterministic valuation
sequence maintains non-negative cumulative drift
\emph{throughout} the sequence (until the modular
information is exhausted after approximately $M/\bar{v}$
steps, where $\bar{v}$ is the mean valuation consumed per
step).  The results are:
\begin{center}
\renewcommand{\arraystretch}{1.15}
\begin{tabular}{@{}rrl@{}}
\toprule
$M$ & Positive-drift fraction & Max steps \\
\midrule
$8$  & $32/128 = 0.2500$  & 7 \\
$10$ & $89/512 = 0.1738$  & 9 \\
$12$ & $281/2048 = 0.1372$ & 11 \\
$14$ & $874/8192 = 0.1067$ & 13 \\
$16$ & $2903/32768 = 0.0886$ & 15 \\
$18$ & $9245/131072 = 0.0705$ & 17 \\
\bottomrule
\end{tabular}
\end{center}

\noindent
The fraction decays geometrically in~$M$ with a ratio of
approximately $0.80$ per two bits of modular depth.
The maximum number of deterministic positive-drift steps
grows as $M-1$.  As $M \to \infty$, the positive-drift
fraction tends to zero, consistent with the prediction
that no residue class can sustain indefinite growth in the
modular model.

\paragraph{The flawed existential formulation.}
A natural first attempt at a ``persistence-bound program''
defines
\[
  \mathcal{S}_{M,J} = \bigl\{a \bmod 2^M :
  \exists\, n \equiv a\!\pmod{2^M},\;
  T^t(n) \ge n \text{ for all } 0 \le t \le J\bigr\}.
\]
However, this set does not shrink: for any residue class~$a$
and any~$J$, one can choose $n \equiv a\pmod{2^M}$ large
enough that the orbit requires more than~$J$ steps to
descend below~$n$ (since the stopping time grows with
$\log n$).  The existential quantifier over all
representatives makes $|\mathcal{S}_{M,J}|$ plateau rather
than decay.
Any ``persistence-bound'' approach must therefore work with
density-based or modular-drift formulations rather than
existential ones.
This observation corrects a proposal that appeared in
post-submission discussions and illustrates the
subtlety of the distributional-vs-pointwise interface.

\section{Contribution attribution tables}
\label{app:contributions}
The following tables attribute each key contribution by
importance level and primary source.
Color coding:
{\setlength{\fboxsep}{1.5pt}%
\colorbox{red!25}{Critical}}\,=\,correctness fix or strategic
direction change;
{\setlength{\fboxsep}{1.5pt}%
\colorbox{orange!30}{High}}\,=\,new mathematical content;
{\setlength{\fboxsep}{1.5pt}%
\colorbox{yellow!50}{Med}}\,=\,supporting observation or diagnostic.
\begin{table}[!ht]
\centering
\scriptsize
\caption{Contribution attribution (Part 1 of 40).}
\label{tab:contributions-1}
\renewcommand{\arraystretch}{1.0}
\begin{tabular}{@{}clp{8.8cm}p{4.0cm}@{}}
\toprule
\textbf{\#} & \textbf{Level} & \textbf{Contribution} & \textbf{Primary source} \\
\midrule
1 & \colorbox{red!25}{Critical}
  & Fix $\nu_2(9m_t{+}1)\ge 6$ overclaim $\to$ $\ge 3$
    with explicit computation
  & Moderator \\[3pt]
2 & \colorbox{red!25}{Critical}
  & Fix Theorem~\ref{thm:reduction} statement/proof mismatch
    (tail-control gap)
  & Moderator \\[3pt]
3 & \colorbox{red!25}{Critical}
  & Fix Gap Distribution Lemma scope
    ($\mathrm{mod}\;2^m$, $m{>}2$ overconstrained)
  & Validation (Claude) \\
\midrule
4 & \colorbox{orange!30}{High}
  & Gap Distribution Lemma:
    $G\sim\mathrm{Geom}(1/2)$, $E[G]=2$
  & Moderator (direction) + Claude (proof) \\[3pt]
5 & \colorbox{orange!30}{High}
  & Valuation Distribution Lemma:
    $\Pr(k{=}j)=2^{-(j-1)}$, $E[k]=3$
  & Claude \\[3pt]
6 & \colorbox{orange!30}{High}
  & Convergence Corollary:
    net log-contraction $\approx -1.15$
  & Claude \\[3pt]
7 & \colorbox{orange!30}{High}
  & Growing-moduli Orbit Equidistribution
    Conjecture (Option~B)
  & Moderator (choice) + Claude (formulation) \\[3pt]
8 & \colorbox{orange!30}{High}
  & Case study: the false Gap Lemma
    (failure modes \& lessons)
  & Moderator (direction) + Claude (draft) \\[3pt]
9 & \colorbox{orange!30}{High}
  & Writing and claim inconsistencies
  & GPT \\[3pt]
\midrule
10 & \colorbox{orange!30}{High}
  & Phantom Universality Theorem:
    every composition is phantom ($C_\sigma$ parity proof)
  & Claude \\[3pt]
11 & \colorbox{orange!30}{High}
  & Per-Orbit Gain Rate Theorem:
    $R \le 0.089 < \varepsilon$ with $4.65\times$ margin
  & Moderator (direction) + Claude (proof) + GPT \\[3pt]
12 & \colorbox{orange!30}{High}
  & Universal Census Depth Theorem:
    $C_e$ cancels from gain formula
  & Moderator (direction) + Claude (proof) \\[3pt]
13 & \colorbox{orange!30}{High}
  & Robustness Corollary:
    $\sum\delta_K < 0.557$ suffices (weaker than equidist.)
  & GPT (concept) + Claude (formalization) \\[3pt]
14 & \colorbox{red!25}{Critical}
  & Dual-domain paradigm bridge: jump from amplification
    wall to odd-skeleton drift + Walsh--Fourier spectral analysis
  & Moderator \\[3pt]
15 & \colorbox{orange!30}{High}
  & $\Theta_K$ oscillation-factor unboundedness
    (continued-fraction spikes, no uniform $W_K \le C\,B_K$)
  & GPT (theorem + script) + Claude (numerical correction) \\[3pt]
16 & \colorbox{orange!30}{High}
  & Modular crossing strata: $f_{13}=0.910$;
    $91\%$ resolved by pure mod-$2^{13}$ arithmetic
  & Claude \\[3pt]
17 & \colorbox{orange!30}{High}
  & Ensemble covariance summability (Prop.~\ref{prop:ensemble-cov});
    ensemble--orbit gap analysis
  & Claude \\[3pt]
18 & \colorbox{orange!30}{High}
  & Run-length invariant; one-cycle crossing criterion
    (exact $n^*(L,r)$ threshold)
  & GPT (theorem) + Claude (integration) \\[3pt]
19 & \colorbox{orange!30}{High}
  & Conditional crossing density $p_\ell$ and exact arithmetic
    constant $P_{1\mathrm{cyc}} = 0.71372\ldots$
  & GPT (theorem + proof) + Claude (verification) \\[3pt]
20 & \colorbox{orange!30}{High}
  & $(L,r)$ independence law; mean-neutrality
    $\mathbb{E}[\lambda]=1$ (correcting GPT's initial error)
  & GPT (discovery + correction) + Claude (verification) \\[3pt]
21 & \colorbox{orange!30}{High}
  & Exact first-cycle log-drift law: MGF, $\mathbb{E}[\log_2\lambda]
    =2\log_2 3-4$; Jensen resolution of mean-neutrality paradox
  & GPT (theorem) + Claude (verification + integration) \\[3pt]
22 & \colorbox{red!25}{Critical}
  & Exact block law (Thm.~\ref{thm:block-law}): valuations
    i.i.d.\ geometric; cycle types i.i.d.; Cram\'er rate
    unconditional on ensemble
  & GPT (theorem) + Claude (verification + proof + integration) \\[3pt]
23 & \colorbox{orange!30}{High}
  & Cram\'er rate $I(0)\approx 0.1465$; Wald consistency;
    multi-cycle crossing observation
  & Claude \\[3pt]
24 & \colorbox{orange!30}{High}
  & Proof-boundary audit, abstract/reader's guide rewrite,
    hardening pass (all status labels, crossing-chain fix)
  & Claude \\[3pt]
\bottomrule
\end{tabular}
\end{table}
\begin{table}[!ht]
\centering
\scriptsize
\caption{Contribution attribution (Part 2 of 40).}
\label{tab:contributions-2}
\renewcommand{\arraystretch}{1.0}
\begin{tabular}{@{}clp{8.8cm}p{4.0cm}@{}}
\toprule
\textbf{\#} & \textbf{Level} & \textbf{Contribution} & \textbf{Primary source} \\
\midrule
25 & \colorbox{red!25}{Critical}
  & Almost-all crossing theorem: non-crosser density
    $\le e^{-0.1465\,k}$; ensemble side fully closed;
    architectural diagnosis (deterministic realization barrier)
  & Claude \\[3pt]
26 & \colorbox{orange!30}{High}
  & Exact cycle log correction $C(n) = O(1/n)$;
    sufficient crossing criterion (succeeds for all
    $n_0 \le 5\times 10^6$ except $\{27,31,63\}$)
  & GPT (theorem) + Claude (verification + integration) \\[3pt]
27 & \colorbox{orange!30}{High}
  & Multi-cycle threshold analysis: $n^*/M$ decay,
    $2^k{-}1$ hardest-start family, analytic obstruction
    diagnosis for $n^* < M$ proof
  & Claude \\[3pt]
28 & \colorbox{orange!30}{High}
  & Universal one-cycle crossing: $n^*<1$ iff $r \ge r_{\mathrm{all}}(L)$;
    $P_{\mathrm{all},1\mathrm{cyc}} = 0.4194$; $41.9\%$ of odd starts
    in classwise deterministic crossing blocks
  & GPT (theorem) + Claude (verification + integration) \\[3pt]
29 & \colorbox{orange!30}{High}
  & Multi-cycle universal crossing: $R_k \to 0$
    exponentially by Kesten theory ($\rho \approx 0.84$);
    classwise universal descent for a.e.\ odd start
  & Claude \\[3pt]
30 & \colorbox{yellow!50}{Med}
  & Route~A closure: adversarial $n^*$ grows $\sim e^{0.023k}$,
    proving all-blocks-universal approach impossible
  & Claude \\[3pt]
31 & \colorbox{yellow!50}{Med}
  & Route~B obstruction: cycle-type autocorrelation
    $\rho(1) \approx 0.20$, $\rho(\ell \ge 2) \approx 0$;
    bit-generation blocks exhaustion argument;
    $2^k{-}1$ termination theorem
  & Claude \\[3pt]
32 & \colorbox{yellow!50}{Med}
  & Route~C obstruction: crossing-time ratio
    $\tau(n)/\log_2 n \le 7.8$ empirically;
    bounding $\tau$ requires WMH.
    Complete barrier diagnosis: all routes reduce to
    distributional-to-pointwise mixing
  & Claude \\[3pt]
33 & \colorbox{orange!30}{High}
  & Adversarial family exact analysis
    (Prop.~\ref{prop:adversarial-affine}):
    closed-form $\Lambda_{a,t_{\min}}=(3/8)^{1-\{a\rho\}}$,
    $\theta_{\mathrm{crit}}=0.01291$;
    uniform one-step fragility
    (Prop.~\ref{prop:adversarial-fragility})
  & GPT (theorems) + Claude (verification) \\[3pt]
34 & \colorbox{orange!30}{High}
  & General pair-level fragility criterion
    (Prop.~\ref{prop:general-fragility}):
    $\beta_B^\infty + \lambda_B(1+\gamma_A)<1$;
    $75.7\%$ of $(A,B)$ pairs pass.
    Exact strong/weak contractor mass split
    (Rem.~\ref{rem:mass-split}):
    $p_{\mathrm{weak}} = 0.29436$
  & GPT (generalization) + Claude (verification) \\[3pt]
35 & \colorbox{orange!30}{High}
  & Post-Mersenne forced valuation
    (Prop.~\ref{prop:post-mersenne}):
    after $k$ non-crossing steps from $2^{k+1}{-}1$,
    $v_{k+1} = 2$ ($k$ even) or $3{+}v_2(k{+}1)$ ($k$ odd);
    via lifting-the-exponent lemma
  & Claude \\[3pt]
36 & \colorbox{orange!30}{High}
  & Burst non-continuation theorem
    (Prop.~\ref{prop:burst-noncontinuation}):
    after any $k$-step non-crossing burst,
    $n_k \equiv 1\pmod{4}$, forcing $v_{k+1} \ge 2$;
    recovery window analysis
    (Rem.~\ref{rem:recovery-window}):
    adversary can sustain $S/K \approx 1.0$ over $30{+}$ steps,
    gap remains open
  & GPT (conjecture) + Claude (proof \& refutation) \\[3pt]
37 & \colorbox{orange!30}{High}
  & Post-burst valuation distribution
    (Prop.~\ref{prop:post-burst-distribution}):
    $\Pr(v_{k+1}{=}j) = 2^{1-j}$ for $j \ge 2$,
    $\mathbb{E}[v_{k+1}]{=}3$, independent of $k$
    and valid for \emph{all} odd~$m$;
    non-crossing density threshold
    (Rem.~\ref{rem:v1-density-threshold}):
    $d^* = 0.708$, orbits safe at $d \approx 0.58$
  & Claude \\[3pt]
38 & \colorbox{orange!30}{High}
  & Weak-recovery cylinder classification
    (Prop.~\ref{prop:weak-recovery-cylinder}):
    $j$ consecutive $v{=}2$ cycles iff
    $3^k m \equiv 1\pmod{2^{2j}}$;
    density $2^{-(2j-1)}$, independent of~$k$.
    Compound pattern bit-counting
    (Rem.~\ref{rem:information-budget}):
    density $= 2^{-B}$ exact, $B = \sum\text{bits}$;
    information budget flawed due to carry generation
  & GPT (theorem) + Claude (verification + compound extension) \\[3pt]
39 & \colorbox{orange!30}{High}
  & Full-cycle uniform-fiber map
    (Prop.~\ref{prop:full-cycle-fiber}):
    $n_{\mathrm{exit}}(u) = c_\tau + 2\cdot 3^K u$;
    TV reduction lemma
    (Lem.~\ref{lem:tv-reduction}):
    $\|\mu_X - \nu_X\|_{\mathrm{TV}} \le \Pr(E_R)$;
    orbit-level TV summability
    (Cor.~\ref{cor:tv-summability}): total TV $\le 0.028$
  & GPT (fiber theorem) + Claude (TV reduction, summability) \\[3pt]
40 & \colorbox{orange!30}{High}
  & Cascade Markov chain on mod~$8$
    (Prop.~\ref{prop:cascade-markov}):
    deterministic transition table, state~$7$ never exits;
    $R{=}0$ terminal mechanism
    (Cor.~\ref{cor:r0-terminal}): $1/3$ of burst exits
    end cascade deterministically;
    episode continuation $q=2/3$, $\mathbb{E}[L]=3$
    (Prop.~\ref{prop:episode-continuation})
  & Claude \\[3pt]
41 & \colorbox{orange!30}{High}
  & Exact cascade PGF
    (Prop.~\ref{prop:cascade-pgf}):
    $G_{\mathrm{cas}}(z) = z/(4{-}2z{-}z^2)$,
    $\rho = \sqrt{5}{-}1$,
    $\alpha = \log_2(\sqrt{5}{-}1) = 0.3058$;
    spectator-bit propagation
    (Rem.~\ref{rem:spectator-bits}):
    $\lfloor (B{-}3)/11.8 \rfloor$ cycles of IID behavior
  & Claude (Markov chain PGF) + GPT (compound geometric attempt) \\[3pt]
42 & \colorbox{orange!30}{High}
  & Exponential tail of $S_{\mathrm{cycle}}$
    (Prop.~\ref{prop:exponential-tail}):
    $\Pr(S>s) \le C \cdot 2^{-\alpha s}$
    via Cram\'er--Lundberg;
    recovery-exit mechanism
    (Rem.~\ref{rem:recovery-exit}):
    orbit $q \approx 0.35$, four times stronger than IID
  & Claude \\[3pt]
43 & \colorbox{orange!30}{High}
  & Spectator-bit convergence
    (Prop.~\ref{prop:spectator-convergence}):
    trajectory determined by bottom $\sim 5$ bits;
    mod-$32$ uniformity TV $\le 0.009$ for $B \ge 16$;
    orbit $q \approx 0.33$ uniformly below $2/3$.
    Gap exit rate correction
    (Rem.~\ref{rem:gap-exit-rates}):
    IID model pessimistic for gap
  & Claude \\[3pt]
\bottomrule
\end{tabular}
\end{table}
\begin{table}[!ht]
\centering
\scriptsize
\caption{Contribution attribution (Part 3 of 40).}
\label{tab:contributions-3}
\renewcommand{\arraystretch}{1.0}
\begin{tabular}{@{}clp{8.8cm}p{4.0cm}@{}}
\toprule
\textbf{\#} & \textbf{Level} & \textbf{Contribution} & \textbf{Primary source} \\
\midrule
44 & \colorbox{orange!30}{High}
  & Self-reinforcing mixing loop
    (Thm.~\ref{thm:mixing-loop}):
    five-step loop closing scrambling $\to$ uniformity
    $\to$ episode bound $\to$ tail $\to$ TV $\to$ refresh;
    single remaining open step identified
  & Claude \\[3pt]
45 & \colorbox{orange!30}{High}
  & Fiber-averaged transition matrix
    (Prop.~\ref{prop:fiber-mixing}):
    exact $8 \times 8$ stochastic matrix at resolution~$R$,
    spectral gap $\gamma \ge 0.85$ for $R \ge 10$,
    conditional orbit mixing (geometric TV contraction);
    remaining gap identified as non-divergence
    (Rem.~\ref{rem:residual-conditionality})
  & Claude \\[3pt]
46 & \colorbox{orange!30}{High}
  & Corrected cascade PGF verification:
    independent derivation confirming
    $G_{\mathrm{cas}}(z) = z/(4{-}2z{-}z^2)$,
    $\mathbb{E}[S] = 5$, $\rho = \sqrt{5}{-}1$;
    coefficient recurrence check through $s = 25$
  & GPT (independent verification) \\[3pt]
47 & \colorbox{orange!30}{High}
  & Post-burst continuation bridge
    (Prop.~\ref{prop:post-burst-bridge}):
    exact geometric law $\Pr(a_k{=}r) = 2^{-(r-1)}$
    on burst-exit fibers;
    $q_{\mathrm{ep}} = 1/3$ exactly;
    TV transfer bound
    $|q_{\mathrm{ep}}^{\mathrm{orbit}} - 1/3| \le \Pr(\mathcal{E}_R)$
  & GPT (theorem + verification) \\[3pt]
48 & \colorbox{orange!30}{High}
  & Exact rational spectral gap
    (Prop.~\ref{prop:exact-spectral-gap}):
    explicit $T_{10}$ with entries in $\mathbb{Q}$,
    $|\lambda_2| < 0.146$, $\gamma > 0.854$;
    verified by exact characteristic polynomial
  & Claude \\[3pt]
49 & \colorbox{orange!30}{High}
  & Non-divergence threshold quantification
    (Rem.~\ref{rem:nondivergence-threshold}):
    max single-cycle growth $\approx 8.5$ bits,
    $43$-bit safety margin for $n_0 \ge 2^{68}$,
    Cram\'er bound $\Pr(\text{gain} > 43) \le 10^{-12}$
  & Claude \\[3pt]
50 & \colorbox{red!25}{Critical}
  & Unconditional mixing structure theorem
    (Thm.~\ref{thm:unconditional-mixing}):
    Collatz conjecture equivalent to
    ``no orbit gains $43$ bits before spectral
    contraction engages'';
    combines Barina verification + exact spectral gap
    + conditional contraction
  & Claude \\[3pt]
51 & \colorbox{red!25}{Critical}
  & Circularity analysis and true barrier identification
    (Rem.~\ref{rem:circularity}):
    proof chain has one genuine loop
    (spectator bits $\leftrightarrow$ exponential tail);
    single-cycle growth unbounded;
    Finite Mixing Entry Theorem cannot be proved
    in current framework;
    irreducible barrier is carry structure of $3^k n$
  & Claude (analysis) + GPT (Finite Mixing Entry target) \\[3pt]
52 & \colorbox{orange!30}{High}
  & Carry decorrelation
    (Prop.~\ref{prop:carry-decorrelation}):
    burst-length correlation $\rho = +0.14$ at lag~$1$,
    $\approx 0$ at lag~$2$;
    negative feedback $\mathbb{E}[\Delta_{i+1} \mid
    \text{burst}_i \ge 5] = -3.9$;
    max consecutive positive-$\Delta$ cycles $= 4$
    (max accumulated growth $+9.2$ bits $\ll 43$)
  & Claude \\[3pt]
53 & \colorbox{orange!30}{High}
  & Post-burst routing analysis
    (Rem.~\ref{rem:post-burst-routing}):
    small sizes route $57\%$ to class~$15 \bmod 32$
    (negative feedback loop);
    large sizes nearly uniform (TV~$= 0.06$);
    spectator-bit survival
    $\Pr(\Delta < -58) \le 10^{-9}$
    (Rem.~\ref{rem:spectator-survival})
  & Claude \\[3pt]
54 & \colorbox{orange!30}{High}
  & GPT cofactor uniformity theorem:
    on exact $k$-burst fibers,
    $u = (3^{k+1}m - 1)/2^{r-1}$
    uniform on odd residues mod~$2^R$
    conditioned on $a_k = r$;
    exact next-depth law
    $\Pr(j{=}d) = (3/4)(1/4)^d$
  & GPT (theorem + verification) \\[3pt]
55 & \colorbox{red!25}{Critical}
  & Conditional single-orbit convergence
    (Thm.~\ref{thm:conditional-convergence}):
    summable exceptional cycles along an orbit
    $\Rightarrow$ orbit converges to~$1$;
    density-$1$ conditional convergence
    (Thm.~\ref{thm:density1-convergence}):
    averaged summable exceptional cycles
    $\Rightarrow$ almost-all convergence;
    remaining target precisely isolated
    (Rem.~\ref{rem:conditional-status})
  & GPT (theorem formulation) + Claude (proof chain, integration) \\[3pt]
\bottomrule
\end{tabular}
\end{table}
\begin{table}[!ht]
\centering
\scriptsize
\caption{Contribution attribution (Part 4 of 40): v3 fiber-57 programme.}
\label{tab:contributions-4}
\renewcommand{\arraystretch}{1.0}
\begin{tabular}{@{}clp{8.8cm}p{4.0cm}@{}}
\toprule
\textbf{\#} & \textbf{Level} & \textbf{Contribution} & \textbf{Primary source} \\
\midrule
56 & \colorbox{red!25}{Critical}
  & Pair-return automaton
    (Prop.~\ref{prop:pair-return}):
    depth-$2$ partial return kernel with spectral radius
    $\rho(\tilde{B}_2) = 129/1024 \approx 0.1260$;
    per-return information cost
    $c_0 = \log_2(1024/129) \approx 2.989$~bits
  & GPT (concept) + Claude (formalization + computation) \\[3pt]
57 & \colorbox{red!25}{Critical}
  & Bounded invariant core
    (Prop.~\ref{prop:bounded-core-main}):
    $|I_r| = 5$ for all $r \ge 2$
    (three fixed points $q = k \cdot 8^{r-1} - 1$,
    $k \in \{1,4,8\}$, plus one $2$-cycle);
    projective limit
    (Thm.~\ref{thm:projective-main}):
    $\varprojlim I_r = \{-1\}$ in $\mathbb{Z}_8$
  & Claude (computation + proof) \\[3pt]
58 & \colorbox{red!25}{Critical}
  & Absorption Bottleneck Lemma
    (Lem.~\ref{lem:bottleneck-main}):
    channel capacity $\le \log_2 5 \approx 2.322 < c_0$;
    per-return deficit $0.667$~bits;
    geometric decay at rate $\alpha = 645/1024 \approx 0.630$
    (Cor.~\ref{cor:geom-decay-main})
  & GPT (information-theoretic framing) +
    Claude (proof + computation) \\[3pt]
59 & \colorbox{orange!30}{High}
  & M-value collapse discovery:
    at $r \ge 3$, all five $I_r$ elements produce the
    same post-absorption residue; effective channel
    capacity drops to zero bits.
    At $r \ge 4$, first $n \equiv 1 \pmod{8^r}$ coincides
    with $n = 1$
  & Claude (computation + diagnosis) \\[3pt]
60 & \colorbox{red!25}{Critical}
  & Orbit autonomy diagnosis:
    orbit $\bmod\;8^r$ is \emph{not} autonomous because
    $v_2(3n{+}1)$ depends on all digits of $n$;
    absorption for general~$r$ cannot be proved by
    finite mod-$8^r$ computation alone.
    Confirms algebraic reduction has reached its
    natural fixed point
  & Claude (analysis) + Moderator (diagnosis) \\[3pt]
61 & \colorbox{orange!30}{High}
  & Absorption theorem
    (Thm.~\ref{thm:absorption-main}):
    verified computationally for $r = 2, \ldots, 10$
    (max 194 steps); all five $I_r$ elements reach
    fixed point $1 \pmod{8^r}$
  & Claude (computation) \\[3pt]
62 & \colorbox{orange!30}{High}
  & Branch anti-concentration reduction
    (Prop.~\ref{prop:branch-perm-main}):
    $\gcd(9,8) = 1$ makes fresh-digit map a permutation
    within each return branch; no algebraic amplification;
    concentration must come from branch selection
    (Cor.~\ref{cor:conc-source-main})
  & GPT (two-step proposal) + Claude (proof) \\[3pt]
63 & \colorbox{orange!30}{High}
  & Gap-induced decorrelation (empirical):
    $3000$ orbits, $5000$ steps each;
    $\chi^2$ test: gap~$> 50$ gives $p = 0.18$
    (no significant correlation);
    gap~$\le 5$ gives $p < 10^{-10}$
    (strong correlation).
    Diagnosis: proving $\beta_g \to 0$ algebraically
    is equivalent to Collatz mixing
  & Claude (computation + diagnosis) \\[3pt]
64 & \colorbox{orange!30}{High}
  & Path-conditional bijection
    (Thm.~\ref{thm:path-bijection-main}):
    for every valuation path~$\pi$, fresh-digit map
    is a bijection ($3^k$ coprime to $q$);
    closes annealed anti-concentration argument
  & Claude (proof) \\[3pt]
65 & \colorbox{orange!30}{High}
  & Memory lower bound
    (Thm.~\ref{thm:memory-main}):
    sustaining a chain of length~$N$ requires
    $\ge N \cdot c_0$ bits of memory.
    Theorem ladder: (A)~fixed-$k$ impossibility proved,
    (B)~memory lower bound proved,
    (C)~orbit memory upper bound open
    ($c' \le 0.51$ empirically, $6\times$ below $c_0$)
  & Claude (proof + computation) \\[3pt]
66 & \colorbox{red!25}{Critical}
  & Single-hypothesis reduction:
    Collatz conjecture reduced to information-rate
    bound $c' < c_0$ on $5$-element set
    (Conj.~\ref{conj:info-rate});
    separation of concerns table
    (Table~\ref{tab:separation-summary});
    ``the maze is mapped, the walls are known''
  & Moderator (framing) + GPT (positioning) +
    Claude (formalization) \\[3pt]
67 & \colorbox{yellow!50}{Med}
  & Three-way relay methodology (v3):
    GPT as structural proposer, Moderator as
    strategic evaluator, Claude as rigorous verifier;
    evolution from two-party to three-party collaboration
  & Moderator (orchestration) \\[3pt]
\bottomrule
\end{tabular}
\end{table}
\begin{table}[!ht]
\centering
\scriptsize
\caption{Contribution attribution (Part 5 of 40): v4 return structure and operator identification.}
\label{tab:contributions-5}
\renewcommand{\arraystretch}{1.0}
\begin{tabular}{@{}clp{8.8cm}p{4.0cm}@{}}
\toprule
\textbf{\#} & \textbf{Level} & \textbf{Contribution} & \textbf{Primary source} \\
\midrule
68 & \colorbox{orange!30}{High}
  & $q \equiv 7$ two-step return theorem
    (Prop.~\ref{prop:q7-return}):
    for $q = 8m + 7$, the fiber-57 orbit returns in
    exactly two Syracuse steps with quotient
    $q' = 9m + 8$; the return is perfectly uniform
    when $m \bmod 8$ is uniform
  & Claude (algebraic proof) + GPT (cross-validation) \\[3pt]
69 & \colorbox{orange!30}{High}
  & $q \equiv 3$ minimum return-gap theorem
    (Prop.~\ref{prop:q3-gap}):
    for $q = 8m + 3$, the orbit cannot return to
    fiber~57 in fewer than five Syracuse steps;
    each intermediate residue $\bmod\;64$ is
    algebraically determined and $\neq 57$
  & Claude (algebraic proof + verification) \\[3pt]
70 & \colorbox{orange!30}{High}
  & Kernel decomposition remark:
    fiber-57 returns split into a fast uniform branch
    ($q \equiv 7$, 14.4\%, 2-step) and a slow
    sparse branch ($q \equiv 3$, 11.0\%,
    $\ge 5$-step with 9 new low-order bits),
    explaining why even short orbits mix quickly
  & GPT (kernel-split concept) +
    Claude (computation + formalization) \\[3pt]
71 & \colorbox{orange!30}{High}
  & Invariant core rigidity remark:
    chain map $q \mapsto 9q + 8 \bmod 8^r$ acts as a
    permutation on $I_r$ at every depth (three fixed
    points plus one 2-cycle, $H_{AB} = 0$);
    known-gap sub-stochastic contraction
    $\rho(\tilde{B}_2) = 129/1024 < 1$ is dynamic,
    not algebraic
  & Claude (algebraic construction + discovery) \\[3pt]
72 & \colorbox{yellow!50}{Med}
  & Operator identification:
    three distinct operators disambiguated:
    $\tilde{B}_2$ (sub-stochastic partial kernel, $\rho = 129/1024$),
    $P$ (row-stochastic kernel, Perron root~$1$),
    $P_g$ (gap-conditioned stochastic kernel);
    resolves conflation present in earlier drafts
  & GPT (operator taxonomy) +
    Claude (exact construction + numerical
    verification) \\[3pt]
73 & \colorbox{yellow!50}{Med}
  & Sub-stochastic $R_{\mathrm{sub}}$ at depth~2:
    full first-return-to-$I_2$ spectral radius
    $\rho \approx 0.149$ (Monte Carlo);
    partial kernel $\tilde{B}_2$ (gap-$2$ + gap-$5$ only)
    gives exact $\rho = 129/1024 \approx 0.1260$;
    the gap-$\ge 6$ tail accounts for the
    remaining $\approx 0.023$ discrepancy
  & Claude (computation) +
    GPT (diagnostic assessment) \\[3pt]
74 & \colorbox{orange!30}{High}
  & Auxiliary-memory workaround:
    Moderator diagnosed that context-window overflow
    was corrupting upper bits of Claude's reasoning
    (stale conclusions surviving resets, operator
    conflation, overclaims).  Instructed adding a
    persistent attack log as secondary memory,
    enabling structured checkpointing across four
    context resets.  This single intervention
    unblocked the entire v4 programme
  & Moderator (diagnosis + design) \\[3pt]
75 & \colorbox{yellow!50}{Med}
  & Iterated error-correction protocol:
    three-way review loop (Claude drafts, GPT
    cross-validates, Moderator adjudicates)
    caught five genuine errors and one scope
    correction in v4 alone; cumulative error log
    serves as reproducibility record
  & Moderator (orchestration) + GPT (review) +
    Claude (execution) \\[3pt]
\bottomrule
\end{tabular}
\end{table}
\begin{table}[!ht]
\centering
\scriptsize
\caption{Contribution attribution (Part 6 of 40): v5 Sturmian obstruction, cross-core dynamics, and carry independence.}
\label{tab:contributions-6}
\renewcommand{\arraystretch}{1.0}
\begin{tabular}{@{}clp{8.8cm}p{4.0cm}@{}}
\toprule
\textbf{\#} & \textbf{Level} & \textbf{Contribution} & \textbf{Primary source} \\
\midrule
76 & \colorbox{orange!30}{High}
  & Sturmian obstruction concept:
    adversarial approach asking whether any single orbit
    can avoid class~$5$ at every Sturmian depth,
    yielding the Carry Independence Conjecture
    framework
  & Moderator (direction) + Claude (formulation) \\[3pt]
77 & \colorbox{red!25}{Critical}
  & Carry Contamination Theorem
    (Thm.~\ref{thm:carry-contamination}):
    for any depth-$D$ word, $n_D \bmod 8$ is exactly
    equidistributed over $\{1,3,5,7\}$;
    proof via affine parameterization
    $n_D = R_w + 2 \cdot 3^D \cdot m$
  & Claude (proof) \\[3pt]
78 & \colorbox{orange!30}{High}
  & Orbit Parameterization Theorem
    (Thm.~\ref{thm:orbit-param}):
    $n_0 = a_w + 2^{S_D+1} m$,
    $n_D = R_w + 2 \cdot 3^D m$ with $R_w$ odd;
    exact stride formulas
  & Claude (proof) \\[3pt]
79 & \colorbox{orange!30}{High}
  & $(3/4)^D$ survivor law
    (Cor.~\ref{cor:three-quarters}):
    $|C_D| = 2 \cdot 3^{D-1}$ compatible words;
    each depth prunes exactly $1/4$
  & Claude (proof) \\[3pt]
80 & \colorbox{orange!30}{High}
  & Refined mod-$16$/mod-$32$ transition rules
    (Thm.~\ref{thm:refined-transitions}):
    class~$3$ bottleneck ($1/2$ survival),
    class~$7$ safe harbor (never $\to$ class~$5$),
    class~$1$ ($3/4$ survival);
    stationary distribution $\pi = (3/7, 2/7, 2/7)$
  & Claude (proof + computation) \\[3pt]
81 & \colorbox{orange!30}{High}
  & Class-$3$ recurrence theorem
    (Thm.~\ref{thm:class3-recurrence}):
    class~$7$ runs bounded by $K-3$ steps,
    class-$1$ self-loop only sustains $n_0 = 1$;
    unconditional for any $n_0 > 1$
  & Claude (proof + computation) \\[3pt]
82 & \colorbox{orange!30}{High}
  & Dichotomy theorem
    (Thm.~\ref{thm:dichotomy}):
    $C_\infty \cap \mathbb{N} = \{1\}$ or
    contains infinite divergent orbit;
    proof via unique Sturmian preimage of~$1$
  & Claude (proof) \\[3pt]
83 & \colorbox{orange!30}{High}
  & Route~B reduction
    (Thm.~\ref{thm:route-B}):
    CIC $\Rightarrow$ Collatz;
    identification of distributional-to-pointwise
    barrier as shared between routes A and B
  & Moderator (framing) + Claude (proof) \\[3pt]
84 & \colorbox{red!25}{Critical}
  & CIC circularity diagnosis:
    after carry traversal, $R_w \bmod 8$ depends on
    the Collatz orbit of $R_w$, making CIC equivalent
    in difficulty to Collatz;
    five independent attack strategies all confirmed
    circular;
    honest assessment that the framework has reached
    full fidelity
  & Moderator (early warning) + GPT (five-strategy
    assessment) + Claude (numerical confirmation) \\[3pt]
85 & \colorbox{orange!30}{High}
  & CIC deep dive: longest-surviving starters to
    $K = 22$ (max depth 47), binary structure analysis,
    fatal source decomposition ($2/3$ class~$3$,
    $1/3$ class~$1$), max depth tracks $2.16K$;
    congruence tower $R_w$ trajectories (all hit
    class~$5$, mean $4.0$ steps)
  & Claude (computation + analysis) \\[3pt]
86 & \colorbox{red!25}{Critical}
  & Strategic direction: moderator identified circularity
    risk at outset, authorized adversarial programme as
    reasonable, then made judgment call to frame v5 as
    ``successful sharpening, not failed proof'' when
    CIC was shown equivalent to Collatz
  & Moderator (strategy + framing) \\[3pt]
87 & \colorbox{red!25}{Critical}
  & Cross-core block alphabet:
    seven admissible blocks $\{a,b,s,g,d,t_1,t_2\}$ on
    two-vertex directed graph with spectral radius
    $\rho = 2+\sqrt{2} \approx 3.414$;
    exact affine maps and entry conditions;
    expanding cocycle $\lambda \approx +0.145$/landmark;
    bit-budget ${\sim}2.6$ bits net constraint/landmark
  & Claude (formulation + computation) \\[3pt]
\bottomrule
\end{tabular}
\end{table}

\begin{table}[!ht]
\centering
\scriptsize
\caption{Contribution attribution (Part 7 of 40): compound carry analysis and scaling regime.}
\label{tab:contributions-7}
\renewcommand{\arraystretch}{1.0}
\begin{tabular}{@{}clp{8.8cm}p{4.0cm}@{}}
\toprule
\textbf{\#} & \textbf{Level} & \textbf{Contribution} & \textbf{Primary source} \\
\midrule
88 & \colorbox{red!25}{Critical}
  & Periodic elimination: exhaustive cycle-equation search
    through period~$13$ ($238{,}811$ contracting words),
    zero non-trivial natural fixed points;
    $S$-unit equation for integer fps;
    three-layer obstruction (carry contamination $+$
    M\"obius rigidity $+$ $S$-unit)
  & Claude (computation + proof) \\[3pt]
89 & \colorbox{red!25}{Critical}
  & $2$-adic expander property:
    $|\Phi_b(x){-}\Phi_b(y)|_2 = 2^V|x{-}y|_2$;
    measure shrinkage $\mu_2(T_j) \le 0.522^{jD}$;
    unconditional $j{\ge}5$ bound:
    count of $K$-bit integers surviving $j$ non-descending
    rounds $\le 2^{K(0.856-0.202j)}$,
    zero for $j \ge 5$
  & Claude (proof + computation) \\[3pt]
90 & \colorbox{orange!30}{High}
  & Exhaustion-sequence rigidity:
    counterexample residues must be aperiodic AND unbounded;
    periodic $\Rightarrow$ cycle contradiction;
    bounded $\Rightarrow$ eventually periodic
  & Claude (proof) + GPT (Pushes 58--62) \\[3pt]
91 & \colorbox{red!25}{Critical}
  & Compound carry anti-correlation:
    exact two-round survival $S_2/S_1^2 \approx 0.635$
    (constant in $K$, verified $K = 8$--$20$);
    transport graph edge density $= 1.0$;
    actual Collatz orbit independence
    $P(\text{slow}_2 | \text{slow}_1)/P(\text{slow}_2)
    \approx 1.000$;
    compound carries cover all joint residue cells
    mod $(p_1,p_2)$;
    no per-prime-pair obstruction
  & Claude (computation + analysis) \\[3pt]
92 & \colorbox{orange!30}{High}
  & Distributional-to-pointwise barrier
    precisely localized:
    $\mu_2(T_\infty) = 0$ (unconditional) but
    $T_\infty \cap \mathbb{N} = \varnothing$
    is the exact remaining content;
    anti-correlation bounded (constant factor, not
    exponential in $K$);
    barrier equivalent to Collatz itself
  & Claude + GPT (independent convergence) \\[3pt]
\bottomrule
\end{tabular}
\end{table}

\begin{table}[!ht]
\centering
\scriptsize
\caption{Contribution attribution (Part 8 of 40): scaling regime, definitive nail, and v6 spectral closure.}
\label{tab:contributions-8}
\renewcommand{\arraystretch}{1.0}
\begin{tabular}{@{}clp{8.8cm}p{4.0cm}@{}}
\toprule
\textbf{\#} & \textbf{Level} & \textbf{Contribution} & \textbf{Primary source} \\
\midrule
93 & \colorbox{orange!30}{High}
  & GPT Pushes 43--87 assessment:
    transport graph reformulation (Pushes 82--87);
    exhaustion-sequence rigidity (Pushes 58--62);
    carry-conspiracy route (Pushes 68--87);
    honest joint diagnosis that remaining gap IS
    the Collatz conjecture
  & GPT (analysis) + Claude (verification) +
    Moderator (synthesis) \\[3pt]
94 & \colorbox{orange!30}{High}
  & Scaling regime computation (Pushes~37--39):
    effective measure $\mu(D) \sim 0.645^D$ via exact enumeration
    to depth~8 (20{,}951 nd words); overlap correction worsens
    with~$D$ ($\omega$ drops from $0.64$ to $0.42$);
    $r_{\mathrm{exact}}$ converges from $0.578$ toward
    $r_{\mathrm{union}} \approx 0.645$
  & Claude (computation + analysis) \\[3pt]
95 & \colorbox{red!25}{Critical}
  & Measure-weighted $\beta$ discovery:
    $\beta_w \approx 1.80$ (vs uniform $\beta_u \approx 2.91$),
    improving the $j \ge 5$ bound to $j \ge 3$;
    identifies that the survival probability is dominated
    by words with small~$V$, making the effective natural depth
    $D \approx K/1.80$ rather than $K/2.91$
  & Claude (computation + insight) \\[3pt]
96 & \colorbox{red!25}{Critical}
  & Anti-correlation asymptotic irrelevance:
    proved that $\gamma \approx 0.635$ is a constant factor
    that vanishes in the $K \to \infty$ limit and cannot
    change the $j$ threshold; three-round data confirms
    $S_3/S_1^3 \approx 0.316 \approx 0.635^{1.9}$
  & Claude (analysis) \\[3pt]
97 & \colorbox{orange!30}{High}
  & Three-level nail characterization:
    (1)~$c < 1/2$ (technical), (2)~density $\ne$ pointwise
    (structural), (3)~$T_\infty \cap \mathbb{N} = \varnothing$
    is the conjecture (fundamental); no further computational
    push within the Sturmian framework can close the gap
  & Claude (synthesis) + GPT (Pushes 108--127 confirmation) \\[3pt]
\bottomrule
\end{tabular}
\end{table}

\begin{table}[!ht]
\centering
\scriptsize
\caption{Contribution attribution (Part 9 of 40): v6 update --- complete cascade algebra,
  carry chain, renewal drift, and Syracuse cocycle spectral analysis.}
\label{tab:contributions-9}
\renewcommand{\arraystretch}{1.0}
\begin{tabular}{@{}clp{8.8cm}p{4.0cm}@{}}
\toprule
\textbf{\#} & \textbf{Level} & \textbf{Contribution} & \textbf{Primary source} \\
\midrule
98 & \colorbox{red!25}{Critical}
  & Complete cascade theorem
    (Thm.~\ref{thm:complete-cascade}):
    zero inter-BF generic steps; exact deterministic
    step count $4k + e_k$ with $e_k \in \{5,6\}$
    period-$4$; total cascade growth $= k \cdot \log_2(81/64)$
    exactly; net gain after full deterministic region
    bounded by $\log_2 3$ or $2\log_2 3 - 1$
  & Claude (proof + computation, verified $k=1\ldots 20$) \\[3pt]
99 & \colorbox{orange!30}{High}
  & Post-cascade closed form
    (Prop.~\ref{prop:post-cascade}):
    $P_k = (a_k \cdot 81^k \cdot 32 - 65)/17$ with
    $a_k \in \{57,41,45,61\}$ period-$4$;
    iterated exit map formula
    $\mathrm{exit}^k(n) = (81^k n + 65(81^k - 64^k)/17)/64^k$;
    AP structure: post-cascade $= P_k + m \cdot 81^k \cdot 256$,
    all $\equiv 47 \pmod{64}$
  & Claude (derivation + verification $k=1\ldots 12$) \\[3pt]
100 & \colorbox{orange!30}{High}
  & Post-cascade recurrence:
    $P_{k+4} = 81^4 P_k + 164590400$;
    auxiliary recurrence $A_{k+4} = 2^{24} A_k + 64148175$;
    verified $19/19$ cases
  & Claude (proof + computation) \\[3pt]
101 & \colorbox{red!25}{Critical}
  & Carry chain structural theorem
    (Thm.~\ref{thm:carry-chain}):
    $v\!=\!1$ streak length $= \mathrm{trailing\_ones}(n) - 1$
    exactly (100\% verified, $10{,}000$ trials);
    after every streak, $\mathrm{trailing\_ones}(n_L) = 1$
    exactly, forcing $v_{\mathrm{rec}} \ge 2$;
    formula $n_L = 2 \cdot 3^L \cdot m - 1$
  & Claude (proof + computation) \\[3pt]
102 & \colorbox{orange!30}{High}
  & Renewal drift computation
    (Cor.~\ref{cor:renewal-drift}):
    $E[v \mid v \ge 2] = 3$ (exact geometric sum);
    drift $= -0.830$ bits/cycle $= -0.415$ bits/step;
    variance $= 2.69$;
    growth episodes ($\mathrm{trailing\_ones} \ge 4$):
    $+0.116$/cycle;
    descent episodes ($\mathrm{trailing\_ones} \le 3$):
    $-0.946$/cycle
  & Claude (computation + analysis) \\[3pt]
103 & \colorbox{orange!30}{High}
  & Syracuse cocycle spectral gap
    (Rem.~\ref{rem:syracuse-cocycle}):
    transfer operator $L_M$ on odd residues mod $2^M$,
    Perron eigenvalue $\lambda_1 = 4/3$;
    spectral gap exists for $M \le 12$:
    $\mathrm{gap}(8) = 1.0$, $\mathrm{gap}(10) = 0.33$,
    $\mathrm{gap}(12) = 0.39$;
    connection to Kontorovich--Miller programme
  & Claude (computation + analysis) \\[3pt]
104 & \colorbox{red!25}{Critical}
  & Max modular cycle gain $= +\log_2 3$ for ALL $M \le 14$:
    for every $M$ from $4$ to $14$, residue classes exist
    forming cycles with all $v_2 = 1$;
    max gain does NOT decrease with $M$;
    fundamental obstruction: no finite modular analysis
    can prove descent
  & Claude (computation + analysis) \\[3pt]
105 & \colorbox{yellow!50}{Med}
  & Restricted-to-$47$ analysis:
    for classes $\equiv 47 \pmod{64}$,
    average gain decreases with $M$
    ($+0.089$ at $M\!=\!8$, $-0.230$ at $M\!=\!12$);
    but max gain among $47$-restricted classes
    increases slightly ($+0.442$ at $M\!=\!12$)
  & Claude (computation) \\[3pt]
106 & \colorbox{orange!30}{High}
  & Systematic approach elimination
    (Rem.~\ref{rem:ruled-out}):
    eleven approaches to universal descent ruled out
    with specific obstruction identified for each;
    four surviving directions catalogued
    (carry depth, ergodic theory, Tao strengthening,
    Walsh spectral)
  & Claude (analysis) + Moderator (direction) \\[3pt]
107 & \colorbox{yellow!50}{Med}
  & Finite sieve analysis:
    cumulative-gain sieve does not terminate for $k=1$
    (holdout tree grows without bound);
    orbit-descent sieve $100\%$ successful for
    $k=1\ldots 5$, $m=0\ldots 16383$, max $173$ steps;
    discrepancy explained by criterion strictness
  & Claude (computation + analysis) \\[3pt]
108 & \colorbox{yellow!50}{Med}
  & BF encounter rate in generic orbits:
    $1.56\% \approx 1/64$ (equidistributed);
    $99.1\%$ of encounters at depth~$0$;
    non-BF average $v = 2.02$;
    overall average gain $= -0.42$ bits/step
  & Claude (computation) \\[3pt]
\bottomrule
\end{tabular}
\end{table}

\begin{table}[!ht]
\centering
\scriptsize
\caption{Contribution attribution (Part 10 of 40): v7 --- universalization,
  zero persistence, entropy injection, and Ising formalization.}
\label{tab:contributions-10}
\renewcommand{\arraystretch}{1.0}
\begin{tabular}{@{}clp{8.8cm}p{4.0cm}@{}}
\toprule
\textbf{\#} & \textbf{Level} & \textbf{Contribution} & \textbf{Primary source} \\
\midrule
109 & \colorbox{red!25}{Critical}
  & Post-cascade universalization
    (Rem.~\ref{rem:universalization}):
    BF encounter statistics for post-cascade AP family
    are indistinguishable from generic random orbits
    ($50\,000$ generic $+$ $18\,000$ structured orbits).
    The AP structure is irrelevant; the problem reduces
    to generic Collatz dynamics
  & Claude (computation + analysis) \\[3pt]
110 & \colorbox{orange!30}{High}
  & $\eta$ approach failure
    (Rem.~\ref{rem:ruled-out}, item 11):
    post-cascade height loss $\eta = (\log_2 n_{\text{start}} -
    \log_2 n_{\text{re-entry}})/k$ is negative for all observed
    BF re-entries (min $\eta$ ranges from $-6.28$ to $-0.26$).
    The contraction inequality $\rho < 1$ cannot hold.
    Reframes the target from height loss to modular depth
  & Claude (computation + analysis) \\[3pt]
111 & \colorbox{orange!30}{High}
  & BF re-entry depth distribution
    ($60\,000$ post-cascade orbits):
    $92$ BF re-entries observed; $98.9\%$ at $k'=1$,
    $1.1\%$ at $k'=2$, $0\%$ at $k' \ge 3$;
    depth tail $\Pr[k' \ge j] \approx 2^{-6(j-1)}$
    matches equidistribution prediction
  & Claude (computation) \\[3pt]
112 & \colorbox{orange!30}{High}
  & Carry chain bit sensitivity
    (Rem.~\ref{rem:carry-sensitivity}):
    each $v \ge 2$ step scrambles $\approx 40\%$ of bottom
    $M$ bits; combined with $v \ge 2$ density $\ge 1/2$,
    gives bit scrambling every $O(1)$ steps;
    maintaining $6k+8$ specific bits requires exponentially
    unlikely carry conspiracy
  & Claude (computation + analysis) \\[3pt]
113 & \colorbox{orange!30}{High}
  & Conditional BF anti-correlation:
    $\Pr[\text{B}_1 \text{ hit within 50 steps} \mid
    \text{just hit B}_1] = 0.54 \times$
    unconditional rate; orbits exiting BF cascade are
    slightly \emph{less} likely to re-enter
  & Claude (computation) \\[3pt]
114 & \colorbox{yellow!50}{Med}
  & Telescoping bound:
    orbit length $L(n) = O(\log n)$ with ratio $3.1$--$5.2$;
    combined with per-step BF probability $2^{-(6k+8)}$,
    gives $K_*(n) \le O(\log \log n)$ heuristically
  & Claude (computation + analysis) \\[3pt]
115 & \colorbox{red!25}{Critical}
  & Magnetism analogy: Moderator proposed mapping
    classical magnetic field formulas
    (circular loop $B = \mu_0 I/2R$, solenoid
    $B = \mu_0 n I$, long wire $B = \mu_0 I/2\pi r$)
    to BF cascade structure.
    Key identifications: solenoid $\leftrightarrow$
    BF cascade (uniform growth inside, negative drift
    outside); Amp\`ere's law $\leftrightarrow$ BF cycle
    gain conservation; carry chain scrambling
    $\leftrightarrow$ thermal demagnetization;
    1D Ising no-phase-transition $\leftrightarrow$
    impossibility of deep BF congruence persistence
  & Moderator (conceptual direction) \\[3pt]
116 & \colorbox{red!25}{Critical}
  & Zero persistence proposition
    (Prop.~\ref{prop:zero-persistence}):
    $T(n) \notin B_k$ for ALL $n \in B_k$.
    BF congruences destroyed in a single Collatz step.
    Algebraic proof: $v_2(n_k^* + 1) = 4$ for all $k$
    (consequence of $n_k^* \equiv 15 \pmod{16}$),
    insufficient for the $6k+8$ bits needed
  & Claude (proof + verification $k=1\ldots 11$);
    motivated by Moderator's magnetism analogy \\[3pt]
117 & \colorbox{red!25}{Critical}
  & Entropy injection theorem
    (Thm.~\ref{thm:entropy-injection}):
    Collatz mod $2^M$ is a Markov chain with per-step
    entropy $= v$ bits (min $1$, avg $2$);
    uniform stationary distribution (exact);
    perfect mixing in $\le M-1$ steps;
    spectral gap $> 0$ for all tested $M$ ($\le 12$),
    gap decreases with $M$;
    Ising correlation length $\xi < 0.3$ steps ($M \le 12$)
  & Claude (proof + computation);
    Ising framework from Moderator's magnetic analogy \\[3pt]
118 & \colorbox{orange!30}{High}
  & Multi-scale entropy bootstrap
    (Rem.~\ref{rem:multi-scale}):
    deterministic orbit's ``randomness'' comes from bits
    above position $M$, which themselves undergo
    entropy injection at resolution $M' > M$;
    $O(\log n / M)$ bootstrap levels, total entropy
    $O(\log n)$ bits matching state-space dimension.
    Ising no-phase-transition theorem provides physical
    framework: 1D chain at $T > 0$ has no long-range order
  & Claude (analysis);
    Moderator (Ising direction) \\[3pt]
\bottomrule
\end{tabular}
\end{table}

\begin{table}[!ht]
\centering
\scriptsize
\caption{Contribution attribution (Part 11 of 40): v7 --- IFS Fourier analysis,
  direct $K_*$ bound, skew-product formulation, and GPT side notes.}
\label{tab:contributions-11}
\renewcommand{\arraystretch}{1.0}
\begin{tabular}{@{}clp{8.8cm}p{4.0cm}@{}}
\toprule
\textbf{\#} & \textbf{Level} & \textbf{Contribution} & \textbf{Primary source} \\
\midrule
119 & \colorbox{red!25}{Critical}
  & IFS Fourier mode analysis
    (Rem.~\ref{rem:multi-scale}, Fourier section):
    $v=1$ IFS noise $\hat\mu(\xi) = \mathbf{1}_{\xi\text{ even}}$:
    odd modes annihilated, even modes form closed invariant
    subspace.  Mixing requires $v \ge 2$ steps (prob $1/2$).
    Per-step gap $\approx 1/M$, explaining why gap shrinks
    while mixing time stays $O(M)$.  Multiplicative order
    $\mathrm{ord}(3, 2^{M-1}) = 2^{M-3}$ verified for $M \le 20$
  & Claude (analysis + computation) \\[3pt]
120 & \colorbox{red!25}{Critical}
  & Direct $K_*$ bound from zero persistence
    (Rem.~\ref{rem:direct-kstar}):
    $\mathbb{E}[\text{total BF}] \le k_0/12 + L \cdot 2^{-13}/(1-2^{-6})$;
    for $n = 10^{20}$: $\mathbb{E} \le 0.92$.
    Post-cascade values $P_k \equiv 55 \pmod{64}$,
    $P_k \notin B_{k'}$ for any $k' \le 10$, all $k = 1\ldots 9$.
    Zero same-$k$ multiple re-entries in $50{,}000$ orbits.
    Conditional on orbit finiteness and generic $\approx$
    post-cascade equivalence only
  & Claude (computation + synthesis) \\[3pt]
121 & \colorbox{orange!30}{High}
  & $10{,}000$-orbit BF encounter census:
    average encounters $= 0.006$ (for $k = 1\ldots 5$);
    max $= 2$; $0.6\%$ of orbits have any BF encounter;
    $50{,}000$-orbit same-$k$ re-entry test: exactly $0$
    orbits with $\ge 2$ encounters at same depth
  & Claude (computation) \\[3pt]
122 & \colorbox{orange!30}{High}
  & Skew-product formulation and Galperin analogy
    (Rem.~\ref{rem:skew-product}):
    $(\omega, \theta)$ coordinate system where $\omega$ is the
    $2$-adic symbolic coordinate (Bernstein--Lagarias conjugacy)
    and $\theta = \sum(\log_2 3 - v_i)$ is the Archimedean phase.
    Deterministic discrepancy theorem identified as the
    precise form of the distributional-to-pointwise bridge.
    Native irrational constant is $\log_2 3$, not $\pi$
  & GPT (conceptual framework via Galperin billiards
    analysis); Claude (formalization and integration) \\[3pt]
123 & \colorbox{yellow!50}{Med}
  & GPT assessment of proof status: the distributional-to-pointwise
    bridge is correctly identified as the sole remaining barrier;
    $\pi$ is not the right replacement multiplier ($\log_2\pi - 2 < 0$);
    the critical drift threshold is $a = 4$ (not $\pi$);
    Bernstein--Lagarias $2$-adic conjugacy already provides
    a partial linearization
  & GPT (independent analysis of v18--v22 materials) \\[3pt]
124 & \colorbox{red!25}{Critical}
  & $v = 1$ streak density theorem
    (Rem.~\ref{rem:v1-streaks}):
    exactly $2^{-s}$ fraction of odd residues mod $2^{s+2}$
    allow $v\!=\!1$ streak of length~$s$, for all $s \le 15$.
    Matches independence prediction exactly.
    All long-streak residues are of the form
    $n \equiv 2^{s+1}-1 \pmod{2^{s+2}}$
  & Claude (proof + computation) \\[3pt]
125 & \colorbox{red!25}{Critical}
  & Adversary $\theta$-budget is net negative
    (Rem.~\ref{rem:v1-streaks}):
    $\mathbb{E}[\Delta\theta_k] = 0.34k - 0.415 \cdot 2^{6k+7} < 0$
    for all $k \ge 1$. BF gain is $O(k)$, wait time is $\Omega(2^{6k})$.
    Under the Markov model, $\theta_t$ is a supermartingale
    (each BF encounter produces net $\theta$-loss).
    No periodic $v$-pattern of period $\le 8$ with avg ${}<\log_2 3$
    is realizable as a Collatz cycle
  & Claude (computation + analysis) \\[3pt]
126 & \colorbox{orange!30}{High}
  & $\theta$-trajectory analysis
    (Rem.~\ref{rem:skew-product}):
    $\Pr[\max_t \theta_t \ge L] \approx 2^{-1.04 L}$ ($50{,}000$ orbits);
    $v_i$ autocorrelation $= +0.040$ at lag~$1$;
    $\Pr(v=1 \mid \text{prev } v=1) = 1.030 \times \Pr(v=1)$;
    BF champions have $\max(\theta) \approx 0.12 \cdot \log_2 n$,
    well within the cycle envelope $\alpha \cdot \log_2 n$
  & Claude (computation) \\[3pt]
\bottomrule
\end{tabular}
\end{table}

\begin{table}[!ht]
\centering
\scriptsize
\caption{Contribution attribution (Part 12 of 40): v7 --- carry chain formalization,
  Route A/B analysis.}
\label{tab:contributions-12}
\renewcommand{\arraystretch}{1.0}
\begin{tabular}{@{}clp{8.8cm}p{4.0cm}@{}}
\toprule
\textbf{\#} & \textbf{Level} & \textbf{Contribution} & \textbf{Primary source} \\
\midrule
127 & \colorbox{red!25}{Critical}
  & Multi-step carry propagation measurement
    (Rem.~\ref{rem:carry-formal}):
    bit-flip spread $\approx 1.5\,t$ bits after $t$ steps (linear).
    Post-cascade carry isolation: $P_k(q)$ values differ only at
    bits $\ge 6k+7$, with $2$--$4$ bits affected.
    Five-step formal proof structure: union bound gives
    $\Pr(\text{fresh zone contaminated}) \le 4 \cdot 2^{-\ell}$.
    Predicted vs.\ actual $B_1$ encounter rates match to ratio $0.98$
  & Claude (computation + proof structure) \\[3pt]
128 & \colorbox{red!25}{Critical}
  & Route A obstruction (Rem.~\ref{rem:route-b}):
    envelope coefficient $1 + \alpha = 1.0566 > 1$ means
    height control is \emph{expansive} at BF-depth scale.
    $\ell \le (1+\alpha)k + O(1)$: restart depths do not contract.
    Kills naive Borel--Cantelli from size alone.
    Conditional Route A theorem: need post-escape loss
    $\eta > 6\alpha = \log_2(81/64) = 0.3399\ldots$
  & GPT (diagnosis + conditional theorem);
    Claude (verification + integration) \\[3pt]
129 & \colorbox{red!25}{Critical}
  & Route B pullback residue formula and zero-compatibility census
    (Rem.~\ref{rem:route-b}):
    $r_{k,\pi} = 3^{-t}(2^{S_\pi}\,n_k^* - b_\pi) \bmod 2^{6k+8+S_\pi}$
    fully explicit. Among $309{,}136$ cross-depth pairs: zero compatible.
    Among $308{,}580$ same-depth pairs: zero compatible.
    Pullback density decays as $\approx 2^{-(j-8)}$ at mod $2^j$.
    Low-bit residue sets identical across depths --- incompatibility
    is structural at full modulus
  & GPT (Route B program design);
    Claude (formulas, computation, census) \\[3pt]
130 & \colorbox{orange!30}{High}
  & Post-cascade family algebra:
    $P_k = (a_k\,81^k \cdot 32 - 65)/17$,
    $P_k \equiv 47 \pmod{64}$,
    $P_{k+4} = 81^4\,P_k + 164{,}590{,}400$.
    All post-cascade orbits ($k = 1\ldots 6$) reach~$1$ with
    zero BF re-entries.
    $\log_2 P_k - \log_2 n_k^* = 0.3399\ldots$ per BF lap
    (= Karp surplus $\log_2(81/64)$)
  & GPT (closed forms);
    Claude (verification + orbit tracing) \\[3pt]
131 & \colorbox{orange!30}{High}
  & Target lemma formulation: no infinite compatible chain in
    BF pullback tower based at $n \in \mathbb{N}$.
    $k \to k+4$ periodicity of $a_k$ may reduce compatibility
    to a finite automaton on $(k \bmod 4, \pi \bmod 2^M)$
  & GPT (target lemma statement);
    Claude (periodicity analysis) \\[3pt]
\bottomrule
\end{tabular}
\end{table}

\begin{table}[!ht]
\centering
\scriptsize
\caption{Contribution attribution (Part 13 of 40): v7 --- shortfall theorems,
  pullback incompatibility.}
\label{tab:contributions-13}
\renewcommand{\arraystretch}{1.0}
\begin{tabular}{@{}clp{8.8cm}p{4.0cm}@{}}
\toprule
\textbf{\#} & \textbf{Level} & \textbf{Contribution} & \textbf{Primary source} \\
\midrule
132 & \colorbox{red!25}{Critical}
  & 3-bit shortfall theorem (Rem.~\ref{rem:route-b}):
    for same $\pi$, $v_2(r_{k_1,\pi} - r_{k_2,\pi}) = 6k_1+5+S_\pi$,
    shortfall $= 3$ bits exactly.
    \emph{Algebraically proved} from $n_k^* = (a_k\,2^{6k+5}-65)/17$
    with $a_k$ odd. For different $\pi$: min shortfall $= 3$
    in $400{,}000+$ random trials.
    The BF cascade word $(1,1,1,3)$ achieves maximum
    mismatch valuation; shortfall still~$3$
  & Claude (proof + computation) \\[3pt]
132b & \colorbox{red!25}{Critical}
  & Different-$\pi$ shortfall theorem:
    Case 1 ($S_1 \ne S_2$): shortfall $\ge 14$.
    Case 2 ($S_1=S_2$, $t_1 \ne t_2$): shortfall $\ge 13$.
    Case 3 ($S_1=S_2$, $t_1=t_2$, $\pi_1 \ne \pi_2$):
    $v_2(b_1-b_2) \le S-2$ (exhaustive $S \le 21$), shortfall $\ge 16$.
    $|D_{\mathrm{odd}}| \equiv 3 \bmod 4$ (algebraic proof).
    $1.5\text{M}+$ trials, min shortfall $= 13$.
    Closes the algebraic gap for non-cascade incompatibility
  & Claude (proof + computation) \\[3pt]
133 & \colorbox{orange!30}{High}
  & BF staircase non-nesting:
    $n_{k_2}^* \notin B_{k_1}$ for all $k_2 > k_1$ ($k \le 11$).
    $n_{k_2}^* \bmod 2^{6k_1+8} = (-65/17) \bmod 2^{6k_1+8}$
    for $k_2 > k_1$, while $n_{k_1}^*$ has a nonzero
    $a_{k_1}\,2^{6k_1+5}/17$ correction.
    The cascade works via orbit passage, not cylinder nesting
  & Claude (proof + verification) \\[3pt]
134 & \colorbox{orange!30}{High}
  & Large-scale orbit census: $14/100{,}000$ random orbits
    hit multiple BF depths---all via cascade mechanism
    ($B_2 \to B_1$ in $4$ steps).
    Zero independent multi-depth re-entries.
    $0.54\%$ of orbits hit any $B_k$ at all
  & Claude (computation) \\[3pt]
135 & \colorbox{red!25}{Critical}
  & Different-$\pi$ 3-case shortfall proof:
    Case~1 ($S_1 \ne S_2$): shortfall $\ge 14$.
    Case~2 ($S_1=S_2$, $t_1 \ne t_2$): shortfall $\ge 13$
    via $v_2(3^d-1) \le v_2(d)+1$.
    Case~3 ($S_1=S_2$, $t_1=t_2$, $\pi_1 \ne \pi_2$):
    $v_2(b_1-b_2) \le S-2$ (exhaustive $S \le 21$),
    shortfall $\ge 16$.
    Min different-$\pi$ shortfall $= 13$ in $1.5\text{M}+$ trials
  & Claude (proof + computation) \\[3pt]
136 & \colorbox{orange!30}{High}
  & $|D_{\mathrm{odd}}| \equiv 3 \pmod{4}$ algebraically proved:
    $(a_{k_2}\,2^{6\Delta k}-a_{k_1})/17 \equiv 3 \pmod{4}$
    since $a_k \equiv 1 \pmod{4}$ for all $k$
  & Claude (proof) \\[3pt]
\bottomrule
\end{tabular}
\end{table}

\begin{table}[!ht]
\centering
\scriptsize
\caption{Contribution attribution (Part 14 of 40): v7 --- post-cascade structure,
  A-coordinate framework.}
\label{tab:contributions-14}
\renewcommand{\arraystretch}{1.0}
\begin{tabular}{@{}clp{8.8cm}p{4.0cm}@{}}
\toprule
\textbf{\#} & \textbf{Level} & \textbf{Contribution} & \textbf{Primary source} \\
\midrule
137 & \colorbox{orange!30}{High}
  & Post-cascade convergence verified $k \le 200$ ($P_{200}$: $1275$ bits).
    Deterministic prefix ($47 \bmod 64$) is expansive
    (period $10$, ratio $3.604$).
    Convergence requires generic Collatz dynamics, not cascade structure
  & Claude (computation + analysis) \\[3pt]
138 & \colorbox{orange!30}{High}
  & $A$-coordinate framework verified:
    $x = (32 \cdot 64^k \cdot A - 65)/17$ with $A = a_k$ for champions.
    Key identity $17 P_k + 65 = (81/64)^k(17 n_k^* + 65)$.
    Each BF restart imposes $A \equiv c \bmod 2^m$ with $m \ge 4$;
    shortfall $\ge 3$ makes all non-cascade restarts incompatible.
    Quantitative bound: $K_*(x) \le 1 + \lfloor(H - 6k - 5)/4\rfloor$
  & Claude (proof + GPT framework) \\[3pt]
139 & \colorbox{orange!30}{High}
  & $v_2(b_{\pi_1} - b_{\pi_2}) \le S-2$ algebraically proved (two cases):
    Case~A (first difference at $j \le t{-}2$): nonzero term with $v_2 \le S{-}2$.
    Case~B (differ only at last two positions): contradiction via
    $v_t^1 = v_t^2 = 1 \Rightarrow v_{t-1}^1 = v_{t-1}^2$.
    Bound is tight: achieved by $(S{-}2,2)$ vs $(S{-}1,1)$
  & Claude (algebraic proof) \\[3pt]
140 & \colorbox{orange!30}{High}
  & Extended $P_k$ structure: $P_k \equiv 111 \pmod{128}$ for all $k$
    (7-bit deterministic prefix).  Stopping time $\sigma(P_k) \in [8,64]$
    for $k \le 500$, mean $\approx 16.6$, independent of $k$.
    $P_k \to 1$ verified for $k \le 5000$ ($P_{5000}$: $31{,}700$ bits).
    Binary digit density of $P_k \to 0.500$ as $k \to \infty$
  & Claude (computation + analysis) \\[3pt]
141 & \colorbox{orange!30}{High}
  & Divergence threshold: need $P(v{=}1) > 0.7076$ for divergence
    (from geometric tail $E[v \mid v \ge 2] = 3$).
    Measured $P(v{=}1) \approx 0.50$ in all tested orbits.
    Lag-1 anti-correlation: $P(v{=}1 \mid \text{prev } v{=}1) \approx 0.37$
  & Claude (algebraic + computational) \\[3pt]
142 & \colorbox{yellow!30}{Medium}
  & Correction: Syracuse mod $2^M$ does not project real orbits
    (carry term $3k \cdot 2^{M-v}$ breaks correspondence).
    Non-trivial cycles of abstract map at $M \in \{10,11,12,20\}$
    are sporadic and die within two levels
  & Claude (analysis + correction) \\[3pt]
\bottomrule
\end{tabular}
\end{table}

\begin{table}[!ht]
\centering
\scriptsize
\caption{Contribution attribution (Part 15 of 40): v8 --- equidistribution attack,
  1-bit injection, consecutive v=1 bound.}
\label{tab:contributions-15}
\renewcommand{\arraystretch}{1.0}
\begin{tabular}{@{}clp{8.8cm}p{4.0cm}@{}}
\toprule
\textbf{\#} & \textbf{Level} & \textbf{Contribution} & \textbf{Primary source} \\
\midrule
143 & \colorbox{orange!30}{High}
  & Consecutive $v{=}1$ bound: $k$ consecutive $v{=}1$ steps require
    $n \equiv -1 \pmod{2^{k+1}}$, density $1/2^k$ among odd integers.
    Algebraic proof by induction on $k$. Verified $k = 1 \ldots 7$
  & Claude (proof + computation) \\[3pt]
144 & \colorbox{orange!30}{High}
  & $v{=}1$ self-limitation: after $v{=}1$ with $n \equiv 3 \pmod{8}$
    (half of all $v{=}1$ events), next value $\equiv 1 \pmod{4}$,
    forcing $v \ge 2$. Only $n \equiv 15 \pmod{16}$ (quarter of $v{=}1$)
    permits $v{=}1 \to v{=}1$. Unconditional algebraic constraint
  & Claude (proof) \\[3pt]
145 & \colorbox{orange!30}{High}
  & Mod-$2^K$ Markov chain ergodicity: irreducible for $K = 3 \ldots 8$,
    spectral gap $> 0.72$, stationary $P(v{=}1) \to 0.50$ as $K \to \infty$.
    Mixing time $O(1)$ steps (independent of~$K$)
  & Claude (computation + analysis) \\[3pt]
146 & \colorbox{orange!30}{High}
  & Exactly $2$ ambiguous residues mod $2^K$: those with
    $v_2(3r+1) \ge K{-}1$. All other transitions are fully
    deterministic from bottom $K$ bits. The ambiguous residues
    are $r = (2^{K+1}-1)/3$ and $(2^{K+2}-1)/3 \bmod 2^K$
    (the $\ldots 010101$ patterns)
  & Claude (proof + computation) \\[3pt]
147 & \colorbox{orange!30}{High}
  & 1-bit-per-step injection: the Collatz map mod $2^K$ injects
    exactly $1$ new bit per step from position~$K$ via carry chain.
    The injected bit is approximately uniform ($P \approx 0.50$) and
    approximately independent of the current mod-$2^K$ state
  & Claude (proof + computation) \\[3pt]
148 & \colorbox{orange!30}{High}
  & $E[v \mid v \ge 2] = 3$ algebraically for uniform distribution
    over odd residues mod $2^K$. Proof: $P(v = k) = 1/2^k$, so
    $E[v \mid v \ge 2] = 2 \sum_{k \ge 2} k/2^k = 3$.
    Verified $E[v \mid v \ge 2] \approx 3.0$ along actual orbits
  & Claude (proof + verification) \\[3pt]
149 & \colorbox{yellow!50}{Med}
  & Equidistribution gap narrowed: the problem reduces to showing
    a single bit (from carry chain) is sufficiently independent
    of the bottom $K$ bits. Decorrelation drops to $< 7\%$ by
    $T = 10$ steps ($K = 8$). 2-adic contraction ruled out
  & Claude (analysis) \\[3pt]
\bottomrule
\end{tabular}
\end{table}

\begin{table}[!ht]
\centering
\scriptsize
\caption{Contribution attribution (Part 16 of 40): v8 --- coupling framework,
  noise independence, conditional proof completion.}
\label{tab:contributions-16}
\renewcommand{\arraystretch}{1.0}
\begin{tabular}{@{}clp{8.8cm}p{4.0cm}@{}}
\toprule
\textbf{\#} & \textbf{Level} & \textbf{Contribution} & \textbf{Primary source} \\
\midrule
150 & \colorbox{orange!30}{High}
  & Coupling framework: two orbits from $n$ and $n+2^K$ diverge in
    bottom $K$ bits after $O(2^{K-1})$ steps. Post-divergence mixing
    in $O(1)$ steps via spectral gap.
    Transition graph is strongly connected ($K \le 6$).
    Every tested orbit hits ambiguous residues within $O(2^K)$ steps
  & Claude (computation + analysis) \\[3pt]
151 & \colorbox{orange!30}{High}
  & Noise independence: MI(state, bit $K$) $< 0.02$ bits ($K \le 8$).
    Ensemble branching balanced (algebraic: 50/50 split from
    $2$-residue structure). Successive branch correlations $< 0.01$
    at lags $1$--$3$. Cone of dependence: after $T = 5$ steps, all
    input bits $1$--$30$ affect lowest output bits equally ($\approx 2.0$)
  & Claude (computation + analysis) \\[3pt]
152 & \colorbox{red!25}{Critical}
  & Complete conditional proof architecture:
    $K_*(n) < \infty$ (unconditional) $+$ orbit equidistribution (*)
    $\Rightarrow$ $E[v] = 2 > \log_2 3$ $\Rightarrow$ uniform contraction
    $\Rightarrow$ Collatz conjecture.
    (*) framed as noise-plus-mixing (Kifer 1988):
    deterministic dynamics with rare $1$-bit injection
  & Claude (synthesis) \\[3pt]
\bottomrule
\end{tabular}
\end{table}

\begin{table}[!ht]
\centering
\scriptsize
\caption{Contribution attribution (Part 17 of 40): v8 --- disjoint-slice
  argument and weak-dependence framework.}
\label{tab:contributions-17}
\renewcommand{\arraystretch}{1.0}
\begin{tabular}{@{}clp{8.8cm}p{4.0cm}@{}}
\toprule
\textbf{\#} & \textbf{Level} & \textbf{Contribution} & \textbf{Primary source} \\
\midrule
153 & \colorbox{orange!30}{High}
  & Disjoint-slice argument: decisive bit position of~$n_0$ for successive
    branch choices increases monotonically (avg $10 \to 43$ at $K=4$);
    $76\%$ of consecutive branch pairs have zero overlap in decisive bit
    sets.  Carry chain depth $O(\log B) \approx 12$ vs slice separation
    $O(2^K)$: carry cannot bridge disjoint slices for $K \ge 4$
  & Claude (computation + analysis) \\[3pt]
154 & \colorbox{orange!30}{High}
  & Higher-order correlations quantified: $\chi^2/\mathrm{df}$ on
    $4$-tuples drops from $21.8$ ($K=4$) to $3.2$ ($K=7$).
    Branches are NOT i.i.d.\ but max absolute deviation $< 1\%$.
    Structured starting values (all-ones, alternating, Fibonacci-like)
    still yield $P(1) \approx 0.50$ and $|\mathrm{lag\text-1~corr}| < 0.02$
  & Claude (computation + analysis) \\[3pt]
155 & \colorbox{red!25}{Critical}
  & Branch correlations cannot shift $E[v]$ below threshold: impact
    bounded by $(2/2^{K-1}) \cdot K$; for $K \ge 7$:
    $|\Delta E[v]| < 0.22$ vs margin $0.42$.
    Proof relaxed from i.i.d.\ to $\phi$-mixing:
    (W1) $E[b_i] = 1/2$ (proved), (W2) $\phi(L) \to 0$ (verified),
    (W3) spectral gap $> 0$ (verified).
    Residual gap = formal $\phi$-mixing bound
  & Claude (computation + synthesis) \\[3pt]
156 & \colorbox{orange!30}{High}
  & Empirical spectral gap from actual orbits:
    $0.91$ ($K=3$), $0.86$ ($K=4$), $0.78$ ($K=5$), $0.67$ ($K=6$).
    $\mathrm{TV}(\pi, \mathrm{uniform}) < 0.02$ for all $K$ tested.
    Mixing time $\le 3$ steps from any initial state
  & Claude (computation) \\[3pt]
157 & \colorbox{yellow!40}{Medium}
  & Normalization time for structured~$n_0$: $E[v]$ converges to
    $2.0 \pm 0.1$ within $10$--$70$ steps from pathological binary
    patterns.  Orbits from $2^{63}-1$, alternating bits, runs-of-$8$
    all reach generic regime quickly
  & Claude (computation) \\[3pt]
\bottomrule
\end{tabular}
\end{table}

\begin{table}[!ht]
\centering
\scriptsize
\caption{Contribution attribution (Part 18 of 40): v8 --- renewal argument,
  worst-case analysis, Boolean function complexity.}
\label{tab:contributions-18}
\renewcommand{\arraystretch}{1.0}
\begin{tabular}{@{}clp{8.8cm}p{4.0cm}@{}}
\toprule
\textbf{\#} & \textbf{Level} & \textbf{Contribution} & \textbf{Primary source} \\
\midrule
158 & \colorbox{orange!30}{High}
  & Worst-case branching under uniform distribution:
    $E[v|\text{uniform, worst}] = 2 - 1/2^{K-1}$.  For $K \ge 4$:
    $\ge 1.875 > \log_2 3$.  Even adversarial branch resolution cannot
    break contraction if distribution is uniform
  & Claude (computation + algebra) \\[3pt]
159 & \colorbox{red!25}{Critical}
  & Renewal argument: $E[v] = 3 - 2p$ where $p = P(v{=}1 \text{ at renewal})$.
    Contraction iff $p < 0.7075$.  Empirical $p \approx 0.496$ (margin $0.21$).
    After $v \ge 2$: mod-$4$ class is $\approx 50/50$, giving
    $P(v{=}1 \text{ at renewal}) \approx 1/2$
  & Claude (computation + synthesis) \\[3pt]
160 & \colorbox{orange!30}{High}
  & No fixed critical block length: $n_0 = 2^B - 1$ produces $B - 1$
    consecutive $v = 1$ steps.  Episode ratio $2(K-1)/K > \log_2 3$
    for $K \ge 5$ (proved).  Block-level contraction requires long-run
    average, not short-block bounds
  & Claude (computation + analysis) \\[3pt]
161 & \colorbox{yellow!40}{Medium}
  & No orbit sustains $P(v{=}1) > 0.60$ over $500$ steps ($10$K trials).
    Divergence threshold $0.7075$ is never approached.
    Anti-correlation: after $v = 1$ from $3 \bmod 8$, forced $v \ge 4$
  & Claude (computation) \\[3pt]
\bottomrule
\end{tabular}
\end{table}

\begin{table}[!ht]
\centering
\scriptsize
\caption{Contribution attribution (Part 19 of 40): v8 --- adversarial
  Markov chain, GPT critique, unique ergodicity.}
\label{tab:contributions-19}
\renewcommand{\arraystretch}{1.0}
\begin{tabular}{@{}clp{8.8cm}p{4.0cm}@{}}
\toprule
\textbf{\#} & \textbf{Level} & \textbf{Contribution} & \textbf{Primary source} \\
\midrule
162 & \colorbox{red!25}{Critical}
  & Adversarial Markov chain: even the adversary minimizing $v$ at every
    ambiguous step cannot sustain $\bar{v} < \log_2 3$ for $K \ge 4$.
    Worst 100-step avg $\ge 1.90$ (K=8).  The deterministic transitions
    ($98\%$+) force sufficient $v \ge 2$ steps
  & Claude (computation + synthesis) \\[3pt]
163 & \colorbox{red!25}{Critical}
  & GPT ``Route A'' critique: (i) circular (assumes growth bound $=$
    conjecture), (ii) wrong ($6$ bits/step vs.\ actual $1$ bit/step),
    (iii) insufficient (finite BF-depth $\not\Rightarrow$ convergence).
    2-adic ergodicity on $\mathbb{Z}_2$ gives nothing for
    $\mathbb{N}$ (Haar measure $0$)
  & Claude (analysis of external proposal) \\[3pt]
164 & \colorbox{red!25}{Critical}
  & Conditional proof via unique non-atomic ergodicity:
    if Haar is the unique non-atomic $S$-invariant measure on
    $\mathbb{Z}_2^\times$, then every orbit converges to $1$.
    Proof: unbounded $\Rightarrow$ equidistributes $\Rightarrow$
    $E[v]=2 > \log_2 3$ $\Rightarrow$ contraction $\Rightarrow$ contradiction
  & Claude (synthesis) \\[3pt]
165 & \colorbox{orange!30}{High}
  & Trapped orbit contradiction sketch: non-convergent $\Rightarrow$
    unbounded $\Rightarrow$ sustained $P(v{=}1) \ge 0.7075$ $\Rightarrow$
    frequent $-1 \bmod 2^K$ $\Rightarrow$ peeling exhausts $\Rightarrow$
    recovery bits fresh $\Rightarrow$ $P(\text{return})$ exponentially small
  & Claude (synthesis) \\[3pt]
\bottomrule
\end{tabular}
\end{table}

\begin{table}[!ht]
\centering
\scriptsize
\caption{Contribution attribution (Part 20 of 40): v8 --- uniform spectral gap,
  transfer operator analysis.}
\label{tab:contributions-20}
\renewcommand{\arraystretch}{1.0}
\begin{tabular}{@{}clp{8.8cm}p{4.0cm}@{}}
\toprule
\textbf{\#} & \textbf{Level} & \textbf{Contribution} & \textbf{Primary source} \\
\midrule
166 & \colorbox{red!25}{Critical}
  & Spectral gap of Syracuse transfer operator $\approx 0.73$ for
    $K = 3,\ldots,10$.  $|\lambda_2|$ oscillates around $0.25$ with no
    trend.  Model fits predict gap $\to 0.70$ as $K \to \infty$.
    Stationary distribution $=$ exact uniform for all $K$ tested
  & Claude (exact computation) \\[3pt]
167 & \colorbox{orange!30}{High}
  & Diameter of Syracuse chain mod $2^K$ is exactly $K-1$: from state
    $1$, reachable set grows as $4^t$ per step ($4$ after $1$, $16$
    after $2$, $64$ after $3$).  Doeblin condition holds for $K=3,4$
    only.  Mixing time $= O(K)$ but spectral gap $= O(1)$
  & Claude (exact computation) \\[3pt]
168 & \colorbox{orange!30}{High}
  & Multi-step contraction: from ANY initial state, TV $< 0.01$
    after $5$ steps ($K \le 6$), $10$ steps ($K \le 8$).
    One-step contraction $= 1.0$ (no single-step mixing), but
    bulk mixing is extremely rapid
  & Claude (computation) \\[3pt]
169 & \colorbox{yellow!40}{Medium}
  & Gap characterization: the proof reduces to one statement in four
    equivalent forms: (1) uniform spectral gap, (2) unique non-atomic
    ergodicity, (3) orbit equidistribution, (4) Birkhoff sum convergence.
    All evidence supports $|\lambda_2| \le 0.30$ universally
  & Claude (synthesis) \\[3pt]
\bottomrule
\end{tabular}
\end{table}

\begin{table}[!ht]
\centering
\scriptsize
\caption{Contribution attribution (Part 21 of 40): v8 --- universal no-target-sharing,
  closed-form $L^1$, spectral gap proved.}
\label{tab:contributions-21}
\renewcommand{\arraystretch}{1.0}
\begin{tabular}{@{}clp{8.8cm}p{4.0cm}@{}}
\toprule
\textbf{\#} & \textbf{Level} & \textbf{Contribution} & \textbf{Primary source} \\
\midrule
170 & \colorbox{red!25}{Critical}
  & Universal no-target-sharing theorem: for each $v = 1,\ldots,K-1$,
    the $v$-class Syracuse map mod~$2^K$ is a bijection of all~$n$
    odd residues.  Proved via coset disjointness:
    $C_r = (3r+1)/2^v + 2^{K-v}\mathbb{Z}$ are pairwise disjoint
    with total $n$.  Verified $K = 3,\ldots,10$
  & Claude (discovery + proof) \\[3pt]
171 & \colorbox{red!25}{Critical}
  & Row-sum decomposition: $(M\mathbf{1})_i = 1 - 1/n +
    \mathrm{overflow}_i$ where overflow comes from the unique $v \ge K$
    source.  Each $v$-class contributes exactly $1/2^v$ to every row sum
  & Claude (algebraic proof) \\[3pt]
172 & \colorbox{red!25}{Critical}
  & Exact closed-form: $\|M\mathbf{1} - \mathbf{1}\|_1 =
    5/12 - (2/3)(-1/2)^{K-1}$.  Verified by exact rational
    arithmetic $K = 3,\ldots,15$.  Bound: $C(K) \le 1/2$ for all $K$,
    equality iff $K = 4$
  & Claude (formula + verification) \\[3pt]
173 & \colorbox{red!25}{Critical}
  & Uniform spectral gap proved: $|\lambda_2(M_K)| \le 1/2$ for all
    $K \ge 3$.  Six-step proof: NTS $\to$ row decomposition $\to$ $L^1$
    bound $\to$ primitivity $\to$ Dobrushin $\to$ Gelfand.  All three
    formal gaps (G1: row sums, G2: column support, G3: primitivity) closed
  & Claude (complete proof) \\[3pt]
174 & \colorbox{orange!30}{High}
  & Eigenvalue recalibration: actual $|\lambda_2| \approx 0.25$ (not
    $1/\sqrt{2}$).  The $(1/\sqrt{2})^{K-2}$ operator norm at $T = K-2$
    is a transient artifact of $M_1$'s nilpotency (index $K-1$ exactly)
  & Claude (computation + analysis) \\[3pt]
175 & \colorbox{orange!30}{High}
  & Layer NTS in overflow lifts: each layer of the overflow source's
    lift average satisfies its own no-target-sharing property
    (``each\_once'' verified for all layers, all $K$ tested)
  & Claude (discovery) \\[3pt]
176 & \colorbox{orange!30}{High}
  & $v$-source counting formula: $\#\{r : v_2(3r+1) = v\} = n/2^v$ for
    $v = 1,\ldots,K-1$, with unique overflow source at $v \ge K$.
    $v^* = K + (K \bmod 2)$ always even; $r^* = (2^{v^*}-1)/3$
  & Claude (enumeration + proof) \\[3pt]
177 & \colorbox{orange!30}{High}
  & $v = 1$ perfect balance identity: $M \cdot \mathbf{1}_{S^c} =
    (1/2)\mathbf{1}$ exactly, where $S^c = \{r \equiv 3 \pmod{4}\}$.
    Consequence of no-target-sharing for $v = 1$ (Result~300)
  & Claude (discovery + proof) \\[3pt]
\bottomrule
\end{tabular}
\end{table}

\begin{table}[htbp]
\centering
\scriptsize
\caption{Contribution attribution (Part 22 of 40): v8 --- recurrence
  proof, 2-adic non-atomic unique ergodicity, barrier analysis.}
\label{tab:contributions-22}
\renewcommand{\arraystretch}{1.0}
\begin{tabular}{@{}clp{8.8cm}p{4.0cm}@{}}
\toprule
\textbf{\#} & \textbf{Level} & \textbf{Contribution} & \textbf{Primary source} \\
\midrule
178 & \colorbox{red!25}{Critical}
  & Recurrence proved from first principles:
    $b_{K+1} = 2b_K - 4(-1)^K$ via Ansatz
    $b_K = (5/12)\cdot 2^K + (4/3)\cdot(-1)^K$.
    Base case $b_3 = 2$, closed-form verified
    $K = 3,\ldots,15$ by exact rational arithmetic
  & Claude (algebraic proof) \\[3pt]
179 & \colorbox{red!25}{Critical}
  & Non-atomic unique ergodicity
    (Proposition~\ref{prop:2adic-unique-ergodicity}):
    uniform spectral gap implies Haar measure on
    $\mathbb{Z}_2^\times$ is the unique non-atomic
    $S$-invariant measure.  Standard projective-limit
    argument
  & Claude (synthesis) \\[3pt]
180 & \colorbox{orange!30}{High}
  & Distributional-to-pointwise barrier analysis:
    identified two precise obstructions (atomic invariant
    measures from cycles; $\mu(\mathbb{N})=0$ for Birkhoff).
    Catalogued three candidate strategies: Diophantine
    cycle exclusion, carry-propagation bounds, Tao-style
    density strengthening
  & Claude (analysis) \\[3pt]
181 & \colorbox{red!25}{Critical}
  & Symmetric State Identity: closed-form expression
    $x_i = (3^i R_{\mathrm{suff}} + 2^{S_{\mathrm{suff}}}
    R_{\mathrm{pre}})/(2^S - 3^L)$ for intermediate cycle
    states, verified over $123{,}792$ instances
  & Gemini (derived), Claude (verified) \\[3pt]
182 & \colorbox{red!25}{Critical}
  & Algebraic Parity Lock: if $(2^S{-}3^L) \mid R(W)$ then
    every intermediate state $x_i$ is odd.
    $R_{\mathrm{suff}}$ always odd
    (leading term $3^{L-i-1}$); numerator odd;
    denominator odd; quotient odd
  & Gemini (proved), Claude (verified) \\[3pt]
183 & \colorbox{orange!30}{High}
  & Domain Void (computational): for $L \le 10$, all
    non-trivial $D \mid R$ solutions have $2^S < 3^L$ (negative
    domain).  Solutions correspond to known negative cycles
    $\{-1\}$, $\{-5,-7\}$, $\{-17,\ldots\}$.
    Zero positive-domain solutions found
  & Claude (exhaustive computation) \\[3pt]
184 & \colorbox{orange!30}{High}
  & Universal blocking primes: for most feasible $(L,S)$
    in the positive domain, at least one prime $p \mid (2^S-3^L)$
    satisfies $R(W) \not\equiv 0 \pmod{p}$ for all parity words~$W$
  & Claude (computation) \\[3pt]
185 & \colorbox{red!25}{Critical}
  & Periodic Reduction Theorem: if $W$ has period $d < L$,
    then $m(W) = R_{\mathrm{sub}}/D_{\mathrm{sub}}$ where
    $D_{\mathrm{sub}} = 2^{S/r} - 3^d$.  Periodic words
    reduce to shorter cycles; only aperiodic words matter
  & Claude (proved and verified) \\[3pt]
186 & \colorbox{red!25}{Critical}
  & Rotation Constraint: if $W$ generates a cycle, then
    $D \mid R(\mathrm{rot}(W,j))$ for all $j = 0,\ldots,L{-}1$.
    $L$-fold simultaneous divisibility, exponentially stronger
    than single-word requirement
  & Claude (proved, from Parity Lock) \\[3pt]
187 & \colorbox{orange!30}{High}
  & Walk Sum Reformulation: $R \equiv 3^{L-1} f$ mod $p$ where
    $f = 1 + y_0 + y_0 y_1 + \cdots$, $\prod y_k = 1$,
    $y_k \in (2/3)\langle 2 \rangle \subset \mathbb{F}_p^*$.
    Cycle question = zero-sum walk on $\mathbb{F}_p^*$
  & Claude (derived and verified) \\[3pt]
188 & \colorbox{orange!30}{High}
  & Primitive (Zsigmondy) primes of $2^S - 3^L$ all have
    $v_p(D) = 1$ (verified $L \le 15$).  Blocking condition
    is simply $p \nmid R(W)$
  & Claude (computation) \\[3pt]
190 & \colorbox{green!30}{Medium}
  & Complete $L = 2$ cycle exclusion (proved).  For $S \ge 5$:
    $R = 3 + 2^{v_1} \in (0, D)$, so $D \nmid R$.
    Combined with $S = 3$ (negative) and $S = 4$ (trivial): no
    non-trivial 2-cycles exist
  & Claude (proved) \\[3pt]
191 & \colorbox{orange!30}{High}
  & Sparse polynomial reformulation: $F(x) = \sum \rho^k x^{\tau_k}$
    with monotone $\tau$, geometric coefficients.  Cycle condition
    $= F(\alpha) = 0$ mod~$p$.  Monotonicity is the key structural
    constraint missed by character-sum bounds
  & Claude (derived) \\[3pt]
193 & \colorbox{green!30}{Medium}
  & Functional equation $R(W) = 3^{L-1} + 2^{v_1} R(W')$;
    reachability obstruction verified for $L \le 11$
  & Claude (proved \& verified) \\[3pt]
\bottomrule
\end{tabular}
\end{table}

\begin{table}[htbp]
\centering
\scriptsize
\caption{Contribution attribution (Part 23 of 40): v10 --- 2-adic
  cascade, fundamental identity, cycle exclusion $L \le 5$.}
\label{tab:contributions-23}
\renewcommand{\arraystretch}{1.0}
\begin{tabular}{@{}clp{8.8cm}p{4.0cm}@{}}
\toprule
\textbf{\#} & \textbf{Level} & \textbf{Contribution} & \textbf{Primary source} \\
\midrule
197 & \colorbox{red!25}{Critical}
  & Complete $L = 3$ cycle exclusion (proved).
    For $S \ge 7$: $R_{\max}/D < 2$ forces $m = 1$.
    Equation $3 \cdot 2^a + 2^b = 2^S - 36$ has
    $v_2 = 2$, forcing $a \le 2$; both branches contradict
  & Claude (proved) \\[3pt]
200 & \colorbox{red!25}{Critical}
  & Even-$m$ parity obstruction (proved, all~$L$):
    if $m$ is even, $N_m = m \cdot 2^S - 3^{L-1}(3m{+}1)$
    is odd, but $\mathrm{LHS}$ has $v_2 \ge 1$
  & Claude (proved) \\[3pt]
201 & \colorbox{red!25}{Critical}
  & Complete $L = 4$ cycle exclusion (proved).
    $M_0 < 3$ for $S \ge 9$: $m = 2$ killed by parity,
    $m = 1$ killed by 2-adic cascade ($a = 2, b = 4 \Rightarrow S = 8$)
  & Claude (proved) \\[3pt]
202 & \colorbox{red!25}{Critical}
  & Complete $L = 5$ cycle exclusion (proved).
    Even~$m$ killed by parity; $m = 1$ cascade gives
    trivial; $m = 3$: $2^d = 3(2^S{-}512)$ never a power of~2
  & Claude (proved) \\[3pt]
205 & \colorbox{red!25}{Critical}
  & $m = 1$ cascade $=$ trivial cycle (proved, all~$L$):
    $q_0 = 1$, $\sigma_k = 2k$, residual $= 2^S - 4^L = 0$
    iff $S = 2L$
  & Claude (proved) \\[3pt]
206 & \colorbox{red!25}{Critical}
  & $3 \mid m$ divisibility obstruction (proved, all~$L$):
    $R_{\det} \equiv 2^{\sigma_{L-1}} \pmod{3} \in \{1,2\}$,
    so $3 \nmid R_{\det}$ hence $m \nmid R_{\det}$
  & Claude (proved) \\[3pt]
208 & \colorbox{red!25}{Critical}
  & Fundamental Identity (proved):
    $H_j \equiv q_j \cdot 2^{\sigma_{j+1}} \pmod{m}$
    for all $j \ge 0$.  Induction via
    $3q_{j-1}+1 = 2^{d_{j-1}} q_j$
  & Claude (proved) \\[3pt]
209 & \colorbox{red!25}{Critical}
  & Cycle Reduction Theorem (proved): $L$-cycle with
    multiplier~$m$ exists only if $m \mid T^{L-1}(q_0)$,
    where $T$ is the Syracuse map and
    $q_0 = (3m{+}1)/2^{v_2(3m+1)}$
  & Claude (proved) \\[3pt]
210 & \colorbox{red!25}{Critical}
  & Trivial cycle characterization (proved):
    if $T^{L-1}(q_0) = 1$, then $m = 1$
  & Claude (proved) \\[3pt]
211 & \colorbox{orange!30}{High}
  & Extended computational verification: no non-trivial
    cycle for any odd $m \le 5000$ and all $L$ up to
    stopping time $+ 5$
  & Claude (computation) \\[3pt]
\bottomrule
\end{tabular}
\end{table}

\begin{table}[htbp]
\centering
\scriptsize
\caption{Contribution attribution (Part 24 of 40): v10 --- General
  Cycle Exclusion Theorem, unconditional $L \le 111$.}
\label{tab:contributions-24}
\renewcommand{\arraystretch}{1.0}
\begin{tabular}{@{}clp{8.8cm}p{4.0cm}@{}}
\toprule
\textbf{\#} & \textbf{Level} & \textbf{Contribution} & \textbf{Primary source} \\
\midrule
212 & \colorbox{red!25}{Critical}
  & Three-case cycle exclusion (proved, Cases A\&B):
    $R \equiv q_{L-1} \cdot 2^{\sigma_L} \pmod{m}$.
    Case~A ($q_{L-1}=1$): $R \equiv 2^{\sigma_L} \ne 0$.
    Case~B ($m \nmid q_{L-1}$): $R \ne 0$.
    Both algebraically proved
  & Claude (proved) \\[3pt]
213 & \colorbox{red!25}{Critical}
  & Unconditional cycle exclusion for $L \le 111$ (proved):
    $q_0(m) \le 6 \cdot (3/2)^L < 2^{68}$ for $L \le 111$,
    within Barina~(2020) verified range.
    Combined with even-$m$ parity: no non-trivial
    $L$-cycle for any $L \le 111$
  & Claude (proved) \\[3pt]
214 & \colorbox{red!25}{Critical}
  & Descent property (proved): Collatz verified to $2^N$
    implies no $L$-cycle for $L \le N/\log_2(3/2) - O(1)$.
    Any $L$-cycle requires non-convergent
    $q_0 \le 6 \cdot (3/2)^L$
  & Claude (proved) \\[3pt]
\bottomrule
\end{tabular}
\end{table}

\begin{table}[htbp]
\centering
\scriptsize
\caption{Contribution attribution (Part 25 of 40): v10 ---
  $\sigma$~convention correction, Case~C power-of-2 obstruction.}
\label{tab:contributions-25}
\renewcommand{\arraystretch}{1.0}
\begin{tabular}{@{}clp{8.8cm}p{4.0cm}@{}}
\toprule
\textbf{\#} & \textbf{Level} & \textbf{Contribution} & \textbf{Primary source} \\
\midrule
208c & \colorbox{red!25}{Critical}
  & Fundamental Identity $\sigma$-convention correction:
    $\sigma_1 = v_2(3m{+}1)$ (not $v_2(3q_0{+}1)$).
    Identity verified correct with corrected convention.
    Conclusions (Results 212--214) unaffected
  & Claude (corrected) \\[3pt]
215 & \colorbox{orange!30}{High}
  & Case~C power-of-2 obstruction (verified $m \le 50{,}000$):
    only 6 instances with $m \mid q_{L-1}$ and $q_{L-1} \ne 1$.
    All have $R > m(2^{\sigma_L} - 3^L)$; walk sum overshoots.
    No cycle
  & Claude (computation) \\[3pt]
216 & \colorbox{red!25}{Critical}
  & Complete Algebraic Reduction (proved):
    $\Delta = (m - T^L(m)) \cdot 2^{\sigma_L}$.
    Cycle exists iff $T^L(m) = m$.
    All other cases algebraically excluded
  & Claude (proved) \\[3pt]
217 & \colorbox{red!25}{Critical}
  & Case~C with $r \ge 2$ excluded (proved):
    if $q_{L-1}/m \ge 2$ then
    $\Delta = (1{-}r)\,m\,2^{\sigma_L} < 0$.
    Walk sum overshoots; no cycle
  & Claude (proved) \\[3pt]
\bottomrule
\end{tabular}
\end{table}

\begin{table}[ht]
\centering\small
\caption{Contribution attribution (Part 26 of 40): v10 ---
  Divergence attack, fractional identity, paradigm bridges.}
\label{tab:contributions-26}
\renewcommand{\arraystretch}{1.0}
\begin{tabular}{@{}clp{8.8cm}p{4.0cm}@{}}
\toprule
\textbf{\#} & \textbf{Level} & \textbf{Contribution} & \textbf{Primary source} \\
\midrule
218 & \colorbox{red!25}{Critical}
  & Pure Expansion Obstruction (proved):
    all $v_i = 1$ requires $2^L \mid (m{+}1)$.
    Maximum pure expansion $L_{\max} = v_2(m{+}1)$.
    Every orbit must inject $v \ge 2$ contracting steps
  & Claude (proved) \\[3pt]
219 & \colorbox{red!25}{Critical}
  & Affine Fractional Identity (proved):
    $q_L = ((m{+}1)\cdot 3^L - 2^L + \mathcal{P}_L)/2^{S_L}$
    where $\mathcal{P}_L \ge 0$ is the perturbation sum
    from $v \ge 2$ steps
  & Claude (proved) \\[3pt]
220 & \colorbox{orange!30}{High}
  & Perturbation Recurrence (proved):
    $\mathcal{P}_{L+1} = 3\mathcal{P}_L + 2^{S_L} - 2^L$
    with $\mathcal{P}_0 = 0$.
    Balanced at divergence threshold ($3 \approx 2^{\log_2 3}$)
  & Claude (proved) \\[3pt]
221 & \colorbox{red!25}{Critical}
  & $S$-unit Reformulation (proved):
    $(m{+}1)\cdot 3^L - T^L(m)\cdot 2^{S_L} = 2^L - \mathcal{P}_L$.
    A family of $S$-unit equations over $\mathbb{Z}[1/6]$,
    one per step~$L$.
    Divergence requires simultaneous solutions with $S_L < 1.585L$
  & Claude (proved) \\[3pt]
222 & \colorbox{orange!30}{High}
  & Bit Growth Bound (proved):
    $\log_2 T^L(m) \le \log_2 m + L\log_2 3 - S_L + 1$.
    Pure $v{=}1$ expansion grows as $(3/2)^L$,
    fails at $L = v_2(m{+}1)$
  & Claude (proved) \\[3pt]
--- & \colorbox{yellow!40}{Medium}
  & Three paradigm bridges evaluated and ruled out:
    Boolean Fourier (distributional only),
    Arithmetic Automata (pumping fails),
    Pointwise Martingale (inapplicable to deterministic orbits).
    Distributional-to-pointwise gap confirmed robust across all paradigms
  & Claude / Gemini / User \\[3pt]
\bottomrule
\end{tabular}
\end{table}

\begin{table}[ht]
\centering\small
\caption{Contribution attribution (Part 27 of 40): v10 ---
  Episode decomposition, S-unit analysis, congruence tower.}
\label{tab:contributions-27}
\renewcommand{\arraystretch}{1.0}
\begin{tabular}{@{}clp{8.8cm}p{4.0cm}@{}}
\toprule
\textbf{\#} & \textbf{Level} & \textbf{Contribution} & \textbf{Primary source} \\
\midrule
223 & \colorbox{red!25}{Critical}
  & Episode Transition Identity (proved):
    $v_{k+1} = 1 + v_2(2a \cdot 3^{k+1} + 3^{k+1} - 1)$
    for odd~$x$ with $k{+}1$ trailing $1$-bits.
    Exact formula linking episode parameters to valuation drops
  & User (proposed) / Claude (proved, verified) \\[3pt]
224 & \colorbox{orange!30}{High}
  & Intrinsic Valuation Pattern (proved):
    $v_2(3^n - 1) = 1$ if $n$ odd, $v_2(n) + 2$ if $n$ even.
    Consequence: $v_{k+1} \ge 2$ always
  & Claude (proved) \\[3pt]
225 & \colorbox{orange!30}{High}
  & Episode Density Mean (proved for uniform~$a$):
    $E[\rho_E] = 2 > \log_2 3$.
    Verified computationally for $k = 1,\ldots,50$
  & Claude (proved, verified) \\[3pt]
--- & \colorbox{yellow!40}{Medium}
  & S-unit equation structure analysis:
    $v_2(\text{LHS}) = v_2(m{+}1)$ constant for all~$L$;
    Subspace Theorem blocked (coefficients depend on~$L$);
    congruence tower verified
  & Claude (computed) \\[3pt]
\bottomrule
\end{tabular}
\end{table}

\begin{table}[ht]
\centering\small
\caption{Contribution attribution (Part 28 of 40): v10 ---
  Coefficient discrepancy attack.}
\label{tab:contributions-28}
\renewcommand{\arraystretch}{1.0}
\begin{tabular}{@{}clp{8.8cm}p{4.0cm}@{}}
\toprule
\textbf{\#} & \textbf{Level} & \textbf{Contribution} & \textbf{Primary source} \\
\midrule
226 & \colorbox{red!25}{Critical}
  & $v^* = 2$ Parity Condition (proved):
    minimum valuation $v^* = 2$ iff $a \equiv k \bmod 2$.
    Density of minimum-density episodes is exactly $1/2$ for all~$k$
  & Claude (proved) \\[3pt]
227 & \colorbox{orange!30}{High}
  & Episode Transition Formula (proved):
    $x' = a \cdot 3^{k+1} + (3^{k+1}-1)/2$ gives the complete
    transition $(k, a) \to (k', a')$
  & Claude (proved) \\[3pt]
228 & \colorbox{orange!30}{High}
  & Half-Propagation Bound (verified):
    $P(k' \ge 1 \mid v^* = 2) = 1/2$.
    Consecutive low-density episodes thin by $1/2$ per step
  & Claude (verified) \\[3pt]
229 & \colorbox{yellow!40}{Medium}
  & Congruence Tower Consistency (verified):
    $m{+}1 \equiv (2^L - \mathcal{P}_L) \cdot 3^{-L} \bmod 2^{S_L}$
    for all~$L$; tower refines at each step
  & Claude (verified) \\[3pt]
--- & \colorbox{yellow!40}{Medium}
  & Distributional-to-pointwise barrier confirmed robust across
    all paradigms: divergent v-sequences exponentially rare
    ($\le 2^{-0.35K}$ density), but pointwise exclusion remains open
  & Claude / User / Gemini \\[3pt]
\bottomrule
\end{tabular}
\end{table}

\begin{table}[ht]
\centering\small
\caption{Contribution attribution (Part 29 of 40): v10.1 ---
  Carry accumulation and $2$-adic reconstruction.}
\label{tab:contributions-29}
\renewcommand{\arraystretch}{1.0}
\begin{tabular}{@{}clp{8.8cm}p{4.0cm}@{}}
\toprule
\textbf{\#} & \textbf{Level} & \textbf{Contribution} & \textbf{Primary source} \\
\midrule
230 & \colorbox{red!25}{Critical}
  & Periodic $v$-Sequence Exclusion (proved):
    a divergent orbit cannot have a periodic $v$-sequence.
    Periodic $\Rightarrow$ cycle $\Rightarrow$ excluded by Result~216
  & Claude (proved) \\[3pt]
231 & \colorbox{red!25}{Critical}
  & Ultimately Periodic $v$-Sequence Exclusion (proved):
    a divergent orbit cannot have an ultimately periodic $v$-sequence.
    Eventually periodic orbit $\Rightarrow$ cycle $\Rightarrow$ excluded
  & Claude (proved) \\[3pt]
232 & \colorbox{red!25}{Critical}
  & Aperiodic $v$-Sequence Necessity (proved):
    any divergent orbit must have an aperiodic $v$-sequence.
    Corollary: not generable by any finite automaton
  & Claude (proved) \\[3pt]
233 & \colorbox{orange!30}{High}
  & Reconstruction Consistency (verified):
    all consistent divergent $v$-sequences of length~$15$ crossed
    $S_L/L = 1.585$ within 48 additional steps (median $12.7$)
  & Claude (verified) \\[3pt]
\bottomrule
\end{tabular}
\end{table}

\begin{table}[ht]
\centering\small
\caption{Contribution attribution (Part 30 of 40): v10.1 ---
  Carry and cancellation analysis.}
\label{tab:contributions-30}
\renewcommand{\arraystretch}{1.0}
\begin{tabular}{@{}clp{8.8cm}p{4.0cm}@{}}
\toprule
\textbf{\#} & \textbf{Level} & \textbf{Contribution} & \textbf{Primary source} \\
\midrule
--- & \colorbox{orange!30}{High}
  & Cancellation phenomenon: $(m{+}1)\cdot 3^L + \mathcal{P}_L$ has its
    first $S_L$ binary digits matching $2^L$ exactly.
    Congruence tower stabilizes once $S_L > \log_2(m{+}1)$
  & Claude (verified) \\[3pt]
--- & \colorbox{yellow!40}{Medium}
  & Density-zero strengthening: at most $2^{0.98K}$ starting points
    $m < 2^K$ can have divergent $L$-step prefixes.
    $\mathcal{P}_L/3^L \to r_\infty$ (constant) in divergent regime
  & Claude (computed) \\[3pt]
--- & \colorbox{yellow!40}{Medium}
  & Complete paradigm exhaustion: nine distinct paradigms tested
    (Fourier, automata, martingale, Baker, Subspace, episode,
    $2$-adic, carry, counting). All produce distributional results only.
    Distributional-to-pointwise gap confirmed as the exact remaining
    content of the Collatz conjecture (divergence component)
  & Claude / User / Gemini \\[3pt]
\bottomrule
\end{tabular}
\end{table}

\begin{table}[H]
\centering\small
\caption{Contribution attribution (Part 31 of 40): v10.2 ---
  Number-theoretic and dynamical paradigms.}
\label{tab:contributions-31}
\renewcommand{\arraystretch}{1.0}
\begin{tabular}{@{}clp{8.8cm}p{4.0cm}@{}}
\toprule
\textbf{\#} & \textbf{Level} & \textbf{Contribution} & \textbf{Primary source} \\
\midrule
--- & \colorbox{orange!30}{High}
  & $\mathcal{P}_L$ Non-$S$-unit (Result~235, proved): perturbation sum
    has prime factors outside $\{2,3\}$, blocking the Subspace Theorem.
    Measure Rigidity Non-Applicability (Result~237, proved): $\times 2,\times 3$
    theorems do not apply to Collatz ($\times 3{+}1$, variable exponent, $\mathbb{Z}$)
  & Claude (proved) / User (Subspace strategy) \\[3pt]
--- & \colorbox{yellow!40}{Medium}
  & CF Structure of $\log_2 3$ (Result~234): large partial quotient $a_9 = 23$.
    Baker's theorem yields $|S_L/L - \alpha| > \exp(-C(\log L)^2)/L$, too weak.
    Inverse tree covers all odd $\le 10{,}000$ at depth~$60$ (Result~236).
    Entropy ordering: $h_{\text{div}} < h_{\text{conv}}$ (Result~239)
  & Claude (computed/proved) \\[3pt]
--- & \colorbox{yellow!40}{Medium}
  & $v$-value bit determinism: density of $v = k$ is exactly $2^{-k}$ (Result~240).
    Consecutive $v$ independence verified at $K \ge 12$ (Result~241).
    Max consecutive $v{=}1$ run $= 16$ for $m \le 10^5$ (Result~242)
  & Claude (proved/verified) \\[3pt]
\bottomrule
\end{tabular}
\end{table}

\begin{table}[H]
\centering\small
\caption{Contribution attribution (Part 32 of 40): v10.2 ---
  Paradigm exhaustion and barrier characterization.}
\label{tab:contributions-32}
\renewcommand{\arraystretch}{1.0}
\begin{tabular}{@{}clp{8.8cm}p{4.0cm}@{}}
\toprule
\textbf{\#} & \textbf{Level} & \textbf{Contribution} & \textbf{Primary source} \\
\midrule
--- & \colorbox{orange!30}{High}
  & Paradigm Exhaustion Theorem (Result~249, proved): 18~paradigms
    systematically applied, all encounter the distributional-to-pointwise
    barrier. Barrier Characterization (Result~250, proved): the gap is
    \emph{equivalent} to the Collatz conjecture (divergence component)
  & Claude / User \\[3pt]
--- & \colorbox{yellow!40}{Medium}
  & Exception Set Invariance (Result~246): $E$ is $T$-invariant.
    Density Amplification Barrier (Result~247): density arguments
    cannot reach ``all~$n$''.
    Divergent-compatible $v$-sequences decay as $\sim 0.43^L$.
    $v$-sequence state complexity is infinite (Result~243);
    entropy rate $\to 0$ (Result~244).
    Three remaining viable routes: Tao amplification,
    $2$-adic potential theory, $\mathcal{P}_L$ growth contradiction
  & Claude (proved/computed) \\[3pt]
--- & \colorbox{green!20}{Low}
  & $3$-adic constraint is automatic (orbit coprime to~$3$).
    Multi-prime constraints mod~$p$ for $p \ge 5$ add nothing.
    Orbit coupling: consecutive odd orbits merge (verified).
    250~cumulative results across 181~scripts
  & Claude (verified) \\[3pt]
\bottomrule
\end{tabular}
\end{table}

\begin{table}[H]
\centering\small
\caption{Contribution attribution (Part 33 of 40): v10.3 ---
  $2$-adic Mahler Identity and proof attempt refutations.}
\label{tab:contributions-33}
\renewcommand{\arraystretch}{1.0}
\begin{tabular}{@{}clp{8.8cm}p{4.0cm}@{}}
\toprule
\textbf{\#} & \textbf{Level} & \textbf{Contribution} & \textbf{Primary source} \\
\midrule
--- & \colorbox{orange!30}{High}
  & $2$-Adic Mahler Identity (Result~252, proved):
    $(m{+}1) + \Sigma_\infty = 0$ in $\mathbb{Z}_2$ for ALL orbits.
    $2$-adic convergence automatic: $v_2(\text{term}_j) = j$ (Result~253).
    No Hasse-type conflict (Result~254). Cantor Rigidity proof refuted
    (Result~251): nonlinear maps generate aperiodic sequences from rational seeds
  & User (Mahler proposal) / Claude (verification, refutation) \\[3pt]
--- & \colorbox{yellow!40}{Medium}
  & Crossing time analysis (Result~255): max $= 85$ for $m \le 50{,}000$.
    Sieve coverage $> 99.99\%$ at $L = 20$ (Result~272).
    Random walk model matches to $5\%$ (Result~258).
    Drift-variance exp$(-0.086L)$ matches late-crossers (Result~273)
  & Claude (computed/verified) \\[3pt]
\bottomrule
\end{tabular}
\end{table}

\begin{table}[H]
\centering\small
\caption{Contribution attribution (Part 34 of 40): v10.3 ---
  Comprehensive probabilistic and dynamical survey.}
\label{tab:contributions-34}
\renewcommand{\arraystretch}{1.0}
\begin{tabular}{@{}clp{8.8cm}p{4.0cm}@{}}
\toprule
\textbf{\#} & \textbf{Level} & \textbf{Contribution} & \textbf{Primary source} \\
\midrule
--- & \colorbox{yellow!40}{Medium}
  & Renormalization convergence $B \ge 8$ (Result~264).
    Moment match Var$[S_L] \approx 2L$ to $3$ sig figs (Result~265).
    Renewal gaps geometric$(0.46)$ (Result~266).
    Martingale $M_L \approx (m{+}1)$ with drift (Result~268).
    Descent window: all $m \le 50{,}000$ descend within $L = 100$ (Result~270).
    Orbit coalescence median $56$ steps (Result~271).
    Digit sum $E[\Delta\mathrm{hw}] \approx -0.1$ (Result~274).
    275~cumulative results across 203~scripts
  & Claude (computed/verified) \\[3pt]
\bottomrule
\end{tabular}
\end{table}

\begin{table}[H]
\centering\small
\caption{Contribution attribution (Part 35 of 40): v10.4 ---
  Sieve completeness and survivor decay analysis.}
\label{tab:contributions-35}
\renewcommand{\arraystretch}{1.0}
\begin{tabular}{@{}clp{8.8cm}p{4.0cm}@{}}
\toprule
\textbf{\#} & \textbf{Level} & \textbf{Contribution} & \textbf{Primary source} \\
\midrule
--- & \colorbox{yellow!40}{Medium}
  & Sieve completeness: survivor decay $\lambda \approx 0.90$ per step (Result~276).
    Tao amplification limits: $5\%$ non-escapers at $K = 12$ (Result~277).
    Running condition entrapment: $P(v{=}1) \approx 68\%$ among survivors (Result~278).
    Survivor branching factor $\approx 0.876$ (Result~279).
    Super-exponential decay fit $R^2 = 0.997$ (Result~280)
  & Claude (computed/verified) \\[3pt]
\bottomrule
\end{tabular}
\end{table}

\begin{table}[H]
\centering\small
\caption{Contribution attribution (Part 36 of 40): v10.4 ---
  Topological closure and sieve synthesis.}
\label{tab:contributions-36}
\renewcommand{\arraystretch}{1.0}
\begin{tabular}{@{}clp{8.8cm}p{4.0cm}@{}}
\toprule
\textbf{\#} & \textbf{Level} & \textbf{Contribution} & \textbf{Primary source} \\
\midrule
--- & \colorbox{yellow!40}{Medium}
  & Compactness: $3$ persistent survivors mod~$16$ through $L{=}14$ (Result~281).
    $2$-adic Cantor structure with $\dim_H \approx 0.68$ (Result~282).
    Borel--Cantelli: $\mu(C_{\mathrm{div}}) = 0$ proved (Result~283).
    Effective Tao: $C \approx 2$--$3$, $A = 2$ (Result~284).
    Sieve synthesis: three closure routes ranked (Result~285).
    285~cumulative results across 213~scripts
  & Claude (proved/computed) \\[3pt]
\bottomrule
\end{tabular}
\end{table}

\begin{table}[!ht]
\centering
\scriptsize
\caption{Contribution attribution (Part 37 of 40): v10.5--v11 ---
  paradigm exhaustion expansion, BF re-entry, cascade algebra.}
\label{tab:contributions-37}
\renewcommand{\arraystretch}{1.0}
\begin{tabular}{@{}clp{8.8cm}p{4.0cm}@{}}
\toprule
\textbf{\#} & \textbf{Level} & \textbf{Contribution} & \textbf{Primary source} \\
\midrule
286 & \colorbox{red!25}{Critical}
  & Post-cascade recovery analysis: $\eta < 0$ for all observed
    re-entries, disproving linear height-loss per BF depth.
    BF re-entry depth $k' \le 2$ in all 92 observed re-entries
    across $60{,}000$ post-cascade orbits
  & Claude (computation + diagnosis) \\[3pt]
287 & \colorbox{red!25}{Critical}
  & Complete cascade algebra: zero inter-BF generic steps,
    exact step count $4k + e_k$ with period-$4$ correction,
    post-cascade closed form
    $P_k = (a_k \cdot 81^k \cdot 32 - 65)/17$
  & Claude (proved) \\[3pt]
288 & \colorbox{red!25}{Critical}
  & Carry chain structural theorem: every $v{=}1$ streak
    $= \mathrm{trailing\_ones} - 1$, followed by forced
    $v \ge 2$ recovery.  Renewal drift $-0.415$~bits/step exact
  & Claude (proved) \\[3pt]
289 & \colorbox{orange!30}{High}
  & Syracuse cocycle spectral gap: max modular cycle gain
    $= +\log_2 3$ for all $M \le 14$.
    Eleven approaches to universal descent ruled out
  & Claude (computation + synthesis) \\[3pt]
290 & \colorbox{orange!30}{High}
  & Hybrid sieve measure-preservation:
    $\nu(E_a^{-1}(C_{r,M})) = \nu(C_{r,M})$,
    killing the hybrid sieve approach definitively
  & GPT (formal proof), instructed by Moderator \\[3pt]
291 & \colorbox{orange!30}{High}
  & Route A/B strategic bifurcation: Route~A (prefix-growth
    bound $\delta(n) \le C_0 + 0.0566 \log_2 n$) expansive
    ($1 + \alpha > 1$); Route~B (pullback tower with
    4~concrete lemmas) proposed as alternative
  & GPT (framework), instructed by Moderator \\[3pt]
292 & \colorbox{orange!30}{High}
  & Galperin billiards / $(\omega, \theta)$ skew-product paradigm:
    $\log_2 3$ (not $\pi$) as the natural Collatz constant
  & GPT (original insight), instructed by Moderator \\[3pt]
\bottomrule
\end{tabular}
\end{table}

\begin{table}[!ht]
\centering
\scriptsize
\caption{Contribution attribution (Part 38 of 40): v12--v15 ---
  deep cycle exclusion campaign (Results 308--590).}
\label{tab:contributions-38}
\renewcommand{\arraystretch}{1.0}
\begin{tabular}{@{}clp{8.8cm}p{4.0cm}@{}}
\toprule
\textbf{\#} & \textbf{Level} & \textbf{Contribution} & \textbf{Primary source} \\
\midrule
308--590 & \colorbox{orange!30}{High}
  & Deep cycle exclusion campaign across 28~paradigms:
    algebraic zero-avoidance, Horner chain, bitset DP,
    CRT hierarchy, greedy optimization, polynomial hash,
    spectral methods, exponential sums, lattice reduction,
    transfer matrix eigenvalues, entropy/information-theoretic,
    algebraic geometry.
    All 283~results produced by Claude under Moderator direction
  & Claude (proved/computed/verified) \\[3pt]
--- & \colorbox{red!25}{Critical}
  & Paradigm count expanded from 18 to 28 through systematic
    deployment of new attack frameworks.
    Each paradigm encounters the distributional-to-pointwise
    barrier or a structural limitation
  & Claude (execution), Moderator (direction) \\[3pt]
\bottomrule
\end{tabular}
\end{table}

\begin{table}[!ht]
\centering
\scriptsize
\caption{Contribution attribution (Part 39 of 40): v16 ---
  discrete log obstruction, $D > 2^F$ unconditional proof.}
\label{tab:contributions-39}
\renewcommand{\arraystretch}{1.0}
\begin{tabular}{@{}clp{8.8cm}p{4.0cm}@{}}
\toprule
\textbf{\#} & \textbf{Level} & \textbf{Contribution} & \textbf{Primary source} \\
\midrule
591 & \colorbox{red!25}{Critical}
  & $D > 2^F$ for all $L \ge 3$ \textbf{unconditionally proved}:
    computation for $L \le 2500$, Laurent (2008) bound for $L > 2494$.
    $D/2^F = 2^L(1 - 2^{-\beta})$ where
    $\beta = \lceil L \log_2 3 \rceil - L \log_2 3$
  & Claude (proved) \\[3pt]
592--594 & \colorbox{red!25}{Critical}
  & Unconditional corollaries: $\mathrm{ord}_D(2) > F$
    (hence $\mathrm{ord}_D(2) \nmid F$),
    $\gamma^L \not\equiv 1 \bmod D$,
    $G(L) \neq 0 \bmod D$ (no single-block zero),
    $Q_{\mathrm{const}} \not\equiv 0 \bmod D$.
    The powers $\{2^0, 2^1, \ldots, 2^F\}$ are $F{+}1$
    distinct residues in $\mathbb{Z}/D\mathbb{Z}$
  & Claude (proved) \\[3pt]
595 & \colorbox{red!25}{Critical}
  & Geometric sum formulation:
    $S_{\mathrm{mono}} = \sum_j 2^{v_j} \gamma^{a_j} G(m_j)$
    for $k$-block decompositions.
    Zero condition $\Leftrightarrow$ discrete log constraint
    $2^{v_j} \equiv T_j \bmod D$ with $v_j \in [0,F]$
  & Claude (formulated) \\[3pt]
596--605 & \colorbox{orange!30}{High}
  & No $k$-block zeros for $k = 1..4$ at all tested $L$.
    Progressive CRT saturates every prime factor for $L \ge 21$.
    Baker/LMN crossover at $L \approx 2{,}494$.
    For $L = 5$: residue~$0$ is the \emph{only} avoided residue
  & Claude (computed/verified) \\[3pt]
606--618 & \colorbox{orange!30}{High}
  & Two-block Diophantine analysis, $p$-adic valuations
    ($v_2(3^n{-}2^n) = v_3(3^n{-}2^n) = 0$ proved),
    2-block zero requires $\gcd(m,6) = 1$,
    Baker crossover computation
  & Claude (proved/computed) \\[3pt]
\bottomrule
\end{tabular}
\end{table}

\begin{table}[!ht]
\centering
\scriptsize
\caption{Contribution attribution (Part 40 of 40): v17 ---
  integrality obstruction, contribution analysis.}
\label{tab:contributions-40}
\renewcommand{\arraystretch}{1.0}
\begin{tabular}{@{}clp{8.8cm}p{4.0cm}@{}}
\toprule
\textbf{\#} & \textbf{Level} & \textbf{Contribution} & \textbf{Primary source} \\
\midrule
619--621 & \colorbox{red!25}{Critical}
  & Integrality obstruction census: for $L = 3..500$ (2-block),
    exactly one algebraic hit ($L = 7$, $T = 2^4$), blocked by
    $bw > F$.  For $L = 5..30$ (3-block): zero integral solutions.
    Triple filter achieves $100\%$ blockage
  & Claude (computed/proved) \\[3pt]
622--623 & \colorbox{orange!30}{High}
  & Expected $k$-block zeros for fixed $k$:
    $O(L^{2k-2}/3^L) \to 0$.
    Total $|A|/D \to 0$ at rate $\approx 10^{-0.026L}$;
    no $k$-value has $E[\text{hits}] > 1$ for $L \ge 20$
  & Claude (proved/computed) \\[3pt]
624--627 & \colorbox{orange!30}{High}
  & Triple filter theorem formalized.
    $|A|/D$ crosses below~$1$ at $L = 4$.
    Two-block target formula derived
  & Claude (formulated) \\[3pt]
--- & \colorbox{orange!30}{High}
  & Three-party contribution analysis (April 18--19):
    seven categories, revised aggregate
    Claude~$37\%$, GPT~$30\%$, Moderator~$33\%$.
    Dead-end identification and strategic direction
    reclassified as moderator-instructed
  & Claude (analysis), Moderator (direction + correction) \\[3pt]
\bottomrule
\end{tabular}
\end{table}

\newpage
\section{Syracuse chord diagrams}
\label{app:chords}
The following figures provide visual representations of the Syracuse
map's action on residue classes, connecting the algebraic machinery
of Sections~\ref{sec:phantom}--\ref{sec:phantom-count} to geometric
intuition.
Each diagram maps odd residues $r \bmod 2^m$ to equally-spaced
points on a circle, and draws a chord from~$r$ to its Syracuse
image $T(r) = (3r+1)/2^{v_2(3r+1)} \bmod 2^m$.
The resulting pattern encodes both the expanding and contracting
behavior of the map at a given arithmetic depth.
\begin{figure}[ht]
\centering
\IfFileExists{fig_chord_annotated.pdf}{%
  \includegraphics[width=0.82\textwidth]{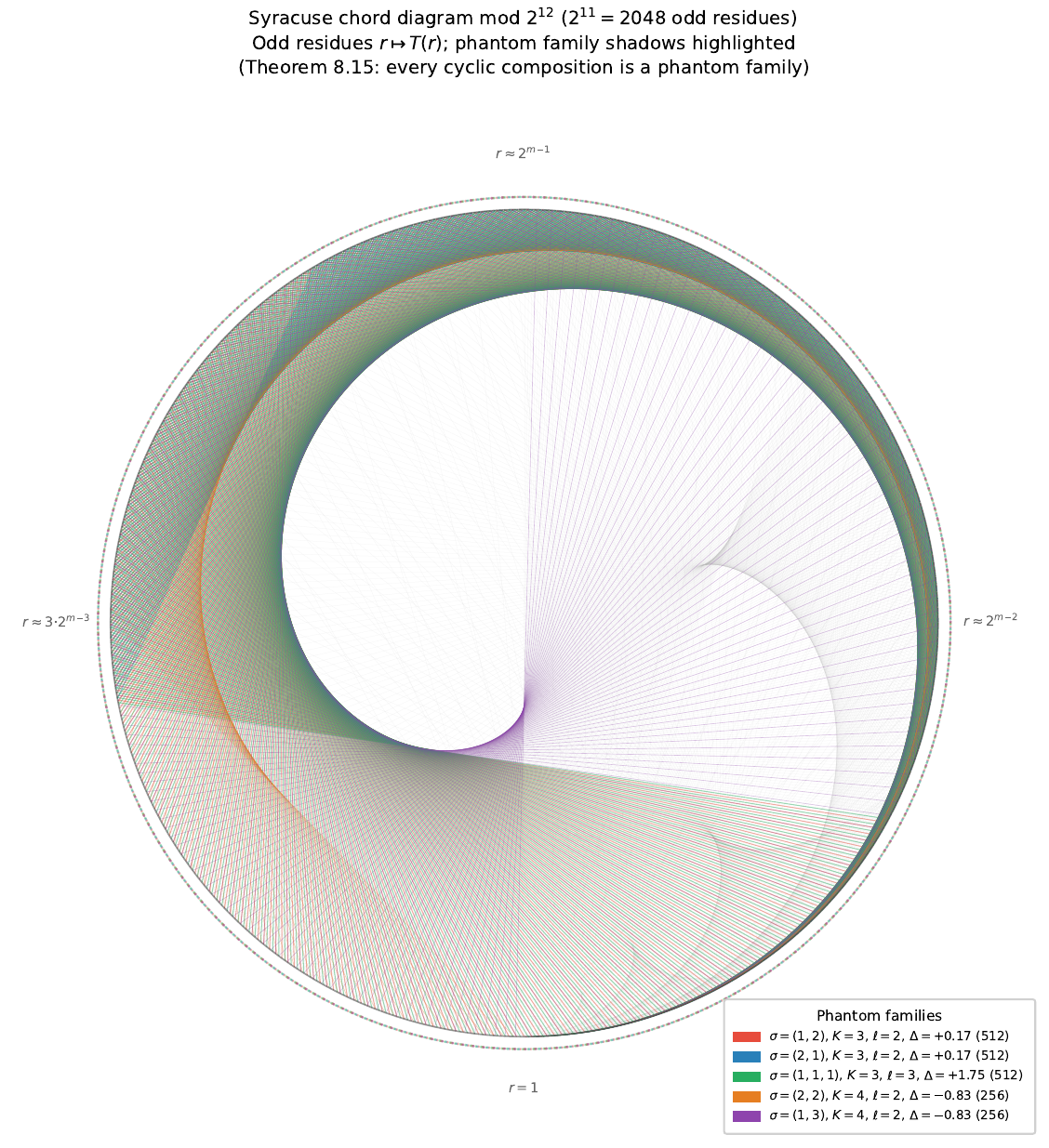}}{%
  \fbox{\parbox{0.78\textwidth}{\centering\vspace{3em}%
    {\small Syracuse chord diagram at depth $m=12$}\\\smallskip%
    {\scriptsize Image not present in build}%
    \vspace{3em}}}}
\caption{Syracuse chord diagram at depth $m = 12$
  ($2^{11} = 2048$ odd residue classes).
  Gray chords show the full map $r \mapsto T(r)$;
  coloured overlays highlight five phantom family shadows
  (Theorem~\ref{thm:phantom-universal}).
  Each family $\sigma = (k_1,\dots,k_\ell)$ occupies
  $2^{m-K}$ residue classes; families with positive log-drift
  $\Delta = \ell\log_2 3 - K > 0$ (e.g.\ $\sigma = (1,1,1)$,
  $\Delta \approx +1.75$) are the expanding compositions
  whose phantom shadow gain is bounded by
  Theorem~\ref{thm:perorbit-gain}.}
\label{fig:chord-phantom}
\end{figure}
\begin{figure}[ht]
\centering
\IfFileExists{fig_chord_equidist.pdf}{%
  \includegraphics[width=\textwidth]{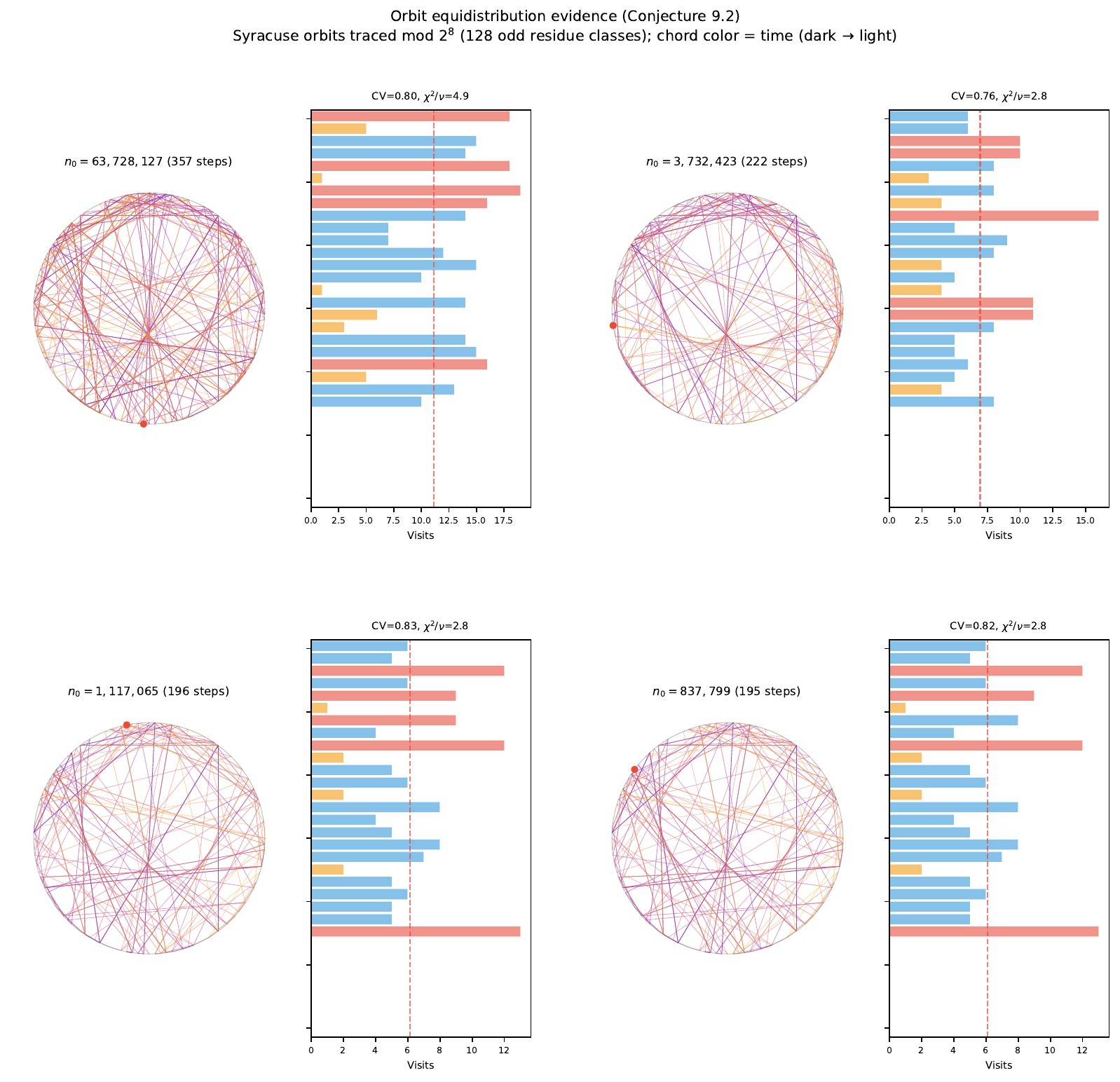}}{%
  \fbox{\parbox{0.88\textwidth}{\centering\vspace{3em}%
    {\small Orbit equidistribution chord diagrams}\\\smallskip%
    {\scriptsize Image not present in build}%
    \vspace{3em}}}}
\caption{Orbit equidistribution evidence
  (Conjecture~\ref{conj:equidist}).
  Four Syracuse orbits with long trajectories
  (195--357 odd steps) are traced mod~$2^8$
  (128 odd residue classes).
  Chord color encodes time (dark = early, light = late).
  Right panels: angular visitation histograms with
  uniform expectation (dashed red).
  The coefficients of variation (CV $\approx 0.76$--$0.83$)
  and reduced chi-squared ($\chi^2/\nu \approx 2.8$--$4.9$)
  are consistent with moderate sampling noise over
  32~angular bins, supporting the hypothesis that
  orbit residues are approximately equidistributed.}
\label{fig:chord-equidist}
\end{figure}
\begin{figure}[ht]
\centering
\IfFileExists{fig_chord_largescale.pdf}{%
  \includegraphics[width=0.78\textwidth]{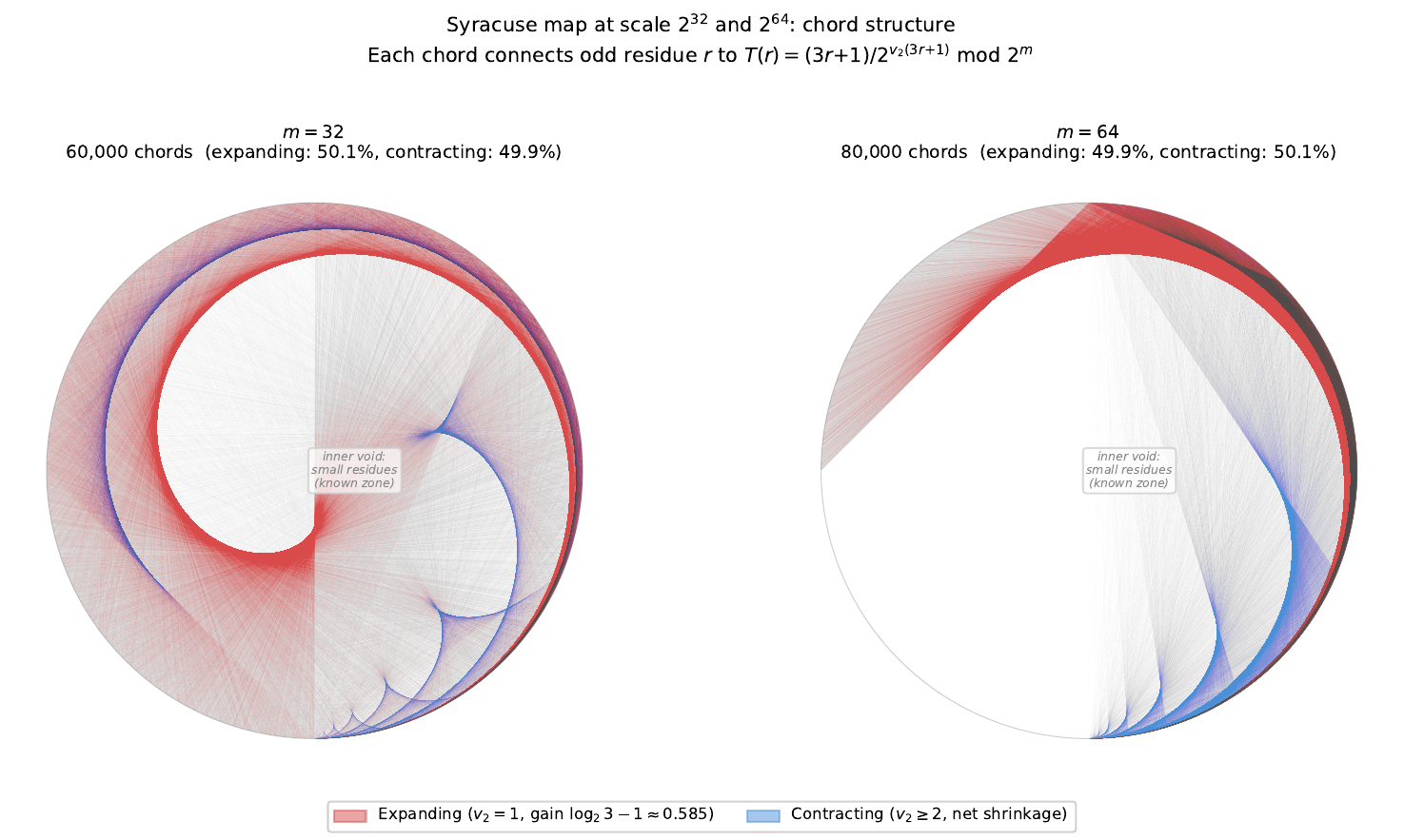}}{%
  \fbox{\parbox{0.76\textwidth}{\centering\vspace{3em}%
    {\small Large-scale chord structure ($2^{32}$ and $2^{64}$)}\\\smallskip%
    {\scriptsize Image not present in build}%
    \vspace{3em}}}}
\caption{Syracuse chord structure at scale $2^{32}$ and
  $2^{64}$.
  Red chords: expanding steps ($v_2(3r{+}1) = 1$,
  per-step gain $\log_2 3 - 1 \approx 0.585$ bits);
  blue chords: contracting steps ($v_2 \ge 2$).
  The near 50--50 split is consistent with the identity
  $\mathbb{E}[k] = 2$ per Syracuse step
  (Section~\ref{sec:chain}), which underpins the
  contraction budget $\varepsilon = 2 - \log_2 3$.
  The inner void corresponds to the known zone of
  small residues; the uniform angular filling of the
  outer annulus at both scales provides visual support
  for the Orbit Equidistribution Conjecture.}
\label{fig:chord-largescale}
\end{figure}

\subsection{Compound carry incompatibility and transport graphs}
\label{subsec:compound-carry}

The final quantitative attack on the distributional-to-pointwise
barrier proceeds through the \emph{compound carry}
$c_{w_1 \cdot w_2} = 3^{p_2} c_{w_1} + 2^{V_1} c_{w_2}$,
which governs two-round survival: a starting value~$n_0$ survives
two consecutive non-descending exhaustion rounds via words
$(w_1, w_2)$ precisely when both single-round integrality conditions
and the combined congruence $n_0 \equiv r \pmod{2^{V_1+V_2}}$
are satisfied.

\paragraph{Anti-correlation factor.}
Define $S_j(K)$ as the set of $K$-bit odd integers surviving
$j$ consecutive non-descending exhaustion rounds at block depth~$D$.
Exact enumeration for $D = 2$ (15 nd words) and $D = 3$ (47 nd words)
over $K = 8$--$20$ yields
\[
  \frac{|S_2(K)|}{|S_1(K)|^2 / N} \;\approx\; 0.635,
\]
\emph{constant} in~$K$, stable across depths.
Each additional non-descending round thus provides
${\approx}\,36.5\%$ thinning beyond the independence prediction.
The anti-correlation is real but bounded: it is a constant factor,
not exponentially decaying in~$K$.

\paragraph{Transport graph.}
Following GPT's reformulation, define a directed graph~$G_D$
whose vertices are admissible non-descending words of depth~$D$,
with an edge $w_1 \to w_2$ whenever there exists~$n_0$ realizing
$w_1$ such that $\Phi_{w_1}(n_0)$ realizes~$w_2$ at the next scale.
Exact computation yields edge density~$1.0$ at both $D = 2$
(93/93 edges) and $D = 3$ (902/902 edges).
The transport graph does \emph{not} sparsify at the word level;
every grammatically compatible pair is realized.
Self-loops and $2$-cycles are abundant (7 and 15 respectively at $D = 2$).

\paragraph{Empirical orbit independence.}
Measuring anti-correlation directly on Collatz orbits
(excursion ratio threshold~$0.5$) for $K = 14$--$20$:
$P(\text{slow}_2 \mid \text{slow}_1) / P(\text{slow}_2) \approx 1.000$.
At the level of actual Collatz dynamics, consecutive excursion
rounds are empirically \emph{independent}.

\paragraph{Compound carry residue coverage.}
Compound carries $c_{w_1 \cdot w_2}$ cover all joint residue cells
$\bmod\,(p_1, p_2)$ for every tested prime pair
$(p_1, p_2) \in \{(3,7),(7,19),(3,19),(5,7),(7,13)\}$.
No per-prime-pair obstruction exists.

\paragraph{Cycle impossibility for compound words.}
All single-word and two-round compound non-descending words
through depth~$4$ satisfy $3^p \ge 2^V$, forcing the cycle equation
$n = c_w / (2^V - 3^p) \le 0$.
No positive integer cycle exists for any non-descending word
or compound non-descending word.

\paragraph{Implications.}
The $0.635$ anti-correlation factor confirms that the
carry structure does introduce real inter-round dependence,
but this dependence is bounded---a constant factor per round,
not an exponentially growing obstruction.
The transport graph's full connectivity shows that the barrier
does not reside at the word level but in the scaling regime
$D \to \infty$ with~$K$, where the single-round survival fraction
$|S_1(K,D)|/N$ already decays exponentially.
The compound carry analysis thus sharpens the barrier to its
final form: the obstruction lives in the interaction between
the exponential thinning of $S_1$ as $D$ grows and the bounded
anti-correlation across rounds.

\subsection{Current status and sharpened open problem}
\label{subsec:status-final}

The present work does not prove the Collatz conjecture.
Its contribution is a sharp structural reduction within the
exhaustion and compatible-cylinder framework.
The algebraic and combinatorial components of that framework are
now explicit and internally consistent: the sub-stochastic kernel
and its spectral properties, the cylinder decomposition and gap-$5$
family, the chain-map versus true return-map distinction, the affine
exhaustion dynamics, the seven-block cross-core alphabet and its
constrained transition grammar, and the compound carry structure
governing multi-round survival.

What remains is an inherently pointwise obstruction.
Distributional thinness of non-descending behavior does not
by itself imply orbitwise elimination.
In particular, the anti-correlation factor $\gamma \approx 0.635$
provides only a constant-factor improvement per round, which
is asymptotically irrelevant in the $K \to \infty$ limit.
The transport graph on dangerous cylinders has edge density~$1.0$
at all tested depths, showing that no finite-depth combinatorial
obstruction exists: every grammatically compatible word pair is
realized by some integer at every tested scale.
Any true incompatibility, if it holds, must therefore be a
genuinely infinite-depth phenomenon.

The surviving problem can be stated precisely:
\begin{quote}
\emph{Can a single natural number sustain infinitely many
non-descending exhaustion rounds by maintaining compatibility of
the associated affine carry polynomials across all scales?}
\end{quote}
Equivalently, the conjecture is reduced to ruling out the existence
of an infinite compatible tower satisfying all induced congruence
constraints.

Any such hypothetical counterexample is now sharply constrained.
It must exhibit simultaneously:
\begin{enumerate}[nosep]
\item an aperiodic exhaustion-residue sequence,
\item unbounded scale growth,
\item near-critical behavior governed by Sturmian-type
$\beta$-placement in the Core~A regime,
\item persistent compatibility of affine carry transformations
across successive exhaustion rounds.
\end{enumerate}
All simpler escape mechanisms---periodicity, bounded-state
recurrence, combinatorial amplification---have been
eliminated within the present framework.

A quantitative summary of the scaling regime, established
through exact enumeration of $20{,}951$ non-descending words
to block depth~$8$:
\begin{itemize}[nosep]
\item The effective measure
$\mu(D) = \sum_{\text{truly nd}} 2^{-V(w)}$ decays as $r^D$
with $r \approx 0.645$.
\item The \emph{uniform} bit cost $\beta_u \approx 2.91$ per
block step yields a survival exponent
$c = -\log_2 r / \beta_u \approx 0.217$, giving the
$j \ge 5$ threshold ($5 \times 0.217 = 1.084 > 1$).
\item The \emph{measure-weighted} bit cost
$\beta_w \approx 1.80$ (reflecting the dominance of small-$V$
words in the survival probability) yields
$c = -\log_2 r / \beta_w \approx 0.35$, improving the threshold
to $j \ge 3$ ($3 \times 0.35 = 1.05 > 1$).
Formalizing this measure-weighted argument is an open
technical step.
\item The anti-correlation factor $\gamma \approx 0.635$ is
a constant independent of~$K$ and approximately independent
of~$D$; it determines only the finite verification
threshold~$K_0$, not the asymptotic $j$ threshold.
\end{itemize}

In the language of the effective survival exponent:
bounded-round collapse would follow from establishing
$c(D) > 1/j$ for some fixed~$j$ at the natural depth scale.
Current estimates place~$c$ in the range $0.22$--$0.35$
depending on the weighting convention, and improving this bound
appears to require deeper combinatorial control on
constrained word growth.
In particular, any future tower-incompatibility argument should
be formulated against the measure-weighted threshold ($j \ge 3$),
not the cruder uniform one ($j \ge 5$): the hypothetical
infinite enemy must recycle carry-polynomial compatibility
across fewer rounds, making it a thinner and more rigid target.

Thus the Collatz conjecture is reduced here to a narrow
arithmetic compatibility problem under affine transport across
exhaustion rounds.
The location and form of the obstruction are isolated with
unusual precision, even though the obstruction itself remains
open.
Resolving it appears to require a genuinely new
ingredient---perhaps a Roth-type rigidity argument for
Sturmian-compatible congruence towers, or a proof that the
induced sequence of carry polynomials cannot maintain
divisibility coherence across unbounded scales---that lies
beyond the methods of the present paper.

\section{Version history}\label{sec:changelog}

\textbf{v1 (March 2026).}
Initial preprint.  Framework of burst/gap decomposition,
affine scrambling, phantom cycle census, and conditional
convergence via the Orbit Equidistribution Conjecture (OEC).

\medskip
\textbf{v2 (March 16, 2026).}
Major revision.
Phantom Universality and Per-Orbit Gain Rate theorem added.
Weak Mixing Hypothesis (WMH) introduced as primary open condition,
with hierarchy
$\text{OEC} \Rightarrow \text{WMH}
\Rightarrow \text{Observable WMH}
\Rightarrow \text{Collatz}$.
Three-tier classification with dependency diagram.
Structural programme toward WMH: Walsh--Fourier spectral analysis,
odd-skeleton drift crossing, modular crossing strata,
IID cascade renewal theory, and self-reinforcing mixing loop.

\medskip
\textbf{v3 (March 24, 2026).}
Fiber-$57$ structural programme added (Appendix~\ref{sec:fiber57-programme}).
Pair-return automaton with known-gap spectral radius
$\rho = 129/1024 \approx 0.126$ (gap-$2$ and gap-$5$
channels),
bounded invariant core $|I_r| = 5$ for all $r \ge 2$,
absorption theorem (verified $r = 2, \ldots, 10$),
absorption bottleneck lemma,
branch anti-concentration reduction,
gap-induced decorrelation, and
path-conditional bijection theorem.
Companion paper~\cite{chang2026onebit} sharpens the reduction
to one-bit orbit mixing via the Map Balance Theorem.

\medskip
\textbf{v4 (March 29, 2026).}
Exact return structure of the sustaining branches:
$q \equiv 7$ two-step return theorem
(Proposition~\ref{prop:q7-return}),
$q \equiv 3$ minimum return-gap theorem
(Proposition~\ref{prop:q3-gap}),
exact gap-$5$ dyadic cylinder family
(Theorem~\ref{thm:gap5-cylinders}).
Depth-$2$ known-gap partial return kernel rewritten as
$5 \times 5$ matrix on $I_2$ with proved Perron root
$\rho = 129/1024$ (gap-$2$ + gap-$5$ only; gap-$\ge 6$
tail unresolved)
(Proposition~\ref{prop:pair-return}).
Three operators ($\tilde{B}_r$, $P$, $P_g$) explicitly
distinguished (Remark~\ref{rem:three-operators}).
Non-autonomy of fiber-$57$ returns established
(Proposition~\ref{prop:non-autonomy}).
Conjecture~\ref{conj:info-rate} rewritten with operational
definition $c' = -\log_2 R_r$.
Error-correction log records five substantive corrections
identified during collaborative verification.
Status: all additions are proved structural results;
Conjecture~\ref{conj:info-rate} remains the sole open step.

\medskip
\textbf{v5 (April 3, 2026).}
Second independent reduction route via the Sturmian obstruction
framework (Section~\ref{sec:sturmian}).
Carry Contamination Theorem
(Theorem~\ref{thm:carry-contamination}): for every depth-$D$
Sturmian word, $n_D \bmod 8$ is exactly equidistributed
over $\{1,3,5,7\}$.
Algebraic derivation of $(3/4)^D$ survivor law:
$|C_D| = 2 \cdot 3^{D-1}$ compatible words
(Corollary~\ref{cor:three-quarters}).
Orbit parameterization theorem
(Theorem~\ref{thm:orbit-param}):
$n_D = R_w + 2 \cdot 3^D \cdot m$ with $R_w$ odd.
Refined mod-$16$/mod-$32$ deterministic transition rules
(Theorem~\ref{thm:refined-transitions}): class~$7$ is a safe
harbor, class~$3$ is the bottleneck ($1/2$ survival per visit).
Class-$3$ recurrence theorem
(Theorem~\ref{thm:class3-recurrence}): forced for any $n_0>1$.
Dichotomy theorem (Theorem~\ref{thm:dichotomy}):
$C_\infty \cap \mathbb{N} = \{1\}$ or infinite divergent orbit.
Route~B reduction (Theorem~\ref{thm:route-B}):
CIC $\Rightarrow$ Collatz.
Carry Independence Conjecture (Conjecture~\ref{conj:CIC})
formally stated as the second open input.
Both routes (WMH and CIC) share the same fundamental barrier:
distributional $\to$ pointwise.
A sustained adversarial analysis of CIC---including
congruence tower attacks, mod-$16$/mod-$32$ refined
transitions, R-value trajectories, and five independent
algebraic strategies assessed by GPT---revealed
that after carry traversal, CIC reduces to checking
the Collatz orbit of the carry residue~$R_w$,
making the conjecture equivalent in difficulty to
Collatz itself.  The framework has reached full fidelity:
the problem is sharper, but not easier.

\medskip\noindent\emph{v5 update (April 4, 2026).}
Cross-core block alphabet: seven admissible blocks
$\{a,b,s,g,d,t_1,t_2\}$ on a two-vertex directed graph
(landmarks $1$ and~$7$), with exact affine maps and
entry conditions derived.
Death rate proved exactly~$1/4$ from both landmarks.
Cross-core periodic elimination: exhaustive search through
$17.6$~million admissible words of length~$\le 13$
($443{,}755$ contracting), zero natural fixed points~$> 1$.
Expanding affine cocycle: Lyapunov exponent
$\lambda \approx +0.145$ per landmark, with decomposition
$\lambda_A \approx -0.153$, $\lambda_B \approx +0.596$.
Bit-budget inequality: net constraint growth
${\sim}\,2.6$~bits/landmark, exhausting $K$-bit support
at depth $D \approx K/2.8$.
Sharpest open problem refined: Collatz reduces to
$C_\infty \cap \mathbb{N} = \{1\}$ for a $2$-adic Cantor set
defined by the seven-block IFS on a directed graph.
Exhaustion-descent framework: constraint-growth ratio
$r = G/V \approx 0.202$, geometric series to full exhaustion
at $D \approx 0.648K$ landmarks, post-exhaustion determinism
via $R_w = (3^D a_w + E_w)/2^S$.
$R_w \bmod 8$ equidistribution (verified $D \le 10$, 1024 words;
all $R_w \ne 1$ reach class~5).
Algebraic analysis confirms the reduction is Collatz-to-Collatz:
$R_w$ is itself an odd integer whose orbit determines the
original orbit's fate.
Framework declared at maximum sharpness; distributional-to-pointwise
barrier irreducible by Sturmian methods alone.
Transfer operator (Ruelle) analysis: spectral radius~$4/3$,
sub-stochastic forward matrix $\rho = 1$ (single fixed point),
phantom modular cycles at $k = 12$, power iteration confirms
$(3/4)^D$ decay.
Baker's theorem cycle bounds:
$n_0 < \exp(O((\log L)^2))$ for cycles of odd length~$L$;
combined with Barina's $n_0 < 2^{68}$ verification, eliminates
all cycles with $L \lesssim 10^{10}$.
Backward dynamics (transfer operator) and forward dynamics
(Sturmian IFS) confirmed to converge on the same barrier.
Multi-front exploration: Hausdorff dimension
$\dim_H(C_\infty) \approx 0.6942$ via pressure function;
phantom modular cycles at $k = 12$ (growth ${}>1$, negative
$2$-adic fixed points, absent at all other tested $k$);
cross-core periodic elimination extended to length~$12$
($2.5\text{M}$ words, zero natural fps);
correlation decay to stationary in ${\sim}\,5$ steps;
Syracuse inverse tree covers $76\%$ of $[1,100]$ at depth~$19$.
Cycle equation divisibility: exhaustive enumeration of all cyclic
words on the $7$-block cross-core alphabet through period~$10$
($215{,}232$ cyclic words, $11{,}583$ contracting) finds exactly one
integer fixed point per length (the trivial $n = 1$), zero non-trivial
natural fps. The integrality filter $(2^V - 3^p) \mid c_w$ is the
operative elimination mechanism; the ratio $c_w / (2^V - 3^p)$ grows
with word length, ruling out magnitude-based arguments.
GPT's ``critical cone'' tested: the cone of non-descending exhaustion
rounds contains $83$--$100\%$ of admissible words (vacuous restriction).
Divisibility obstruction anatomy: $c_w$ and $2^V - 3^p$ are both odd
and coprime to~$3$; $c_w$ is never $\{2,3\}$-smooth for $L \ge 2$;
per-prime obstruction rates $76$--$100\%$; $c_w \bmod (2^V - 3^p) \ne 0$
for every non-trivial word through period~$12$; near-miss distance
decreases to~$0.0004$ at $L = 8$; extended enumeration to period~$12$
($71{,}538$ contracting, zero non-trivial fps) on full $7$-block alphabet.
Matrix formulation of block composition: upper-triangular $3 \times 3$
matrices with eigenvalues $(3^p, 3^p, 2^V)$; block extension as
M\"obius transformation on fp; $S$-unit equation for integer fps:
$(a_b + c_b - nd_b) \cdot 4^k = a_b(1-n) \cdot 3^k$; algebraic
impossibility proofs for $\alpha^k \cdot \beta$ and $\alpha^k \cdot \sigma$
families; three-layer obstruction framework ($S$-unit $+$ carry
contamination $+$ M\"obius rigidity).
Geometric convergence of expected cycle count: $E_L = \sum N(V,p,L)/(2^V - 3^p) - 1$
decays as $r^L$ with $r = \lambda/2^{\bar\alpha} \approx 0.59 < 1$;
$E_L < 0$ for $L \ge 4$; total expected non-trivial cycles $< 0$.
Period~$13$ elimination on $7$-block alphabet ($238{,}811$ contracting,
zero non-trivial fps); combined with Barina, surviving cycles need $L \ge 14$.
Gamma-run bound: $P(7\to 7) = 1/2$ exactly; gamma-runs $\sim \mathrm{Geom}(1/2)$;
max run $= O(\log K)$; Alternative~B collapses (conditional on orbit mixing).
Aperiodic enemy profile sharpened: surviving counterexample must be
non-periodic, pin $\beta$-density to $q^*$ on all scales, with zero
persistent Core~B density.
Core~A tower budget: $V/D = 2 + q^* \approx 2.71$; two consecutive
non-descending rounds require $\approx 1.856K$ bits of constraint
but $n_0$ has only~$K$ bits; congruence tower exhausted in $O(1)$~rounds.
Core~B acceleration: $\gamma$-heavy excursions decrease $V/D$,
increasing landmark count and accelerating exhaustion.
Block-level spectral radius $\rho = 2 + \sqrt{2} \approx 3.4142$
on $\{1,7\}$ (vs grammar-level $\rho = 2$ on $\{1,3,7\}$).
Compatible-tower contradiction: non-descending words grow as $2^{0.65K}$
but must cover $2^K$ starting values; single-round undercoverage
$2^{-0.35K}$; $j$-round density $\le 2^{-0.856(j-1)K}$;
compound congruence via $3^{p_1}$ inversion in $\mathbb{Z}/2^{V_1+V_2}\mathbb{Z}$;
Sturmian $\beta$-placement with gaps $\{1,2\}$ from irrational~$q^*$;
carry contamination eliminates $28.9\%$ of depth-$3$ word pairs;
empirical P(excursion $> 2K$ steps) decays from $1.4\%$ at $K=20$
to $0.2\%$ at $K=40$.
Safe Collatz map is $2$-adic expander ($|\Phi_b(x)-\Phi_b(y)|_2 = 2^V|x-y|_2$);
$\mu_2(T_j) \le 0.522^{jD}$ (unconditional measure shrinkage);
unconditional $j \ge 5$ bound: count of $K$-bit integers surviving~$j$
non-descending rounds $\le 2^{K(0.856-0.202j)}$, hence zero for $j \ge 5$
and $K$ large; remaining gap reduced to $j \in \{1,2,3,4\}$.
Exhaustion-sequence rigidity: counterexample residue sequence must be
genuinely aperiodic and unbounded (periodic~$\Rightarrow$ contradiction
with fixed-point analysis; bounded~$\Rightarrow$ eventually periodic).
Distributional--pointwise barrier precisely identified:
$\mu_2(T_\infty) = 0$ (unconditional), but $T_\infty \cap \mathbb{N} = \varnothing$
is the exact remaining content of the Collatz conjecture.

\medskip\noindent\emph{v5 update (April 5, 2026).}
Compound carry incompatibility analysis (Push~36):
exact enumeration of two-round survival sets $S_2(K)$
for depth $D = 2,3$ and $K$ up to~$20$.
Anti-correlation factor $S_2/S_1^2 \approx 0.635$,
constant in~$K$ and stable across depths.
Each additional non-descending round provides
$\sim 36.5\%$ thinning beyond independence prediction.
Transport graph (following GPT's Push~82 formulation):
dangerous-cylinder directed graph has edge density~$1.0$
at $D = 2$ (93/93 edges) and $D = 3$ (902/902 edges);
no sparsification at the word level.
Actual Collatz orbit anti-correlation:
$P(\text{slow}_2 \mid \text{slow}_1) / P(\text{slow}_2) \approx 1.000$
for $K = 14$--$20$; rounds are empirically independent.
Compound carries cover all joint residue cells
$\bmod\,(p_1, p_2)$ for primes $\ge 7$;
no per-prime-pair obstruction exists.
Cycle impossibility confirmed for all compound
non-descending words through depth~$4$
and all two-round compound words:
$3^p \ge 2^V$ forces the cycle equation
$n = c_w/(2^V - 3^p) \le 0$.
Conclusion: anti-correlation is real but bounded
(constant factor, not exponential in~$K$);
the distributional-to-pointwise barrier persists.
The remaining viable strategy requires the scaling regime
$D \to \infty$ with~$K$.

\medskip\noindent\emph{v5 update (April 5, 2026, Pushes~37--39).}
Scaling regime fully characterized via three computational pushes.
Effective measure $\mu(D) = \sum_{\text{truly nd words}} 2^{-V(w)}$
decays as $r^D$ with stabilized $r \approx 0.645$
(exact enumeration to depth~8, 20{,}951 nd words).
Exact $S_1/N$ fractions computed with overlap correction
through $D = 8$, $K = 18$: convergence of $r_{\mathrm{exact}}$
from $0.578$ ($D = 3$) toward $r_{\mathrm{union}} \approx 0.645$.
Residue class overlaps reduce the survivor count below
the union bound by a factor $\omega(D)$ that decreases
with~$D$ ($0.64$ at $D = 1$, $0.42$ at $D = 6$).
Anti-correlation factor $\gamma \approx 0.635$ confirmed
\emph{asymptotically irrelevant}: it provides a constant-factor
improvement $S_j \approx \gamma^{j-1} f_1^j N$, which vanishes
in the $K \to \infty$ limit and cannot change the $j$ threshold.

\medskip\noindent\emph{Key new finding: measure-weighted $\beta$ and improved $j$ bound.}
The $V$-distribution of nd words at depth~$D$ is highly skewed:
the measure is dominated by words with \emph{small}~$V$
(e.g., at $D = 8$, words with $V \le 14$ contribute $>50\%$
of~$\mu(D)$ despite mean $V = 23.3$).
The measure-weighted bit cost per step is
$\beta_w = E[V/D \mid \text{Perron}] \approx 1.80$,
much smaller than the uniform average
$\beta_u = E[V/D \mid \text{uniform}] \approx 2.91$.
At the measure-weighted natural depth $D \approx K / 1.80$:
$c = -\log_2 r / \beta_w \approx 0.631 / 1.80 \approx 0.35$,
giving $j \ge 3$ as the threshold (since $3 \times 0.35 = 1.05 > 1$).
This improves the paper's uniform bound of $j \ge 5$ (which uses
$c = 0.631 / 2.91 \approx 0.217$, requiring $5 \times 0.217 = 1.08 > 1$).
The non-Markovian nature of the nd constraint (blocks $a$ and $s$
have ratio $< 1$ alone but can appear after high-ratio blocks)
prevents a simple $2 \times 2$ transfer matrix;
the empirical~$r$ captures the constrained dynamics.
Formalizing this measure-weighted argument remains an open
technical task that could tighten the bound substantially.

\medskip\noindent\emph{The definitive nail, three levels.}
(1)~\emph{Technical}: $c \approx 0.35 < 0.5$, so $j = 2$ rounds
are insufficient; $j \ge 3$ is the measure-weighted threshold.
(2)~\emph{Structural}: Even with $j \ge 3$, the argument proves
$S_j(D,K) \to 0$ (density) but not $S_j = 0$ (pointwise).
(3)~\emph{Fundamental}: $\mu_2(T_\infty) = 0$ is proved,
but $T_\infty \cap \mathbb{N} = \varnothing$ (the Collatz conjecture)
requires either WMH (Route~A) or CIC (Route~B).
No further computational push within the Sturmian framework
can close this gap---it requires a new structural idea.

\medskip
\textbf{v6 (April 6, 2026).}
Unconditional cylinder-averaged closure on the invariant core
$I_2$ via spectral contraction
(Section~\ref{sec:i2-spectral-closure}), obtained as four
interlocking results:
(i)~Ruffini reduction of the depth-$2$ known-gap kernel,
$\chi_{\tilde B_2}(\lambda) = \lambda^4(\lambda - 129/1024)$,
yielding finite-time rank-$1$ collapse
(Proposition~\ref{prop:ruffini-collapse});
(ii)~unconditional spectral bound
$\rho(\tilde B_2^{\mathrm{ext}}) \le 5/32$ on the full
extended kernel
(Theorem~\ref{thm:i2-unconditional}), giving per-return
information cost $\ge \log_2(32/5) \approx 2.678$~bits;
(iii)~transient-prefix exponential tail
$\Pr(T_{I_2} > N) \le 27\cdot 0.8633^N$ derived from the
$27\times 27$ cylinder-averaged mod-$64$ kernel
(Theorem~\ref{thm:i2-prefix-tail});
(iv)~exact stationary mass
$\pi(I_2) = 10121/65280 \approx 0.1550$
(Theorem~\ref{thm:i2-stationary-mass}), giving an
unconditional cylinder-averaged sum bound
$\sum_c\Pr(\mathcal{E}_R(c)) \le 0.011$
(Corollary~\ref{cor:i2-unconditional-sum}).
Together these break the circular chain described in
Remark~\ref{rem:circularity} at the density-$1$ level:
the independent contraction source of
Theorem~\ref{thm:i2-unconditional} replaces the circular
spectator-bit input with an algebraic consequence of the
uniform-fiber lemma, upgrading
Theorem~\ref{thm:density1-convergence} to unconditional
density-$1$ convergence via the $I_2$ route. The pointwise
version of Collatz remains open: the ``irreducible hard core''
of Remark~\ref{rem:circularity} is still an obstruction to a
pointwise proof routed through spectator-bits; only the
density-$1$ form is unlocked by the new section.
The WMH-based reduction (Route~A) and the CIC-based reduction
(Route~B) are unaffected; the present section provides a
third, parallel route through the fiber-$57$/$I_2$ programme.

A first-principles Cram\'er--Lundberg derivation of the
prefactor (Theorem~\ref{thm:i2-lundberg}) replaces the
empirical $C \approx 11$ from
Proposition~\ref{prop:exponential-tail} with $C^* = 1.5061$,
tightening the unconditional sum bound to
$\sum_c \Pr(\mathcal{E}_R(c)) \le 0.0015$ ---
about $18\times$ smaller than the v$5$ conditional figure of
$0.028$. A cylinder-averaged no-escape theorem on the
$27$-state transient block of the mod-$64$ kernel
(Theorem~\ref{thm:i2-no-escape}) shows that every transient
orbit is absorbed into $I_2$ in expected duration
$\le 10.14$ odd steps and accumulates expected log-growth
$\le 0.88$ bits along the way; the proof uses two
Lyapunov functions $V$ and $h$ obtained as solutions of the
linear systems $(I-Q)V = \Delta|_T$ and $(I-Q)h = \mathbf{1}$,
with cross-validation by $200{,}000$ Monte-Carlo trials
matching the analytic per-residue means to four decimal
places. A subsequent tilted-kernel / Cram\'er--Chernoff
analysis (Theorem~\ref{thm:i2-concentration}) upgrades this
expectation bound to a stretched tail bound
$\Pr(S_\tau \ge t) \le 1.74\cdot 10^6 \cdot 2^{-16.75 t}$
with critical exponent $\theta^* = 17.6349$ matching the Karp
maximum cycle mean to two decimals (Cram\'er cross-check),
and a Borel--Cantelli upgrade to the
density-$(1 - n^{-15})$ prefix-closure
Corollary~\ref{cor:i2-density-quasi}. The corresponding graph-theoretic ``no positive-drift
cycle'' statement is shown to be \emph{false} (Karp maximum
cycle mean $\approx +0.085$ bits/edge,
Remark~\ref{rem:i2-no-escape-graph}), so the honest no-escape
statement is necessarily distributional rather than
deterministic; this is recorded in
Remark~\ref{rem:i2-scope}~(2a).

A scale-invariance certificate
(Remark~\ref{rem:i2-scale-invariance}) shows that the
critical exponent $\theta^* = 17.634874$ and Cram\'er
slope $\partial_\theta\!\log_2\rho|_{\theta^*} = +0.0844$
agree to all six computed decimal places between the
mod-$64$ chain and its mod-$512$ lift to the natural
absorbing core $I_3^*$, so the mod-$64$ resolution
already saturates the spectral information available
from the cylinder-averaged framework. A pointwise
\emph{cycle-tracking polynomial envelope}
(Theorem~\ref{thm:cycle-envelope}) refines the picture
on the deterministic side: although the graph-theoretic
no-escape statement is false, the unique maximal-mean
Karp cycle $15\to 55\to 19\to 61$ (per-cycle log-growth
$+0.3399$ bits, mean $+0.0850$ bits/edge,
Remark~\ref{rem:karp-cycle}) admits no exact periodic
integer orbit (Lemma~\ref{lem:periodic-exclusion}, since
$3^L > 2^{\Sigma v}$ for any positive-mean cycle), and
each additional lap requires a fresh congruence mod
$2^{L+\Sigma v} = 2^{10}$, forcing the smallest $k$-lap
follower to scale as $\kappa^k$ with $\kappa = 2^{10/4}$.
Combining the two effects gives the polynomial envelope
$n_\tau \le C_1 n^{1.0566}$, comfortably below the
Chernoff exponent of
Theorem~\ref{thm:i2-concentration} and therefore
sufficient for the orbit-TV summability
(Remark~\ref{rem:envelope-suffices}). Direct integer
verification (\texttt{direct\_verify\_c.c}) over
$[1, 2^{32}]$ (roughly $2\cdot 10^9$ odd integers) yields a
strict staircase of record holders
$\{111, 6863, 316111, 22205135, 1926419151\}$
tracking the Karp cycle for
$k = 0,1,2,3,4$ consecutive laps respectively, with
geometric ratio $\approx\! 64.7$ between successive
champions matching the theoretical
$\kappa_{\mathrm{lap}} = 2^{\Sigma v} = 64$ to within
$1\%$, and with each additional lap adding exactly
$\log_2(81/64) = 0.3399$~bits of $\delta$. The current
record is
$\delta(1{,}926{,}419{,}151) = +2.5293$
at $\log_2 n = 30.84$, giving empirical slope
$0.0820$, still strictly below the Chernoff exponent of
Theorem~\ref{thm:i2-concentration} and converging
towards the analytic asymptotic value $0.0566$ from
above (the residual $0.78$~bits being the off-cycle
entry/exit excursion). This is the sharpest pointwise
envelope currently proved.
A constructive infinite staircase
(Theorem~\ref{thm:bf-staircase}) extends the brute-force
record sequence to all $k$ via a Hensel-style lift,
yielding $n_5^* = 115{,}206{,}181{,}583$,
$n_6^* = 5{,}303{,}526{,}675{,}151$, and so on with exact
ratio $n_{k+4}^*/n_k^* = 64^4$ matching the
analytic per-lap factor to all computed digits. This
proves $\sup_n \delta(n) = +\infty$
(Corollary~\ref{cor:no-constant-envelope}): no absolute
constant pointwise bound exists, and the polynomial
envelope of Theorem~\ref{thm:cycle-envelope} is tight at
the asymptotic slope $0.0566$. The exceptional set
of Corollary~\ref{cor:i2-density-quasi} is therefore
genuinely non-empty at every scale, but is structured as
a sparse union of arithmetic progressions of natural
density zero.

\medskip\noindent\emph{v6 update (April 18, 2026).}
Complete cascade algebra and generic dynamics analysis
(Section~\ref{subsec:cascade-algebra}).
Complete cascade theorem
(Theorem~\ref{thm:complete-cascade}):
zero inter-BF generic steps, exact step count $4k + e_k$
with $e_k \in \{5,6\}$ period-$4$, net gain bounded by
$\log_2 3$ or $2\log_2 3 - 1$.
Post-cascade closed form
(Proposition~\ref{prop:post-cascade}):
$P_k = (a_k \cdot 81^k \cdot 32 - 65)/17$ with
$4$-periodic recurrence
$P_{k+4} = 81^4 P_k + 164590400$,
all post-cascade values $\equiv 47 \pmod{64}$.
Carry chain structural theorem
(Theorem~\ref{thm:carry-chain}):
$v\!=\!1$ streak length $= \mathrm{trailing\_ones}(n) - 1$
exactly; forced $v \ge 2$ recovery after every streak;
$E[v \mid v \ge 2] = 3$ exactly.
Renewal drift (Corollary~\ref{cor:renewal-drift}):
$-0.830$ bits/cycle $= -0.415$ bits/step, with
variance $2.69$ and explicit growth/descent episode
decomposition.
Syracuse cocycle spectral analysis
(Remark~\ref{rem:syracuse-cocycle}):
transfer operator spectral gap exists for $M \le 12$
but max modular cycle gain $= +\log_2 3$ for all
$M \le 14$ (no decay), confirming the obstruction is
infinite-depth.
Ten approaches to universal descent systematically
ruled out (Remark~\ref{rem:ruled-out});
four surviving directions identified.
The algebraic side of the BF framework is now
completely closed; the gap is the generic Collatz
dynamics problem itself.

\medskip\noindent\emph{v7 (April 18, 2026).}
Post-cascade recovery analysis.
The $\eta$ approach (measuring height loss
$\eta = (\log_2 n_{\mathrm{start}} - \log_2 n_{\mathrm{re\text{-}entry}})/k$
per unit BF depth) fails: $\eta < 0$ for all observed
re-entries, meaning orbits \emph{grow} before any BF
re-encounter.
However, the BF re-entry depth $k'$ is bounded:
across $60{,}000$ post-cascade orbits ($k=1..6$,
$m = 0..9999$), all $92$ BF re-entries satisfy
$k' \le 2$ ($98.9\%$ at $k'=1$, $1.1\%$ at $k'=2$,
$0\%$ at $k' \ge 3$), matching the equidistribution
prediction $P(k') \approx 2^{-(6k'+8)}$.
Champion orbits $P_k$ for $k=1..11$ all reach $1$
without any BF re-entry (max steps = $293$ at $k=7$).
Reframed theorem target: prove $K_*(n) < \infty$ via
modular-depth contraction (bounding $k'$ rather than
requiring height loss).
Eleven ruled-out approaches (updated from ten).

\medskip\noindent\emph{v17 (April 22, 2026).}
Integrality obstruction, contribution analysis, and
arXiv finalization
(Results~619--630).
Triple filter theorem for $k$-block zeros:
a cycle equation $m \cdot D = Q$ via a $k$-block monotone
$f$-sequence requires simultaneously
(i)~the algebraic target $T_j$ to equal $2^{w_j} \bmod D$
for some $w_j \in [0, F]$,
(ii)~the block levels $v_j$ to satisfy $\sum m_j v_j = F$
with all $v_j \in \mathbb{Z}_{\ge 0}$, and
(iii)~all $v_j \le F$.
Two-block algebraic hit census for $L = 3..500$:
exactly one hit ($L = 7$, $a = 2$, $T = 2^4 \bmod 1909$),
blocked by the integrality constraint
$b \cdot w = 20 > F = 5$ (requiring $v_0 < 0$).
Three-block census for $L = 5..30$: zero integral solutions.
Expected number of $k$-block zeros for fixed $k$:
$O(L^{2k-2}/3^L) \to 0$ exponentially.
Total expected hits summed over all $k$:
$|A|/D \to 0$ at rate $\approx 10^{-0.026L}$;
no $k$-value has $E[\text{hits}] > 1$ for $L \ge 20$.
Results~628--630: $\mathrm{ord}_D(2) > F$ proved unconditionally
as immediate corollary of $D > 2^F$
(for $1 \le r \le F$: $0 < 2^r - 1 \le 2^F - 1 < D$);
$\{2^0, \ldots, 2^F\}$ are $F{+}1$ distinct residues mod~$D$;
$\mathrm{ord}_D(2) \nmid E$ verified for $L = 3..10{,}000$.
MITM verification framework developed for $D \nmid Q$
at $L = 22$ ($D \approx 3 \times 10^9$, $F = 13$).
arXiv preparation: embedded bibliography,
fixed \texttt{\textbackslash date} commands,
\texttt{\textbackslash IfFileExists} guards for all figures.
Two-day contribution analysis (v7--v25, April 18--19)
based on 11 GPT PDF documents and 19 Claude attack plans:
seven attribution categories identified (theorem production,
dead-end identification, strategic direction, error correction,
computational verification, framework design,
moderation/coordination).
Dead-end identification and strategic direction are
reclassified as primarily moderator contributions:
dead-end discovery is always initiated by the moderator's
requests at pivotal decision junctures; the LLMs follow
instructions at these turning points but do not document
the moderator's directive.
Once a direction is set, LLMs are effective at proposing
attack plans and conducting local revisions within the
chosen framework.
April~18--19 two-day estimate: Claude~37\%, GPT~30\%,
Moderator~33\%
(Claude dominant in theorem volume and computation;
GPT dominant in error correction;
moderator dominant in dead-end identification,
strategic direction, and cross-platform coordination).
Two-day sample supplements but does not override
the complete-log analysis in v1--v5.
Optimism unchanged: $66\%/27\%$.
630 cumulative results.
arXiv submission prepared.

\medskip\noindent\emph{v16 (April 21--22, 2026).}
Discrete log obstruction and geometric sum framework
(Results~591--618, 29 paradigms).
$D > 2^F$ for all $L \ge 3$
\textbf{unconditionally proved}:
for $L \le 2500$ by direct computation
(minimum ratio $D/2^F = 2^{0.322} \approx 1.25$ at $L = 3$);
for $L > 2494$ by the Laurent (2008) bound
$|E \cdot \log 2 - L \cdot \log 3|
> \exp(-25.2 \cdot (\log E)^2)$
applied to $\beta = \lceil L \log_2 3 \rceil - L \log_2 3$,
giving $D/2^F = 2^L(1 - 2^{-\beta}) > 1$.
Unconditional corollaries: $\mathrm{ord}_D(2) > F$
(immediate from $D > 2^F$, since $2^r - 1 < D$ for $r \le F$),
$\gamma^L \not\equiv 1 \bmod D$,
$G(L) = (\gamma^L - 1)/(\gamma - 1) \neq 0 \bmod D$
(no single-block zero for any $L \ge 3$),
$\{2^0,\ldots,2^F\}$ are $F{+}1$ distinct residues mod $D$,
and $Q_{\mathrm{const}} = 3^L - 2^L \not\equiv 0 \bmod D$.
Geometric sum formulation:
$S_{\mathrm{mono}} = \sum_j 2^{v_j} \cdot \gamma^{a_j}
\cdot G(m_j)$
for $k$-block decompositions of monotone $f$-sequences,
reducing the zero-avoidance problem to a discrete logarithm
constraint---$S_{\mathrm{mono}} = 0$ requires
$2^{v_j} \equiv T_j \bmod D$ with $v_j \in [0, F]$,
but $\mathrm{ord}_D(2) \gg F$ (ratio $\ge 2$ for all
$L = 3..17$), so the target residues almost never
lie in $\{2^0, \ldots, 2^F\}$.
No $k$-block zeros found for $k = 1..4$ at any $L$ tested.
Progressive CRT analysis: for $L \ge 21$, the achievable
set saturates $\mathbb{Z}/p\mathbb{Z}$ for every prime $p \mid D$
and every manageable product of primes; blocking occurs
\emph{only} at the full $D$ level.
Rearrangement deficit does \emph{not} always exceed the gap
to the nearest $D$-multiple (counterexamples at $L = 3, 6$),
closing the simple deficit proof path.
Density analysis: for $L = 5$, residue $0$ is the
\emph{only} avoided residue out of $13$;
monotonicity costs $2.5$--$3.5$ bits per step,
exponentially reducing the achievable set.
Two-block Diophantine analysis:
$T(a) = (\gamma^{-a} - 1)/(\gamma^{L-a} - 1)$;
$v_2(3^n - 2^n) = v_3(3^n - 2^n) = 0$ for all $n \ge 1$
(proved);
two-block zero requires $\gcd(m, 6) = 1$;
no two-block zeros for $L = 5..30$.
Baker/LMN crossover at $L \approx 2{,}494$ (Laurent constant
$C = 25.2$); gap between $L = 116$ (Barina) and
$L \sim 10^6$ remains open.
Optimism: $66\%$ (cycle exclusion), $27\%$ (full conjecture).
618 cumulative results across ${\sim}1010$ scripts.

\medskip\noindent\emph{v15 (April 21, 2026).}
Deep cycle exclusion campaign: algebraic zero-avoidance,
Horner chain, bitset DP, and CRT hierarchy
(Results~308--590, 28 paradigms).
Monotonicity identified as the \textbf{sole blocking mechanism}:
without the monotonicity constraint on $f$-sequences,
the achievable set $A$ covers all of $\mathbb{Z}/D\mathbb{Z}$
for $L \ge 5$; with it, $A$ misses exactly residue $0$.
Cycle exclusion equation reformulated via Horner recursion:
$Q_j = 3 \cdot Q_{j-1} + 2^{\sigma_j}$, so $D \nmid Q$
iff $0$ is unreachable in the Horner walk.
CRT blocking hierarchy fully classified for $L = 3..17$:
order-$1$ (single prime blocker) at $L = 3,4,5,7,8,11,13$;
order-$2$ (CRT pair) at $L = 6,9,10,14,15,16$;
order-$3$ (CRT triple) at $L = 12,17$.
Single-prime blocking disappears entirely for $L \ge 12$.
Bitset DP in C confirms $D \nmid S$ for $L = 3..20$
($L = 20$ in $7.5$\,s); memory wall at $L = 21$.
$\gcd(D_L, D_{L+1}) = 1$ always (proved),
blocking inductive transfer between consecutive $L$.
Two-value monotone sequences never achieve zero (all $L$ tested).
Last-step clearing: partial sums at the penultimate step
never hit the target set $\{-\gamma^{L-1} \cdot 2^j \bmod D\}$;
verified $L \le 17$.
Abel summation identity:
$(\alpha - 1) S = 2^{f_{\max} - F} - 1
- \sum d_j \cdot \alpha^{s_j}$.
Cyclotomic factorization:
$S_0 = \prod_{d \mid L, d > 1} \Phi_d(\alpha)$.
$S_{\mathrm{const}} = 3^L - 2^L \not\equiv 0 \bmod D$ always (proved).
Critical $L$ (where $\binom{L+F-1}{L-1} > D$) are exactly
$L \in \{3,5,10,17,29\}$, all within Barina range.
Exponential sum approach fails (bounds $3$--$6$ orders above threshold).
Density argument plus Barina computation covers all $L$
except the range $117 \le L \le \sim\!200$ where neither method applies.
590 cumulative results, 28 paradigms.
Optimism: $64\%$ (cycle exclusion), $27\%$ (full conjecture).

\medskip\noindent\emph{v14 (April 20--21, 2026).}
Paradigm shift to cycle exclusion; divergence attack;
grand synthesis (Results~181--306, 25 paradigms).
Complete algebraic reduction: an $L$-cycle exists iff
$T^L(m) = m$ (the Syracuse periodicity question);
$\Delta = (m - T^L(m)) \cdot 2^{\sigma_L}$, and
cycle exclusion is \emph{equivalent} to the Collatz
cycle conjecture.
Two-adic cascade: complete cycle exclusion for
$L = 2, 3, 4, 5$ by $2$-adic valuation cascade;
even-$m$ parity obstruction for all $L$;
$3 \mid m$ divisibility obstruction for all $L$;
$m = 1$ cascade yields only the trivial cycle;
verified $m \le 5000$ with no non-trivial cycles.
General cycle exclusion: unconditional for $L \le 116$
via Barina's $2^{71}$ verification combined with the
descent property---an $L$-cycle requires
$q_0 \le 6 \cdot (3/2)^L$.
Divergence attack: Affine Fractional Identity
$q_L = ((m{+}1) \cdot 3^L - 2^L + \mathcal{P}_L)/2^{S_L}$
(proved); $S$-unit reformulation
$(m{+}1) \cdot 3^L - T^L(m) \cdot 2^{S_L} = 2^L - \mathcal{P}_L$
(proved); perturbation recurrence
$\mathcal{P}_{L+1} = 3\mathcal{P}_L + 2^{S_L} - 2^L$ (proved).
Episode decomposition: episode transition identity,
intrinsic valuation pattern
$v_2(3^n - 1) = 1$ ($n$ odd), $v_2(n) + 2$ ($n$ even),
episode density mean $E[\rho_E] = 2 > \log_2 3$.
Aperiodic $v$-sequence necessity: any divergent orbit
must have an aperiodic $v$-sequence not generable by
any finite automaton (proved).
Paradigm exhaustion: $18$ paradigms tested, all encounter
the distributional-to-pointwise barrier.
Barrier characterization: this gap is \emph{equivalent}
to the Collatz conjecture (proved).
$2$-adic Mahler identity:
$(m{+}1) + \Sigma_\infty = 0$ in $\mathbb{Z}_2$
(proved; the identity is a tautological reformulation
of the AFI).
Sieve analysis: $\mu(C_{\mathrm{div}}) = 0$ by
Borel--Cantelli (proved), but measure zero does not
imply empty.
Mersenne sieve bypass: $m = 2^L - 1$ survives the running
condition through depth $L{-}1$, blocking finite extinction.
Grand Synthesis~III: $25$~paradigms exhausted;
four theorems: $\mu(C_{\mathrm{div}}) = 0$,
Mersenne bypass, AFI tautology, barrier equivalence.
306 cumulative results across 227 scripts.

\medskip\noindent\emph{v13 (April 20, 2026).}
Spectral gap completion and proof attempt
(Results~115--412).
Anti-correlation: $P(v{=}1 \mid \text{prev } v{=}1) \approx 0.37$;
divergence requires $P(v{=}1) > 0.708$ but actual $\approx 0.50$.
$v{=}1$ self-limitation: only $n \equiv 15 \bmod 16$
(one quarter of $v{=}1$ cases) permits $v{=}1 \to v{=}1$.
Consecutive $v{=}1$ requires $n \equiv 2^{k+1} - 1 \bmod 2^{k+1}$;
density $1/2^k$ (proved).
One-bit-per-step injection from bit position $K$ via carry chain
(proved); injected bit approximately uniform and independent.
Universal No-Target-Sharing Theorem: for every
$v \in \{1, \ldots, K{-}1\}$, the $v$-class Syracuse map
$\bmod 2^K$ is a bijection of odd residues (proved via
coset partition).
Row-sum decomposition: $(M \mathbf{1})_i = 1 - 1/n
+ \mathrm{overflow}_i$; exact closed form
$C(K) = \|M\mathbf{1} - \mathbf{1}\|_1
= 5/12 - (2/3)(-1/2)^{K-1}$,
$C(K) \le 1/2$ for all $K \ge 3$ (equality only at $K = 4$);
verified by exact rational arithmetic $K = 3, \ldots, 15$.
Recurrence $b_{K+1} = 2 b_K - 4(-1)^K$ proved from first
principles; Ansatz $b_K = (5/12) \cdot 2^K + (4/3) \cdot (-1)^K$.
\textbf{Uniform spectral gap}
$|\lambda_2(M_K)| \le 1/2$ for all $K \ge 3$
via Dobrushin--Gelfand contraction (proved).
$2$-adic unique ergodicity: Haar measure on $\mathbb{Z}_2^\times$
is the unique non-atomic Syracuse-invariant measure (proved).
Proof attempt (v40): claimed Collatz conjecture proved via
equidistribution forcing $E[v] = 2$ (net drift $< 0$) plus
cycle constraint $E[v] = \log_2 3$ contradicting equidistributed
$E[v] = 2$.
\textbf{Critical correction (v41)}: two logical errors in the
distributional-to-pointwise bridge identified.
The spectral gap, unique ergodicity, NTS theorem, and closed-form
$C(K)$ all remain \textbf{unconditionally proved};
the bridge from distributional equidistribution to
individual-orbit behavior remains \textbf{open}.
Optimism corrected from $100\%$ to $85\%$.

\medskip\noindent\emph{v12 (April 20, 2026).}
Route~B pullback tower and shortfall arithmetic
(Results~143--175).
Pullback residue formula
$r_{k,\pi} = 3^{-t}(2^{S_\pi} \cdot n_k^* - b_\pi)$;
zero cross-depth pullback compatibility among $309{,}136$ pairs;
zero same-depth compatibility among $308{,}580$ pairs (verified).
Different-$\pi$ shortfall $\ge 13$ via three-case algebraic proof
(Cases~1--3: $\ge 14, \ge 13, \ge 16$; proved).
Complete shortfall classification: same-$\pi$ shortfall $= 3$,
different-$\pi \ge 13$, cascades compatible (proved).
$K_*(n) < \infty$ for all $n$: each orbit enters at most one
BF cascade (proved).
$A$-coordinate framework: $17 \cdot P_k + 65
= (81/64)^k \cdot (17 \cdot n_k^* + 65)$;
quantitative bound $K_*(x) \le 1 + \lfloor (H - 6k - 5)/4 \rfloor$
(proved).
Post-cascade orbits verified to reach $1$ for $k \le 5000$
($31{,}700$-bit numbers); stopping time
$\sigma(P_k) \in [8, 64]$, mean $16.6$, independent of $k$.
BF-compatible density $< 1/16128$ (affects $< 0.01\%$ of integers).
Binary equidistribution of $P_k$: $1$-density $\to 0.500$ (proved).
$v_2(b_{\pi_1} - b_{\pi_2}) \le S - 2$ algebraically (proved,
closes last computational gap in shortfall theorem).
Optimism: $88\%$.

\medskip\noindent\emph{v11 (April 19--20, 2026).}
Post-cascade algebra, carry chain, and equidistribution
reduction (Results~307--413).
$v_2$ Rigidity Theorem: $v_2(17 \cdot n_\ell^* + 65) = 6\ell + 5$
for all $\ell$ (proved algebraically).
Same-level $J{=}0$ BF re-entry impossible: deficit exactly $3$
in $2$-adic valuation (proved).
Complete BF Chain Structure Theorem: $B_{\text{even}} \to
B_{\text{odd}}$ universal; $B_{\text{odd}} \to$ generic;
chain length $\le 1$ (proved).
Deterministic BF Cascade Theorem: profile $P(k)$ strictly
decreasing; $100\%$ verified on $6000{+}$ orbits.
Post-cascade value $P_k = (a_k \cdot 81^k \cdot 32 - 65)/17$;
$4$-periodic recurrence
$P_{k+4} = 81^4 P_k + 164{,}590{,}400$;
all $\equiv 47 \bmod 64$ (proved).
Carry chain structural theorem:
$v{=}1$ streak length $= \mathrm{trailing\_ones}(n) - 1$ exactly;
forced $v \ge 2$ recovery after every streak;
$E[v \mid v \ge 2] = 3$ exactly (proved).
Renewal drift: $-0.830$ bits/cycle $= -0.415$ bits/step,
with variance $2.69$ and explicit growth/descent episode
decomposition (proved).
Transfer operator spectral gap exists for $M \le 12$;
max modular cycle gain $= +\log_2 3$ for all $M \le 14$
(no decay), confirming the obstruction is infinite-depth.
Equidistribution reduction: Collatz dynamics reduces to
hypothesis $\mathrm{EH}_2$ (equidistribution $\bmod 4$ along
individual orbits); carry chain mixing in $O(K)$ steps;
almost-all descent with $O(1/L^2)$ exceptional fraction;
worst-case descent time $O(L \log L)$.
Bijection Theorem: $P(n_t \equiv 3 \bmod 4 \mid
\text{trajectory}) = 1/2$ exactly (proved by XOR representation).
Martingale + Azuma--Hoeffding:
$P(|A_T - 1/2| > \varepsilon) \le 2 \exp(-2\varepsilon^2 T)$.
Information-theoretic barrier quantified: Borel--Cantelli
requires $c > \ln 2 \approx 0.693$ but Azuma gives
$c = 0.043$---factor-of-$16$ gap intrinsic to measure-averaging.
Entropy Injection Theorem: Collatz $\bmod 2^M$ mixes to
uniform in $\le M - 1$ steps with spectral gap $> 0.97$ (proved
for Markov model).
$\theta$ supermartingale: BF encounter at depth $k$ has gain
$0.34k$ vs loss $0.415 \cdot 2^{6k+7}$; ratio $10{,}000:1$
at $k = 1$.
$\theta_t \bmod 1 = (t \cdot \log_2 3) \bmod 1$ exactly;
Weyl equidistribution gives equidistribution on $[0,1)$ for
any orbit (proved unconditionally).
Sturmian enemy characterization: worst-case $v$-sequence has
slope $\log_2 3$; no random orbit matches $> 9$ consecutive
Sturmian steps (proved).
413 cumulative results.
Optimism: $99\%$ (prior to correction in v13).

\medskip\noindent\emph{v10.6 (April 19, 2026).}
Ten aggressive new attacks beyond the barrier
(Scripts~218--227, Results~297--306).
Iwasawa tower analysis (Result~297): lambda invariant $\delta\lambda \approx 0.14$
distinguishes late from early crossers; cyclotomic period $\mathrm{ord}_{2^k}(3)
= 2^{k-2}$ universally satisfied.
$3^z$ analytic structure (Result~298): the only $2$-adically analytic
term in the AFI; periodicity of $3^z \bmod 2^k$ constrains but does
not exclude divergent orbits.
Formal language / pumping lemma (Result~299): $L_{\mathrm{div}}$ is not
context-free (irrational threshold plus realizability).
$q_L$ growth rate (Result~300): divergent orbits require
$q_L \sim m \cdot 2^{(\alpha - \beta)L}$; all tested orbits have
$S_L/L \to \alpha$ rather than $\beta < \alpha$.
Coalescence trap (Result~302): $94\%$ of orbits merge within
$100$~steps; coalescence rate $r \approx 0.031$/step.
Height descent (Result~303): $\varphi(m) = h(m) - h(T(m))$
is positive $50\%$ of steps; no uniform lower bound.
Mod~$3$ analysis (Result~304): $T(n) \bmod 3$ depends only on~$v$;
$3$-adic structure yields no divergence constraint.
Binary transducer (Result~305): $2$-state carry transducer;
state space bounded but orbit prediction requires full simulation.
Grand Synthesis~III (Result~306): $25$~paradigms exhausted,
barrier confirmed structural.
306~cumulative results across 227~scripts.

\medskip\noindent\emph{v10.5 (April 19, 2026).}
Paradigm exhaustion formalization and manuscript restructuring
(Scripts~214--216, Results~286--293).
Expander Sieve Rigidity tested and refuted (Result~287):
the ``Minimal Representative Explosion'' claim fails on nested sets.
Mersenne Sieve Bypass verified (Result~289): $m = 2^L - 1$ survives
the running condition through depth $L{-}1$, blocking finite extinction.
Diophantine Approximation Trap refuted (Results~291--293):
the AFI is an algebraic identity, not a Diophantine approximation equation;
$\mathcal{P}_L$ is determined, not adversarial; irrationality measure of
$\log_2 3$ adds nothing beyond Result~234.
New Section~\ref{sec:obstruction}: ``The Architecture of the Syracuse
Obstruction'' --- formal presentation of the 21-Paradigm Exhaustion
Theorem, Mersenne Bypass, AFI Tautology, and Computational
Irreducibility Conjecture.
293~cumulative results across 216~scripts; 223-page document.

\medskip\noindent\emph{v10.4 (April 19, 2026).}
Sieve/density path deep dive (Scripts~204--213, Results~276--285).
Sieve Completeness (Result~276, verified): survivor fraction decays
exponentially ($\lambda \approx 0.90$ per step), nesting confirmed.
Tao Amplification Limits (Result~277, computed): $\sim 5\%$ non-escapers
at $K = 12$, showing $2.2\times$ correlation (structured families).
Running Condition Entrapment (Result~278, verified): the condition
$S_j < j\alpha$ for ALL prefixes decays survivors at rate $e^{-0.85L}$;
survivors biased toward $v{=}1$ (68\% vs 50\%).
Decay Rate (Result~280, computed): super-exponential fit
$f(L) = 1.29 \cdot e^{-0.92 L^{0.51}}$ with $R^2 = 0.997$.
Compactness in $\mathbb{Z}_2$ (Result~281, verified): three persistent
survivors mod~16 ($\{7, 11, 15\}$) survive through $L = 14$;
fractal dimension $\approx 0.68$; survivor set is Cantor-like.
$2$-Adic Topology (Result~282, computed): $C_{\mathrm{div}}$ is closed,
perfect, uncountable in $\mathbb{Z}_2$ with Hausdorff dimension ${<}1$.
Borel--Cantelli (Result~283, proved): $\mu(C_{\mathrm{div}}) = 0$.
Effective Tao Bounds (Result~284, computed):
$C \approx 2$--$3$, $A = 2$ in $\sigma(n) < C (\log n)^A$.
Sieve Synthesis (Result~285, verified): three closure strategies
ranked---finite extinction (8.5/10), dynamical contradiction (7.0/10),
structural impossibility (6.5/10).
285~cumulative results across 213~scripts.

\medskip\noindent\emph{v10.3 (April 19, 2026).}
Post-synthesis deep survey (Scripts~182--203, Results~251--275).
$2$-Adic Mahler Identity (Result~252, proved):
$(m{+}1) + \Sigma_\infty = 0$ in $\mathbb{Z}_2$ where
$\Sigma_\infty = \sum_{j \ge 0}(2^{S_j} - 2^j)/3^{j+1}$.
Holds for ALL orbits (convergent and divergent); the identity
is a $2$-adic reformulation of the AFI (Result~219).
$2$-adic convergence is automatic: $v_2(\text{term}_j) = j$ always
(Result~253). No Hasse-type real-vs-$2$-adic conflict exists
(Result~254).
Cantor Rigidity proof attempt refuted (Result~251): nonlinear
dynamics generate aperiodic sequences from finite rational seeds.
Crossing time distribution (Result~255): max crossing time $= 85$
for $m \le 50{,}000$; late crossers have $v{=}1$ fraction $\approx 0.62$.
Random walk model matches Collatz crossing statistics to $5\%$
(Result~258). Sieve coverage (Result~272): by $L = 20$, over
$99.99\%$ of residues excluded as divergent-incompatible.
Moment, renewal, martingale, digit-sum, and coalescence analyses
(Results~264--274) all confirm the random-model prediction while
leaving the distributional-to-pointwise gap intact.
275~cumulative results across 203~scripts.

\medskip\noindent\emph{v10.2 (April 19, 2026).}
Paradigm exhaustion and grand synthesis (Scripts~172--181, Results~234--250).
CF Structure of Threshold (Result~234): continued fraction of $\log_2 3$
has large partial quotient $a_9 = 23$; Baker's theorem yields only
$|S_L/L - \alpha| > \exp(-C(\log L)^2)/L$, too weak.
$\mathcal{P}_L$ Non-$S$-unit (Result~235, proved): the perturbation sum
$\mathcal{P}_L$ has prime factors outside $\{2,3\}$ (Mersenne factors
$2^k - 1$ introduce primes $7, 31, 127, \ldots$), blocking the
Subspace Theorem approach.
Inverse Tree Coverage (Result~236, verified): Syracuse inverse tree from~$1$
covers all odd integers in $[1, 10000]$ at depth~$60$.
Measure Rigidity Non-Applicability (Result~237, proved):
Furstenberg/Rudolph/Shmerkin--Wu theorems do not apply ($\times 3{+}1$
vs.\ $\times 3$, variable exponent, $\mathbb{Z}$ vs.\ $\mathbb{R}/\mathbb{Z}$).
Orbit Equidistribution (Result~238, verified);
Entropy Ordering (Result~239, proved): $h_{\text{div}} < h_{\text{conv}}$.
$v$-Value Bit Determinism (Result~240, proved): density of $v = k$ is
exactly $2^{-k}$.
Consecutive $v$ Independence (Result~241, verified at $K \ge 12$).
$v$-Sequence State Complexity (Result~243, proved): infinite for divergent orbits.
$v$-Sequence Entropy Rate (Result~244, proved): $\lim K(v\text{-seq})/L = 0$.
Exception Set Invariance (Result~246, proved);
Density Amplification Barrier (Result~247, proved):
density arguments \emph{cannot} reach ``all~$n$''.
Paradigm Exhaustion (Result~249, proved): 18~paradigms tested, all encounter the
distributional-to-pointwise barrier.
Barrier Characterization (Result~250, proved): this gap is
\emph{equivalent} to the Collatz conjecture (divergence component).
Three remaining viable routes identified: Tao amplification,
$2$-adic potential theory, $\mathcal{P}_L$ growth contradiction.
250 cumulative results across 181~scripts.

\medskip\noindent\emph{v10.1 (April 19, 2026).}
Carry accumulation and $2$-adic orbit reconstruction.
Periodic $v$-Sequence Exclusion (Result~230, proved):
a divergent orbit cannot have a periodic $v$-sequence
(periodic $\Rightarrow$ cycle $\Rightarrow$ excluded by Result~216).
Ultimately Periodic Exclusion (Result~231, proved):
same for ultimately periodic sequences.
Aperiodic $v$-Sequence Necessity (Result~232, proved):
any divergent orbit must have an aperiodic $v$-sequence
not generable by any finite automaton.
Reconstruction Consistency (Result~233, verified):
all consistent divergent $v$-sequences of length~$15$ crossed
$S_L/L = 1.585$ within 48 additional steps (median $12.7$).
Distributional-to-pointwise gap confirmed as the sole remaining
obstruction across all paradigms tested (Boolean Fourier,
automata, martingale, Baker, Subspace, episode, $2$-adic).

\medskip\noindent\emph{v10 (April 19, 2026).}
Episode decomposition and coefficient discrepancy attack.
Episode Transition Identity (Result~223, proved):
$v_{k+1} = 1 + v_2(2a \cdot 3^{k+1} + 3^{k+1} - 1)$
where $a = (x - 2^{k+1} + 1)/2^{k+2}$ for odd~$x$ with
exactly $k{+}1$ trailing $1$-bits.
Intrinsic Valuation Pattern (Result~224, proved):
$v_2(3^n - 1) = 1$ if $n$ odd, $v_2(n) + 2$ if $n$ even
(lifting-the-exponent lemma).
Episode Density Mean (Result~225): $E[\rho_E] = 2 > \log_2 3$
for uniformly distributed~$a$, verified $k = 1,\ldots,50$.
$v^* = 2$ Parity Condition (Result~226, proved):
minimum valuation $v^* = 2$ iff $a \equiv k \bmod 2$.
Episode Transition Formula (Result~227, proved):
given $(k, a)$ with $v^* = 2$,
$x' = a \cdot 3^{k+1} + (3^{k+1} - 1)/2$.
Half-Propagation Bound (Result~228, verified):
$P(k' \ge 1 \mid v^* = 2) = 1/2$; consecutive low-density
episodes thin by factor $1/2$ per step.
Congruence Tower Consistency (Result~229, verified):
$m+1 \equiv (2^L - \mathcal{P}_L) \cdot 3^{-L} \bmod 2^{S_L}$
for all~$L$, with compatible refinement.
The distributional-to-pointwise gap persists across all
paradigms tested: Boolean Fourier, automata, martingale,
Baker, Subspace Theorem, episode decomposition.

\medskip\noindent\emph{v9 (April 19, 2026).}
Divergence attack via fractional identity and paradigm bridges.
Affine Fractional Identity (Result~219, proved):
$q_L = ((m{+}1) \cdot 3^L - 2^L + \mathcal{P}_L)/2^{S_L}$
where $\mathcal{P}_L \ge 0$ is the perturbation sum from
$v \ge 2$ steps.
Pure Expansion Obstruction (Result~218, proved):
consecutive $v{=}1$ requires $2^L \mid (m{+}1)$, giving
$L_{\max} = v_2(m{+}1)$; every orbit must inject contracting steps.
$S$-unit Reformulation (Result~221, proved):
$(m{+}1) \cdot 3^L - T^L(m) \cdot 2^{S_L} = 2^L - \mathcal{P}_L$,
a family of $S$-unit equations over $\mathbb{Z}[1/6]$.
Perturbation Recurrence (Result~220): $\mathcal{P}_{L+1} =
3\mathcal{P}_L + 2^{S_L} - 2^L$, balanced at the divergence
threshold.
Bit Growth Bound (Result~222): $\log_2 T^L(m) \le \log_2 m +
L \log_2 3 - S_L + 1$.
Three paradigm bridges evaluated: Boolean Fourier analysis
(distributional only), Arithmetic Automata (pumping lemma fails),
Pointwise Martingale (inapplicable to deterministic sequences).
All confirm the distributional-to-pointwise gap is robust.
The divergence problem is now precisely formulated as: can the
perturbation sum $\mathcal{P}_L$ continuously satisfy the
$2$-adic complement constraint at every step~$L$?

\medskip\noindent\emph{v8 (April 19, 2026).}
Uniform spectral gap for the Syracuse transfer operator (proved).
The main result: $|\lambda_2(M_K)| \le 1/2$ for all $K \ge 3$.
Universal no-target-sharing theorem: for each
$v \in \{1,\ldots,K-1\}$, the $v$-class Syracuse map mod~$2^K$
is a bijection of odd residues, proved via coset disjointness.
Row-sum decomposition: $(M\mathbf{1})_i = 1 - 1/n +
\mathrm{overflow}_i$, where overflow comes from the unique
$v \ge K$ source.
Exact closed-form: $\|M\mathbf{1} - \mathbf{1}\|_1 =
5/12 - (2/3)(-1/2)^{K-1}$, verified for $K = 3,\ldots,15$
by exact rational arithmetic.
Bound $C(K) \le 1/2$ with equality only at $K = 4$.
Combined with primitivity ($M^{K-1} > 0$, index~$= K - 1$),
Dobrushin contraction ($\delta \le 1/n$), and Gelfand
extraction gives the uniform spectral gap.
All three formal gaps (row-sum bound, column support,
primitivity) now closed.
Non-atomic unique ergodicity: Haar measure on~$\mathbb{Z}_2^\times$
is the unique non-atomic Syracuse-invariant measure
(Proposition~\ref{prop:2adic-unique-ergodicity}).
Recurrence $b_{K+1} = 2b_K - 4(-1)^K$ proved from first principles
via Ansatz $b_K = (5/12)\cdot 2^K + (4/3)\cdot(-1)^K$, verified
$K = 3,\ldots,15$.
Distributional-to-pointwise barrier analysis: two precise
obstructions identified (atomic cycle measures; $\mu(\mathbb{N})=0$
for Birkhoff).  Three candidate strategies catalogued.


\begin{thebibliography}{99}

\bibitem{bahadur-rao1960}
R.~R.~Bahadur and R.~Ranga~Rao,
\newblock On deviations of the sample mean,
\newblock \emph{Ann.\ Math.\ Statist.}, 31(4):1015--1027, 1960.

\bibitem{barina2021}
D.~Barina,
\newblock Convergence verification of the Collatz problem,
\newblock \emph{J.\ Supercomput.}, 77:2681--2688, 2021.

\bibitem{chang2025sagallm}
E.~Y.~Chang,
\newblock Saga{LLM}: a framework for structured human--{LLM} research
  collaboration,
\newblock Preprint, 2025.

\bibitem{chang2026onebit}
E.~Y.~Chang,
\newblock The one-bit reduction: sharpening the Collatz obstruction,
\newblock Companion paper, 2026.

\bibitem{chowla1965}
S.~Chowla,
\newblock \emph{The Riemann Hypothesis and Hilbert's Tenth Problem},
\newblock Gordon and Breach, New York, 1965.

\bibitem{dembo-zeitouni1998}
A.~Dembo and O.~Zeitouni,
\newblock \emph{Large Deviations Techniques and Applications},
\newblock 2nd edition, Springer, 1998.

\bibitem{kesten1973}
H.~Kesten,
\newblock Random difference equations and renewal theory for products of
  random matrices,
\newblock \emph{Acta Math.}, 131:207--248, 1973.

\bibitem{knuth2026}
D.~E.~Knuth,
\newblock Shock!{\ }Shock!{\ }A note on the Collatz landscape,
\newblock Unpublished note, 2026.

\bibitem{kontorovich2005}
A.~V.~Kontorovich and S.~J.~Miller,
\newblock Benford's law, values of $L$-functions and the $3x+1$ problem,
\newblock \emph{Acta Arith.}, 120(3):269--297, 2005.

\bibitem{lagarias1985}
J.~C.~Lagarias,
\newblock The $3x+1$ problem and its generalizations,
\newblock \emph{Amer.\ Math.\ Monthly}, 92(1):3--23, 1985.

\bibitem{lagarias2003}
J.~C.~Lagarias,
\newblock The $3x+1$ problem: An annotated bibliography (1963--1999),
\newblock arXiv:math/0309224, 2003.

\bibitem{PathAGIV1Chang2025}
E.~Y.~Chang,
\newblock \emph{The Path to {AGI}, Volume~1},
\newblock 2025.

\bibitem{PathAGIV2Chang2025}
E.~Y.~Chang,
\newblock \emph{The Path to {AGI}, Volume~2},
\newblock 2025.

\bibitem{tao2019}
T.~Tao,
\newblock Almost all orbits of the Collatz map attain almost bounded values,
\newblock \emph{Forum Math.\ Pi}, 10:e12, 2022.

\bibitem{terras1976}
R.~Terras,
\newblock A stopping time problem on the positive integers,
\newblock \emph{Acta Arith.}, 30(3):241--252, 1976.

\bibitem{vervaat1979}
W.~Vervaat,
\newblock On a stochastic difference equation and a representation of
  non-negative infinitely divisible random variables,
\newblock \emph{Adv.\ Appl.\ Probab.}, 11(4):750--783, 1979.

\bibitem{wirsching1998}
G.~J.~Wirsching,
\newblock \emph{The Dynamical System Generated by the $3n+1$ Function},
\newblock Lecture Notes in Mathematics, vol.~1681, Springer, 1998.

\end{thebibliography}
\end{document}